\documentclass[10pt, reqno]{amsart}

\makeatletter
\def\@tocline#1#2#3#4#5#6#7{\relax
  \ifnum #1>\c@tocdepth
  \else
    \par \addpenalty\@secpenalty\addvspace{#2}
    \begingroup \hyphenpenalty\@M
    \@ifempty{#4}{
      \@tempdima\csname r@tocindent\number#1\endcsname\relax
    }{
      \@tempdima#4\relax
    }
    \parindent\z@ \leftskip#3\relax \advance\leftskip\@tempdima\relax
    \rightskip\@pnumwidth plus4em \parfillskip-\@pnumwidth
    #5\leavevmode\hskip-\@tempdima
      \ifcase #1
       \or\or \hskip 1em \or \hskip 2em \else \hskip 3em \fi
      #6\nobreak\relax
    \dotfill\hbox to\@pnumwidth{\@tocpagenum{#7}}\par
    \nobreak
    \endgroup
  \fi}
\makeatother

\usepackage[latin1]{inputenc}
\usepackage{amsmath,amsthm,amssymb,color,mathtools,esint,comment}
\definecolor{darkgreen}{rgb}{0.0, 0.6, 0.13}
\usepackage{graphicx}
\usepackage{caption}
\usepackage[shortlabels]{enumitem}
\setlist{  
  listparindent=\parindent,
  parsep=0pt,
}
\usepackage{xcolor}
\usepackage{makecell}
\usepackage{array}

\usepackage{hyperref}

\usepackage{float}

\usepackage{mathrsfs}

\extrafloats{500}
\usepackage{import}  

\usepackage{geometry}
\usepackage{bbm}

\usepackage{longtable}
\usepackage{autobreak}

\usepackage{bm}

\usepackage[normalem]{ulem}

\newcommand{\newterm}[1]{{\textbf{#1}}}

\newcommand{\pb}{{\boldsymbol{p}}}

\newcommand{\supp}{\mathrm{supp}\,}

\newcommand{\Dirac}{\boldsymbol{\Delta}}
\newcommand{\dirac}{\boldsymbol{\delta}}

\newcommand{\Cb}{\mathbb C}

\newcommand{\Eb}{\mathbb E}

\newcommand{\Kb}{\mathbb K}

\newcommand{\Mb}{\mathbb M}
\newcommand{\Nb}{\mathbb N}

\newcommand{\Pb}{\mathbb P}
\newcommand{\Qb}{\mathbb Q}
\newcommand{\Rb}{\mathbb R}
\newcommand{\Sb}{\mathbb S}
\newcommand{\Tb}{\mathbb T}

\newcommand{\Vb}{\mathbb V}

\newcommand{\Xb}{\mathbb X}
\newcommand{\Yb}{\mathbb Y}
\newcommand{\Zb}{\mathbb Z}

\newcommand{\oneb}{\mathbbm{1}}

\newcommand{\Ac}{\mathcal A}
\newcommand{\Bc}{\mathcal B}
\newcommand{\Cc}{\mathcal C}
\newcommand{\Dc}{\mathcal D}
\newcommand{\Ec}{\mathcal E}
\newcommand{\Fc}{\mathcal F}
\newcommand{\Gc}{\mathcal G}
\newcommand{\Hc}{\mathcal H}
\newcommand{\Ic}{\mathcal I}

\newcommand{\Lc}{\mathcal L}
\renewcommand{\Mc}{\mathcal M}
\newcommand{\Nc}{\mathcal N}

\newcommand{\Pc}{\mathcal P}
\newcommand{\Qc}{\mathcal Q}
\newcommand{\Rc}{\mathcal R}
\newcommand{\Sc}{\mathcal S}
\newcommand{\Tc}{\mathcal T}

\newcommand{\Vc}{\mathcal V}

\newcommand{\Xc}{\mathcal X}
\newcommand{\Yc}{\mathcal Y}
\newcommand{\Zc}{\mathcal Z}

\newcommand{\Sf}{\mathfrak S}
\newcommand{\Bf}{\mathfrak B}

\newcommand{\Rf}{\mathfrak R}

\newcommand{\mf}{\mathfrak m}
\newcommand{\pf}{\mathfrak p}

\newcommand{\rf}{\mathfrak r}
\newcommand{\of}{\mathfrak o}
\newcommand{\Af}{\mathfrak A}
\newcommand{\lf}{\mathfrak l}
\newcommand{\tf}{\mathfrak t}
\newcommand{\qf}{\mathfrak q}

\newcommand{\Lf}{\mathfrak L}

\newcommand{\Cf}{\mathfrak C}

\newcommand{\Mf}{\mathfrak M}

\newcommand{\nf}{\mathfrak n}

\newcommand{\Fs}{\mathscr F}
\newcommand{\Ts}{\mathscr T}

\newcommand{\vz}{\boldsymbol{z}}

\newcommand{\vLambda}{\boldsymbol{\Lambda}}

\newcommand{\vt}{\boldsymbol{t}}

\theoremstyle{definition}
\newtheorem{mainthm}{Theorem}
\newtheorem{theorem}{Theorem}[section]

\newtheorem{lemma}[theorem]{Lemma}
\newtheorem{proposition}[theorem]{Proposition}

\newtheorem{definition}[theorem]{Definition}
\newtheorem{reduct}{Reduction}

\theoremstyle{remark}
\newtheorem{remark}[theorem]{Remark}
\numberwithin{equation}{section}

\begin{document}

\title{Long time derivation of the Boltzmann equation from hard sphere dynamics}

\author{Yu Deng}
\address{\textsc{Department of Mathematics, University of Chicago, Chicago, IL, USA}}
\email{\texttt{yudeng@usc.edu}}
\author{Zaher Hani}
\address{\textsc{Department of Mathematics, University of Michigan, Ann Arbor, MI, USA}}
\email{\texttt{zhani@umich.edu}}
\author{Xiao Ma}
\address{\textsc{Department of Mathematics, University of Michigan, Ann Arbor, MI, USA}}
\email{\texttt{mxiao@umich.edu}}
\date{}

\begin{abstract}
We provide a rigorous derivation of Boltzmann's kinetic equation from the hard sphere system for rarefied gas, which is valid for arbitrarily long times, as long as the solution to the Boltzmann equation exists. This extends Lanford's landmark theorem \cite{Lan75}, which justifies this derivation for a sufficiently short time. In a companion paper \cite{DHM25}, we connect this derivation to existing literature on hydrodynamic limits. This completes the resolution of Hilbert's Sixth Problem pertaining to the derivation of fluid equations from Newton's laws, in the case of a rarefied, hard sphere gas.

The general strategy follows the paradigm introduced by the first two authors for the long-time derivation of the wave kinetic equation in wave turbulence theory. This is based on propagating a long-time cumulant ansatz, which keeps memory of the full collision history of the relevant particles, by an important partial time expansion. The heart of the matter is proving the smallness of these cumulants in $L^1$, which can be reduced to  combinatorial properties for the associated diagrams which we call \emph{molecules}. These properties are then proved by devising an elaborate \emph{cutting algorithm}, which is a major novelty of this work.

\end{abstract}
\maketitle
\tableofcontents

\section{Introduction}\label{sec.intro}

Dating back to the year 1872, Boltzmann's kinetic theory is one of the most profound and revolutionary scientific ideas ever introduced. At the heart of this theory lies Boltzmann's kinetic equation which gives an effective \emph{macroscopic} description of the statistics of a large system of \emph{microscopic} interacting particles. A significant feature of this theory, which caused some controversy at the time of its inception, is the fact that the macroscopic description is given by a \emph{time-irreversible} equation (i.e. with an entropy functional that increases forward in time), which is derived from a microscopic system satisfying the \emph{time-reversible} Newton's laws of motion. This is known as the emergence of the \emph{arrow of time}.

Unsurprisingly, this profound physical theory prompted the deep mathematical question, namely to rigorously justify the derivation of Boltzmann's equation starting from Newton's laws. This question was asked by Hilbert, in his address at the International Congress of Mathematicians in 1900, as the \emph{sixth problem}\footnote{This was one of the two questions proposed in Hilbert's follow-up article \cite{H01} expanding on his ICM address. The other question was to establish the mathematical foundation of probability theory, which was done in the 1930--40s.} in his celebrated  list of outstanding mathematical problems for the twentieth century. More precisely, he proposed the question of giving a mathematically rigorous derivation for the macroscopic equations of fluids, such as the Euler and Navier-Stokes equations, from the Newtonian laws governing microscopic particle interactions, with the derivation of Boltzmann's kinetic equation as a crucial intermediate step. This problem requires justifying two limits: (i) the \emph{kinetic limit} where we pass from the (microscopic) Newtonian dynamics for an $N$-particle system to Boltzmann's kinetic equation for the one-particle density function in the limit $N\to \infty$, and (ii) the \emph{hydrodynamic limit} where we pass from Boltzmann's kinetic equation to the (macroscopic) equations of fluid motion in the limit where the collision rate goes to infinity. 

Despite being the subject of extensive study in the last century, this problem remained open for more than 120 years after it was first announced. While the second hydrodynamic limit is much better understood (see the discussions in Section \ref{sec.future_2} below), the main difficulty in the resolution of Hilbert's sixth problem lies in justifying the first kinetic limit, which is the subject of this manuscript. In this direction, the state of art can be briefly summarized as follows, see Section \ref{sec.past} for a more extensive review. A landmark result of Lanford \cite{Lan75} in 1975 provides the rigorous derivation of the Boltzmann equation for \emph{sufficiently short time}. This result has since been completed, improved, and further developed by different groups of authors (see for example Illner-Pulvirenti \cite{IP86}, Gallagher-Saint Raymond-Texier \cite{GST14} and more recent works \cite{PS17,BGSS20,BGSS20_2,BGSS22,BGSS22_2} concerning fluctuations). However, all these works are still restricted to perturbative regimes, either under short-time or near-vacuum conditions or at the exact Maxwellian equilibrium.

The main result of this paper fills this gap, and thus completes the justification of the first limit needed for the resolution of Hilbert's sixth problem. More precisely, we provide a rigorous derivation of Boltzmann's equation from hard sphere dynamics, for arbitrarily long time that covers the full lifespan of the Boltzmann solution. This solution may come from any given initial data, which can also be arbitrarily large. In particular, if a solution to the Boltzmann equation exists globally, then our derivation is valid for all finite times, and times going (slowly) to infinity. In the companion work \cite{DHM25}, we will connect this to the second hydrodynamic limit, thus giving a full resolution to Hilbert's sixth problem (see Section \ref{sec.future_2}).

\subsection{Setup and the main result}\label{sec.setup} In this section we state our main result, Theorem \ref{th.main}. We start by describing the microscopic dynamics of the interacting particles, the kinetic formalism involving $s$-particle density functions, and the Boltzmann equation governing the effective dynamics of those functions. The microscopic system here is given by the standard hard sphere particle model, under the Boltzmann-Grad scaling law. We shall discuss some other scenarios in Section \ref{sec.future_1}.
\subsubsection{The hard sphere dynamics}
Consider $N$ hard sphere particles with diameter $\varepsilon$.
\begin{definition}[The hard sphere dynamics \cite{Alexander}]\label{def.hard_sphere} We define the dynamics for the $N$-particle hard sphere system with diameter $\varepsilon$ in dimension $d\geq 2$. This will be referred to as the \textbf{Original} or \textbf{O-dynamics} when necessary, to distinguish it with the various modified dynamics in Definitions \ref{def.modified_dynamics}, \ref{def.molecule_truncated_dynamics}, \ref{def.E_dynamics} and \ref{def.T_dynamics} below.
\begin{enumerate}
\item\label {it.O1}\emph{State vectors $z_j$ and $\vz_N$}. We define $z_j=(x_j,v_j)$ and $\vz_N=(z_j)_{1\leq j\leq N}$, where $x_j$ and $v_j$ are the center of mass and velocity of particle $j$. These $z_j$ and $\vz_N$ are called the \textbf{state vector} of the $j$-th particle and the collection of all $N$ particles, respectively.
\item\label{it.O2} \emph{The domain $\Dc_N$}. We define the \textbf{non-overlapping domain} $\Dc_N$ as
\begin{equation}\label{eq.setDN}\Dc_N:=\big\{\vz_N=(z_1,\cdots,z_N)\in \Rb^{2dN}:|x_i-x_j|\geq\varepsilon\,\,(\forall i\neq j)\big\}.
\end{equation} 
\item\label{it.O3} \emph{The hard sphere dynamics $\vz_N(t)$}. Given initial configuration $\vz_N^0=(z_1^0,\cdots,z_N^0)\in\Dc_N$, we define the \textbf{hard sphere dynamics} $\vz_N(t)=(z_j(t))_{1\leq j\leq N}$ as the time evolution of the following system:
\begin{enumerate}
\item\label{it.O31} We have $\vz_N(0)=\vz_N^0$.
\item\label{it.O32} Suppose $\vz_N(t')$ is known for $t'\in[0,t]$. If $|x_i(t)-x_j(t)|=\varepsilon$ for some $(i,j)$, then we have
\begin{equation}\label{eq.hardsphere}
\left\{
\begin{aligned}x_i(t^+)&=x_i(t),\quad x_j(t^+)=x_j(t),\\
v_i(t^+)&=v_i(t)-\big((v_i(t)-v_j(t))\cdot\omega\big)\omega,\quad \omega:=(x_i(t)-x_j(t))/\varepsilon\in \Sb^{d-1};\\
v_j(t^+)&=v_j(t)+\big((v_i(t)-v_j(t))\cdot\omega\big)\omega.
\end{aligned}
\right.
\end{equation} where $t^+$ indicates right limit at time $t$, see {\color{blue}Figure \ref{fig.1}}. Note that $x_j(t)$ are always continuous in $t$; in this definition (and similarly for all the modified dynamics below) we also always require $v_j(t)$ to be left continuous in $t$, so $v_j(t)=v_j(t^-)$.
\begin{figure}[h!]
\centering
\includegraphics[width=0.18\linewidth]{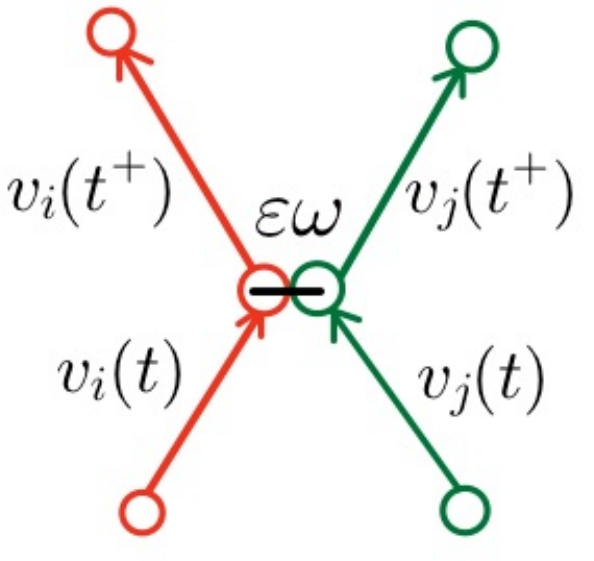}
\caption{An illustration of the hard sphere dynamics (Definition \ref{def.hard_sphere} (\ref{it.O32})): $(v_{i}(t),v_j(t))=(v_{i}(t^-),v_j(t^-))$ and $(v_{i}(t^+),v_j(t^+))$ are incoming and outgoing velocities, and $\varepsilon\omega$ is the vector connecting the centers of the two colliding particles, which has length $\varepsilon$.}
\label{fig.1}
\end{figure}
\item \label{it.O33} If for certain $i$ there is no $j$ such that the scenario in (\ref{it.O32}) happens, then we have
\begin{equation}\label{eq.hardsphere2}\frac{\mathrm{d}}{\mathrm{d}t}(x_i,v_i)=(v_i,0).
\end{equation}
\end{enumerate}
\item \label{it.O4}\emph{The flow map $\Hc_N(t)$ and flow operator $\Sc_N(t)$}. We define the \textbf{flow map} $\Hc_N(t)$ by
\begin{equation}\label{eq.defH^O}\Hc_N(t):\Dc_N\to\Dc_N,\quad\Hc_N(t)(\vz_N^0)=\vz_N(t)\end{equation} where $\vz_N(t)$ is defined by (\ref{it.O3}). Define also the \textbf{flow operator} $\Sc_N(t)$, for functions $f=f(\vz_N)$, by
\begin{equation}\label{eq.defS^O}
(\Sc_N(t) f)(\vz_N)=f\big((\Hc_N(t))^{-1}\vz_N\big).
\end{equation}
\item \label{it.O5} \emph{Collision}. We define a \textbf{collision} to be an instance that (\ref{eq.hardsphere}) in (\ref{it.O32}) occurs. Note that these collisions exactly correspond to the discontinuities of $v_j(t)$.
\end{enumerate}
\end{definition}
We list some basic properties of the hard sphere dynamics, which hold up to Lebesgue zero sets:
\begin{proposition}\label{prop.hardsphere} There exists a closed subset $Z\subseteq \Dc_N$ with Lebesgue measure $|Z|=0$, such that the time evolution defined by the hard sphere system (Definition \ref{def.hard_sphere} (\ref{it.O3})) exists and is unique on $\Dc_N\backslash Z$, and satisfies
\begin{enumerate}
\item No two collisions happen at the same time $t$;
\item The total number of collisions has an absolute upper bound that depends only on $d$ and $N$;
\item The flow maps $\Hc_N(t)$ defined in Definition \ref{def.hard_sphere} (\ref{it.O4}) are measure preserving diffeomorphisms from $\Dc_N\backslash Z$ to itself, and satisfy the flow or semi-group property $\Hc_N(t+s)=\Hc_N(t)\Hc_N(s)$ for $t,s\geq 0$.
\end{enumerate}
\end{proposition}
\begin{proof} This is proved in Alexander \cite{Alexander}.
\end{proof}

From now on we will always ignore the Lebesgue zero set $Z$ in Proposition \ref{prop.hardsphere} (as they do not contribute to Lebesgue integrals), and will not distinguish the set $\Dc_N$ from $\Dc_N\backslash Z$. We will do the same for all the modified dynamics below (Definitions \ref{def.modified_dynamics}, \ref{def.molecule_truncated_dynamics}, \ref{def.E_dynamics} and \ref{def.T_dynamics}).

\subsubsection{The grand canonical ensemble}\label{sec.grand_canon} In this paper, we study the hard sphere system with a random number of particles which have random initial states, whose law is given by the grand canonical ensemble. 

The grand canonical ensemble is defined by a density function $W_0:\cup_N \Dc_N\to\Rb_{\geq 0}$ (or equivalently a sequence of density functions $W_{0,N}:\Dc_N\to\Rb_{\geq 0}$) that determines the law of the random initial data for the hard sphere system. We state the definition as follows:
\begin{definition}[The grand canonical ensemble] \label{def.grand_canon} Fix $0<\varepsilon\ll 1$ and a nonnegative function $f_0=f_0(z)$ with $\int f_0(z)\,\mathrm{d}z=1$. We define the \textbf{grand canonical ensemble} as follows.
\begin{enumerate}
\item\label{it.Can0} \emph{The grand canonical domain $\Dc$}. We define $\Dc:=\bigcup_{N}\Dc_N$ to be the \textbf{grand canonical domain}, so $\vz\in\Dc$ always represents $\vz=\vz_N\in\Dc_N$ for some $N$. We also define the hard sphere dynamics $\vz(t)$ on $\Dc$, and the flow map $\Hc(t)$ on $\Dc$, by specifying them to be the $\vz_N(t)$ in Definition \ref{def.hard_sphere} (\ref{it.O3}), and the $\Hc_N(t)$ in Definition \ref{def.hard_sphere} (\ref{it.O4}), for $\vz^0=\vz_N^0\in\Dc_N$.
\item \label{it.Can1} \emph{Random data and initial density}. Assume $\vz^0\in \Dc$ is a random variable, whose law is given by the \textbf{initial density functions} $(W_{0,N})$, in the sense that \begin{equation}\label{eq.inti_density}
\Pb\big(\vz^0=\vz_N^0\in A\subseteq\Dc_N\big)=\frac{1}{N!}\int_{A} W_{0,N}(\vz_N)\,\mathrm{d}\vz_N
\end{equation}
for any $N$ and any $A\subseteq\Dc_N$, where $(W_{0,N})$ is given by
\begin{equation}\label{eq.N_par_ensemble} 
    \frac{1}{N!}W_{0,N}(\vz_N):=\frac{1}{\Zc}\frac{\varepsilon^{-(d-1)N}}{N!}\prod_{j=1}^N f_0(z_j)\cdot \mathbbm{1}_{\Dc_N}(\vz_N).
\end{equation}
Here, $\mathbbm{1}_{\Dc_N}$ is the indicator function of the set $\Dc_N$, and the partition function $\Zc$ is defined as follows
\begin{equation}\label{eq.partition_func}
    \Zc:=1+\sum_{N=1}^\infty\frac{\varepsilon^{-(d-1)N}}{N!}\int_{\Rb^{2dN}}\prod_{j=1}^N f_0(z_j)\cdot \mathbbm{1}_{\Dc_N}(\vz_N)\,\mathrm{d}\vz_N.
\end{equation}
\item \label{it.Can2}\emph{Evolution of random data}. Let $\vz^0\in \Dc$ be the random variable defined in (\ref{it.Can1}) above, and let $\vz(t)=\Hc(t)\vz^0$ be the evolution of initial data $\vz^0$ by hard sphere dynamics, defined as in Definition \ref{def.hard_sphere} (\ref{it.O3})--(\ref{it.O4}), then $\vz(t)$ is also a $\Dc$-valued random variable. We define the \textbf{density functions} $W_N(t,\vz_N)$ for the law of the random variable $\vz(t)$ by the relation
\begin{equation}\label{eq.inti_density_t}
\Pb\big(\vz(t)=\vz_N(t)\in A\subseteq\Dc_N\big)=\frac{1}{N!}\int_{A} W_{N}(t,\vz_N)\,\mathrm{d}\vz_N.
\end{equation} Then it is easy to see (using the volume preserving property in Proposition \ref{prop.hardsphere}) that
\begin{equation}\label{eq.time_t_ensemble}W_N(t,\vz_N):=(\Sc_N(t)W_{0,N})(\vz_N).
\end{equation}
\item\label{it.Can3} \emph{The correlation function $f_s(t)$}. Given $s\in\Nb$, we define the $s$-particle (rescaled) \textbf{correlation function} $f_s(t,\vz_s)$ by the following formula
\begin{equation}\label{eq.s_par_cor}
f_s(t,\vz_s):=\varepsilon^{(d-1)s}\sum_{n=0}^\infty\frac{1}{n!}\int_{\Rb^{2dn}}W_{s+n}(t,\vz_{s+n})\,\mathrm{d}z_{s+1}\cdots \mathrm{d}z_{s+n},
\end{equation} 
where we abbreviate $\vz_s=(z_1,\cdots,z_s)$ and $\vz_{s+n}=(z_1,\cdots,z_{s+n})$ similar to $\vz_N$.
\item\label{it.Can4} \emph{The Boltzmann-Grad scaling law}. The formula (\ref{eq.N_par_ensemble}), particularly the choice of the power $\varepsilon^{-(d-1)N}$, implies that
\begin{equation}\label{eq.Bol-Grad}\Eb(N)\cdot \varepsilon^{d-1}\approx 1
\end{equation} up to error $O(\varepsilon)$, where $\vz^0\in\Dc$ is the random initial data defined in (\ref{it.Can1}) and $N$ is the random variable (i.e. number of particles) determined by $\vz^0=\vz_N^0\in\Dc_N$. In other words, we are considering the hard sphere dynamics with $\approx \varepsilon^{-(d-1)}$ particles (on average). This is referred to as the \textbf{Boltzmann-Grad scaling law}, and will be fixed throughout the paper.
\end{enumerate}
\end{definition}

\subsubsection{The Boltzmann equation} Starting from works of Maxwell and Boltzmann in 1860--1870s, it is predicted that the one-particle correlation function $f_1(t,z)=f_1(t,x,v)$ should solve the Boltzmann equation in the kinetic limit $\varepsilon\to 0$, which is defined as follows.
\begin{definition}[The hard sphere Boltzmann equation]\label{def.boltzmann} We define the Cauchy problem for the \textbf{Boltzmann equation} for hard sphere collisions, with initial data $f(0,z)=f_0(z)$ where $z=(x,v)\in\Rb^{2d}$, as follows:
\begin{equation}\label{eq.boltzmann}
\left\{
\begin{aligned}
(\partial_t+v\cdot\nabla_x)f&=\int_{\Rb^d}\int_{\Sb^{d-1}}\big((v-v_1)\cdot\omega\big)_+\cdot(f'f_1'-ff_1)\,\mathrm{d}\omega\mathrm{d}v_1;\\
f(0,x,v)&=f_0(x,v).
\end{aligned}
\right.
\end{equation} The right hand side of (\ref{eq.boltzmann}) is referred to as the \textbf{collision operator}, where we denote
\[f=f(x,v),\quad f_1=f(x,v_1),\quad f'=f(x,v'),\quad f_1'=f(x,v_1')\] and
\begin{equation}\label{eq.boltzmann2} 
    v'=v-\big((v-v_1)\cdot\omega\big)\omega,\quad v_1'=v_1+\big((v-v_1)\cdot\omega\big)\omega,
\end{equation} 
and $\big((v-v_1)\cdot\omega\big)_+$ is the positive part of $(v-v_1)\cdot\omega$.
\end{definition}
\subsubsection{The main result} We now state the main result of this paper. First we introduce the following norm for functions $f: \Rb^d \times \Rb^d \to \Cb$,
\begin{equation}\label{eq.boltzmann_decay_2}
\|f(x,v)\|_{\mathrm{Bol}^\beta}:=\sum_{k\in\Zb^d}\sup_{|x-k|\leq 1,v\in\Rb^d}e^{\beta|v|^2}|f(x,v)|.
\end{equation} 
\begin{mainthm}\label{th.main} Fix $d\geq 2, \beta>0$, and a nonnegative function $f_0=f_0(z)\geq 0$, with $\int f_0(z)\,\mathrm{d}z=1$. Suppose the solution $f(t,z)$ to the Boltzmann equation (\ref{eq.boltzmann}) exists on the time interval $[0,t_{\mathrm{fin}}]$, such that
\begin{equation}\label{eq.boltzmann_decay_1}
\big\|e^{2\beta  |v|^2}f(t)\big\|_{L_{x,v}^\infty}\leq A<\infty,\quad \forall t\in[0,t_{\mathrm{fin}}].
\end{equation}
Note that $t_{\mathrm{fin}}$ can be arbitrarily large. Suppose also that the initial data $f_0$ satisfies
\begin{equation}\label{eq.boltzmann_decay_3}\|f_0\|_{\mathrm{Bol}^{2\beta}}+\|\nabla_x f_0\|_{\mathrm{Bol}^{2\beta}}\leq B_0<\infty.
\end{equation}

Consider the $d$ dimensional hard sphere system of diameter $\varepsilon$ particles (Definition \ref{def.hard_sphere}), with random initial data $\vz^0$ given by the grand canonical ensemble (Definition \ref{def.grand_canon}), under the Boltzmann-Grad scaling law \eqref{eq.Bol-Grad}. Let $\varepsilon$ be small enough depending on $(d,t_{\mathrm{fin}},\beta,A, B_0)$. Then, uniformly in $t\in[0,t_{\mathrm{fin}}]$ and in $s\leq |\log\varepsilon|$, the $s$-particle correlation functions $f_s(t)$ defined as in \eqref{eq.N_par_ensemble}--\eqref{eq.s_par_cor}, satisfy
\begin{equation}\label{eq.maineqn}
\bigg\|f_s(t,\vz_s)-\prod_{j=1}^s f(t,z_j)\cdot \mathbbm{1}_{\Dc_s}(\vz_s)\bigg\|_{L^1(\Rb^{2ds})}\leq \varepsilon^\theta,
\end{equation} 
where $\theta>0$ is an absolute constant depending only on the dimension $d$. 
\end{mainthm}
    
\begin{remark}
\label{rem.theorem} We make a few comments on the scope of Theorem \ref{th.main}.
\begin{enumerate}
\item \emph{Assumptions on the Boltzmann solution}: the Maxwellian decay of $f(t,z)$ in $v$ is standard for solutions to the Boltzmann equation, and in principle can be replaced by bounds in any space in which the Boltzmann equation is locally well-posed by fixed point arguments. On the other hand, the required integrability in $x$ in (\ref{eq.boltzmann_decay_2}) for the initial data $f_0$ is basically at the level of $L_x^1$, which is the minimal requirement to define and normalize the ensemble in \eqref{eq.N_par_ensemble}--(\ref{eq.partition_func}).
\item\label{it.intro_rem_2} \emph{The time interval $[0,t_{\mathrm{fin}}]$}: if the solution $f(t,z)$ to the Boltzmann equation (\ref{eq.boltzmann}) develops singularity at some finite time $t_{\mathrm{sing}}$, then the $t_{\mathrm{fin}}$ in Theorem \ref{th.main} can be chosen as anything strictly less than $t_{\mathrm{sing}}$. If the solution exists and is uniformly bounded for all time $t>0$, then the $t_{\mathrm{fin}}$ in Theorem \ref{th.main} can be chosen as a function of $\varepsilon$ that diverges as $\varepsilon\to 0$; the current proof would lead to $t_{\mathrm{fin}}\sim\log|\log\varepsilon|$ (or a power of it). We do not expect this to be optimal, and it is an interesting question whether this can be improved (say to $t_{\mathrm{fin}}\sim|\log\varepsilon|$).
\item \emph{Parameters and convergence rate}:  the restriction $s\leq|\log\varepsilon|$ and the convergence rate $\varepsilon^\theta$ in (\ref{eq.maineqn}) are not expected to be optimal, and are chosen primarily for convenience. The method of our proof would cover the case $s\leq \varepsilon^{-\kappa}$ for sufficiently small $\kappa>0$ depending on $t_{\mathrm{fin}}$; on the other hand, if $s$ is as large as $\varepsilon^{-1}$ then (\ref{eq.maineqn}) might not hold even for initial data, see \cite{PS17}.

Next, for $s=1$, the optimal value of $\theta$ should be $1-$ and can be obtained by refining our proof; as for the non-factorized parts or cumulants for $s\geq 2$, their optimal decay rate is linked to the derivation of the fluctuation equations about the Boltzmann limit, which we will discuss in Section \ref{sec.future_3} below. Finally, the value $\Eb(N)\cdot\varepsilon^{d-1}$ (inverse of the Knudsen number) occurring in the Boltzmann-Grad scaling law (\ref{eq.Bol-Grad}) can be replaced by any constant $\alpha$ independent of $\varepsilon$ (by replacing the $\varepsilon^{-(d-1)N}$ in (\ref{eq.N_par_ensemble})--(\ref{eq.partition_func}) with $(\alpha\varepsilon^{-(d-1)})^N$), or (a suitable power of) $\log|\log\varepsilon|$ as in (\ref{it.intro_rem_2}).
\item \emph{The grand canonical ensemble}: we choose to work with the grand canonical ensemble (Definition \ref{def.grand_canon}), as opposed to the canonical ensemble where the particle number $N$ is fixed, only to simplify some calculations with cumulants. This simplification was observed in several recent works on this subject \cite{PS17,BGSS20, BGSS20_2, BGSS22_2}. Theorem \ref{th.main} remains true for the canonical ensemble, but the proof would be much longer due to the extra manipulations needed to reduce the expressions of cumulants. 
\end{enumerate}
\end{remark}
\subsection{Background and history}\label{sec.past} We briefly discuss the background and history of the rigorous derivation of Boltzmann's equation, and the closely related problem of deriving its wave counterpart, i.e. the wave kinetic equation in wave turbulence theory.
\subsubsection{Particle systems and the Boltzmann equation}\label{sec.past_1} The study of large particle systems has a long history, starting from the beginning of the theory of statistical physics. The kinetic theory, which lies at the heart of nonequilibrium statistical physics, started with the following genius observation of Boltzmann \cite{Bol72,Bol02}: assuming independence of states between different particles, the full statistics of the system can be described by a single closed effective equation, namely the \emph{kinetic equation}, in the limit of infinite degrees of freedom.

The first instance of what is now known as the Boltzmann equation was written down in the 1870s; since then, the mathematical justification of this limiting process has been a major challenge. The main difficulty comes from justifying the assumption of (asymptotic) independence of different particles, which is known as the \emph{Stosszahlansatz} or \emph{molecular chaos}. The formulation of this problem as known today was first put forward by Grad, and the works of Cercignani and Grad in the 1950s introduced the important idea of the BBGKY hierarchy to this subject. In the mean time, a simplified model was also introduced and studied by Kac (\cite{Kac56}, see also subsequent and related developments \cite{OVY93,QY98,MS11,MT12,MM13}). Motivated by these contributions, Lanford \cite{Lan75} (see also King \cite{King75}) was able to set a first milestone in this subject, by proving the molecular chaos assumption and deriving the Boltzmann equation for sufficiently small time.

Since \cite{Lan75}, there have been numerous results devoted to the derivation of the Boltzmann equation and related problems. We do not attempt to be exhaustive in recounting the vast literature, but only list some of the most important ones below; see the surveys \cite{SpohnBook,BGSS18, BGSS23} and references therein for the other contributions.
\begin{itemize}
\item Illner-Pulvirenti \cite{IP86} (see also \cite{IP89}), which obtained a long-time result but for the near-vacuum case (i.e. the collision rate $\Eb(N)\cdot\varepsilon^{d-1}\ll 1$ in (\ref{eq.Bol-Grad}));
\item Gallagher-Saint Raymond-Texier \cite{GST14}, which completed and refined the result of \cite{Lan75}, and obtained explicit convergence rates;
\item Pulvirenti-Simonella \cite{PS17}, which obtained explicit decay rates for cumulants associated with the hard sphere system, with the exponent proportional to the number of involved particles, again in the short-time regime;
\item Bodineau-Gallagher-Saint Raymond \cite{BGS16, BGS17,BGS18}, which studied the related problems of deriving Brownian motion, Ornstein-Uhlenbeck processes, and the Stokes-Fourier system from the hard sphere dynamics, by passing through a linearized Boltzmann equation;
\item Bodineau-Gallagher-Saint Raymond-Simonella \cite{BGSS20,BGSS20_2,BGSS22,BGSS22_2}, which investigated the further question of deriving the equation for fluctuations, either for short time, or for longer times but at the exact Maxwellian equilibrium. 
\end{itemize}

Finally, apart from the question of deriving the kinetic equations, there has also been a huge literature concerning the solution theory of these equations (Boltzmann, Landau etc.). In general strong solutions exist in local and perturbative settings, but only weak solutions are known globally (see however \cite{GS23,ISV24} on recent global existence results in the homogeneous setting). We will not elaborate on this, but refer to the reviews in \cite{UT06, Villani} and references therein.
\subsubsection{Wave turbulence and the wave kinetic equation}\label{sec.past_2} One of the most important developments in the physics of the twentieth century was the emergence of quantum physics as a theory to describe subatomic particle dynamics. A central idea here is the \emph{particle-wave} duality, which contributed to an exponential rise in the physical and mathematical interest in nonlinear wave systems. This led to the recasting of many classical particle theories in the context of interacting wave systems, including the theory of nonequilibrium statistical physics for waves, known as the wave turbulence or wave kinetic theory.

Same as the particle case, the \emph{wave kinetic theory} aims at describing the statistics of a wave system using a closed equation for the lowest order correlation; usually this is the second moment of the wave field, as opposed to the one-particle density function in the kinetic theory for particles. The role of the Boltzmann equation is then played by the \emph{wave kinetic equation} (WKE). For instance, in dimension $d=3$, when the wave system is modeled by the cubic nonlinear Schr\"odinger equation, the corresponding wave kinetic equation takes the form 
\begin{equation}\label{eq.wke}
(\partial_t+v\cdot\nabla_x)f=\int_{\Rb^3}\int_{\Sb^{2}}\big((v-v_1)\cdot\omega\big)_+\cdot\big(f'f_1'(f+f_1)-ff_1(f'+f_1')\big)\,\mathrm{d}\omega\mathrm{d}v_1
\end{equation} (including also the homogeneous variant where the $v\cdot\nabla_x$ term in (\ref{eq.wke}) is absent), where we use the same notation as in \eqref{eq.boltzmann2}. The wave kinetic equations, similar to \eqref{eq.wke} but for different wave systems, were first introduced in the work of Peierls \cite{Pei29} in 1929 for anharmonic crystals. This was shortly followed by works of Nordheim \cite{Nor28} and Uehling and Uhlenbeck \cite{UU33} in the context of quantum gases. Since then, the theory has attracted a great deal of attention from mathematicians and physicists, and has developed into a rich and broad subject. In particular, these contributions included the influential works of Hasselmann \cite{Has62,Has63} in 1962--1963, which applied the wave kinetic theory to ocean waves and led to the modern methods of ocean surface forecasting, as well as the work of Zakharov \cite{Zak65} in 1965, which provided a systematic approach to the theory and how it can be used to understand turbulent phenomena for waves. In \cite{Zak65}, Zakharov also introduced the notion of Kolmogorov-Zakharov spectra, along the same lines as Kolmogorov's theory of hydrodynamic turbulence. It is for this reason that the theory is also commonly known as wave turbulence theory. We refer the reader to \cite{Naz11} for a textbook treatment.

Mathematically, the study of the derivation of wave kinetic equations like (\ref{eq.wke}) started with the work of Spohn \cite{Spo77}. For linear models, this program was completed in Erd\"{o}s-Yau \cite{EY00} and Erd\"{o}s-Salmhofer-Yau \cite{ESY08}. For the full nonlinear problem, an important early milestone was due to Lukkarinen-Spohn \cite{LS11} and subsequent partial results were obtained in \cite{BGHS19,DH19,CG19,CG20}. We refer the reader to the introduction of \cite{DH23} for a more exhaustive historical recount. The full derivation of (\ref{eq.wke}) in the short-time regime was obtained by the first two authors \cite{DH21,DH21_2,DH23}, serving as the wave-counterpart of Lanford's theorem.

Recently, in \cite{DH23_2}, the first two authors have extended the derivation of (\ref{eq.wke}) to arbitrarily long times that cover the full lifespan of solutions to (\ref{eq.wke}). We remark that the result of \cite{DH23_2} is parallel to Theorem \ref{th.main}, not only conceptually but also on a very concrete level. In fact, while the proof in the current paper is self-contained, many of the key ideas in this proof are obtained by adapting the ideas from the proof in \cite{DH23_2} to the particle setting. In Section \ref{subsec.idea.compare} we will make a detailed comparison between the proofs in \cite{DH23_2} and in the current paper, including both similarities and differences.

\subsection{Main ideas of the proof}\label{sec.strategy} In this subsection, we present a mostly philosophical description of the proof, its strategy, and its most important ideas. A more elaborate and concrete discussion of these ideas can be found in Section \ref{sec.main_idea}. 

To a large extent, the framework of the proof in this paper is the same as that in \cite{DH23_2}. Below we will explain this framework in the current context of particle systems.

\subsubsection{Time layering}\label{sec.intro.layers} The well-known fundamental difficulty in extending Lanford's theorem to longer times, which is the same difficulty in the long-time derivation of the wave kinetic equation, is the \emph{divergence of time  expansions of solutions}. In fact, in both settings, the proof of the short-time result is based on matching two different time expansions: expansion I, given by the Duhamel expansion of the dynamics of the microscopic system (or the BBGKY hierarchy satisfied by the correlation functions $f_s$), and expansion II, given by the Duhamel expansion of the kinetic equation. The bulk of the proof\footnote{In the wave case, there is an extra level of complication coming from the fact that the series does not converge absolutely, but only conditionally.}  goes into identifying the leading terms in expansion I, matching them order-by-order with the terms in expansion II, and proving that all the remaining non-leading terms vanish in the limit. At a basic level, both expansions are bounded like $\sum_n C^n t^n$, which is convergent only for sufficiently small times $|t|<C^{-1}$ given by the radius of convergence. Beyond this time both expansions will diverge, and this methods fails in a fundamental way.

To go beyond this threshold and reach an arbitrarily large time $t= t_{\mathrm{fin}}$, a natural idea is to subdivide the long time interval $[0,t]$ into $\Lf$ subintervals (called time layers) $[(\ell-1) \tau, \ell \tau]$ with $1\leq \ell \leq \Lf$ and $\tau$ sufficiently small, such that on each subinterval, the solution to the kinetic equation is given by a convergent series expansion with base time $(\ell-1) \tau$. This \textbf{time layering} is the first key idea in the proof, which allows to resolve the divergence problem for expansion II mentioned above, regarding solutions to the Boltzmann equation. This leaves us with the task of controlling the microscopic system (expansion I) on $[(\ell-1) \tau, \ell \tau]$ using the data at time $(\ell-1)\tau$. An obvious difficulty arises with this approach, as the data at $(\ell-1)\tau$ will no longer have exact independence between different particles. Such independence (or equivalently the complete factorization of the correlation functions $f_s$) is granted\footnote{Strictly speaking, we do not have exact independence due to the $\mathbbm 1_{\Dc_N}(\vz_N)$ factor in \eqref{eq.N_par_ensemble}, but the cumulants caused by this factor are easy to deal with (see Proposition \ref{prop.initial_cumulant}), and will be ignored here.} at time $0$, and plays a fundamental role in controlling the expansion I of the microscopic system in the short-time proof.

\subsubsection{Cumulants}\label{subsec.cumexp_intro} In order to go from time $(\ell-1)\tau$ to $\ell\tau$, we need to propagate some information for the ensemble at time $\ell\tau$ that describes departure from independence. This will be given by the {\bf cumulants}. Following \cite{PS17}, the cumulants can be defined as the $E_H$ occurring in the following expansion of the $s$-point density function:
\begin{equation}\label{eq.cumulant_intro}
    f_s(t, \vz_s) = \prod_{j \in [s]} f^{\Ac}(t, z_j) + 
    \underbrace{\sum_{\varnothing \neq H \subseteq [s]} 
    \bigg(\prod_{j \in [s] \setminus H} f^{\Ac}(t, z_j)\bigg) 
    \cdot E_H(t, \vz_H)}_{O_{L^1}(\varepsilon^{\theta})},
\end{equation}
where $[s] = \{1, \ldots, s\}$ and $\vz_H = (z_j)_{j \in H}$. Here the first term on the right-hand side represents the fully factorized part of the correlation, with $f^{\Ac}(t, z)$ being a one-particle density function, and the remaining terms involving cumulants $E_H$ describe the deviation from independence. Similar cumulant expansions have also played a central role in the wave setting \cite{DH23_2} and in a series of recent works on particle systems \cite{PS17, BGSS22_2, BGSS20, BGSS23}. 

The physical interpretation of \eqref{eq.cumulant_intro} is as follows: the left-hand side denotes the probability density of finding $s$-particles, say $1, \cdots, s$, with states $z_1,\cdots, z_s$, at time $t$. This depends on the collision history of these particles among themselves and with remaining particles. The first term on the right hand side comes from collision histories where the $s$ particles are disconnected, i.e. they do not interact with each other and each of them interacts with a disjoint set of the remaining particles. The second term comes from collision histories where the particles in $H$ are connected to each other via collisions. 

By (\ref{eq.cumulant_intro}), we see that Theorem \ref{th.main} would follow if we can show that (a) $f^\Ac$ converges to the Boltzmann solution, and (b) $\|E_H\|_{L^1}$ converges to zero in the limit as $\varepsilon \to 0$.

\subsubsection{The cumulant ansatz}\label{subsec.cumansatz_intro} As discussed in Section \ref{subsec.cumexp_intro}, a natural route of proving Theorem \ref{th.main} is to approximate $f^\Ac$ by the Boltzmann solution (Proposition \ref{prop.est_fa}; this is relatively easy, see Section \ref{subsec.idea.fA} for further discussions, and Section \ref{sec.fa} for the proof), and prove that
\begin{itemize}
\item For each $1\leq \ell\leq \Lf$, we have the cumulant expansion (\ref{eq.cumulant_intro}) at $t=\ell\tau$, with some constant $\alpha>0$ such that
\begin{equation}\label{eq.cumuest_intro}
\|E_H(\ell\tau)\|_{L^1}\leq \varepsilon^{\alpha|H|},\quad \forall H\neq\varnothing.
\end{equation}
\end{itemize}

This $E_H$ estimate is stated in Proposition \ref{prop.cumulant_est}. Note that this is a reasonable conclusion to aim for. In fact, for $\ell=0$ (\ref{eq.cumuest_intro}) is trivial, while for $\ell=1$ this follows from Lanford's theorem  and the proof in \cite{PS17}. The next natural question is: how to proceed on the intervals $[(\ell-1)\tau,\ell\tau]$, to prove (\ref{eq.cumuest_intro}) at time $\ell\tau$.

In \cite{PS17}, the authors show that if the cumulants at $(\ell-1)\tau$ satisfy the $L^\infty$ estimate 
$$\|E_H((\ell-1)\tau)\|_{L^\infty}\leq \varepsilon^{\beta |H|},$$
then one can recover the $L^1$ estimate \eqref{eq.cumuest_intro} at time $\ell \tau$. The crucial point here is that one has to start with the stronger $L^\infty$ bound, and can only recover the weaker $L^1$ bound, therefore this estimate cannot be iterated.

There is something deep behind this phenomenon: for \emph{any} Banach space $X$, it is \emph{not possible} to obtain smallness of $\|E_H(\ell\tau)\|_X$ (in the form of (\ref{eq.cumuest_intro})) using only smallness of the same norm $\|E_H((\ell-1)\tau)\|_{X}$, precisely because of the \emph{time irreversibility} of Boltzmann's equation.

In fact, suppose the contrary, that such an inductive proof of the $\|\cdot\|_X$ bound was possible. Take $\ell=2$. Then, by the \emph{time reversibility} of the Newtonian dynamics (Definition \ref{def.hard_sphere}), we could either go forward from time $\tau$ to time $2\tau$, or backward to time $0$, which just differ by flipping the signs of the velocities $v_j$ at time $\tau$. Since $\|\cdot\|_X$ is invariant under such sign flipping (which is true for any natural Banach norm such as $L^p$ norms), we would conclude from \eqref{eq.cumulant_intro} (assuming the approximation $f^\Ac(t,z)\approx f(t, z)$) that
\[
f_s(2\tau, z_1, \cdots, z_s)=f(2\tau,z_1) \cdots f(2\tau, z_s)+O_{X}(\varepsilon^\upsilon)
\]
as well as
\[
f_s(0, z_1, \cdots, z_s)=f(2\tau,z_1) \cdots f(2\tau, z_s)+O_{X}(\varepsilon^\upsilon).
\]
This is impossible, as the leading term in $f_s(0,\cdots)$ is given by tensor product of $f(0,\cdot)$ and not $f(2\tau,\cdot)$.

By the above discussion, we know that \eqref{eq.cumuest_intro} cannot be proved by direct induction in $\ell$. In fact, in order to prove \eqref{eq.cumuest_intro}, we need to propagate \textbf{structural information} (which goes beyond any Banach norms) for the cumulants $E_H(\ell\tau)$ by induction in $\ell$. Such structural information will be provided by a \textbf{partial time expansion} argument, a second key ingredient of the proof, which we discuss below.

Start from a cumulant $E_H(\ell\tau)$, and assume $f^\Ac$ has already been approximated by the Boltzmann solution. By performing a cluster expansion (see Section \ref{sec.proof_cumulant_formula} or \cite{BGSS22_2}) on $[(\ell-1)\tau,\ell\tau]$, we can express this cumulant as a sum of Duhamel time integrals involving $f^\Ac((\ell-1) \tau)$ and $E_H((\ell-1) \tau)$. Now, the crucial point in the partial time expansion is that, we \emph{refrain from} expanding $f^\Ac((\ell-1) \tau)$ any further in time (and use instead its proximity to the Boltzmann solution), and only expand the cumulant terms $E_H((\ell-1) \tau)$ in terms of $f^\Ac((\ell-2)\tau)$ and $E_H((\ell-2) \tau)$, and so on. This is very important in our proof as it allows us to avoid the $t^n$ divergence (discussed in Section \ref{sec.intro.layers}); indeed, such divergence only comes from the aggregation if we were to expand the leading terms $f^\Ac$ all the way to time $t=0$, which we are avoiding in this partial time expansion. On the other hand, the cumulants come with power gains (of size $\varepsilon^{\alpha|H|}$, see (\ref{eq.cumuest_intro})), which offsets the $t^n$ divergence even if they are expanded all the way to time $t=0$.

This structural information constructed above will be further discussed in Section \ref{sec.molecule_representation_intro} below. As far as we know, such partial expansion construction has not appeared in previous literature, and is crucial in obtaining a time irreversible equation from time reversible dynamics.

\subsubsection{Representing cumulants by molecules}\label{sec.molecule_representation_intro} Recall the physical interpretation of the cumulants in Section \ref{subsec.cumexp_intro} and the partial time expansion in Section \ref{subsec.cumansatz_intro}. With these, we know that the cumulant $E_H(\ell \tau)$ can be written as a sum of contributions corresponding to probability densities of possible collision histories on $[0,\ell\tau]$ involving particles in $H$. These collision histories are constructed by the partial time expansions, and they satisfy that the particles in $H$ are connected to each other via collisions (cf. Section \ref{subsec.cumexp_intro}).

If we disregard the precise positions, velocities and collision times in the particle trajectories, and focus only on their combinatorial structure, then each collision history will \textbf{topologically reduce} to an abstract diagram, called a \textbf{molecule}. Such diagrams encode the possible patterns of collisions among the particles in $H$, and disregard the precise geometric information. This process of abstraction is illustrated in {\color{blue}Figure \ref{fig.topological_reduction}} (Section \ref{sec.molecule}), where particle trajectories are depicted as sets of edges sharing the same color, and collisions are represented by nodes, called \textbf{atoms}. In addition, the molecule also retains the ordering between different collisions of the same particle, and the layer of each collision (represented by horizontal lines in {\color{blue}Figure \ref{fig.topological_reduction}}).

Now, with the molecular representation of collision histories, we can formally state the structural information for the cumulant $E_H(\ell\tau)$ in Section \ref{subsec.cumansatz_intro}, as follows:
\begin{equation}\label{eq.molecule_representation_intro}
    |E_H(\ell \tau)| \le \sum_{\Mb} |\mathcal{I}\mathcal{N}_{\Mb}|,
\end{equation}
where the sum is taken over molecules $\Mb$ that are constructed from the partial expansion process in Section \ref{subsec.cumansatz_intro} (see Definition \ref{def.set_T_F} and Proposition \ref{prop.cumulant_formula}), and $|\Ic\Nc_\Mb|$ is an explicit expression representing the normalized probability density that the collision history described by $\Mb$ happens (see Definition \ref{def.associated_int}). 

The formulas \eqref{eq.cumulant_intro} and \eqref{eq.molecule_representation_intro} will be stated in Proposition \ref{prop.cumulant_formula} and proved in Section \ref{sec.proof_cumulant_formula}. For reasons that will be discussed in Section \ref{sec.counting_molecule_intro} below, in this proof we also need a truncation on the dynamics, which is defined in Section \ref{sec.truncation_large_molecule}. For a more detailed discussion of all the above ingredients, see Section \ref{subsec.idea.cumulants}.

\subsubsection{Reduction to molecule estimates and truncation of the dynamics}\label{sec.counting_molecule_intro} To obtain the bound on $\|E_H\|_{L^1}$ in \eqref{eq.cumuest_intro} using \eqref{eq.molecule_representation_intro}, one has to obtain effective estimates for the $L^1$ norm of the normalized probability densities $|\Ic\Nc_\Mb|$. Here, the central quantity in understanding such estimates is the number of \textbf{recollisions} $\rho$ in the molecule $\Mb$. This is defined as the \emph{circuit rank} of $\Mb$, viewed as an undirected graph, which is the number of independent cycles in $\Mb$. Heuristically, one expects that as the recollision number $\rho$ increases, the number of molecules $\Mb$ in the sum \eqref{eq.molecule_representation_intro} also increases. This divergence has to be balanced by gains in the estimates on $|\Ic\Nc_\Mb|$ for molecules with larger $\rho$. This heuristic is quantified by establishing the following two estimates. Given $\rho\geq 1$, we have:
\begin{itemize}
    \item \emph{Bounding the number of terms in \eqref{eq.molecule_representation_intro}.} The number of molecules $\Mb$ with $\rho$ recollisions satisfies the bound
    \begin{equation}\label{eq.molecule_number_upp_intro}
        \# \Mb \le C^{|\Mb|} \cdot |\log \varepsilon|^{C\rho};
    \end{equation}
    \item \emph{Bounding the contribution of each term.} The contribution of each molecule $\Mb$ is bounded by
    \begin{equation}\label{eq.molecule_term_upp_intro}
         \big\||\mathcal{I}\mathcal{N}_{\Mb}|\big\|_{L^1} \le \tau^{|\Mb|} \cdot \varepsilon^{\upsilon \rho}\,(\upsilon>0).
    \end{equation}
\end{itemize}

Note how the above two bounds exhibit the effect of recollisions: each more recollision leads to $|\log\varepsilon|^C$ loss in the number of molecules in \eqref{eq.molecule_number_upp_intro}, which is offset by the $\varepsilon^{\upsilon}$ gain in the renormalized probabilities in \eqref{eq.molecule_term_upp_intro}. The cumulant estimate \eqref{eq.cumuest_intro} then follows directly from \eqref{eq.molecule_representation_intro}--\eqref{eq.molecule_term_upp_intro}, provided that $\varepsilon\ll\tau\ll 1$.
 
From our construction, the molecule $\Mb$ consists of $\ell$ layers, and can be seen as a concatenation of $\ell$ sub-molecules $\Mb_{\ell'}\,(\ell'\leq\ell)$, where $\Mb_{\ell'}$ contains the atoms in layer $\ell'$. To prove \eqref{eq.molecule_number_upp_intro} and \eqref{eq.molecule_term_upp_intro}, we need to impose an upper bound on the number of atoms and recollisions in each of the sub-molecule $\Mb_{\ell'}$. To enforce the above restrictions, we introduce the notion of \textbf{truncated dynamics} (see Definition \ref{def.T_dynamics}) within each time layer. More precisely, in the truncated dynamics, at each time where a collision would occur, we check whether or not this collision would cause a violation to the following two conditions:
\begin{enumerate}
\item\label{it.intro_it_1} The number of particles in each cluster is at most $\Lambda:=|\log\varepsilon|^{O(1)}$, and 
\item\label{it.intro_it_2} The number of recollisions (i.e. independent cycles) in each cluster is at most $\Gamma=O(1)$,
\end{enumerate} 
where \textbf{clusters} are the subsets of particles that are connected by collisions. If this collision would cause a violation, then we ``turn off'' this collision and allow the particles to cross each other; otherwise we allow this collision just as in the original hard sphere dynamics. 

Of course, this truncation also leads to an error term $f_s^{\mathrm{err}}$ which occurs when the original dynamics (or in fact the extended version of it, see Definition \ref{def.E_dynamics}) does not coincide with the truncated dynamics. We will control this error term in Proposition \ref{prop.cumulant_error}. In fact, it can be expressed in terms of $|\mathcal{I}\mathcal{N}_\Mb|$, and bounded using similar arguments, in a much similar way as $E_H$, see Section \ref{subsec.idea.fA} and Section \ref{sec.error}.

Given \eqref{it.intro_it_1} and \eqref{it.intro_it_2}, it follows that the number of atoms in any molecule $\Mb$ can be bounded by $|\log \varepsilon|^{O(1)}$. The estimate \eqref{eq.molecule_number_upp_intro} then follows from the well-known combinatorial fact that the number of graphs whose vertices have maximum degree 4, with at most $m$ vertices, and $\rho$ independent cycles  is bounded by $C^m \cdot m^\rho$ for some constant $C > 0$. Applying this to molecules yields:
\[
    \#\Mb \le C^{|\Mb|} \cdot |\Mb|^{\rho} \le C^{|\Mb|} \cdot |\log \varepsilon|^{C \rho},
\]
as claimed in \eqref{eq.molecule_number_upp_intro} (Proposition \ref{prop.layerrec3}; see Section \ref{subsec.idea.molecule} for more details). As a result, we are left with proving the estimate \eqref{eq.molecule_term_upp_intro}, which is the most delicate part of the proof. 

\subsubsection{The cutting algorithm}\label{sec.cutting_alg_intro} In the proof, we will reduce the estimate of $\big\||\mathcal{I}\mathcal{N}_{\Mb}|\big\|_{L^1}$ to a purely combinatorial problem for the molecule $\Mb$. To this end, we first rewrite it as $\big\||\mathcal{I}\mathcal{N}_{\Mb}|\big\|_{L^1} = \varepsilon^{(d-1)|H|}\cdot\mathcal{I}_{\Mb}(Q)$, where $\mathcal{I}_{\Mb}$ is an operator involving integration over all intermediate positions and velocities (associated with the edges of $\Mb$) and all collision times (associated with the atoms of $\Mb$), and $Q$ is a specific input function. These integrals involve Dirac $\dirac$ functions, which encode the relations between incoming and outgoing velocities at each collision. The reduction to this form is contained in Section \ref{sec.local_associated}, see Proposition \ref{prop.local_int}.

Then, to get optimal bounds for $\Ic_\Mb(Q)$, we need to find the correct order of variables in which to integrate. For this, we introduce the key notion of \textbf{cutting}, see Definition \ref{def.cutting}, which provides a systematic way of deciding the order of integration in $\Ic_\Mb$ (i.e. applying the Fubini theorem). More precisely, each cutting divides the molecule $\Mb$ into two separate molecules $\Mb_1$ and $\Mb_2$; by first integrating in those variables associated with the molecule $\Mb_2$ (with those variables associated with the molecule $\Mb_1$ regarded as fixed) and then integrating in those variables associated with the molecule $\Mb_1$, we get that
\begin{equation}\label{eq.molcut_intro}
    \Ic_\Mb = \Ic_{\Mb_1}\circ \Ic_{\Mb_2}.
\end{equation}

Now, using (\ref{eq.molcut_intro}), our strategy is to devise a \textbf{cutting algorithm} which cuts $\Mb$ into a number of ``small" molecules $\Mb_j$ with one or two atoms, called \textbf{elementary molecules} (see Definition \ref{def.elementary}). The operators $\Ic_{\Mb_j}$ can be directly calculated and estimated (say in the $L^\infty\to L^\infty$ norm), so that applying (\ref{eq.molcut_intro}) then yields an upper bound for $\Ic_{\Mb}$. The two most important types of elementary molecules are the $\{3\}$ molecules which represent collisions (i.e. a particle $\pb$ colliding with a given particle $\pb'$), and the $\{33\}$ molecules which represent collisions followed by recollisions (i.e. a particle $\pb$ colliding with a given particle $\pb'$, and then colliding again with another given particle $\pb''$). It is known that $\{3\}$ molecules are \textbf{normal} in the sense that integrating it does not gain any power, while $\{33\}$ molecules are \textbf{good} in the sense that integrating it gains a power $\varepsilon^\upsilon$ for some constant $\upsilon>0$ (there are also $\{4\}$-molecules which are usually viewed as \textbf{bad} due to ``wasting" of dimensions of integration). The discussion of elementary molecules, and associated weight and volume bounds, are proved in Section \ref{sec.elem_int}, see Definition \ref{def.good_normal}.

With the above discussions, it is now clear that the main ingredient in the proof of \eqref{eq.molecule_term_upp_intro} is to construct a cutting algorithm (see Proposition \ref{prop.comb_est}) that produces sufficiently many good molecules. This is a major novelty and main technical ingredient, and is discussed in more details in Section \ref{subsec.idea.alg}. Roughly speaking, the cutting algorithm consists of two big steps: (A) reducing the multi-layer case to the two-layer case by performing a layer refining and a layer selection process; (B) treating the two-layer case by exhibiting a dichotomy between two main scenarios and applying several basic algorithms.

In Section \ref{sec.toy}, we start with a simplified case, illustrating the main ideas of the algorithm; then we extend these ideas to full generality, and present the proofs in steps (A) and (B) in Sections \ref{sec.layer} and \ref{sec.maincr} respectively.

\subsection{Future horizons}\label{sec.future} Here we list some future directions, now within reach, with the proof of Theorem \ref{th.main}.
\subsubsection{The hydrodynamic limit}\label{sec.future_2} After deriving the Boltzmann equation from particle dynamics in the kinetic limit, we would need to match it with the hydrodynamic limit and derive the fluid equations (Euler and Navier-Stokes). Formally, this requires considering densities that are \emph{local Maxwellians}, i.e. those functions that are Maxwellian in $v$ but with mass, center and variance depending on $(t,x)$:
\begin{equation}\label{eq.localmaxw}f(t,x,v)=\frac{\rho(t,x)}{(2\pi T(t,x))^{d/2}}\exp\bigg(-\frac{|v-u(t,x)|^2}{2T(t,x)}\bigg).\end{equation} By considering perturbation expansions (i.e. Hilbert and Chapman-Enskog expansions) around (\ref{eq.localmaxw}), we can derive the compressible Euler equation for $(\rho,u,T)$ and the incompressible Euler and Navier-Stokes equations for subsequent terms, in a suitable scaling limit.

Compared to the derivation of the Boltzmann equation, the derivation of fluid equations is much better understood, see for example \cite{DEL89,GS04,Sai09,GT20}. In \cite{DHM25}, by combining these results with an analog of Theorem \ref{th.main} on the torus $\Tb^d$, we can finally execute both steps of Hilbert's program in tandem, and derive the fluid equations from Newtonian particle particle systems in the Boltzmann-Grad limit.

\subsubsection{The short-range potential case}\label{sec.future_1} With the methods in the current paper, a natural problem that should be within reach, is the short range potential case studied in \cite{GST14}, namely
\begin{equation}\label{eq.potential}\frac{\mathrm{d}x_i}{\mathrm{d}t}=v_i,\quad \frac{\mathrm{d}v_i}{\mathrm{d}t}=-\frac{1}{\varepsilon}\sum_{j\neq i}\nabla\Phi\bigg(\frac{x_i-x_j}{\varepsilon}\bigg),
\end{equation} where $\Phi$ is a suitable potential function. We believe that Theorem \ref{th.main} should remain true in this case, but the proof would need to involve (among other things) extra technical ingredients to deal with the possibility of simultaneous interactions between more than two particles via (\ref{eq.potential}).

\subsubsection{Fluctuations and large deviations}
\label{sec.future_3} The rigorous study of this direction was initiated in the recent works of Bodineau-Gallagher-Saint Raymond-Simonella \cite{BGSS20,BGSS20_2,BGSS22}. Here one considers the random variable (i.e. observable) $\zeta:=\frac{1}{N}\sum_{j=1}^N F(z_j(t))$, where $F=F(z)$ is a fixed function and $\vz_N(t)=(z_j(t))_{1\leq j\leq N}$ is the solution to the hard sphere system in Definition \ref{def.hard_sphere} (\ref{it.O3}). Then up to errors that vanish when $\varepsilon\to 0$, we have
\[\zeta\approx \Eb(\zeta)\approx\int_{\Rb^{2d}}F(z)f(t,z)\,\mathrm{d}z,\] which can be viewed as a ``law of large numbers" for the random variables $F(z_j(t))$. The fluctuation, which is given by
\begin{equation}\eta=\sqrt{N}\bigg(\frac{1}{N}\sum_{j=1}^N F(z_j(t))-\int_{\Rb^{2d}}F(z)f(t,z)\,\mathrm{d}z\bigg),
\end{equation} should then be described by the ``central limit theorem". Such results have been established in the short-time regime in \cite{BGSS20}, and for longer times at the exact global Maxwellian equilibrium (i.e. (\ref{eq.localmaxw}) with $\rho=T=1$ and $u=0$) in \cite{BGSS20_2,BGSS22}.

With the proof of Theorem \ref{th.main}, we expect that these results in \cite{BGSS20,BGSS20_2,BGSS22} can be extended to the long-time and off-equilibrium setting, where the initial independence assumption in the short time setting is substituted by the precise description of $f_s$ and $E_H$ at any time $\ell\tau$ given by \eqref{eq.cumulant_intro} and (\ref{eq.molecule_representation_intro}).

\subsection{Parameters and notations}\label{sec.notation}  Before starting the proof, we first note that, by the assumptions (\ref{eq.boltzmann_decay_1}) and (\ref{eq.boltzmann_decay_3}) in Theorem \ref{th.main}, and using Proposition \ref{prop.pres_of_decay}, it follows that
\begin{equation}\label{eq.boltzmann_decay_4}
    \|f(t)\|_{\mathrm{Bol}^\beta}+\|\nabla_xf(t)\|_{\mathrm{Bol}^\beta}\leq B<\infty,\quad \forall t\in[0,t_{\mathrm{fin}}],
\end{equation} where $B$ is a constant depending only on $(d,t_{\mathrm{fin}},\beta,A,B_0)$. In the proof below we will always use (\ref{eq.boltzmann_decay_4}) instead of (\ref{eq.boltzmann_decay_1}) and (\ref{eq.boltzmann_decay_3}). In this section we fix some parameters and define some notations, which will be used throughout the proof.
\begin{definition}[Parameters and notations]
\label{def.notation} We define the following parameters and notations.
\begin{enumerate}
\item\label{it.time_layer} \emph{Time layering}. Recall that $(d,t_{\mathrm{fin}},\beta,A,B_0)$ has been fixed as in Theorem \ref{th.main}, and $B$ has been fixed as in (\ref{eq.boltzmann_decay_4}) depending only on these quantities. Throughout this paper we will assume $t=t_{\mathrm{fin}}$ (otherwise replace $t_{\mathrm{fin}}$ by $t$). Fix a large positive integer $\Lf$ depending on $(d,t_{\mathrm{fin}},\beta,A,B)$, and \textbf{subdivide the time interval} $[0,t_{\mathrm{fin}}]$ as follows:
\begin{equation}\label{eq.def_short_time}
\tau:=\frac{t_{\mathrm{fin}}}{\Lf};\qquad [0,t_{\mathrm{fin}}]=\bigcup_{\ell=1}^\Lf [(\ell-1)\tau,\ell\tau].
\end{equation}
\item\label{it.defC*} \emph{Convention for constants}. Throughout the proof, we will distinguish \textbf{two types of (large) constants}: those depending only on $(d,t_{\mathrm{fin}},\beta,A,B)$ and not on $\Lf$ will be denoted by $C$, while those depending on $\Lf$ will be denoted by $C^*$. The notions of $\lesssim$ and $\lesssim_*$ etc. are defined accordingly. We assume $\varepsilon$ is small enough depending on $C^*$.

We also fix an \textbf{explicit sequence of large constants} $(C_1^*,C_2^*,\cdots, C_{15}^*)$ such that \footnote{In practice, only in the layer selection part (cf. Definition \ref{def.layer_select}) do we need growths like $C_8^*\gg (C_7^*)^{10\Lf}$ and $C_6^*\gg (C_5^*)^{10\Gamma\Lf}$; for other $j$ we can choose (say) $C_{j+1}^*=(C_j^*)^{100}$.}
 \begin{equation}\label{eq.defC*}C_{j+1}^*\gg(\textrm{a suitable function of }C_j^*)\end{equation} all depending on $\Lf$. We understand that (i) for any large constant $C$ under the above convention, we always assume $\Lf^C\ll C_j^*$ for each $j$, and (ii) any large constant which we denote by $C^*$ will always be $\gg C_j^*$ for each $j$.
\item\label{it.defnot3} \emph{Other parameters}. \textbf{We fix} $\Gamma=\Gamma(d)$ sufficiently large but depending only on $d$, and fix a \textbf{sequence of parameters} $(A_\ell,\Lambda_\ell)\,(0\leq \ell\leq\Lf)$ by
\begin{equation}\label{eq.defLambdaseq}
A_{\Lf}=|\log\varepsilon|;\quad \Lambda_{\ell}=A_{\ell}^{10d},\quad A_{\ell-1}=\Lambda_{\ell}^{10d}.
\end{equation} They will be used in the truncation of correlation functions (Section \ref{sec.trunc_cor}). It is obvious that \begin{equation}\label{eq.ineq_Al}A_\ell,\Lambda_\ell\leq |\log\varepsilon|^{C^*},\quad \Lambda_{\ell-1}\gg A_{\ell-1}\gg \Lambda_\ell^2A_\ell.\end{equation} We also fix \textbf{two sequences} (which will be used in the analysis of $f^\Ac$ in Section \ref{sec.fa})
\begin{equation}
\label{eq.decayseq}\beta_\ell:=\bigg(\frac{9}{10}+\frac{\Lf-\ell}{10\Lf}\bigg)\frac{\beta}{2};\quad \theta_0:=3^{-d-8},\quad \theta_{\ell}:=\bigg(\frac{9}{10}+\frac{\Lf-\ell}{10\Lf}\bigg)\theta_0
\end{equation} for $1\leq\ell\leq\Lf$, and set $\theta=\theta_0/2$ in Theorem \ref{th.main}. For the proof in Sections \ref{sec.treat_integral}--\ref{sec.error},we also fix $\upsilon:=3^{-d-1}$.
\item\label{it.vector} \emph{Set, vector and integral notations}. For any positive integers $a$ and $b$, define \[[a:b]=\{a,a+1,\cdots,b\},\quad[a]:=[1:a].\] For any set $A$, as already done in Section \ref{sec.setup}, we shall adopt the \textbf{vector notation} $\vz_A:=(z_j)_{j\in A}$, and $\vz_{a:b}:=\vz_{[a:b]}$ and $\vz_a:=\vz_{[a]}$, and $\mathrm{d}\vz_A=\prod_{j\in A}\mathrm{d}z_j$. The same notion applies for other variables, such as $\boldsymbol{t}_A=(t_j)_{j\in A}$ etc. Define also the \textbf{tensor product} by $f^{\otimes A}(\vz_A)=\prod_{j\in A}f(z_j)$. The notion $\vz$ without any subscripts will be used to denote some $\vz_N$ where $N$ depends on the context.
\item \emph{Indicator functions and linear operators}. In the proof we will need to use a lot of \textbf{indicator functions}. These functions will be denoted by $\mathbbm{1}_{(\cdots)}$, which are defined by certain properties represented by $(\cdots)$ depending on the context, and take values in $\{0,1\}$. Define similarly the functions $\chi_{(\cdots)}$, except that those may take values in $\{-1,0,1\}$. Moreover, in Definition \ref{def.associated_op_nonlocal} (and several other places related to initial cumulants) we will define the notion $\mathbbm{1}_\Lc^{\varepsilon}$, which equals a positive power of $\varepsilon$ at points where the usual indicator function would be $0$. In practice they have the same behavior as usual indicator functions, and serve to simplify notations by unifying certain cases.

We will view any function $F$ (including the above indicator functions) as a \textbf{linear operator} defined by multiplication. For any linear operators $\Sc_1$ and $\Sc_2$ we define $\Sc_1\circ\Sc_2$ as their composition. For any linear operator $\Sc$ and function $F$, the notion $\Sc\circ F$ means the composition of $\Sc$ with the multiplication operator defined by $F$. This is different from the result of applying $\Sc$ to the function $F$, which we denote by $\Sc (F)$.
\item\label{it.sum_over_equiv} \emph{Summation over equivalence classes}. Let $\sim$ be an equivalence class on a set $X$, define $X / \sim$ be the set of equivalence classes, with the equivalent class containing $x\in X$ denoted by $[x]$. Let $\mathrm{Fun}$ be a function on $X$; if it is constant on each equivalence class, then the sum
 \[\sum_{[x] \in X / \sim} \mathrm{Fun}(x)\] is well-defined by taking $x$ as an arbitrary representative from each $[x]$. Otherwise, we define \[\sum_{[x]\in X/\sim}\mathrm{Fun}(x) := \sum_{[x]\in X/\sim} \frac{1}{|[x]|}\sum_{x'\in [x]}\mathrm{Fun}(x'),\] where $|[x]|$ is the cardinality of $[x]$. We remark that the case when $\mathrm{Fun}$ is not constant on each equivalence class, occurs only in Proposition \ref{prop.molecule_representation} and its proof.
\item\label{it.misc}\emph{Miscellaneous}. We use the notion $|A|$ or $\# A$ to denote the cardinality of finite sets $A$, and use $|A|$ to denote the Lebesgue measure of Borel sets $A\subseteq\Rb^p$ for some $p$. Define $X_+=\max(X,0)$ and $X_-=-\min(X,0)$. In addition to the standard $O(\cdot)$ notation, we also use $O_1(A)$ to denote quantities $X$ that satisfy the exact inequality $|X|\leq A$.
\end{enumerate}
\end{definition}

\subsubsection*{Acknowledgements} The authors are partly supported by a Simons Collaboration Grant on Wave Turbulence. The first author is supported in part by NSF grant DMS-2531437 and a Sloan Fellowship. The second author is supported in part by NSF grants DMS-1936640 and DMS-2350242. The authors would like to thank Thierry Bodineau, Isabelle Gallagher, Ning Jiang, Toan Nguyen, Mario Pulvirenti, Chenjiayue Qi and Laure Saint-Raymond for helpful discussions (including suggestions for the revision), and Yan Guo, Lingbing He, Clement Mouhot, Luis Silvestre, Sergio Simonella, Herbert Spohn, Nikolay Tzvetkov, and Huafei Yan for helpful discussions.

\section{A detailed overview of the proof}\label{sec.main_idea}
In this section, we give a more detailed sketch of the proof of Theorem \ref{th.main}. As explained in Section \ref{sec.strategy}, we will divide $[0, t_{\mathrm{fin}}]$ into smaller intervals $[(\ell-1)\tau,\ell\tau]$ where $\tau\ll 1$, see (\ref{eq.def_short_time}). Recall also the notion $\vz_A:=(z_j)_{j\in A}$ and $f^{\otimes A}$ in Definition \ref{def.notation} (\ref{it.vector}).

\subsection{The cumulant expansion formula}\label{subsec.idea.cumulants} As discussed in  Section \ref{sec.strategy}, the first step of the proof is to develop the cumulant expansion formulas \eqref{eq.cumulant_intro} and \eqref{eq.molecule_representation_intro}. This is inspired by the previous works \cite{PS17, BGSS22_2}, but with the key new idea of multi-layer, partial time expansion (see Section \ref{subsec.cumansatz_intro}) and several extra twists.

\subsubsection{Statement of the formula}\label{sec.statement_formula_idea} We start by stating the  precise version of \eqref{eq.cumulant_intro} and \eqref{eq.molecule_representation_intro} as follows:
\begin{equation}\label{eq.cumulant_idea}
\begin{gathered}
    f(\ell\tau, \vz_s) = \sum_{H \subseteq [s]} 
    \left[f^{\Ac}(\ell\tau, \cdot)\right]^{\otimes([s] \setminus H)} 
    \cdot E_H(\ell\tau, \vz_H) + \mathrm{Err}(\ell\tau, \vz_s), \\
    |E_H(\ell\tau)| \le  \sum_{\Mb \in \mathcal{F},\, r(\Mb)=H} |\mathcal{I}\mathcal{N}_{\Mb}|.
\end{gathered}
\end{equation}
Here $\mathrm{Err}(\ell\tau, \vz_s)$, omitted in \eqref{eq.cumulant_intro}, denotes an extremely small error term, which we will neglect in the remainder of the analysis. The sets $\mathcal{F}$, $r(\Mb)$ and the quantity $|\mathcal{I}\mathcal{N}_{\Mb}|$ will be defined and discussed below.

As explained in Section \ref{sec.molecule_representation_intro}, a \textbf{molecule} $\Mb$ (Definition \ref{def.molecule}) is a diagram, illustrated in {\color{blue}Figure \ref{fig.molecule}} (Section \ref{sec.molecule}), that encodes the combinatorial structure of the collision history. It is obtained by disregarding the precise spatial and temporal details of particle trajectories, retaining only the topological pattern of interactions, in a process called \textbf{topological reduction} (Definition \ref{def.top_reduction}). It consists of several key components listed below, which are illustrated in {\color{blue}Figure \ref{fig.molecule}}:
\begin{itemize}
    \item \textbf{C/O atoms}, which represent individual collisions/overlaps, are depicted as nodes like $\circ$/$\bullet$.
    \item \textbf{Edges}, which represent segments of free transport between collisions along a particle trajectory, are depicted as edges between nodes.
    \item \textbf{Particle lines}, which represent the full trajectory of each particle, are visualized in {\color{blue}Figure \ref{fig.molecule}} as chains of edges sharing the same color. They will be identified with particles in Sections \ref{sec.formula_cumulant}--\ref{sec.proof_cumulant_formula}. Each particle line has two \textbf{ends}, corresponding to the initial and final trajectories of the particle.
    \item \textbf{Layers} are the subsets $\Mb_{\ell'}\subseteq\Mb$ in {\color{blue}Figure \ref{fig.molecule}} between the horizontal lines representing times $\ell'\tau$ and $(\ell'+1)\tau$. Dividing $\Mb$ into layers $\Mb_{\ell'}$ corresponds to assigning a layer $\ell'$ to each atom $\nf$ which corresponds to a collision that happens in the time interval $[\ell' \tau, (\ell'+1)\tau]$.
    \item Each particle line $\pb$ has bottom end $e$, and a corresponding \textbf{start layer} $\ell_1[\pb]=\ell_1[e]$; this indicates the time $(\ell_1[\pb]-1)\tau$ at which the particle starts interacting with the system. Similarly, it also has a top end $e'$ and a \textbf{finish} layer $\ell_2[\pb]=\ell_2[e']$, indicating the time $\ell_2[\pb]\cdot\tau$ at which the particle stops interacting with the system.
    \item \textbf{Root particle lines} are the particle lines in $H$ (occurring in the cumulant $E_H$), i.e. those penetrating the top horizontal line $\ell\tau$ in {\color{blue}Figure \ref{fig.molecule}}.
\end{itemize}

We now explain the definitions of $|\mathcal{I}\mathcal{N}_{\Mb}|$ and $\mathcal{F}$ in \eqref{eq.cumulant_idea}. The expression $|\mathcal{I}\mathcal{N}_{\Mb}|$ is given by
\begin{equation}\label{eq.associated_integral_molecule_abs_idea}
    |\mathcal{I}\mathcal{N}_{\Mb}|(\vz_{r(\Mb)}) \coloneqq \varepsilon^{-(d-1)|\Mb|_{p \backslash r}} 
    \int_{\mathbb{R}^{2d|\Mb|_{p \backslash r}}} 
    (\mathcal{S}_{\Mb}\circ\mathbbm{1}_{\Mb}) \bigg(\mathbbm{1}_{\mathcal{L}}^\varepsilon\cdot
    \prod_{\pb \in p(\Mb)} \left| f^{\Ac}((\ell_1[\pb]-1)\tau) \right| \bigg)
    \mathrm{d}\vz_{(p \backslash r)(\Mb)},
\end{equation}
where the notations are explained as follows:

\begin{itemize}
    \item $r(\Mb)$, $p(\Mb)$, and $(p \backslash r)(\Mb)$ denote the sets of root particle lines, all particle lines, and non-root particle lines respectively, and $|\Mb|_{p \backslash r} \coloneqq |(p \backslash r)(\Mb)|$. Each particle in $p(\Mb)$ is connected (via collisions and overlaps) to some root particle in $r(\Mb)$. Moreover, $\ell_1[\pb]$ is the start layer of $\pb$ as above, i.e. the first time that $\pb$ is involved in the system.

    \item $\mathcal{S}_{\Mb}$ denotes the dynamics involving only the particles in $p(\Mb)$. The indicator function $\mathbbm{1}_{\Mb}$ restricts the integration to trajectories that topologically reduce to the molecule $\Mb$. The term $\mathbbm{1}_{\mathcal{L}}^\varepsilon$ enforces additional constraints that certain particles are initially close to each other, which comes from initial cumulants.

    \item The expression in \eqref{eq.associated_integral_molecule_abs_idea} expresses $|\mathcal{I}\mathcal{N}_{\Mb}|$ as the renormalized probability density contributed by the collision history encoded by $\Mb$ (with additional information contained in $\mathbbm{1}_{\mathcal{L}}$). Here $\vz_{(p\backslash r)(\Mb)}$ represent the final states of non-root particles, and integration in them leaves the result being a density function of only the root variables $\vz_{r(\Mb)}$. 
    
The input factors $f^\Ac ((\ell_1[\pb]-1)\tau)$ come from the start time $(\ell_1[\pb]-1)\tau$ where particle $\pb$ enters the system; their variables represent the initial states of particles, which are related to the variables $\vz_{p(\Mb)}$ by the transformation $\Sc_\Mb$ (strictly speaking, there should be a linear transport in the variables of $f^\Ac$ which transforms between time $0$ and $(\ell_1[\pb]-1)\tau$, see (\ref{eq.extraQ})).
\end{itemize}

The sum in \eqref{eq.cumulant_idea} is taken over molecules $\Mb \in \mathcal{F}$, which is defined as follows. Recall $\Mb_{\ell'}\subseteq\Mb$ is formed by atoms with layer $\ell'$. Define $r(\Mb_{\ell'})$ to be the set of particle lines $\pb$ that cross into layer $\ell'$ from higher layers (i.e. when the start layer $\ell_1[\pb]\leq\ell'$, and the finish layer $\ell_2[\pb]>\ell'$), which is also the the root particle set of $\Mb_{\ell'}$. Then, roughly speaking, the set $\mathcal{F}$ consists of molecules $\Mb$  that satisfy the following conditions:
\begin{enumerate}
    \item\label{it.set_F_l_1_intro} Every connected component of $\Mb_{\ell'}$ must contain at least one particle line from $r(\Mb_{\ell'})$.

    \item\label{it.set_F_l_2_intro} The number of atoms in $\Mb_{\ell'}$ is bounded above by $|\log \varepsilon|^{O(1)}$.

    \item\label{it.set_F_l_3_intro} Each layer $\Mb_{\ell'}$ is ``almost" a forest; i.e. we have proper control on its number of cycles.

    \item\label{it.set_F_l_4_intro} For all $1 \leq \ell' \leq \ell$, every particle line in $r(\Mb_{\ell'})$ is connected either to another particle line in $r(\Mb_{\ell'})$, or to a particle line in $r(\Mb_{\ell'-1})$, via edges within $\Mb_{\ell'}$.
\end{enumerate}

In the above conditions, (\ref{it.set_F_l_1_intro}) states that in layer $\Mb_{\ell'}$, we only track the dynamics of those particles that are connected to the root particles in $r(\Mb_{\ell'})$, through a chain of collisions and overlaps within $\Mb_{\ell'}$. Next,  (\ref{it.set_F_l_2_intro})--(\ref{it.set_F_l_3_intro}) come from the truncation of dynamics (described in Section \ref{sec.counting_molecule_intro}), which limits the number of particles and recollisions in each layer. 

Finally, \eqref{it.set_F_l_4_intro} states that every particle in $r(\Mb_{\ell'})$ is either connected (via collisions and overlaps within $\Mb_{\ell'}$) to another particle line in $r(\Mb_{\ell'})$, or connected to a particle that occurs in the cumulant $E_H$ at time $(\ell'-1)\tau$ (i.e. crosses into layers $<\ell'$ from layer $\ell'$ and belongs to $r(\Mb_{\ell'-1})$). This is because we are consideing the collision histories that contribute to the cumulant $E_H$, and because of our partial expansion (Section \ref{subsec.cumansatz_intro}), where only the cumulants $E_H((\ell'-1)\tau)$ are further expanded into previous time layers.

The precise definitions of $|\mathcal{I}\mathcal{N}_{\Mb}|$ and $\mathcal{F}$ will be given in Definitions \ref{def.associated_int} and \ref{def.set_T_F}, respectively. The cumulant expansion formula \eqref{eq.cumulant_idea} will be formally stated in Proposition \ref{prop.cumulant_formula}. The conditions (\ref{it.set_F_l_1_intro})--(\ref{it.set_F_l_4_intro}) listed above correspond, with minor modifications, to those in Definition \ref{def.set_T_F} \eqref{it.set_F_l_1}--\eqref{it.set_F_l_4}.

\subsubsection{Proof of the formula}\label{sec.proof_formula_idea} In this section, we sketch the proof of the cumulant formula \eqref{eq.cumulant_idea}. By \eqref{eq.time_t_ensemble} and \eqref{eq.s_par_cor}, we have
\begin{equation}\label{eq.proof_formula_idea_1}
\begin{gathered}
    f(\ell\tau,\vz_s) = \varepsilon^{(d-1)s} \sum_{n=0}^\infty \frac{1}{n!} \int_{\Rb^{2dn}} W_{s+n}(\ell\tau,\vz_{s+n})\,\mathrm{d}z_{s+1} \cdots \mathrm{d}z_{s+n}, \\
    W_{N}(\ell\tau) = \Sc_{N}(\tau)(W_{N}((\ell-1)\tau)).
\end{gathered}
\end{equation}

\emph{Prerequisite: The truncated dynamics.} Before starting, we note that, in order that the molecule $\Mb$ satisfies (\ref{it.set_F_l_1_intro})--(\ref{it.set_F_l_3_intro}) in Section \ref{sec.statement_formula_idea} (in particular, $\Mb_{\ell'}$ being almost a forest means we need to control the number of recollisions), we need to perform a truncation on the dynamics (see Definition \ref{def.T_dynamics}) within each layer.

In practice, the truncation is done by monitoring each potential collision event. When a collision is about to occur (i.e. a pre-collisional configuration forms), we check whether it would violate the constraints in \eqref{it.set_F_l_2_intro} or \eqref{it.set_F_l_3_intro}. If such a violation would happen, we ``turn-off'' the collision by allowing the particles to pass through one another without a collision; otherwise, we allow the collision to happen as in the original dynamics.

\emph{Step 1: The induction setup.} The general strategy is to verify \eqref{eq.cumulant_idea} by induction. For the base case $\ell=0$, we rely on the exact factorization at $t=0$ in (\ref{eq.N_par_ensemble}), except the factor $\mathbbm 1_{\Dc_N}$. This latter factor leads to the initial cumulants, which is the reason of the factors $\mathbbm{1}_\Lc^{\varepsilon}$ in (\ref{eq.associated_integral_molecule_abs_idea}). These initial cumulants can be treated in the standard way (see Proposition \ref{prop.initial_cumulant}), and will be ignored in this section.

Then, we need the induction step, where we derive \eqref{eq.cumulant_idea} for $\ell$ from the same formula for $\ell-1$. In order to identify the factorized and cumulant parts at time $\ell\tau$ based on those at $(\ell-1)\tau$, we need to perform the \textbf{cluster expansion}, which is detailed in \emph{Step 2} below. Here, the clusters associated with the trajectories on $[(\ell-1)\tau,\ell\tau]$ (and the corresponding transport operators $\Sc_N(\tau)$), are defined 
as the maximal subsets of particles that are connected by collisions within $[(\ell-1)\tau,\ell\tau]$. Denote by $\Mb$ the molecule (without O-atoms) which is the result of topological reduction of the trajectory on $[(\ell-1)\tau,\ell\tau]$, and let $\Mb_j$ be the connected components of $\Mb$, corresponding to each cluster.

\emph{Step 2: The cluster expansion.} Let $p(\Mb)$ be the set of all particles in $\Mb$ etc. Let $\mathbbm{1}_{\Mb}$ be the indicator function that the trajectory of the particles in $p(\Mb)$  topologically reduces to $\Mb$, and let $\mathbbm{1}_{\Mb_j \not\sim \Mb_{j'}}$ be the indicator function that the particles in $p(\Mb_j)$ never collide with those in $p(\Mb_{j'})$. Then we have the following cluster expansion identity:
\begin{equation}\label{eq.POU}
1 = \sum_{k \leq N} \sum_{\{\Mb_1, \ldots, \Mb_k\}} \bigg( \prod_{1 \leq j < j' \leq k} \mathbbm{1}_{\Mb_j \not\sim \Mb_{j'}} \bigg) \bigg( \prod_{j=1}^{k} \mathbbm{1}_{\Mb_j} \bigg),
\end{equation}
which gives a partition of unity over all non-overlapping cluster decompositions.

Correspondingly, we have the following cluster expansion formula for the dynamics $\Sc_N(\tau)$:
\begin{equation}\label{eq.W_Nexp}
\Sc_{N}(\tau) = \sum_{k \leq N} \sum_{\{\Mb_1, \ldots, \Mb_k\}} \bigg( \prod_{j=1}^{k} \Sc_{\Mb_j}(\tau) \cdot \mathbbm{1}_{\Mb_j} \bigg)\bigg( \prod_{1 \leq j < j' \leq k} \mathbbm{1}_{\Mb_j \not\sim \Mb_{j'}} \bigg),
\end{equation}
where $\Sc_{\Mb_j}(\tau)$ denotes the dynamics restricted to each cluster. By applying (\ref{eq.W_Nexp}) to $W_N((\ell-1)\tau)$ and using the first equation in (\ref{eq.proof_formula_idea_1}), this leads to an expression of $f_s(\ell\tau)$ in terms of cluster expansions and $f_s((\ell-1)\tau)$, which in particular contains the indicator functions $\prod \mathbbm{1}_{\Mb_j \not\sim \Mb_{j'}}$.

\emph{Step 3: Penrose argument and overlaps.} The indicator function $\mathbbm{1}_{\Mb_j \not\sim \Mb_{j'}}$ in (\ref{eq.W_Nexp}) indicates there is no overlap between the trajectories of particles in $p(\Mb_j)$ and $p(\Mb_{j'})$, which implies the further factorization $\mathbbm{1}_{\Mb_j \not\sim \Mb_{j'}}=\prod_{(e,e')}\mathbbm{1}_{(\Mb_j,e) \not\sim (\Mb_{j'},e')}$, where $(e,e')$ runs over pairs of edges $e\in\Mb_j$ and $e'\in\Mb_{j'}$, and $\mathbbm{1}_{(\Mb_j,e) \not\sim (\Mb_{j'},e')}$ indicates that (the trajectory piece corresponding to) the edge $e$ and the edge $e'$ do not overlap.

Then, using the idendity $\mathbbm{1}_{(\Mb_j,e) \not\sim (\Mb_{j'},e')}=1-\mathbbm{1}_{(\Mb_j,e)\sim (\Mb_{j'},e')}$, we can expand the product $\prod \mathbbm{1}_{(\Mb_j,e) \not\sim (\Mb_{j'},e')}$ into a summation of terms, using the standard \emph{Penrose argument} similar to those in \cite{BGSS20}. Here, compared to \cite{BGSS20}, we need to perform an extra truncation in the Penrose argument (in order to bound the number of clusters by $\Lambda=|\log\varepsilon|^C$), which leads to extra complications.

In the end, we can recast (\ref{eq.W_Nexp}) as a summation over single layer molecules $\Mb$ which satisfy (\ref{it.set_F_l_1_intro})--(\ref{it.set_F_l_3_intro}) in Section \ref{sec.statement_formula_idea}. This $\Mb$ is formed by joining the clusters $\Mb_j$ (which contain only C-atoms) by O-atoms, which correspond to \textbf{overlaps} between edges $e\in\Mb_j$ and $e'\in\Mb_{j'}$, which come from the indicator functions $\mathbbm{1}_{(\Mb_j,e)\sim (\Mb_{j'},e')}$. Then, once $\Mb$ is fixed, we can insert the induction hypothesis \eqref{eq.cumulant_idea} for $(\ell-1)\tau$, and factorize the resulting term, to get an expression of form \eqref{eq.cumulant_idea} for $\ell\tau$, thus completing the inductive step.

\emph{Step 4: Formation of multi-layer molecules.} With \emph{Steps 1--3} above, we can then induct on $\ell$ to prove the first line in \eqref{eq.cumulant_idea}, as well as an expression (upper bound) of $|E_H(\ell\tau)|$ in terms of summation over sequences $(\Mb_{\ell},\cdots,\Mb_1)$ of molecules. Then, by merging the corresponding top and bottom ends for the same particle in each different $\Mb_{\ell'}$, we can construct a multi-layer molecule $\Mb$ from the sequence $(\Mb_{\ell'})$, so the summation over this sequence is rewritten as a summation over multi-layered molecule $\Mb$, which leads to the second line in \eqref{eq.cumulant_idea} and completes the proof.

The prerequisite, i.e. truncation of dynamics, will be discussed in Section \ref{sec.truncation_large_molecule}. Then, \emph{Steps 1--4} will be carried out in Section \ref{sec.proof_cumulant_formula}. The cluster expansion in \emph{Step 2} is contained in Section \ref{sec.cluster_expansion} (Lemma \ref{lem.dynamics_factorization}), and the Penrose argument in \emph{Step 2} is contained in Section \ref{sec.penrose} (Proposition \ref{prop.S_N_decomposition}). Then we complete the inductive step in Section \ref{sec.f_to_f}--\ref{sec.E_to_E} (Proposition \ref{prop.cumulant_single}), and the base case (initial cumulant) in Section \ref{sec.initial_cumulant} (Proposition \ref{prop.initial_cumulant}). Finally, in Section \ref{sec.proof_cumulant_formula_final}, we carry out \emph{Step 4} and put everything together to complete the proof.

\subsection{Reduction to local integral expressions}\label{subsec.idea.molecule} Now we have proved \eqref{eq.cumulant_intro} and \eqref{eq.molecule_representation_intro}. As discussed in Section \ref{sec.counting_molecule_intro}, in order to prove Theorem \ref{th.main}, it suffices to bound $\|E_H\|_{L^1}$. In this subsection, we apply \eqref{eq.cumulant_idea} to reduce the estimate for $\|E_H\|_{L^1}$ to individual estimates for a certain quantity $\Ic_\Mb(Q_\Mb)$ given by local integral expressions, which is equivalent to $\big\||\Ic\Nc_\Mb|\big\|_{L^1}$. This is done in two steps.
\subsubsection{Counting the number of molecules} The first step is to prove the upper bound for the number of molecules, i.e. \eqref{eq.molecule_number_upp_intro}, where the quantity $\rho$ is defined by $\rho=\sum_{\ell'} (s_{\ell'} + \Rf_{\ell'})$ \footnote{This $\rho$ is equivalent (up to multiplicative constants) to the number of recollisions defined in Section \ref{sec.counting_molecule_intro}.}. Here, $\Rf_{\ell'}$ is the number of the independent cycles (i.e., the circuit rank) in $\Mb_{\ell'}$, and $s_{\ell'} = |r(\Mb_{\ell'})|$ is the number of particle lines crossing into layer $\ell'$ from higher layers.

The proof of \eqref{eq.molecule_number_upp_intro} is contained in Proposition \ref{prop.layerrec3} in Section \ref{sec.int_dia_exp}. By property \eqref{it.set_F_l_2_intro} in Section \ref{sec.statement_formula_idea}, we know that the number of atoms in each layer $\Mb_{\ell'}$ can be bounded by $|\log \varepsilon|^{C}$. Then \eqref{eq.molecule_number_upp_intro} follows from the combinatorial fact that the number of graphs with all degrees $\leq 4$, at most $m$ vertices, and at most $n$ independent cycles, is bounded by $C^m \cdot m^n$ for some constant $C > 0$. Applying this to molecules yields:
\[
\begin{aligned}
    \#\Mb \le& \prod_{\ell'} \#\Mb_{\ell'}\cdot \#\{\textrm{choices of $s_{\ell'}$ particle lines going to the layer below}\}
    \\
    \le& \prod_{\ell'} C^{|\Mb_{\ell'}|} \cdot |\Mb_{\ell'}|^{\Rf_{\ell'}} \cdot |\log \varepsilon|^{Cs_{\ell'}} \le C^{|\Mb|} \cdot |\log \varepsilon|^{C \rho},
\end{aligned}
\]
which proves \eqref{eq.molecule_number_upp_intro}. Here we used the fact that the number of particle lines in $\Mb_{\ell'}$ is $\le |\log \varepsilon|^{O(1)}$, and $\#\{\textrm{choices of $s_{\ell'}$ particle lines going to the layer below}\}\le |\log \varepsilon|^{Cs_{\ell'}}$. 

\subsubsection{Local representation formula}\label{sec.local_rep_idea} It now suffices to prove (\ref{eq.molecule_term_upp_intro}). In this subsection, we reduce the quantity $\big\||\Ic\Nc_\Mb|\big\|_{L^1}$ to another quantity $\Ic_\Mb(Q_\Mb)$ which is given by a local integral expression of $\dirac$ functions. This $\Ic_\Mb(Q_\Mb)$ is defined as follows.

Consider a molecule $\Mb$. We define the \textbf{associated variables} $z_e$ and $t_\nf$ for each edge $e$ and atom $\nf$, forming $\vz_\Ec = (z_e : e \in \Ec)$ and $\vt_\Mb = (t_\nf : \nf \in \Mb)$, where $\Ec$ is the edge set of $\Mb$ (and $\Mb$ viewed as the atom set of itself). For a C-atom $\nf \in \Mb$ (there is a similar and simpler definition for O-atoms), let $(e_1, e_2)$ and $(e_1', e_2')$ denote its bottom and top edges respectively. We define the \newterm{associated distribution} $\Dirac_\nf = \boldsymbol{\Delta} = \boldsymbol{\Delta}(z_{e_1}, z_{e_2}, z_{e_1'}, z_{e_2'}, t_\nf)$ by
\begin{equation}\label{eq.associated_dist_C_idea}
\begin{aligned}
    \boldsymbol{\Delta}(z_{e_1}, z_{e_2}, z_{e_1'}, z_{e_2'}, t_\nf) &:= 
    \boldsymbol{\delta}\big(x_{e_1'} - x_{e_1} + t_\nf (v_{e_1'} - v_{e_1})\big) \cdot 
    \boldsymbol{\delta}\big(x_{e_2'} - x_{e_2} + t_\nf (v_{e_2'} - v_{e_2})\big) \\
    &\quad \cdot \boldsymbol{\delta}\big(|x_{e_1} - x_{e_2} + t_\nf (v_{e_1} - v_{e_2})| - \varepsilon\big) \cdot 
    \left[(v_{e_1} - v_{e_2}) \cdot \omega\right]_- \\
    &\quad \cdot \boldsymbol{\delta}\big(v_{e_1'} - v_{e_1} + [(v_{e_1} - v_{e_2}) \cdot \omega] \omega\big) \cdot 
    \boldsymbol{\delta}\big(v_{e_2'} - v_{e_2} - [(v_{e_1} - v_{e_2}) \cdot \omega] \omega\big),
\end{aligned}
\end{equation}
where $\omega := \varepsilon^{-1}(x_{e_1} - x_{e_2} + t_\nf (v_{e_1} - v_{e_2}))$ is a unit vector, and $z_- := -\min(z, 0)$. On the support of $\Dirac_\nf$, if $(x_{e_1}, v_{e_1})$ and $(x_{e_2}, v_{e_2})$ represent the configurations of two particles immediately before the collision (pre-collisional configurations), then $(x_{e_1'}, v_{e_1'})$ and $(x_{e_2'}, v_{e_2'})$ correspond to their configurations immediately after the collision (post-collisional configurations).

We define the \textbf{associated domain} $\Dc$ by
\begin{equation}\label{eq.associated_domain_idea}
    \Dc := \left\{ \vt_\Mb = (t_\nf) \in (\mathbb{R}^+)^{|\Mb|} \,\middle|\, 
    (\ell' - 1)\tau < t_\nf < \ell'\tau \text{ if } \ell[\nf] = \ell', \quad 
    t_\nf < t_{\nf^+} \text{ if } \nf^+ \text{ is the parent of } \nf \right\}.
\end{equation}

Then, we define the \newterm{associated operator} $\Ic_\Mb$ by
\begin{equation}\label{eq.associated_int_op_idea}
    \Ic_\Mb(Q):=\varepsilon^{-(d-1)(|\Ec|-2|\Mb|)}\int_{\mathbb{R}^{2d|\Ec|}\times \Rb^{|\Mb|}}\bigg(\prod_{\nf\in\Mb}\boldsymbol{\Delta}_\nf\bigg)\cdot Q \,\mathrm{d}\vz_{\Ec}\,\mathrm{d}\vt_\Mb.
\end{equation}

The main result of this subsection is then stated as follows: for any full molecule $\Mb$ with $r(\Mb) = H$, we have the following equality
\begin{equation}\label{eq.local_int_idea}
    \left\| |\Ic\Nc_\Mb|(\vz_H) \right\|_{L^1} = 
    \varepsilon^{(d-1)|H|} \cdot \Ic_\Mb(Q_\Mb),
\end{equation}
where $Q_\Mb$ is a bounded function with certain support properties, whose precise form is less important here.

The proof of (\ref{eq.local_int_idea}) is presented in Proposition \ref{prop.local_int} in Section \ref{sec.int_reduce}. The main point here is that the factors in the integral (\ref{eq.local_int_idea}) are \emph{local} i.e. each individual $\Dirac_\nf$ contains only those variables at atom $\nf$ and edges at $\nf$, and \emph{do not contain arbitrarily long chains of substitutions as in $\Ic\Nc_\Mb$}. This is very important, as it allows us to apply the \textbf{cutting algorithm} below.

\subsection{The cutting algorithm}\label{subsec.idea.alg} Recall that we need to prove (\ref{eq.molecule_term_upp_intro}). Using (\ref{eq.local_int_idea}), we only need to establish an upper bound for $\Ic_\Mb(Q_\Mb)$. This is a multiple integral; involving many integrated variables; in order to estimate it, we need to reduce it to an \emph{iterated integral} by choosing an appropriate order of variables in which we integrate. This process is implemented through the \textbf{cutting operation} and \textbf{cutting algorithm} (as sketched in Section \ref{sec.cutting_alg_intro}), which we now introduce.

\subsubsection{Cutting operations, and elementary molecules}\label{sec.cutting_op_idea} 
\begin{figure}[h!]
    \centering
     \includegraphics[width=0.4\linewidth]{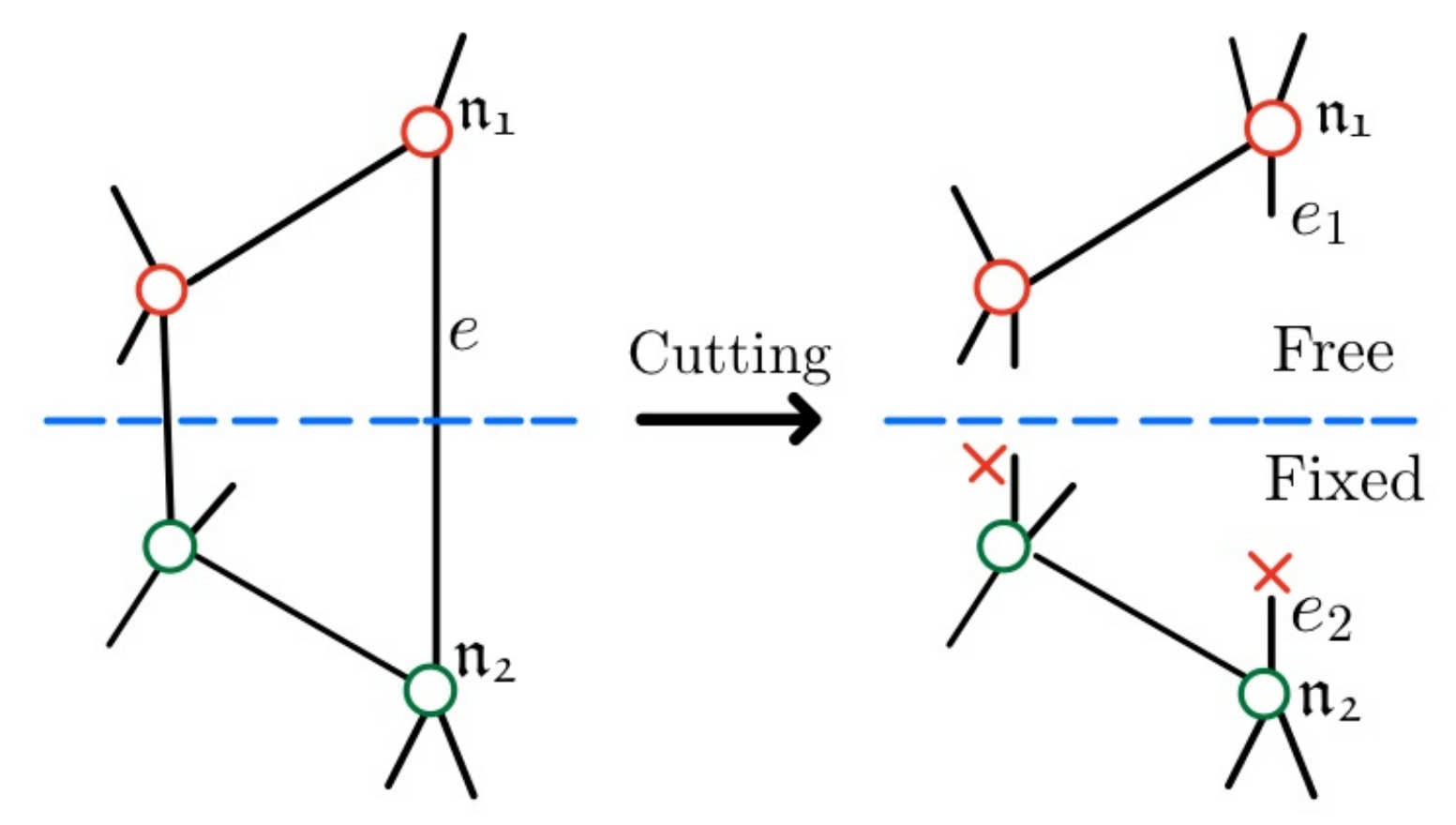}
    \caption{Cutting a molecule $\Mb$ into smaller molecules $\Mb_1$ and $\Mb_2$ (where red atoms in $\Mb_1$ are cut as free, and green atoms in $\Mb_2$ cut as fixed). The fixed ends are marked with red $\times$ symbols.}
    \label{fig.cutting}
\end{figure}

Given any molecule $\Mb$, \textbf{cutting} is the operation that divides $\Mb$ into two disjoint sub-molecules $\Mb_1$ and $\Mb_2$, where each bond connecting an atom $\nf_1\in\Mb_1$ to an atom $\nf_2\in \Mb_2$ is turned into a \textbf{free end} at $\nf_1$ and a \textbf{fixed end} at $\nf_2$. We say $\Mb_1$ is \textbf{cut as free} and $\Mb_2$ is \textbf{cut as fixed}. See {\color{blue}Figure \ref{fig.cutting}} for an illustration, and Definition \ref{def.cutting} for the precise definition. 

In terms of the integral in $\Ic_\Mb$, cutting $\Mb_1$ as free (and $\Mb_2=\Mb\backslash \Mb_1$ as fixed) simply means that we first integrate in the $(z_e)$ and $(t_\nf)$ variables associated with $\Mb_2$ (treating the variables associated with $\Mb_1$ as fixed) and then integrate in those variables associated with $\Mb_1$. Here when a bond $e$ between $\Mb_1$ and $\Mb_2$ is cut into a free end $e_1$ at $\Mb_1$ and a fixed end $e_2$ at $\Mb_2$, we always impose the condition $z_e=z_{e_1}=z_{e_2}$; see {\color{blue}Figure \ref{fig.cutting}}. This is merely a diagrammatic representation of Fubini's theorem.

With the above conventions, it is straightforward to verify that
\begin{equation} \label{eq.comp_intro}
    \Ic_\Mb = \Ic_{\Mb_1}\circ\Ic_{\Mb_2}    
\end{equation}
where in $\Ic_{\Mb_2}$, we are integrating in the variables of free ends and bonds of $\Mb_2$, treating the variables of fixed ends as frozen in the integral. In fact, (\ref{eq.comp_intro}) follows from the calculation\begin{equation}
\begin{split}
\Ic_\Mb(Q) & = \int_{\Rb^{2d|\Ec|}\times\Rb^{|\Mb|}}\bigg(\prod_{\nf\in\Mb}\boldsymbol{\Delta}_\nf\bigg)\cdot Q \,\mathrm{d}\vz_{\Ec}\mathrm{d}\vt_\Mb
\\
&=\underbrace{\int_{\Rb^{2d|\Ec_1|}\times\Rb^{|\Mb_1|}} \bigg(\prod_{\nf\in\Mb_1}\boldsymbol{\Delta}_\nf\bigg)}_{\Ic_{\Mb_1}}\cdot \underbrace{\bigg(\int_{\Rb^{2d|\Ec_2|}\times\Rb^{|\Mb_2|}} \bigg(\prod_{\nf\in\Mb_2}\boldsymbol{\Delta}_\nf\bigg)\cdot Q \,\mathrm{d}\vz_{\Ec_2}\mathrm{d}\vt_{\Mb_2}\bigg)}_{\Ic_{\Mb_2}(Q)}\,\mathrm{d}\vz_{\Ec_1}\mathrm{d}\vt_{\Mb_1}
\\
&=\Ic_{\Mb_1}\circ\Ic_{\Mb_2}(Q).
\end{split}
\end{equation}

In general, in (\ref{eq.comp_intro}), $\Mb$ is allowed to have fixed ends, and the definitions of $\Ic_\Mb$ (and $\Ic_{\Mb_1}$) in \eqref{eq.associated_int_op_idea} is modified by replacing $\Ec$ with $\Ec_* \coloneqq \{\textrm{bonds and free ends}\} $. The proof of (\ref{eq.comp_intro}) in full generality is presented in Proposition \ref{prop.cutting}.

In our cutting algorithm, presented in Section \ref{sec.intro_alg} below, we will repeatedly cut $\Mb$ into smaller molecules until all of them becomes \textbf{elementary molecules} $\Mb_j$, such that each of them contains only one or two atoms, so the integral $\Ic_{\Mb_j}$ can be calculated explicitly. We can then control $\Ic_\Mb(Q_\Mb)$ by $\Ic_\Mb(Q_\Mb)\leq \|Q_\Mb\|_{L^\infty}\cdot \prod_j \|\Ic_{\Mb_j}\|_{L^\infty\to L^\infty}$.

There are four major types of elementary molecules, denoted by \{2\}, \{3\}, \{4\} and \{33\}, indicated by the \textbf{degree} (i.e number of bonds plus free ends, \emph{not} counting fixed ends) of each atom. See {\color{blue}Figure \ref{fig.elementary}} (Section \ref{sec.oper_int}) for a picture of elementary molecules and Definition \ref{def.elementary} for its rigorous definition. Note that \{33\} molecules are divided into types \{33A\} and \{33B\}, and we also allow \{44\} molecules which play minor roles and will not be discussed here.

The key observation is that, for each \{33\} molecule $\Mb_j$ that occurs in the cutting process (see $\Mb_2$ in {\color{blue}Figure \ref{fig.cutting}}), we have the operator norm bound $\|\Ic_{\Mb_j}\|_{L^\infty\to L^\infty}\leq \varepsilon^{\upsilon}$ for some positive constant $\upsilon$, at least in typical cases. For this reason, we define each \{33\} molecule to be \textbf{good}, and each \{3\} and \{2\} molecule to be \textbf{normal}. In fact, a \{3\} molecule corresponds to a normal collision, a \{2\} molecule corresponds to a recollision, and a \{33\} molecule (which can be further decomposed to a \{3\}  and a \{2\} molecule) corresponds to a recollision that follows a normal collision, which constitutes the basic unit of gain, as already observed in \cite{DH21} in the context of wave turbulence. Finally, we define a \{4\} molecule to be \textbf{bad} as it ``wastes" one integral and leads to losses $\|\Ic_{\Mb_j}\|_{L^\infty\to L^\infty}\leq \varepsilon^{-(d-1)}$. See Definition \ref{def.good_normal} for the precise version of classification of elementary molecules.

Returning to the estimate of $\Ic_\Mb(Q_\Mb)$, suppose there exists a cutting sequence that decomposes a molecule $\Mb$ into a collection of molecules $(\Mb_j)$. Then, by the discussion above and the fact that $Q_\Mb$ is bounded (say $\|Q_\Mb\| \le 1$), we obtain
\begin{equation}\label{eq.apply_cutting_idea}
    \left\| |\Ic\Nc_\Mb|(\vz_H) \right\|_{L^1} = \varepsilon^{(d-1)|H|}\Ic_\Mb(Q_\Mb) 
    \leq \prod_j \|\Ic_{\Mb_j}\|_{L^\infty \to L^\infty}
    \le \varepsilon^{\upsilon \cdot \#_{\mathrm{good}} + (d-1) (|H|- \#_{\mathrm{bad}})},
\end{equation}
where we gain a factor of $\varepsilon^\upsilon$ for each good molecule and lose a factor of $\varepsilon^{-(d-1)}$ for each bad molecule. By (\ref{eq.apply_cutting_idea}), it is clear that (\ref{eq.molecule_term_upp_intro}) would follow if we can construct a cutting sequence satisfying
\begin{equation}\label{eq.comb_intro1}
    \upsilon \cdot \#_{\mathrm{good}} +(d-1) \cdot (|H|-\#_{\mathrm{bad}}) \geq c \cdot \rho.
\end{equation} This inequality is stated in Proposition \ref{prop.comb_est}, and is the main goal of the combinatorial part of this paper. We discuss the algorithm of constructing the cutting sequence, in Section \ref{sec.intro_alg} below.

The discussion and estimates of elementary molecules are contained in Section \ref{sec.elem_int}, see Propositions \ref{prop.intmini}--\ref{prop.intmini_2} and Definition \ref{def.good_normal}. The proof of (\ref{eq.molecule_term_upp_intro}) (as part of Proposition \ref{prop.cumulant_est}, see (\ref{eq.proof_int_est_from_alg_step0_*})) is contained in Section \ref{sec.summary}, using also the auxiliary results (Propositions \ref{prop.weight}--\ref{prop.volume}) in Section \ref{sec.aux}.

\subsubsection{The cutting algorithm}\label{sec.intro_alg} Finally we discuss the proof of (\ref{eq.comb_intro1}) stated as Proposition \ref{prop.comb_est}. This is the main technical ingredient of the proof, and will occupy Sections \ref{sec.toy}--\ref{sec.maincr}. The idea is to maximize the number of good (i.e. \{33\}) molecules and minimize the number of bad (i.e. \{4\}) molecules, and balance between them if needed. Below we discuss the main ideas and a few important points in this construction.

\textbf{Toy models.} The full algorithm, see Sections \ref{sec.layer}--\ref{sec.maincr}, is quite complicated. In Section \ref{sec.toy}, we will neglect a few technicalities by making the simplifying assumptions that (i)  
$\Mb$ has no O-atoms, and (ii) each layer $\Mb_{\ell'}$ is precisely a forest. Moreover, we define the \textbf{toy models} (I, I+, II and III, see Definitions \ref{def.toy}, \ref{def.toy1+}, \ref{def.toy2}, \ref{def.toy3} and {\color{blue}Figure \ref{fig.toy models}}) to be specific cases of \emph{two layer molecules} $\Mb$: it has an upper layer $\Mb_U$ and a lower layer $\Mb_D$, each of which is a tree. The $\rho$ in (\ref{eq.comb_intro1}) is comparable to the number of bonds between $\Mb_U$ and $\Mb_D$.

The proof in Section \ref{sec.toy} (under the simplifications (i)--(ii)) will consist of two parts: reducing general multi-layer molecules $\Mb$ to one of the toy models (via \emph{layer selection}), and designing a cutting algorithm for each individual toy model.

\textbf{The UP algorithm.} We start with a simple algorithm called \textbf{UP}, see Definitions \ref{def.up_toy}. Basically, in this algorithm we select a \emph{lowest degree 3} atom $\nf$, and cut the atoms in $S_\nf$ (the set of descendants of $\nf$) from high to low, see {\color{blue}Figure \ref{fig.up}} (Section \ref{sec.toy1+}).

We can show that all the molecules resulting from \textbf{UP} are elementary by preserving the key monotonicity property: each atom either has \emph{no top fixed end} or \emph{no bottom fixed end}, see Proposition \ref{prop.up_toy} (\ref{it.up_toypf_1}). Moreover, this algorithm \emph{minimizes} $\#_{\{4\}}$ (the number of \{4\} molecules) by allowing only one \{4\} molecule per connected component of $\Mb$. 

The drawback is that, the \textbf{UP} algorithm does \emph{not necessarily maximize} $\#_{\{33\}}$ (the number of \{33\} molecules), so it is not directly useful in our proof. However, it is still useful in ``finishing off" unimportant parts of the molecule where we do not expect any gain. Moreover, in some important cases, one can actually make special adaptations to \textbf{UP} to guarantee the abundance of \{33\} molecules; this leads to the important algorithms \textbf{2CONNUP} and \textbf{3COMPUP} below.

\textbf{Toy model I.} The feature of toy model I, see Definition \ref{def.toy}, is that each atom in $\Mb_U$ has exactly one bond with $\Mb_D$. In this case, the cutting algorithm is defined in Definition \ref{def.toy_alg}. It combines a trivial algorithm in $\Mb_U$ (always cutting lowest) with a greedy algorithm in $\Mb_D$ (cutting a \{33\} molecule, and the related \{343\} molecule, whenever possible). See {\color{blue}Figure \ref{fig.toy model 1 alg}} (Section \ref{sec.toy1}).

In this algorithm, one can show that the number of \{4\} molecules is minimized, i.e. $\#_{\{4\}}\leq 1$; moreover the greedy algorithm in $\Mb_D$ guarantees the abundance of \{33\} molecules, i.e. $\#_{\{33\}}\gtrsim\rho$, see Proposition \ref{prop.toy}. This relies on the fact that we never get any degree 2 atom in $\Mb_D$ (Lemma \ref{lem.toy_1}), which crucially uses that $\Mb_D$ is a forest.

\textbf{Toy model I+.} Toy model I+ (Definition \ref{def.toy1+}) is similar to toy model I, but only a subset $X\subseteq\Mb_U$ of atoms has one bond with $\Mb_D$, and the others have none. Moreover, in toy model I+ we assume the number of connected components of $X$ satisfies $\#_{\mathrm{comp}(X)}\ll\rho$.

In this case the algorithm (which is essentially the \textbf{MAINUD} algorithm) is defined in Definition \ref{def.toy1+_alg}, and is similar to Definition \ref{def.toy_alg}, see {\color{blue}Figures \ref{fig.toy model 1+ alg1}--\ref{fig.toy model 1+ alg2}} (Section \ref{sec.toy1+}). The only difference is that we only apply the naive ``cutting lowest" algorithm in $X$, and instead apply \textbf{UP} in the rest of $\Mb_U$, to minimize the number of \{4\} molecules. In fact this guarantees that $\#_{\{4\}}\lesssim\#_{\mathrm{comp}(X)}$, which is negligible compared to the gain $\#_{\{33\}}\gtrsim\rho$, see Proposition \ref{prop.toy1+_alg}.

\textbf{Toy model II.} In toy model II (Definition \ref{def.toy2}), we have the same assumptions as toy model I, but with $\#_{\mathrm{comp}(X)}\gtrsim\rho$. In this case the algorithm is defined in Definition \ref{def.toy2_alg}. Basically we first remove $\Mb_D$ and cut it using \textbf{DOWN} the dual of \textbf{UP}, and then cut $\Mb_U$ using a special variant of \textbf{UP} which is called the \textbf{3COMPUP} algorithm, see {\color{blue}Figure \ref{fig.toy model 2 alg}} (Section \ref{sec.toy2}).

The idea here is that, after removing $\Mb_D$, the degree 3 atoms in $\Mb_U$ are precisely those in $X$. Then, $\Mb_U$ will satisfy the key property that the degree 3 atoms are \emph{separated by degree 4 atoms}, i.e. the number of connected components of degree 3 atoms satisfies $\#_{\mathrm{3comp}(\Mb_U)}=\#_{\mathrm{comp}(X)}\gtrsim\rho$. Under this assumption, one can design a special variant of \textbf{UP} (basically choosing \emph{two} lowest deg 3 atoms $\nf$ instead of one), namely the \textbf{3COMPUP}, to guarantee the abundance of \{33\} molecules, see Proposition \ref{prop.toy2}.

\textbf{Toy model III.} In toy model III (Definition \ref{def.toy3}), we assume there exists $\gtrsim\rho$ atoms in $\Mb_U$, each having \emph{two bonds} with $\Mb_D$; they are called \emph{2-connections}. In this case, we simply apply the \textbf{DOWN} algorithm (dual of \textbf{UP} in $\Mb_D$), see Definition \ref{def.toy3_alg} and {\color{blue}Figure \ref{fig.toy model 3 alg}} (Section \ref{sec.toy3}). In the process, each 2-connection atom $\nf$ will naturally occur as a \{33\} molecule, which we then cut it to obtain a good component. Consequently we get $\#_{\{33\}}\gtrsim\rho$ and $\#_{\{4\}}=1$, see Proposition \ref{prop.toy3}.

\textbf{Reduction to toy models.} With the algorithm constructed for each toy model, it now remains to reduce general multi-layer molecules $\Mb$ (under the assumptions (i)--(ii)) to one of the toy models. The first step is to reduce $\Mb$ to a \emph{two-layer} molecule by a process called \textbf{layer selection}.

Since we have ignored the recollisions \emph{within each layer} $\Mb_{\ell'}$ by assuming each layer is a forest, we know that the recollisions in $\Mb$ must happen \emph{across different layers}. These recollisions come from the condition (\ref{it.set_F_l_4_intro}) in Section \ref{sec.statement_formula_idea}: this condition implies that the particle lines in $r(\Mb_{\ell'})$ are connected both \emph{within layer $\Mb_{\ell'+1}$} and \emph{within layers $\Mb_{\leq\ell'}:=\cup_{\ell''\leq\ell'}\Mb_{\ell''}$}, resulting in cycles.

The idea of layer selection is as follows. First choose a layer $\ell_1$ such that $s_{\ell_1-1}=|r(\Mb_{\ell_1-1})|\gtrsim\rho$, and that $s_{\ell_1-1}\gg s_{\ell_1}$; this means that the particle lines in $r(\Mb_{\ell_1-1})$ are connected to each other (with negligible exceptions) within layer $\Mb_{\ell_1}$. Now, choose the \emph{highest} layer $\ell_2<\ell_1$ such that sufficiently many particle lines in $r(\Mb_{\ell_1-1})$ are connected to each other \emph{within $\Mb_{[\ell_2:\ell_1-1]}:=\cup_{\ell_2\leq\ell'\leq\ell_1-1}\Mb_{\ell'}$}. Then, by putting together layers $\Mb_{\ell'}\,(\ell_2<\ell'\leq\ell_1)$ into $\Mb_U$ (and removing a negligible number of exceptions), and defining $\Mb_{\ell_2}:=\Mb_D$, we obtain a two layer molecule where each of $\Mb_U$ and $\Mb_D$ is a tree (or forest), and there exist $\gtrsim\rho$ bonds between them. For details see Definition \ref{def.layer_select_toy} and Proposition \ref{prop.layer_select_toy}, and {\color{blue}Figure \ref{fig.layerselect}} (Section \ref{sec.toy_multi}).

Once we are in the two layer setting, we then consider the \emph{dichotomy} determined by the quantities $\#_{\mathrm{2conn}}$ (the number of 2-connections as defined above) and $\#_{\mathrm{comp}(X)}$ (the number of connected components of $X$), where $X$ is the set of atoms in $\Mb_U$ that have exactly one bone with $\Mb_D$. By definition, if $\#_{\mathrm{2conn}}\gtrsim\rho$ then $\Mb$ is toy model III. If $\#_{\mathrm{2conn}}\ll\rho$, then the 2-connections will be negligible and we may assume$\#_{\mathrm{2conn}}=0$. Then, $\Mb$ will be toy model I+ or II, depending on whether $\#_{\mathrm{comp}(X)}\ll\rho$ or $\#_{\mathrm{comp}(X)}\gtrsim\rho$. This completes the proof in Section \ref{sec.toy}.

\textbf{The full algorithm.} Finally we discuss the extra ingredients in the full algorithm in Sections \ref{sec.layer}--\ref{sec.maincr}, in addition to those in Section \ref{sec.toy}. These are summarized in Section \ref{sec.toy_reduce} and include the followings:

\begin{itemize}
\item The O-atoms: in the general case we need to replace the relevant notions (adjacency, connectedness etc.) by the ones associated with ov-segments, and modify the cutting algorithms accordingly. For example, the algorithms \textbf{UP} (Definition \ref{def.up_toy}), \textbf{2CONNUP} (Definition \ref{def.toy3_alg}), \textbf{3COMPUP} (Definition \ref{def.toy2_alg}) and \textbf{MAINUD} (Definition \ref{def.toy1+_alg}), have their full versions defined in Definitions \ref{def.alg_up}, \ref{def.alg_2connup}, \ref{def.3comp_alg} and \ref{def.alg_maincr} respectively.
\item Layer refinement: this is discussed in Definition \ref{def.layer_refine} (Section \ref{sec.layer_refine}), and is needed to treat the possible recollisions \emph{within each layer} $\Mb_{\ell'}$. We further sub-divide each $\Mb_{\ell'}$ into finitely many subsets, such that each of them is precisely a forest, which recovers assumption (ii). This also leads to an extra step in the layer selection process, see Definition \ref{def.layer_select}.
\item The strong an weak resonances: these are exceptional cases where \{33\} molecules do not gain powers, see Definitions \ref{def.strdeg} and \ref{def.weadeg}. In reality, if this effect becomes significant, we shall perform alternative cutting algorithms that exploit such degeneracies to gain powers of $\varepsilon$ in a different way, see Propositions \ref{prop.case2} and \ref{prop.comb_est_case4}.
\item The pre-processing step: this is discussed in Definition \ref{def.func_select} (Section \ref{sec.preprocess}). This additional step is needed to deal with possible fixed ends in $\Mb_D$ for the two-layer molecule $\Mb$, which may cause trouble to the greedy algorithm in the toy model I+.
\end{itemize}

For a summary of the ingredients of the proofs in Sections \ref{sec.layer}--\ref{sec.maincr} and their logical relations, see {\color{blue} Figure \ref{fig.flowchart}} (Section \ref{sec.toy_reduce}).
\subsection{The Boltzmann approximation and truncation error}\label{subsec.idea.fA} In Sections \ref{subsec.idea.cumulants}--\ref{subsec.idea.alg} above, we have focused on the proof of the $E_H$ estimates in (\ref{eq.cumuest_intro}) (see Proposition \ref{prop.cumulant_est}). Now it remains to prove the approximation of $f^\Ac$ by the Boltzmann solution, see Proposition \ref{prop.est_fa}.

By discussions in Section \ref{subsec.idea.cumulants}, we know that $f^\Ac(\ell\tau)$ is given by an explicit series expansion of iterated time integrals that \emph{involve only $f^\Ac((\ell-1)\tau)$}, therefore this part follows from essentially the same arguments in the short-time theory. We identify those leading terms that do not contain recollision, and those sub-leading terms that contain recollision. We then match each leading term with the corresponding term that occurs in the Duhamel expansion of $f(\ell\tau)$ in terms of $f((\ell-1)\tau)$; note that the cluster expansion (see (\ref{eq.fAterm})) actually leads to extra terms in addition to those occurring in Boltzmann Duhamel expansions, but these extra terms turn out to cancel each other, as shown in Proposition \ref{prop.cancel}. As for the sub-leading terms, they can be proved to be vanishing by applying a very simple special case of the cutting algorithm in Sections \ref{sec.toy}--\ref{sec.maincr} (see Section \ref{sec.fa_subleading}, \textbf{Point 5}), since we only have \emph{one single layer} and the number of recollisions in each cluster is bounded by $\Gamma$. All these are contained in Section \ref{sec.fa}.

Finally we discuss the estimate for the error term $f_s^{\mathrm{err}}$ caused by the truncated dynamics. This estimate is stated in Proposition \ref{prop.cumulant_error} and is proved in Section \ref{sec.error}. Here again we reduce the estimate to the construction of a cutting algorithm. However, as we now have many recollisions, one new difficulty arises due to the possibility that a small group of particles within distance $O(\varepsilon)$ of each other producing many recollisions within time $O(\varepsilon)$. To treat this, we rely on \emph{a new ingredient}, namely the absolute upper bound on the number of collisions given the number of particles, which is proved in Burago-Ferleger-Kononenko \cite{BFK98}. Basically, if $\Gamma$ is sufficiently large, then these $\Gamma$ recollisions must involve sufficiently many particles, with each new particle providing a \{33\} molecule in a suitable algorithm.

The algorithm in Section \ref{sec.error} is similar to the one in Sections \ref{sec.toy}--\ref{sec.maincr}, and is based on a two-layer molecule formed by $A^+$ and $A$ where $A$ is (almost) a tree and $A^+$ is arbitrary; see Proposition \ref{prop.trans} (which actually involves three layers $A^+$, $A$ and $A^-$, but in reality we only exploit two of them). The main algorithm, called  \textbf{MAINTRUP}, is similar to the \textbf{MAINUD} algorithm in Section \ref{sec.mainud}, where each time we cut a lowest atom in $A^+$, and keep cutting any \{33\} or \{343\} structure in $A$ that occur. This will then guarantee a lower bound for the number of \{33\} molecules, provided we have a lower bound for the number $\#_{\mathrm{conn}}^+$ of bonds connecting an atom in $A$ to an atom in a specific subset $X^+(A)\subseteq A^+$. This last lower bound is then provided by (the inverse function of) the Burago-Ferleger-Kononenko upper bound in \cite{BFK98}, see Proposition \ref{prop.upper_bound_col}. For more details of the algorithm, see Definition \ref{def.alg_maintrup} and Proposition \ref{prop.alg_maintrup}.
\subsection{Comparison between wave and particle systems}\label{subsec.idea.compare} As mentioned in Section \ref{sec.past}, the proof in this paper is deeply inspired by that in the wave setting by the first two authors \cite{DH23_2}. Here we make a brief comparison between the two proofs, in terms of their similarities and differences. In the case of wave kinetic theory, the microscopic system in \cite{DH23_2} is given by the nonlinear Schr\"odinger (NLS) equation on a torus $\Tb^d_L$ of size $L$, with a nonlinearity of strength $\alpha$. The lowest-order nontrivial statistical quantity is the second moment of the Fourier modes $\mathbb E|\widehat u(t, k)|^2$, where $\widehat u(t,k)$ is the Fourier transform of the solution $u(t,x)$ to the NLS equation, which is to be approximated by the solution of the (homogeneous) wave kinetic equation \eqref{eq.wke}, in the limit $L\to \infty$ and $\alpha\to 0$. This kinetic limit is taken under an appropriate choice of scaling law between $L$ and $\alpha$, namely $\alpha=L^{-\gamma}$, where $\gamma\in(0,1)$ in the homogeneous setting, see \cite{DH23_2}.

At a conceptual level, the two proofs are very similar. The fundamental difficulty represented by the divergence of the series expansion on the long time interval $[0, t_{\mathrm{fin}}]$ is the same, and the same goes for the need to create an arrow of time where, in both cases, the microscopic system is time-reversible but the limiting kinetic equation is irreversible. The idea of performing a partial expansion, where one refrains from expanding the lowest order statistic (the $f^\Ac$ here, and the second moments in the wave case) but expands the cumulants all the way to time $t=0$, is also the same. This keeps track of all the relevant collisions or interactions on the long time interval $[0, t_{\mathrm{fin}}]$ that affect the cumulants at time $t_{\mathrm{fin}}$.

Even at a more technical level, there are uncanny similarities between the two proofs. First, both proofs start with a truncation of the dynamics in order to control the size of the diagrams appearing in the cumulant expansions. Also, in both proofs, the needed cumulant estimate is reduced to a combinatorial problem that is resolved through an elaborate algorithm. This combinatorial problem, in both cases, involves the combinatorial structure of \emph{molecules}, which are formed by \emph{atoms} that represent collisions or interactions. While the algorithms are vastly different in the wave case and particle case, there is an important similarity in the gain mechanism. As shown in Section \ref{subsec.idea.alg}, the basic unit of gain here is the special structure of \{33\} molecules; similarly, in the wave case, the basic unit of gain also comes from the \{33\} molecules and the associated five-vector counting problem.

All this being said, the analogy between the particle and wave kinetic theories has its limitations, and despite the above similarities, each proof carries its own distinct set of challenges. In the wave case, a large box limit needs to be performed to go from the discrete Fourier space of $\Tb^d_L$ to the continuum setting of the wave kinetic equation \eqref{eq.wke}. This links the wave kinetic theory to some deep number theory problems, which seems to be absent in the particle setting. Also, the Feynman diagram expansion for the NLS equation turns out to converge \emph{only conditionally and not absolutely}, and requires highly nontrivial cancellations between arbitrarily large diagrams, which is another distinctive feature and novelty of  \cite{DH23_2}. On the other hand, the cumulant expansion in the particle setting is technically more involved, due to the application of inclusion-exclusion principle and the resulting notion of overlaps, which does not seem to be have an analog in the wave case. Also the main algorithm in the particle case seems to be conceptually more complicated than the corresponding algorithm in the wave case, which may be due to the singular nature of the $\dirac$ integrals in $\Ic_\Mb(Q_\Mb)$ compared to the counting problems in \cite{DH23_2}.
\subsection{The rest of this paper} The rest of this paper is organized as follows. Note that, although the notion of cutting and cutting algorithm is the main novelty and central to the proof, they have to wait until Section \ref{sec.cutting} (cutting) and Section \ref{sec.toy} (algorithm) to appear, due to the necessary preparations in Sections \ref{sec.formula_cumulant}--\ref{sec.local_associated}.
\begin{itemize}
\item Part I: Sections \ref{sec.formula_cumulant}--\ref{sec.proof_cumulant_formula} are devoted the statement and proof of Proposition \ref{prop.cumulant_formula}. 
\begin{itemize}
\item[\labelitemi] In Section \ref{sec.formula_cumulant} we define all the relevant notions (molecules, topological reduction etc.) in Sections \ref{sec.molecule}--\ref{sec.molecule_prescribed}, and then state Proposition \ref{prop.cumulant_formula} in Section \ref{sec.formula_cumulant_statement}. 
\item[\labelitemi] In Section \ref{sec.truncation_large_molecule} we define the truncated dynamics (Definition \ref{def.T_dynamics}) and its properties, which are needed in the proof of Proposition \ref{prop.cumulant_formula}.
\item[\labelitemi] In Section \ref{sec.proof_cumulant_formula} we prove Proposition \ref{prop.cumulant_formula}.
\end{itemize}
\item Interlude I: Section \ref{sec.reduction} contains the statement of the main propositions, Propositions \ref{prop.est_fa}--\ref{prop.cumulant_error}, of this paper. We also prove that these propositions imply Theorem \ref{th.main}.
\item Part II: Sections \ref{sec.local_associated}--\ref{sec.treat_integral} are devoted to the first half of the proof of Proposition \ref{prop.cumulant_est}, where we reduce it to Proposition \ref{prop.comb_est}.
\begin{itemize}
\item[\labelitemi] In Section \ref{sec.local_associated} we make the necessary preparations, and reduce the the expression $|\Ic\Nc_\Mb|$ occurring in the upper bound of $E_H$ to the local integration $\Ic_\Mb(Q_\Mb)$.
\item[\labelitemi] In Section \ref{sec.cutting} we define the notion of cutting (Definition \ref{def.cutting}), study its properties, and how the integral $\Ic_\Mb(Q_\Mb)$ behaves under cutting (Proposition \ref{prop.cutting}).
\item[\labelitemi] In Section \ref{sec.treat_integral} we state the main result concerning cutting algorithms, i.e. Proposition \ref{prop.comb_est}, and prove that it implies Proposition \ref{prop.cumulant_est}.
\end{itemize}
\item Interlude II: Section \ref{sec.prepare} contains a few simple reductions and reiterations of properties of the molecule $\Mb$, in preparation of the proof of Proposition \ref{prop.comb_est}.
\item Part III: Sections \ref{sec.toy}--\ref{sec.maincr} are devoted to the second half of the proof of Proposition \ref{prop.cumulant_est}, i.e. the proof of Proposition \ref{prop.comb_est}.
\begin{itemize}
\item[\labelitemi] In Section \ref{sec.toy} we study the two-layer toy models, namely toy models I, I+, II and III (in Sections \ref{sec.toy1}--\ref{sec.toy3} respectively). In Section \ref{sec.toy_multi} we define layer selection in the toy model case, and in Section \ref{sec.toy_reduce} we sketch the additional ingredients needed in the full algorithm.
 \item[\labelitemi] In Section \ref{sec.layer} we present the first half of the proof of Proposition \ref{prop.comb_est}, where we reduce it to Proposition \ref{prop.case5}. This section contains the \textbf{UP} algorithm (Definition \ref{def.alg_up}), strong and weak degeneracies (Definitions \ref{def.strdeg}, \ref{def.weadeg}), layer refinement (Section \ref{sec.layer_refine}) and layer selection (Section \ref{sec.layer_select}).
  \item[\labelitemi] In Section \ref{sec.maincr} we finish the proof of Proposition \ref{prop.case5}. In Sections \ref{sec.2connup}--\ref{sec.mainud} we introduce the algorithms \textbf{2CONNUP}, \textbf{3COMPUP}, pre-processing and \textbf{MAINUD} respectively. Finally, in Section \ref{sec.finish} we finish the proof of Proposition \ref{prop.case5} by invoking the dichotomy and reducing to one of the toy model cases in Section \ref{sec.toy}. 
\end{itemize}
\item The Finale: Sections \ref{sec.fa}--\ref{sec.error} are devoted to the proof of Propositions \ref{prop.est_fa} and \ref{prop.cumulant_error}.
\item Appendices: Appendix \ref{app.aux}--\ref{app.glossary} contain auxiliary results and tables summarizing important notions.
\end{itemize}

\section{The cumulant expansion formula}\label{sec.formula_cumulant} In this section, we will state the main cumulant formula in Proposition \ref{prop.cumulant_formula}, which has form
\begin{equation}\label{eq.cumulant_expansion_sec_intro}
    f_s(\ell\tau,\vz_s) = \sum_{H\subseteq [s]} (f^\Ac(\ell\tau))^{\otimes([s]\backslash H)} \cdot E_H(\ell\tau, \vz_H) + (\mathrm{error\ terms}),
\end{equation} 
\begin{equation}\label{eq.formula_cumulant_sec_intro} 
    |E_H(\ell\tau, \vz_H)|\leq\sum_{[\Mb]\in \Fc_{\vLambda_{\ell}}}|\Ic\Nc_\Mb|(\vz_H).
\end{equation} 
Here $f^\Ac$ is an approximate solution to the Boltzmann equation, and $\Mb$ is a molecule (Definition \ref{def.molecule}). Also $\Fc_{\vLambda_{\ell}}$ is a collection of molecules (Definition \ref{def.set_T_F}) and $|\Ic\Nc_\Mb|$ is an integral expression associated with the molecule $\Mb$ (Definition \ref{def.associated_int}), which represents the densities given by trajectories that topologically reduce to $\Mb$ (Definition \ref{def.top_reduction}).

The statement of Proposition \ref{prop.cumulant_formula} involves many concepts that need to be defined. We define these notions in Sections \ref{sec.molecule}--\ref{sec.molecule_prescribed}, which are related to molecules $\Mb$ and topological reductions. Then, in Section \ref{sec.formula_cumulant_statement} we define the set $\mathcal{F}_{\boldsymbol{\Lambda}_\ell}$ and the associated integral $|\Ic\Nc_\Mb|(\vz_H)$ in (\ref{eq.formula_cumulant_sec_intro}), and then state Proposition \ref{prop.cumulant_formula}. This will be the main goal in the first part of the paper, whose proof will occupy Sections \ref{sec.truncation_large_molecule}--\ref{sec.proof_cumulant_formula}.

\subsection{Molecules}\label{sec.molecule} We start with the definition of molecules. As seen in {\color{blue}Figure \ref{fig.topological_reduction}} (Section \ref{sec.topological_reduction}), molecules are topological reductions of actual particle trajectories (see Definition \ref{def.top_reduction} for precise formulation).

\begin{definition}[Molecules]\label{def.molecule} We define a (layered) molecule to be a diagram demonstrated by {\color{blue}Figure \ref{fig.molecule}}.
More rigorously, we define a \newterm{(layered) molecule} $\Mb=(\Mc,\Ec,\Pc,\Lc)$ to be a diagram as follows. For the correspondence between molecule pictures and physical particle pictures, see {\color{blue}Table \ref{tab.trans}} in Appendix \ref{app.trans}.
\begin{figure}[h!]
    \includegraphics[width=0.5\linewidth]{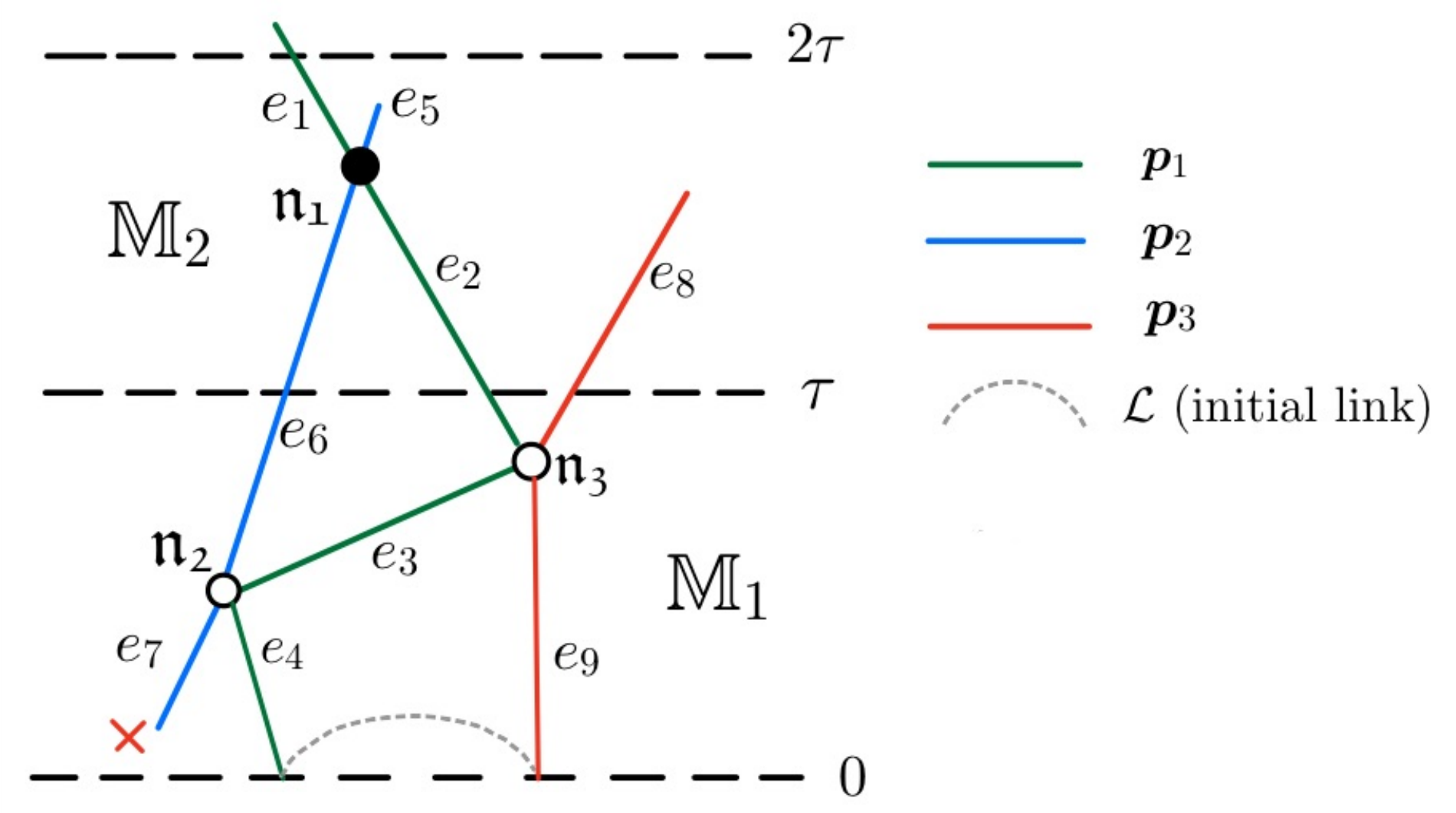}
    \caption{An example of molecules.}
    \label{fig.molecule}
\end{figure}
\begin{enumerate}
    \item \textit{Atoms, edges, bonds and ends.} $\Mc$ is a set of nodes referred to as \newterm{atoms}, and $\Ec$ is a set of \newterm{edges} consisting of \newterm{bonds} and \newterm{ends} to be specified below. For convenience, below we will abuse notation and use $\Mb$ instead of $\Mc$ to denote the set of atoms in $\Mb$.
    \begin{enumerate}
        \item\label{it.molecule_atom} An atom $\nf\in \Mb$ is either a \newterm{C-atom} (collision) or an \newterm{O-atom} (overlap), denoted respectively by $\circ$ and $\bullet$ as in {\color{blue}Figure \ref{fig.molecule}}. Each atom is assigned with a unique \newterm{layer} $\ell[\nf]$ which is an integer in $[\underline{\ell}:\ell]$, where $\ell$ is the highest layer and $\underline{\ell}$ is the lowest layer of $\Mb$. Define $\Mb_{\ell'}$ to be the set of atoms $\nf\in\Mb$ of layer $\ell[\nf]=\ell'$. The boundaries between layers are indicated by horizontal lines, as in {\color{blue}Figure \ref{fig.molecule}}.

        \item\label{it.molecule_edge} Each edge in $\Ec$ is either a \newterm{bond} between two atoms $\nf$ and $\nf'$, or an \textbf{end} emanating from one atom $\nf$, or an \textbf{empty end} that does not involve any atom. Each edge is either \textbf{top} or \textbf{bottom} when viewed at each involved atom $\nf$; a top (resp. bottom) edge is drawn above (resp. below) the atom, as in {\color{blue}Figure \ref{fig.molecule}}. Moreover, each end is also either \newterm{free} or \newterm{fixed}, with fixed ends indicated by a $\times$, as shown in {\color{blue}Figure \ref{fig.molecule}}. We understand by convention that empty ends are always free, and count as both top and bottom.
        \item\label{it.molecule_parent} For two atoms $\nf$ and $\nf'$, if a bond exists between $\nf$ and $\nf'$ that is bottom at $\nf$ and top at $\nf'$, we say $\nf$ is \textbf{parent} of $\nf'$ and $\nf'$ is \textbf{child} of $\nf$.
    \end{enumerate}
    \item\label{it.particle_line_same} \textit{Particle lines.} $\Pc$ is a collection of paths: each path starts from a bottom end at an atom, then keeps connecting to a parent by a bond, until reaching a top end at another atom. Each chain $\pb\in\Pc$ is called a \textbf{particle line}. We require that each edge belongs to a unique particle line; in particular, each empty end is also viewed as a particle line by convention. For a top edge and a bottom edge at an atom $\nf$, we say they are \textbf{serial} if they belong to the same particle line. Note that each particle line $\pb$ contains a unique bottom end and a unique top end.  We will identify particle lines with physical particles, see {\color{blue}Table \ref{tab.trans}}.
    
    Recall $\ell$ defined in (\ref{it.molecule_atom}). We may define a \textbf{start layer} $\ell_1[\pb]\in[\underline{\ell}:\ell]$ and \textbf{finish layer} $\ell_2[\pb]\in[\underline{\ell}:\ell+1]$ for each particle line $\pb$. As $\pb$ contains a unique bottom (and unique top) end, we may also equivalently write $\ell_1[e]$ (resp. $\ell_2[e]$) for each bottom (resp. top) end $e$ instead of $\ell_1[\pb]$ and $\ell_2[\pb]$. We refer to the interval $[\ell_1[\pb]:\ell_2[\pb]]$ as the \textbf{lifespan} of the particle line $\pb$. In {\color{blue}Figure \ref{fig.molecule}}, the lifespan of $\pb_1$ and $\pb_2$ are $[1: 3]$ and $[1: 2]$, respectively. For $\pb_1$, its finish layer is $\ell_2[\pb_1] = \ell + 1 = 3$, which is above the highest layer $\ell$. To indicate this, we draw it penetrating the top horizontal line.
    \item\label{it.initial_link} \textit{Initial links.} Finally, $\Lc$ is a collection of pairs of bottom ends $(e,e')$ with $\ell_1[e]=\ell_1[e']=1$ (provided $\underline{\ell}=1$). Since each particle line has a unique bottom end by (\ref{it.particle_line_same}), we can equivalently view $\Lc$ as a collection of pairs of particle lines with start layer $1$. Each pair in $\Lc$ is called an \textbf{initial link}. In {\color{blue}Figure \ref{fig.molecule}}, $\Lc$ contains a single pair $(e_4, e_9)$, which is indicated by a dashed line connecting them.
    \item\label{it.ov_segment} \textit{Ov-segments.} For later use, we also define an \textbf{ov-segment} to be a path between two atoms (or ends) such that all intermediate atoms in this path are O-atoms (except the possible endpoints), and this path is contained in a single particle line. We say an ov-segment is \textbf{maximal} if it is not contained in other ov-segments. In {\color{blue}Figure \ref{fig.molecule}}, the path consisting of edges $e_1$ and $e_2$, and atoms $\nf_1$ and $\nf_3$, is an example of a maximal ov-segment. Note that each edge in $\Mb$ belongs to a unique maximal ov-segment. Also each bond in $\Mb$ is a special case of ov-segment (not necessarily maximal).
    \item The objects $(\Mc, \Ec, \Pc, \Lc)$ satisfies the following requirements. 
    \begin{enumerate}
        \item Each atom has exactly 2 top edges and exactly 2 bottom edges.
        \item Starting from any atom $\nf$, it is impossible to get back to $\nf$ by iteratively taking parent or by iteratively taking child.
        \item If $\nf$ is a parent of $\nf'$ then we must have $\ell[\nf]\geq \ell[\nf']$. If $\nf$ belongs to a particle line $\pb$ then we must have $\ell_1[\pb]\leq\ell[\nf]\leq\ell_2[\pb]$.
    \end{enumerate}
\end{enumerate}
\end{definition}
\begin{definition}[Full molecules]\label{def.full_molecule}
    If a molecule $\Mb$ has no fixed end, we say $\Mb$ is a \newterm{full molecule}. 
\end{definition}

In the statement of the cumulant formula, we will need to label different particles to distinguish them. Also the particles corresponding to $\vz_s$ or $\vz_H$ in $f_s(\ell\tau, \vz_s)$ or cumulant $E_H(\ell\tau, \vz_H)$, i.e. the root particles, will play a special role.

\begin{definition}[Label and root particle lines]\label{def.root_particle_line} Given a molecule $\Mb$, we may assign \textbf{labels} (which takes values in $\Nb$) to each particle line $\pb$ of $\Mb$. This makes $\Mb$ a \textbf{labeled molecule}. Note that two labeled molecules with the same molecule structure (Definition \ref{def.molecule}) but different labelings are regarded as different.

For any molecule $\Mb$, define the \textbf{set of all particle lines} in $\Mb$ to be $p(\Mb)$. Define a particle line $\pb$ to be a \textbf{root particle line}, if its finish layer $\ell_2[\pb]=\ell+1$, where $\ell$ is the highest layer in $\Mb$. Define the set of root particle lines to be $r(\Mb)$, and the set of non-root particles to be $(p\backslash r)(\Mb):=p(\Mb)\backslash r(\Mb)$.

For two labeled molecules $\Mb$ and $\Mb'$ with $\ell$ layers, we define them to be \textbf{equivalent}, denoted $\Mb\sim\Mb'$, if they have the same molecule structure (Definition \ref{def.molecule}) and the same labels for root particle lines. Define $[\Mb]$ to be the equivalence class of $\Mb$, which is obtained by relabeling non-root particle lines in $\Mb$. 

In {\color{blue}Figure \ref{fig.molecule}}, we can label $\pb_1$, $\pb_2$ and $\pb_3$ by $1$, $2$ and $3$ respectively, where $\pb_1$ is the root particle line.
\end{definition}
\begin{remark}\label{rem.label} Clearly an equivalence class $[\Mb]$ is uniquely determined by an unlabeled molecule $\Mb$, and a fixed labeling of all root particle lines. In the proof of the cumulant expansion formula (Proposition \ref{prop.cumulant_formula}) in Sections \ref{sec.truncation_large_molecule}--\ref{sec.proof_cumulant_formula}, it is important that we use $\Mb$ to represent labeled molecules and sum over equivalence classes $[\Mb]$ (see relevant notations in Definition \ref{def.notation} (\ref{it.sum_over_equiv})). Once this proof is finished, in subsequent analysis of the $|\Ic\Nc_\Mb|$ terms, the labeling will not play any role, so we will ignore it and use $\Mb$ to denote unlabeled molecules in Sections \ref{sec.reduction}--\ref{sec.error} (we still need to label the root particles, which will not play a major role, see Remark \ref{rem.label_root}).
\end{remark}
We introduce some additional notions associated with molecules, which will be used in the proof below.

\begin{definition}[Ordering of atoms]\label{def.molecule_order} Let $\Mb$ be a molecule and $\nf,\nf'\in\Mb$ be atoms.
\begin{enumerate}
    \item \emph{Descendants and ancestors.} If $\nf$ is obtained from $\nf'$ by iteratively taking parents, then we say $\nf$ is an \newterm{ancestor} of $\nf'$, and $\nf'$ is a \newterm{descendant} of $\nf$. We use $S_\nf$ (resp. $Z_\nf$) to denote the set of descendants (resp. ancestors) of $\nf$. By convention we understand that $\nf\in S_\nf$ (and $\nf\in Z_\nf$).
    
    \item\label{it.partial_order} \emph{The partial ordering of atoms and edges.} We also define a \newterm{partial ordering} between all atoms, by defining $\nf'\prec\nf$ if and only if $\nf'\neq\nf$ and is a descendant of $\nf$ (equivalently $\nf$ is ancestor of $\nf'$). In particular $\nf'\prec\nf$ if $\nf'$ is a child of $\nf$ (equivalently $\nf$ is parent of $\nf'$). 
    \item\label{it.minimax} Given any subset $A\subseteq\Mb$, we can define a \newterm{lowest} atom in $A$ to be any atom $\nf\in A$ that satisfies $\nf'\not\in A$ for each $\nf'\prec\nf$ (i.e. $\nf'\neq\nf$ is descendant of $\nf$). Define \newterm{highest} atoms in the same way. These are just minimal and maximal elements of $A$ under the partial ordering $\prec$, and note that they need not be unique for general $A$.
\end{enumerate}
\end{definition}

\begin{definition}[Various sets associated with molecules]
\label{def.sets_molecule} Let $\Mb=(\Mc,\Ec,\Pc,\Lc)$ be a molecule.
\begin{enumerate}
    \item\label{it.subsets_Ec} \emph{Sets $\Ec_*$ and $\Ec_{\mathrm{end}}$.} Define $\Ec_*$ to be the set of all edges that are either bonds or free ends (i.e. not counting fixed ends). Define also the sets $\Ec_{\mathrm{end}}$ (resp. $\Ec_{\mathrm{end}}^-$, $\Ec_{\mathrm{end}}^+$) to be the set of all ends (resp. bottom ends, top ends). Note that empty ends belong to $\Ec_{\mathrm{end}}^\pm$ by convention.
    
    Note that the sets $\Ec_{\mathrm{end}}^-$ and $\Ec_{\mathrm{end}}^+$ are in bijection with the set $p(\Mb)$ of all particle lines, see Definition \ref{def.molecule} (\ref{it.particle_line_same}).
    \item \emph{Number of atoms and particle lines.} We will use $|\Mb|$, $|\Mb|_{\mathrm{O}}$, $|\Mb|_p$, $|\Mb|_r$ and $|\Mb|_{p\backslash r}$ to denote the number of atoms, O-atoms, particle lines (i.e. $|p(\Mb)|$), root particle lines (i.e. $|r(\Mb)|$), and non-root particle lines (i.e. $|(p\backslash r)(\Mb)|$) of $\Mb$ respectively.
\end{enumerate}
\end{definition}

\subsection{Topological reduction}\label{sec.topological_reduction} In this section, we rigorously formulate the notion of topological reduction. We start with the notion of admissible trajectories, which generally occur from modifications of hard sphere dynamics where certain collisions are neglected according to certain rules (see Definition \ref{def.modified_dynamics}).
\begin{definition}[Admissible trajectories]\label{def.admissible}
Let $\vz(t)=\vz_N(t)=(z_j(t))_j$ be any particle trajectory (a priori with no equation or restriction). We say it is \textbf{admissible} (see {\color{blue}Figure \ref{fig.3}}), if
\begin{figure}[h!]
\includegraphics[scale=0.25]{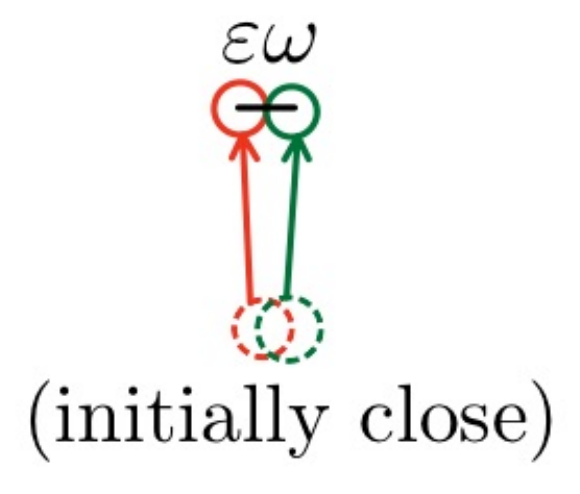}
\caption{Admissible trajectories (Definition \ref{def.admissible}), where it is allowed to have $|x_i-x_j|<\varepsilon$, but two particles reaching distance $\varepsilon$ from ``inside" will not collide and not change velocity, as they do not satisfy the pre-collisional condition.}
\label{fig.3}
\end{figure}
\begin{enumerate}
\item Each $z_j(t)=(x_j(t),v_j(t))$ is piecewise linear, $x_j(t)$ is continuous and $v_j(t)$ is left continuous, with a discrete set of times of velocity discontinuity referred to as \textbf{collisions};
\item\label{it.precol} At each collision time $t$, there exist two and only two particles whose velocities are discontinuous, their positions $(x_i,x_j)$ and velocities $(v_i,v_j)$ satisfy (\ref{eq.hardsphere}) with $x_i(t)-x_j(t)=\varepsilon\omega\,(\omega\in\Sb^{d-1})$, and $(v_i(t)-v_j(t))\cdot\omega\leq 0$. We refer to this set of conditions as the \textbf{pre-collisional condition}.
\item \label{it.unique}Additionally, we define an \textbf{overlap} to be a time when $x_i(t)-x_j(t)=\varepsilon\omega\,(\omega\in\Sb^{d-1})$ and $(v_i(t)-v_j(t))\cdot\omega\leq 0$, but the velocities $v_i$ and $v_j$ are continuous at time $t$.
\end{enumerate}
\end{definition}
\begin{definition}[Topological reduction]\label{def.top_reduction} Let $\vz(t)$ be an admissible trajectory for $t\in[0,\ell\tau]$, and $\Mb$ be a molecule with $\ell$ layers (Definition \ref{def.molecule}). We say $\vz(t)$ \textbf{topologically reduces} to $\Mb$, if there is a bijection between maximal ov-segments of $\Mb$ (Definition \ref{def.molecule} (\ref{it.ov_segment})) and maximal linear pieces of the trajectory $\vz(t)$, see {\color{blue}Figure \ref{fig.topological_reduction}}, which satisfies the following requirements:
\begin{enumerate}
    \item If two maximal ov-segments intersect at a C-atom, then the corresponding linear pieces of $\vz(t)$ have a collision. If they intersect at an O-atom, then the corresponding linear pieces have an overlap.
    \item The induced correspondence between atoms and collisions/overlaps preserves order (i.e. if $\nf$ is parent of $\nf'$ then the corresponding collision/overlap happens later in time), and \emph{is bijective on collisions} (it is injective on overlaps by time ordering, but might not be surjective).
    \item If two maximal ov-segments belong to the same particle line in $\Mb$, then the corresponding linear pieces belong to the trajectory of the same particle.
    \item\label{it.reduct_layer} The collision/overlap corresponding to an atom in layer $\ell$ must occur in the time interval $[(\ell-1)\tau, \ell\tau]$.
\end{enumerate}

\noindent For a summary of correspondence of concepts in molecules and physical trajectories, see {\color{blue}Table \ref{tab.trans}}. In {\color{blue}Figure \ref{fig.topological_reduction}}, the overlap in $[\tau, 2\tau]$ of the red and green particles corresponds to the O-atom on layer $2$, and the collision in $[0, \tau]$ of the red and green particles corresponds to the C-atom on layer $1$.
\end{definition}

\begin{figure}[h!]
    \centering
    \includegraphics[width=0.7\linewidth]{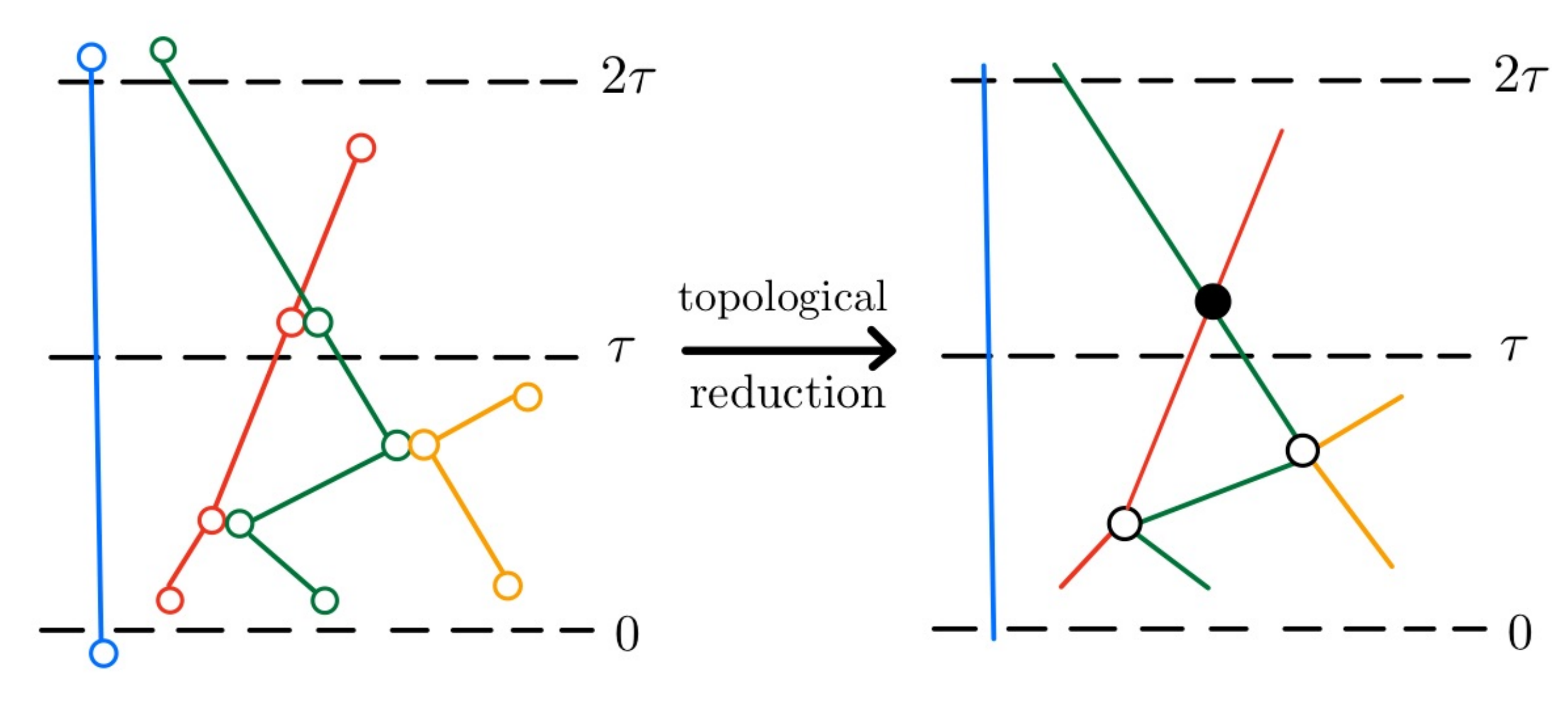}
    \caption{An illustration of topological reduction.}
\label{fig.topological_reduction}
\end{figure}
\begin{remark}\label{rem.top_reduction_uniqueness} 
The topological reduction $\Mb$ of a given trajectory $\vz(t)$ may not be unique, since we do not require a bijective correspondence between O-atoms and the actual set of overlaps. The O-atoms may represent only an arbitrary subset of all overlaps.  However there is always a \emph{unique} topological reduction that is C-molecule (i.e. contains only C-atoms). 
\end{remark}

\subsection{Clusters}\label{sec.cluster} As seen in Definition \ref{def.top_reduction}, O-atoms correspond to overlaps that do not affect the dynamics. Ignoring all such overlaps leads to the definition of clusters.

\begin{definition}[Create and delete an O-atom]\label{def.create_delete} Given a molecule $\Mb$, we define the following operations, see {\color{blue} Figure \ref{fig.deletion}}. The particle lines and ov-segments will not be broken after these operations.
\begin{figure}[h!]
    \centering
    \includegraphics[width=0.4\linewidth]{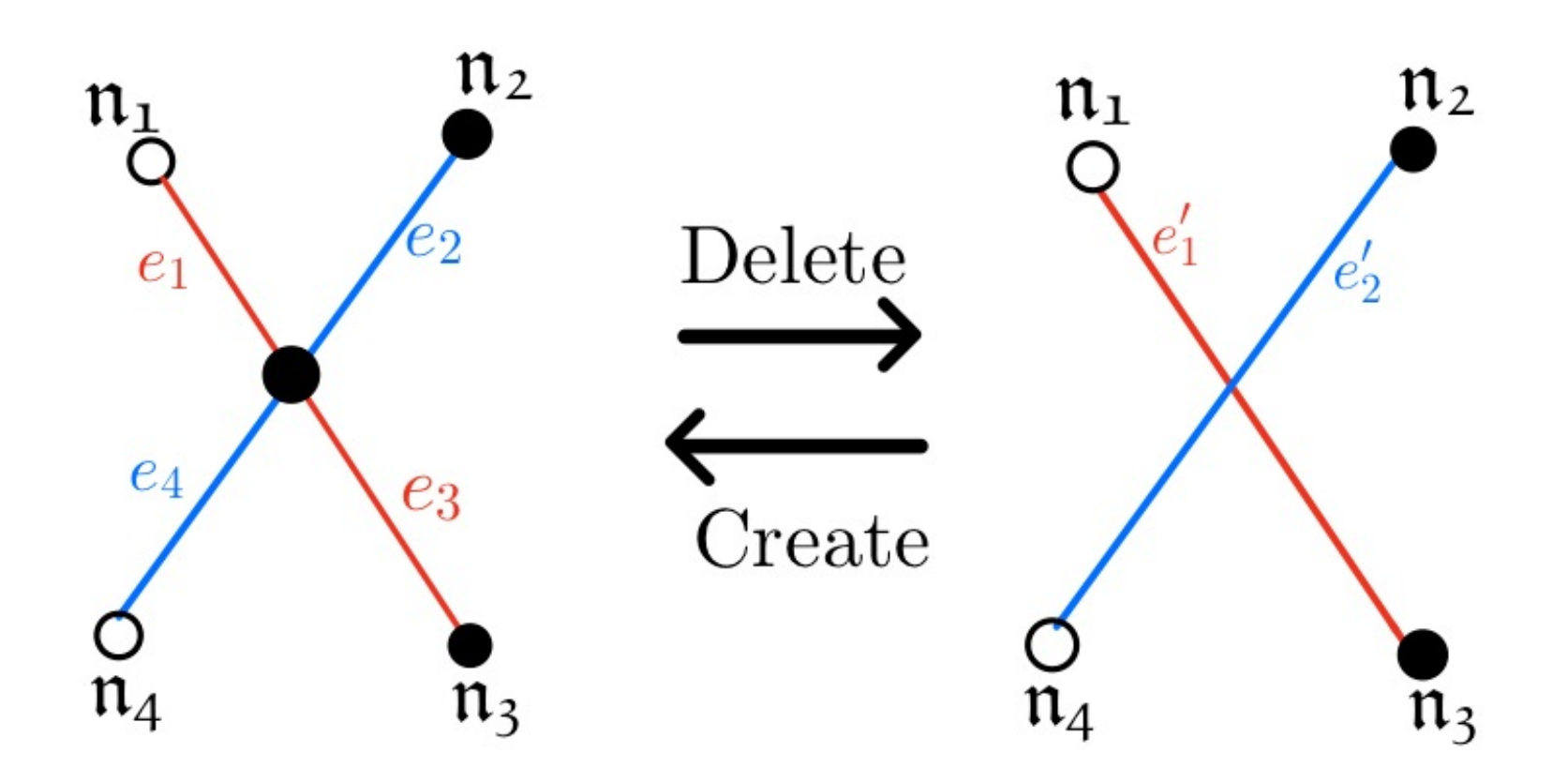}
    \caption{Creating and deleting O-atoms as in Definition \ref{def.create_delete}. Note that each atom $\nf_j$ can be either C- or O-atoms, or absent.}
    \label{fig.deletion}
\end{figure}
\begin{enumerate}
    \item \emph{Delete an O-atom.} Let $\of$ be an O-atom with 4 edges $e_i\,(1\leq i\leq 4)$ connected to other atoms $\nf_i$, where $e_1$ and $e_3$ are serial and same for $e_2$ and $e_4$. By \textbf{deleting} the O-atom $\of$, we remove $\of$ and the 4 edges $e_j$ from $\Mb$, and add two bonds $e_1'$ between $\nf_1$ and $\nf_3$, and $e_2'$ between $\nf_2$ and $\nf_4$, with the same top/bottom assignment as $(e_1,e_3)$ and $(e_2,e_4)$. Here if some $\nf_i$ is absent, say $e_3$ is a free end at $\of$, then the bond $(e_1,e_3)$ should be replaced by a free end at $\nf_1$, etc.
    \item \emph{Create an O-atom.} Given two edges $e_1'$ between $\nf_1$ and $\nf_3$ and $e_2'$ between $\nf_2$ and $\nf_4$, we can \textbf{create} an O-atom $\of$ by reverting the above process: add a new O-atom $\of$, remove $e_1'$ and $e_2'$, add the bonds $e_j$ between $\of$ and $\nf_j$ with the same top/bottom assignment as $e_1'$ and $e_2'$, and make $e_1$ serial with $e_3$ and $e_2$ serial with $e_4$. We also refer to this as \textbf{joining} $e_1'$ and $e_2'$ by $\of$.
\end{enumerate}
\end{definition}
\begin{definition}[Clusters]\label{def.cluster} Given a molecule $\Mb$, we obtain a new molecule $\Mb'$ by deleting all O-atoms from $\Mb$ as in Definition \ref{def.create_delete}, leaving only C-atoms. The connected components $(\Mb_1, \ldots, \Mb_k)$ of $\Mb'$ are called the \textbf{clusters} of $\Mb$. We say that an O-atom $\of$ \textbf{joins} $\Mb_j$ and $\Mb_{j'}$ if the two particle lines at $\of$ connects the C-atoms in $\Mb_j$ and in $\Mb_{j'}$ respectively. We also define the auxiliary \textbf{cluster graph} to be the graph with nodes being clusters $\Mb_j$, and an edge is drawn between $\Mb_j$ and $\Mb_{j'}$ for each O-atom joining them. See {\color{blue} Figure \ref{fig.overlap3}} for an illustration.
\end{definition}
\begin{figure}[h!]
    \centering
    \includegraphics[width=0.75\linewidth]{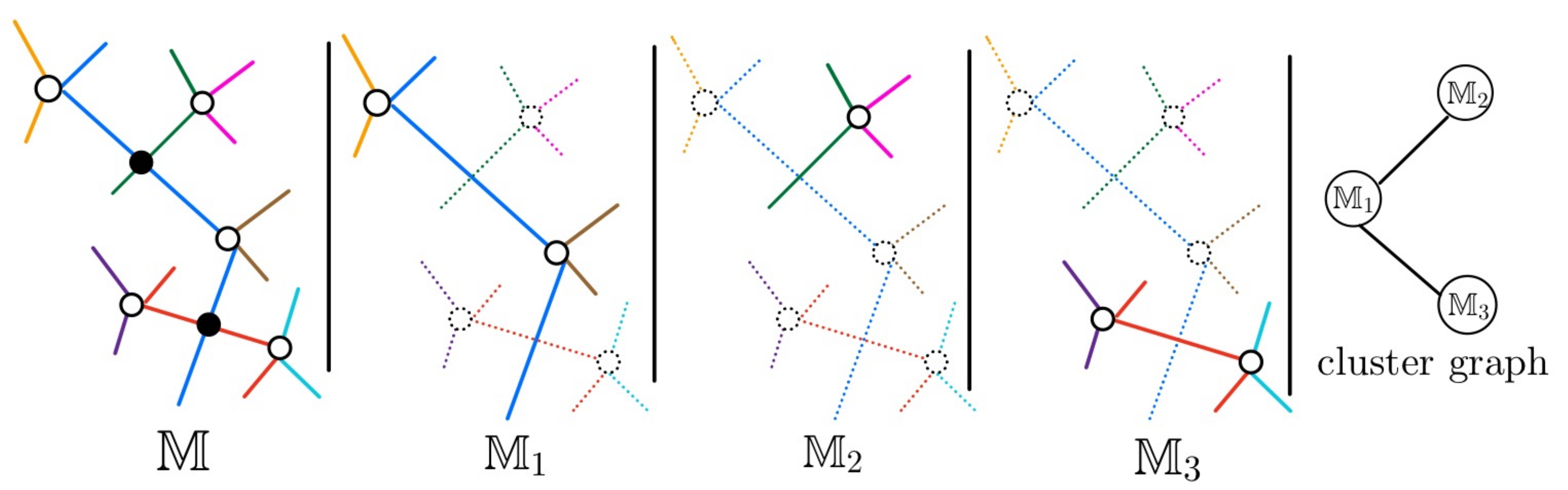}
    \caption{An example of clusters and O-atoms joining them as in Definition \ref{def.cluster}, and the corresponding cluster graph. Here we have 3 clusters $\Mb_j$, with $(\Mb_1,\Mb_2)$ and $(\Mb_1,\Mb_3)$ each joined by an O-atom.}
    \label{fig.overlap3}
\end{figure}

\begin{definition}[C-molecule]
    A molecule consists of only C-atoms are referred to as \textbf{C-molecules}. Clusters are always C-molecules.
\end{definition}
\begin{definition}[Recollision number]\label{def.recollision_number} Given a subset $A\subseteq\Mb$, we define its recollision number by $\rho(A)=\#_{\mathrm{bond}(A)}-|A|+\#_{\mathrm{comp}(A)}$, where $\#_{\mathrm{bond}(A)}$ and $\#_{\mathrm{comp}(A)}$ are respectively the number of bonds between atoms in $A$, and the number of connected components of $A$. We say a set $A$ is \textbf{recollisional} (or equivalently \textbf{cyclic}) if $\rho(A)>0$.
\end{definition}
\begin{remark} Recollision number is an important notion and will be used frequently later. It equals the \emph{circuit rank} which is the \emph{number of (independent) cycles} in $A$. It is easy to see that $\rho(A)\geq 0$ for any $A$.
\end{remark}

\subsection{Molecule prescribed dynamics}\label{sec.molecule_prescribed} In Proposition \ref{prop.cumulant_formula}, the quantity $\Ic\Nc_\Mb$ basically represents an integral over the non-root particle variables, under the restriction that the particle trajectories under certain dynamics topologically reduces to $\Mb$. In order to exclude the possible effects of other collisions, we need to introduce the molecule prescribed dynamics (Definition \ref{def.molecule_truncated_dynamics}), which ignores any possible collision that does not come from the molecule. This is an example of general modified dynamics, which we introduce below.

\begin{definition}[Modified dynamics]\label{def.modified_dynamics} We define the \textbf{modified dynamics}, which is governed by a given property $\mathtt{Pro}$, as follows.    
\begin{enumerate}
\item\label{it.property} \emph{Governing property.} $\mathtt{Pro}$ is a property, which can be expressed in terms of an admissible trajectory $\boldsymbol{\gamma}$ and two particles $\pb,\pb'$. This means that, given any admissible trajectory $\boldsymbol{\gamma}$ and particles $\pb,\pb'$, we can uniquely determine whether $\mathtt{Pro}(\boldsymbol{\gamma},\pb,\pb')$ equals $\mathtt{True}$ or $\mathtt{False}$.
\item\label{it.modified_dynamics} \emph{Modified dynamics.} Given initial configuration $\vz^0=\vz_N^0$, we define the dynamics $\vz^{\mathtt{Pro}}(t)=\vz^{\mathtt{Pro}}(t,\vz^0)$ governed by $\mathtt{Pro}$, as follows:
\begin{enumerate}
\item\label{it.modified_dynamics_0} We have $\vz^{\mathtt{Pro}}(0)=\vz^0$.
\item\label{it.modified_dynamics_1} Suppose $\vz^{\mathtt{Pro}}(t')$ is known for $t'\in[0,t]$. Let $\boldsymbol{\gamma}$ be the admissible trajectory $\vz^{\mathtt{Pro}}(t')\,(t'\in[0,t])$. If two particles satisfy the pre-collisional condition (Definition \ref{def.admissible} (\ref{it.precol})) at time $t$, \emph{and also $\mathtt{Pro}(\boldsymbol{\gamma},\pb,\pb')=\mathtt{True}$} for the molecule $\Mb'$ and the particle lines $\pb,\pb'$ corresponding to these two particles, then they collide as in (\ref{eq.hardsphere}).
\item\label{it.modified_dynamics_2} If a particle does not satisfy (\ref{it.modified_dynamics_1}) with any other particle, then it does not have any collision and moves linearly by (\ref{eq.hardsphere2}).
\end{enumerate}
\item \emph{Flow map and flow operator.} We define the \textbf{flow map} $\Hc^{\mathtt{Pro}}(t)=\Hc_N^{\mathtt{Pro}}(t)$ by
\begin{equation}\label{eq.defH^M}
\Hc^{\mathtt{Pro}}(t):\Rb^{2dN}\to\Rb^{2dN},\quad\Hc^{\mathtt{Pro}}(t)(\vz^0)=\vz^{\mathtt{Pro}}(t)
\end{equation} 
where $\vz^{\mathtt{Pro}}(t)=\vz_N^{\mathtt{Pro}}(t)$ is defined in (\ref{it.modified_dynamics}). The mapping $\Hc_N(t)$ need not be injective; however, we can still define the \textbf{flow operator} $\Sc^{\mathtt{Pro}}(t)$ by
\begin{equation}\label{eq.defS^M}
    \Sc^{\mathtt{Pro}}(t)f(\vz)=\sum_{\vz^0:\Hc^{\mathtt{Pro}}(t)(\vz^0)=\vz}f(\vz^0).
\end{equation}
\item \emph{Collision and overlap.} The notions of \textbf{collisions} and \textbf{overlaps} for any trajectory of modified dynamics, are defined as in Definition \ref{def.admissible}.
\item \emph{Label of particles.} In this definition, we assume that particles are labeled by $[N]=\{1,\cdots,N\}$ (see the convention in Definition \ref{def.notation} (\ref{it.vector})). For a general labeling set $\lambda$, the formulation remains the same, except that the subscript $N$ in expressions like $\vz_N(t)$, $\Hc_N(t)$, $\Sc_N(t)$ etc. is replaced by $\lambda$.
\end{enumerate}
\end{definition}
\begin{remark}\label{rem.modified} We make a few remarks regarding Definition \ref{def.modified_dynamics}. First, up to Lebesgue zero sets, the modified dynamics is well-defined just as Definition \ref{def.hard_sphere}, and produces admissible trajectories in Definition \ref{def.admissible}. Also, we allow $|x_i-x_j|<\varepsilon$ and do not restrict to $\Dc_N$, and only allow collisions under pre-collisional condition.

Second, the modified dynamics can be described as ``collision $\Leftrightarrow$ pre-collisional and property $\mathtt{Pro}$ holds", in particular the dynamics is completely determined by the condition $\mathtt{Pro}$. Examples of this $\mathtt{Pro}$ include the molecule prescribed dynamics (Definition \ref{def.molecule_truncated_dynamics}), extended dynamics (Definition \ref{def.E_dynamics}) and truncated dynamics (Definition \ref{def.T_dynamics}).
\end{remark}

\begin{figure}[h!]
\includegraphics[width=0.6\linewidth]{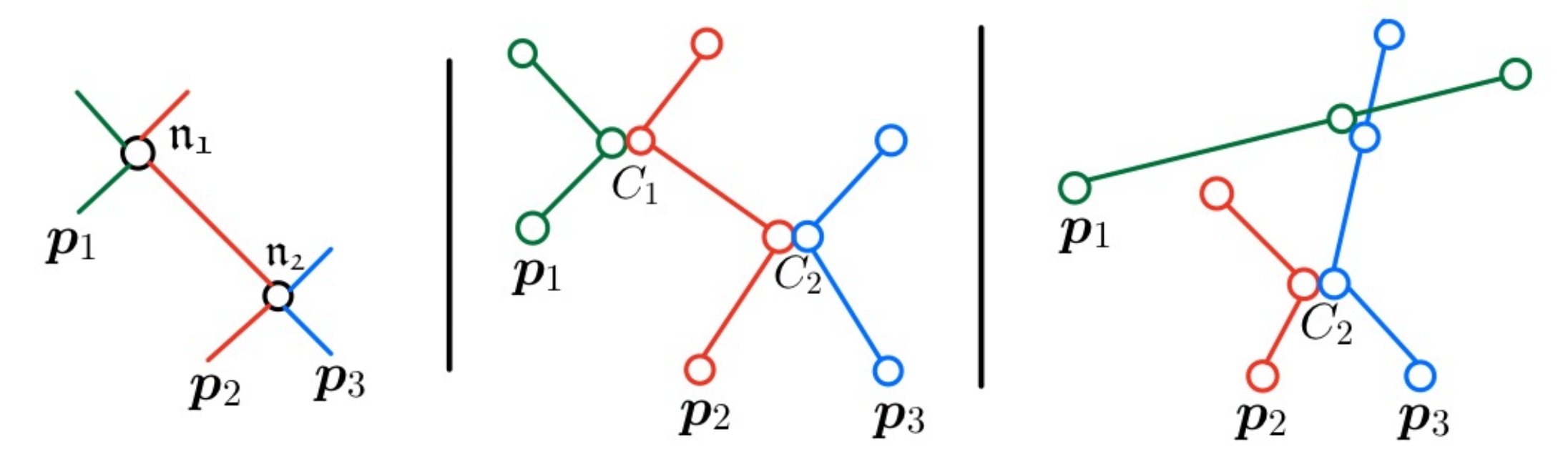}
\caption{Left: a C-molecule $\Mb$. Middle and right: two examples of the $\Mb$-prescribed dynamics. In both cases, the pre-collisional configuration between $\pb_2$ and $\pb_3$ results in a collision $C_2$ corresponding to the C-atom $\nf_2$. In the middle picture we have another collision $C_1$ corresponding to $\nf_1$; however, in the right picture, the pre-collisional configuration of $\pb_1$ and $\pb_3$ does not allow for a collision, as it does not fit in $\Mb$ (according to $\Mb$, the $\pb_1$ should collide with $\pb_2$ and not $\pb_3$).}
\label{fig.prescribe}
\end{figure}

\begin{definition}[Molecule prescribed dynamics]\label{def.molecule_truncated_dynamics} Given a C-molecule $\Mb$, the dynamics \textbf{prescribed by $\Mb$} is a modified dynamics as in Definition \ref{def.modified_dynamics}, with particle set being $p(\Mb)$, and governing property given by
\[\mathtt{Pro}(\boldsymbol{\gamma},\pb,\pb')=\mathtt{True}\Leftrightarrow \textrm{$\Mb'\subseteq\Mb$, and $\pb$ and $\pb'$ intersects at a C-atom $\nf\in\Mb\backslash\Mb'$, and $S_\nf\backslash\{\nf\}\subseteq\Mb'$}.\] Here $\Mb'$ is the unique C-molecule which is topological reduction\footnote{In this topological reduction, the requirement in Definition \ref{def.top_reduction} (\ref{it.reduct_layer}) will be adjusted accordingly: when $t\in[\ell\tau,(\ell+1)\tau]$ then $\Mb'$ will have $\ell$ layers, and those collisions and overlaps corresponding to an atom in layer $\Mb'_{\ell'}$ must occur in $[(\ell'-1)\tau,\ell'\tau]$ if $\ell'<\ell$, and in $[(\ell-1)\tau,t]$ if $\ell'=\ell$.} of $\boldsymbol{\gamma}$ as in Definition \ref{def.top_reduction} (the uniqueness is guaranteed by Remark \ref{rem.top_reduction_uniqueness}), $S_\nf$ is the set of descendants of $\nf$ (Definition \ref{def.molecule_order}), and $\pb$ and $\pb'$ are viewed as particle lines in $\Mb$. See {\color{blue}Figure \ref{fig.prescribe}.}

If $\Mb$ contains O-atoms, the the dynamics prescribed by $\Mb$ is defined to be the dynamics prescribed by $\widetilde{\Mb}$, which is the C-molecule formed by deleting all the O-atoms in $\Mb$.
\end{definition}

\subsection{Statement of the cumulant formula}\label{sec.formula_cumulant_statement} 
\begin{figure}[h!]
\includegraphics[width=0.55\linewidth]{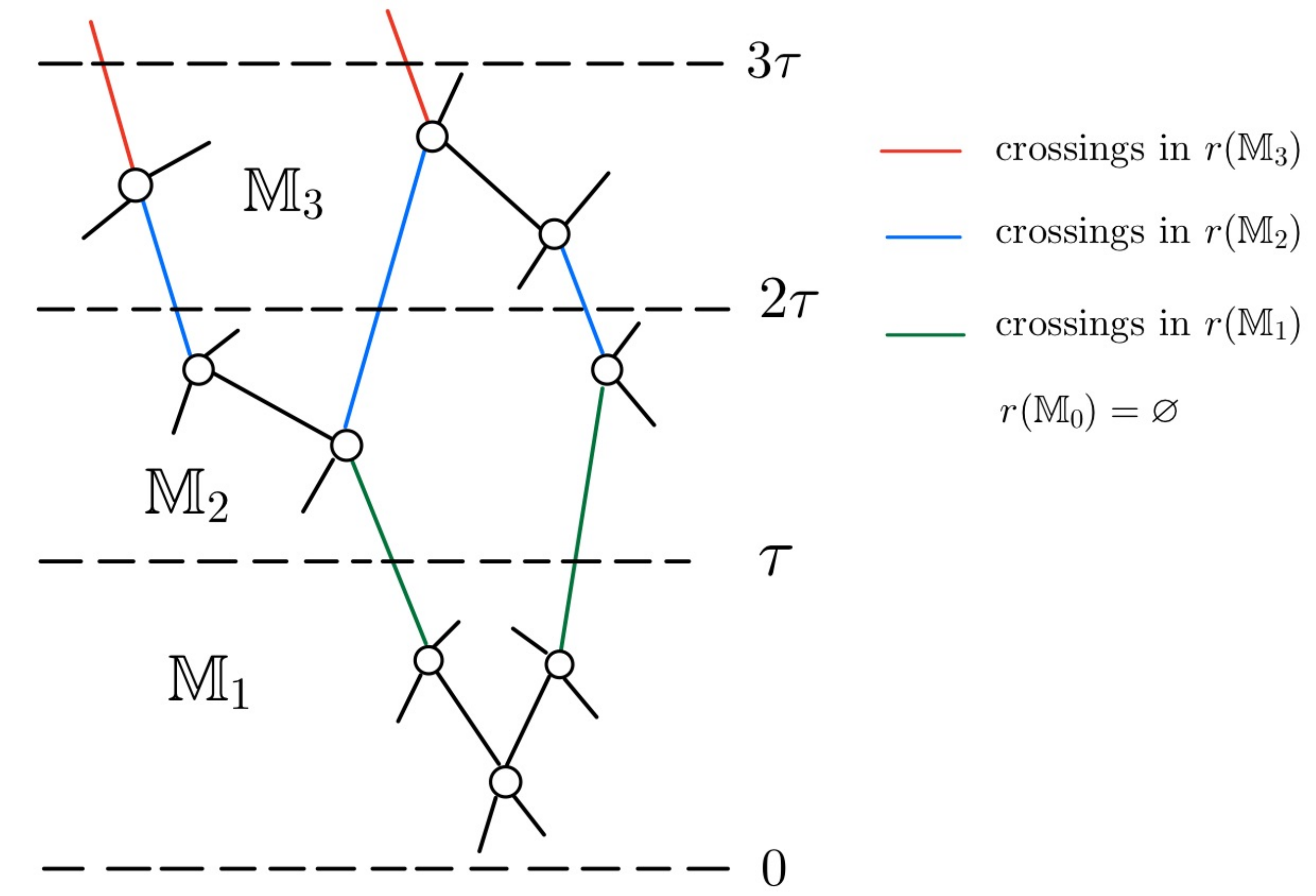}
\caption{Examples of the $p(\Mb_{\ell'})$ and $r(\Mb_{\ell'})$ sets in Definition \ref{def.prsets}. The red edges belong to those particle lines in $r(\Mb_3)$, and intuitively represent that they cross into layer $\Mb_3$; same for blue and green edges
(corresponding to $r(\Mb_2)$ and $r(\Mb_1)$; note that different $r(\Mb_{\ell'})$ may have common particle lines). The $p(\Mb_{\ell'})$ consists of those particle lines that are present in layer $\Mb_{\ell'}$ (including those in $r(\Mb_{\ell'})$). In this picture, Definition \ref{def.set_T_F} (\ref{it.set_F_l_2}) and (\ref{it.set_F_l_4}), and Proposition \ref{prop.mol_axiom} (\ref{it.axiom5}) and (\ref{it.axiom6}) (without showing initial links), are verified.
}
\label{fig.moleculeprop}
\end{figure}

In this section, we provide a rigorous statement of the cumulant formulas (\ref{eq.cumulant_expansion_sec_intro})--(\ref{eq.formula_cumulant_sec_intro}), see Proposition \ref{prop.cumulant_formula}. We start with the definition of the sets $p(\Mb_{\ell'})$ and $r(\Mb_{\ell'})$, of particles and crossings in each layer $\Mb_{\ell'}$.
\begin{definition}\label{def.prsets} Let $\Mb$ be a molecule as in Definition \ref{def.molecule} with $\ell$ layers. For $\ell'\in[\underline{\ell}:\ell]$, define the sets
\begin{align}\label{eq.pset}
p(\Mb_{\ell'})&:=\{\pb:\ell_1[\pb]\leq\ell'\leq\ell_2[\pb]\},\\
\label{eq.rset}
r(\Mb_{\ell'})&:=\{\pb:\ell_1[\pb]\leq\ell'<\ell'+1\leq\ell_2[\pb]\}.
\end{align} 
We refer to the particle lines in $r(\Mb_{\ell'})$ as the \textbf{root particles of layer $\ell'$}. We refer to $p(\Mb_{\ell'})$ as the set of particle lines \textbf{occurring in layer $\ell'$}, and $r(\Mb_{\ell'})$ as the set of particle lines \textbf{crossing into layer $\ell'$}, see {\color{blue}Figure \ref{fig.molecule}} where $r(\Mb_2)=\{1\}$ and $r(\Mb_1)=\{1,2,3\}$ (assume that $p_1$, $p_2$, $p_3$ are labeled by $1$, $2$, $3$). Define also $|\Mb_{\ell'}|_p=|p(\Mb_{\ell'})|$ and $|\Mb_{\ell'}|_r=|r(\Mb_{\ell'})|$ etc. similar to Definition \ref{def.sets_molecule}.
\end{definition}
Before proceeding, we need one more notion of two particle lines being \emph{connected} via a set (including a layer). This is needed in the properties of molecules $\Mb$ in Definition \ref{def.set_T_F}.
\begin{definition}\label{def.connectedvia} Let $\pb_1$ and $\pb_2$ be two particle lines, and $A\subseteq\Mb$ be a set. We say:
\begin{itemize}
\item $\pb_1$ \textbf{intersects} $A$, if an atom exists that belongs to both $\pb_1$ and $A$.
\item $\pb_1$ and $\pb_2$ are \textbf{connected via $A$}, if $\pb_1$ and $\pb_2$ intersect the same component of $A$;
\item $\pb_1$ is \textbf{connected to a cycle within $A$}, if $\pb_1$ intersects a component of $A$ that contains a cycle;
\item $\pb_1$ \textbf{forms an initial link within $A$}, if there exists a particle line $\pb_1'$ which either equals $\pb_1$ or is connected to $\pb_1$ via $A$, such that the bottom end in $\pb_1'$ is involved in an initial link.
\end{itemize}
\end{definition}
Now we can define the objects needed in Proposition \ref{prop.cumulant_formula}. First is the set $\Fc_{\vLambda_{\ell}}$ in \eqref{eq.formula_cumulant_sec_intro}, as well as related sets $\Tc_{\Lambda_{\ell}}$ occurring in the definition of $f^\Ac$ etc.
\begin{definition}[Set of molecules $\Tc$ and $\Fc$]\label{def.set_T_F} Given $\ell\leq\Lf$, $\vLambda_{\ell}:=(\Lambda_{\ell'})_{\ell'\leq\ell}$ and $\Gamma$, we define the following sets of molecules. Recall the notions in Definition \ref{def.molecule} and the sets $p(\Mb_{\ell'})$ and $r(\Mb_{\ell'})$ in Definition \ref{def.prsets}.
\begin{enumerate}
    \item\label{it.set_F} \emph{The set $\mathcal{F}_{\boldsymbol{\Lambda}_\ell}$.} Define $\mathcal{F}_{\boldsymbol{\Lambda}_\ell}$ to be the \textbf{set of equivalence class} $[\Mb]$ of labeled molecules $\Mb$ that satisfy the following properties:
    \begin{enumerate}
        \item\label{it.set_F_l_1} $\Mb$ has lowest layer $\underline{\ell}=1$ and highest layer $\ell$, with $|r(\Mb_{\ell'})|\leq A_{\ell'}$ for each $\ell'\in [1:\ell]$. Each component of $\Mb_{\ell'}$ contains $<\Lambda_{\ell'}$ clusters (see Definition \ref{def.cluster}), and for each cluster $\Mb_j$ we have $|\Mb_j|_p\le \Lambda_{\ell'}$ and $\rho(\Mb_j)\le \Gamma$ (Definitions \ref{def.sets_molecule} and \ref{def.recollision_number}).
        \item\label{it.set_F_l_2} For each $\ell'\in[1:\ell]$, each particle line in $p(\Mb_{\ell'})$ is either connected to a particle line in $r(\Mb_{\ell'})$ via $\Mb_{\ell'}$, or belongs to $r(\Mb_{\ell'})$. Each component of $\Mb_{\ell'}$ contains a particle line in $r(\Mb_{\ell'})$. See {\color{blue}Figure \ref{fig.moleculeprop}.}
        \item\label{it.set_F_l_3} The cluster graph (Definition \ref{def.cluster}) of each $\Mb_{\ell'}$ is a forest. The initial link graph, formed by initial link conditions between particle lines $\pb\in p(\Mb_1)$, is also a forest.
        \item\label{it.set_F_l_4} For $\ell'\in[2:\ell]$, each particle line in $r(\Mb_{\ell'})$ is either connected to another particle line in $r(\Mb_{\ell'})$ via $\Mb_{\ell'}$, or connected to a particle line in $r(\Mb_{\ell'-1})$ via $\Mb_{\ell'}$ (in the sense of Definition  \ref{def.connectedvia}), or belongs to $r(\Mb_{\ell'-1})$. For $\ell'=1$, the set $r(\Mb_{0})$ in the above definition should be replaced by the set of particle lines that are involved in an initial link. See {\color{blue}Figure \ref{fig.moleculeprop}.}
    \end{enumerate}
        \item\label{it.set_F_err} \emph{The set $\mathcal{F}_{\boldsymbol{\Lambda}_\ell}^{\mathrm{err}}$.} $\mathcal{F}_{\boldsymbol{\Lambda}_\ell}^{\mathrm{err}}$ is the set of equivalence classes $[\Mb]$ as in (\ref{it.set_F}) that satisfy conditions \eqref{it.set_F_l_1}--\eqref{it.set_F_l_4}, except that \eqref{it.set_F_l_4} is assumed only for $\ell'\leq\ell-1$. Moreover, for $\ell'=\ell$, instead of  \eqref{it.set_F_l_1} we assume that each component of $\Mb_{\ell}$ contains $\leq\Lambda_{\ell}$ clusters and there exists a unique component that contains \emph{exactly} $\Lambda_{\ell}$ clusters.
    \item\label{it.set_T} \emph{The set $\mathcal{T}_{\Lambda_\ell}$ and $\mathcal{T}_{\Lambda_\ell}^{\mathrm{err}}$.} Given $\Lambda_\ell$ (which is a single number), define $\mathcal{T}_{\Lambda_\ell}$ to be the set of equivalence class $[\Mb]$ of labeled molecules $\Mb$ such that it has only one layer $\ell$ (i.e. $\underline{\ell}=\ell$) and no initial link, has only one root particle labeled by $1$, and satisfies \eqref{it.set_F_l_1}--\eqref{it.set_F_l_3}.
    
    Define also $\mathcal{T}_{\Lambda_\ell}^{\mathrm{err}}$ to be the set of $[\Mb]$ such that  $\Mb$ is obtained from a molecule in $\Tc_{\Lambda_{\ell}}$ by joining two edges $e_1$ and $e_2$ with an extra O-atom.
    \item\label{it.set_F_err_2} \emph{The set $\mathcal{F}_{\boldsymbol{\Lambda}_\ell}^{\mathrm{trc.err}}$.} Define $\mathcal{F}_{\boldsymbol{\Lambda}_\ell}^{\mathrm{trc.err}}$ to be the set of equivalence class $[\Mb]$ of labeled molecules $\Mb$ such that it has $\ell$ layers and satisfies (\ref{it.set_F_l_1})--(\ref{it.set_F_l_4}) for $\ell'\leq \ell-1$. Moreover the layer $\Mb_{\ell}$ is connected (apart from possible empty ends which are root particle lines) and has only C-atoms, and satisfies that 
    \begin{equation}\label{eq.trc_err}|\Mb_{\ell}|_p\leq 3\Lambda_{\ell},\ \mathrm{and}\ \rho(\Mb_{\ell})\leq 2\Gamma;\qquad \mathrm{moreover, \ either\ }|\Mb_{\ell}|_p>\Lambda_{\ell} \mathrm{\ or\ }\rho(\Mb_{\ell})>\Gamma.
    \end{equation}
\end{enumerate}
\end{definition}

\begin{remark}\label{rem.mol_property} We briefly explain the ideas behind the items in Definition \ref{def.set_T_F}. Here (\ref{it.set_F_l_1}) comes from the truncation of dynamics (Definition \ref{def.T_dynamics}) and the truncation in the Penrose argument, and (\ref{it.set_F_l_3}) also comes from the construction in the Penrose argument. Next, (\ref{it.set_F_l_2}) follows from the fact that we only expand the particles that are connected to root particles, and (\ref{it.set_F_l_4}) follows from the condition for the cumulant $E_H$. The other conditions (\ref{it.set_F_err})--(\ref{it.set_F_err_2}) are similar but with modifications associated with single layer or error terms.
\end{remark}

Next, we define the integral expression $|\Ic\Nc_\Mb|$ in \eqref{eq.formula_cumulant_sec_intro}, and its signed variant $\Ic\Nc_\Mb$, in Definitions \ref{def.associated_op_nonlocal}--\ref{def.associated_int}.
\begin{definition}[Associated transport operators]\label{def.associated_op_nonlocal} Let $\Mb$ be a molecule. We define the following objects:
\begin{enumerate}
    \item\label{it.prescribed_flow} \emph{Prescribed dynamics and flow operators.} We use $\vz_\Mb(t)$ to denote the $\Mb$-prescribed dynamics in Definition \ref{def.molecule_truncated_dynamics}, with $\Hc_\Mb(t)$ and $\Sc_\Mb(t)$ etc. defined accordingly as in Definition \ref{def.modified_dynamics}. Note that the subscript $N$ is omitted, as $N=|\Mb|_p$ can be inferred from $\Mb$.
    \item\label{it.associated_ind} \emph{Indicator functions $\mathbbm{1}_{\Mb}$ and $\mathbbm{1}_{\Lc}^\varepsilon$.} We define the indicator function $\mathbbm{1}_{\Mb}(\vz^0)$ such that
    \begin{equation}\label{eq.ind_1M}
    \mathbbm{1}_{\Mb}(\vz^0)=1\Leftrightarrow \textrm{the trajectory $\vz_\Mb(t)$, with initial configuration $\vz^0$, topologically reduces to $\Mb$.}
    \end{equation}
    
    Define also the indicator function $\mathbbm{1}_{\Lc}^\varepsilon$, associated with the initial link structure $\Lc$, by 
    \begin{equation}\label{eq.1_L}
        \mathbbm{1}_\Lc^{\varepsilon} = \prod_{(e_1, e_2)\in \Lc} \mathbbm{1}_{e_1\sim_{\mathrm{in}} e_2}^\varepsilon\qquad\textrm{and}\qquad \mathbbm{1}_{e_1\sim_{\mathrm{in}} e_2}^\varepsilon = \left\{\begin{aligned}
            &1 \quad &&|x_{e_1}-x_{e_2}|\le \varepsilon,
            \\
            &\varepsilon^\upsilon \quad &&|x_{e_1}-x_{e_2}|\ge \varepsilon,
        \end{aligned}\right.
    \end{equation} 
    where $\upsilon = 3^{-d-1}$ (Definition \ref{def.notation} (\ref{it.defnot3})).
    \item \emph{The operator $(\Sc\circ\mathbbm{1})_{\Mb}$ and $|(\Sc\circ\mathbbm{1})_{\Mb}|$.} Define the \textbf{associated transport operators} $(\Sc\circ\mathbbm{1})_\Mb$ and $|(\Sc\circ\mathbbm{1})_\Mb|$ by 
    \begin{equation}\label{eq.associated_op}
        (\Sc\circ\mathbbm{1})_\Mb = (-1)^{|\Mb|_{\mathrm{O}}}\cdot\Sc_\Mb((\ell-\underline{\ell}+1)\tau)\circ\mathbbm{1}_\Mb\quad\text{and}\quad|(\Sc\circ\mathbbm{1})_\Mb| = \Sc_\Mb((\ell-\underline{\ell}+1)\tau)\circ \mathbbm{1}_\Mb,
    \end{equation} where $\ell$ is the highest layer and $\underline{\ell}$ is the lowest layer of $\Mb$, and recall $|\Mb|_{\mathrm{O}}$ is the number of O-atoms in $\Mb$, see Definition \ref{def.sets_molecule}.
\end{enumerate}

\end{definition}

\begin{definition}[Associated integrals $\Ic\Nc_\Mb$ and $|\Ic\Nc_\Mb|$]\label{def.associated_int} Let $\Mb$ be a molecule, let the lowest layer of $\Mb$ be $\underline{\ell}$. Define the \textbf{associated integrals} $\Ic\Nc_\Mb$ and $|\Ic\Nc_\Mb|$ as follows:
    \begin{align}\label{eq.associated_integral_molecule}
        \Ic\Nc_\Mb(\vz_{r(\Mb)}) &\coloneqq \varepsilon^{-(d-1)|\Mb|_{p\backslash r}}\int_{\Rb^{2d|\Mb|_{p\backslash r}}} (\Sc\circ\mathbbm{1})_{\Mb}\, Q\,\mathrm{d}\vz_{(p\backslash r)(\Mb)},\\\label{eq.associated_integral_molecule_abs}
        |\Ic\Nc_\Mb|(\vz_{r(\Mb)}) &\coloneqq \varepsilon^{-(d-1)|\Mb|_{p\backslash r}}\int_{\Rb^{2d|\Mb|_{p\backslash r}}} |(\Sc\circ\mathbbm{1})_{\Mb}|\, |Q|\,\mathrm{d}\vz_{(p\backslash r)(\Mb)},
    \end{align}
    where  $\vz_{(p\backslash r)(\Mb)}:=(\vz_\pb)_{\pb\in (p\backslash r)(\Mb)}$, and
    \begin{equation}\label{eq.extraQ}
    Q:=\mathbbm{1}_\Lc^{\varepsilon}\cdot\prod_{\pb\in p(\Mb)} f^\Ac((\ell_1[\pb]-1)\tau,x_\pb'+(\ell_1[\pb]-\underline{\ell})\tau\cdot v_\pb',v_\pb').
    \end{equation}
\end{definition}

\begin{remark}\label{rem.cumulant_formula_objects} We make a few remarks on Definitions \ref{def.associated_op_nonlocal}--\ref{def.associated_int}.
\begin{enumerate}

\item\label{it.trans_oper} The time $(\ell-\underline{\ell}+1)\tau$ is because $\Mb$ has $\ell-\underline{\ell}+1$ layers, and each layer corresponding to evolution of time $\tau$. Note that, in the support of $\mathbbm{1}_\Mb$, the $\Mb$-prescribed dynamics contains exactly the collisions in $\Mb$. Thus $(\Sc\circ\mathbbm{1})_\Mb$ is given by a fixed composition of explicitly defined changes of variables of form (\ref{eq.hardsphere}).
\item In (\ref{eq.associated_integral_molecule})--(\ref{eq.associated_integral_molecule_abs}) we have two sets of variables, namely $z_\pb$ and $z_\pb'$ for each particle $\pb$. They represent the states of particle $\pb$ at the final time $\ell\tau$ and initial time $(\underline{\ell}-1)\tau$ respectively, in the $\Mb$-prescribed dynamics. In particular they are related by the explicit change of variables in $(\Sc\circ\mathbbm{1})_\Mb$ as in (\ref{it.trans_oper}) (under the support of $\mathbbm{1}_\Mb$ contained in $(\Sc\circ\mathbbm{1})_\Mb$).
\item The integrals $\Ic\Nc_\Mb$ and $|\Ic\Nc_\Mb|$ in (\ref{eq.associated_integral_molecule})--(\ref{eq.associated_integral_molecule_abs}) represent the probability density under the assumption that all collisions and overlaps in $\Mb$ (and initial links) occur as expected. This is clear from the transport operator and indicator function in $\Sc\circ\mathbbm{1}$ as described above. The integral in all \emph{non-root} particles $\vz_{(p\backslash r)(\Mb)}$ leaves the density as a function of the root particles $\vz_{r(\Mb)}$.

Finally, the input function $Q$ in (\ref{eq.extraQ}) contains the initial link/initial cumulant factor $\mathbbm{1}_\Lc^{\varepsilon}$. It also contains the $f^\Ac$ factors (approximate Boltzmann solution) for each particle $\pb$, at time $(\ell_1[\pb]-1)\tau$ corresponding to the start layer $\ell_1[\pb]$. This is the first time $\pb$ becomes relevant in the dynamics prescribed by $\Mb$. Note the linear transport $x_\pb'+(\ell_1[\pb]-\underline{\ell})\tau\cdot v_\pb'$ in (\ref{eq.extraQ}); this is due to the fact that $f^\Ac$ is measures at time $(\ell_1[\pb]-1)\tau$, while $x_\pb'$ represents the state of $\pb$ at time $(\underline{\ell}-1)\tau$.
\end{enumerate}
\end{remark}
Now, we are finally ready to state the main cumulant expansion formula in this section, in Proposition \ref{prop.cumulant_formula} below.
\begin{proposition}\label{prop.cumulant_formula} Recall $\Lambda_\ell$ and $A_\ell$ defined in \eqref{eq.defLambdaseq}, and $t_{\mathrm{fin}}=\Lf\tau$. Assume the initial density function $W_{0,N}(\vz_N)$ is given by \eqref{eq.N_par_ensemble}. Then we can define suitable quantities $\widetilde{f}_s(\ell\tau,\vz_s)$ (for all $\ell\in[\Lf]$ and $s\leq A_\ell$) and $f_s^{\mathrm{err}}(t_{\mathrm{fin}},\vz_s)$ (for $s\leq A_\Lf$), see (\ref{eq.truncated_correlation})--(\ref{eq.truncated_error}), such that we have the decomposition
\begin{equation}\label{eq.f_s_dec}
    f_s(t_{\mathrm{fin}}, \boldsymbol{z}_s) = \widetilde{f}_s(t_{\mathrm{fin}}, \boldsymbol{z}_s) + f_s^{\mathrm{err}}(t_{\mathrm{fin}}, \boldsymbol{z}_s),
\end{equation}
and for each $\ell\in[0:\Lf]$ and $s\leq A_\ell$, we have
\begin{equation}\label{eq.cumulant_expansion}
    \widetilde{f}_{s}(\ell\tau,\vz_s)=\sum_{H\subseteq [s]}(f^\Ac(\ell\tau))^{\otimes([s]\backslash H)}\cdot E_H(\ell\tau,\vz_H)+\mathrm{Err}(\ell\tau, \vz_s).
\end{equation} 
In (\ref{eq.f_s_dec})--\eqref{eq.cumulant_expansion} we have the following estimates:
\begin{enumerate}
    \item\label{it.cumulant_formula_1} $(f^\Ac(\ell\tau))^{\otimes([s]\backslash H)}$ is a function of $\vz_{[s]\backslash H}$, which is the tensor product of the single variable function $f^\Ac(\ell\tau, z)$. This function is defined inductively in $\ell$ by $f^\Ac(0)=f_0$ and
    \begin{equation}\label{eq.fAterm}
        f^{\mathcal{A}}(\ell\tau,z) = \sum_{[\Mb]\in \mathcal{T}_{\Lambda_\ell}}\Ic\Nc_\Mb(z) + f^{\mathcal{A},\mathrm{err}}(\ell\tau,z),\quad |f^{\mathcal{A},\mathrm{err}}(\ell\tau,z)| \le|\log\varepsilon|^{C^*} \sum_{[\Mb]\in \mathcal{T}_{\Lambda_\ell}^{\mathrm{err}}} |\Ic\Nc_\Mb|(z).
    \end{equation}
    \item The cumulant $E_H(\ell \tau)$ is bounded by (where $\ell\geq 1$, for $\ell=0$ see Proposition \ref{prop.initial_cumulant}) 
    \begin{equation}\label{eq.cumulant_formula}
        |E_H(\ell\tau,\vz_H)| \le 
        \sum_{[\Mb]\in \mathcal{F}_{\boldsymbol{\Lambda}_\ell},r(\Mb) = H} |\Ic\Nc_\Mb|(\vz_H).
    \end{equation}
    \item We have $\mathrm{Err}=\mathrm{Err}^1+\mathrm{Err}^2$ (where $\ell\geq 1$, for $\ell=0$ we have $\mathrm{Err}(0)=0$), and
    \begin{align}\label{eq.cumulant_formula_err}
        |\mathrm{Err}^1(\ell\tau, \vz_s)|&\leq \sum_{[\Mb]\in \mathcal{F}_{\boldsymbol{\Lambda}_\ell}^{\mathrm{err}},r(\Mb) = [s]}|\Ic\Nc_\Mb|(\vz_s),\\
\label{eq.cumulant_formula_err2}
        \|\mathrm{Err}^2(\ell\tau)\|_{L^1} &\le \varepsilon^{-2(d-1)\Lambda_{\ell}^2A_{\ell}}\cdot \|\mathrm{Err}((\ell-1)\tau)\|_{L^1}.
    \end{align}
    \item The truncation error $f^{\mathrm{err}}(t_{\mathrm{fin}}, \vz_s)$ in (\ref{eq.f_s_dec}) satisfies
    \begin{equation}\label{eq.cumulant_formula_err_trunc}
        \|f^{\mathrm{err}}(t_{\mathrm{fin}}, \vz_s)\|_{L^1}\leq \sum_{\ell=1}^\Lf\bigg(\sum_{[\Mb]\in \mathcal{F}_{\boldsymbol{\Lambda}_\ell}^{\mathrm{trc.err}},r(\Mb) = [s]}\big\||\Ic\Nc_\Mb|(\vz_s)\big\|_{L^1}+\varepsilon^{-8(d-1)\Lambda_\ell}\cdot\|\mathrm{Err}((\ell-1)\tau)\|_{L^1}\bigg).
    \end{equation}
\end{enumerate}
\end{proposition}
\begin{proof}
    The proof will occupy Section \ref{sec.truncation_large_molecule}--\ref{sec.proof_cumulant_formula}, and will be finished in Section \ref{sec.proof_cumulant_formula_final}.
\end{proof}
\begin{remark}\label{rem.label_root} The quantity $\widetilde{f}_s$ is defined by the truncated dynamics, see Section \ref{sec.truncation_large_molecule}. Note that the equation \eqref{eq.fAterm} defines a recurrence relation for $f^{\mathcal{A}}$, because both $\Ic\Nc_\Mb$ and $|\Ic\Nc_\Mb|$ are expressed as integrals involving products of $f^{\mathcal{A}}$. These expressions are also constant on equivalence classes $[\Mb]$.

The summation (\ref{eq.cumulant_formula}) (and similarly (\ref{eq.cumulant_formula_err}) and (\ref{eq.cumulant_formula_err_trunc})) involves the condition $r(\Mb)=H$. Fixing the labels of these root particles may lead to an additional $|H|!$ factor, which is usually harmless (but requires some extra discussion in Section \ref{sec.error}). Apart from this, the equivalence classes $[\Mb]$ of labeled molecules $\Mb$ are then the same as unlabeled molecules, see Remark \ref{rem.label}.
\end{remark}

\section{Truncation of dynamics}\label{sec.truncation_large_molecule} We now start the proof of Proposition \ref{prop.cumulant_formula}. In this section we will define the truncated dynamics and establish the decomposition (\ref{eq.f_s_dec}). This is needed in order to control the number of particles and recollisions in each cluster in each time layer, that occur in the molecule expansion (for example in Definition \ref{def.set_T_F} \eqref{it.set_F_l_1}, the number of particles in layer $\ell'$ is bounded by $\Lambda_{\ell'}^2A_{\ell'}$). 

Recall from \eqref{eq.time_t_ensemble} and \eqref{eq.s_par_cor} that $f(t, \boldsymbol{z}_s)$ is given by 
\begin{equation}\label{eq.def_f_s'}
    f(t, \boldsymbol{z}_s) = \varepsilon^{(d-1)s}\sum_{n = 0}^{\infty} \frac{1}{n!} \int (\mathcal{S}_{s+n}(t)W_{0, s+n})(\boldsymbol{z}_{s+n})
    \, \mathrm{d}\boldsymbol{z}_{[s+1:s+n]},
\end{equation} 
so the correlation $f(t,\vz_s)$ is determined by the collisions in the dynamics $\Sc_{s+n}(t)$, for which we will define the truncation.

\subsection{The extended and truncated dynamics}\label{sec.E_T_dynamics} In this subsection, we define the extended and truncated dynamics as two instances of modified dynamics in Definition \ref{def.modified_dynamics}.

\begin{definition}[The extended dynamics]\label{def.E_dynamics} We define the \textbf{extended dynamics} (or \textbf{E-dynamics}), denoted by $\vz^E(t)$, to be the modified dynamics in Definition \ref{def.modified_dynamics}, with the governing property $\mathtt{Pro}(\boldsymbol{\gamma},\pb,\pb')\equiv\mathtt{True}$. 
\end{definition}
\begin{remark} In the E-dynamics, we have a collision whenever the pre-collisional condition is satisfied. In particular, if the distance of any two particles becomes $\geq\varepsilon$, then this separation is preserved thereafter.

The E-dynamics is a natural extension of the original dynamics (Definition \ref{def.hard_sphere}) from $\Dc_N$ to the whole space $\Rb^{2dN}$. If we restrict the initial data of the E-dynamics to $\Dc_N$, then the trajectory stays in $\Dc_N$ forever, and the E-dynamics coincides with the O-dynamics (see Proposition \ref{prop.dynamics_E_T_property} below).
\end{remark}

\begin{definition}[The truncated dynamics]\label{def.T_dynamics} Let $(\Lambda,\Gamma)$ be fixed. We define the \textbf{truncated dynamics}  (or $(\Lambda,\Gamma)$-, or \textbf{T-dynamics}), denoted by $\vz^{\Lambda,\Gamma}(t)$, to be the modified dynamics in Definition \ref{def.modified_dynamics}, with the governing property $\mathtt{Pro}(\boldsymbol{\gamma},\pb,\pb')$ given by
\begin{multline}\label{eq.trunc_dym}
\mathtt{Pro}(\boldsymbol{\gamma},\pb,\pb')=\mathtt{True}\Leftrightarrow\textrm{either ($\Mb_i=\Mb_j$ and $\rho(\Mb_i)\leq\Gamma-1$),}\\
\textrm{or ($\Mb_i\neq \Mb_j$ and $|\Mb_i|+|\Mb_j|\leq\Lambda$ and $\rho(\Mb_i)+\rho(\Mb_j)\leq\Gamma$)},
\end{multline} where $\Mb_i$ (resp. $\Mb_j$) is the connected components of $\Mb$ in which the particle line $\pb$ (resp. $\pb'$) occurs, and $\Mb$ is the unique C-molecule which is topological reduction of $\boldsymbol{\gamma}$.
\end{definition}
\begin{remark}\label{rem.T_dyn} In the $(\Lambda,\Gamma)$-truncated dynamics, we force the property that each cluster has size $\leq\Lambda$ and recollision number $\leq\Gamma$ (Proposition \ref{prop.dynamics_E_T_property} (\ref{it.dynamics_E_T_property_4})), by neglecting any collisions that may violate this property.
\end{remark}
\begin{proposition}\label{prop.dynamics_E_T_property} The E- and $(\Lambda,\Gamma)$-dynamics (Definitions \ref{def.E_dynamics} and \ref{def.T_dynamics}) satisfy the following properties.
    \begin{enumerate}
        \item\label{it.dynamics_E_T_property_1} Up to some Lebesgue zero set of $\vz^0$, the E- and T- dynamics in Definitions \ref{def.E_dynamics} and \ref{def.T_dynamics} are well-defined and produces admissible trajectories in Definition \ref{def.admissible}. 
        \item\label{it.dynamics_E_T_property_2} Define $(\Hc^E,\Sc^E)$ and $(\Hc^{\Lambda,\Gamma},\Sc^{\Lambda,\Gamma})$ as special cases of (\ref{eq.defH^M}) and (\ref{eq.defS^M}). Then, when the initial data $\vz^0\in \Dc_N$, the E-dynamics coincides with the O-dynamics in Definition \ref{def.hard_sphere}, i.e. $\vz^E(t) = \vz(t)$ and $\Hc_N^{E}(t) (\vz^0) = \Hc_N(t) (\vz^0)$. Moreover, we have $\Sc_N^{E}(t)f = \Sc_N(t)f$ if $f$ is supported in $\Dc_N$.
        \item\label{it.dynamics_E_T_property_3} The $\Hc^E_N(t)$ and $\mathcal{S}_N^E(t)$ defined in (\ref{it.dynamics_E_T_property_2}) satisfy the following flow or semi-group property for $t,s\geq 0$:
        \begin{equation}\label{eq.semigroup}
            \mathcal{H}_N^E(t+s)=\mathcal{H}_N^E(t)\mathcal{H}_N^E(s),\qquad\mathcal{S}_N^E(t+s)=\mathcal{S}_N^E(t)\mathcal{S}_N^E(s).
        \end{equation}
 In general, the $\Hc^{\Lambda,\Gamma}_N(t)$ and $\mathcal{S}^{\Lambda,\Gamma}_N(t)$ might not satisfy the flow property.
 \item\label{it.dynamics_E_T_property_3.5} Both transport operators $\mathcal{S}_N^E(t)$ and $\mathcal{S}^{\Lambda,\Gamma}_N(t)$ preserve the $L^1$ integral:
 \begin{equation}\label{eq.L1pres}\int (\Sc_N^E(t)f)\,\mathrm{d}x\mathrm{d}v=\int (\Sc_N^{\Lambda,\Gamma}(t)f)\,\mathrm{d}x\mathrm{d}v=\int f\,\mathrm{d}x\mathrm{d}v.\end{equation}
        \item \label{it.dynamics_E_T_property_4} In the $(\Lambda,\Gamma)$-dynamics, suppose a trajectory topologically reduces to a C-molecule $\Mb$. Then each connected component $\Mb_j$ of $\Mb$ satisfies $|\Mb_j|_p\leq\Lambda$ and $\rho(\Mb_j)\leq\Gamma$, (recall the notions of) Moreover, for O-, E- and $(\Lambda,\Gamma)$-dynamics, the number of collisions $|\Mb|$ has an absolute upper bound that depends only on $d$ and $N$.
    \end{enumerate}
\end{proposition}
\begin{proof} See Appendix \ref{app.aux}.
\end{proof}

\subsection{The truncated domain}\label{sec.truncated_domain} In this subsection, we introduce the domain $\mathcal D^{\Lambda,\Gamma}_N$, which characterizes the initial data such that the truncated dynamics coincide with the extended dynamics.

\begin{definition}[The truncated and error domain $\mathcal D^{\Lambda,\Gamma}_N$ and $\Ec_N^{\Lambda,\Gamma}$]\label{def.truncated_domain} Given $(\Lambda,\Gamma)$, we define the \textbf{truncated domain} $\mathcal D^{\Lambda,\Gamma}_N$ as follows:
\begin{equation}\label{eq.truncated_domain}
    \mathcal D^{\Lambda,\Gamma}_N = \bigcup_{k = 1}^N\bigcup_{\substack{\{\lambda_1, \cdots, \lambda_k\} \in \mathcal{P}^k_N\\|\lambda_j|\le \Lambda}}\bigg\{\vz^0\in \Rb^{2dN}: \textrm{$p(\Mb_j)=\lambda_j$ and $\rho(\Mb_j)\leq\Gamma$ for each $1\leq j\leq k$}\bigg\},
\end{equation} where $\Mb_j\,(1\leq j\leq k)$ are the components of the (unique) C-molecule topological reduction of the admissible trajectory $\vz^E(t)\,(t\in[0,\tau])$ with initial configuration $\vz^0$. We also refer to these $\Mb_j$ as \textbf{clusters}.

We also define the \textbf{error domain} $\Ec_N^{\Lambda,\Gamma}$ as 
\begin{equation}\label{eq.error_domain}
    \Ec_N^{\Lambda,\Gamma}:=\Rb^{2dN}\backslash\Dc_N^{\Lambda,\Gamma}.
\end{equation}
Their indicator functions are denoted by $\mathbbm{1}_{\Dc_N^{\Lambda,\Gamma}}$ and $\mathbbm{1}_{\Ec_N^{\Lambda,\Gamma}}$, and we have $\mathbbm{1}_{\Ec_N^{\Lambda,\Gamma}} = 1-\mathbbm{1}_{\Dc_N^{\Lambda,\Gamma}}$.
\end{definition}

The following proposition ensures that the E- and $(\Lambda,\Gamma)$-dynamics coincide on $\Dc_N^{\Lambda,\Gamma}$.

\begin{proposition}\label{prop.dynamics_comparison} For each $\vz^0\in \Dc_N^{\Lambda,\Gamma}$ and $t\in[0,\tau]$ we have $\Hc_N^{\Lambda,\Gamma}(t)(\vz^0)=\Hc_N^{E}(t)(\vz^0)$. Moreover, if $f$ is supported in $\Dc_N^{\Lambda,\Gamma}$, then we have $\Sc_N^{\Lambda,\Gamma}(t)f=\Sc_N^{E}(t)f$. For general $f$ we have
\begin{equation}\label{eq.S_trunc_equal_to_S}
    \Sc_N^{\Lambda,\Gamma}(t)\big(\mathbbm{1}_{\Dc_N^{\Lambda,\Gamma}}\cdot f\big)=\Sc^{E}_N(t) \big(\mathbbm{1}_{\Dc_N^{\Lambda,\Gamma}}\cdot f\big).
\end{equation}
\end{proposition}
\begin{proof}
Suppose $\vz^0\in \Dc_N^{\Lambda,\Gamma}$. By definition of $\Dc_N^{\Lambda,\Gamma}$, in the whole E-dynamics $\vz_N^E(t)$ for $t\in[0,\tau]$, the size $|\Mb_j|$ of each cluster $\Mb_j$ (under the E-dynamics) never exceeds $\Lambda$, and the recollision number $\rho(\Mb_j)$ never exceeds $\Gamma$. Therefore, whenever a collision happens in the E-dynamics, we must also have a collision in the $(\Lambda,\Gamma)$-dynamics. This implies that the E-dynamics $\vz_N^E(t)$ and $(\Lambda,\Gamma)$-dynamics $\vz_N^{\Lambda,\Gamma}(t)$ have exactly the same collisions, and thus coincide because they have the same initial data.

Now assume $f$ is supported in $\Dc_N^{\Lambda,\Gamma}$, then 
\begin{equation} \label{eq.dynamics_comparison}
    (\Sc_N^{\Lambda,\Gamma}(t)f)(\vz)=\sum_{\vz^0\in \Dc_{N}^{\Lambda,\Gamma}:\Hc_N^{\Lambda,\Gamma}(t)(\vz^0)=\vz}f(\vz^0)=\sum_{\vz^0\in \Dc_{N}^{\Lambda,\Gamma}:\Hc^{E}_N(t)(\vz^0)=\vz}f(\vz^0)=(\Sc_N^{E}(t)f)(\vz).
\end{equation} 
Then \eqref{eq.S_trunc_equal_to_S} follows because $\mathbbm{1}_{\Dc_N^{\Lambda,\Gamma}}\cdot f$ is supported in $\Dc_N^{\Lambda,\Gamma}$.
\end{proof}

\subsection{The truncated correlation functions}\label{sec.trunc_cor}
In this subsection we establish the decomposition (\ref{eq.f_s_dec}) by applying the $(\Lambda_\ell,\Gamma)$-truncated dynamics on each time interval $[(\ell-1)\tau,\ell\tau]$. 

Recall the $\Lf$ and $\tau$ defined in Definition \ref{def.notation} (\ref{it.time_layer}). By the semi-group property of $\Sc_N^E(t)$ (Proposition \ref{prop.dynamics_E_T_property} \eqref{it.dynamics_E_T_property_3}), we have $\mathcal{S}_N^{E}(t_{\mathrm{fin}}) = (\Sc^{E}_N(\tau))^\Lf$. Moreover, since $W_{0,N}$ is supported in $\Dc_N$ (see \eqref{eq.N_par_ensemble}), by Proposition \ref{prop.dynamics_E_T_property} \eqref{it.dynamics_E_T_property_2}, we have $\mathcal{S}_N(t_{\mathrm{fin}}) W_{0,N} = \mathcal{S}_N^{E}(t_{\mathrm{fin}}) W_{0,N}$. Combining the above two facts, we get
\begin{equation}\label{eq.dynamics_dec}
        W_N(t_{\mathrm{fin}}) = \mathcal{S}_N(t_{\mathrm{fin}}) W_{0,N} = \mathcal{S}_N^{E}(t_{\mathrm{fin}}) W_{0,N} =\big((\Sc^{E}_N(\tau))^\Lf\big) W_{0,N}.
\end{equation}  

Recall the sequence $(A_\ell)$ and $(\Lambda_\ell)$ defined in \eqref{eq.defLambdaseq} and $\Gamma$ introduced in Definition \ref{def.notation}. We then define $\Dc_N^\ell:=\Dc_N^{\Lambda_\ell,\Gamma}$ and $\Ec_N^\ell:=\Ec_N^{\Lambda_\ell,\Gamma}$. By definitions of the two indicator functions $\mathbbm{1}_{\Dc_N^\ell}$ and $\mathbbm{1}_{\Ec_N^\ell}$ (Definition \ref{def.truncated_domain}), we have
\begin{equation}\label{eq.1_D^l_E^l_dec}
   1=\mathbbm{1}_{\Dc_N^\ell}+\mathbbm{1}_{\Ec_N^\ell}. 
\end{equation} 

Define $\Hc_N^\ell:=\Hc_N^{\Lambda_\ell,\Gamma}$ and $\Sc_N^\ell:=\Sc_N^{\Lambda_\ell,\Gamma}$, which are given by the $(\Lambda,\Gamma)$-dynamics with $\Lambda=\Lambda_{\ell}$. By Proposition \ref{prop.dynamics_comparison}, we know that $\mathcal{S}^{E}_{N}(\tau) \circ\mathbbm{1}_{\mathcal D_N^{\ell}} = \Sc_N^\ell(\tau) \circ\mathbbm{1}_{\Dc_N^\ell}$ which implies that $(\mathcal{S}_{N}^{E}(\tau) -\Sc_N^\ell(\tau))\circ\mathbbm{1}_{\mathcal D_N^{\ell}} = 0 $. We then have
\begin{equation}\label{eq.S1_split}
\begin{aligned}
        \mathcal{S}_{N}^{E}(\tau)&= \mathcal{S}_{N}^\ell(\tau)+(\Sc^{E}_N(\tau)-\Sc_N^\ell(\tau))\circ(\mathbbm{1}_{\Dc_N^\ell}+\mathbbm{1}_{\Ec_N^\ell})\\
        &= \Sc_N^\ell(\tau) + (\Sc^{E}_N(\tau)-\Sc_N^\ell(\tau))\circ\mathbbm{1}_{\Ec_N^\ell}=\Sc_N^\ell(\tau)+\Sc_{\ell,N}^{\mathrm{err}}(\tau),
        \end{aligned}
\end{equation}
where we define the operator $\mathcal{S}^{\mathrm{err}}_{\ell,N}(\tau)$ as 
\begin{equation}\label{eq.S^err}
    \Sc_{\ell,N}^{\mathrm{err}}(\tau):=(\Sc^{E}_N(\tau)-\Sc_N^\ell(\tau))\circ\mathbbm{1}_{\Ec_N^\ell}.
\end{equation}

Now, by plugging \eqref{eq.S1_split} into \eqref{eq.dynamics_dec}, we get
\begin{equation}\label{eq.iterate_truncation}
    \begin{split}
        W_N(t_{\mathrm{fin}})&=\big(\Sc_N^E(\tau)\big)^\Lf W_{0,N}
        \\ 
        &=\bigg(\prod_{\ell=1}^\Lf \Sc_N^\ell(\tau)\bigg)W_{0,N}+\sum_{\ell=1}^{\Lf}\bigg[(\Sc^{E}_N(\tau))^{\Lf-\ell+1}\prod_{\ell'=1}^{\ell-1}\Sc_N^{\ell'}(\tau)-(\Sc^{E}_N(\tau))^{\Lf-\ell}\prod_{\ell'=1}^{\ell}\Sc_N^{\ell'}(\tau)\bigg]W_{0,N}
        \\
        &=\bigg(\prod_{\ell=1}^\Lf\Sc_N^{\ell}(\tau)\bigg)W_{0,N}+\sum_{\ell=1}^\Lf\bigg((\Sc^{E}_N(\tau))^{\Lf-\ell}\circ\mathcal{S}^{\mathrm{err}}_{\ell,N}(\tau)\circ\prod_{\ell'=1}^{\ell-1}\Sc_N^{\ell'}(\tau)\bigg)W_{0,N}.
    \end{split}
\end{equation} 
Here in the products $\prod_{\ell=1}^\Lf$ the factors with \emph{larger $\ell$ should appear on the left}. For each $\ell$, we define 
\begin{equation}\label{eq.W_widetilde_l}
    \widetilde{W}_N(\ell\tau) = \bigg(\prod_{\ell'=1}^{\ell}\mathcal{S}^{\ell'}_N(\tau)\bigg) W_{0,N},
\end{equation}
then \eqref{eq.iterate_truncation} we be written as $W_N(t_{\mathrm{fin}})=\widetilde{W}_N(t_{\mathrm{fin}})+W_N^{\mathrm{err}}(t_{\mathrm{fin}})$, where
\begin{equation}\label{eq.W_widetilde}
    \widetilde{W}_N(t_{\mathrm{fin}}) = \widetilde{W}_N(\mathfrak{L}\tau),
\end{equation}
\begin{equation}\label{eq.W_err}
    W^{\mathrm{err}}_N(t_{\mathrm{fin}}) = \sum_{\ell = 1}^{\mathfrak{L}}\bigg((\Sc^{E}_N(\tau))^{\Lf-\ell}\circ\mathcal{S}^{\mathrm{err}}_{\ell,N}(\tau)\bigg)\widetilde{W}_N((\ell-1)\tau).
\end{equation}
\begin{definition}[Truncated correlation function]\label{def.ftrunc} We define the \textbf{truncated correlation function} $\widetilde{f}_s$, and the \textbf{error correlation function} $f_s^{\mathrm{err}}$, by
\begin{align}\label{eq.truncated_correlation}
    \widetilde{f}_s(\ell\tau, \boldsymbol{z}_s) &= \varepsilon^{s(d-1)}\sum_{n = 0}^{\infty} \frac{1}{n!} \int_{\Rb^{2dn}} \widetilde{W}_{s+n}(\ell\tau,\boldsymbol{z}_{s+n})\,\mathrm{d}\vz_{[s+1:s+n]},\\
    \label{eq.truncated_error}
    f_s^{\mathrm{err}}(t_{\mathrm{fin}}, \boldsymbol{z}_s) &= \varepsilon^{s(d-1)}\sum_{n = 0}^{\infty} \frac{1}{n!} \int_{\Rb^{2dn}} W^{\mathrm{err}}_{s+n}(t_{\mathrm{fin}},\boldsymbol{z}_{s+n})  \,\mathrm{d}\vz_{[s+1:s+n]},
\end{align} where $\widetilde{W}_{s+n}$ and $W_{s+n}^{\mathrm{err}}$ are defined in (\ref{eq.W_widetilde_l}) and (\ref{eq.W_err}). This completes the proof of (\ref{eq.f_s_dec}).
\end{definition}

\section{Proof of the cumulant formula}\label{sec.proof_cumulant_formula} We already have (\ref{eq.f_s_dec}), with $\widetilde{f}_s$ and $f_s^{\mathrm{err}}$ defined in (\ref{eq.truncated_correlation}) and (\ref{eq.truncated_error}) respectively. To complete Proposition \ref{prop.cumulant_formula}, in this section, we prove the upper bound for $f_s^{\mathrm{err}}(t_{\mathrm{fin}}, \boldsymbol{z}_s)$ in \eqref{eq.cumulant_formula_err_trunc}, and the cumulant expansion formula \eqref{eq.cumulant_expansion} for $\widetilde{f}_s(\ell\tau)$, with bounds (\ref{eq.fAterm})--(\ref{eq.cumulant_formula_err2}). This is done in several steps.

Recall that the definition of $\widetilde{W}_N$ and $\widetilde{f}_s$ involves the transport operator $\Sc_N^{\Lambda_{\ell'},\Gamma}(\tau)$ defined by the truncated dynamics. In Sections \ref{sec.cluster_expansion}--\ref{sec.penrose} we will prove a crucial decomposition formula for this operator, namely Proposition \ref{prop.S_N_decomposition}. Then in Sections \ref{sec.f_to_f}--\ref{sec.E_to_E}, we will use this decomposition to prove Proposition \ref{prop.cumulant_single}. This is the key inductive step, which establishes the cumulant decomposition (\ref{eq.cumulant_expansion}) at time $\ell\tau$ provided the same expansion holds at time $(\ell-1)\tau$, and expresses the cumulant $E_H(\ell\tau)$ in terms of $E_H((\ell-1)\tau)$ and a molecule $\Mb_{\ell}$ with single layer $\ell$. Then in Sections \ref{sec.initial_cumulant}--\ref{sec.proof_cumulant_formula_final}, by induction in $\ell$ using Proposition \ref{prop.cumulant_single}, and using the initial cumulant formula (time $t=0$, see Proposition \ref{prop.initial_cumulant}), we complete the proof of Proposition \ref{prop.cumulant_formula}.

\subsection{Cluster expansion}\label{sec.cluster_expansion}

In this section, we perform the cluster expansion for the $(\Lambda,\Gamma)$-truncated dynamics and prove Lemma \ref{lem.dynamics_factorization} by introducing the C-molecules, which allow us to write the transport operator $\Sc_N^{\Lambda,\Gamma}(t)$ as a summation of tensor products $\prod_{j=1}^k\Sc_{\Mb_j}(t)$ associated with dynamics prescribed by clusters $\Mb_j$.

\begin{definition}[Partial dynamics]\label{def.partial_dynamics} Recall $\vz=\vz_N=(z_1,\cdots,z_N)$ as in Definition \ref{def.notation} (\ref{it.vector}). Let $\lambda\subseteq[N]$ and $\vz_\lambda=(z_j)_{j\in\lambda}$. Given initial configuration $\vz(0)=\vz^0$, let $\vz(t)$ be a modified dynamics (for example $\vz^E(t)$ or $\vz^{\Lambda,\Gamma}(t)$), then we can define the following notions.
\begin{enumerate}
    \item\label{it.dynamics_partial} \emph{The partial dynamics $\vz_\lambda(t)$.} Consider $\lambda\subseteq[N]$ as a set of particles. We define the \textbf{partial dynamics} $\vz_\lambda(t)$ to be the dynamics following the same evolution rule but only taking into account particles in $\lambda$, and with initial configuration $\vz_\lambda^0=(z_j^0)_{j\in\lambda}$. When $\lambda = [N]$, we refer to it as the \textbf{full dynamics}.
    \item\label{it.partial_flow_map} \emph{The partial flow maps $\Hc_\lambda(t)$ and $\Sc_\lambda(t)$.} We define the \textbf{partial flow maps} $\Hc_\lambda(t)$ by $\Hc_\lambda(t)(\vz_\lambda^0)=\vz_\lambda(t)$. It can be extended to $\vz^0 = \vz^0_N$ by  $\Hc_\lambda(t)\vz^0 = \Hc_\lambda(t)(\vz_\lambda^0,\vz^0_{\lambda^c})=(\vz_\lambda(t),\vz^0_{\lambda^c})$. We define the \textbf{partial flow operators} $\Sc_\lambda(t)$ as follows, 
    \begin{equation}\label{eq.transport_op_partial}
    \begin{gathered}
        (\Sc_\lambda(t)f)(\vz_\lambda) = \sum_{\vz_\lambda^0:\Hc_\lambda(t)(\vz_\lambda^0)=\vz_\lambda}f(\vz_\lambda^0),
        \\
        (\Sc_\lambda(t)f)(\vz) = (\Sc_\lambda(t)f)(\vz_\lambda,\vz_{\lambda^c})=\sum_{\vz_\lambda^0:\Hc_\lambda(t)(\vz_\lambda^0)=\vz_\lambda}f(\vz_\lambda^0,\,\vz_{\lambda^c})
    \end{gathered}
    \end{equation}
    Finally, when $\lambda=[s]$, we will denote $\vz_\lambda=\vz_s$, and similarly $\Hc_\lambda=\Hc_s$ and $\Sc_\lambda=\Sc_s$.
    \item \emph{Composition of partial dynamics.} Given disjoint subsets $(\lambda_1, \cdots, \lambda_k)$, the \textbf{composition} of the partial dynamics of these subsets is the  dynamics specified by the flow map $\prod_{j=1}^k\Hc_{\lambda_j}(t)$.
\end{enumerate}
\end{definition}
\begin{remark}[]
Note that in general $\vz_\lambda(t)\neq (\vz_N(t))_\lambda$, where the left hand side is defined by the partial dynamics, and the right hand side is defined by selecting the particles in $\lambda$ from the full dynamics. This is due to the possible collisions between particles in $\lambda$ and $\lambda^c$. 

Note also that, in the composition of partial dynamics of $(\lambda_1, \ldots, \lambda_k)$, collisions occur \emph{only} between two particles in the same subset $\lambda_j$; any pre-collisional configuration between two particles in different subsets $\lambda_j$ are always ignored.
\end{remark}

\begin{definition}[Interaction of partial dynamics]\label{def.indicator_cluster} Let $\Mb$ and $\Mb'$ be \emph{C-molecules} with only one layer $\ell$, and particle sets $p(\Mb)$ and $p(\Mb')$ respectively. Let $\vz^{\Lambda,\Gamma}(t)$ be the truncated dynamics on $[(\ell-1)\tau, \ell\tau]$ and $\vz^0 = \vz^{\Lambda,\Gamma}|_{t = (\ell-1)\tau}$ be the initial configuration (with particle set being $p(\Mb)$ or $p(\Mb')$ or $p(\Mb)\cup p(\Mb')$).
\begin{enumerate}
    \item\label{it.indicator_cluster_1} \emph{The indicator functions $\mathbbm{1}^{\Lambda,\Gamma}_{\Mb}$.} We define the \textbf{indicator function $\mathbbm{1}^{\Lambda,\Gamma}_{\Mb}(\vz^0)$} to be $1$ if and only if the trajectory $\vz_{p(\Mb)}^{\Lambda,\Gamma}(t)$ topologically reduces to $\Mb$, and $0$ otherwise. Note that this is different\footnote{Here $\mathbbm{1}_\Mb$ only requires all collisions in $\Mb$ happen as expected, while $\mathbbm{1}^{\Lambda,\Gamma}_{\Mb}$ additionally prohibits other pre-collisional configurations not in $\Mb$.} from $\mathbbm{1}_\Mb$ in Definition \ref{def.associated_op_nonlocal} \eqref{it.associated_ind}, as the underlying dynamics is $\vz^{\Lambda,\Gamma}(t)$ instead of $\vz_\Mb(t)$.

    \item\label{it.indicator_cluster_2} \emph{The indicator functions $\mathbbm{1}_{\Mb\not\sim \Mb'}$.} We write $\Mb\not\sim \Mb'$ if and only if in the trajectory $\vz_{p(\Mb)\cup p(\Mb')}^{\Lambda,\Gamma}(t)$, there is no collisions between any atom in $p(\Mb)$ and any atom in $p(\Mb')$. Define the \textbf{indicator function $\mathbbm{1}_{\Mb \not\sim \Mb'}(\vz^0)$} to be $1$ if $\Mb \not\sim \Mb'$, and $0$ otherwise.
    \item\label{it.indicator_cluster_3} \emph{The indicator functions $\mathbbm{1}_{(\Mb,e)\not\sim(\Mb',e')}$.} Given two molecules $(\Mb, \Mb')$ and two edges $e\in \Mb$ and $e'\in\Mb'$, we write $(\Mb,e)\not\sim(\Mb',e')$ if any of the following two situations occur.
    \begin{enumerate}
        \item\label{it.indicator_cluster_3_a} \emph{No overlap:} under the composition dynamics $\Hc_{p(\Mb)}^{\Lambda, \Gamma}\circ \Hc_{p(\Mb')}^{\Lambda, \Gamma}$, the linear segments in the trajectory corresponding to $e$ and $e'$ never enter a pre-collisional configuration at any time $t_*$.
        \item\label{it.indicator_cluster_3_b} \emph{Truncated overlap:} $e$ and $e'$ do enter a pre-collisional configuration under the dynamics $\Hc_{p(\Mb)}^{\Lambda, \Gamma}\circ \Hc_{p(\Mb')}^{\Lambda, \Gamma}$ at time $t_*$, but $\mathtt{Pro}(\boldsymbol{\gamma},\pb,\pb')=\mathtt{False}$. Here $\pb,\pb'$ are the particles containing edges $e$ and $e'$ respectively, $\boldsymbol{\gamma}$ is the trajectory of the $\Hc_{p(\Mb)}^{\Lambda, \Gamma}\circ \Hc_{p(\Mb')}^{\Lambda, \Gamma}$ dynamics up to time $t_*$, and $\mathtt{Pro}$ is the condition in the $(\Lambda,\Gamma)$-truncated dynamics (Definition \ref{def.T_dynamics}) with particle set $p(\Mb)\cup p(\Mb')$.
    \end{enumerate}
    We define the \textbf{indicator function $\mathbbm{1}_{(\Mb,e)\not\sim(\Mb',e')}(\vz^0)$} to be $1$ if $(\Mb,e)\not\sim(\Mb',e')$, and $0$ otherwise. Note that these objects are well-defined only when the trajectory of particles in $p(\Mb)$ and $p(\Mb')$ topologically reduce to $\Mb$ and $\Mb'$.
\end{enumerate}
\end{definition}
\begin{remark}\label{rem.trunc_overlap} We explain the idea behind Definition \ref{def.indicator_cluster} (\ref{it.indicator_cluster_3}). The statement $(\Mb,e)\not\sim(\Mb',e')$ here means that, either the edges $e$ and $e'$ do not overlap, or their overlap \emph{does not cause a collision} in the $\Hc_{p(\Mb)\cup p(\Mb')}^{\Lambda,\Gamma}$ dynamics due to the $(\Lambda,\Gamma)$-truncation. We have to define it in this way (instead of just saying $e$ and $e'$ do not overlap), in order to ensure that the equality (\ref{eq.dynamics_factorization_1_M_notsim_M}) is true. Note that $(\Mb,e)\not\sim(\Mb',e')$ is \emph{less restrictive} than saying $e$ and $e'$ do not overlap, therefore the opposite statement $(\Mb,e)\sim(\Mb',e')$ is \emph{more restrictive} than saying $e$ and $e'$ do overlap. This implies an \emph{upper bound} for the corresponding indicator function $\mathbbm{1}_{(\Mb,e)\sim(\Mb',e')}$, which is all we need in all subsequent proofs. We emphasize that below \emph{it is always safe} to treat $(\Mb,e)\not\sim(\Mb',e')$ in the same way as the natural non-overlapping condition.
\end{remark}

Now we state the main result of this subsection, Lemma \ref{lem.dynamics_factorization}. The central formula here is (\ref{eq.dynamics_factorization}), which factorizes the transport operator $\Sc_N^{\Lambda,\Gamma}$ into $\Sc_{\Mb_j}(t)$. The product of indicator functions $\mathbbm{1}_{(\Mb_j,e)\not\sim(\Mb_{j'},e')}$ will be treated in Section \ref{sec.penrose} using the Penrose argument.

\begin{lemma}[]\label{lem.dynamics_factorization} We have the following decompositions.
\begin{enumerate}
    \item\label{it.dynamics_factorization_1} We have
    \begin{equation}\label{eq.dynamics_factorization_1}
        1 = \sum_{\Mb} \mathbbm{1}_\Mb^{\Lambda,\Gamma},
    \end{equation} where the summation is taken over all labeled C-molecules $\Mb$ (of single layer) with $p(\Mb)=[N]$. In this sum we may also assume that each component $\Mb_j$ of $\Mb$ satisfies $|\Mb_j|_p\leq\Lambda$ and $\rho(\Mb_j)\leq\Gamma$.
    \item\label{it.dynamics_factorization_2} If $\{\Mb_1, \cdots, \Mb_k\}$ are the connected components of $\Mb$, and $f$ is supported on $\mathbbm{1}_\Mb^{\Lambda,\Gamma}$, then we have
    \begin{align}\label{eq.dynamics_factorization_1_M}
        \mathbbm{1}_\Mb^{\Lambda,\Gamma} &= \bigg(\prod_{j= 1}^{k} \mathbbm{1}_{\Mb_j}^{\Lambda,\Gamma} \bigg) \bigg( \prod_{1\le j < j'\le k} \mathbbm{1}_{\Mb_j \not\sim \Mb_{j'}} \bigg),\\
\label{eq.dynamics_factorization_S}
        \mathcal{S}_N^{\Lambda,\Gamma}(t) f &= \prod_{j= 1}^{k} \mathcal{S}_{\Mb_j}(t) f.
    \end{align}
    \item\label{it.dynamics_factorization_3} In the set $(\supp\mathbbm{1}_{\Mb_j}^{\Lambda,\Gamma})\cap (\supp\mathbbm{1}_{\Mb_{j'}}^{\Lambda,\Gamma})$, we have
    \begin{equation}\label{eq.dynamics_factorization_1_M_notsim_M}
        \mathbbm{1}_{\Mb_j \not\sim \Mb_{j'}} = \prod_{e\in \Mb_j,\, e'\in \Mb_{j'}} \mathbbm{1}_{(\Mb_j,e) \not\sim(\Mb_{j'},e')}
    \end{equation}
    \item\label{it.dynamics_factorization_4} By putting together (\ref{it.dynamics_factorization_1})--(\ref{it.dynamics_factorization_3}), we get
    \begin{equation}\label{eq.dynamics_factorization}
        \mathcal{S}_N^{\Lambda,\Gamma}(t)=\sum_{k \leq N} \sum_{\{\Mb_1, \cdots, \Mb_k\}}
        \bigg(\prod_{j= 1}^{k} \mathcal{S}_{\Mb_j}(t) \circ\mathbbm{1}_{\Mb_j}^{\Lambda,\Gamma} \bigg) \circ\bigg( \prod_{1\le j < j'\le k} \prod_{e\in \Mb_j,\, e'\in \Mb_{j'}} \mathbbm{1}_{(\Mb_j,e) \not\sim(\Mb_{j'},e')} \bigg).
    \end{equation} Here the summation is taken over $\{\Mb_1,\cdots,\Mb_k\}$ which are connected components of a labeled C-molecule $\Mb$ as in (\ref{it.dynamics_factorization_1}) (equivalently this can also be taken over $\Mb$). Note that each $\Mb_j$ is a labeled C-molecule itself, and they satisfy that $p(\Mb_1)\cup \cdots\cup p(\Mb_k) = [N]$ and $|\Mb_j|_p\le \Lambda$, $|\chi_{\Mb_j}|\le \Gamma$.
\end{enumerate}    
\end{lemma}
Before proving Lemma \ref{lem.dynamics_factorization}, we first state and prove the following Lemma \ref{lem.clusters_property}, which contains the core parts of Lemma \ref{lem.dynamics_factorization}, i.e. (\ref{eq.dynamics_factorization_1_M})--(\ref{eq.dynamics_factorization_S}).
\begin{lemma}\label{lem.clusters_property}
Let $\Mb$ be a molecule with connected components $\{\Mb_1, \cdots, \Mb_k\}$ and $p(\Mb) = [N]$. Then, up to Lebesgue zero sets, the trajectory $\{\vz^{\Lambda,\Gamma}(t)\}_{t\in [0, \tau]}$ with initial data $\vz^0 = \vz^{\Lambda,\Gamma}(0)$ topologically reduces to $\Mb$, if and only if the following conditions hold:
\begin{enumerate}
    \item\label{it.clusters_property_1} For $1\le j\le k$, the partial dynamics $\vz^{\Lambda,\Gamma}_{p(\Mb_j)}(t)$ topologically reduces to $\Mb_j$, and
    \item\label{it.clusters_property_2} For $1\le j < j'\le k$, $\Mb_j\not\sim \Mb_{j'}$ in the partial dynamics $\vz^{\Lambda,\Gamma}_{p(\Mb_j)\cup p(\Mb_{j'})}(t)$.
\end{enumerate}
Let $\Sf$ be the set of initial configurations $\vz^0$ defined by the above conditions. If $\vz^0\in \Sf$ and $f$ is supported in $\Sf$, then for $t\in[0,\tau]$ we have
\begin{equation}\label{eq.dynamics_factorize}
    \Hc_N^{\Lambda,\Gamma}(t)\vz_N^0=\bigg(\prod_{j=1}^k\Hc_{\Mb_j}(t)\bigg)\vz_N^0,\quad\textrm{and}\quad \mathcal{S}_N^{\Lambda,\Gamma}(t)f = \bigg(\prod_{j=1}^k \mathcal{S}_{\Mb_j}(t)\bigg)f.
\end{equation} 
\end{lemma}

\begin{proof} We divide the proof into 4 parts. Recall $\vz_\lambda^{\Lambda,\Gamma}(t)$ is the partial dynamics, and $(\vz^{\Lambda,\Gamma}(t))_\lambda$ is obtained from selecting the particles in $\lambda$ from the full dynamics $\vz^{\Lambda,\Gamma}(t)$; in general they are not equal.

\textbf{Proof part 1.} In this part we prove that, if no particle in $\lambda$ collides with particle in $\lambda^c$ in the full dynamics $\vz^{\Lambda,\Gamma}(t)$, then we have $(\vz^{\Lambda,\Gamma}(t))_\lambda = \vz_\lambda^{\Lambda,\Gamma}(t)$.

To prove this fact, define \[t_* = \sup \big\{t: (\vz^{\Lambda,\Gamma}(t))_\lambda = \vz_\lambda^{\Lambda,\Gamma}(t)\big\},\] it suffices to show that $t_* = \tau$. If $t_*<\tau$, by left continuity we have $(\vz^{\Lambda,\Gamma}(t))_\lambda = \vz_\lambda^{\Lambda,\Gamma}(t)$ for all $t\leq t_*$. Let the C-molecule topological reduction of $\vz^{\Lambda,\Gamma}(t)$ and $\vz_\lambda^{\Lambda,\Gamma}(t)$ (up to time $t_*$) be $\Mb$ and $\Mb_{\lambda}$, since no collision occurs between particles in $\lambda$ and $\lambda^c$, we know that $\Mb_\lambda$ is just a connected component of $\Mb$ consisting particle lines in $\lambda$. Consider the possible collision in either $(\vz^{\Lambda,\Gamma}(t))_\lambda$ or $\vz_\lambda^{\Lambda,\Gamma}(t)$ at time $t_*$. By Definition \ref{def.T_dynamics} of the $(\Lambda,\Gamma)$-dynamics, we know that, whether this collision is allowed is determined only by properties of the molecule $\Mb_\lambda$, in both cases. As such, the two dynamics must follow the same scattering rule at $t_*$. Moreover, up to a zero measure, no other collisions occur for $t>t_*$ sufficiently close to $t_*$. Therefore,  the two dynamics remain identical for such $t$, contradicting the maximality of $t_*$. This contradiction shows that we must have $t_*=\tau$, as desired.

\textbf{Proof part 2.} In this part, we show that \eqref{it.clusters_property_1} and \eqref{it.clusters_property_2} hold under the assumption that $\vz^{\Lambda,\Gamma}(t)$ topologically reduces to $\Mb$. In fact, this assumption implies that no collision happens between any particle in $p(\Mb_j)$ and any particle in $p(\Mb_{j'})$ in the full dynamics $\vz^{\Lambda,\Gamma}(t)$. By Part 1, we know that 
$$\vz_{p(\Mb_j)}^{\Lambda,\Gamma}(t) = (\vz^{\Lambda,\Gamma}(t))_{p(\Mb_j)}\quad\textrm{and}\quad\vz_{p(\Mb_j)\cup p(\Mb_{j'})}^{\Lambda,\Gamma}(t) = (\vz^{\Lambda,\Gamma}(t))_{p(\Mb_j)\cup p(\Mb_{j'})}.$$ 
Clearly the first equation above implies \eqref{it.clusters_property_1}. The second equation implies that there is also no collision between $p(\Mb_j)$ and $p(\Mb_{j'})$ in the partial dynamics $\vz_{p(\Mb_j)\cup p(\Mb_{j'})}^{\Lambda,\Gamma}(t)$, which is \eqref{it.clusters_property_2}.

\textbf{Proof part 3.} In this part, we show that \eqref{it.clusters_property_1} and \eqref{it.clusters_property_2} implies that $\vz^{\Lambda,\Gamma}(t)$ topologically reduces to $\Mb$. Let $\lambda_j = p(\Mb_j)$ and $\lambda_{j'} = p(\Mb_{j'})$. We first show that under \eqref{it.clusters_property_1} and \eqref{it.clusters_property_2}, there is no collision between any $\lambda_j$ and $\lambda_{j'}$ in $\vz^{\Lambda,\Gamma}(t)$. Assume the contrary, there must exist a first collision time $t_*$ between some $\lambda_j$ and $\lambda_{j'}$. Let the next collision time in $\vz^{\Lambda,\Gamma}(t)$ be $t_*+\delta$, then in the time interval $[0, t_*+\delta]$, there is no collision between particles in $\lambda_j\cup\lambda_{j'}$ and $(\lambda_j\cup\lambda_{j'})^c$. By Part 1, we know that $(\vz^{\Lambda,\Gamma}(t))_{\lambda_j\cup \lambda_{j'}} = \vz_{\lambda_j\cup \lambda_{j}}^{\Lambda,\Gamma}(t)$ for $t\in [0, t_*+\delta]$, which is impossible because there is a collision between $\lambda_j$ and $\lambda_{j'}$ in $(\vz^{\Lambda,\Gamma}(t))_{\lambda_j\cup \lambda_{j'}}$ while there is no collision between $\lambda_j$ and $\lambda_{j'}$ in $\vz_{\lambda_j\cup \lambda_{j'}}^{\Lambda,\Gamma}(t)$ by \eqref{it.clusters_property_2}.

Now assume $\vz^{\Lambda,\Gamma}(t)$ topologically reduces to $\Mb'$, we only need to show $\Mb'=\Mb$. Since there is no collision between any $\lambda_j$ and $\lambda_{j'}$ in $\vz^{\Lambda,\Gamma}(t)$, by Part 1 we know that $\vz_{p(\Mb_j)}^{\Lambda,\Gamma}(t) = (\vz^{\Lambda,\Gamma}(t))_{p(\Mb_j)}$. As $\vz_{p(\Mb_j)}^{\Lambda,\Gamma}(t)$ topologically reduces to $\Mb_j$ by \eqref{it.clusters_property_1}, we know that the components of $\Mb'$ are identical to $\Mb_j$, hence $\Mb' = \Mb$.

\textbf{Proof part 4.} Finally we prove \eqref{eq.dynamics_factorize}. If $\vz^{\Lambda,\Gamma}$ topologically reduces to $\Mb$ with components $\{\Mb_{1}, \cdots,\Mb_{k}\}$, there is no collision between any $p(\Mb_j)$ and $p(\Mb_{j'})$, so $(\vz^{\Lambda,\Gamma}(t))_{p(\Mb_j)} = \vz_{p(\Mb_{j'})}^{\Lambda,\Gamma}(t)$ by Part 1. Let $\lambda_j=p(\Mb_j)$, we may identify $\vz=\vz_N$ with $(\vz_{\lambda_1},\cdots, \vz_{\lambda_k})$, so we have
\begin{equation}\label{eq.proof_cluster_lem_1}
    \begin{split}
        \Hc_N^{\Lambda,\Gamma}(t)\vz_N^0 &= \vz^{\Lambda,\Gamma}(t) = \left((\vz^{\Lambda,\Gamma}(t))_{\lambda_1},\cdots, (\vz^{\Lambda,\Gamma}(t))_{\lambda_k}\right) = \left(\vz_{\lambda_1}^{\Lambda,\Gamma}(t),\cdots, \vz_{\lambda_k}^{\Lambda,\Gamma}(t)\right) 
        \\
        &= \left(\Hc_{\lambda_1}^{\Lambda,\Gamma}(t)\vz^0_{\lambda_1},\cdots, \Hc_{\lambda_k}^{\Lambda,\Gamma}(t)\vz^0_{\lambda_k}\right) = \bigg(\prod_{j=1}^k\Hc_{\lambda_j}^{\Lambda,\Gamma}(t)\bigg)\vz_N^0 = \bigg(\prod_{j=1}^k\Hc_{\Mb_j}(t)\bigg)\vz_N^0.
    \end{split}
\end{equation} 
In the last step, we have applied the fact that $\Hc_{\lambda_j}^{\Lambda,\Gamma}(t,\vz^0) = \Hc_{\Mb_j}(t,\vz^0)$ (where $\Hc_{\Mb_j}$ is the molecule prescribed dynamics in Definition \ref{def.molecule_truncated_dynamics}), which is because the trajectory $\vz_{\lambda_j}^{\Lambda,\Gamma}(t)$ topologically reduces to $\Mb_j$. Finally, the $\Sc$ equality in \eqref{eq.dynamics_factorize} follows from the $\Hc$ equality, in the same way as (\ref{eq.dynamics_comparison}). This proves \eqref{eq.dynamics_factorize}.
\end{proof}

\begin{proof}[Proof of Lemma \ref{lem.dynamics_factorization}] All the molecules in this proof are C-molecules. First, the decomposition \eqref{it.dynamics_factorization_1} is obvious, because every trajectory topologically reduces to a unique C-molecule (Remark \ref{rem.top_reduction_uniqueness}). The property that $|\Mb_j|_p\leq\Lambda$ and $\rho(\Mb_j)\leq\Gamma$ follows from the $(\Lambda,\Gamma)$-truncation. Second, the equalities \eqref{eq.dynamics_factorization_1_M}--\eqref{eq.dynamics_factorization_S} in \eqref{it.dynamics_factorization_2} follow directly from Lemma \ref{lem.clusters_property} in view of the various indicator functions defined in Definition \ref{def.indicator_cluster}. As \eqref{it.dynamics_factorization_4} follows from simply putting together \eqref{it.dynamics_factorization_1}--\eqref{it.dynamics_factorization_3}, now we only need to prove \eqref{it.dynamics_factorization_3}.

To prove \eqref{it.dynamics_factorization_3}, it suffices to show that 
\begin{equation}\label{eq.proof_dynamics_factorization_1}
    \Mb \not\sim \Mb'\Leftrightarrow (\Mb,e)\not\sim(\Mb',e')\textrm{ for \emph{all} }e\in \Mb\textrm{ and }e'\in\Mb'.    
\end{equation} In fact, suppose the $\not\sim$ is replaced by $\sim$ on either side of (\ref{eq.proof_dynamics_factorization_1}), then there exists a pre-collisional configuration between a particle $\pb\in\Mb$ and a particle $\pb'\in\Mb'$, either in the composition dynamics $\Hc_{p(\Mb)}^{\Lambda, \Gamma}\circ \Hc_{p(\Mb')}^{\Lambda, \Gamma}$ or in the full dynamics $\Hc_{p(\Mb)\cup p(\Mb')}^{\Lambda, \Gamma}$. Consider the \emph{first time} $t_*$ when such a pre-collisional configuration occurs, then the two trajectories \emph{must coincide before} $t_*$, call it $\boldsymbol{\gamma}$. Now, if $\Mb\sim\Mb'$ due to a collision in the full dynamics, then we must have $\mathtt{Pro}(\boldsymbol{\gamma},\pb,\pb')=\mathtt{True}$ by in Definition \ref{def.T_dynamics}, hence $(\Mb,e)\sim(\Mb',e')$ for some edges $e$ and $e'$ in the particle lines $\pb$ and $\pb'$ by Definition \ref{def.indicator_cluster} (\ref{it.indicator_cluster_3}). Conversely, if $(\Mb,e)\sim(\Mb',e')$, then we also have $\mathtt{Pro}(\boldsymbol{\gamma},\pb,\pb')=\mathtt{True}$ by Definition \ref{def.indicator_cluster} (\ref{it.indicator_cluster_3}), hence $\Mb\sim\Mb'$ by Definition \ref{def.T_dynamics}.
\end{proof}

\subsection{Penrose argument and molecule representation}\label{sec.penrose} In this subsection, we start from (\ref{eq.dynamics_factorization}) in Lemma \ref{lem.dynamics_factorization}, and treat the product $\prod\mathbbm{1}_{(\Mb_j,e)\not\sim(\Mb_{j'},e')}$ using the Penrose argument which introduces \textbf{the O-atoms}, with the goal of proving (\ref{eq.S_N_decomposition}) in Proposition \ref{prop.S_N_decomposition}. This is motivated by similar arguments in  \cite{BGSS20}, except that we start from the product $\prod\mathbbm{1}_{(\Mb_j,e)\not\sim(\Mb_{j'},e')}$, instead of $\prod\mathbbm{1}_{\Mb_j\not\sim\Mb_{j'}}$ in \cite{BGSS20}.

\begin{definition}[]\label{def.indicator_general} Let $\Lambda$ and $\Gamma$ be fixed. The molecules $\Mb$ and $\Mb'$ in this definition have only one layer $\ell$, but may contain O-atoms. Let $\vz^{\Lambda,\Gamma}(t)$ be the truncated dynamics on $[(\ell-1)\tau, \ell\tau]$ (with suitable particle set), and let $|\Mb|_{\mathrm{O}}$ be the number of O-atoms in $\Mb$. We define the following objects:
\begin{enumerate}
    \item\label{it.indicator_cluster_T} \emph{The indicator function $\mathbbm{1}^{\Lambda,\Gamma}_{\Mb}$.} Define the \textbf{indicator function} $\mathbbm{1}^{\Lambda,\Gamma}_{\Mb}$, which extends the definition in Definition \ref{def.indicator_cluster} to allow O-atoms in $\Mb$, by
    \begin{equation}\label{eq.1_notsim}
        \mathbbm{1}^{\Lambda,\Gamma}_{\Mb} = \prod_{j = 1}^{k}  \mathbbm{1}^{\Lambda,\Gamma}_{\Mb_j}\cdot \prod_{1\le j < j' \le k} \prod_{\substack{e\in \Mb_j,\, e'\in \Mb_{j'} \\ \textrm{joined by O-atom}}}  \mathbbm{1}_{(\Mb_j,e)\sim(\Mb_{j'},e')} \cdot \mathbbm{1}_\prec.
    \end{equation}
    Here $\{\Mb_1, \cdots, \Mb_k\}$ are clusters of $\Mb$ obtained by deleting all O-atoms (Definition \ref{def.cluster}), and $e$ and $e'$ are edges in $\Mb_j$ and $\Mb_{j'}$ respectively (obtained from deleting O-atoms) that are joined by an O-atom in $\Mb$. The indicator function $\mathbbm{1}_\prec$ indicates that if $\of\prec \of'$ are O-atoms of $\Mb$ that belong to the same ov-segment, then the overlap corresponding to $\of$ occurs earlier than that of $\of'$.
    \item \emph{The operator $(\Sc\circ\mathbbm{1})_\Mb^{\mathrm{Pen}}$.} Recall the operators $(\Sc\circ\mathbbm{1})_\Mb$ and $|(\Sc\circ\mathbbm{1})_\Mb|$ defined in \eqref{eq.associated_op} (with $(\ell-\underline{\ell}+1)\tau$ replaced by $\tau$ as $\Mb$ has only one layer):
    \begin{equation}\label{eq.associated_op_single_layer}
        (\Sc\circ\mathbbm{1})_\Mb = (-1)^{|\Mb|_{\mathrm{O}}}\cdot\Sc_\Mb(\tau)\circ\mathbbm{1}_\Mb,\quad|(\Sc\circ\mathbbm{1})_\Mb| = \Sc_\Mb(\tau)\circ\mathbbm{1}_\Mb.
    \end{equation} We now define the operator $(\Sc\circ\mathbbm{1})_\Mb^{\mathrm{Pen}}$ which is variant of $(\Sc\circ\mathbbm{1})_\Mb$. The exact definition is very complicated and can be found in \eqref{eq.proof_S_N_decomposition_25_1} (which is never used in the rest of this paper); in the proof we will only rely on its several natural properties, which we list as follows.
    \begin{enumerate}
        \item\label{it.S1_Mb_sim_1} $(\Sc\circ\mathbbm{1})_\Mb^{\mathrm{Pen}}$ is an operator depending only on $\vz_{p(\Mb)}$ and satisfies
        \begin{equation}\label{eq.S1_Mb_sim_upper}
            |(\Sc\circ\mathbbm{1})_\Mb^{\mathrm{Pen}} f|\le |(\Sc\circ\mathbbm{1})_\Mb| (|f|),
        \end{equation}
        \begin{equation}\label{eq.S1_Mb_sim_asym}
            \big|(\Sc\circ\mathbbm{1})_\Mb^{\mathrm{Pen}} f - (\Sc\circ\mathbbm{1})_\Mb f \big|\le \sum_{e, e'\in \Mb} \big|(\Sc\circ\mathbbm{1})_{\Mb_{e,e'}} \big| (|f|),
        \end{equation}
        where $\Mb_{e,e'}$ is the molecule obtained by creating an O-atom joining $e$ and $e'$ in $\Mb$. 
        \item\label{it.S1_Mb_sim_2} If $\Mb$ can be decomposed into disconnected parts $\Mb = \Mb^1\sqcup \cdots\sqcup \Mb^k$, then 
        \begin{equation}\label{eq.S1_Mb_sim_asym_factorization}
            (\Sc\circ\mathbbm{1})_\Mb^{\mathrm{Pen}} = (\Sc\circ\mathbbm{1})_{\Mb^1}^{\mathrm{Pen}}\circ(\Sc\circ\mathbbm{1})_{\Mb^2}^{\mathrm{Pen}}\circ\cdots\circ (\Sc\circ\mathbbm{1})_{\Mb^k}^{\mathrm{Pen}}.
        \end{equation}
        \item\label{it.S1_Mb_sim_3} If $\Mb'$ is a obtained from relabeling non-root particles in $\Mb$ while preserving\footnote{This means that if two particles $\pb$ and $\pb'$ are in the same component, and $\pb$ has larger label than $\pb'$ before relabeling (recall labels are positive integers), then the same is true after relabeling.} the order of particle labels within each component of $\Mb$, then we have 
        \begin{equation}\label{eq.S1_Mb_sim_asym_relabelling}
            (\Sc\circ\mathbbm{1})_{\Mb'}^{\mathrm{Pen}} f' = (\Sc\circ\mathbbm{1})_\Mb^{\mathrm{Pen}} f
        \end{equation} 
        where $f'$ is the function obtained by renaming the variables according to this relabeling.
    \end{enumerate} 
    \item \emph{The set of labeled molecules $\Fs_{\Lambda}$ and $\Fs_{\Lambda}^{\mathrm{err}}$.}\label{it.setF} Define $\Fs_{\Lambda}$ (resp. $\Ts_{\Lambda}$ if $s=1$) and $\Fs_{\Lambda}^{\mathrm{err}}$ to be the set of single layer, labeled molecules $\Mb$ (\emph{not} equivalence classes) satisfying the following conditions:
    \begin{enumerate}
        \item\label{it.setF1} $r(\Mb) = [s]$ and every component of $\Mb$ contains at least one particle line $r\in r(\Mb)$. Define $\Mb(r)$ to be the component containing $r$; note that $\Mb(r)$ may be the same for different $r$.
        \item\label{it.setF2} For $\Mb\in \Fs_\Lambda$, each component of $\Mb$ contains $<\Lambda$ clusters. For $\Mb\in \Fs_\Lambda^{\mathrm{err}}$, each component contains $\leq \Lambda$ clusters, and there exists a unique component that contains exactly $\Lambda$ clusters. In either case, each cluster $\Mb_j$ satisfies $|\Mb_j|_p\leq\Lambda$ and $\rho(\Mb_j)\leq\Gamma$.
        \item\label{it.setF3} The cluster graph (Definition \ref{def.cluster}) of $\Mb$ is a forest.
    \end{enumerate}
\end{enumerate}
\end{definition}
\begin{remark} The Penrose argument we will perform in the proof of Proposition \ref{prop.S_N_decomposition} below involves a lot of artificial choices, which makes the definition of $(\Sc\circ\mathbbm{1})_\Mb^{\mathrm{Pen}}$ also highly artificial. As already pointed out, the exact formula is not important and we only need the properties (\ref{it.S1_Mb_sim_1})--(\ref{it.S1_Mb_sim_3}) in Definition \ref{def.indicator_general} (in practice $(\Sc\circ\mathbbm{1})_\Mb^{\mathrm{Pen}}$ will be $(\Sc\circ\mathbbm{1})_\Mb$ times some artificial indicator function, so these properties are natural).

The sets of molecules $\Fs_\Lambda$ and $\Fs_\Lambda^{\mathrm{err}}$ in (\ref{it.setF})
\label{rem.penrose_indicator} come from certain truncations we perform in the Penrose argument. Basically, if the number of clusters in each component is $<\Lambda$ then we are in (\ref{it.setF1}) in Definition \ref{def.indicator_general}; if not, we choose a first component with $\geq\Lambda$ clusters, truncate it at exactly $\Lambda$ clusters and treat the remaining clusters trivially, which leads to (\ref{it.setF2}) in Definition \ref{def.indicator_general}.
\end{remark}
\begin{lemma}[]\label{lem.1_M_notsim_le_1_M}
    For any molecule $\Mb$ we have $\mathbbm{1}^{\Lambda,\Gamma}_\Mb\le \mathbbm{1}_\Mb$.
\end{lemma}
\begin{proof} First assume $\Mb$ is C-molecule (so $\mathbbm{1}^{\Lambda,\Gamma}_\Mb$ is as in Definition \ref{def.indicator_cluster}). Then $\mathbbm{1}_\Mb$ only requires all collisions in $\Mb$ happen as expected, while $\mathbbm{1}^{\Lambda,\Gamma}_{\Mb}$ additionally prohibits other pre-collisional configurations not in $\Mb$, so we have $\mathbbm{1}^{\Lambda,\Gamma}_\Mb\le \mathbbm{1}_\Mb$.

If $\Mb$ contains O-atoms, by \eqref{eq.1_notsim} and Definition \ref{def.indicator_cluster} \eqref{it.indicator_cluster_3} (see Remark \ref{rem.trunc_overlap}) and the above proof, we have
\begin{equation}
    \mathbbm{1}^{\Lambda,\Gamma}_{\Mb} \le \prod_{j = 1}^{k}  \mathbbm{1}_{\Mb_j}\cdot \prod_{1\le j < j' \le k} \prod_{\substack{e\in \Mb_j,\, e'\in \Mb_{j'} \\ \textrm{joined by O-atom}}}  \mathbbm{1}_{e\textrm{ overlaps with }e'} \cdot \mathbbm{1}_\prec,
\end{equation}
and the right hand side exactly equals $\mathbbm{1}_\Mb$ (the $\mathbbm{1}_{\prec}$ specifies the time order of O-atoms in $\Mb$ which is not included in other indicator functions). This completes the proof.
\end{proof}

Now we state the main conclusion of this section, namely Proposition \ref{prop.S_N_decomposition} below.
\begin{proposition}[]\label{prop.S_N_decomposition} Let $(\Lambda,\Gamma)$ and $N$ and $s$ be fixed. We have
\begin{equation}\label{eq.S_N_decomposition}
    \mathcal{S}_N^{\Lambda,\Gamma}(\tau) = \sum_{\substack{\Mb\in \Fs_{\Lambda}\\ p(\Mb)\subseteq[N],\, r(\Mb)=[s]}}(\Sc\circ\mathbbm{1})_{\Mb}^{\mathrm{Pen}}\circ\Sc_{[N]\backslash p(\Mb)}^{\Lambda,\Gamma} + \sum_{\substack{\Mb\in \Fs_{\Lambda}^{\mathrm{err}}\\ p(\Mb)\subseteq[N],\, r(\Mb)=[s]}} O_1\bigg(|(\Sc\circ\mathbbm{1})_{\Mb}|\circ\Sc_{[N]\backslash p(\Mb)}^{\Lambda,\Gamma}\bigg).
\end{equation}
where $|O_1(A)|\leq A$ (see Definition \ref{def.notation} (\ref{it.misc})).
\end{proposition}
\begin{proof} Start from \eqref{eq.dynamics_factorization}. We have
\begin{equation}\label{eq.proof_S_N_decomposition_0_2}
\begin{gathered}
    \mathcal{S}_N^{\Lambda,\Gamma}(\tau) =\sum_{\overline{\Mb}}\bigg(\prod_{j=1}^k \Sc_{\Mb_j}(\tau)\circ \mathbbm{1}_{\Mb_j}^{\Lambda,\Gamma} \bigg)\circ\bigg(\prod_{1\le j < j'\le k} \prod_{e\in \Mb_j,\, e'\in \Mb_{j'}} \mathbbm{1}_{(\Mb_j,e) \not\sim(\Mb_{j'},e')}\bigg),
    \\
    \textrm{where }\overline{\Mb}\textrm{ is C-molecule, and }\{\Mb_1, \cdots, \Mb_k\}\textrm{ are connected components of }\overline{\Mb}.
\end{gathered}
\end{equation} We divide the proof into 4 parts.

\textbf{Proof part 1: preparations.} Let the labeled C-molecule $\overline{\Mb}$, and the clusters $\Mb_j$, be fixed for now. Note that $p(\overline{\Mb})=[N]$, and also fix the root particles set to be $[s]$. A cluster $\Mb_j$ is called a root cluster if it contains a root particle.

The $\Mb$ on the right hand side of (\ref{eq.S_N_decomposition}) will be formed by selecting certain clusters $\Mb_j$ and joining them by certain O-atoms. In particular this $\Mb$ is uniquely determined by (i) the set $p(\Mb)$ and the O-atoms joining edges $e\in\Mb_j$ and $e'\in\Mb_{j'}$, and (ii) the parent-child orderings between O-atoms in the same ov-segment.

In the proof below, the ordering in (ii) plays no role in the main arguments below: once (i) is fixed, we can arbitrarily set the ordering $\prec$ between O-atoms in the same ov-segment, and each choice of ordering leads to a unique molecule $\Mb$ in (\ref{eq.S_N_decomposition}). Correspondingly, we decompose $1=\sum_{\prec}\mathbbm{1}_\prec$ in (\ref{eq.1_notsim}), which occurs in $(\Sc\circ\mathbbm{1})_{\Mb}^{\mathrm{Pen}}$ in (\ref{eq.S_N_decomposition}).

We also introduce the following notations, which will be used in the proof below:
\begin{itemize}
\item Let $E$ be the set of edge pairs $(e,e')$ such that $e\in\Mb_j$ and $e'\in\Mb_{j'}$ for some $j\neq j'$. Denote elements of $E$ by $\of$ (corresponding to possible overlaps), and denote the event $(\Mb_j,e) \sim(\Mb_{j'},e')$ by $A_\of$.
\item For any $G\subseteq E$, define $\overline{\Mb}_G$ to be the result of creating an O-atom at each $\of=(e,e')\in G$ (Definition \ref{def.create_delete}) starting from $\overline{\Mb}$, which is in one-to-one correspondence with $G$. Define also $\mathrm{Cl}_G$ to be the cluster graph (Definition \ref{def.cluster}) of $\overline{\Mb}_G$.
\item Let $\Fs$ and $\Fs^{\mathrm{err}}$ be collections of subsets $F\subseteq E$ satisfying the following properties:
\begin{itemize}
\item Each component of the graph $\overline{\Mb}_F$ that contains an O-atom must contain a root cluster.
\item For $F\in\Fs$, each component of $\overline{\Mb}_F$ contains $<\Lambda$ clusters. For $F\in\Fs^{\mathrm{err}}$, each component of $\overline{\Mb}_F$ contains $\leq\Lambda$ clusters, and there exists a unique component that contains exactly $\Lambda$ clusters.
\item The cluster graph $\mathrm{Cl}_F$ is a forest.
\end{itemize}
\item For any $F\in\Fs\cup\Fs^{\mathrm{err}}$, define $\Mb=\Mb_F\subseteq\overline{\Mb}_F$ which is formed by the components of $\overline{\Mb}_F$ that contain at least one root cluster, and $\Mb'=\Mb_F':=\overline{\Mb}_F\backslash \Mb$ (which is a C-atom with each cluster being a component, due to the above properties of $F$).
\end{itemize} 

\textbf{Proof part 2: Penrose argument.} To treat the product $\prod \mathbbm{1}_{(\Mb_j,e) \not\sim(\Mb_{j'},e')}$ in (\ref{eq.proof_S_N_decomposition_0_2}), we apply the Penrose argument in the spirit of \cite{BGSS20}. By writing $\mathbbm{1}_{(\Mb_j,e) \not\sim(\Mb_{j'},e')}=1-\mathbbm{1}_{(\Mb_j,e) \sim(\Mb_{j'},e')}$ and expanding, we get
\begin{equation}\label{eq.incl_excl}\prod \mathbbm{1}_{(\Mb_j,e) \not\sim(\Mb_{j'},e')}=\sum_{G\in 2^E}(-1)^{|G|}\mathbbm{1}_G,\quad \mathbbm{1}_G=\prod_{(e,e')\in G}\mathbbm{1}_{(\Mb_j,e) \sim(\Mb_{j'},e')}.\end{equation} Recall also that $(e,e')=\of$ and $\mathbbm{1}_{(\Mb_j,e) \sim(\Mb_{j'},e')}=\mathbbm{1}_{A_\of}$.

The main issue with (\ref{eq.incl_excl}) is that the number of terms $\mathbbm{1}_G$ is too large. We will suitably group these terms and exploit the cancellation coming from the signs $(-1)^{|G|}$, to reduce the sum $\sum_{G\in 2^E}$ in (\ref{eq.incl_excl}) to a much smaller sum $\sum_{F\in\Fs\cup\Fs^{\mathrm{err}}}$, by the following argument:
\begin{enumerate}
\item\label{it.penrose1} We construct a suitable mapping $\pi:2^E\to\Fs\cup\Fs^{\mathrm{err}}$.
\item\label{it.penrose2} We show that, for each $F\in\Fs\cup\Fs^{\mathrm{err}}$ that belongs to the range of $\pi$, one can define a \emph{forbidden set} $FB(F)$, such that for any $G\in 2^E$ we have
\begin{equation}\label{eq.pvfbset}\pi(G)=F\Leftrightarrow F\subseteq G\textrm{ and }FB(F)\cap G=\varnothing.\end{equation}
\item\label{it.penrose3} Define also the \emph{allowed set} $AL(F):=E\backslash (F\cup FB(F))$. Let $(\Mb,\Mb')=(\Mb_F,\Mb_F')$ be defined as above, we also show that
\begin{equation}\label{eq.alset}AL(F)=\{(e,e'):e,e'\in\Mb'\}\cup AL'(F);
\end{equation} if $F\in\Fs$, we also show that
\begin{equation}\label{eq.alset2} AL'(F)\subseteq\{(e,e'):e,e'\textrm{ in the same component of }\Mb\}.
\end{equation} Finally, for $F\in \Fs$, the set $\mathrm{AL}(F)$ is invariant if we relabel the non-root particles in $\Mb=\Mb_F$ while preserving the order of particle labels within each component of $\Mb$.
\end{enumerate}

The mapping $\pi$ will be constructed (with properties (\ref{it.penrose2})--(\ref{it.penrose3}) proved) in Part 4 of the proof below. In this and the next part, we use \emph{only these properties} to finish the proof of (\ref{eq.S_N_decomposition}). The point is that we rewrite the sum in (\ref{eq.incl_excl}) as 
\[\sum_{G\in 2^E}(\cdots)=\sum_{F\in \Fs}\sum_{G:\pi(G)=F}(\cdots)+\sum_{F\in \Fs^{\mathrm{err}}}\sum_{G:\pi(G)=F}(\cdots)\] (where we also restrict $F$ to the range of $\pi$) and apply (\ref{eq.pvfbset}) to factorize the inner sum in $G$. More precisely:
\begin{equation}\label{eq.incl_excl_2}
\begin{aligned}
\prod \mathbbm{1}_{(\Mb_j,e) \not\sim(\Mb_{j'},e')}=\sum_{G\in 2^E}(-1)^{|G|}\mathbbm{1}_G&=\sum_{F\in\Fs\cup \Fs^{\mathrm{err}}}\sum_{G:\pi(G)=F}\prod_{\of=(e,e')\in G}(-\mathbbm{1}_{A_\of})\\
&=\sum_{F\in\Fs\cup \Fs^{\mathrm{err}}}\sum_{G=F\cup H,\,H\in 2^{AL(F)}}\prod_{\of\in F\cup H}(-\mathbbm{1}_{A_\of})\\
&=\sum_{F\in\Fs\cup \Fs^{\mathrm{err}}}\prod_{\of\in F}(-\mathbbm{1}_{A_\of})\cdot\sum_{H\in 2^{AL(F)}}(-1)^{|H|}\prod_{\of\in H}\mathbbm{1}_{A_\of}\\
&=\sum_{F\in\Fs\cup \Fs^{\mathrm{err}}}\prod_{\of\in F}(-\mathbbm{1}_{A_\of})\cdot\prod_{\of\in AL(F)}(1-\mathbbm{1}_{A_\of})\\
&=\sum_{F\in\Fs\cup \Fs^{\mathrm{err}}}(-1)^{|F|}\prod_{\of\in F}\mathbbm{1}_{A_\of}\cdot\prod_{\of\in AL(F)}\mathbbm{1}_{A_\of^c}.
\end{aligned}
\end{equation} Expanding the abbreviations, we get from (\ref{eq.incl_excl_2}) that
\begin{equation}\label{eq.proof_S_N_decomposition_22_1}
    \prod_{1\le j < j'\le k} \prod_{e\in \Mb_j,\, e'\in \Mb_{j'}} \mathbbm{1}_{(\Mb_j,e)\not\sim(\Mb_{j'},e')} = \sum_{F\in \Fs\cup \Fs^{\mathrm{err}}} (-1)^{|F|} \prod_{(e, e')\in F} \mathbbm{1}_{(\Mb_j,e)\sim(\Mb_{j'},e')} \prod_{(e, e')\in AL(F)} \mathbbm{1}_{(\Mb_j,e)\not\sim(\Mb_{j'},e')}.
\end{equation}
Now, by inserting \eqref{eq.proof_S_N_decomposition_22_1} into \eqref{eq.proof_S_N_decomposition_0_2}, and decomposing $1 = \sum_{\prec}\mathbbm{1}_\prec$ as described in Part 1 (with $\prec$ corresponding to the ordering between O-atoms that constitutes part of the molecule), we get
\begin{align}
    \mathcal{S}_N^{\Lambda,\Gamma}(\tau) &= \sum_{\overline{\Mb}}\bigg(\prod_{j=1}^k \Sc_{\Mb_j} \circ\mathbbm{1}_{\Mb_j}^{\Lambda,\Gamma}\bigg)\circ\bigg( \prod_{1\le j < j'\le k} \prod_{e\in \Mb_j,\, e'\in \Mb_{j'}} \mathbbm{1}_{(\Mb_j,e)\not\sim(\Mb_{j'},e')} \bigg)
    \nonumber\\
    &=\sum_{(\overline{\Mb},F\in \Fs\cup \Fs^{\mathrm{err}},\prec)} (-1)^{|F|}\, \underbrace{\bigg(\prod_{j=1}^k \Sc_{\Mb_j} \circ\mathbbm{1}_{\Mb_j}^{\Lambda,\Gamma}\bigg)\circ\bigg( \prod_{(e, e')\in F} \mathbbm{1}_{(\Mb_j,e)\sim(\Mb_{j'},e')} \mathbbm{1}_\prec\bigg)}_{=\Sc_{\overline{\Mb}_F} \circ\mathbbm{1}_{\overline{\Mb}_F}^{\Lambda,\Gamma}\textrm{ by \eqref{eq.1_notsim}}} \circ\bigg(\prod_{(e, e')\in AL(F)} \mathbbm{1}_{(\Mb_j,e)\not\sim(\Mb_{j'},e')}\bigg)
    \nonumber\\
    \label{eq.proof_S_N_decomposition_22_2}&=\sum_{\overline{\Mb}_F\,(F\in \Fs\cup \Fs^{\mathrm{err}})} (-1)^{|F|}\, (\Sc_{\overline{\Mb}_F} \circ\mathbbm{1}_{\overline{\Mb}_F}^{\Lambda,\Gamma})\circ\bigg( \prod_{(e, e')\in AL(F)} \mathbbm{1}_{(\Mb_j,e)\not\sim(\Mb_{j'},e')}\bigg).
\end{align}Here in the last line, we have used that the triple $(\overline{\Mb},F,\prec)$ uniquely corresponds to the molecule $\overline{\Mb}_F$, and that $\Sc_{\overline{\Mb}_F} = \prod_{j=1}^k \Sc_{\Mb_j}$ (since different clusters $\Mb_j$ are connected by overlaps). 

Moreover, if we define $\Mb=\Mb_F$ and $\Mb'=\Mb_F'$ in the notion of Part 1, then it is easy to see that $\overline{\Mb}_F$ is in one-to-one correspondence with the pair $(\Mb,\Mb')$; in addition, we have $F\in\Fs$ (resp. $F\in\Fs^{\mathrm{err}}$) if and only if $\Mb\in \Fs_\Lambda$ (resp. $\Mb\in \Fs_\Lambda^{\mathrm{err}}$) due to Definition \ref{def.indicator_general} (\ref{it.setF}) and the definitions in Part 1. Note also that the O-atoms in $\Mb$ are in bijection with elements in $F$ (so $|F|=|\Mb|_{\mathrm{O}}$), we can then rewrite (\ref{eq.proof_S_N_decomposition_22_2}) as
\begin{equation}\label{eq.proof_S_N_decomposition_23_1}
    \mathcal{S}_N^{\Lambda,\Gamma}(\tau) = \sum_{\substack{\Mb\in \Fs_{\ell, \Lambda}\cup\Fs_{\ell, \Lambda}^{\mathrm{err}}\\ p(\Mb)\subseteq[N],\, r(\Mb)=[s]}} \sum_{\Mb' : p(\Mb') = [N]\backslash p(\Mb)} (-1)^{|\Mb|_{\mathrm{O}}} \cdot(\Sc_{\overline{\Mb}_F} \circ\mathbbm{1}_{\overline{\Mb}_F}^{\Lambda,\Gamma})\circ\bigg( \prod_{(e, e')\in \mathrm{AL}(F)} \mathbbm{1}_{(\Mb_j,e)\not\sim(\Mb_{j'},e')}\bigg).
\end{equation} Note also that $\Mb'$ is a C-molecule as shown in Part 1.

\textbf{Proof part 3: the proof of (\ref{eq.S_N_decomposition}).} Now, starting from (\ref{eq.proof_S_N_decomposition_23_1}), we further split the product over $AL(F)$. By (\ref{eq.alset}) we know
\begin{equation}\label{eq.proof_S_N_decomposition_24_2}
    \prod_{(e, e')\in AL(F)} \mathbbm{1}_{(\Mb_j,e)\not\sim(\Mb_{j'},e')} 
    =\prod_{(e, e')\in AL'(F)} \mathbbm{1}_{(\Mb_j,e)\not\sim(\Mb_{j'},e')}\cdot\prod_{e, e'\in \Mb'} \mathbbm{1}_{(\Mb_j,e)\not\sim(\Mb_{j'},e')}.
\end{equation} Note that we will use the product $\prod_{(e, e')} \mathbbm{1}_{(\Mb_j,e)\not\sim(\Mb_{j'},e')}$ in (\ref{eq.proof_S_N_decomposition_24_2}) only when $F\in \Fs$; when $F\in\Fs^{\mathrm{err}}$ we will trivially bound it by $O_1(1)$. Also by \eqref{eq.1_notsim}, we have
\begin{equation}\label{eq.proof_S_N_decomposition_25_3}
\begin{aligned}
    \Sc_{\overline{\Mb}_F} \circ\mathbbm{1}_{\overline{\Mb}_F}^{\Lambda,\Gamma} &= \bigg(\prod_{j=1}^k \Sc_{\Mb_j}\bigg)\circ\bigg(\prod_{j=1}^k \mathbbm{1}_{\Mb_j}^{\Lambda,\Gamma} \prod_{(e, e')\in F} \mathbbm{1}_{(\Mb_j,e)\sim(\Mb_{j'},e')} \cdot\mathbbm{1}_\prec\bigg)
    \\
    &=\underbrace{\bigg(\prod_{\Mb_j\subseteq\Mb} \Sc_{\Mb_j} \circ\mathbbm{1}_{\Mb_j}^{\Lambda,\Gamma}\bigg)\circ\bigg( \prod_{(e, e')\in F} \mathbbm{1}_{(\Mb_j,e)\sim(\Mb_{j'},e')} \mathbbm{1}_\prec\bigg)}_{= \Sc_\Mb \circ\mathbbm{1}_\Mb^{\Lambda,\Gamma}\textrm{ by \eqref{eq.1_notsim}}}\circ \bigg(\prod_{\Mb_j\subseteq\Mb'} \Sc_{\Mb_j} \circ\mathbbm{1}_{\Mb_j}^{\Lambda,\Gamma}\bigg)
    \\ 
    &= \Sc_\Mb \circ\mathbbm{1}_\Mb^{\Lambda,\Gamma} \circ \bigg(\prod_{\Mb_j\subseteq\Mb'} \Sc_{\Mb_j} \circ\mathbbm{1}_{\Mb_j}^{\Lambda,\Gamma}\bigg).
\end{aligned}
\end{equation}
Now we define the operator $(\Sc\circ\mathbbm{1})_\Mb^{\mathrm{Pen}}$ by 
\begin{equation}\label{eq.proof_S_N_decomposition_25_1}
    (\Sc\circ\mathbbm{1})_\Mb^{\mathrm{Pen}} = (-1)^{|\Mb|_{\mathrm{O}}}\cdot\Sc_\Mb \circ \mathbbm{1}_\Mb^{\Lambda,\Gamma} \circ\bigg(\prod_{(e, e')\in AL'(F)} \mathbbm{1}_{(\Mb_j,e)\not\sim(\Mb_{j'},e')}\bigg).
\end{equation} Then, by inserting (\ref{eq.proof_S_N_decomposition_24_2})--(\ref{eq.proof_S_N_decomposition_25_1}) into (\ref{eq.proof_S_N_decomposition_23_1}), we get that

\begin{align}
    \mathcal{S}_N^{\Lambda,\Gamma}(\tau) 
    &=\sum_{\substack{\Mb\in \Fs_{\ell, \Lambda}\cup\Fs_{\ell, \Lambda}^{\mathrm{err}}\\ p(\Mb)\subseteq[N],\, r(\Mb)=[s]}} \sum_{\Mb' : p(\Mb') = [N]\backslash p(\Mb)} (-1)^{|\Mb|_{\mathrm{O}}}\cdot (\Sc_{\overline{\Mb}_F} \circ\mathbbm{1}_{\overline{\Mb}_F}^{\Lambda,\Gamma})\circ\bigg( \prod_{(e, e')\in AL(F)} \mathbbm{1}_{(\Mb_j,e)\not\sim(\Mb_{j'},e')}\bigg)
        \nonumber\\
    &=\sum_{\substack{\Mb\in \Fs_{\Lambda}\\ p(\Mb)\subseteq[N],\, r(\Mb)=[s]}} \sum_{\Mb' : p(\Mb') = [N]\backslash p(\Mb)} \underbrace{(-1)^{|\Mb|_{\mathrm{O}}}\cdot\Sc_\Mb \circ\mathbbm{1}_\Mb^{\Lambda,\Gamma}\circ \bigg(\prod_{(e, e')\in AL'(F)} \mathbbm{1}_{(\Mb_j,e)\not\sim(\Mb_{j'},e')}\bigg)}_{=(\Sc\circ\mathbbm{1})_\Mb^{\mathrm{Pen}}, \textrm{ by \eqref{eq.proof_S_N_decomposition_25_1}}}\nonumber\\
    &\phantom{=}\circ\bigg(\prod_{\Mb_j\subseteq\Mb'} \Sc_{\Mb_j} \circ\mathbbm{1}_{\Mb_j}^{\Lambda,\Gamma}\bigg)\circ\bigg(\prod_{e, e'\in \Mb'} \mathbbm{1}_{(\Mb_j,e)\not\sim(\Mb_{j'},e')}\bigg)
    \nonumber\\
    &+\sum_{\substack{\Mb\in \Fs_{\Lambda}^{\mathrm{err}}\\ p(\Mb)\subseteq[N],\, r(\Mb)=[s]}} \sum_{\Mb' : p(\Mb') = [N]\backslash p(\Mb)} \Sc_\Mb \circ\mathbbm{1}_\Mb^{\Lambda,\Gamma}\circ O_1(1)\circ\bigg(\prod_{\Mb_j\subseteq\Mb'} \Sc_{\Mb_j} \circ\mathbbm{1}_{\Mb_j}^{\Lambda,\Gamma}\bigg)\circ\bigg(\prod_{e, e'\in \Mb'} \mathbbm{1}_{(\Mb_j,e)\not\sim(\Mb_{j'},e')}\bigg)\nonumber\\
    &=\sum_{\substack{\Mb\in \Fs_{\Lambda}\\ p(\Mb)\subseteq[N],\, r(\Mb)=[s]}}(\Sc\circ\mathbbm{1})_\Mb^{\mathrm{Pen}}\circ\bigg[\underbrace{\sum_{\Mb' : p(\Mb') = [N]\backslash p(\Mb)}\bigg(\prod_{\Mb_j\subseteq\Mb'} \Sc_{\Mb_j} \circ\mathbbm{1}_{\Mb_j}^{\Lambda,\Gamma}\bigg)\circ\bigg(\prod_{e, e'\in \Mb'} \mathbbm{1}_{(\Mb_j,e)\not\sim(\Mb_{j'},e')}\bigg)}_{=\Sc_{[N]\backslash p(\Mb)}^{\Lambda,\Gamma}\textrm{ by (\ref{eq.dynamics_factorization})}}\bigg]\nonumber\\
    &+\sum_{\substack{\Mb\in \Fs_{\Lambda}\\ p(\Mb)\subseteq[N],\, r(\Mb)=[s]}}\Sc_\Mb \circ\mathbbm{1}_\Mb^{\Lambda,\Gamma}\circ O_1(1)\circ\bigg[\underbrace{\sum_{\Mb' : p(\Mb') = [N]\backslash p(\Mb)}\bigg(\prod_{\Mb_j\subseteq\Mb'} \Sc_{\Mb_j} \circ\mathbbm{1}_{\Mb_j}^{\Lambda,\Gamma}\bigg)\circ\bigg(\prod_{e, e'\in \Mb'} \mathbbm{1}_{(\Mb_j,e)\not\sim(\Mb_{j'},e')}\bigg)}_{=\Sc_{[N]\backslash p(\Mb)}^{\Lambda,\Gamma}\textrm{ by (\ref{eq.dynamics_factorization})}}\bigg]\nonumber\\
    \label{eq.proof_S_N_decomposition_26_1}&=\sum_{\substack{\Mb\in \Fs_{\Lambda}\\ p(\Mb)\subseteq[N],\, r(\Mb)=[s]}}(\Sc\circ\mathbbm{1})_{\Mb}^{\mathrm{Pen}}\circ\Sc_{[N]\backslash p(\Mb)}^{\Lambda,\Gamma} + \sum_{\substack{\Mb\in \Fs_{\Lambda}^{\mathrm{err}}\\ p(\Mb)\subseteq[N],\, r(\Mb)=[s]}} O_1\bigg(|(\Sc\circ\mathbbm{1})_{\Mb}|\circ\Sc_{[N]\backslash p(\Mb)}^{\Lambda,\Gamma}\bigg),
\end{align} where the time parameter in all $\Sc$ operators are $\tau$. This proves (\ref{eq.S_N_decomposition}). 

It remains to show that the $(\Sc\circ\mathbbm{1})_{\Mb}^{\mathrm{Pen}}$ constructed in \eqref{eq.proof_S_N_decomposition_25_1} satisfies the properties in Definition \ref{def.indicator_general} \eqref{it.S1_Mb_sim_1}-\eqref{it.S1_Mb_sim_3}. In fact, by \eqref{eq.proof_S_N_decomposition_25_1}, we know that 
\begin{equation}\label{eq.proof_S_N_decomposition_4_1}
    (\Sc\circ\mathbbm{1})_\Mb^{\mathrm{Pen}} = (-1)^{|\Mb|_{\mathrm{O}}}\cdot\Sc_\Mb \circ\mathbbm{1}_\Mb^{\Lambda,\Gamma} \circ\bigg(\prod_{(e, e')\in AL'(F)} \mathbbm{1}_{(\Mb_j,e)\not\sim(\Mb_{j'},e')}\bigg).
\end{equation} To prove Definition \ref{def.indicator_general} \eqref{it.S1_Mb_sim_1}, we need to verify \eqref{eq.S1_Mb_sim_upper} and \eqref{eq.S1_Mb_sim_asym}. Here \eqref{eq.S1_Mb_sim_upper} is a simple corollary of the fact that $\mathbbm{1}_\Mb^{\Lambda,\Gamma}\le \mathbbm{1}_\Mb$ (by Lemma \ref{lem.1_M_notsim_le_1_M}) and $\mathbbm{1}_{(\Mb_j,e)\not\sim(\Mb_{j'},e')} \le 1$. Moreover \eqref{eq.S1_Mb_sim_asym} is a corollary of the fact that \[\bigg|\prod_{(e, e')\in A} \mathbbm{1}_{(\Mb_j,e)\not\sim(\Mb_{j'},e')} - 1\bigg|\le \sum_{(e, e')\in A}\mathbbm{1}_{(\Mb_j,e)\sim(\Mb_{j'},e')};\qquad A:=AL'(F).\]

Now we prove Definition \ref{def.indicator_general} \eqref{it.S1_Mb_sim_2}. Assume $\Mb$ is the disconnected union $\Mb = \Mb^1\sqcup \cdots\sqcup \Mb^k$, it is easy to check that $\Sc_\Mb \circ\mathbbm{1}_\Mb^{\Lambda,\Gamma} = (\Sc_{\Mb^1} \circ\mathbbm{1}_{\Mb^1}^{\Lambda,\Gamma})\circ\cdots\circ(\Sc_{\Mb^k} \circ\mathbbm{1}_{\Mb^k}^{\Lambda,\Gamma})$. For the product over $AL'(F)$, a similar decomposition also holds true:
$$
\prod_{(e, e')\in AL'(F)} \mathbbm{1}_{(\Mb_j,e)\not\sim(\Mb_{j'},e')} = \prod_{i=1}^k\prod_{\substack{(e,e')\in AL'(F)\\e,e'\in\Mb^i}}\mathbbm{1}_{(\Mb_j,e)\not\sim(\Mb_{j'},e')},
$$ where we have used the fact that if $(e,e')\in AL'(F)$ then they must be in the same component of $\Mb$ (hence the same $\Mb^i$), due to (\ref{eq.alset2}). The property in Definition \ref{def.indicator_general} \eqref{it.S1_Mb_sim_2} then follows. Finally Definition \ref{def.indicator_general} \eqref{it.S1_Mb_sim_3} follows from the property of $AL(F)$ under relabeling, see property (\ref{it.penrose3}) in Part 2.

This completes the proof of Proposition \ref{prop.S_N_decomposition}, pending the construction of the mapping $\pi$. We will construct this $\pi$ in Part 4 below.

\textbf{Proof part 4: construction of $\pi$.} In this part, we define the map $\pi$ and prove properties (\ref{it.penrose2})--(\ref{it.penrose3}) in Part 2. The idea is similar to Proposition 2.3.3 in \cite{BGSS20}, but the arguments are more complicated due to truncation.

Before defining $\pi$, we first make some preparations. Recall $E$ is the set of edge pairs $(e,e')$; given a labeling of particle lines, we can define a linear ordering on the set $E$ as follows: suppose $\of=(e,e')\in E$, where $e$ and $e'$ belong to particle lines $\pb$ and $\pb'$ respectively, and $\pb<\pb'$ (in terms of labels). Then the linear ordering is defined by first comparing $\pb$ (in terms of labels), then comparing $\pb'$, and then comparing the order in the same particle line (descendant $<$ ancestor). Clearly this linear ordering is invariant if we relabel the non-root particles, but keeping the order of particle labels within each component of $\Mb$.

Now let $G\in 2^E$ be fixed. In the definition below, we will consider the graph $\overline{\Mb}_G$ defined in Part 1 (by creating an O-atom at each $\of=(e,e')\in G$) and its cluster graph $\mathrm{Cl}_G$. For each root particle line $r\in[s]$, let $\Qb_r$ be the (root) cluster containing $r$, and $\overline{\Mb}_G(r)$ be the component of $\overline{\Mb}_G$ containing $r$ (and $\Qb_r$), and $\mathrm{Cl}_G(r)$ be the corresponding component of $\mathrm{Cl}_G$. For each cluster $\Qb\in\mathrm{Cl}_G(r)$, consider the graph theoretic distance between $\Qb$ and $\Qb_r$ (i.e. the length of shortest path connecting them), and define $\mathrm{Cl}_G^q(r)$ to be the set of $\Qb\in\mathrm{Cl}_G(r)$ with the above distance being $q\geq 0$.

\textbf{Algorithm defining $\pi$.} For each root particle $r\in[s]$, let $\overline{\Mb}_G(r)$ and $\mathrm{Cl}_G(r)$ be defined above. We distinguish two cases:

\emph{Case 1: if $|\mathrm{Cl}_G(r)|<\Lambda$ for each $r$.} We perform the following steps:
\begin{enumerate}
\item\label{it.penrose_alg_1} Start from the first root particle $r=1$, and work in the graph $\mathrm{Cl}_G(1)$ (note that this graph may contain multiple edges, but this does not affect the arguments below). 
\item\label{it.penrose_alg_2} Consider all the elements $\of=(e,e')\in G$ (i.e. O-atoms in $\overline{\Mb}_G$) that join a cluster in $\mathrm{Cl}_G^q(1)$ and a cluster in $\mathrm{Cl}_G^{q+1}(1)$ for some $q\geq 0$. We order them first in the increasing order in $q$, and then in the increasing order in $\of$ (following the linear order defined above).
\item\label{it.penrose_alg_3} Start with the set $S:=\{\Qb_1\}$. Each time, select the smallest $\of=(e,e')$ in (\ref{it.penrose_alg_2}) that still remains in $\overline{\Mb}_G(1)$. Let $\Kb$ be the cluster in $\mathrm{Cl}_G^{q+1}(1)$ that is joined by $\of$ to a cluster in $\mathrm{Cl}_G^{q}(1)$, we add $\{\Kb\}$ to the set $S$, and delete all other elements of $G$ (i.e. O-atoms) that join $\Kb$ to another cluster in $S$. Repeat until $S$ exhausts all clusters in $\mathrm{Cl}_G(1)$, then stop working in this component.
\item\label{it.penrose_alg_4} Then, select the next smallest root particle $r$ that does not belong to any previously operated component, and repeat (\ref{it.penrose_alg_1})--(\ref{it.penrose_alg_3}) until all root clusters have been accounted for. Finally, delete all O-atoms that do not belong to any $\overline{\Mb}_G(r)$. Define the result to be $F=\pi(G)$.
\end{enumerate}

\emph{Case 2: if $|\mathrm{Cl}_G(r)|\geq\Lambda$ for some $r$.} We perform the same steps as in \emph{Case 1}, except at the first time when $|\mathrm{Cl}_G(r)|\geq\Lambda$. At this time we perform (\ref{it.penrose_alg_1})--(\ref{it.penrose_alg_3}) as in \emph{Case 1}, but terminate early when $|S|=\Lambda$ (without adding the next atom to $S$). Then we delete all the $\of$ that still remains in $\overline{\Mb}_G(r)$ and do not join two clusters in $S$. After this, stop exploring (even if there are still root clusters) and delete all O-atoms that do not belong to any previously operated component. Define the result to be $F=\pi(G)$.

It is easy to verify the following properties of $\pi(G)=F$; in particular, they imply that $F\in\Fs\cup\Fs^{\mathrm{err}}$ by definitions in Part 1.
\begin{itemize}
\item The cluster graph $\mathrm{Cl}_F$ is a forest. This is because in our construction (\ref{it.penrose_alg_3}), each new cluster in $S$ has only one overlap (after deletion) joining it to earlier clusters in $S$.
\item If $|\mathrm{Cl}_G(r)|<\Lambda$, then $\mathrm{Cl}_F(r)=\mathrm{Cl}_G(r)$ and $\mathrm{Cl}_F^q(r)=\mathrm{Cl}_G^q(r)$. This is because in our construction (\ref{it.penrose_alg_3}), the set $S$ remains connected in the cluster graph, and each cluster in $S\cap \mathrm{Cl}_G^q(r)$ is always connected to $\Qb_r$ by a path of length $q$ even after the deletions.
\item In \emph{Case 2}, $\mathrm{Cl}_F$ has a unique component $\mathrm{Cl}_F(r_0)$ that contains $\Lambda$ clusters; each component $\mathrm{Cl}_F(r')$ that is operated before $\mathrm{Cl}_F(r_0)$ contains $<\Lambda$ clusters, and each component other than $\mathrm{Cl}_F(r_0)$ and $\mathrm{Cl}_F(r')$ contains only one cluster. This is clear from the construction in \emph{Case 2}.
\end{itemize}

\textbf{Definition of $FB(F)$.} Now we define the set $FB(F)$ as follows.

\emph{Case 1: when $|\mathrm{Cl}_F(r)|<\Lambda$ for each $r$.} We define $\of\in FB(F)$ if and only if
    \begin{enumerate}
        \item\label{it.forbidden_1_1} $\of$ joins two clusters in two different components $\mathrm{Cl}_F(r)$;
        \item\label{it.forbidden_1_2} $\of$ joins a cluster in $\mathrm{Cl}_F^q(r)$ and a cluster in $\mathrm{Cl}_F^{q'}(r)$ with $|q-q'|\geq 2$; or
        \item\label{it.forbidden_1_3} $\of\not\in F$ joins a cluster $\Qb \in \mathrm{Cl}_F^q(r)$ and a cluster $\Kb \in \mathrm{Cl}_F^{q+1}(r)$, and the unique $\of'\in F$ joining $\Kb$ to a cluster in $\mathrm{Cl}_F^q(r)$ satisfies $\of<\of'$ in the linear ordering; or
        \item\label{it.forbidden_1_4} $\of$ joins a cluster in $\mathrm{Cl}_F(r)$ and a cluster in $\mathrm{Cl}_F$ not in any $\mathrm{Cl}_F(r')$.
    \end{enumerate}
    
\emph{Case 2: when $|\mathrm{Cl}_F(r_0)|=\Lambda$.} We define the components $\mathrm{Cl}_F(r_0)$ and previously operated components $\mathrm{Cl}_F(r')$ as above, and let their union be $\mathrm{Cl}_F(\leq r_0)$. Consider also the last cluster $\Kb$ reaching $|S|=\Lambda$ in \emph{Case 2} of the algorithm; assume it belongs to $\mathrm{Cl}_{F}^{q_0}(r_0)$ and the associated smallest $\of$ is $\of_0$. Then, we define $\of\in FB(F)$ if and only if
    \begin{enumerate}[resume]
        \item $\of$ satisfies one of \label{it.forbidden_2_1}\eqref{it.forbidden_1_1}--\eqref{it.forbidden_1_3} holds within $\mathrm{Cl}_F(\leq r_0)$; or
        \item\label{it.forbidden_2_2} $\of$ joins one cluster $\Kb$ outside $\mathrm{Cl}_{\leq r_0}$ and another cluster $\Qb$, which either belongs to $\mathrm{Cl}_F(r')$ for some $r'$ as above, or belongs to $\mathrm{Cl}_F^q(r_0)$ for some $q<q_0-1$, or belongs to $\mathrm{Cl}_F^{q_0-1}(r_0)$ and $\of<\of_0$.
    \end{enumerate}

\textbf{Proof of (\ref{it.penrose2})--(\ref{it.penrose3}) in Part 2.} Now, with the mapping $\pi$ and the set $FB(F)$ defined as above, we can prove the properties (\ref{it.penrose2})--(\ref{it.penrose3}) in Part 2. 

The main point is to show (\ref{eq.pvfbset}). Recall the properties \eqref{it.forbidden_1_1}--\eqref{it.forbidden_2_2} in the definition of $FB(F)$ above. We first show that if $\pi(G)= F$ then $FB(F)\cap G = \varnothing$ (as $F\subseteq G$ is obvious). Assume $\of\in FB(F)\cap G$ for contradiction, and consider each of the cases \eqref{it.forbidden_1_1}--\eqref{it.forbidden_2_2}.

\eqref{it.forbidden_1_1} If $\of$ joins a clusters in $\mathrm{Cl}_F(r)$ and a cluster in $\mathrm{Cl}_F(r')\neq \mathrm{Cl}_F(r)$, then we have $\mathrm{Cl}_G(r)= \mathrm{Cl}_G(r')$ because of $\of$, but this is contradiction because $\mathrm{Cl}_F(r)=\mathrm{Cl}_G(r)$ and $\mathrm{Cl}_F(r')=\mathrm{Cl}_G(r')$.

\eqref{it.forbidden_1_2} If $\of$ joins a clusters in $\mathrm{Cl}_F^q(r)$ and a cluster in $\mathrm{Cl}_F^{q'}(r)$ with $q'\geq q+2$, then this latter cluster belongs to $\mathrm{Cl}_G^{q+1}(r)$ or lower because of $\of$, but this is contradiction because $\mathrm{Cl}_G^{q+1}(r)=\mathrm{Cl}_F^{q+1}(r)$ is disjoint with $\mathrm{Cl}_F^{q'}(r)$.

\eqref{it.forbidden_1_3} If $\of$ joins $\Qb \in \mathrm{Cl}_F^q(r)$ and $\Kb \in \mathrm{Cl}_F^{q+1}(r)$, and satisfies $\of < \of'$ where $\of'\in F$, then consider the time when we add $\{\Kb\}$ into $S$. By the construction in (\ref{it.penrose_alg_3}) of the algorithm, we would select $\of$ and delete $\of'$, so the result will not be $F$, contradiction.

\eqref{it.forbidden_1_4} If $\of$ joins one cluster in $\mathrm{Cl}_F(r)$ to a cluster $\Kb$ not in any $\mathrm{Cl}_F(r')$, then $\Kb\in \mathrm{Cl}_G(r)$ because of $\of$, but this implies $\Kb\in \mathrm{Cl}_F(r)$, contradiction.

\eqref{it.forbidden_2_1} This case is basically the same as \eqref{it.forbidden_1_1}--\eqref{it.forbidden_1_3} but restricted to $\mathrm{Cl}_F(r_0)$ and $\mathrm{Cl}_F(r')$.

\eqref{it.forbidden_2_2} Assume $\of$ joins $\Kb\not\in\mathrm{Cl}_F(\leq r_0)$ and $\Qb$ as specified. If $\Qb\in\mathrm{Cl}_F(r')$ for some $r'$ in the definition of $FB(F)$ above, then we have $\Kb\in\mathrm{Cl}_G(r')=\mathrm{Cl}_F(r')$ due to $\of$, which is a contradiction. If $\Qb\in\mathrm{Cl}_F^q(r_0)$ for $q<q_0-1$, then similarly $\Kb\in\mathrm{Cl}_G^{q+1}(r_0)\subseteq\mathrm{Cl}_F(r_0)$ (the subset relation is because $q+1<q_0$, so all clusters in $\mathrm{Cl}_G^{q+1}(r_0)$ must have been exhausted in the construction in (\ref{it.penrose3}) of the algorithm before reaching $|S|=\Lambda$), which is also contradiction. Finally if $\Qb\in\mathrm{Cl}_F^{q_0-1}(r_0)$ and $\of<\of_0$, then again the in construction we would select $\of$ at some point before selecting $\of_0$, so the result will not be $F$, contradiction.

We next show that if $ F\subseteq G$ and $FB(F)\cap G = \varnothing$, then we must have $\pi(G)= F$. This is basically done by reverting the proof in cases \eqref{it.forbidden_1_1}--\eqref{it.forbidden_2_2} above, which we carry out below.

Assume first that $|\mathrm{Cl}_F(r)|<\Lambda$ for each $r$. Then, using the absence of elements in \eqref{it.forbidden_1_1} and \eqref{it.forbidden_1_4}, we deduce that $\mathrm{Cl}_F(r)=\mathrm{Cl}_G(r)$ for each $r$; using the absence of elements in \eqref{it.forbidden_1_2}, we can inductively show that $\mathrm{Cl}_F^q(r)=\mathrm{Cl}_G^q(r)$ for each $r$ and $q$. Finally, using the assumptions for \eqref{it.forbidden_1_3}, we see that the choice of edges in the construction of $\pi(G)$ exactly leads to those edges already in $F$. 

The case when some $|F_{r_0}|=\Lambda$ is similar. The same arguments above implies that $\mathrm{Cl}_F(r')=\mathrm{Cl}_G(r')$ for those $r'$ in the definition of $FB(F)$ above, and $\mathrm{Cl}_F^q(r_0)=\mathrm{Cl}_G^q(r_0)$ for $q<q_0$, and the result up to this point coincides with $F$. Then, when we execute (\ref{it.penrose3}) of the algorithm and each time choosing a smallest remaining $\of$, the absence of elements in \eqref{it.forbidden_2_2} guarantees that we cannot add anything outside $\mathrm{Cl}_F(\leq r_0)$ into $S$. So the result coincides with $F$, as desired.

Finally, we prove the properties of $AL(F)$ in (\ref{it.penrose3}) in Part 2. In fact, the first two properties (\ref{eq.alset})--(\ref{eq.alset2}) follow immediately from the definition of $FB(F)$ (and thus $AL(F)$) above. The last property concerning relabeling is also obvious, because the set $AL(F)$ is constructed using the linear ordering for $\of$, which is invariant under the said relabelings.
\end{proof}

\subsection{Recurrence formula for $\widetilde{f}_s$}\label{sec.f_to_f} In this subsection we apply the decomposition of $\Sc_N^{\Lambda,\Gamma}$, i.e. (\ref{eq.S_N_decomposition}) in Proposition \ref{prop.S_N_decomposition}, to obtain a representation formula of $\widetilde{f}_s(\ell\tau)$ at time $\ell\tau$, in terms of single-layer molecules $\Mb$ and the same quantity $\widetilde{f}_s$ at time $(\ell-1)\tau$. See Proposition \ref{prop.molecule_representation}.

\begin{definition}[Sets $\Fc_{\Lambda_\ell}$ and $\Fc_{\Lambda_\ell}^{\mathrm{err}}$]\label{def.set_F_single_layer} Given $\Lambda_\ell$, define $\Fc_{\Lambda_\ell}$ (resp. $\Fc_{\Lambda_\ell}^{\mathrm{err}}$) to be the set of \textbf{equivalence class} $[\Mb]$ of labeled molecules $\Mb$ of single layer $\ell$, such that $\Mb\in\Fs_{\Lambda_{\ell}}$ (resp. $\Mb\in\Fs_{\Lambda_\ell}^{\mathrm{err}}$) as in Definition \ref{def.indicator_general} (\ref{it.setF}). Note that the $\Tc_{\Lambda_\ell}$ in Definition \ref{def.set_T_F} (\ref{it.set_T}) is a special case of this $\Fc_{\Lambda_\ell}$ with only one root particle (the $\Tc_{\Lambda_\ell}^{\mathrm{err}}$ is unrelated to $\Fc_{\Lambda_\ell}^{\mathrm{err}}$, though.)
\end{definition}

\begin{proposition}[Recurrence formula for $\widetilde{f}_s$]\label{prop.molecule_representation} Recall $\widetilde{f}_s(\ell\tau,\vz_s)$ defined in (\ref{eq.W_widetilde_l}) and (\ref{eq.truncated_correlation}). For $\Mb$ with one layer $\ell$, define 
\begin{equation}\label{eq.molecule_terms}
\begin{gathered}
    \widetilde{f}_{\Mb}(\ell\tau, \vz_{r(\Mb)}) = \varepsilon^{-(d-1)|\Mb|_{p\backslash r}} \int (\Sc\circ\mathbbm{1})^{\mathrm{Pen}}_\Mb\big( \widetilde{f}_{|\Mb|_p}((\ell-1)\tau, \vz_{p(\Mb)}')\big) d\vz_{(p\backslash r)(\Mb)},
    \\
    \widetilde{f}_{\Mb, \mathrm{upp}}(\ell\tau, \vz_{r(\Mb)}) = \varepsilon^{-(d-1)|\Mb|_{p\backslash r}} \int |(\Sc\circ\mathbbm{1})_\Mb| \big(|\widetilde{f}_{|\Mb|_p}((\ell-1)\tau, \vz_{p(\Mb)}')|\big) d\vz_{(p\backslash r)(\Mb)}.
\end{gathered}
\end{equation}
Then we have the following formula
\begin{equation}\label{eq.molecule_representation}
\begin{aligned}
    \widetilde{f}(\ell\tau, \vz_s) &= \sum_{[\Mb]\in \Fc_{\Lambda_\ell},\,r(\Mb) = [s]} \widetilde{f}_{\Mb}(\ell\tau, \vz_s) + \mathrm{Err}(\ell\tau, \vz_s),
    \\ 
    |\mathrm{Err}(\ell\tau, \vz_s)|& \le \sum_{[\Mb]\in \Fc_{\Lambda_\ell}^{\mathrm{err}},\,r(\Mb) = [s]} \widetilde{f}_{\Mb, \mathrm{upp}}(\ell\tau, \vz_s).
\end{aligned}
\end{equation}
\end{proposition}
\begin{remark}\label{rem.f_to_f} The expressions in (\ref{eq.molecule_terms}) are just like those in (\ref{eq.associated_integral_molecule})--(\ref{eq.associated_integral_molecule_abs}) in Definition \ref{def.associated_int}, with the transport operator $(\Sc\circ\mathbbm{1})_\Mb$ (and its variant) associated with the molecule $\Mb$, and integration in the non-root particles. Here the molecule $\Mb$ has single layer $\ell$, as we evolve from time $(\ell-1)\tau$ to time $\ell\tau$. The input function is the density $\widetilde{f}_{|\Mb|_p}$ at time $(\ell-1)\tau$, with $\vz_{p(\Mb)}'$ being the states of the particles at time $(\ell-1)\tau$. Note also that $\ell_1[\pb]-\underline{\ell}=0$ in (\ref{eq.associated_integral_molecule})--(\ref{eq.associated_integral_molecule_abs}) because $\Mb$ has only one layer $\ell$.
  
  We remark that the $\widetilde{f}_{\Mb}$ in (\ref{eq.molecule_terms}) is \emph{not} constant on each equivalence class $[\Mb]$ due to the $(\Sc\circ\mathbbm{1})^{\mathrm{Pen}}_\Mb$ factor (in fact this is the \emph{only} instance of non-constant cases throughout this paper). By Definition \ref{def.notation} \eqref{it.sum_over_equiv}, we should interpret the sum as
  \begin{equation}\label{eq.sum_equiv}\sum_{[\Mb]\in \Fc_{\Lambda_\ell},\,r(\Mb) = [s]}(\cdots):=\sum_{\substack{\Mb\in \Fs_{\Lambda_\ell}\\p(\Mb)=[|\Mb|_p],r(\Mb) = [s]}}\frac{1}{(|\Mb|_p-s)!}(\cdots).
  \end{equation} Here note that (i) we may always restrict $p(\Mb)=[|\Mb|_p]$ in each equivalence class $[\Mb]$, without affecting the expression $\widetilde{f}_{\Mb}$ or the average in $[\Mb]$ (thanks to Definition \ref{def.indicator_general} \eqref{it.S1_Mb_sim_3}), and (ii) once $p(\Mb)=[|\Mb|_p]$ is fixed, the size of the equivalence class equals $|[\Mb]|=(|\Mb|_p-s)!$ due to permuting the $|\Mb|_p-s$ non-root labels.
\end{remark}
\begin{proof} In the summation below we only consider the sum over $\Mb\in \Fs_{\Lambda_{\ell}}$ and $[\Mb]\in\Fc_{\Lambda_\ell}$, as the the remaining sum over $\Fs_{\Lambda_{\ell}}^{\mathrm{err}}$ and $\Fc_{\Lambda_{\ell}}^{\mathrm{err}}$ is similar (in fact easier as $\widetilde{f}_{\Mb,\mathrm{upp}}$ is invariant under equivalence relation).

By \eqref{eq.W_widetilde_l} and \eqref{eq.truncated_correlation}, we have
\begin{equation}\label{eq.proof_molecule_representation_1}
    \widetilde{f}(\ell\tau, \boldsymbol{z}_s) = \varepsilon^{s(d-1)}\sum_{n=0}^\infty \frac{1}{n!} \int_{\Rb^{2dn}} \mathcal{S}^{\Lambda_\ell, \Gamma}_{s+n}(\tau) \widetilde{W}_{s+n}((\ell-1)\tau,\boldsymbol{z}_{s+n}')\,\mathrm{d}\vz_{[s+1:s+n]}.
\end{equation} Here, for distinction, we used a different set of variables $\vz_{s+n}'$ for the states at $(\ell-1)\tau$. Inserting \eqref{eq.S_N_decomposition}, we get the decomposition
\begin{equation}\label{eq.proof_molecule_representation_1-}
\begin{aligned}
    \widetilde{f}(\ell\tau, \boldsymbol{z}_s) 
    &= \varepsilon^{s(d-1)}\sum_{n = 0}^{\infty} \sum_{\substack{\Mb\in \Fs_{ \Lambda_\ell}\\ p(\Mb)\subseteq[s+n],\, r(\Mb)=[s]}} \frac{1}{n!} \int_{\Rb^{2dn}} (\Sc\circ\mathbbm{1})_{\Mb}^{\mathrm{Pen}}\circ\Sc_{[s+n]\backslash p(\Mb)} \widetilde{W}_{s+n}((\ell-1)\tau,\boldsymbol{z}_{s+n}')\,\mathrm{d}\vz_{[s+1:s+n]}\\
     &+ (\mathrm{Err}),
    \end{aligned}
\end{equation} where $(\mathrm{Err})$ represents the sum over $\Fs_{\Lambda_\ell}^{\mathrm{err}}$. By Definition \ref{def.indicator_general} \eqref{it.S1_Mb_sim_3}, in the summation (\ref{eq.proof_molecule_representation_1-}), we may always restrict $p(\Mb)=[|\Mb|_p]$, with an additional factor $\binom{n}{|\Mb|_p-s}$ (this fixes the \emph{set of labels} of non-root particles, but not the label of each \emph{individual} particles). This leads to the expression
\begin{equation}\label{eq.proof_molecule_representation_1+}
\begin{aligned}
\widetilde{f}(\ell\tau, \boldsymbol{z}_s)&=\varepsilon^{s(d-1)} \sum_{n=0}^\infty
\sum_{\substack{\Mb\in \Fs_{ \Lambda_\ell}\\ p(\Mb) = [|\Mb|_p], \,r(\Mb) = [s]}} \frac{1}{(|\Mb|_p-s)!(n-|\Mb|_p+s)!}\\
&\times\int_{\Rb^{2dn}} (\Sc\circ\mathbbm{1})_{\Mb}^{\mathrm{Pen}}\circ\Sc_{[s+n]\backslash p(\Mb)}^{\Lambda_\ell,\Gamma} \widetilde{W}_{s+n}((\ell-1)\tau,\boldsymbol{z}_{s+n}')\,\mathrm{d}\vz_{[s+1:s+n]}  + (\mathrm{Err}).
\end{aligned}
\end{equation}

Upon switching the order of summation, we get
\begin{equation}\label{eq.proof_molecule_representation_2}
\begin{aligned}
    &\widetilde{f}(\ell\tau, \boldsymbol{z}_s) 
    = \varepsilon^{s(d-1)}\sum_{\substack{\Mb\in \Fs_{\Lambda_\ell}\\ p(\Mb) = [|\Mb|_p], \, r(\Mb) = [s]}} \frac{1}{(|\Mb|_p-s)!} \int (\Sc\circ\mathbbm{1})_{\Mb}^{\mathrm{Pen}} 
    \\ 
    &
    \bigg(\sum_{n = |\Mb|_p-s}^{\infty} \frac{1}{(n-|\Mb|_p+s)!} \int \Sc_{[|\Mb|_p+1:s+n]}^{\Lambda_\ell,\Gamma} \widetilde{W}_{s+n}((\ell-1)\tau,\boldsymbol{z}_{s+n}')\,\mathrm{d}\vz_{[|\Mb|_p+1:s+n]}\bigg) \mathrm{d}\vz_{(p\backslash r)(\Mb)} + (\mathrm{Err}).
\end{aligned}
\end{equation} Consider the parenthesis in the second line in (\ref{eq.proof_molecule_representation_2}). Using the $L^1$ invariance of $\Sc^{\Lambda,\Gamma}$ (Proposition \ref{prop.dynamics_E_T_property} (\ref{it.dynamics_E_T_property_3.5})) and \eqref{eq.s_par_cor}, we conclude that this parenthesis equals $\varepsilon^{-(d-1)|\Mb|_p}\widetilde{f}_{|\Mb|_p}((\ell-1)\tau, \vz_{p(\Mb)}')$. This implies that
\begin{equation}\label{eq.proof_molecule_representation_3}
    \widetilde{f}(\ell\tau, \boldsymbol{z}_s) 
    = \sum_{\substack{\Mb\in \Fs_{\Lambda_\ell}\\ p(\Mb) = [|\Mb|_p],\,r(\Mb) = [s]}} \frac{\varepsilon^{-(d-1)(|\Mb|_p-s)}}{(|\Mb|_p-s)!} \int (\Sc\circ\mathbbm{1})_{\Mb}^{\mathrm{Pen}} 
    \widetilde{f}_{|\Mb|_p}((\ell-1)\tau, \vz_{p(\Mb)}') \,\mathrm{d}\vz_{p(\Mb)} + (\mathrm{Err}).
\end{equation} Now, using (\ref{eq.sum_equiv}) and note that $s=|\Mb|_r$, we see that the main term in (\ref{eq.proof_molecule_representation_3}) exactly matches that in (\ref{eq.molecule_representation}). Since the error term $(\mathrm{Err})$ can be treated in the same way, this completes the proof of Proposition \ref{prop.molecule_representation}.
\end{proof}

\subsection{Recurrence formula for cumulants}\label{sec.E_to_E} In this subsection, starting from (\ref{eq.molecule_representation}) in Proposition \ref{prop.molecule_representation}, we shall insert the ansatz (\ref{eq.cumulant_expansion}) for $\widetilde{f}_s((\ell-1)\tau)$ and subsequently derive the same expansion for $\widetilde{f}_s(\ell\tau)$. This also includes a recurrence formula for both $f^\Ac$ and $E_H$ in (\ref{eq.cumulant_expansion}), see Proposition \ref{prop.cumulant_single}. This is the main step in the proof of Proposition \ref{prop.cumulant_formula}; with it, we simply proceed by induction in $\ell$, and use the initial cumulant expansion for $\ell=0$ (Proposition \ref{prop.initial_cumulant} in Section \ref{sec.initial_cumulant}), which is done in Section \ref{sec.proof_cumulant_formula_final}.

\begin{definition}[]\label{def.r_conn_H}
    Given a molecule $\Mb$ and $H\subseteq p(\Mb)$, we define $r(\Mb)\searrow H$, if each particle line $r\in r(\Mb)$ belongs to a component of $\Mb$, which contains either another particle line in $r(\Mb)$, or a particle line in $H$.
\end{definition}
\begin{remark} Using Definition \ref{def.connectedvia}, the condition $r(\Mb)\searrow H$ can also be expressed as ``each particle line in $r(\Mb)$ is connected to either another particle line in $r(\Mb)$ via $\Mb$, or connected to a particle line in $H$ via $\Mb$, or belongs to $H$", in the same manner as Definition \ref{def.set_T_F} (\ref{it.set_F_l_4}).

In reality, $r(\Mb)=r(\Mb_{\ell'})$ will be the root particle set of a certain layer $\Mb_{\ell'}$, and $H=r(\Mb_{\ell'-1})$ will be the root particle set of the previous layer. The symbol $r(\Mb)\searrow H$ means that $r(\Mb)$ is either connected ``horizontally" to itself, or ``vertically" to $H$ in the previous layer.
\end{remark}

\begin{proposition}\label{prop.cumulant_single} Fix $\ell$ and $(A_\ell,\Lambda_\ell)$ defined in \eqref{eq.defLambdaseq}. Assume for each set $Q_{\ell-1}$ with $s^+=|Q_{\ell-1}|\leq \Lambda_\ell^2A_\ell$, we have cumulant expansion at $t = (\ell - 1)\tau$,
\begin{equation}\label{eq.cumulant_l-1}
    \widetilde{f}_{s^+}((\ell-1)\tau,\vz_{Q_{\ell-1}})=\sum_{H_{\ell-1}\subseteq Q_{\ell-1}}(f^\Ac((\ell-1)\tau))^{\otimes(Q_{\ell-1}\backslash H_{\ell-1})}\cdot E_{H_{\ell-1}}((\ell-1)\tau)+\mathrm{Err}((\ell-1)\tau,\vz_{Q_{\ell-1}}),
\end{equation} 
where $E_{H_{\ell-1}}((\ell-1)\tau)$ is a function of $\vz_{H_{\ell-1}}$ (with $E_{\varnothing}=1$ by convention) and is symmetric in its variables, $f^\Ac((\ell-1)\tau)$ is a single variable function of $z = (x, v)$, and the tensor product $(f^\Ac((\ell-1)\tau))^{\otimes(Q_{\ell-1}\backslash H_{\ell-1})}$ is a function of $\vz_{Q_{\ell-1}\backslash H_{\ell-1}}$. 

Then for each $Q_\ell$ with $s=|Q_\ell|\leq A_\ell$, we have the cumulant expansion at $t = \ell\tau$:
\begin{equation}\label{eq.cumulant_l}
\widetilde{f}_s(\ell\tau,\vz_{Q_\ell})=\sum_{H_\ell\subseteq Q_\ell}(f^{\Ac}(\ell\tau))^{\otimes(Q_\ell\backslash H_\ell)}\cdot E_{H_\ell}(\ell\tau,\vz_{H_\ell})+\mathrm{Err}(\ell\tau,\vz_{Q_\ell}).
\end{equation} 
Here, if we define
\begin{equation}\label{eq.associated_int_notsim_single}
\begin{gathered}
    \Ic\Nc_{\Mb, H}^{\mathrm{Pen}}(\vz_{r(\Mb)}) \coloneqq \varepsilon^{-(d-1)|\Mb|_{p\backslash r}}\int_{\Rb^{2d|\Mb|_{p\backslash r}}} (\Sc\circ\mathbbm{1})_{\Mb}^{\mathrm{Pen}} \left[\left(f^\Ac((\ell-1)\tau)\right)^{\otimes (p(\Mb)\backslash H)}\cdot E_H((\ell-1)\tau,\vz_H')\right]\mathrm{d}\vz_{(p\backslash r)(\Mb)},
    \\
    |\Ic\Nc_{\Mb, H}|(\vz_{r(\Mb)}) \coloneqq \varepsilon^{-(d-1)|\Mb|_{p\backslash r}}\int_{\Rb^{2d|\Mb|_{p\backslash r}}} |(\Sc\circ\mathbbm{1})_{\Mb}|\left[\left|f^\Ac((\ell-1)\tau)\right|^{\otimes (p(\Mb)\backslash H)}\cdot |E_H((\ell-1)\tau,\vz_H')|\right]\mathrm{d}\vz_{(p\backslash r)(\Mb)},
\end{gathered}
\end{equation}
then $f^{\Ac}$, $E_{H_\ell}$ and $\mathrm{Err}$ in \eqref{eq.cumulant_l} are explicitly given by the following formulas:
\begin{enumerate}
    \item\label{it.cumulant_single_1} $(f^{\Ac}(\ell\tau))^{\otimes(Q_\ell\backslash H_\ell)}$ is a function of $\vz_{Q_\ell\backslash H_\ell}$ which is a tensor product of the single variable function $f^{\Ac}(\ell\tau)$ given by
    \begin{equation}\label{eq.fAterm_single}
        f^{\mathcal{A}}(\ell\tau,z) = \sum_{[\Mb]\in \mathcal{T}_{\Lambda_\ell}}\Ic\Nc_{\Mb, \varnothing}^{\mathrm{Pen}}(\vz_{r(\Mb)}).
    \end{equation}
    Note that the right hand side of (\ref{eq.fAterm_single}) is also a function of a single variable $z$ as $|r(\Mb)|=1$.
    \item\label{it.cumulant_single_2} The cumulant $E_{H_\ell}(\ell\tau,\vz_{H_\ell})$ satisfies
    \begin{equation}\label{eq.EHterm_single}
        |E_{H_\ell}(\ell\tau,\vz_{H_\ell})| \le \sum_{\substack{[\Mb_\ell]\in \mathcal{F}_{ \Lambda_\ell},\, r(\Mb_\ell) = H_\ell\\ H_{\ell-1}\subseteq p(\Mb_\ell),\,r(\Mb_\ell)\searrow H_{\ell-1}}}|\Ic\Nc_{\Mb_\ell, H_{\ell-1}}|(\vz_{r(\Mb_\ell)}).
    \end{equation}  
    \item\label{it.cumulant_single_3} The error term $\mathrm{Err}(\ell\tau,\vz_{Q_\ell})$ satisfies $|\mathrm{Err}| \le |\mathrm{Err}^1|+|\mathrm{Err}^2|$, where $\mathrm{Err}^1$ and $\mathrm{Err}^2$ satisfy that
    \begin{align}\label{eq.Err1term}
        |\mathrm{Err}^1(\ell\tau,\vz_{Q_\ell})| &\le \sum_{[\Mb_\ell]\in \mathcal{F}_{ \Lambda_\ell}^{\mathrm{err}},\, r(\Mb_\ell) = Q_\ell,\, H_{\ell-1}\subseteq p(\Mb_\ell)}|\Ic\Nc_{\Mb_\ell, H_{\ell-1}}|(\vz_{r(\Mb_\ell)}),\\
        \label{eq.Err2term}
        \|\mathrm{Err}^2(\ell\tau)\|_{L^1} &\le \varepsilon^{-2(d-1)\Lambda_\ell^2A_\ell}\cdot \|\mathrm{Err}((\ell-1)\tau)\|_{L^1}.
    \end{align}
\end{enumerate}
\end{proposition}
\begin{remark}\label{rem.single_layer} We make a few remarks regarding Proposition \ref{prop.cumulant_single}. First, the function $E_H(\ell \tau)$ depends only on $|H|$ and is symmetric in its variables, and the functions $\mathrm{Err}$, $ \mathrm{Err}^1$ and $\mathrm{Err}^2$ are symmetric in all their variables. These will be guaranteed for the rest of this paper.

Second, the key condition $r(\Mb_\ell)\searrow H_{\ell-1}$ in (\ref{eq.EHterm_single}) comes from the requirement that no root particle in $E_{H_\ell}$ should be completely independent from the rest (and thus be factorized out). This is equivalent to the condition in Definition \ref{def.set_T_F} (\ref{it.set_F_l_4}), and is crucial in the subsequent algorithms.
\end{remark}
\begin{proof} Start with \eqref{eq.molecule_terms}--\eqref{eq.molecule_representation}, where we replace the variable $\vz_s$ by $\vz_{Q_\ell}$. We have
\begin{equation}\label{eq.proof_cumulant_single_1}
    \widetilde{f}(\ell\tau, \vz_{Q_\ell}) = \textrm{Main}^0(\ell\tau, \vz_{Q_\ell}) + \mathrm{Err}^0(\ell\tau, \vz_{Q_\ell}),
\end{equation}
where 
\begin{equation}\label{eq.proof_cumulant_single_2}
\begin{aligned}
    \textrm{Main}^0(\ell\tau, \vz_{Q_\ell}) &= \sum_{[\Mb_\ell]\in \Fc_{\Lambda_\ell}, \, r(\Mb_\ell) = Q_\ell} \widetilde{f}_{\Mb_\ell}(\ell\tau, \vz_{Q_\ell}),
    \\
    \widetilde{f}_{\Mb_\ell}(\ell\tau, \vz_{Q_\ell}) &= \varepsilon^{-(d-1)|\Mb_\ell|_{p\backslash r}} \int (\Sc\circ\mathbbm{1})^{\mathrm{Pen}}_{\Mb_\ell} \widetilde{f}_{|\Mb_\ell|_p}((\ell-1)\tau, \vz_{p(\Mb_\ell)}') \,\mathrm{d}\vz_{(p\backslash r)(\Mb_\ell)},
\end{aligned}
\end{equation}
and
\begin{equation}\label{eq.proof_cumulant_single_3}
\begin{aligned}
    |\mathrm{Err}^0(\ell\tau, \vz_{Q_\ell})| &\le \sum_{[\Mb_\ell]\in \Fc_{\Lambda_\ell}^{\mathrm{err}}, \, r(\Mb_\ell) = Q_\ell} \widetilde{f}_{\Mb_\ell, \mathrm{upp}}(\ell\tau, \vz_{Q_\ell})
    \\
    \widetilde{f}_{\Mb_\ell, \mathrm{upp}}(\ell\tau, \vz_{Q_\ell}) &= \varepsilon^{-(d-1)|\Mb_\ell|_{p\backslash r}} \int |(\Sc\circ\mathbbm{1})_{\Mb_\ell}| |\widetilde{f}_{|\Mb_\ell|_p}((\ell-1)\tau, \vz_{p(\Mb_\ell)}')| \,\mathrm{d}\vz_{(p\backslash r)(\Mb_\ell)}.
\end{aligned}
\end{equation}

In the proof below, we will show that the $\textrm{Err}^0$ term in (\ref{eq.proof_cumulant_single_1}) contributes part of $\mathrm{Err}^1$ and part of $\mathrm{Err}^2$, which satisfy \eqref{it.cumulant_single_3}. The $\textrm{Main}^0$ term in (\ref{eq.proof_cumulant_single_1}) contributes part of $\mathrm{Err}^2$, as well as the main expression
\begin{equation}\label{eq.proof_cumulant_single_*}
    \sum_{H_\ell\subseteq Q_\ell}(f^{\Ac}(\ell\tau))^{\otimes(Q_\ell\backslash H_\ell)}\cdot E_{H_\ell}(\ell\tau,\vz_{H_\ell})
\end{equation} in (\ref{eq.cumulant_l}), with $f^\Ac(\ell\tau)$ and $E_H(\ell\tau)$ satisfying \eqref{it.cumulant_single_1}--\eqref{it.cumulant_single_2}. We divide the proof into 3 parts.

\textbf{Proof part 1.} We first construct the $\mathrm{Err}^1$ and $\mathrm{Err}^2$ terms and prove \eqref{it.cumulant_single_3}. Start with $\mathrm{Err}^0$ in (\ref{eq.proof_cumulant_single_1}). Note that the number of root clusters in $\Mb_{\ell}$ does not exceed $s=Q_\ell\leq A$. By Definition \ref{def.indicator_general} (\ref{it.setF}) and Definition \ref{def.set_F_single_layer}, we see that the number of particle lines in $\Mb_\ell$ does not exceed $\Lambda_\ell^2A_\ell$. Thus, we may insert \eqref{eq.cumulant_l-1} into the second line of \eqref{eq.proof_cumulant_single_3}, and use the definitions (\ref{eq.associated_int_notsim_single}), to get
\begin{equation}\label{eq.proof_cumulant_single_2_1}
\begin{aligned}
    \widetilde{f}_{\Mb_\ell, \mathrm{upp}}(\ell\tau, \vz_{Q_\ell}) &= \varepsilon^{-(d-1)|\Mb_\ell|_{p\backslash r}} \int |(\Sc\circ\mathbbm{1})_{\Mb_\ell}| |\widetilde{f}_{|\Mb_\ell|_p}((\ell-1)\tau, \vz_{p(\Mb_\ell)}')| \,\mathrm{d}\vz_{(p\backslash r)(\Mb_\ell)}
    \\
    &=\sum_{H_{\ell-1}\subseteq p(\Mb_\ell)} |\Ic\Nc_{\Mb_\ell, H_{\ell-1}}|(\vz_{Q_\ell})+\mathrm{Err}_{\Mb_\ell}^2(\ell\tau,\vz_{Q_\ell}),
\end{aligned}
\end{equation}
where $\mathrm{Err}_{\Mb_{\ell}}^{2}$ satisfies that
\begin{equation}\label{eq.proof_cumulant_single_1_1'}
    |\mathrm{Err}_{\Mb_{\ell}}^2(\ell\tau, \vz_{Q_\ell})|\leq\varepsilon^{-(d-1)|\Mb_\ell|_{p\backslash r}} \int |(\Sc\circ\mathbbm{1})_{\Mb_\ell}| |\mathrm{Err}((\ell-1)\tau, \vz_{p(\Mb_\ell)}')| \,\mathrm{d}\vz_{(p\backslash r)(\Mb_\ell)}.
\end{equation}

Inserting (\ref{eq.proof_cumulant_single_2_1})--(\ref{eq.proof_cumulant_single_1_1'}) into the first line of \eqref{eq.proof_cumulant_single_3}, we get
\begin{equation}\label{eq.proof_cumulant_single_2_2}
    |\textrm{Err}^0(\ell\tau, \vz_{Q_\ell})| \le \sum_{[\Mb_\ell]\in \Fc_{ \Lambda_\ell}^{\mathrm{err}}, \, r(\Mb_\ell) = Q_\ell,\, H_{\ell-1}\subseteq p(\Mb_\ell)} |\Ic\Nc_{\Mb_\ell, H_{\ell-1}}|(\vz_{Q_\ell}) + |\textrm{Err}^{2,2}(\ell\tau, \vz_{Q_\ell})|,
\end{equation}
where 
\begin{equation}\label{eq.proof_cumulant_single_2_3}
\begin{aligned}
    |\Ic\Nc_{\Mb_\ell, H_{\ell-1}}|(\vz_{Q_\ell}) &= \varepsilon^{-(d-1)|\Mb_\ell|_{p\backslash r}} \int |(\Sc\circ\mathbbm{1})_{\Mb_\ell}| \left(|f^\Ac((\ell-1)\tau)|^{\otimes(p(\Mb_\ell)\backslash H_{\ell-1})}\cdot |E_{H_{\ell-1}}((\ell-1)\tau)|\right),\\
    |\textrm{Err}^{2,2}(\ell\tau, \vz_{Q_\ell})| &\le \sum_{[\Mb_\ell]\in\Fc_{ \Lambda_\ell}^{\mathrm{err}}, \, r(\Mb_\ell) = Q_\ell}\varepsilon^{-(d-1)|\Mb_\ell|_{p\backslash r}} \int |(\Sc\circ\mathbbm{1})_{\Mb_\ell}| |\mathrm{Err}((\ell-1)\tau, \vz_{p(\Mb_\ell)})| \,\mathrm{d}\vz_{(p\backslash r)(\Mb_\ell)}.
    \end{aligned}
\end{equation}

Now we define the $\sum(\cdots)$ term on the right hand side of (\ref{eq.proof_cumulant_single_2_2}) to be $\mathrm{Err}^1$, then \eqref{eq.Err1term} follows immediately from this definition. As for the term $\mathrm{Err}^{2,2}$, we shall combine it with the term $\mathrm{Err}^{2,1}$ below (which comes from $\mathrm{Main}^0$ in (\ref{eq.proof_cumulant_single_1})) to form $\mathrm{Err}^2$, and prove it satisfies \eqref{eq.Err2term}, which would then complete the proof of \eqref{it.cumulant_single_3}.

Now consider $\mathrm{Main}^0$ in (\ref{eq.proof_cumulant_single_1}). By inserting \eqref{eq.cumulant_l-1} into the second line of \eqref{eq.proof_cumulant_single_2}, and using the definitions (\ref{eq.associated_int_notsim_single}), we get
\begin{equation}\label{eq.proof_cumulant_single_1_1}
\begin{aligned}
    \widetilde{f}_{\Mb_\ell}(\ell\tau, \vz_{Q_\ell}) &= \varepsilon^{-(d-1)|\Mb_\ell|_{p\backslash r}} \int (\Sc\circ\mathbbm{1})^{\mathrm{Pen}}_{\Mb_\ell} \widetilde{f}_{|\Mb_\ell|_p}((\ell-1)\tau, \vz_{p(\Mb_\ell)}') \,\mathrm{d}\vz_{(p\backslash r)(\Mb_\ell)}
    \\
   &=\sum_{H_{\ell-1}\subseteq p(\Mb_\ell)} \Ic\Nc_{\Mb_\ell, H_{\ell-1}}^{\mathrm{Pen}}(\vz_{Q_\ell})+\mathrm{Err}_{\Mb_{\ell}}^{2}(\ell\tau, \vz_{Q_\ell}),
\end{aligned}
\end{equation} where $\mathrm{Err}_{\Mb_{\ell}}^{2}$ satisfies the same bound (\ref{eq.proof_cumulant_single_1_1'}).
Inserting \eqref{eq.proof_cumulant_single_1_1} into the first line of \eqref{eq.proof_cumulant_single_2}, we get
\begin{equation}\label{eq.proof_cumulant_single_1_2}
    \textrm{Main}^0(\ell\tau, \vz_{Q_\ell}) = \sum_{[\Mb_\ell]\in \Fc_{ \Lambda_\ell}, \, r(\Mb_\ell) = Q_\ell,\, H_{\ell-1}\subseteq p(\Mb_\ell)} \Ic\Nc_{\Mb_\ell, H_{\ell-1}}^{\mathrm{Pen}}(\vz_{Q_\ell})+\mathrm{Err}^{2,1}(\ell\tau, \vz_{Q_\ell}),
\end{equation}
where
\begin{equation}\label{eq.proof_cumulant_single_1_2+}
    \Ic\Nc_{\Mb_\ell, H_{\ell-1}}^{\mathrm{Pen}}(\vz_{Q_\ell}) = \varepsilon^{-(d-1)|\Mb_\ell|_{p\backslash r}} \int (\Sc\circ\mathbbm{1})^{\mathrm{Pen}}_{\Mb_\ell} \left((f^\Ac((\ell-1)\tau))^{\otimes(p(\Mb_\ell)\backslash H_{\ell-1})}\cdot E_{H_{\ell-1}}((\ell-1)\tau)\right),
    \end{equation} and $\mathrm{Err}^{2,1}$ is the same as $\mathrm{Err}^{2,2}$ in (\ref{eq.proof_cumulant_single_2_3}) but with summation taken over $[\Mb_\ell]\in\Fc_{\Lambda_\ell}$. By combining this with $\mathrm{Err}^{2,2}$ we get the term $\mathrm{Err}^2$, which the same as in (\ref{eq.proof_cumulant_single_2_3}) but with summation taken over $[\Mb_\ell]\in\Fc_{\Lambda_\ell}\cup \Fc_{\Lambda_\ell}^{\mathrm{err}}$.

Next we prove \eqref{eq.Err2term} for $\mathrm{Err}^2$. Note that $|\Mb_\ell|_{p\backslash r}\leq \Lambda_\ell^2A_\ell$ as proved before. Using $L^1$ conservation, it is easy to see that the $L^1$ norm of each individual term in the summation on the right hand side of (\ref{eq.proof_cumulant_single_2_3}) is bounded by $\varepsilon^{-(d-1)\Lambda_\ell^2A_\ell}\cdot\|\mathrm{Err}((\ell-1)\tau)\|_{L^1}$. Moreover, consider all possible choices of labeled molecules $[\Mb_{\ell}]\in\Fc_{\Lambda_\ell}\cup\Fc_{\Lambda_\ell}^{\mathrm{err}}$ (equivalently $\Mb_\ell\in \Fs_{\Lambda_\ell}\cup\Fs_{\Lambda_\ell}^{\mathrm{err}}$). Since the total number of particles is $\leq\Lambda_\ell^2A_\ell\leq |\log\varepsilon|^{C^*}$, and each cluster has at most $\Gamma$ recollisions due to truncation, we may simply fix all the collisions and overlaps and their relative ordering, and the combinatorial factor caused by this is bounded by $(|\log\varepsilon|^{C^*})^{\Lambda_\ell^2A_\ell}$. This proves (\ref{eq.Err2term}) and finishes the proof of \eqref{it.cumulant_single_3}.

\textbf{Proof part 2.} It remains to construct the expression (\ref{eq.proof_cumulant_single_*}) and prove \eqref{it.cumulant_single_1}--\eqref{it.cumulant_single_2}. This will come from the $\sum(\cdots)$ term in (\ref{eq.proof_cumulant_single_1_2}). Note that this involves a molecule $\Mb_\ell$ and a set $H_{\ell-1}\subseteq p(\Mb_\ell)$. Consider the components of $\Mb_{\ell}$. We say a component is \emph{single}, if it contains a unique root particle $r\in r(\Mb_\ell)$ (we also say this $r$ is \emph{single}), and does not contain any particle in $H_{\ell-1}$.

Note that the single components are labeled by single root particles $r$, which we denote by $\Mb_\ell(r)$. Also define $\Mb_{\ell}'$ to be the union of all \emph{non-single} components of $\Mb_{\ell}$, and denote $r(\Mb_\ell')=H_\ell\subseteq Q_\ell$ to be the set of non-single particles. By definition of single, it is easy to see that $r(\Mb_{\ell}')\searrow H_{\ell-1}$. See {\color{blue}Figure \ref{fig.ehset}}.
\begin{figure}[h!]
\includegraphics[width=0.5\linewidth]{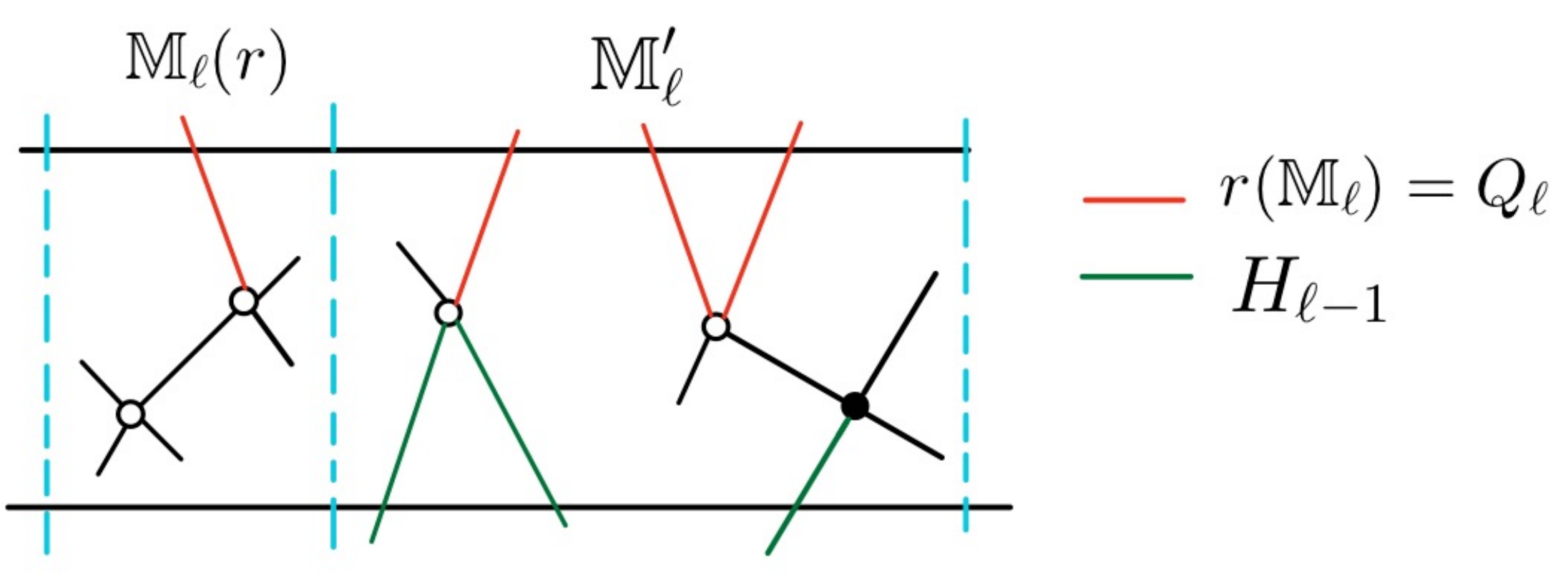}
\caption{The molecule $\Mb_{\ell}$ and its single component $\Mb_\ell(r)$, and non-single component $\Mb_{\ell}'$. Here the particle lines in $r(\Mb_\ell)=Q_\ell$ are marked as red, and those in $H_{\ell-1}$ marked as green. In particular, $H_\ell=r(\Mb_{\ell}')$ is formed by the 3 red particle lines in $\Mb_\ell'$.}
\label{fig.ehset}
\end{figure}

Note that, as \emph{unlabeled} molecules, $\Mb_\ell$ is in one-to-one correspondence with $(\Mb_\ell(r))_{r\not\in H_\ell}$ and $\Mb_{\ell}'$. Next we show that
\begin{equation}\label{eq.proof_cumulant_single_1_5}
    \Ic\Nc_{\Mb_\ell, H_{\ell-1}}^{\mathrm{Pen}}(\vz_{Q_\ell}) = \bigg(\prod_{r\notin H_\ell} \Ic\Nc_{\Mb_\ell(r), \varnothing}^{\mathrm{Pen}}(z_r) \bigg) \cdot \Ic\Nc_{\Mb_\ell', H_{\ell-1}}^{\mathrm{Pen}}(\vz_{H_\ell}).
\end{equation} To prove this, simply note that
\begin{equation}\label{eq.proof_cumulant_single_3_1}
    (\Sc\circ\mathbbm{1})_{\Mb_\ell}^{\mathrm{Pen}} = \prod_{r\notin H_\ell}(\Sc\circ\mathbbm{1})_{\Mb_\ell(r)}^{\mathrm{Pen}}\circ (\Sc\circ\mathbbm{1})_{\Mb_\ell'}^{\mathrm{Pen}}
\end{equation}
 by \eqref{eq.S1_Mb_sim_asym_factorization}. Inserting this into the the first line of \eqref{eq.associated_int_notsim_single}, we get
\begin{equation}\label{eq.proof_cumulant_single_3_2}
\begin{aligned}
    \Ic\Nc_{\Mb_\ell, H_{\ell-1}}^{\mathrm{Pen}}(\vz_{Q_\ell}) &= \varepsilon^{-(d-1)|\Mb_\ell|_{p\backslash r}} \int (\Sc\circ\mathbbm{1})^{\mathrm{Pen}}_{\Mb_\ell} \left((f^\Ac((\ell-1)\tau))^{\otimes(p(\Mb_\ell)\backslash H_{\ell-1})}\cdot E_{H_{\ell-1}}((\ell-1)\tau)\right)
    \\
    &= \varepsilon^{-(d-1)|\Mb_\ell|_{p\backslash r}} \int \prod_{r\notin H_\ell}\left((\Sc\circ\mathbbm{1})_{\Mb_\ell(r)}^{\mathrm{Pen}} (f^\Ac((\ell-1)\tau))^{\otimes(p(\Mb_\ell(r)))}\right) \times
    \\
    &\phantom{=\varepsilon^{-(d-1)|\Mb_\ell|_{p\backslash r}}} \int(\Sc\circ\mathbbm{1})_{\Mb_\ell'}^{\mathrm{Pen}} \left((f^\Ac((\ell-1)\tau))^{\otimes(p(\Mb_\ell')\backslash H_{\ell-1})}\cdot E_{H_{\ell-1}}((\ell-1)\tau)\right)
    \\
    &=\bigg(\prod_{r\notin H_\ell} \Ic\Nc_{\Mb_\ell(r), \varnothing}^{\mathrm{Pen}}(z_r) \bigg) \cdot \Ic\Nc_{\Mb_\ell', H_{\ell-1}}^{\mathrm{Pen}}(\vz_{H_\ell})
\end{aligned}
\end{equation}
which proves \eqref{eq.proof_cumulant_single_1_5}.

Now, note that the right hand side of \eqref{eq.proof_cumulant_single_1_5} has the structure \[\prod_{r\not\in H_\ell}\mathrm{Fun}(\Mb_\ell(r))\cdot \mathrm{Fun}(\Mb_\ell',H_{\ell-1}),\] where $\mathrm{Fun}(\cdot)$ is an arbitrary expression depending only on the object $(\cdot)$. We claim that using the convention in Definition \ref{def.notation} (\ref{it.sum_over_equiv}), we have
\begin{equation}\label{eq.proof_cumulant_single_1_6}
\begin{aligned}
    &\quad\ \ \sum_{[\Mb_\ell]\in \Fc_{\Lambda_\ell}, \, r(\Mb_\ell) = Q_\ell,\, H_{\ell-1}\subseteq p(\Mb_\ell)} \prod_{r\not\in H_\ell}\mathrm{Fun}(\Mb_\ell(r))\cdot \mathrm{Fun}(\Mb_\ell',H_{\ell-1})\\
    &= \sum_{H_\ell\subseteq Q_\ell} \bigg(\prod_{r\notin H_\ell} \sum_{[\Mb_\ell(r)]\in \Tc_{ \Lambda_\ell}}\mathrm{Fun}(\Mb_\ell(r))\bigg)\cdot\sum_{\substack{[\Mb_\ell']\in \Fc_{\Lambda_\ell}, \, r(\Mb_\ell') = H_\ell\\ H_{\ell-1}\subseteq p(\Mb_\ell'),\,r(\Mb_\ell')\searrow H_{\ell-1}}}\mathrm{Fun}(\Mb_\ell',H_{\ell-1}).
    \end{aligned}
\end{equation} Note that this equality is not trivial (even though $\Mb_\ell$ is in one-to-one correspondence with $(\Mb_\ell(r))_{r\not\in H_\ell}$ and $\Mb_{\ell}'$), because we need to account for labelings of these molecules, and in general $\mathrm{Fun}$ is not constant on equivalence classes.

The proof of (\ref{eq.proof_cumulant_single_1_6}) involves some complicated calculation of multinomial coefficients, and will be left to Part 4 below. For now we use it to finish the proof of Proposition \ref{prop.cumulant_single}. Indeed, by inserting \eqref{eq.proof_cumulant_single_1_5} and \eqref{eq.proof_cumulant_single_1_6} into the first line of \eqref{eq.proof_cumulant_single_1_2}, we get
\begin{equation}\label{eq.proof_cumulant_single_1_7}
\begin{aligned}
    \textrm{$\sum(\cdots)$ in (\ref{eq.proof_cumulant_single_1_2})} &= \sum_{[\Mb_\ell]\in \Fc_{ \Lambda_\ell}, \, r(\Mb_\ell) = Q_\ell,\, H_{\ell-1}\subseteq p(\Mb_\ell)} \Ic\Nc_{\Mb_\ell, H_{\ell-1}}^{\mathrm{Pen}}(\vz_{Q_\ell})
    \\
    &=\sum_{H_\ell\subseteq Q_\ell} \bigg(\prod_{r\notin H_\ell} \sum_{[\Mb_\ell(r)]\in \Tc_{\Lambda_\ell}} \Ic\Nc_{\Mb_\ell(r), \varnothing}^{\mathrm{Pen}}(z_r)\bigg)\cdot\sum_{\substack{[\Mb_\ell']\in \Fc_{\Lambda_\ell}, \, r(\Mb_\ell') = H_\ell,\\ H_{\ell-1}\subseteq p(\Mb_\ell'),\,r(\Mb_\ell')\searrow H_{\ell-1}}} \Ic\Nc_{\Mb_\ell', H_{\ell-1}}^{\mathrm{Pen}}(\vz_{H_\ell})
    \\
    &= \sum_{H_\ell\subseteq Q_\ell}(f^{\Ac}(\ell\tau))^{\otimes(Q_\ell\backslash H_\ell)}\cdot E_{H_\ell}(\ell\tau,\vz_{H_\ell}),
\end{aligned}
\end{equation}
where 
\begin{equation}\label{eq.proof_cumulant_single_1_8}
\begin{aligned}
    f^{\Ac}(\ell\tau) &= \sum_{[\Mb_\ell(r)]\in \Tc_{ \Lambda_\ell}} \Ic\Nc_{\Mb_\ell(r), \varnothing}^{\mathrm{Pen}}(z_r)
    \\ 
    E_{H_\ell}(\ell\tau,\vz_{H_\ell}) &= \sum_{\substack{[\Mb_\ell']\in \Fc_{\Lambda_\ell}, \, r(\Mb_\ell') = H_\ell\\ H_{\ell-1}\subseteq p(\Mb_\ell'),\,r(\Mb_\ell')\searrow H_{\ell-1}}} \Ic\Nc_{\Mb_\ell', H_{\ell-1}}^{\mathrm{Pen}}(\vz_{H_\ell}).
\end{aligned}
\end{equation} This proves \eqref{eq.proof_cumulant_single_*}, together with the formulas \eqref{eq.fAterm_single}--\eqref{eq.EHterm_single} in \eqref{it.cumulant_single_1}--\eqref{it.cumulant_single_2} (note that trivially $\Ic\Nc_{\Mb,H}^{\mathrm{Pen}}\leq|\Ic\Nc_{\Mb,H}|$). Therefore, the proof of Proposition \ref{prop.cumulant_single} is complete, pending the proof of (\ref{eq.proof_cumulant_single_1_6}).

\textbf{Proof part 4.} Finally we prove (\ref{eq.proof_cumulant_single_1_6}). Without loss of generality, we assume that $r(\Mb_\ell) = Q_\ell = [s]$. In the expressions below, for convenience of writing, we will replace $\mathrm{Fun}(\Mb_\ell(r))$ by $\bigcirc_r$, and $\mathrm{Fun}(\Mb_\ell',H_{\ell-1})$ by $\bigcirc'$, and their product (in the first line of (\ref{eq.proof_cumulant_single_1_6})) by $\bigcirc$. By (\ref{eq.sum_equiv}), we have
\begin{equation}\label{eq.proof_cumulant_single_4_1}
\begin{aligned}
    \sum_{\substack{[\Mb_\ell]\in \Fc_{\Lambda_\ell}, \, r(\Mb_\ell) = [s]\\ H_{\ell-1}\subseteq p(\Mb_\ell)}} \bigcirc&= \sum_{\substack{\Mb_\ell\in \Fs_{\Lambda_\ell}, \, r(\Mb_\ell) = [s]\\ H_{\ell-1}\subseteq p(\Mb_\ell),\,p(\Mb_\ell) = [|\Mb_\ell|_p]}}\frac{1}{(|\Mb_\ell|_p-s)!}\cdot\bigcirc
    \\
    &= \sum_{n_\ell=s}^\infty \frac{1}{(n_\ell-s)!} \sum_{\substack{\Mb_\ell\in \Fs_{\Lambda_\ell}, \, r(\Mb_\ell) = [s]\\ H_{\ell-1}\subseteq p(\Mb_\ell),\,p(\Mb_\ell) = [n_\ell]}}\bigcirc
    \\
    &= \sum_{n_\ell=s}^\infty \frac{1}{(n_\ell-s)!}\sum_{H_\ell\subseteq [s]}\sum_{\substack{(\Mb_\ell(r))_{r\not\in H_\ell}\\ \Mb_\ell(r)\in \Ts_{\Lambda_\ell},\,r(\Mb_\ell(r))=\{r\}}}\prod_r\bigcirc_r\cdot\sum_{\substack{\Mb_\ell'\in\Fs_{\Lambda_\ell}, H_{\ell-1}\subseteq p(\Mb_\ell')\\p(\Mb_\ell)=[n_\ell],\,r(\Mb_\ell')=H_\ell\\ r(\Mb_\ell')\searrow H_{\ell-1}}}\bigcirc'.
\end{aligned}
\end{equation}

By Definition \ref{def.indicator_general} \eqref{it.S1_Mb_sim_3}, we know that the summand $\bigcirc$ remains invariant under any relabeling that preserves the order of labels within each component. Now consider the summation in (\ref{eq.proof_cumulant_single_4_1}) with $n_\ell$, $H_\ell$ and $|\Mb_\ell'|_{p\backslash r}:=m'$ fixed, and with $|\Mb_\ell(r)|_{p\backslash r}:=m_r$ fixed for each $r$. The number of choices for the \emph{non-root particle sets} of each $\Mb_\ell(r)$ and $\Mb_\ell'$ is given by the multinomial coefficient $\binom{n_\ell-s}{m_1,\cdots,m_r,m'}$. For each fixed choice of \emph{these sets}, we have
\begin{equation}\label{eq.inv1}\sum_{\substack{(\Mb_\ell(r))_{r\not\in H_\ell}:\,p(\Mb_\ell(r))\,\textrm{fixed}\\ \Mb_\ell(r)\in \Ts_{\Lambda_\ell},\,r(\Mb_\ell(r))=\{r\}}}\prod_r\bigcirc_r=\prod_{r\not\in H_\ell}\sum_{\substack{\Mb_\ell(r)\in\Ts_{\Lambda_\ell}:\,p(\Mb_\ell(r))\,\textrm{fixed}\\r(\Mb_\ell(r))=\{r\}}}\bigcirc_r=\prod_{r\not\in H_\ell}m_r!\sum_{\substack{[\Mb_\ell(r)]\in\Tc_{\Lambda_\ell}:|\Mb_\ell(r)|_{p\backslash r}=m_r\\r(\Mb_\ell(r))=\{r\}}}\bigcirc_r,\end{equation} where the last equality is because different choices of $p(\Mb_\ell(r))$ do not affect the expression $\bigcirc_r$ or the average in the equivalence class $[\Mb_\ell(r)]$ by the above invariance property. Similarly we have
\begin{equation}\label{eq.inv2}\sum_{\substack{\Mb_\ell'\in\Fs_{\Lambda_\ell}, H_{\ell-1}\subseteq p(\Mb_\ell')\\p(\Mb_\ell')\,\textrm{fixed},\,r(\Mb_\ell')=H_\ell\\ r(\Mb_\ell')\searrow H_{\ell-1}}}\bigcirc'=(m')!\sum_{\substack{[\Mb_\ell']\in \Fc_{\Lambda_\ell}:|\Mb_\ell'|_{p\backslash r}=m', \, r(\Mb_\ell') = H_\ell\\ H_{\ell-1}\subseteq p(\Mb_\ell'),\,r(\Mb_\ell')\searrow H_{\ell-1}}}\bigcirc'.\end{equation}

By inserting (\ref{eq.inv1})--(\ref{eq.inv2}) into (\ref{eq.proof_cumulant_single_4_1}), and first summing in $\Mb_\ell(r)$ and $\Mb_\ell'$ and then in the non-root particle sets and finally in $m_r$ and $m'$, we get
\begin{equation}\label{eq.proof_cumulant_single_fin}
\begin{aligned}(\ref{eq.proof_cumulant_single_4_1})&=\sum_{n_\ell=s}^\infty\frac{1}{(n_\ell-s)!}\sum_{H_\ell\subseteq [s]}\sum_{\substack{(m_r)_{r\not\in H_\ell},m'\\m'+\sum_r m_r=n_\ell-s}}\binom{n_\ell-s}{m_1,\cdots,m_r,m'}\\
&\times \prod_{r\not\in H_\ell}m_r!\sum_{\substack{[\Mb_\ell(r)]\in\Tc_{\Lambda_\ell}:|\Mb_\ell(r)|_{p\backslash r}=m_r\\r(\Mb_\ell(r))=\{r\}}}\bigcirc_r\cdot (m')!\sum_{\substack{[\Mb_\ell']\in \Fc_{\Lambda_\ell}:|\Mb_\ell'|_{p\backslash r}=m', \, r(\Mb_\ell') = H_\ell\\ H_{\ell-1}\subseteq p(\Mb_\ell'),\,r(\Mb_\ell')\searrow H_{\ell-1}}}\bigcirc'.
\end{aligned}
\end{equation} Using the multinomial formula, we see that the combinatorial factors in (\ref{eq.proof_cumulant_single_fin}) exactly cancel. Let the summations in the second line of (\ref{eq.proof_cumulant_single_fin}) be $X_{m_r}$ and $Y_{m'}$ respectively, then
\begin{equation}\label{eq.proof_cumulant_single_fin2}
\begin{aligned}(\ref{eq.proof_cumulant_single_4_1})&=\sum_{H_\ell\subseteq [s]}\sum_{n_\ell=s}^\infty\sum_{\substack{(m_r)_{r\not\in H_\ell},m'\\m'+\sum_r m_r=n_\ell-s}}\prod_{r\not\in H_\ell} X_{m_r}\cdot Y_{m'}\\
&=\sum_{H_\ell\subseteq [s]}\sum_{(m_r)_{r\not\in H_\ell},m'}\prod_{r\not\in H_\ell} X_{m_r}\cdot Y_{m'}=\sum_{H_\ell\subseteq [s]}\prod_{r\not\in H_\ell}\bigg(\sum_{m_r} X_{m_r}\bigg)\cdot\bigg(\sum_{m'}Y_{m'}\bigg)
\end{aligned}
\end{equation} where we have changed the order of summation. Finally, by definition of $X_{m_r}$ and $Y_{m'}$ we have
\begin{equation}\label{eq.proof_cumulant_single_fin3}\sum_{m_r}X_{m_r}=\sum_{\substack{[\Mb_\ell(r)]\in \Tc_{ \Lambda_\ell}\\ r(\Mb_\ell(r))=\{r\}}}\bigcirc_r=\sum_{[\Mb_\ell(r)]\in \Tc_{ \Lambda_\ell}}\bigcirc_r,\qquad \sum_{m'}Y_{m'}=\sum_{\substack{[\Mb_\ell']\in \Fc_{\Lambda_\ell}, \, r(\Mb_\ell') = H_\ell\\ H_{\ell-1}\subseteq p(\Mb_\ell'),\,r(\Mb_\ell')\searrow H_{\ell-1}}}\bigcirc'\end{equation} (note that $\Mb_\ell(r)$ has only one root particle, and renaming it into $r$ or anything else does not affect the expression $\bigcirc_r$ or the average in the equivalence class). By putting together (\ref{eq.proof_cumulant_single_fin2}) and (\ref{eq.proof_cumulant_single_fin3}), we have proved (\ref{eq.proof_cumulant_single_1_6}). This completes the proof of Proposition \ref{prop.cumulant_single}.
\end{proof}

\subsection{The initial cumulant expansion}\label{sec.initial_cumulant} In this subsection we prove the initial cumulant expansion (\ref{eq.cumulant_l}) for $\ell=0$, without the $\mathrm{Err}$ term, see Proposition \ref{prop.initial_cumulant}. The initial cumulants $E_{H_0}(0)$ correspond to the closeness scenarios $|x_{p_1}-x_{p_2}|\leq\varepsilon$ which, stem from the $\mathbbm{1}_{\Dc_N}$ indicator functions in (\ref{eq.N_par_ensemble}).

\begin{proposition}[Initial cumulant expansion]\label{prop.initial_cumulant} Assume that the initial density function $W_{0,N}(\vz_N)$ is given by \eqref{eq.N_par_ensemble}. Then for each set $Q_0$ with $|Q_0|\leq A_0$ ($A_0$ is defined in \eqref{eq.defLambdaseq}), we have initial cumulant expansion at $t = 0$:
\begin{equation}\label{eq.initial_cumulant_multilayer}
    f(0,\vz_{Q_0})=\sum_{H_0\subseteq Q_0}(f^\Ac(0))^{\otimes (Q_0\backslash H_0)}\cdot E_{H_0}(0,\vz_{H_0}),
\end{equation} 
where $f^\Ac(0,z)=f_0(z)$ and $E_{H_0}(0,\vz_{H_0})$ satisfies the following upper bound
\begin{equation}\label{eq.initial_cumulant_bound}
    |E_{H_0}(0,\vz_{H_0})|\le (f^\Ac(0))^{\otimes H_0}\sum_{\Lc:p(\Lc)=H_0} \mathbbm{1}_\Lc^{\varepsilon}.
\end{equation} 
Here $\Lc$ is a collection of pairs $(p_1,p_2)$ with $p_1,p_2\in H_0$ (cf. Definition \ref{def.molecule} (\ref{it.initial_link}), where bottom ends are identified with particle lines). This $\Lc$ is a forest with no isolated vertex, when viewed as a graph over vertex set $H_0$. Moreover, $\mathbbm{1}_\Lc^{\varepsilon}$ is defined by 
\begin{equation}\label{eq.1_L'}
    \mathbbm{1}_\Lc^{\varepsilon} = \prod_{(p_1, p_2)\in \Lc} \mathbbm{1}_{p_1\sim_{\mathrm{in}} p_2}^\varepsilon\qquad\textrm{and}\qquad \mathbbm{1}_{p_1\sim_{\mathrm{in}} p_2}^{\varepsilon} = \left\{\begin{aligned}
        &1 \quad &&|x_{p_1}-x_{p_2}|\le \varepsilon,
        \\
        &\varepsilon^\upsilon \quad &&|x_{p_1}-x_{p_2}|\ge \varepsilon.
    \end{aligned}\right.
\end{equation} 
where $\upsilon = 3^{-d-1}$, as in (\ref{eq.1_L}).
\end{proposition}
\begin{proof} The following proof is an adaptation of the proof of \cite{PS17}, Proposition A.1. To start with, we introduce the following notations:
\begin{equation}\label{eq.proof_initial_cumulant_notation}
    \mathbbm{1}_{A\not\sim B} = \mathbbm{1}_{\{|z_p-z_{p'}|\ge\varepsilon,\ \forall p\in A\textrm{ and }p\neq p'\in B\}}\quad\textrm{and}\quad\mathbbm{1}_{p\sim B} = \mathbbm{1}_{\{\exists p'\in B,\textrm{ such that }|z_p-z_{p'}|\leq\varepsilon\}}.
\end{equation}
    
Without loss of generality, we assume that $Q_0 = [s]$. By \eqref{eq.N_par_ensemble}, \eqref{eq.partition_func} and \eqref{eq.s_par_cor}, we have
\begin{equation}\label{eq.proof_initial_cumulant_1}
\begin{aligned}
    f(0,\vz_s)&=\frac{1}{\Zc}\sum_{n=0}^\infty\frac{\varepsilon^{-(d-1)n}}{n!}\int_{\Dc_{s+n}} (f^\Ac(0))^{\otimes [s+n]} (\vz_{s+n})\,\mathrm{d}\vz_{[s+1:s+n]},
    \\
    \Zc&=1+\sum_{N=1}^\infty\frac{\varepsilon^{-(d-1)N}}{N!}\int_{\Dc_N} (f^\Ac(0))^{\otimes [N]}(\vz_N)\,\mathrm{d}\vz_N,
\end{aligned}
\end{equation}
which implies that
\begin{equation}\label{eq.proof_initial_cumulant_2}
    f(0,\vz_s)=(f^\Ac(0))^{\otimes [s]}\mathbbm{1}_{[s]\not\sim [s]}\cdot\frac{1}{\Zc}\sum_{n=0}^\infty\frac{\varepsilon^{-(d-1)n}}{n!}\int_{\Dc_{[s+1:s+n]}} \mathbbm{1}_{[s]\not\sim [s+1:s+n]} (f^\Ac(0))^{\otimes [s+1:s+n]}\,\mathrm{d}\vz_{[s+1:s+n]}.
\end{equation} 

Note that
\begin{equation}\label{eq.proof_initial_cumulant_3}
    \mathbbm{1}_{[s]\not\sim [s+1:s+n]} = \prod_{p\in [s]}\mathbbm{1}_{p\not\sim [s+1:s+n]} = \sum_{H_1\subseteq [s]} (-1)^{|H_1|}\cdot\prod_{p\in H_1}\mathbbm{1}_{p\sim [s+1:s+n]}.
\end{equation} We will fix $H_1$ and define $M=[s]\backslash H_1$, then put aside the factors $\mathbbm{1}_{H_1\not\sim H_1}\cdot\mathbbm{1}_{H_1\not\sim M}$ that occur in the $\mathbbm{1}_{[s]\not\sim[s]}$ in (\ref{eq.proof_initial_cumulant_2}), and focus on $\mathbbm{1}_{M\not\sim M}$. Here we apply a Penrose type argument in the setting of Part 2 of the proof of Proposition \ref{prop.S_N_decomposition}, by setting
\begin{itemize}
\item $E = \{(p, p'):p\neq p'\in M\}$, so subsets of $E$ are identified with graphs over vertex set $M$, and denote $\mathbbm{1}_{A_{\of}}:=\mathbbm{1}_{p\sim p'}$ if $\of=(p,p')$.
\item There is no $\Fs^{\mathrm{err}}$, and $\Fs$ is the set of $F\subseteq E$ whose corresponding graph are forests;
\item The mapping $\pi:2^E\to\Fs$ is any algorithm of generating a spanning forest $F=\pi(G)$ from a graph $G$. For example, we can fix a linear ordering of elements of $E$ based on labels of particles, then apply the same procedure that is applied to the graph $\mathrm{Cl}_G$ in (\ref{it.penrose_alg_3}) in the definition of the algorithm in Part 4 of the proof of Proposition \ref{prop.S_N_decomposition}.
\end{itemize} It is easy to see, as in Part 4 of the proof of Proposition \ref{prop.S_N_decomposition}, that (\ref{eq.pvfbset}) is true for the above-defined $\Fs$ and $\pi$. Then we get

\begin{equation}\label{eq.proof_initial_cumulant_4}
    \mathbbm{1}_{M\not\sim M} = \prod_{\of\in E}(1-\mathbbm{1}_{A_\of}) = \sum_{\Lc_2} (-1)^{|\Lc_2|}\mathbbm{1}_{\Lc_2}^{\mathrm{Pen}},
\end{equation}
where $\Lc_2$ is a graph over vertex set $M=[s]\backslash H_1$ that is a forest (corresponding to $F\in\Fs$ above), $|\Lc_2|=|F|$ is the number of edges in the graph $\Lc_2$, $p(\Lc_2)$ is the set of non-isolated vertices in $\Lc_2$, and
\[\mathbbm{1}_{\Lc_2}^{\mathrm{Pen}} := \prod_{(p,p') \in \Lc_2} \mathbbm{1}_{p\sim p'} \prod_{(p,p') \in AL(F)} \mathbbm{1}_{p\not\sim p'};\qquad \mathbbm{1}_{\Lc_2}^{\mathrm{Pen}}\le \mathbbm{1}_{\Lc_2}:=\prod_{(p,p') \in \Lc_2} \mathbbm{1}_{p\sim p'}.\]

Inserting \eqref{eq.proof_initial_cumulant_3} into \eqref{eq.proof_initial_cumulant_2}, and letting $H_1\cup p(\Lc_2)=H_0$, we get
\begin{align}\label{eq.proof_initial_cumulant_5}
    f(0,\vz_s) &= \sum_{H_0\subseteq [s]}(f^\Ac(0))^{\otimes ([s]\backslash H_0)}\cdot E_H(0,\vz_{H_0}), \\
    E_H(0,\vz_{H_0}) &= (f^\Ac(0))^{\otimes H}\frac{1}{\Zc}\sum_{n=0}^\infty\frac{\varepsilon^{-(d-1)n}}{n!}\sum_{(H_1,\Lc_2):H_1\cup p(\Lc_2)=H_0} (-1)^{|\Lc_2|+|H_1|}\mathbbm{1}_{H_1\not\sim H_1}\cdot\mathbbm{1}_{H_1\not\sim ([s]\backslash H_1)}\cdot\mathbbm{1}_{\Lc_2}^{\mathrm{Pen}}\nonumber\\
    \label{eq.proof_initial_cumulant_6}&\times\int_{\Dc_{[s+1:s+n]}} \prod_{p\in H_1}\mathbbm{1}_{p\sim [s+1:s+n]}\cdot (f^\Ac(0))^{\otimes [s+1:s+n]} \,\mathrm{d}\vz_{[s+1:s+n]}.
\end{align}

Recall $\upsilon=3^{-d-1}$. Let $K\subseteq [s+1:s+n]$ be the set of $p'$ such that $p'\sim p$ for some $p\in H_1$. Since each $p'$ corresponds to at most $(3\upsilon)^{-1}=3^d$ different $p$ (the number of non-intersecting unit balls that all intersect a given unit ball cannot exceed $3^d$, due to a simple covering argument), we conclude that $|K|\ge 3\upsilon |H_1|$, thus
\begin{equation}\label{eq.proof_initial_cumulant_7}
    \prod_{p\in H_1}\mathbbm{1}_{p\sim [s+1:s+n]} \le \sum_{K:K\subseteq [s+1:s+n],\, |K|\ge 3\upsilon |H_1|}\, \prod_{p'\in K}\mathbbm{1}_{p'\sim H_1}.
\end{equation}

Denote the integral in (\ref{eq.proof_initial_cumulant_6}) by $(\mathrm{Int})$, then by (\ref{eq.proof_initial_cumulant_7}) we have
\begin{equation}\label{eq.proof_initial_cumulant_8}
\begin{aligned}
    &\mathrm{(Int)} \le \sum_{K:K\subseteq [s+1:s+n],\, |K|\ge 3\upsilon |H_1|}\int_{\Dc_{[s+1:s+n]}} \prod_{p'\in K}\mathbbm{1}_{p'\sim H_1}\cdot (f^\Ac(0))^{\otimes [s+1:s+n]} \,\mathrm{d}\vz_{[s+1:s+n]}
    \\
    &\le \sum_{K:K\subseteq [s+1:s+n],\, |K|\ge 3\upsilon |H_1|} \int_{\Dc_{[s+1:s+n]\backslash K}} (f^\Ac(0))^{\otimes [s+1:s+n]\backslash K} \,\mathrm{d}\vz_{[s+1:s+n]\backslash K}\cdot \prod_{p'\in K} \int \mathbbm{1}_{p'\sim H_1} f^\Ac(0, z_{p'})\,\mathrm{d} z_{p'}
    \\
    &\le \sum_{K:K\subseteq [s+1:s+n],\, |K|\ge 3\upsilon |H_1|} \varepsilon^{(d-\frac{1}{2})|K|} \int_{\Dc_{[s+1:s+n]\backslash K}} (f^\Ac(0))^{\otimes [s+1:s+n]\backslash K} \,\mathrm{d}\vz_{[s+1:s+n]\backslash K}
    \\
    &\le \sum_{3\upsilon |H_1|\le m\le n } \varepsilon^{(d-\frac{1}{2})m}\binom{n}{m} \int_{\Dc_{[n-m]}} (f^\Ac(0))^{\otimes [n-m]} \,\mathrm{d}\vz_{[n-m]}.
\end{aligned}
\end{equation} Here in the third line we have used $\int \mathbbm{1}_{p'\sim H_1} f^\Ac(0, z_{p'})\,\mathrm{d} z_{p'}\le |H_1|\varepsilon^{d-1} \le \varepsilon^{d-\frac{1}{2}}$. In the fourth line we have fixed $|K|=m$ which leaves $\binom{n}{m}$ choices for $K$, and renamed $\Dc_{[s+1:s+n]\backslash K}$ into $\Dc_{[n-m]}$.

Finally, inserting (\ref{eq.proof_initial_cumulant_8}) into \eqref{eq.proof_initial_cumulant_6}, we get
\begin{equation}\label{eq.proof_initial_cumulant_9}
\begin{aligned}
    &|E_{H_0}(0,\vz_{H_0})| \leq (f^\Ac(0))^{\otimes H_0} \frac{1}{\Zc}\sum_{n=0}^\infty\frac{\varepsilon^{-(d-1)n}}{n!}\sum_{(H_1,\Lc_2):H_1\cup p(\Lc_2)=H_0} \mathbbm{1}_{H_1\not\sim H_1}\cdot\mathbbm{1}_{H_1\not\sim ([s]\backslash H_1)}\cdot\mathbbm{1}_{\Lc_2}^{\mathrm{Pen}}
    \\
    &\phantom{|E_{H_0}(0,\vz_{H_0})|}\times\sum_{3\upsilon |H_1|\le m\le n }\varepsilon^{(d-\frac{1}{2})m} \binom{n}{m} \int_{\Dc_{[n-m]}} (f^\Ac(0))^{\otimes [n-m]} \,\mathrm{d}\vz_{[n-m]}
    \\
    &=(f^\Ac(0))^{\otimes H_0}\sum_{H_1\cup p(\Lc_2)=H_0}\mathbbm{1}_{\Lc_2} \sum_{m\ge 3\upsilon |H_1|}\frac{\varepsilon^{\frac{m}{2}}}{m!}\frac{1}{\Zc}\underbrace{\left(\sum_{n=m}^\infty\frac{\varepsilon^{-(d-1)(n-m)}}{(n-m)!} \int_{\Dc_{[n-m]}} (f^\Ac(0))^{\otimes [n-m]} \,\mathrm{d}\vz_{[n-m]}\right)}_{=\Zc,\textrm{ by \eqref{eq.partition_func} and the change of variable }N=n-m}
    \\
    &=(f^\Ac(0))^{\otimes H_0}\sum_{H_1\cup p(\Lc_2)=H_0}\mathbbm{1}_{\Lc_2} \sum_{m\ge 3\upsilon |H_1|}\frac{\varepsilon^{\frac{m}{2}}}{m!} 
    \le \sum_{H_1\cup p(\Lc_2)=H_0}\varepsilon^{\upsilon |H_1|}\cdot\mathbbm{1}_{\Lc_2} 
    \\
    &\le(f^\Ac(0))^{\otimes H_0} \sum_{\Lc: p(\Lc)=H_0}\mathbbm{1}_\Lc^{\varepsilon}.
\end{aligned}
\end{equation}
Here in the last step, we choose $\Lc=\Lc_1\cup\Lc_2$, where $\Lc_1$ is an arbitrary forest with vertex set $H_1$; since $\mathbbm{1}_{p_1\sim_{\mathrm{in}} p_2}\ge \varepsilon^{\upsilon}$, it follows that $\mathbbm{1}_\Lc^{\varepsilon}\ge \varepsilon^{\upsilon |H_1|}\cdot\mathbbm{1}_{\Lc_2}$. This completes the proof.
\end{proof}

\subsection{Proof of Proposition \ref{prop.cumulant_formula}}\label{sec.proof_cumulant_formula_final} In this subsection, we put together Propositions \ref{prop.cumulant_single} and \ref{prop.initial_cumulant}, to prove Proposition \ref{prop.cumulant_formula}.
\begin{proof}[Proof of Proposition \ref{prop.cumulant_formula}] We divide the proof into 5 parts.

\textbf{Proof part 1.} We prove there exists $E_H=E_H(\ell\tau,\vz_H)$ that verifies (\ref{eq.cumulant_expansion}). In fact, this follows immediately by induction on $\ell$, using Proposition \ref{prop.initial_cumulant} (base case $\ell=0$) and Proposition \ref{prop.cumulant_single} (inductive step).

\textbf{Proof part 2.} We prove that $f^\Ac$ in (\ref{eq.cumulant_expansion}) satisfies (\ref{eq.fAterm}). In fact, we have $f^\Ac(0)=f_0$ by Proposition \ref{prop.initial_cumulant}; by \ref{eq.fAterm_single}), we see that $f^\Ac(\ell\tau)$ satisfies
\begin{equation}\label{eq.proof_cumulant_multi_2_1}
    f^{\mathcal{A}}(\ell\tau,z) = \sum_{[\Mb]\in \mathcal{T}_{\Lambda_\ell}}\Ic\Nc_{\Mb, \varnothing}^{\mathrm{Pen}}(z)
\end{equation}
with $z=\vz_{r(\Mb)}$. Moreover, by first line of \eqref{eq.associated_int_notsim_single}, we see that $\Ic\Nc_{\Mb, \varnothing}^{\mathrm{Pen}}(z)$ is given by
\begin{equation}\label{eq.proof_cumulant_multi_2_2}
    \Ic\Nc_{\Mb, \varnothing}^{\mathrm{Pen}}(z) = \varepsilon^{-(d-1)|\Mb|_{p\backslash r}}\int_{\Rb^{2d|\Mb|_{p\backslash r}}} (\Sc\circ\mathbbm{1})_{\Mb}^{\mathrm{Pen}}\left[\left(f^\Ac((\ell-1)\tau)\right)^{\otimes (p(\Mb))}\right]\mathrm{d}\vz_{(p\backslash r)(\Mb)}.
\end{equation} Recall also the quantity $\Ic\Nc_\Mb(z)$ defined in (\ref{eq.associated_integral_molecule}) in Definition \ref{def.associated_int}; using that $\mathbbm{1}_\Lc^{\varepsilon}=1$ ($[\Mb]\in\Tc_{\Lambda_\ell}$ has no initial link) and $\ell_1[\pb]-\underline{\ell}=0$ ($\Mb$ has only one layer $\ell$), we see that $\Ic\Nc_\Mb(z)$ is the same as $\Ic\Nc_{\Mb, \varnothing}^{\mathrm{Pen}}(z)$ but with $(\Sc\circ\mathbbm{1})_{\Mb}^{\mathrm{Pen}}$ in (\ref{eq.proof_cumulant_multi_2_2}) replaced by $(\Sc\circ\mathbbm{1})_{\Mb}$. Similarly, the quantity $|\Ic\Nc_{\Mb,\varnothing}|(z)$ defined in \eqref{eq.associated_int_notsim_single} is identical with the quantity $|\Ic\Nc_\Mb|(z)$ defined in (\ref{eq.associated_integral_molecule_abs}).

Now, by inserting \eqref{eq.S1_Mb_sim_asym} into (\ref{eq.proof_cumulant_multi_2_2}), we get
\begin{equation}\label{eq.proof_cumulant_multi_2_3}
\begin{aligned}
    \Ic\Nc_{\Mb, \varnothing}^{\mathrm{Pen}}(z) &= \varepsilon^{-(d-1)|\Mb|_{p\backslash r}}\int_{\Rb^{2d|\Mb|_{p\backslash r}}} \bigg[(\Sc\circ\mathbbm{1})_{\Mb}\left(f^\Ac((\ell-1)\tau)\right)^{\otimes (p(\Mb))}
    \\
    &\phantom{= \varepsilon^{-(d-1)|\Mb|_{p\backslash r}}\int_{\Rb^{2d|\Mb|_{p\backslash r}}}\bigg[}+ \sum_{(e,e')}O_1\left(|(\Sc\circ\mathbbm{1})_{\widetilde{\Mb}}| |f^\Ac((\ell-1)\tau)|^{\otimes (p(\Mb))}\right)\bigg]\,\mathrm{d}\vz_{(p\backslash r)(\Mb)}
    \\
    &= \Ic\Nc_{\Mb}(z) +\sum_{(e,e')} O_1\big(|\Ic\Nc_{\widetilde{\Mb}}|(z)\big).
\end{aligned}
\end{equation}
Here $|O_1(A)|\leq A$, and $\widetilde{\Mb}=\Mb_{e,e'}$ is obtained from $\Mb$ by adding one more overlap $\of$ joining edges $e$ and $e'$, so that $[\widetilde{\Mb}]\in\Tc_{\Lambda_\ell}^{\mathrm{err}}$ if $[\Mb]\in\Tc_{\Lambda_\ell}$. 

By inserting (\ref{eq.proof_cumulant_multi_2_3}) into \eqref{eq.proof_cumulant_multi_2_1}, we get
\[f^\Ac(\ell\tau,z)=\sum_{[\Mb]\in\Tc_{\Lambda_\ell}}\Ic\Nc_\Mb(z)+\sum_{(\Mb,\widetilde{\Mb})}O_1\big(|\Ic\Nc_{\widetilde{\Mb}}|(z)\big),\] where $[\widetilde{\Mb}]\in\Tc_{\Lambda_\ell}^{\mathrm{err}}$ is obtained from $[\Mb]\in\Tc_{\Lambda_\ell}$ by adding one O-atom as above. This completes the proof of \eqref{eq.fAterm}, since the number of $[\Mb]$ with fixed $[\widetilde{\Mb}]$ is bounded by $|\log\varepsilon|^{C^*}$.

\textbf{Proof part 3.} We prove that $E_H$ in (\ref{eq.cumulant_expansion}) satisfies (\ref{eq.cumulant_formula}). In fact, by \eqref{eq.EHterm_single} and the second line of \eqref{eq.associated_int_notsim_single}, we have the recurrence inequality
\begin{equation}\label{eq.proof_cumulant_multi_3_1}
\begin{aligned}
    |E_{H_\ell}(\ell\tau,\vz_{H_\ell})| &\le \sum_{\substack{[\Mb_\ell]\in \mathcal{F}_{ \Lambda_\ell},\, r(\Mb_\ell) = H_\ell\\H_{\ell-1}\subseteq p(\Mb_\ell),\,r(\Mb_\ell)\searrow H_{\ell-1}}} \varepsilon^{-(d-1)|\Mb_\ell|_{p\backslash r}}\int_{\Rb^{2d|\Mb_\ell|_{p\backslash r}}} |(\Sc\circ\mathbbm{1})_{\Mb_\ell}|
    \\
    &\qquad\qquad\qquad\left[\left|f^\Ac((\ell-1)\tau)\right|^{\otimes (p(\Mb_\ell)\backslash H_{\ell-1})}\cdot |E_{H_{\ell-1}}((\ell-1)\tau,\vz_{H_{\ell-1}})|\right]\mathrm{d}\vz_{(p\backslash r)(\Mb_\ell)}.
\end{aligned}
\end{equation}  
This inequality can be schematically written as $|E_{H_\ell}(\ell\tau,\vz_{H_\ell})|\le \mathrm{Op}_\ell(|E_{H_{\ell-1}}((\ell-1)\tau,\vz_{H_{\ell-1}})|)$. Iterating it then yields that $|E_{H_\ell}(\ell\tau,\vz_{H_\ell})|\le (\mathrm{Op}_\ell\circ\cdots\circ\mathrm{Op}_1)(|E_{H_0}(0,\vz_{H_0})|)$. Using also (\ref{eq.initial_cumulant_bound}), this leads to
\begin{equation}\label{eq.proof_cumulant_multi_3_2}
\begin{aligned}
    |E_{H_\ell}(\ell\tau,\vz_{H_\ell})| &\leq  \sum_{\substack{([\Mb_{\ell'}])_{\ell'\in[1:\ell]}: \,[\Mb_{\ell'}]\in \mathcal{F}_{\Lambda_{\ell'}}\\ r(\Mb_{\ell'}) = H_{\ell'},\,H_{\ell'-1}\subseteq p(\Mb_{\ell'})\\r(\Mb_{\ell'})\searrow H_{\ell'-1}}}\sum_{\Lc} \prod_{\ell' = 1}^\ell \varepsilon^{-(d-1)|\Mb_{\ell'}|_{p\backslash r}}\int_{\Rb^{2d|\Mb_{\ell'}|_{p\backslash r}}}\mathrm{d}\vz_{(p\backslash r)(\Mb_{\ell'})}
    \\
    & \times\prod_{\ell' = 1}^\ell|(\Sc\circ\mathbbm{1})_{\Mb_{\ell'}}|\bigg(\prod_{\ell' = 1}^\ell\left|f^\Ac((\ell'-1)\tau)\right|^{\otimes (p(\Mb_{\ell'})\backslash H_{\ell'-1})}|f^\Ac(0)|^{\otimes(H_0)}\cdot \mathbbm{1}_\Lc^\varepsilon\bigg),
\end{aligned}
\end{equation} where in the product $\prod_{\ell'=1}^\ell$, the higher values $\ell'$ occur on the left; also $\Lc$ is as in Proposition \ref{prop.initial_cumulant}. Note that (\ref{eq.proof_cumulant_multi_3_2}) is a sum over equivalence classed of single layer molecules $[\Mb_{\ell'}]$ for each $\ell'\in[1:\ell]$, and the summand is constant in each equivalence class (which is easy to verify). Then, from now on, we may treat these equivalence classes as unlabeled molecules (cf. Remark \ref{rem.label}).

\textbf{Constructing $\Mb$ from $\Mb_{\ell'}$.} Next, we shall combine all these $\Mb_{\ell'}$ into one multi-layer molecule $\Mb$ by the following procedure: for each $\ell'\in[1:\ell-1]$, note that $H_{\ell'}=r(\Mb_{\ell'})\subseteq p(\Mb_{\ell'})\cap p(\Mb_{\ell'+1})$ in (\ref{eq.proof_cumulant_multi_3_2}). For each particle line in $H_{\ell'}$, it has a unique top end in $\Mb_{\ell'}$ and a unique bottom end in $\Mb_{\ell'+1}$, and we merge them into one bond. See {\color{blue}Figure \ref{fig.mergelayer}}.
\begin{figure}[h!]
\includegraphics[width=0.6\linewidth]{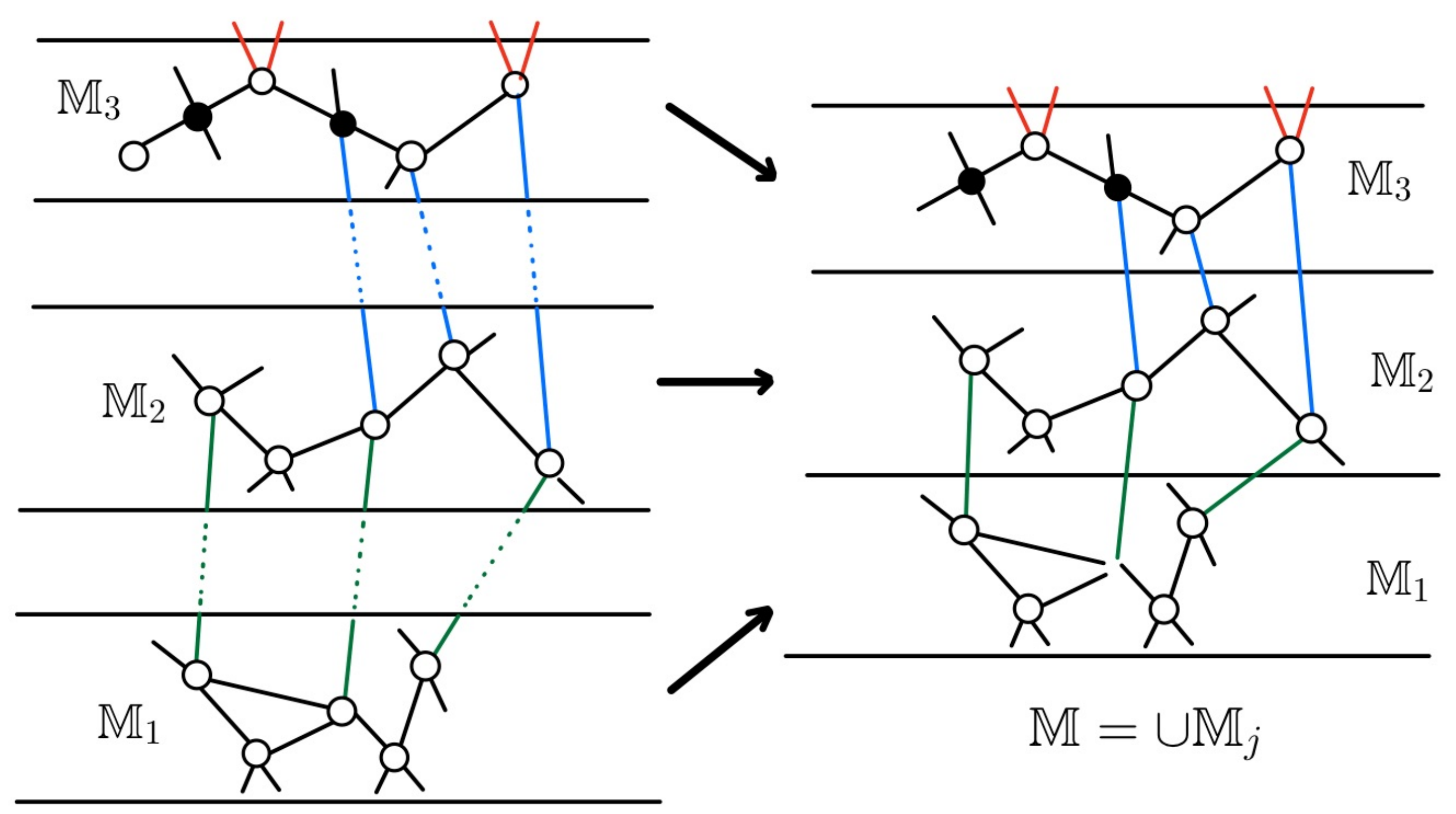}
\caption{Construction of $\Mb$ from single layer molecules $\Mb_{\ell'}$. Here $\ell=3$; the red (resp. blue and green) edges belong to $H_3=r(\Mb_3)$ (resp. $H_2$ and $H_1$). For example, in the left picture, each of the 3 blue top ends in $\Mb_2$ and 3 bottom ends in $\Mb_3$ belongs to a particle line in $H_2=r(\Mb_2)$; in the right picture, we merge the blue ends that belong to the same particle line, to form 3 blue bonds.
}
\label{fig.mergelayer}
\end{figure}

\textbf{Expressing $E_{H_\ell}$ using $\Mb$.} With the above correspondence between $\Mb$ and $\Mb_{\ell'}$, we can verify that
\begin{equation}\label{eq.change_of_sum}\sum_{\substack{([\Mb_{\ell'}])_{\ell'\in[1:\ell]}: \,[\Mb_{\ell'}]\in \mathcal{F}_{\Lambda_{\ell'}}\\ r(\Mb_{\ell'}) = H_{\ell'},\,H_{\ell'-1}\subseteq p(\Mb_{\ell'})\\r(\Mb_{\ell'})\searrow H_{\ell'-1}}}\sum_{\Lc}(\cdots)=\sum_{[\Mb]\in\Fc_{\vLambda_{\ell'}},r(\Mb)=H_\ell}(\cdots),\end{equation}for the summation in (\ref{eq.proof_cumulant_multi_3_2}), where $\Fc_{\vLambda_\ell}$ is defined in Definition \ref{def.set_T_F} (\ref{it.set_F}). Here we note that
\begin{itemize}
\item Properties (\ref{it.set_F_l_1})--(\ref{it.set_F_l_3}) in Definition \ref{def.set_T_F} exactly correspond to the assumptions $[\Mb_{\ell'}]\in \Fc_{\Lambda_{\ell'}}$ in (\ref{eq.proof_cumulant_multi_3_2}), cf. Definition \ref{def.set_F_single_layer} and Definition \ref{def.indicator_general} (\ref{it.setF}); note that the set $[s]$ in Definition \ref{def.indicator_general} (\ref{it.setF1}) should be replaced by $H_\ell$ (which is another fixed set).
\item Property (\ref{it.set_F_l_4}) in Definition \ref{def.set_T_F} exactly correspond to the assumptions $r(\Mb_{\ell'})\searrow H_{\ell'-1}$ in (\ref{eq.proof_cumulant_multi_3_2}), cf. Definition \ref{def.r_conn_H}. Note that for $\ell'=1$, the set $H_0$ is precisely those involved in initial links $\Lc$ (see Proposition \ref{prop.initial_cumulant}), which is also consistent with Definition \ref{def.set_T_F} (\ref{it.set_F_l_4}).
\item Finally, the assumptions $r(\Mb_{\ell'})=H_{\ell'}$ in (\ref{eq.proof_cumulant_multi_3_2}) uniquely links $H_{\ell'}$ to $\Mb$, and the assumptions $H_{\ell'-1}\subseteq p(\Mb_{\ell'})$ allows to define the above merging process to form $\Mb$.
\end{itemize} 

It is also easy to see that
\begin{equation}\label{eq.multi_layer_new_1}(p\backslash r)(\Mb)=\bigsqcup_{\ell'=1}^{\ell}(p\backslash r)(\Mb_{\ell'}),\qquad p(\Mb)=\bigsqcup_{\ell'=1}^{\ell}(p(\Mb_{\ell'})\backslash H_{\ell'-1})\sqcup H_0,\end{equation}using (\ref{eq.pset})--(\ref{eq.rset}) in Definition \ref{def.prsets}. In particular, we may identify the measures and the input factors: \begin{equation}\label{eq.multi_layer_new_2}\prod_{\ell'=1}^\ell \left(\varepsilon^{-(d-1)|\Mb_{\ell'}|_{p\backslash r}}\,\mathrm{d}\vz_{(p\backslash r)(\Mb_{\ell'})}\right)=\varepsilon^{-(d-1)|\Mb|_{p\backslash r}}\,\mathrm{d}\vz_{(p\backslash r)(\Mb)},\end{equation} \begin{equation}\label{eq.multi_layer_new_3}\prod_{\ell'}\left|f^\Ac((\ell'-1)\tau)\right|^{\otimes (p(\Mb_{\ell'})\backslash H_{\ell'-1})}|f^\Ac(0)|^{\otimes(H_0)}=\prod_{\pb\in p(\Mb)}|f^\Ac((\ell_1[\pb]-1)\tau)|,\end{equation} where $\ell_1[\pb]=\ell'$ if $\pb\in p(\Mb_{\ell'})\backslash H_{\ell'-1}$ and $\ell_1[\pb]=1$ if $\pb\in H_0$, which follows from (\ref{eq.pset})--(\ref{eq.rset}) and $H_{\ell'}=r(\Mb_{\ell'})$.

Next we treat the composition operator $\prod_{\ell'=1}^\ell|(\Sc\circ\mathbbm{1})_{\Mb_{\ell'}}|$ in (\ref{eq.proof_cumulant_multi_3_2}). We claim that
\begin{equation}\label{eq.multi_layer_new_4}
\prod_{\ell'=1}^\ell|(\Sc\circ\mathbbm{1})_{\Mb_{\ell'}}|=\Rc_2^{-1}\circ|(\Sc\circ\mathbbm{1})_\Mb|\circ \Rc_1^{-1},
\end{equation} where $\Rc_jf:=f\circ \Tc_j^{-1}$ with
\begin{equation}\label{eq.multi_layer_new_5}
\Tc_1:(x_\pb',v_\pb')\mapsto (x_\pb'+(\ell_1[\pb]-1)\tau\cdot v_\pb',v_\pb'),\quad \Tc_2:(x_\pb,v_\pb)\mapsto (x_\pb+(\ell-\ell_2[\pb])_+\tau\cdot v_\pb,v_\pb).
\end{equation} To prove (\ref{eq.multi_layer_new_4}), recall the flow map $\Hc_\Mb$ defined in Definition \ref{def.associated_op_nonlocal} (\ref{it.prescribed_flow}). Using that $|(\Sc\circ\mathbbm{1})_\Mb| f=(\mathbbm{1}_\Mb\cdot f)\circ \Hc_\Mb^{-1}$ (see (\ref{eq.associated_op}), note that $\Hc_\Mb$ is injective on the support of $\mathbbm{1}_\Mb$) and similarly for $\Mb_{\ell'}$, we only need to prove\begin{equation}\label{eq.multi_layer_new_6}\Hc_\Mb(\ell\tau)=\Tc_2\circ(\Hc_{\Mb_\ell}(\tau)\circ\cdots\circ\Hc_{\Mb_1}(\tau))\circ\Tc_1\end{equation} on the support of $\mathbbm{1}_\Mb$. To prove (\ref{eq.multi_layer_new_6}), simply note that in the prescribed dynamics $\Hc_\Mb$, each particle $\pb$ first moves linearly from time $0$ to time $(\ell_1[\pb]-1)\tau$ (i.e. before its start layer $\ell_1[\pb]$ at which time $\pb$ starts to occur in the $\Mb$-prescribed dynamics), then moves according to the dynamics $\Hc_{\Mb_\ell}(\tau)\circ\cdots\circ\Hc_{\Mb_1}(\tau)$ (with only those in the lifespan of $\pb$ having effect), and then moves linearly from time\footnote{The $\ell_2[\pb]$ may be $\ell+1$ by convention, but the relevant time would still be $\ell\tau$, hence the $(\ell-\ell_2[\pb])_+$ in (\ref{eq.multi_layer_new_5}).} $\ell_2[\pb]\tau$ to time $\ell\tau$ (i.e. after the finish layer $\ell_2[\pb]$ at which time $\pb$ no longer occurs in the $\Mb$-prescribed dynamics). By composing these three parts of dynamics, we obtain  (\ref{eq.multi_layer_new_6}) and hence  (\ref{eq.multi_layer_new_4}).

Now, using (\ref{eq.multi_layer_new_4}), we can replace the $\prod_{\ell'=1}^\ell|(\Sc\circ\mathbbm{1})_{\Mb_{\ell'}}|$ in (\ref{eq.proof_cumulant_multi_3_2}), up to the two linear transports $\Rc_1^{-1}$ and $\Rc_2$. The $\Rc_1^{-1}$ acts on the $f^\Ac$ factors and turns them into $|f^\Ac((\ell_1[\pb]-1)\tau,x_\pb'+(\ell_1[\pb]-1)\tau\cdot v_\pb',v_\pb')|$ as occur in (\ref{eq.extraQ}) with $\underline{\ell}=1$, and does not affect $\mathbbm{1}_\Lc^{\varepsilon}$ as $\ell_1[\pb]=1$ for initial links. The $\Rc_2^{-1}$ has not effect on the integral, as they act only on the non-root variables of $\Mb$ (for root particles $\pb\in r(\Mb)$ we have $r_2[\pb]\in\{\ell,\ell+1\}$ thus $(\ell-\ell_2[\pb])_+=0$) which are integrated, and the linear transport $\Rc_2^{-1}$ preserves integration.

\textbf{Finishing the proof of (\ref{eq.cumulant_formula}).} Let $H=H_\ell$. By putting together all the above and inserting (\ref{eq.change_of_sum}), (\ref{eq.multi_layer_new_2})--(\ref{eq.multi_layer_new_4}) into (\ref{eq.proof_cumulant_multi_3_2}), we get
\begin{multline}\label{eq.proof_cumulant_multi_3_4}
    |E_{H}(\ell\tau,\vz_{H})| \le \sum_{[\Mb]\in \mathcal{F}_{\boldsymbol{\Lambda}_\ell},r(\Mb) = H} \varepsilon^{-(d-1)|\Mb|_{p\backslash r}}\int_{\Rb^{2d|\Mb|_{p\backslash r}}}|(\Sc\circ\mathbbm{1})_\Mb|
    \\
    \bigg(\mathbbm{1}_\Lc^{\varepsilon}\cdot\prod_{\pb\in p(\Mb)} |f^\Ac((\ell_1[\pb]-1)\tau,x_\pb'+(\ell_1[\pb]-\underline{\ell})\tau\cdot v_\pb',v_\pb')|\bigg)\mathrm{d}\vz_{(p\backslash r)(\Mb)}, 
\end{multline} which proves (\ref{eq.cumulant_formula}).

\textbf{Proof part 4.} We prove that $\mathrm{Err}$ in (\ref{eq.cumulant_expansion}) satisfies \eqref{eq.cumulant_formula_err}--\eqref{eq.cumulant_formula_err2}. Note $\mathrm{Err}(0)=0$ by (\ref{eq.initial_cumulant_multilayer}). We define $\mathrm{Err}^1$ and $\mathrm{Err}^2$ as in Proposition \ref{prop.cumulant_single} (\ref{it.cumulant_single_3}), then \eqref{eq.cumulant_formula_err2} follows immediately from \eqref{eq.Err2term}. The estimate \eqref{eq.cumulant_formula_err} for $\mathrm{Err}^1$ follows from the same argument as in the proof of (\ref{eq.cumulant_formula}) in Part 3 above (starting from (\ref{eq.Err1term})), with the only difference being the layer $\ell$ molecule $\Mb_{\ell}$: we have $\Mb_\ell\in\Fc_{\Lambda_\ell}^{\mathrm{err}}$ instead of $\Fc_{\Lambda_\ell}$, that $p(\Mb_\ell)=Q_\ell$ instead of $H_\ell$, and do not require $r(\Mb_\ell)\searrow H_{\ell-1}$. This difference leads to the condition for $\Fc_{\vLambda_\ell}^{\mathrm{err}}$ in Definition \ref{def.set_T_F} (\ref{it.set_F_err}), hence $\Mb\in\Fc_{\vLambda_\ell}^{\mathrm{err}}$ in \eqref{eq.cumulant_formula_err}. This completes the proof of \eqref{eq.cumulant_formula_err}--\eqref{eq.cumulant_formula_err2}.

\textbf{Proof part 5.} Finally we prove that $f_s^{\mathrm{err}}$ in (\ref{eq.f_s_dec}) satisfies (\ref{eq.cumulant_formula_err_trunc}). Recall that this $f_s^{\mathrm{err}}$ is defined in (\ref{eq.truncated_error}). We start with an upper bound for $\mathbbm{1}_{\mathcal{E}_N^{\Lambda,\Gamma}}$ in Definition \ref{def.truncated_domain}.

\textbf{Upper bound for $\mathbbm{1}_{\mathcal{E}_N^{\Lambda,\Gamma}}$.} We first prove the following upper bound:
\begin{equation}\label{eq.error_upper_grammar_rec}\mathbbm{1}_{\Ec_N^{\Lambda,\Gamma}}\leq \sum_{\Mb\in\Fs_{\Lambda_\ell}^{\mathrm{trc.err}},\,r(\Mb)=[s]}\mathbbm{1}_\Mb,
    \end{equation} where $\mathbbm{1}_\Mb$ is defined in (\ref{eq.ind_1M}) in Definition \ref{def.associated_op_nonlocal}. $\Fs_{\Lambda_\ell}^{\mathrm{trc.err}}$ is the collection of all labeled molecules $\Mb$ of a single layer $\ell$ that is connected apart from possible empty ends (which are root particle lines), has only C-atoms, and satisfies (\ref{eq.trc_err}) in Definition \ref{def.set_T_F} (\ref{it.set_F_err_2}). Define $\Fc_{\Lambda_\ell}^{\mathrm{trc.err}}$ to be the set of corresponding equivalence classes.
    
    In fact, to prove (\ref{eq.error_upper_grammar_rec}), we only need to prove that, if $\vz^0\in\Ec_N^{\Lambda,\Gamma}$, then there exists $\Mb\in\Fs_{\Lambda_\ell}^{\mathrm{trc.err}}$ such that $\mathbbm{1}_\Mb(\vz^0)=1$. To prove this, consider the trajectory $\vz^E(t)$ defined by the E-dynamics with initial configuration $\vz^0$. For each time $t\in[0,\tau]$, define $\Mb(t)$ to be the C-molecule topological reduction of the trajectory up to time $t$.

By Definition \ref{def.truncated_domain}, we know that for $t=\tau$, the molecule $\Mb(\tau)$ has a component $\Mb_i$ with either $|\Mb_i|_p>\Lambda$ or $\rho(\Mb_i)>\Gamma$; in this case we say $\Mb_i$ is \emph{large}. Consider the \emph{first time} $t_*$ such that $\Mb(t_*^+)$ has a large component $\Mb_i$; there must be a collision $\nf$ happening at $t_*$, such that $\Mb=\Mb(t_*^+)$ is obtained from $\Mb(t_*^-)$ by adding one more collision $\nf$. Note that this collision either happens within a component $\Mb_j$ of $\Mb(t_*^-)$ (so does not affect $|\Mb_j|_p$ and increases $\rho(\Mb_j)$ by $1$), or merges two components $\Mb_j$ and $\Mb_k$ of $\Mb(t_*^-)$ together (so the new component $\Mb_i$ of $\Mb(t_*^+)$ satisfies $|\Mb_i|_p=|\Mb_j|_p+|\Mb_k|_p$ and $\rho(\Mb_i)=\rho(\Mb_j)+\rho(\Mb_k)$). By definition of $t_*$, we know that no component of $\Mb(t_*^-)$ is large, so $|\Mb_j|_p\leq\Lambda$ and $\rho(\Mb_j)\leq\Gamma$, and the same for $\Mb_k$. In particular, we know that $\Mb_i$ is large, and $|\Mb_i|_p\leq 2\Lambda_\ell$ and $\rho(\Mb_i)\leq2\Gamma$.

Now, let $\Mb$ be the molecule formed by adding to $\Mb_i$ the empty ends corresponding to each root particle in $[s]$ that does not belong to $\Mb_i$ (these particles move linearly without any collision), with $p(\Mb)=[s]\cup p(\Mb_i)$ and $r(\Mb)=[s]$. By the above property of $\Mb_i$, and using that $s\leq\Lambda_\ell$, it is easy to see that $\Mb$ satisfies (\ref{eq.trc_err}). Moreover the fact that $\Mb_i$ is a component of the C-molecule topological reduction $\Mb$ of the E-dynamics, implies that $\mathbbm{1}_{\Mb_i}=1$ (see the proof of Lemma \ref{lem.1_M_notsim_le_1_M}), and thus $\mathbbm{1}_\Mb=1$. This proves (\ref{eq.error_upper_grammar_rec}).

\textbf{Finishing the proof of (\ref{eq.cumulant_formula_err_trunc}).} Recall from (\ref{eq.W_err}) and (\ref{eq.truncated_error}) that
\begin{equation}f_s^{\mathrm{err}}(t_{\mathrm{fin}},\vz_s)=\varepsilon^{s(d-1)}\sum_{n=0}^\infty\frac{1}{n!}W_{s+n}^{\mathrm{err}}(t_{\mathrm{fin}},\vz_{s+n})=\varepsilon^{s(d-1)}\sum_{\ell=1}^{\mathfrak{L}}\sum_{n=0}^\infty\frac{1}{n!}\big[(\Sc^E_{s+n}(\tau))^{{\mathfrak{L}}-\ell}\circ \Sc_{\ell,s+n}^{\mathrm{err}}(\tau)\big](\widetilde{W}_{s+n}((\ell-1)\tau)).
\end{equation} Now, using (\ref{eq.error_upper_grammar_rec}), the definition $\Ec_N^\ell=\Ec_N^{\Lambda_\ell,\Gamma}$, the definition of $\Sc_{\ell,N}^{\mathrm{err}}$ in (\ref{eq.S^err}), and the $L^1$ preserving property of $\Sc_N^E$ and $\Sc_N^\ell$ (Proposition \ref{prop.dynamics_E_T_property} (\ref{it.dynamics_E_T_property_3.5}), we obtain that
\begin{align}
\|f_s^{\mathrm{err}}(t_{\mathrm{fin}},\vz_s)\|_{L_{\vz_s}^1}&=\varepsilon^{s(d-1)}\sum_{\ell=1}^{\mathfrak{L}}\sum_{n=0}^\infty\frac{1}{n!}\big\|\big[(\Sc^E_{s+n}(\tau))^{{\mathfrak{L}}-\ell}\circ \Sc_{\ell,s+n}^{\mathrm{err}}(\tau)\big](\widetilde{W}_{s+n}((\ell-1)\tau))\big\|_{L^1}
\nonumber\\
&\leq 2\varepsilon^{s(d-1)}\sum_{\ell=1}^{\mathfrak{L}}\sum_{n=0}^\infty\frac{1}{n!}\big\|\mathbbm{1}_{\Ec_{s+n}^{\ell}}\cdot\widetilde{W}_{s+n}((\ell-1)\tau)\big\|_{L^1}
\nonumber\\
\label{eq.proof_upper_err_1}&\leq 2\varepsilon^{s(d-1)}\sum_{\ell=1}^{\mathfrak{L}}\sum_{n=0}^\infty\sum_{{\Mb_{\ell}}\in\Fs_{\Lambda_\ell}^{\mathrm{trc.err}},\,r({\Mb_{\ell}})=[s]}\frac{1}{n!}\big\|\mathbbm{1}_{\Mb_{\ell}}\cdot\widetilde{W}_{s+n}((\ell-1)\tau)\big\|_{L^1}.
\end{align}
In the sum (\ref{eq.proof_upper_err_1}), we argue as in the proof of Proposition \ref{prop.molecule_representation}: by symmetry we may always restrict $p({\Mb_{\ell}})=[|{\Mb_{\ell}}|_p]$ with an additional factor $\binom{n}{|{\Mb_{\ell}}|_p-s}$. Using also (\ref{eq.sum_equiv}) which is also true in the current setting, we obtain that
\begin{align}
\|f_s^{\mathrm{err}}(t_{\mathrm{fin}},\vz_s)\|_{L_{\vz_s}^1}&\leq2\varepsilon^{s(d-1)}\sum_{\ell=1}^{\mathfrak{L}}\sum_{n=0}^\infty\sum_{{\Mb_{\ell}}\in\Fs_{\Lambda_\ell}^{\mathrm{trc.err}},\,r({\Mb_{\ell}})=[s]}\frac{1}{n!}\big\|\mathbbm{1}_{\Mb_{\ell}}\cdot\widetilde{W}_{s+n}((\ell-1)\tau)\big\|_{L^1}\nonumber\\
&\leq2\varepsilon^{s(d-1)}\sum_{\ell=1}^{\mathfrak{L}}\sum_{n=0}^\infty\sum_{\substack{{\Mb_{\ell}}\in\Fs_{\Lambda_\ell}^{\mathrm{trc.err}}\\p({\Mb_{\ell}})=[|{\Mb_{\ell}}|_p],\,r({\Mb_{\ell}})=[s]}}\frac{1}{n!}\binom{n}{|{\Mb_{\ell}}|_p-s}\big\|\mathbbm{1}_{\Mb_{\ell}}\cdot\widetilde{W}_{s+n}((\ell-1)\tau)\big\|_{L^1}\nonumber\\
&=2\varepsilon^{s(d-1)}\sum_{\ell=1}^{\mathfrak{L}}\sum_{n=0}^\infty\sum_{[{\Mb_{\ell}}]\in\Fc_{\Lambda_\ell}^{\mathrm{trc.err}},\,r({\Mb_{\ell}})=[s]}\frac{1}{(n-|{\Mb_{\ell}}|_p+s)!}\big\|\mathbbm{1}_{\Mb_{\ell}}\cdot\widetilde{W}_{s+n}((\ell-1)\tau)\big\|_{L^1}\nonumber\\
&=2\varepsilon^{s(d-1)}\sum_{\ell=1}^{\mathfrak{L}}\sum_{[{\Mb_{\ell}}]\in\Fc_{\Lambda_\ell}^{\mathrm{trc.err}},\,r({\Mb_{\ell}})=[s]}\int_{\Rb^{2d|{\Mb_{\ell}}|_p}}\mathbbm{1}_{\Mb_{\ell}}(\vz_{p({\Mb_{\ell}})})\,\mathrm{d}\vz_{p({\Mb_{\ell}})}\nonumber\\
\label{eq.proof_upper_err_new}&\phantom{2\varepsilon^{s(d-1)}}\times\sum_{n=|{\Mb_{\ell}}|_p-s}^\infty \frac{1}{(n-|{\Mb_{\ell}}|_p+s)!}\int_{\Rb^{2d(n-|{\Mb_{\ell}}|_p+s)}}\widetilde{W}_{s+n}((\ell-1)\tau)\,\mathrm{d}\vz_{[s+n]\backslash p({\Mb_{\ell}})}.
\end{align}
Since the summation in (\ref{eq.proof_upper_err_new}) is equal to $\varepsilon^{-(d-2)|{\Mb_{\ell}}|_p}\widetilde{f}_{|{\Mb_{\ell}}|_p}((\ell-1)\tau,\vz_{p({\Mb_{\ell}})})$ by (\ref{eq.truncated_correlation}), and using also the $L^1$ preserving property of $\Sc_{\Mb_{\ell}}$, we conclude that
\begin{equation}\label{eq.proof_upper_err_new_2}
\begin{aligned}\|f_s^{\mathrm{err}}(t_{\mathrm{fin}},\vz_s)\|_{L_{\vz_s}^1}&\leq 2\sum_{\ell=1}^\Lf\sum_{[{\Mb_{\ell}}]\in\Fc_{\Lambda_\ell}^{\mathrm{trc.err}},\,r({\Mb_{\ell}})=[s]}\varepsilon^{-(d-1)|{\Mb_{\ell}}|_{p\backslash r}}\int_{\Rb^{2d|{\Mb_{\ell}}|_p}}\mathbbm{1}_{\Mb_{\ell}}\cdot \widetilde{f}_{|{\Mb_{\ell}}|_p}((\ell-1)\tau)\,\mathrm{d}\vz_{p({\Mb_{\ell}})}\\
&\leq 2\sum_{\ell=1}^\Lf\sum_{[{\Mb_{\ell}}]\in\Fc_{\Lambda_\ell}^{\mathrm{trc.err}},\,r({\Mb_{\ell}})=[s]}\varepsilon^{-(d-1)|{\Mb_{\ell}}|_{p\backslash r}}\int_{\Rb^{2d|{\Mb_{\ell}}|_p}}|(\Sc\circ\mathbbm{1})_{\Mb_{\ell}}||\widetilde{f}_{|{\Mb_{\ell}}|_p}((\ell-1)\tau)\,\mathrm{d}\vz_{p({\Mb_{\ell}})}.
\end{aligned}
\end{equation} 

The $\widetilde{f}((\ell-1)\tau)$ term in (\ref{eq.proof_upper_err_new_2}) contains a main term and an error term $\mathrm{Err}$ as in (\ref{eq.cumulant_expansion}). The contribution of $\mathrm{Err}$ leads to the second term on the right hand side of (\ref{eq.cumulant_formula_err_trunc}), using $|\Mb_\ell|_{p\backslash r}\leq 3\Lambda_\ell$ and following the proof of \eqref{eq.Err2term} in Part 1 of Proposition \ref{prop.cumulant_single}. For the contribution of the main term in (\ref{eq.cumulant_expansion}), we repeat the same arguments as in the proof of (\ref{eq.cumulant_formula}) and (\ref{eq.cumulant_formula_err}) in Parts 3 and 4 above (i.e. expand $\widetilde{f}((\ell-1)\tau)$ into previous layers and putting all different layers into a multi-layer molecule $\Mb$; note we have extra integration in particles $r(\Mb)=[s]$ here). This then leads to the first term on the right hand side of (\ref{eq.cumulant_formula_err_trunc}), where $[\Mb]\in\Fc_{\vLambda_{\ell}}^{\mathrm{trc.err}}$ because $[\Mb_{\ell}]\in\Fc_{\Lambda_\ell}^{\mathrm{trc.err}}$, see Definition \ref{def.set_T_F} (\ref{it.set_F_err_2}). This proves (\ref{eq.cumulant_formula_err_trunc}), and finally finishes the proof of Proposition \ref{prop.cumulant_formula}.
\end{proof}
\begin{remark}\label{rem.unlabeled} With Proposition \ref{prop.cumulant_formula} proved, in the rest of this paper, we need to estimate the quantities $\Ic\Nc_\Mb$ and $|\Ic\Nc_\Mb|$ on the right hand side of (\ref{eq.fAterm})--(\ref{eq.cumulant_formula_err_trunc}). Note that these quantities are constant in equivalence classes $[\Mb]$. Thus, as said in Remark \ref{rem.label}, for the rest of this paper, we will use $\Mb$ to denote unlabeled molecules only (plus labelings of root particles only).
\end{remark}

\section{Statements of cumulant and $f^\Ac$ estimates}\label{sec.reduction} In this section we will use Proposition \ref{prop.cumulant_formula} to reduce Theorem \ref{th.main} to the 3 main estimates (see Propositions Propositions \ref{prop.est_fa}--\ref{prop.cumulant_error} below), about the quantities $f^\Ac$, $E_H$, $\mathrm{Err}^1$, $\mathrm{Err}^2$ and $f^{\mathrm{err}}$ defined in Proposition \ref{prop.cumulant_formula}. These propositions will be the main goal in the rest of the paper, whose proof will occupy Sections \ref{sec.local_associated}--\ref{sec.error}.

In the first proposition, we show that $f^\Ac$ is indeed an approximate solution to the Boltzmann equation. 
\begin{proposition}
\label{prop.est_fa} Let $f^\Ac$ be defined as in Proposition \ref{prop.cumulant_formula} \eqref{it.cumulant_formula_1}, inductively via (\ref{eq.fAterm}), with $f^\Ac(0,z)=f_0(z)$; also recall $\beta_\ell$ and $\theta_\ell$ defined in (\ref{eq.decayseq}) and $\|\cdot\|_{\mathrm{Bol}^{\beta_\ell}}$ defined by \eqref{eq.boltzmann_decay_2}. Then we have
\begin{equation}\label{eq.approxfA_1}
    f^\Ac(\ell\tau,z)=f(\ell\tau,z)+(f^\Ac)_{\mathrm{rem}}(\ell\tau,z),
\end{equation}
 where $f(t,z)$ is the solution to the Boltzmann equation as in Theorem \ref{th.main}, and $(f^\Ac)_{\mathrm{rem}}$ satisfies that
\begin{equation}\label{eq.approxfA_2}
    \|(f^\Ac)_{\mathrm{rem}}(\ell\tau)\|_{\mathrm{Bol}^{\beta_\ell}}\leq\varepsilon^{\theta_\ell}.
\end{equation}
\end{proposition}
\begin{proof} The proof will occupy Section \ref{sec.fa}, and will be finished in Section \ref{sec.fa_subleading}.
\end{proof}

In the second proposition, we prove upper bounds for $E_H$ and $\mathrm{Err}_{\ell}^1$. This is the crucial proposition, and will be the main goal of the second part of the paper.
\begin{proposition}
\label{prop.cumulant_est} Recall the $E_H(\ell\tau,\vz_H)$ and $\mathrm{Err}^1(\ell\tau,\vz_s)$ introduced in Proposition \ref{prop.cumulant_formula}, and also the notion $C_j^*$ in (\ref{eq.defC*}) in Definition \ref{def.notation} (\ref{it.defC*}). Assume that the estimates (\ref{eq.approxfA_1})--(\ref{eq.approxfA_2}) hold with $\ell$ replaced by any $\ell'\leq\ell-1$. Then we have
\begin{equation}\label{eq.cumulant_est_1}
\|E_H(\ell\tau,\vz_H)\|_{L^1}\leq\varepsilon^{3^{-d-2}+(C_{14}^*)^{-1}\cdot |H|}
\end{equation} for any $H\neq\varnothing$, and
\begin{equation}\label{eq.cumulant_est_2}
\|\mathrm{Err}^1(\ell\tau,\vz_s)\|_{L^1}\leq \tau^{\Lambda_{\ell}/10}.
\end{equation}
\end{proposition}
\begin{proof}
    The proof will occupy Sections \ref{sec.local_associated}--\ref{sec.treat_integral}, and will be finished in Section \ref{sec.summary}, pending a crucial combinatorial statement, Proposition \ref{prop.comb_est}. Then, the proof of Proposition \ref{prop.comb_est} will be the main goal of the third part of the paper, and will occupy Sections \ref{sec.prepare}--\ref{sec.maincr}.
\end{proof}

In the last proposition we prove upper bounds for $f_s^{\mathrm{err}}(t_{\mathrm{fin}})$.
\begin{proposition}
\label{prop.cumulant_error} Recall the $f_s^{\mathrm{err}}(t_{\mathrm{fin}},\vz_s)$ introduced in Proposition \ref{prop.cumulant_formula}. Then we have
\begin{equation}\label{eq.cumulant_error}
    \|f_s^{\mathrm{err}}(t_{\mathrm{fin}},\vz_s)\|_{L^1}\leq \varepsilon^{3^{-d-2}}.
\end{equation}
\end{proposition}
\begin{proof} The proof will occupy Section \ref{sec.error}, based on Proposition \ref{prop.comb_est_extra} which will be proved in Section \ref{sec.last_proof}.
\end{proof}

Now we prove Theorem \ref{th.main} assuming Propositions \ref{prop.est_fa}--\ref{prop.cumulant_error}.
\begin{proof}[Proof of Theorem \ref{th.main} assuming Propositions \ref{prop.est_fa}--\ref{prop.cumulant_error}] The goal is to prove \eqref{eq.maineqn} for $f_s(t_{\mathrm{fin}},\vz_s)$ (recall we have set $t=t_{\mathrm{fin}}=\Lf\tau$), with $s\leq|\log\varepsilon|=A_\Lf$. By \eqref{eq.f_s_dec} and \eqref{eq.cumulant_error}, we may replace $f_s(t_{\mathrm{fin}},\vz_s)$ in \eqref{eq.maineqn} by $\widetilde{f}_s(t_{\mathrm{fin}},\vz_s)$ in (\ref{eq.cumulant_expansion}), which is 
\begin{equation}\label{eq.cumulant_est_3}
\widetilde{f}_s(t_{\mathrm{fin}},\vz_s)=\sum_{H\subseteq[s]}(f^\Ac(t_{\mathrm{fin}}))^{\otimes([s]\backslash H)}\cdot E_H(t_{\mathrm{fin}},\vz_H)+\mathrm{Err}(\Lf\tau,\vz_s).
\end{equation} 
By iteratively applying the $\mathrm{Err}^2$ estimate \eqref{eq.cumulant_formula_err2} in Proposition \ref{prop.cumulant_formula} and combining with \eqref{eq.cumulant_est_2}, we get that
\[
    \|\mathrm{Err}(\Lf\tau,\vz_s)\|_{L^1}\leq\sum_{\ell=0}^\Lf\tau^{\Lambda_{\ell}/10}\cdot \prod_{\ell'=\ell+1}^\Lf\varepsilon^{-(d-1)\Lambda_{\ell'}^2A_{\ell'}}
\] 
which is negligible thanks to our choice of the sequences $(\Lambda_\ell)$ and $(A_\ell)$ in (\ref{eq.ineq_Al}) in Definition \ref{def.notation} (\ref{it.defnot3}), which ensures that $\Lambda_{\ell} \gg \Lambda_{\ell'}^2A_{\ell'}$ if $\ell'>\ell$ and $\tau^{\Lambda_{\ell}/10}\cdot \varepsilon^{-(d-1)\Lambda_{\ell'}^2A_{\ell'}}\ll 1$. Therefore, we may omit the $\mathrm{Err}$ term from (\ref{eq.cumulant_est_3}).

Next, (\ref{eq.approxfA_1}) and (\ref{eq.approxfA_2}) implies that $\|f^\Ac(t_{\mathrm{fin}})-f(t_{\mathrm{fin}})\|_{L^1}\leq \varepsilon^{(\theta+\theta_\Lf)/2}$ (with $\theta=\theta_0/2$ and $\theta_\Lf$ defined in Definition \ref{def.notation} (\ref{it.defnot3})), and thus $\|f^\Ac(t_{\mathrm{fin}})\|_{L^1}\leq 1+\varepsilon^\theta$ using $\int f_0(z)\,\mathrm{d}z=1$ and the conservation of $L^1$ norm for the Boltzmann equation (\ref{eq.boltzmann}). Using also the bound (\ref{eq.cumulant_est_1}) for $E_H$ and the fact that $s\leq|\log\varepsilon|$, it follows that the terms corresponding to $H\neq\varnothing$ in (\ref{eq.cumulant_est_3}) are also negligible in $L^1$. We then know that the function $\widetilde{f}_s(t_{\mathrm{fin}})$ is close to the following tensor product up to error $O_{L^1}(\varepsilon^\theta)$:
\[
    (f^\Ac(t_{\mathrm{fin}}))^{\otimes [s]}=\prod_{p\in [s]}f^\Ac(t_{\mathrm{fin}},z_p).
\] 
Once again, by $\|f^\Ac(t_{\mathrm{fin}})-f(t_{\mathrm{fin}})\|_{L^1}\leq \varepsilon^{(\theta+\theta_\Lf)/2}$, we know the function $\widetilde{f}_s(t_{\mathrm{fin}})$ is close to $(f(t_{\mathrm{fin}}))^{\otimes [s]}$ up to error $O_{L^1}(\varepsilon^\theta)$. Then \eqref{eq.maineqn} follows by using the following estimates for the indicator function $\mathbbm{1}_{\Dc_s}(\vz_s)$
\[
\bigg\|\prod_{p\in[s]}f(t_{\mathrm{fin}},z_p)\cdot\big(1-\mathbbm{1}_{\Dc_s}(\vz_s)\big)\bigg\|_{L^1}\leq Cs^2\varepsilon^d\leq\varepsilon^{d/2}
\] 
because $s\leq|\log\varepsilon|$. This proves Theorem \ref{th.main}.
\end{proof}

\section{Preparations for $E_H$ estimates} \label{sec.local_associated}  From this section until Section \ref{sec.maincr}, we will prove Proposition \ref{prop.cumulant_est} (and the associated Proposition \ref{prop.comb_est}). We start from the upper bounds of $E_H$ and $\mathrm{Err}^1$ in (\ref{eq.cumulant_formula})--(\ref{eq.cumulant_formula_err}) in Proposition \ref{prop.cumulant_formula}, which have the form $\sum_{\Mb}|\Ic\Nc_\Mb|$ (here and below we will replace $[\Mb]$ by $\Mb$ as $\Mb$ now represents unlabeled molecules, see Remark \ref{rem.unlabeled}). To bound these terms, it suffices to (i) control the number of terms $\Mb$ in the summation, and (ii) control the size of $|\Ic\Nc_\Mb|$ for each individual $\Mb$.

In this section, first in Section \ref{sec.int_dia_exp} we prove an upper bound (Proposition \ref{prop.layerrec3}) on the number of choices of $\Mb$, and then in Section \ref{sec.int_reduce} we reduce the $L^1$ norm of $|\Ic\Nc_\Mb|$ to an integral $\Ic_\Mb(Q_\Mb)$ involving Dirac $\dirac$ functions that localizes each individual collision and overlap (Proposition \ref{prop.local_int}). The estimates of this quantity $\Ic_\Mb(Q_\Mb)$ will then occupy Sections \ref{sec.cutting}--\ref{sec.treat_integral}.

\subsection{Upper bounds of combinatorial factors}\label{sec.int_dia_exp} In this subsection we prove an upper bound on the number of molecules $\Mb$ in Proposition \ref{prop.layerrec3}. For this we need to introduce the following parameter $\rho$ for molecules $\Mb$:
\begin{definition}[The parameter $\rho$]\label{def.parameter_rho_old} Given a molecule $\Mb$ with highest layer $\ell$ and lowest layer $1$, with layer $\ell'$ denoted by $\Mb_{\ell'}$. Define the parameters $s_{\ell'}$ and $\Rf_{\ell'}$ by $s_{\ell'} \coloneqq |r(\Mb_{\ell'})|$ and $\Rf_{\ell'} \coloneqq \rho(\Mb_{\ell'})$ for $\ell'\in[1:\ell]$; also define $s_0=|H_0|$ where $H_0$ is the initial cumulant set (see (\ref{eq.proof_cumulant_multi_3_2})), i.e. those particles involved in initial links. Now define the parameter $\rho$ by
\begin{equation}\label{eq.def_rho_old}
    \rho = \sum_{\ell'=0}^\ell s_{\ell'}+\sum_{\ell'=1}^\ell\Rf_{\ell'}.
\end{equation}
We interpret the quantity $\rho$ as the sum of $\Rf_{\ell'}$ (the total number of recollisions within each single layer), together with the sum of $s_{\ell'}$ (the total number of particles lines crossing different layers).
\end{definition}

We now state and prove the main combinatorial upper bound in this section.
\begin{proposition}\label{prop.layerrec3} Let $\Mb$ be a molecule and $\rho$ be as in Definition \ref{def.parameter_rho_old}. Assume $H$ is a set and $|H|\le |\log\varepsilon|^{C^*}$. Consider all molecules $\Mb\in \mathcal{F}_{\boldsymbol{\Lambda}_\ell}$ (recall that we are summing over non-labeled molecules $\Mb$) that satisfy $r(\Mb)=H$, with the size $|\Mb|$ of $\Mb$ fixed, and the $\rho$ defined in Definition \ref{def.parameter_rho_old} also fixed. Then the number of choices for $\Mb$, together with the labelings of root particles, is bounded by
\begin{equation}\label{eq.layercomb}
    \sum_{\Mb} 1\leq C^{|\Mb|}\cdot |\log\varepsilon|^{C^*\rho}.
\end{equation} The same result is true if we replace $\mathcal{F}_{\boldsymbol{\Lambda}_\ell}$ by $\mathcal{F}_{\boldsymbol{\Lambda}_\ell}^{\mathrm{err}}$ or $\Tc_{\Lambda_\ell}$ or $\Fc_{\vLambda_\ell}^{\mathrm{trc.err}}$.
\end{proposition}
\begin{proof} We only consider the case $\Mb\in \mathcal{F}_{\boldsymbol{\Lambda}_\ell}$, as the other cases are similar. Note that $\Mb$ is constructed from its layers $\Mb_{\ell'}$ by the following process, see Part 3 of the proof of Proposition \ref{prop.cumulant_formula} and {\color{blue}Figure \ref{fig.mergelayer}}:
\begin{itemize}
\item Label the root particles, and draw a molecule $\Mb_{\ell'}$ with a single layer $\ell'$, for each $\ell'\in[1:\ell]$;
\item For each $\ell'\in[1:\ell-1]$, identify a subset of $p(\Mb_{\ell'})$ with a subset of $p(\Mb_{\ell'+1})$, call this identified subset $H_{\ell'}$ and $|H_{\ell'}|=s_{\ell'}$;
\item For each particle line in $H_{\ell'}$, merge the corresponding top end in $\Mb_{\ell'}$ and bottom end in $\Mb_{\ell'+1}$ to form a bond;
\item Finally choose a set $H_0\subseteq p(\Mb_1)$, and fix the initial links $\Lc$ between particle lines in $H_0$ as in Proposition \ref{prop.initial_cumulant}.
\end{itemize}

We now estimate the number of choices for $\Mb$, by considering each of the above steps.

\textbf{Proof part 1}. In this part we prove that $|\Mb|\leq|\log\varepsilon|^{C^*}$ and $|\Mb|_p\leq |\log\varepsilon|^{C^*}$. In fact, by Definition \ref{def.set_T_F} (\ref{it.set_F_l_1})--(\ref{it.set_F_l_2}) we know that each cluster of $\Mb_{\ell'}$ has at most $\Lambda_{\ell'}$ particle lines and at most $\Lambda_{\ell'}+\Gamma$ collisions (using the upper bound of $\Gamma$ for the recollision number). Then each component has at most $\Lambda_{\ell'}^2$ particle lines and at most $2\Lambda_{\ell'}^2$ collisions (as each component has at most $\Lambda_\ell$ clusters). Then each $\Mb_{\ell'}$ has at most $A_{\ell'}\Lambda_{\ell'}^2$ particle lines and $2A_{\ell'}\Lambda_{\ell'}^2$ collisions (as each component must contain a root particle line, and there are $s_{\ell'}\leq A_{\ell'}$ of them). All these values are bounded as $\leq|\log\varepsilon|^{C^*}$ by Definition \ref{def.notation} (\ref{it.defnot3}), hence the result.

\textbf{Proof part 2.} In this part we consider the first step of labeling and constructing individual $\Mb_{\ell'}$. Note that labelings lead to $|H|!\leq|\log\varepsilon|^{C^*\rho}$ choices as $|H|=|r(\Mb)|=|r(\Mb_\ell)|=s_\ell\leq\rho$. Next we may fix the values of $(|\Mb_{\ell'}|)_{\ell'\in[1:\ell]}$. As these numbers add up to the fixed $|\Mb|$, we have at most $2^{|\Mb|}$ choices which is acceptable in (\ref{eq.layercomb}). Similarly, we can also fix the values of $(s_{\ell'})$ and $(\rho(\Mb_{\ell'}))$ at an acceptable price. Now, note that each $\Mb_{\ell'}$ is a graph with each atom having at most $4$ bonds; moreover it is formed by adding at most $\rho(\Mb_{\ell'})$ edges to a forest. 

It is known that the number of trees (and forests) with $|\Mb_{\ell'}|$ atoms and each atom having at most $4$ bonds can be bounded by $10^{|\Mb_{\ell'}|}$ (the result for forests follows from that for trees; for trees, we may fix a root and turn the tree into a rooted, at most ternary tree, and use the well-known counting upper bound for the number of binary and ternary trees).

Once the tree parts are fixed, we add in the $\leq\rho$ extra edges one by one; each edge has at most $|\Mb_{\ell'}|^2\leq |\log\varepsilon|^{C^*}$ choices, which leads to $|\log\varepsilon|^{C^*\rho(\Mb_{\ell'})}$ choices in total. All these factors are acceptable in (\ref{eq.layercomb}) even after accounting for different $\ell'$ (for example $\prod_{\ell'} 10^{|\Mb_{\ell'}|}=10^{|\Mb|}$).

\textbf{Proof part 3.} In this part we consider the steps of matching and merging particle lines/ends in $\Mb_{\ell'}$ and $\Mb_{\ell'+1}$, see {\color{blue}Figure \ref{fig.mergelayer}}. Recall that $|H_{\ell'}|=s_{\ell'}$; selecting a set of this cardinality in $\Mb_{\ell'}$ and $\Mb_{\ell'+1}$ leads to at most \[\binom{|\Mb_{\ell'}|_p}{s_{\ell'}}\binom{|\Mb_{\ell'+1}|_p}{s_{\ell'}}\leq |\log\varepsilon|^{C^*s_{\ell'}}\] choices. Putting different $\ell'$ together, we again get a factor of $|\log\varepsilon|^{C^*\rho}$ which is acceptable.

\textbf{Proof part 4.} In this part we consider the last step of constructing $H_0$ and $\Lc$. Note that $|H_0|=s_0$, and $|\Lc|$ (i.e. the number of initial links in $\Lc$) is bounded by $2|H_0|$ as $\Lc$ is a forest. The number of choices of $H_0$ and $\Lc$ is then trivially bounded by $|\log\varepsilon|^{C^*s_0}$, which is also acceptable.

By putting Parts 1--4 above together that accounts for all steps in the construction of $\Mb$, we have proved (\ref{eq.layercomb}) and Proposition \ref{prop.layerrec3}.
\end{proof}

\subsection{Reducing to the integral $\Ic_\Mb(Q_\Mb)$}\label{sec.int_reduce}

In this subsection, we rewrite the integrals $|\Ic\Nc_\Mb|$ in Proposition \ref{prop.cumulant_formula} as integrals $\Ic_\Mb(Q_\Mb)$ involving Dirac functions that localize each individual collision and overlap, see Definition \ref{def.associated_op} and Proposition \ref{prop.local_int}. As we will see in Section \ref{sec.cutting}, this representation behaves well under the cutting operation, and allows us to obtain efficient estimates.

\begin{definition}[Associated operator of a molecule]\label{def.associated_op} Consider the molecule $\Mb=(\Mc,\Ec,\Pc,\Lc)$, we define the following objects.
\begin{enumerate}
    \item\label{it.associated_vars} \emph{Associated variables $\vz_\Ec$ and $\vt_\Mb$.} We associate each atom $\nf\in \Mb$ with a time variable $t_\nf$ and each edge $e\in \Ec$ with a position-velocity vector $z_e=(x_e,v_e)$. We refer to them as \newterm{associated variables} and denote the collections of these variables by $\vz_\Ec=(z_e:e\in\Ec)$ and $\vt_\Mb=(t_\nf:\nf\in\Mb)$.
    \item\label{it.associated_dist} \emph{Associated distributions $\boldsymbol{\Delta}_\nf$.} Given an atom $\nf\in\Mb$, let $(e_1,e_2)$ and $(e_1',e_2')$ be the two bottom edges and top edges at $\nf$ respectively, such that $e_1$ and $e_1'$ (resp. $e_2$ and $e_2'$) are serial. We define the \newterm{associated distribution} $\Dirac_\nf=\boldsymbol{\Delta}=\boldsymbol{\Delta}(z_{e_1},z_{e_2},z_{e_1'},z_{e_2'},t_\nf)$ as follows. 
    \begin{enumerate}
        \item If $\nf$ is an C-atom, then
        \begin{equation}\label{eq.associated_dist_C}
        \begin{aligned}
            \boldsymbol{\Delta}(z_{e_1},z_{e_2},z_{e_1'},z_{e_2'},t_\nf)&:=\boldsymbol{\delta}\big(x_{e_1'}-x_{e_1}+t_\nf(v_{e_1'}-v_{e_1})\big)\cdot \boldsymbol{\delta}\big(x_{e_2'}-x_{e_2}+t_\nf(v_{e_2'}-v_{e_2})\big)\\&\times \boldsymbol{\delta}\big(|x_{e_1}-x_{e_2}+t_\nf(v_{e_1}-v_{e_2})|-\varepsilon\big)\cdot \big[(v_{e_1}-v_{e_2})\cdot \omega\big]_-\\&\times\boldsymbol{\delta}\big(v_{e_1'}-v_{e_1}+[(v_{e_1}-v_{e_2})\cdot\omega]\omega\big)\cdot\boldsymbol{\delta}\big(v_{e_2'}-v_{e_2}-[(v_{e_1}-v_{e_2})\cdot\omega]\omega\big),
        \end{aligned}
        \end{equation} 
        where $\omega:=\varepsilon^{-1}(x_{e_1}-x_{e_2}+t_\nf(v_{e_1}-v_{e_2}))$ is a unit vector, and $z_-=-\min(z,0)$.
        \item If $\nf$ is an O-atom, then
        \begin{equation}\label{eq.associated_dist_O}
        \begin{aligned}
            \boldsymbol{\Delta}(z_{e_1},z_{e_2},z_{e_1'},z_{e_2'},t_\nf)&:=\boldsymbol{\delta}(x_{e_1'}-x_{e_1})\cdot \boldsymbol{\delta}(x_{e_2'}-x_{e_2})\\&\times \boldsymbol{\delta}\big(|x_{e_1}-x_{e_2}+t_\nf(v_{e_1}-v_{e_2})|-\varepsilon\big)\cdot \big[(v_{e_1}-v_{e_2})\cdot \omega\big]_-\\&\times\boldsymbol{\delta}\big(v_{e_1'}-v_{e_1}\big)\cdot\boldsymbol{\delta}\big(v_{e_2'}-v_{e_2}\big)
        \end{aligned}
        \end{equation} where $\omega$ is as above. Note that (\ref{eq.associated_dist_O}) forces $z_{e_1}=z_{e_1'}$ and $z_{e_2}=z_{e_2'}$.
    \end{enumerate}
    \item\label{it.associated_domain} \emph{The associated domain $\mathcal{D}$.} The associated domain $\mathcal{D}$ is defined by
    \begin{equation}\label{eq.associated_domain}
        \Dc:=\big\{\vt_\Mb=(t_\nf)\in (\mathbb{R}^+)^{|\Mb|}: (\ell'-1)\tau<t_\nf<\ell'\tau\mathrm{\ if\ }\nf\in\Mb_{\ell'};\ t_{\nf}<t_{\nf^+}\mathrm{\ if\ }\nf^+\mathrm{\ is\ parent\ of\ }\nf \big\}.
    \end{equation}
    \item\label{it.associated_int} \emph{The associated operator $\Ic_\Mb$.}  With the associated variables $\vz_\Ec=(z_e:e\in\Ec)$ and $\vt_\Mb=(t_\nf:\nf\in\Mb)$, we define the \newterm{associated integral operator} $\Ic_\Mb$ by 
    \begin{equation}\label{eq.associated_int_op}
        \begin{aligned}
            \Ic_\Mb(Q):=\varepsilon^{-(d-1)(|\Ec_*|-2|\Mb|)}\int_{\mathbb{R}^{2d|\Ec_*|}\times \Rb^{|\Mb|}} 
            \prod_{\nf\in\Mb}\boldsymbol{\Delta}_\nf \ Q(\vz_\Ec,\vt_\Mb) \,\mathrm{d}\vz_{\Ec_*}\,\mathrm{d}\vt_\Mb. 
        \end{aligned}
    \end{equation}   
    Here $\Ec_*$ is the set of edges of $\Mb$ that are \emph{bonds or free ends} (i.e. not counting fixed ends; for full molecules $\Mb$ this makes no difference). Note that in (\ref{eq.associated_int_op}) we are only integrating in the variables $\vz_{\Ec_*}$ for bonds and free ends, thus $\Ic_\Mb(Q)$ is a function (or distribution) of \emph{all fixed ends of $\Mb$}.
\end{enumerate}
\end{definition}
\begin{remark}\label{rem.delta_int} We make a few remarks concerning Definition \ref{def.associated_op}.
\begin{enumerate}
\item The time variables $t_\nf$ represent the collision or overlap time of $\nf$. The vector variables $z_e = (x_e,v_e)$ are associated with the linear piece of trajectory along edges $e$, where $v_e$ is the velocity (which is constant along this piece), and $x_e$ is the position transported backwards to time $0$ (so the actual position of the particle at time
$t$ in this linear piece is $x_e+tv_e$).
\item The Dirac functions in $\Dirac_\nf$ prescribes the collision or overlap by forcing the elastic collision formula
in (\ref{eq.associated_dist_C}) and forcing the same velocity in (\ref{eq.associated_dist_O}). In (\ref{eq.associated_int_op}), we treat all the variables $z_e$ for \emph{fixed ends} $e$ as fixed, so the result of this integral is a function (or distribution) or these fixed variables.
\item In (\ref{eq.associated_domain}) we have the natural (partial) time ordering, consistent with the atom ordering in Definition \ref{def.molecule_order}: the time $t_\nf$ for $\nf$ must be before the time $t_{\nf^+}$ for $\nf^+$ (its parent, or equivalently the next
collision). Moreover collisions in layer $\Mb_{\ell'}$ must happen in time interval $[(\ell'-1)\tau,\ell'\tau]$,which leads to natural orderings between different layers. These time orderings will occasionally be useful in the proof below. Apart from these, we do not assume any total ordering of collision times (except in an auxiliary step in the proof of Proposition \ref{prop.local_int} below).
\end{enumerate}
\end{remark}
We now state the main result of this subsection.
\begin{proposition}\label{prop.local_int} Let $\Mb$ be a full molecule (so $\Ec_*=\Ec$) with $r(\Mb) = H$. Let $|\Ic\Nc_\Mb|$ be the associated integral of $\Mb$ defined in \eqref{eq.associated_integral_molecule_abs}. Then, we have
\begin{equation}\label{eq.local_int}
    \|\, |\Ic\Nc_\Mb|(\vz_{H})\|_{L^1} = \varepsilon^{(d-1)|H|}\cdot\Ic_\Mb(Q_\Mb).
\end{equation}
where $Q_\Mb$ is the function given by 
\begin{equation}\label{eq.Q_M}
    Q_\Mb = \mathbbm{1}_\Dc(\vt_\Mb)\cdot \mathbbm{1}_\Lc^{\varepsilon}(\vz_\Ec)\cdot \prod_{e\in \Ec_{\mathrm{end}}^-}\big|f^\Ac\big((\ell_1[e]-1)\tau,x_e-(\ell_1[e]-1)\tau\cdot v_e,v_e\big)\big|.
\end{equation}
Here in \eqref{eq.Q_M}, $\Ec_{\mathrm{end}}^-$ is the set of bottom ends of $\Mb$ (including empty ends, see Definition \ref{def.sets_molecule} \eqref{it.subsets_Ec}), and $\ell_1[e]$ is the initial layer for bottom end $e$ in Definition \ref{def.molecule} \eqref{it.particle_line_same}. The indicator functions $\mathbbm{1}_\Dc$ is defined in \eqref{eq.associated_domain}, and $\mathbbm{1}_\Lc^{\varepsilon}$ is defined in \eqref{eq.1_L}.
\end{proposition}
\begin{proof} We divide the proof into 3 parts.

\textbf{Proof part 1.} In this part, we make some simplifications to (\ref{eq.local_int}) and reduce it to (\ref{eq.proof_local_int_1_*}). First we prove $|\Mb|_p = |\Ec| - 2|\Mb|$. This is because each atom in $\Mb$ has exactly $4$ edges; by adding them together and noticing that each bond (between two atoms) is counted twice and each end is counted once, we obtain $4|\Mb| = 2|\Ec| - |\Ec_{\mathrm{end}}^-| - |\Ec_{\mathrm{end}}^+|$ (ignoring empty ends; however this equality remains true even with empty ends). Using also that $|\Ec_{\mathrm{end}}^-| = |\Ec_{\mathrm{end}}^+|=|\Mb|_p$ (as each particle line has a unique top end and a unique bottom end, see Definition \ref{def.molecule} (\ref{it.particle_line_same})), we conclude that $|\Mb|_p = |\Ec| - 2|\Mb|$, as desired.

By $|\Mb|_p = |\Ec| - 2|\Mb|$, and note that $|H|=|\Mb|_r$, we see that the powers of $\varepsilon$ on both sides of (\ref{eq.local_int}) match. Next, in the definition of $|\Ic\Nc_\Mb|$ in (\ref{eq.associated_integral_molecule_abs}), by integrating in the $\vz_{r(\Mb)}$ variables, using the definition $|(\Sc\circ\mathbbm{1})_\Mb|$ in \eqref{eq.associated_op} and the $L^1$ preserving property of $\Sc_\Mb$, we see that
\begin{equation}\label{eq.proof_local_int_*}
\big\||\Ic\Nc_\Mb|\big\|_{L^1}=\varepsilon^{-(d-1)|\Mb|_{p\backslash r}}\int_{\Rb^{2d|\Mb|_p}}\mathbbm{1}_\Mb\cdot|(|Q(\vz_{p(\Mb)}')|)\,\mathrm{d}\vz_{p(\Mb)}',
\end{equation} where $Q$ is as in (\ref{eq.extraQ}). We may identify a particle line $\pb\in p(\Mb)$ with its bottom end $e\in\Ec_{\mathrm{end}}^-$, and note also that $Q_\Mb=|Q|\cdot\mathbbm{1}_\Dc(\vt_\Mb)$; as such, we see that (\ref{eq.local_int}) would follow if we can prove \begin{equation}\label{eq.proof_local_int_1_*}
\int_{\Rb^{2d|\Ec\backslash \Ec_{\mathrm{end}}^-|}\times \Dc}\prod_{\nf\in\Mb}\boldsymbol{\Delta}(z_{e_1},z_{e_2},z_{e_1'},z_{e_2'},t_\nf)\,\mathrm{d}\vz_{\Ec\backslash\Ec_{\mathrm{end}}^-}\mathrm{d}\vt_\Mb= \mathbbm{1}_{\Mb}\big(\vz_{\Ec_{\mathrm{end}}^-}\big),\quad a.e.,
\end{equation}
where we identify $\pb\in p(\Mb)$ with the bottom end $e\in\Ec_{\mathrm{end}}^-$ and identify $\vz_\pb'$ in (\ref{eq.proof_local_int_*}) with $z_e$. (Indeed, once we prove (\ref{eq.proof_local_int_1_*}), we simply integrate both sides against $\varepsilon^{-(d-1)|\Mb|_{p\backslash r}}\cdot|Q|$; the right hand side becomes (\ref{eq.proof_local_int_*}), while the left hand side becomes $\varepsilon^{(d-1)|H|}\cdot \Ic_\Mb(Q_\Mb)$ using the above identification, relation between $|Q|$ and $Q_\Mb$, and the restriction $\vt_\Mb\in \Dc$ in (\ref{eq.proof_local_int_1_*}).)

In the rest of this proof, we focus on proving \eqref{eq.proof_local_int_1_*}.

\textbf{Proof part 2.} In this part, we introduce an extension of \eqref{eq.proof_local_int_1_*} for the purpose of induction. This extension contains two extra parameters $\ell_{\mathrm{min}}$ and $t_{\mathrm{min}}\in [(\ell_{\mathrm{min}}-1)\tau,\ell_{\mathrm{min}}\tau]$, and an (extra) \emph{total ordering} $\prec$ for all atoms in $\Mb$, which is an extension of the partial ordering in Definition \ref{def.molecule_order}. For each such total ordering $\prec$ and each $t_{\mathrm{min}}$, define the extension of \eqref{eq.proof_local_int_1_*} as follows.

In this extension, we require that $\ell[\nf]\geq \ell_{\mathrm{min}}$ for all atoms in $\nf\in\Mb$ (we still view the lowest layer of $\Mb$ as $1$). Moreover, define the new indicator function \begin{equation}\label{eq.proof_local_int_2_1}
    \mathbbm{1}_{\Mb, \prec, t_{\mathrm{min}}}\big(\vz_{\Ec_{\mathrm{end}}^-}\big)
\end{equation}
to be the same indicator function $\mathbbm{1}_\Mb$, but with the additional requirements that (i) each collision and overlap happens after time $t_{\mathrm{min}}$, and (ii) if $\nf_1\prec\nf_2$ in the total ordering, then the collision/overlap corresponding to $\nf_1$ happens before that of $\nf_2$. Similarly, we also define the new domain
\begin{equation}\label{eq.proof_local_int_2_2}
    \Dc_{\mathrm{tot}} \coloneqq \{t_\Mb:\, t_{\mathrm{min}}<t_{\nf_1}<t_{\nf_2}\textrm{ for any }\nf_1\prec \nf_2\}.
\end{equation} Then, the extension is written as \begin{equation}\label{eq.proof_local_int_2_*}
    \int_{\Rb^{2d|\Ec\backslash \Ec_{\mathrm{end}}^-|}\times \Dc_{\mathrm{tot}}}\prod_{\nf\in\Mb}\boldsymbol{\Delta}(z_{e_1},z_{e_2},z_{e_1'},z_{e_2'},t_\nf)\,\mathrm{d}\vz_{\Ec\backslash\Ec_{\mathrm{end}}^-}\mathrm{d}\vt_\Mb= \mathbbm{1}_{\Mb, \prec, t_{\mathrm{min}}}\big(\vz_{\Ec_{\mathrm{end}}^-}\big),\quad a.e.
\end{equation}

If we can prove (\ref{eq.proof_local_int_2_*}) for each fixed $\prec$, then \eqref{eq.proof_local_int_1_*} follows by setting $\ell_{\mathrm{min}}=1$, $t_{\mathrm{min}}=0$ and summing over all different choices of the total ordering $\prec$. In the rest of this proof, we focus on proving (\ref{eq.proof_local_int_2_*}).

\textbf{Proof part 3.} In this part, we prove \eqref{eq.proof_local_int_2_*} by induction on the number of collisions and overlaps. 

If $\Mb$ has no collision or overlap, then $\Mb=\varnothing$ and both sides of \eqref{eq.proof_local_int_2_*} are vacuously $1$. Now suppose \eqref{eq.proof_local_int_2_*} is true for molecules with fewer collisions and overlaps than $\Mb$. For $\Mb$, we choose the smallest atom $\nf\in\Mb$ under the total ordering $\prec$. We will assume $\nf$ is a collision happening between two particles $p_1$ and $p_2$; the overlap case is similar and much easier (where the incoming velocity coincide with the outgoing ones).

Consider the C-atom $\nf\in\Mb$. It has no child as it is the smallest atom, which means it has two bottom ends $(e_1,e_2)$. Let the corresponding serial top edges at $\nf$ be $(e_1',e_2')$, and let $\pb_j$ be the particle line that $e_j$ and $e_j'$ belong to. We will assume $e_1'$ and $e_2'$ are both bonds; the case when one or both of them are ends only requires trivial modifications. Now, if we remove the atom $\nf$ and turn $(e_1',e_2')$ to two bottom ends at the other atoms they connect to, we will reduce $\Mb$ to another molecule $\Mb'$ (with the extra total ordering inherited from $\Mb$).

On the left hand side of \eqref{eq.proof_local_int_2_*}, we may first fix $t_\nf$ and $(x_{e_1'},v_{e_1'},x_{e_2'},v_{e_2'})$, and integrate in the remaining variables. By induction hypothesis for $\Mb'$, the integral in the remaining variables equals 
\begin{equation}\label{eq.defofO'}
    \mathbbm{1}_{\Mb', \prec', t_\nf}\big(x_{e_1'},v_{e_1'},x_{e_2'},v_{e_2'},(z_e)_{e\in\Ec_{\mathrm{end}}^-\backslash \{e_1,e_2\}}\big),
\end{equation} 
where the indicator function $\mathbbm{1}_{\widetilde{\Mb}', \prec', t_\nf}$ is defined for the molecule $\Mb'$ with the total ordering $\prec'$ inherited from $\prec$, and $t_{\mathrm{min}}$ replaced by $t_\nf$. We then integrate in the variables $t_\nf$ and $(x_{e_1'},v_{e_1'},x_{e_2'},v_{e_2'})$; regarding this integral, we claim (for any function $Q\geq 0$) that
\begin{multline}\label{eq.singlecol}
    \int_{\Rb^{4d}\times I}\boldsymbol{\Delta}(z_{e_1},z_{e_2},z_{e_1'},z_{e_2'},t_\nf)\cdot Q(z_{e_1'},z_{e_2'},t_\nf)\,\mathrm{d}z_{e_1'}\mathrm{d}z_{e_2'}\mathrm{d}t_\nf\\=\mathbbm{1}_{\mathrm{col}}(z_{e_1},z_{e_2})\cdot Q\big(x_{\mathrm{out}}^1-t_{\mathrm{col}}\cdot v_{\mathrm{out}}^1,v_{\mathrm{out}}^1,x_{\mathrm{out}}^2-t_{\mathrm{col}}\cdot v_{\mathrm{out}}^2,v_{\mathrm{out}}^2,t_{\mathrm{col}}\big).
\end{multline} Here $I=[(\ell[\nf]-1)\tau,\ell[\nf]\tau]\cap [t_{\mathrm{min}},+\infty)$, $\mathbbm{1}_{\mathrm{col}}$ is the indicator function that a collision happens between particles $p_1$ and $p_2$ with initial position and velocity $(z_{e_1},z_{e_2})$, at collision time $t_{\mathrm{col}}\in I$, and $z_{\mathrm{out}}^j=(x_{\mathrm{out}}^j,v_{\mathrm{out}}^j)$ for $j\in\{1,2\}$ are the position and velocity of the particles $p_1$ and $p_2$ right after the collision; note that $z_{\mathrm{out}}^j$ and $t_{\mathrm{col}}$ are both functions of $(z_{e_1},z_{e_2})$.

The proof of \eqref{eq.singlecol} follows from first integrating in $x_{e_1'}$ and $x_{e_2'}$, then in $v_{e_1'}$ and $v_{e_2'}$, and then in $t_\nf$. Note that:
\begin{enumerate}[{(i)}]
\item The Dirac $\boldsymbol{\delta}$ functions imply that a collision happens between particles $p_1$ and $p_2$;
\item The factor $[(v_{e_1}-v_{e_2})\cdot\omega]_-$ supported in $(v_{e_1}-v_{e_2})\cdot\omega\leq 0$ implies that the particles form a pre-collisional configuration at collision time;
\item The collision time $t_{\mathrm{col}}$ is the unique time $t_{\mathrm{col}}$ satisfying $|(x_{e_1}+t_{\mathrm{col}}v_{e_1})-(x_{e_2}+t_{\mathrm{col}}v_{e_2})|=\varepsilon$ and the pre-collisional condition;
\item The value of $z_{\mathrm{out}}^j=(x_{\mathrm{out}}^j,v_{\mathrm{out}}^j)=(x_{e_j'}+t_{\mathrm{col}}\,v_{e_j'},v_{e_j'})$ equals the outgoing position and velocity of particles $p_j$ immediately after the collision, which is uniquely determined by $(z_{e_1},z_{e_2})$;
\item Once $t_{\mathrm{col}}$ and $v_{e_j'}$ are fixed, the integral in $x_{e_{j}'}$ contains a single $\boldsymbol{\delta}$ factor, and once $t_{\mathrm{col}}$ is fixed, the integral in $v_{e_j'}$ also contains a single $\delta$ factor. Finally it is easy to check that
\[|(v_{e_1}-v_{e_2})\cdot\omega|\cdot\boldsymbol{\delta}\big(|x_{e_1}-x_{e_2}+t_\nf(v_{e_1}-v_{e_2})|-\varepsilon\big)=\boldsymbol{\delta}(t_\nf-t_{\mathrm{col}}),\] so integrating in the order of $x_{e_j'}$ then $v_{e_j'}$ and $t_\nf$ then proves \eqref{eq.singlecol}.
\end{enumerate}

Now, by \eqref{eq.singlecol} and induction hypothesis, we know that the left hand side of \eqref{eq.proof_local_int_2_*} equals
\begin{equation}\label{eq.combinecol}
    \mathbbm{1}_{\mathrm{col}}(z_{e_1},z_{e_2})\cdot \mathbbm{1}_{\Mb', \prec', t_{\mathrm{col}}}\big(x_{\mathrm{out}}^1-t_{\mathrm{col}}\cdot v_{\mathrm{out}}^1,v_{\mathrm{out}}^1,x_{\mathrm{out}}^2-t_{\mathrm{col}}\cdot v_{\mathrm{out}}^2,v_{\mathrm{out}}^2,(z_e)_{e\in\Ec_{\mathrm{end}}^-\backslash \{e_1,e_2\}}\big).
\end{equation} 
The indicator function $\mathbbm{1}_{\mathrm{col}}$ is equivalent to the requirement that the particles $p_1$ and $p_2$ collide in the layer $[(\ell[\nf]-1)\tau,\ell[\nf]\tau]$ and after time $t_{\mathrm{min}}$. The indicator function $\mathbbm{1}_{\Mb', \prec', t_{\mathrm{col}}}$ is equivalent to requiring that, if particle $p_j\,(j\in\{1,2\})$ had initial position $x_{\mathrm{out}}^j-t_{\mathrm{col}}\cdot v_{\mathrm{out}}^j$ and initial velocity $v_{\mathrm{out}}^j$, and other particles had the same initial position and velocity given by the corresponding $z_e$ $(e\in\Ec_{\mathrm{end}}^-\backslash \{e_1,e_2\})$, then the exact collisions and overlaps in the $\Mb'$-prescribed dynamics (with total ordering $\prec'$) would happen, and they would all happen (in the corresponding time layer and) after time $t_{\mathrm{col}}$. 

Note that the position and velocity of $(p_1,p_2)$ at time $t_{\mathrm{col}}$, with the above virtual initial configuration, is exactly $(z_{\mathrm{out}}^1,z_{\mathrm{out}}^2)$, which equals the same outgoing position and velocity of $(p_1,p_2)$ at time $t_{\mathrm{col}}$ right after the first collision $\nf$, with the actual initial configuration $z_e$ $(e\in\Ec_{\mathrm{end}}^-)$. Since we are restricting all collisions and overlaps other than $\nf$ to happen after $t_{\mathrm{col}}$, we know that $\mathbbm{1}_{\Mb', \prec', t_{\mathrm{col}}}$ is equivalent to requiring that the exact collisions and overlaps in the $\Mb$-prescribed dynamics (other than $\nf$, with total ordering $\prec$) all happen (in the corresponding time layer and) after the collision time $t_{\mathrm{col}}$ for $\nf$. Therefore, the product indicator function $\mathbbm{1}_{\mathrm{col}}\cdot\mathbbm{1}_{\Mb', \prec', t_{\mathrm{col}}}$ equals $\mathbbm{1}_{\Mb}$ in \eqref{eq.proof_local_int_2_*}. This proves \eqref{eq.proof_local_int_2_*} and hence Proposition \ref{prop.local_int}.
\end{proof}
\begin{remark}\label{rem.zeroset} By (\ref{eq.local_int}) and (\ref{eq.proof_local_int_*}), we know that changing the value of $Q_\Mb$ on Lebesgue zero sets does not affect the value of $\Ic_\Mb(Q_\Mb)$. Below we will exploit this to make sure that the functions $Q_\Mb$ integrated in \eqref{eq.associated_int_op} are everywhere defined Borel functions, and their support satisfy that no two collisions happen at the same time, and no collision happens at any $\ell\tau$.
\end{remark}
\begin{remark}\label{rem.overlap_property} Note the following difference between C-atoms and O-atoms. In \eqref{eq.associated_dist_C}, the distribution $\Dirac$ is symmetric under permutation $(e_1,e_2,e_1',e_2')\leftrightarrow(e_2,e_1,e_2',e_1')$, and also $(e_1,e_2,e_1',e_2')\leftrightarrow(e_1',e_2',e_1,e_2)$ except for the $[(v_{e_1}-v_{e_2})\cdot\omega]_-$ factor. In other words, for any C-atom we must distinguish top and bottom edges, and also distinguish serial edges, in order for \eqref{eq.associated_dist_C} to be well-defined; however there is also symmetry between ``top" and ``bottom" in the case when the $[(v_{e_1}-v_{e_2})\cdot\omega]_-$ factor does not play a role.

On the other hand, in \eqref{eq.associated_dist_O}, the distribution $\Dirac$ is symmetric under permutation exchanging $e_1$ with $e_1'$, or exchanging $e_2$ with $e_2'$, or exchanging $(e_1,e_1')$ with $(e_2,e_2')$. Therefore, for any O-atom we do not need to distinguish top and bottom edges, as long as the serial edges (or particle lines) are properly distinguished.
\end{remark}

\section{Operations on molecules}\label{sec.cutting} From this section until Section \ref{sec.treat_integral}, we will prove Proposition \ref{prop.cumulant_est} (pending Proposition \ref{prop.comb_est} which will be proved in Sections \ref{sec.prepare}--\ref{sec.maincr}). Thanks to the preparations in Section \ref{sec.local_associated}, we only need to estimate the quantity $\Ic_\Mb(Q_\Mb)$ in (\ref{eq.local_int})--(\ref{eq.Q_M}) for a fixed molecule $\Mb$. 

In Section \ref{sec.oper}, we define several operations on the molecule $\Mb$, which reduce it to smaller and simpler \emph{elementary molecules} (see Definition \ref{def.elementary}). In Section \ref{sec.oper_int} we study the behavior of the quantity $\Ic_\Mb(Q)$ under these operations, see Propositions \ref{prop.deleting}--\ref{prop.cutting}. With these, it then suffices to estimate the $\Ic_\Mb(Q)$ integrals for elementary molecules, which are done in Section \ref{sec.treat_integral}.

\subsection{Cutting, deleting and splitting} \label{sec.oper}

The operations we will define include cutting (see Definition \ref{def.cutting}; this is the main operation that breaks $\Mb$ into two smaller molecules, which corresponds to rewriting the integral $\Ic_\Mb$ in (\ref{eq.associated_int_op}) as an iterated integral), deleting (see Definition \ref{def.delete}; this forgets the overlap condition at an O-atom), and splitting (see Definition \ref{def.splitting}; this restricts to different sub-domains in the integral, or equivalently the support of $Q$ in (\ref{eq.associated_int_op})).

\begin{figure}[h!]
    \centering
    \includegraphics[width=0.35\linewidth]{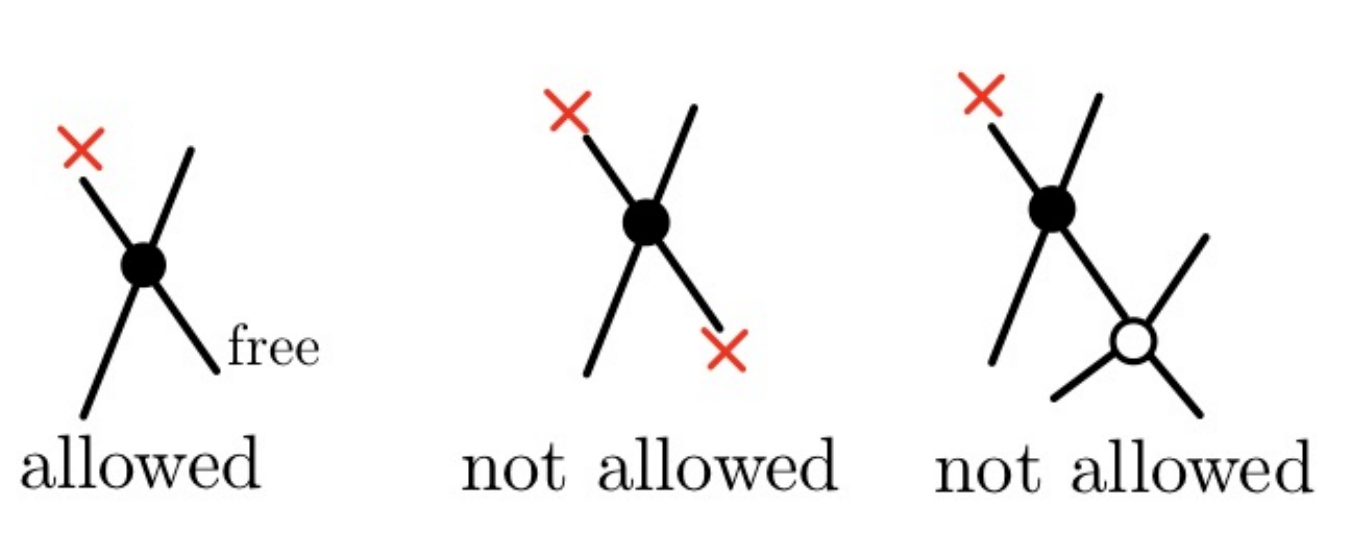}
    \caption{Structures allowed and not allowed at O-atoms in regular molecules. The simple pair is formed by a fixed end and serial free end at an O-atom as shown on the left.}
    \label{fig.regular}
\end{figure}

\begin{figure}[h!]
    \centering
    \includegraphics[width=0.8\linewidth]{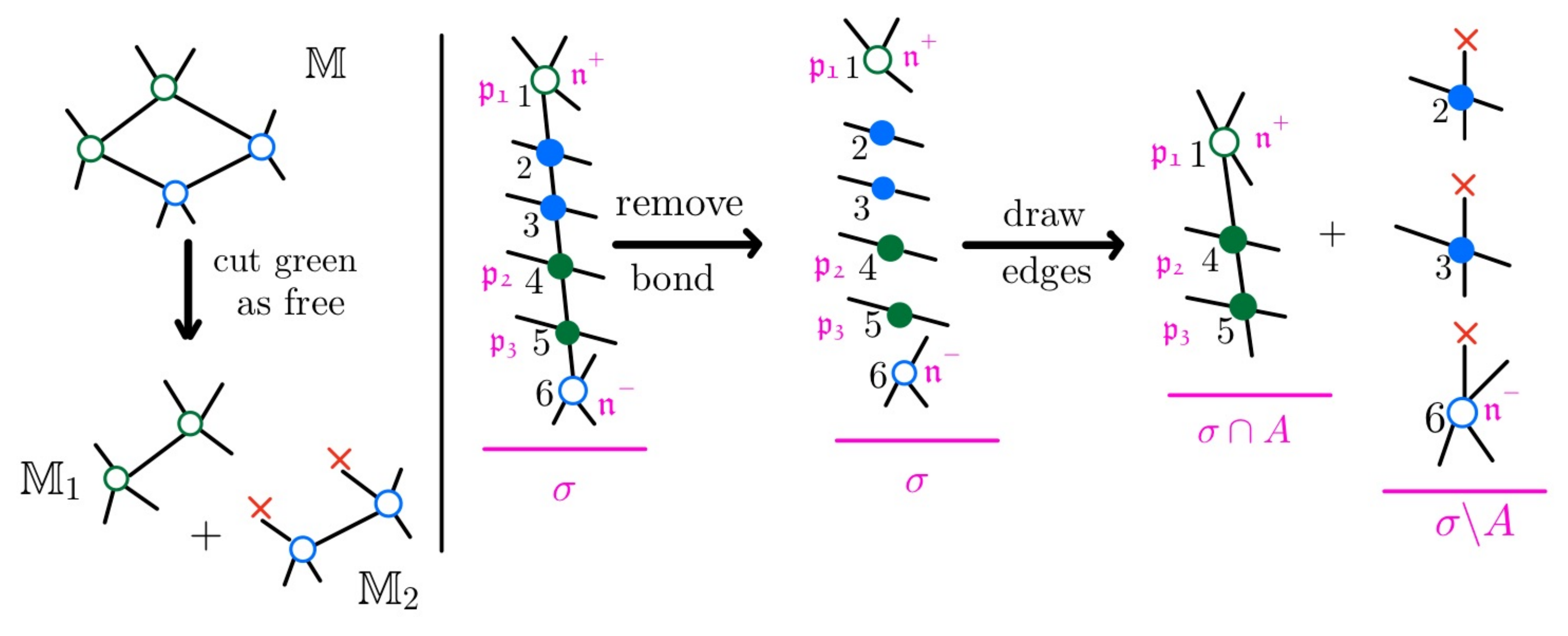}
    \caption{Example of cutting (Definition \ref{def.cutting}). Left: simple case of cutting without O-atoms. Right: cutting in the presence of O-atoms (showing only one maximal ov-segment $\sigma$). Here green atoms belong to $A$ and become $\Mb_1$ after cutting, while blue atoms belong to $\Mb\backslash A$ and become $\Mb_2$ after cutting. The $\nf^\pm$ and $\pf_j$ are as indicated.}
    \label{fig.cutting2}
\end{figure}

The notion of deleting an O-atom is already defined as part of Definition \ref{def.create_delete}. We repeat here for the reader's convenience, together with a few more notions for molecules.
\begin{definition}[Deleting and various other notions]\label{def.delete} We define the following notions for a molecule $\Mb$.
\begin{enumerate}
\item\label{it.delete_new} \emph{Delete an O-atom.} Let $\of$ be an O-atom with 4 edges $e_i\,(1\leq i\leq 4)$ connected to other atoms $\nf_i$, where $e_1$ and $e_3$ are serial and same for $e_2$ and $e_4$. By \textbf{deleting} the O-atom $\of$, we remove $\of$ and the 4 edges $e_j$ from $\Mb$, and add two bonds $e_1'$ between $\nf_1$ and $\nf_3$, and $e_2'$ between $\nf_2$ and $\nf_4$, with the same top/bottom assignment as $(e_1,e_3)$ and $(e_2,e_4)$. Here if some $\nf_i$ is absent, say $e_3$ is a free end at $\of$, then the bond $(e_1,e_3)$ should be replaced by a free end at $\nf_1$, etc.
\item \emph{Degree.} Define the \newterm{degree} of $\nf$ (abbreviated \newterm{deg}) to be the number of edges at $\nf$ that are not fixed ends (counting only bonds and free ends). Note that any atom has deg at most 4, and every atom has deg 4 if and only $\Mb$ is full.
\item \emph{Sub-molecules.} Let $\Mb'\subseteq\Mb$ be a subset of atoms. We will regard $\Mb'$ as a molecule (referred to as a \newterm{sub-molecule} of $\Mb$), by viewing each bond between an atom in $\Mb'$ and an atom in $\Mb\backslash\Mb'$ as a free end at the atom in $\Mb'$.
\end{enumerate}
\end{definition}
Next we define the important notion of regular molecules, which is an important property that governs our treatment of O-atoms and definition of cutting (Definition \ref{def.cutting}) below.
\begin{definition}\label{def.reg} We say $\Mb$ is \textbf{regular}, if for any O-atom $\nf\in\Mb$, any fixed end at $\nf$ must be serial with a free end at $\nf$ (we say they form a \textbf{simple pair}), see {\color{blue}Figure \ref{fig.regular}}. This property will be assumed for all the molecules studied below, and will be preserved throughout the reduction process.
\end{definition}
Now we introduce the operation of cutting. For some intuitions see Remark \ref{rem.cutting}.
\begin{definition}[Cutting]\label{def.cutting} Let $\Mb$ be a regular molecule and $A\subseteq\Mb$ be a set of atoms. The operation of \newterm{cutting} $A$ \newterm{as free} (resp. $\Mb\backslash A$ \newterm{as fixed}) turns $\Mb$ into a new molecule $\Mb'$ which is the disjoint union of two molecules, i.e. $\Mb'=\Mb_1\sqcup\Mb_2$. For convenience, we may also say that $\Mb$ is ``cut into two smaller molecules $\Mb_1$ and $\Mb_2$". We understand that all empty ends belong to $\Mb_1$; as for the atoms and other edges, we describe the procedure as follows, see {\color{blue} Figure \ref{fig.cutting2}}. For some intuitions see also Remark \ref{rem.cutting}.
\begin{enumerate}
\item\label{it.cutting_1} The atom sets of $\Mb_1$ and $\Mb_2$ are $A$ and $\Mb\backslash A$ respectively. If $\Mb$ contains only C-atoms (so each maximal ov-segment is just a bond), then we simply break each bond $e$ between an atom in $A$ and an atom in $\Mb\backslash A$, by removing it and drawing a free end $e_1$ (resp. fixed end $e_2$) at the corresponding atom in $A$ (resp. $\Mb\backslash A$). We form $\Mb_1$ (resp. $\Mb_2$) by atoms in $A$ (resp. $\Mb\backslash A$) and these edges.
\item\label{it.cutting_2} Now, suppose $\Mb$ contains O-atoms. For each maximal ov-segment $\sigma$ that contains an atom in $A$ and an atom in $\Mb\backslash A$ (other maximal ov-segments are not affected), let the highest and lowest atoms of $\sigma$ be $\nf^+$ and $\nf^-$, and let $(\pf_1,\cdots,\pf_s)$ be the atoms in $A\cap\sigma$ (ordered high to low). Then, we take the following steps, which extend the steps in (\ref{it.cutting_1}) from bonds to ov-segments, see {\color{blue}Figure \ref{fig.cutting2}}:
\begin{enumerate}
\item\label{it.cutting_2a} First remove every bond in $\sigma$. Then, for $1\leq j\leq s-1$, draw a bond between $\pf_j$ and $\pf_{j+1}$. Then, if $\pf_1\neq\nf^+$, draw a top free end at $\pf_1$; if $\pf_s\neq\nf^-$, draw a bottom free end at $\pf_s$. Let $\Mb_1$ be formed by atoms in $A$ and these edges (constructed from different $\sigma$).
\item\label{it.cutting_2b} Next, if $\nf^+\in\sigma\backslash A$, draw a bottom fixed end at $\nf^+$; if $\nf^-\in\sigma\backslash A$, draw a top fixed end at $\nf^-$. Then, for any $\nf^\pm\neq \nf\in\sigma\backslash A$, draw a simple pair at $\nf$ consisting of a free and and a (serial) fixed end. Let $\Mb_2$ be formed by atoms in $\Mb\backslash A$ and these edges (constructed from different $\sigma$).
\end{enumerate}
\end{enumerate}
\end{definition}
Note that, deleting O-atoms does not affect the correspondence between molecular and physical pictures under topological reduction. However, this correspondence will break down after cutting. Therefore, we will only talk in the molecule language after performing any cutting. For example, notions like $\pb$ will only be referring to particle lines as in Definition \ref{def.molecule} (\ref{it.particle_line_same}), instead of physical particles.

\begin{figure}[h!]
    \centering
    \includegraphics[width=0.26\linewidth]{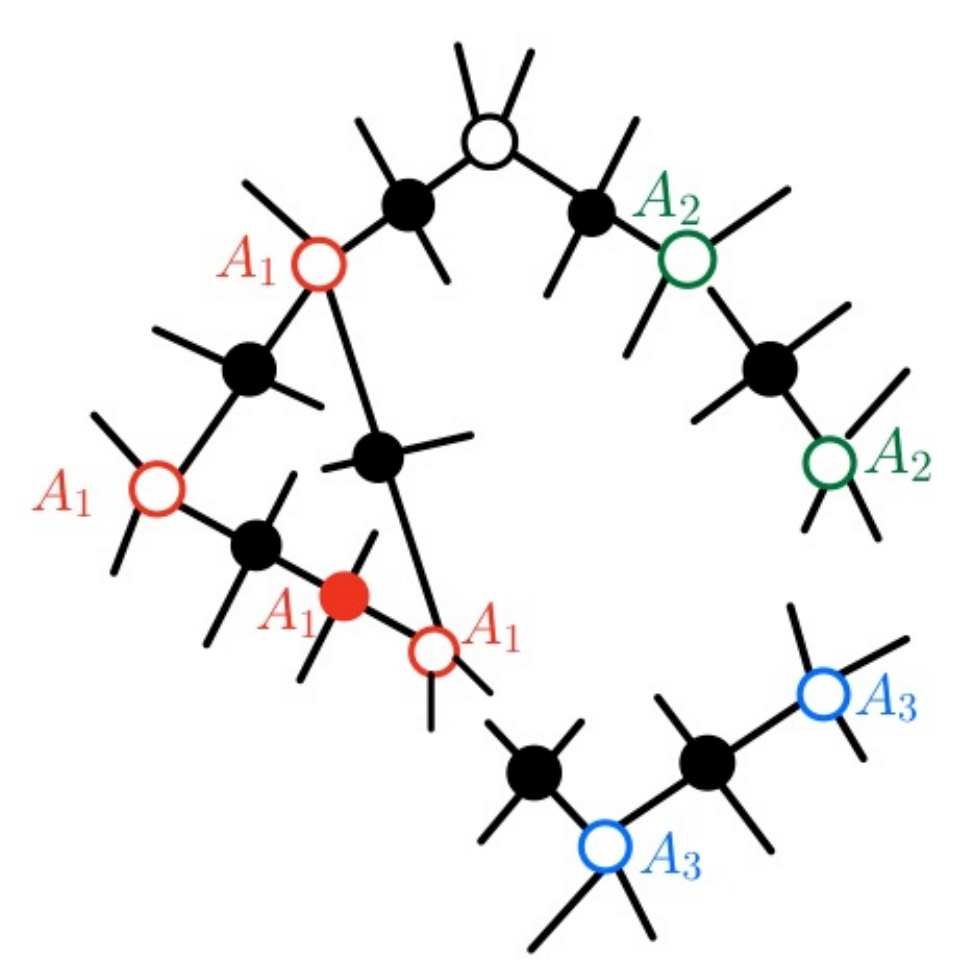}
    \caption{The ov- related notions in Definition \ref{def.ov_connect}. Here $A$ is the set of colored atoms, which has 3 ov-components $A_j$ represented by 3 different colors (as opposed to 7 connected components). The ov-distance between $A_1$ and $A_2$ is $2$, while the ov-distance between $A_1$ and $A_3$ is $\infty$.}
    \label{fig.ov notions}
\end{figure}

In connection with the notion of cutting, we define the notion of ov-adjacency and a few more naturally related notions, see {\color{blue} Figure \ref{fig.ov notions}}:
\begin{definition}[Ov- notions]\label{def.ov_connect} Recall molecules $\Mb$ and ov-segments in Definition \ref{def.molecule}. 
\begin{enumerate}
\item \emph{Ov-adjacency.} We define two atoms $\nf_1$ and $\nf_2$ to be \textbf{ov-adjacent}, if they are connected by an ov-segment; if $\nf_1$ is an ancestor (resp. descendant) of $\nf_2$ along this ov-segment, we say it is an \textbf{ov-parent} (resp. \textbf{ov-child}) of $\nf_2$.
\item \emph{Ov-connectedness.} For a subset $A\subseteq\Mb$, we say it is \textbf{ov-connected} if for any two atoms $\nf,\nf'\in A$ there exist atoms $\nf_j\in \Mb\,(0\leq j\leq q)$ such that $(\nf_0,\nf_q)=(\nf,\nf')$ and $\nf_j$ is ov-adjacent with $\nf_{j+1}$. The maximal ov-connected subsets of $A$ are called \textbf{ov-components} of $A$.
\item\emph{Ov-distance.} Finally, define the \textbf{ov-distance} of two subsets $B,B'\subseteq \Mb$ to be the minimal $q$ such that there exist atoms $\nf_j\in A\,(0\leq j\leq q)$ such that $(\nf_0,\nf_q)\in B\times B'$ and $\nf_j$ is ov-adjacent with $\nf_{j+1}$. Note that the ov-distance is $0$ if $B\cap B'\neq\varnothing$, and is $\infty$ if they are contained in different connected components\footnote{The ov-components of the \emph{whole molecule} $\Mb$ are the same as the components. However the ov-components of \emph{subsets} $A\subseteq\Mb$ might not be the same as its components.} of $\Mb$.
\end{enumerate}  
\end{definition}

\begin{lemma}[]\label{lem.regularity}
    Regularity is preserved after deleting an O-atom or cutting a subset $A$ as fixed or free. 
\end{lemma}
\begin{proof} Simply note that, the forbidden configuration, where a fixed end at an O-atom is serial with another fixed end or a bond, cannot be constructed in the deleting or cutting process. This can be directly verified for each step in Definitions \ref{def.delete} and \ref{def.cutting}.
\end{proof}
\begin{remark}\label{rem.cutting} We make a few remarks concerning Definitions \ref{def.cutting}--\ref{def.ov_connect}.
\begin{enumerate}
\item The way to understand cutting is the following: suppose we cut $A$ as free (and $\Mb\backslash A$ as fixed) and consider a maximal ov-segment $\sigma$ as in Definition \ref{def.cutting} (\ref{it.cutting_2}), see {\color{blue}Figure \ref{fig.cutting2}}. Then, in $\Mb_2$ (which has atom set $\Mb\backslash A$), \emph{every possible bond} in $\sigma$ is broken, and each atom in $\Mb_2$ has a new fixed end. We may simply say that ``the ov-segment $\sigma$ is broken in $\Mb_2$". On the other hand, in $\Mb_1$ (which has atom set $A$), the ov-segments are preserved and the atoms in $\Mb_1\cap\sigma$ are still ov-adjacent along $\sigma$ (with the $\Mb_2$ atoms between them neglected) and they do not have any new fixed end.
\item In particular, if we cut $A$ as free, then this cutting generates (one or more) fixed end at every atom $\nf\in\Mb\backslash A$ that is \emph{ov-adjacent to an atom in $A$}, and not at any other atom. When $\nf$ is C-atom, this fixed end is \emph{bottom} (resp. \emph{top}) if $\nf$ is \emph{ov-parent} (resp. \emph{ov-child}) of an atom in $A$ (as for O-atoms, the bottom/top distinction is less important, but it is important that $\Mb$ is regular, see Remark \ref{rem.reg}).
\end{enumerate}
\end{remark}
\begin{remark}\label{rem.reg} Recall Definition \ref{def.reg} or regular molecules. This regularity will be retained in all algorithms: any simple pair formed by fixed end at O-atom and the serial free end is completely unrelated to the rest of the molecule (and this O-atom gets detached from both its neighbors). Some simple consequences of this will be used frequently below: (i) any O-atom of deg 2 must not be adjacent to any other atom, and cutting it does not affect the rest of the molecule; (ii) for any O-atom $\nf$ of deg 3, it can be connected to other atoms along one ov-segment. All ov-parents and ov-children of $\nf$ (if any) must belong to the same (maximal) ov-segment and are ov-adjacent to each other.
\end{remark}
\begin{remark} Note that, deleting O-atoms (Definition \ref{def.delete}) does not affect the correspondence between molecular and physical pictures under topological reduction. However, after cutting, this correspondence will break down, and we will only talk in the molecule language afterwards (for example notions like $\pb$ will only be referring to particle lines as in Definition \ref{def.molecule} (\ref{it.particle_line_same})), instead of physical particles.
\end{remark}

We now discuss the operation of splitting into sub-cases. This is not an operation on the molecule $\Mb$ but rather an operation on the domain of integration (or equivalently on the support of $Q$) in the integral (\ref{eq.associated_int_op}). In reality, it is used to classify the various degenerate and non-degenerate cases witch will require different strategy in the algorithm (cf. Remark \ref{rem.elementary}).
\begin{definition}\label{def.splitting} In \textbf{splitting}, we do not change $\Mb$ but split the support of $Q$ in (\ref{eq.associated_int_op}) as follows: decompose $1$ into finitely many indicator functions $1=\sum_{j}\mathbbm{1}_{\Sc_j}$, and decompose $Q=\sum_j Q\mathbbm{1}_{\Sc_j}$. Each of these terms corresponds to a \textbf{sub-case} where we make extra \textbf{restrictions} on the support of $Q$ provided by the conditions defining $\Sc_j$. Examples of such defining conditions include $|x_e-x_{e'}|\geq\lambda$ (or oppositely $|x_e-x_{e'}|\leq\lambda$), $|v_e-v_{e'}|\geq\lambda$ or $|t_\nf-t_{\nf'}|\geq\lambda$ (or opposite) etc., where $(e,e')$ are any two given edges, $(\nf,\nf')$ are any two given atoms, and $\lambda$ is any given quantity, cf. the definition of good molecules in Definition \ref{def.good_normal}.
\end{definition}

Now we define the notion of operation sequences, which are sequences consisting of operations of deleting, cutting and splitting:

\begin{definition}[Operation sequence]\label{def.cutting_algorithm} Given a molecule $\Mb$ and $Q$, define an \textbf{operation sequence} to be a sequence of operations acting on $\Mb$ or $Q$, such that each one is either \textbf{cutting} a subset of $\Mb$ as in Definition \ref{def.cutting}, or \textbf{deleting} an O-atom in $\Mb$ as in Definition \ref{def.delete}, or \textbf{splitting} into finitely many sub-cases (by decomposing $1$ into indicator functions) as in Definition \ref{def.splitting}. We require that all the deletings must occur before all the cuttings (splittings can occur anywhere in the sequence). If an operation sequence does not contain splitting, then we call it a \textbf{cutting sequence}.
\end{definition}

\begin{figure}[h!]
    \centering
    \includegraphics[width=0.55\linewidth]{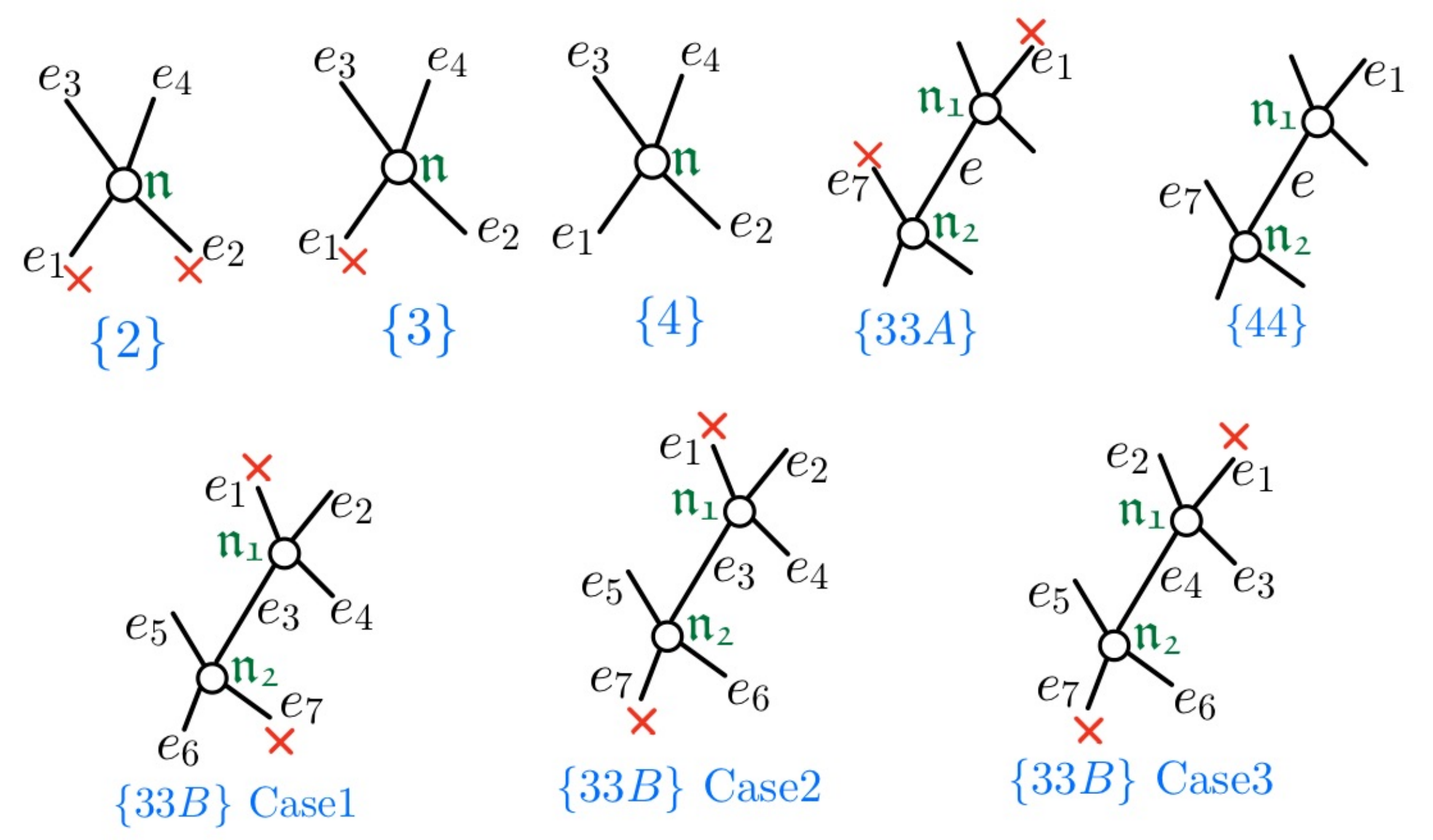}
    \caption{The elementary molecules as in Definition \ref{def.elementary}. Note the good/normal/bad property of this molecule depends on the extra restrictions on the variables due to splitting, see Definition \ref{def.good_normal}; the edges $e_j$ and atoms $\nf_j$ are marked as will be used in the proof of Propositions \ref{prop.intmini}--\ref{prop.intmini_2}.}
    \label{fig.elementary}
\end{figure}

Finally, we define the elementary molecules, which are the result of the cutting sequences we will construct in the algorithms below, and the main object of study in Section \ref{sec.treat_integral}:
\begin{definition}[Elementary molecules]\label{def.elementary} We define a molecule $\Mb$ to be elementary, if it satisfies one of the properties, see {\color{blue} Figure \ref{fig.elementary}}:
\begin{enumerate}
\item \textbf{\{2\} molecules}: $\Mb$ contains one atom of deg 2, with both fixed ends being top (or both being bottom) in the case of C-atom, and with the two fixed ends not serial in the case of O-atom.
\item \textbf{\{3\} molecules}: $\Mb$ contains one atom of deg 3.
\item \textbf{\{4\} molecules}: $\Mb$ contains one atom of deg 4.
\item\textbf{\{33\} molecules}: $\Mb$ contains two atoms connected by a bond, which are both deg 3. We say $\Mb$ is \textbf{\{33A\} molecule} if we can cut one atom as free such that the other atom becomes a \{2\} molecule as defined above. Otherwise we say $\Mb$ is \textbf{\{33B\} molecule}.
\item\textbf{\{44\} molecules}: $\Mb$ contains two atoms connected by a bond, which are both deg 4.
\end{enumerate}
\end{definition}

\subsection{Operations on associated operators}\label{sec.oper_int} In this subsection, we study the behavior of $\Ic_\Mb (Q)$ in (\ref{eq.associated_int_op}) under the cutting, deleting and splitting operations. Namely, we will prove $\Ic_\Mb(Q)\leq \Ic_{\Mb'}(Q)$ after deleting (Propositions \ref{prop.deleting}) and $\Ic_\Mb(Q)=\Ic_{\Mb_1}\circ\Ic_{\Mb_2}(Q)$ after cutting (Proposition \ref{prop.cutting}).

Recall the integral (\ref{eq.associated_int_op}), where $Q$ is a function of $\vt_\Mb$ and $\vz_\Ec$. If we perform deleting or cutting to $\Mb$ and turn it into $\Mb'$ or $\Mb'=\Mb_1\sqcup\Mb_2$, then we will get a new set of vector variables $z_e$ associated with the new molecule (the time variables remain the same as atoms are not affected by cutting). In order to define $\Ic_{\Mb'}(Q)$ and $\Ic_{\Mb_1}\circ\Ic_{\Mb_2}(Q)$ for the same $Q$, we need to match the vector variables for $\Mb'$ with those for $\Mb$. More precisely, we need to define the following two mappings:
\begin{itemize}
\item A mapping from the set of \emph{bonds plus free ends} of $\Mb$ to that of $\Mb'$; this decides how we transform the function $Q$ in (\ref{eq.associated_int_op}) between the integral variables of $\Mb$ and those of $\Mb'$.
\item A mapping from the set of new \emph{fixed ends} of $\Mb_2$ to the set of new \emph{bonds plus free ends} of $\Mb_1$; this decides, while evaluating $\Ic_{\Mb_2}(Q)$, the value that $z_f$ is \emph{fixed to be}, for each fixed end $f$ of $\Mb_2$.
\end{itemize}

\begin{figure}[h!]
    \centering
    \includegraphics[width=0.65\linewidth]{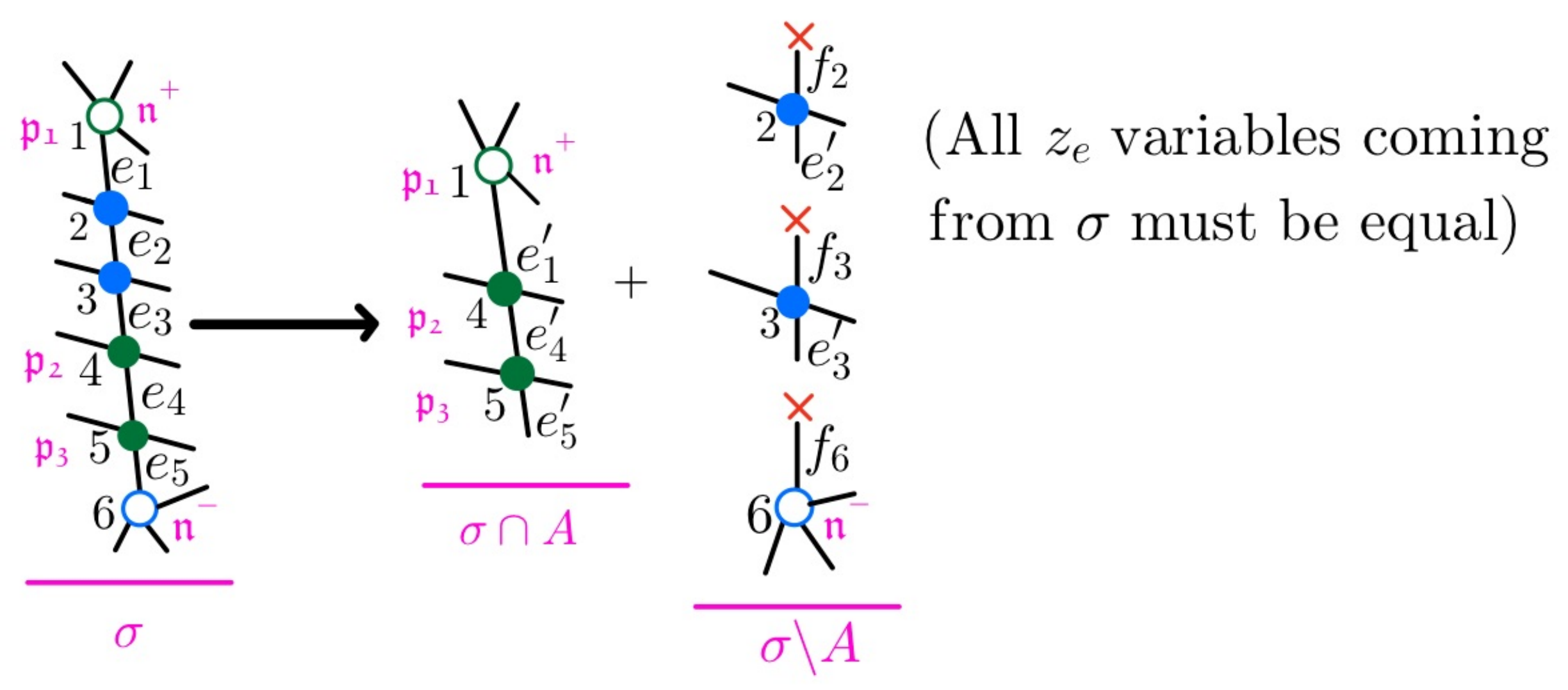}
    \caption{The transformation of associated variables under cutting (Definition \ref{def.cutting}), see Definition \ref{def.associated_vars_op}. This picture shows the case of the example shown in {\color{blue}Figure \ref{fig.cutting2}}, where $e_j\mapsto e_j'$ and $f_j\mapsto e_1'$ in the notation of the mapping in Definition \ref{def.associated_vars_op}. Note that all the $z_e$ variables for edges $e$ coming from this ov-segment are all equal to each other, see the proof of Proposition \ref{prop.cutting}.}
    \label{fig.cutting2+}
\end{figure}

In Definition \ref{def.associated_vars_op} below, we will construct these mappings under deleting and cutting operations, which provides the correspondences between vector variables $z_e$.
\begin{definition}[Transformation of associated variables under cutting and deleting]\label{def.associated_vars_op} Suppose $\Mb'$ is formed from $\Mb$ by deleting or cutting. Define the following mappings:
\begin{enumerate}
    \item\label{it.associated_vars_op_2} \emph{The case of deleting}. Let $e_j$ and $(e_1',e_2')$ be as in Definition \ref{def.delete} (\ref{it.delete_new}), we define the mapping  $e_1,e_3\mapsto e_1'$ and $e_2,e_4\mapsto e_2'$. This means, we identify the variables $z_{e_1}$ and $z_{e_3}$ before deleting, with $z_{e_1'}$ after deleting, and same for $(e_2,e_4)$.
    \item \label{it.associated_vars_op_2.5} \emph{The case of cutting: no O-atoms}. As in Definition \ref{def.cutting} (\ref{it.cutting_1}), for each bond $e$ between $A$ and $\Mb\backslash A$ that is turned into a free end $e_1\in\Mb_1$ and a fixed end $e_2\in \Mb_2$, we define the mapping $e\mapsto e_1$ and $e_2\mapsto e_1$. This means, we identify $z_{e}$ before cutting with $z_{e_1}$ after cutting, and also fix the fixed end variable $z_{e_2}$ to be $z_{e_1}$.
    \item\label{it.associated_vars_op_3} \emph{The case of cutting: general}.  Let $\sigma$ and related notions be as in Definition \ref{def.cutting} (\ref{it.cutting_2}), we define the mapping as follows. For each bond $e\in\sigma$, consider the two atoms of $e$, and the one that is closer to $\pf_1$; we then map $e$ to the free end $e'\in\Mb'$ at this atom. For each fixed end $e\in\Mb_2$, we map it to the same bottom edge $e'\in\Mb_1\cap\sigma$ at $\pf_1$ (it could be a free end or a bond). For example, in {\color{blue}Figure \ref{fig.cutting2+}}, we have $e_j\mapsto e_j'$ and $f_j\mapsto e_1'$. We also identify and fix the relevant variables as in (\ref{it.associated_vars_op_2})--(\ref{it.associated_vars_op_2.5}).
\end{enumerate}
\end{definition}
\begin{proposition}\label{prop.deleting}
Let $\Mb'$ be formed from $\Mb$ by deleting one O-atom (Definition \ref{def.delete}). Then we have 
\begin{equation}\label{eq.deleting}
    \Ic_{\Mb}(Q) \le \Ic_{\Mb'}(Q), \qquad \forall \,Q\ge 0.
\end{equation}
\end{proposition}

\begin{proof} Note that the power $|\Ec_*|-2|\Mb|$ in (\ref{eq.associated_int_op}) is the same for $\Mb$ and $\Mb'$, as we have $1$ less atom and $2$ less bonds/free ends when going from $\Mb$ to $\Mb'$. Now let $e_j$ and $(e_1',e_2')$ be as in Definition \ref{def.delete} (\ref{it.delete_new}), then we have $z_{e_1'}=z_{e_1}=z_{e_3}$ by the mapping in Definition \ref{def.associated_vars_op} (\ref{it.associated_vars_op_2}) and the Dirac factors in (\ref{eq.associated_dist_O}), and similarly $z_{e_2'}=z_{e_2}=z_{e_4}$. Concerning the integral in (\ref{eq.associated_int_op}), using the identification of variables in Definition \ref{def.associated_vars_op}, we see that $\Ic_\Mb(Q)$ is the same as $\Ic_{\Mb'}(Q)$ but with the extra integral
\[\int\dirac(z_{e_3}-z_{e_1})\dirac(z_{e_4}-z_{e_2})\dirac\big(|x_{e_1}-x_{e_2}+t_\nf(v_{e_1}-v_{e_2})|-\varepsilon\big)\cdot[(v_{e_1}-v_{e_2})\cdot\omega]_-\,\mathrm{d}z_{e_3}\mathrm{d}z_{e_4}\mathrm{d}t_\nf.\] By (\ref{eq.singlecol}) we know this integral equals an indicator function, which is $\leq 1$, hence $\Ic_\Mb(Q)\leq \Ic_{\Mb'}(Q)$.
\end{proof}

\begin{proposition}\label{prop.cutting}
Assume that a molecule $\Mb$ is cut into $\Mb_1$ and $\Mb_2$ (Definition \ref{def.cutting}), where $\Mb_1$ and $\Mb_2$ are cut as free and fixed respectively. Then we have
\begin{equation}\label{eq.cutting}
    \Ic_{\Mb}(Q) = \Ic_{\Mb_1}\circ \Ic_{\Mb_2}(Q).
\end{equation}

\end{proposition}
\begin{proof} Recall the formula \eqref{eq.associated_int_op}. Let $\Ec_{*,j}$ and $\Ec_*'$ be the $\Ec_*$ sets for $\Mb_j$ and $\Mb'$ respectively, then it is easy to see that $|\Mb|=|\Mb_1|+|\Mb_2|$ and $|\Ec_*|=|\Ec_{*,1}|+|\Ec_{*,2}|=|\Ec_*'|$, so the powers on both sides of (\ref{eq.cutting}) match. Below we will neglect these powers of $\varepsilon$ for simplicity. Now by \eqref{eq.associated_int_op}, we have
\begin{equation}\label{eq.proof_cutting_lemma_1}
        \Ic_\Mb(Q)=\int_{\mathbb{R}^{2d|\Ec_*|}\times \Rb^{|\Mb|}} 
        \bigg(\prod_{\nf\in\Mb}\boldsymbol{\Delta}_\nf\bigg)\cdot Q\,\mathrm{d}\vz_{\Ec_*}\,\mathrm{d}\vt_\Mb.
        \end{equation}  
          On the other hand, by Fubini we have
          \begin{equation}\label{eq.proof_cutting_lemma_1+}
          \begin{aligned}
        \Ic_{\Mb_1}\circ\Ic_{\Mb_2}(Q)&=\int_{\mathbb{R}^{2d|\Ec_{1,*}|}\times \Rb^{|\Mb_1|}} \bigg(\prod_{\nf_1\in\Mb_1}\boldsymbol{\Delta}_{\nf_1}\bigg)\,\mathrm{d}\vz_{\Ec_{1,*}}\mathrm{d}\vt_{\Mb_1}\int_{\mathbb{R}^{2d|\Ec_{2,*}|}\times \Rb^{|\Mb_2|}} \bigg(\prod_{\nf_2\in\Mb_2}\boldsymbol{\Delta}_{\nf_2}\bigg)\cdot Q\,\mathrm{d}\vz_{\Ec_{2,*}}\mathrm{d}\vt_{\Mb_2}
    \\
    &=\int_{\mathbb{R}^{2d|\Ec_*'|}\times \Rb^{|\Mb|}} 
        \bigg(\prod_{\nf\in\Mb}\boldsymbol{\Delta}_\nf \bigg)\cdot Q \,\mathrm{d}\vz_{\Ec_*'}\mathrm{d}\vt_\Mb.
        \end{aligned}
        \end{equation} Here in the second line in (\ref{eq.proof_cutting_lemma_1+}), we are identifying the variables $\vz_{\Ec_*'}$ with the variables $\vz_{\Ec_*}$, and fixing each fixed end variable $z_e$ in $\Mb_2$ to be the corresponding free end variable $z_{e'}$ in $\Mb_1$, by the mapping in Definition \ref{def.associated_vars_op}. it then suffices to show that (\ref{eq.proof_cutting_lemma_1}) equals (\ref{eq.proof_cutting_lemma_1+}) under the above identification and assignments.
        
\textbf{No O-atom case.} Assume $\Mb$ has no O-atom, then we are in the simpler scenario in Definition \ref{def.cutting} (\ref{it.cutting_1}). In this case, each new free end $e_1\in\Mb_1$ is paired with a new fixed end $e_2\in\Mb_2$ and they come from breaking a bond $e\in\Mb$; by Definition \ref{def.associated_vars_op} (\ref{it.associated_vars_op_2.5}) we see that the integration in (\ref{eq.proof_cutting_lemma_1}) and (\ref{eq.proof_cutting_lemma_1+}) are the same. Moreover, the $\Dirac_\nf$ and $Q$ factors are also the same under the above identification and assignments, so trivially (\ref{eq.proof_cutting_lemma_1}) equals (\ref{eq.proof_cutting_lemma_1+}).

\textbf{General case.} In the general case with O-atoms, the identification and assignments of variables and more complicated. We only need to consider a fixed maximal ov-segment $\sigma$ (along with the relevant notions $\pf_j$) as in Definition \ref{def.cutting} (\ref{it.cutting_2}). Let the bonds in $\sigma$ be $e_j,\,(1\leq j\leq q)$ from top to bottom, then in view of the Dirac factors in (\ref{eq.associated_dist_O}), we know the integral (\ref{eq.proof_cutting_lemma_1}) involves the Dirac factors 
\begin{equation}\label{eq.compare_delta_0}\prod_{j=1}^{q-1}\dirac(z_{e_{j+1}}-z_{e_j}).\end{equation} On the other hand, for each $\pf_i\in A\cap\sigma$, consider the bond connecting $\pf_i$ to its child in $\sigma$, and let it be $e_{j_i}$. Then, by the identification and assignments as in Definition \ref{def.associated_vars_op} (\ref{it.associated_vars_op_3}), we know that the integral (\ref{eq.proof_cutting_lemma_1+}) involves the Dirac factors
\begin{equation}\label{eq.compare_delta_1}\prod_{i=1}^{s-1}\dirac\big(z_{e_{j_{i+1}}}-z_{e_{j_i}}\big)\cdot\prod_{e\not\in\{e_{j_i}\}}\dirac\big(z_e-z_{e_{j_1}}\big)\end{equation} It is easy to see that (\ref{eq.compare_delta_0})=(\ref{eq.compare_delta_1}), as both sides are exactly the Dirac mass supported on the linear subspace $\{z_{e_1}=\cdots =z_{e_q}\}$. For example, in the context of {\color{blue}Figure \ref{fig.cutting2+}}, where $(j_1,j_2,j_3)=(1,4,5)$, we have
\[\dirac(z_{e_1}-z_{e_2})\dirac(z_{e_2}-z_{e_3})\dirac(z_{e_3}-z_{e_4})\dirac(z_{e_4}-z_{e_5})=\dirac(z_{e_1}-z_{e_4})\dirac(z_{e_4}-z_{e_5})\cdot\dirac(z_{e_2}-z_{e_1})\dirac(z_{e_3}-z_{e_1}).\]

Finally, note that all the Dirac functions in (\ref{eq.compare_delta_0}) and (\ref{eq.compare_delta_1}) imply that all the variables $z_{e_j}$ are equal, so by inserting them into the rest of the $\Dirac_\nf$ and $Q$ factors in (\ref{eq.proof_cutting_lemma_1}) and (\ref{eq.proof_cutting_lemma_1+}), we conclude that these remaining factors are also the same for (\ref{eq.proof_cutting_lemma_1}) and (\ref{eq.proof_cutting_lemma_1+}). This then shows that (\ref{eq.proof_cutting_lemma_1}) equals (\ref{eq.proof_cutting_lemma_1+}), and completes the proof of Proposition \ref{prop.cutting}.
\end{proof}
\begin{remark}
In general, \eqref{eq.cutting} holds with both sides being interpreted in the sense of distributions. However, in this paper, we will only encounter the case where $\Ic_\Mb(Q)$ is a function whenever $Q$ is, in which case \eqref{eq.cutting} holds classically. The condition that $\Ic_{\Mb}(Q)$ is a function is really a condition on the molecule $\Mb$, which restricts the number and places of its fixed ends, and will always be satisfies throughout the algorithm. That being said, one can still interpret \eqref{eq.cutting} in bigger generality in the sense of distributions as follows: If one replaces all the delta functions appearing in the definition of $\Ic_\Mb$ by approximations of identity, then the analog of formula \eqref{eq.cutting} holds with no ambiguity. The same result then holds by taking the limits of those approximations of identity (in the appropriate distributional sense).
\end{remark}

Note that the operators $\Ic_{\Mb_1}$ and $\Ic_{\Mb_2}$ in Proposition \ref{prop.cutting} in general do not commute (unless $\Mb$ has no fixed end). Therefore, it is important to order the resulting components after cutting, such that $\Mb_1$ which is cut as free always occurs to the left of $\Mb_2$ which is cut as fixed (i.e. the operator $\Ic_{\Mb_2}$ is applied first, and the corresponding variables are integrated first in the iterated integral). This ordering is defined in the following definition.
\begin{definition}[Ordering components after cutting]\label{def.ordered_disjoint} Suppose a molecule $\Mb$ becomes a new molecule $\Mb'$ after a sequence of cutting operations; we refer to the connected components of $\Mb'$ as final components. We define a total ordering $\prec_{\mathrm{cut}}$ between these final components inductively, as follows (see {\color{blue} Figure \ref{fig.choice_2}} for an example):
\begin{enumerate}
    \item\label{it.disjoint_1} In the first cutting operation, suppose $\Mb$ is cut into $\Mb_1$ (as free) and $\Mb_2$ (as fixed). Then for any final component $\Xc,\Yc$ of $\Mb'$, we require that $\Xc\prec_{\mathrm{cut}}\Yc$ if $\Xc\subseteq\Mb_1$ and $\Yc\subseteq\Mb_2$.
    \item For the final components that are contained in $\Mb_1$ (resp. in $\Mb_2$), order them inductively using $\prec_{\mathrm{cut}}$ defined for the subsequent sequence of cutting operations applied to $\Mb_1$ (resp. to $\Mb_2$).
    \end{enumerate}
    
    If $\Mb'$ is the result of $\Mb$ after (deleting and) cutting, and $\Mb'$ is the disjoint union of its final components, i.e. $\Mb' = \Mb_1 \sqcup \cdots \sqcup \Mb_k$ with $\Mb_1 \prec_{\mathrm{cut}} \cdots \prec_{\mathrm{cut}} \Mb_k$, then we define the \textbf{associated integral operator} by
    \begin{equation}\label{eq.associated_op_disjoint}
        \Ic_{\Mb'} = \Ic_{\Mb_1}\circ \Ic_{\Mb_2}\circ\cdots\circ \Ic_{\Mb_k}.
    \end{equation}

\end{definition}
\begin{remark}\label{rem.cut_order} As explained above, if $\Mb_i\prec_{\mathrm{cut}}\Mb_j$, then the operator $\Ic_{\Mb_j}$ is applied first, and the $\Mb_j$ variables are integrated first. In many cases of the algorithm below, each time we will cut an elementary component as free, and work on the remaining molecule (cf. Remark \ref{rem.aftercut} below); in this case $\prec_{\mathrm{cut}}$ is just the order of being cut, i.e. $\Yc\prec_{\mathrm{cut}}\Xc$ if and only if $\Yc$ is cut before $\Xc$. This ordering will be useful in identifying a good molecule after splitting, see Definition \ref{def.good_normal} (\ref{it.good_1}), (\ref{it.good_3}).

Suppose we have a cutting sequence, and one of the cuttings breaks a maximal ov-segment $\sigma$ as in Definition \ref{def.cutting} and creates free and fixed ends. Then for any such fixed end $e$, there always exists a free end or bond $e'$ such that $z_e=z_{e'}$, and we have $\Yc\prec_{\mathrm{cut}}\Xc$ where $\Xc$ and $\Yc$ are the final components containing $e$ and $e'$ respectively. In fact, we can choose $e'$ to be the image of $e$ under the mapping in Definition \ref{def.associated_vars_op} (\ref{it.associated_vars_op_3}), then $\Yc\prec_{\mathrm{cut}}\Xc$ follows from Definition \ref{def.ordered_disjoint} (\ref{it.disjoint_1}). This will also be helpful in getting good molecules below.
\end{remark}

By putting together Definition \ref{def.splitting}, Proposition \ref{prop.deleting} and Proposition \ref{prop.cutting}, we arrive at the main conclusion of this section, namely
\begin{proposition}[]\label{prop.cutting_algorithm} Start with a full molecule $\Mb$ and $\Ic_\Mb(Q)$ as in (\ref{eq.local_int}), where $Q=Q_\Mb$ as in (\ref{eq.Q_M}). If we apply any cutting sequence (Definition \ref{def.cutting_algorithm}) to $\Mb$ and turn it to $\Mb'$, with the associated variables transformed as in Definition \ref{def.associated_vars_op}, then we always have
\begin{equation}\label{eq.cutting_sequence}
\Ic_\Mb(Q)\leq \Ic_{\Mb'}(Q),
\end{equation}
where $\Ic_\Mb(Q)$ is defined as in (\ref{eq.associated_op_disjoint}) with $\Mb_1\prec_{\mathrm{cut}}\cdots\prec_{\mathrm{cut}}\Mb_k$ being the final components of $\Mb'$.

If we apply any operation sequence (including also splitting which operates on the support of $Q$) to $\Ic_\Mb(Q)$ and get a contribution $\Ic_{\Mb'}(Q')$ for each sub-case, then we have   
\begin{equation}\label{eq.cutting_algorithm}
    \Ic_{\Mb}(Q)\le \sum_{\textrm{sub-cases}}\Ic_{\Mb'}(Q'),
\end{equation}  
where the summation is taken over all sub-cases. Note that $Q' = Q\cdot\mathbbm{1}_\Pc$ for some $\Pc$ representing the restriction in this sub-case (see Definition \ref{def.splitting}).
\end{proposition}
\begin{proof} We only need to prove (\ref{eq.cutting_sequence}). In fact, for any splitting $Q=\sum_{j}Q\mathbbm{1}_{\Sc_j}$ as in Definition \ref{def.splitting}, we trivially have $\Ic_\Mb(Q)=\sum_j\Ic_\Mb(Q\mathbbm{1}_{\Sc_j})$, so (\ref{eq.cutting_algorithm}) follows from (\ref{eq.cutting_sequence}).

Now to prove (\ref{eq.cutting_sequence}), using Proposition \ref{prop.deleting}, we can treat the deleting operation and only need to consider cuttings, in which case we will prove the equality in (\ref{eq.cutting_sequence}). We proceed by induction: the case of one cutting follows from Proposition \ref{prop.cutting}. Now for any $\Mb$, say the first cutting cuts $\Mb_1$ as free and $\Mb_2$ as fixed. In subsequent cuttings, assume $\Mb_1$ is cut into $\Mb_1'=\sqcup_{j}\Mb_{1j}$ (in increasing $\prec_{\mathrm{cut}}$ order) and $\Mb_2$ is cut into $\Mb_2'=\sqcup_j\Mb_{2j}$ (in increasing $\prec_{\mathrm{cut}}$ order), then the final result of the cutting sequence is $\Mb'=(\sqcup_{j}\Mb_{1j})\sqcup (\sqcup_{j}\Mb_{2j})$ in increasing $\prec_{\mathrm{cut}}$ order, due to Definition \ref{def.ordered_disjoint}. Then by induction hypothesis we have
\[\Ic_\Mb(Q)=\Ic_{\Mb_1}\circ\Ic_{\Mb_2}(Q)=\bigg(\prod_{j}\Ic_{\Mb_{1j}}\bigg)\circ\bigg(\prod_{j}\Ic_{\Mb_{2j}}\bigg)(Q)=\Ic_{\Mb'}(Q),\] which completes the proof.
\end{proof}

\section{Estimating the integral $\Ic_\Mb(Q_\Mb)$}\label{sec.treat_integral}

In this section, we complete the proof of Proposition \ref{prop.cumulant_est}, pending the proof of the pure combinatorial Proposition \ref{prop.comb_est} (which is in Sections \ref{sec.prepare}--\ref{sec.maincr}). Starting from Proposition \ref{prop.cumulant_formula}, using Propositions \ref{prop.layerrec3}, \ref{prop.local_int} and \ref{prop.cutting_algorithm}, we can reduce Proposition \ref{prop.cumulant_est} to the estimate of quantities of form $\Ic_{\Mb'}(Q')$, where $\Mb'$ (and $Q'$) are possible results of an operation sequence applied to $\Mb$ (and $Q=Q_\Mb$), corresponding to possible sub-cases. 

The operation sequence is constructed in Proposition \ref{prop.comb_est} using the algorithms in Sections \ref{sec.prepare}--\ref{sec.maincr}. In Proposition \ref{prop.comb_est} we can guarantee that:
\begin{itemize}
\item The number of sub-cases is acceptable (i.e. $\leq C^{|\Mb|}\cdot|\log\varepsilon|^{C^*\rho}$ in Proposition \ref{prop.comb_est});
\item $\Mb'$ is formed by elementary molecules (Definition \ref{def.elementary});
\item We have adequate control of these elementary molecules (and associated support information of $Q'$) in terms of good and bad molecules (Definition \ref{def.good_normal} and (\ref{eq.overall_alg})).
\end{itemize}

The goal of this section is to use the above information to prove an estimate for $\Ic_{\Mb'}(Q')$, namely (\ref{eq.proof_int_est_from_alg_step0_*}). Note that this $\Ic_{\Mb'}$ is a composition of $\Ic_{\Mb_j}$ where $\Mb_j$ are components of $\Mb'$ which are elementary molecules. In Section \ref{sec.elem_int} we will estimate each $\Ic_{\Mb_j}$, and classify them as good, normal, and bad in Definition \ref{def.good_normal}. Then in Section \ref{sec.summary}, we put together these local estimates and apply Proposition \ref{prop.comb_est}, to prove (\ref{eq.proof_int_est_from_alg_step0_*}) and thus finish the proof of Proposition \ref{prop.cumulant_est}. These arguments rely on some technical estimates concerning volumes of $(v_e)$ and cross-sections in (\ref{eq.associated_int_op}), which we leave to Section \ref{sec.aux}.
\subsection{Integrals for elementary molecules}\label{sec.elem_int} In this subsection, we prove upper bounds for $\Ic_{\Mb}$ for an elementary molecule $\Mb$, in Propositions \ref{prop.intmini}--\ref{prop.intmini_2}. We will use $\nf_j$ and $e_j$ to denote atoms and edges in $\Mb$, and abbreviate $t_{\nf_j}:=t_j$ and $z_{e_j}:=z_j$. For future reference, let us explicitly write out the support of $\boldsymbol{\Delta}_{\nf}(z_{e_1}, z_{e_2}, z_{e_1'}, z_{e_2'}, t_\nf) = \boldsymbol{\Delta}_{\nf}(z_1, z_2, z_3, z_4, t_1)$, which is defined in \eqref{eq.associated_dist_C} and \eqref{eq.associated_dist_O}, as follows:
\begin{equation}\label{eq.delta_supp}
\left\{
\begin{aligned}
x_2&=x_1+t_1(v_1-v_2)-\varepsilon\omega,
\\
x_3&=x_1+t_1\big[(v_1-v_2)\cdot\omega\big]\omega,
\\
v_3&=v_1-\big[(v_1-v_2)\cdot\omega\big]\omega,
\\
x_4&=x_2-t_1\big[(v_1-v_2)\cdot\omega\big]\omega,
\\
v_4&=v_2+[(v_1-v_2)\cdot\omega\big]\omega,
\end{aligned}
\right.\ \textrm{if $\nf$ is C-atom;}\quad \left\{
\begin{aligned}
x_2&=x_1+t_1(v_1-v_2)-\varepsilon\omega,\\
x_3&=x_1,\\
v_3&=v_1,\\
x_4&=x_2=x_1+t_1(v_1-v_2)-\varepsilon\omega,\\
v_4&=v_2
\end{aligned}
\right.\ \textrm{if $\nf$ is O-atom.}
\end{equation}

\begin{proposition}[Upper bounds for one-atom $\Ic_\Mb$]\label{prop.intmini} Let $\Ic_\Mb$ be defined as in \eqref{eq.associated_int_op}, where $\Mb$ contains only one atom $\nf$ with 4 edges $e_1,\dots,e_4$, see {\color{blue}Figure \ref{fig.elementary}}. We assume that $e_1$ (resp. $e_2$) is serial with $e_3$ (resp. $e_4$). Let $Q$ be a function supported in $|v_j|\leq|\log\varepsilon|^{C^*}$ and $|x_j-x_j^0|\leq |\log\varepsilon|^{C^*}$, where $(x_j,v_j)=z_{e_j}$ and $x_j^0$ are fixed vectors, and also $\vt_\Mb\in\Dc$ with $\Dc$ as in (\ref{eq.associated_domain}).
\begin{enumerate}
\item\label{it.intmini_1} Suppose $\Mb$ is a \{2\} molecule as in Definition \ref{def.elementary}, i.e. the two fixed ends are either both bottom or both top, by symmetry we may assume $(e_1,e_2)$ are fixed and $(e_3,e_4)$ are free. Then
\begin{equation}\label{eq.intmini_1}
    \Ic_\Mb(Q) = \oneb_{\mathrm{col}}(z_1,z_2)\cdot Q.
\end{equation} 
Here in \eqref{eq.intmini_1}, the indicator $\oneb_{\mathrm{col}}$ is defined to be $1$ if and only if there exists a unique $t_1$ such that 
\begin{equation}\label{eq.intmini_1+}
    |x_1-x_2+t_1(v_1-v_2)|=\varepsilon,\quad (v_1-v_2)\cdot\omega\leq 0,
\end{equation} 
where $\omega:=\varepsilon^{-1}(x_1-x_2+t_1(v_1-v_2))$ is a unit vector. We also have the following estimate
\begin{equation}\label{eq.int_elem_est_deg2}
    \|\Ic_\Mb(Q)\|_{L^\infty}\leq \|Q\|_{L^\infty}
\end{equation}
Moreover, in $\Ic_\Mb(Q)$, the function $Q$ depends on $(z_1, z_2, z_3, z_4, t_1)$. On the right-hand side of \eqref{eq.intmini_1}, $Q$ becomes a function of $(z_1, z_2)$, obtained by a change of variables applied to the original $Q$. This change of variables is determined as follows: $t_1$ is expressed as a function of $(z_1, z_2)$ via \eqref{eq.intmini_1+}, and $(z_3, z_4)$ is defined as a function of $(z_1, z_2)$, either by the last four equations in the left half of \eqref{eq.delta_supp} if $\Mb$ is a C-atom, or by $(z_3, z_4) = (z_1, z_2)$ if $\Mb$ is an O-atom.
\item\label{it.intmini_2} Suppose $\Mb$ is a \{3\} molecule as in Definition \ref{def.elementary}, by symmetry we may assume $e_1$ is fixed and $(e_2,e_3,e_4)$ are free. Then
\begin{equation}\label{eq.intmini_3}
\Ic_\Mb(Q) = \int_{\Rb\times\Sb^{d-1}\times\Rb^d}\big[(v_1-v_2)\cdot\omega\big]_-\cdot Q\,\mathrm{d}t_1\mathrm{d}\omega\mathrm{d}v_2.
\end{equation} 
Here in \eqref{eq.intmini_3}, $(z_1,z_2,z_3,z_4,t_1)$ in the input variables of $Q$ are functions of the integral variables $(t, \omega, v_2)$ given by \eqref{eq.delta_supp}, with $(x_1,v_1)$ fixed. In the integral \eqref{eq.intmini_3} the domain of integration can be restricted to $t_1\in [(\ell[\nf]-1)\tau,\ell[\nf]\tau]$ since $\supp (Q)\subseteq \Dc$. 

\item\label{it.intmini_3} Suppose $\Mb$ is a \{4\} molecule as in Definition \ref{def.elementary}. Then
\begin{equation}\label{eq.intmini_6}
    \Ic_\Mb(Q)=\varepsilon^{-(d-1)}\int_{\Rb\times\Sb^{d-1}\times(\Rb^d)^3}\big[(v_1-v_2)\cdot\omega\big]_-\cdot Q\,\mathrm{d}t_1\mathrm{d}\omega\mathrm{d}x_1\mathrm{d}v_1\mathrm{d}v_2.
\end{equation} 
Here in \eqref{eq.intmini_6}, $(z_1,z_2,z_3,z_4,t_1)$ in the input variables of $Q$ are functions of the integral variables $(t, \omega, x_1,v_2,v_2)$ given by \eqref{eq.delta_supp}. In the integral \eqref{eq.intmini_6} the domain of integration can be restricted to $t_1\in [(\ell[\nf]-1)\tau,\ell[\nf]\tau]$ since $\supp Q\subseteq \Dc$.
We also have the following estimate
\begin{equation}\label{eq.int_elem_est_deg4}
    \|\Ic_\Mb(Q)\|_{L^\infty}\leq \varepsilon^{-(d-1)}|\log\varepsilon|^{C^*}\cdot\|Q\|_{L^\infty}.
\end{equation}
\item\label{it.intmini_empty_line} Suppose $\Mb$ is an empty end. Then
\begin{equation}\label{eq.intmini_empty_line}
    \Ic_\Mb(Q):=\varepsilon^{-(d-1)}\int_{\mathbb{R}^{2d}} Q \,\mathrm{d}z_{e}.
\end{equation} 
We also have the following estimate
\begin{equation}\label{eq.int_elem_est_empty_line}
    \|\Ic_\Mb(Q)\|_{L^\infty}\leq \varepsilon^{-(d-1)}|\log\varepsilon|^{C^*}\cdot\|Q\|_{L^\infty}.
\end{equation}
\end{enumerate}
\end{proposition}
\begin{proof} We may assume $\|Q\|_{L^\infty}=1$. We divide the proof into several parts.
    
\textbf{Proof of \eqref{it.intmini_1}.} We start with \eqref{eq.intmini_1} and consider the C-atom case (the O-atom case is similar and much easier). Note that here we are fixing $(x_1,v_1,x_2,v_2)$ and integrating in $(x_3,v_3,x_4,v_4)$ and $t_1$, which is the same as in \eqref{eq.singlecol} in the proof of Proposition \ref{prop.local_int}. The equality \eqref{eq.intmini_1} then follows from the same proof, where we first integrate in $(x_3,x_4)$, then in $(v_3,v_4)$, and then in $t_1$. Then (\ref{eq.int_elem_est_deg2}) follows directly from \eqref{eq.intmini_1}.

\textbf{Proof of \eqref{it.intmini_2}.} Now consider \eqref{eq.intmini_3}; again consider the C-atom case, and note that here we are fixing $(x_1,v_1)$ and integrating in $(t_1,x_2,v_2,x_3,v_3,x_4,v_4)$, i.e.
\begin{equation}\label{eq.intmini_proof_1}
\begin{aligned}
\Ic_\Mb(Q)&=\varepsilon^{-(d-1)}\int_\Rb\int_{\Rb^{6d}}\dirac\big(x_3-x_1+t_1(v_3-v_1)\big)\cdot\dirac\big(x_4-x_2+t_1(v_4-v_2)\big)\\
&\times\dirac\big(v_3-v_1+[(v_1-v_2)\cdot\omega]\omega\big)\dirac\big(v_4-v_2-[(v_1-v_2)\cdot\omega]\omega\big)\\
&\times\dirac\big(|x_1-x_2+t_1(v_1-v_2)|-\varepsilon\big)\cdot\big[(v_1-v_2)\cdot\omega\big]_-\cdot Q(x_1,v_1,\cdots,x_4,v_4,t_1)\,\mathrm{d}t_1\prod_{j=2}^4\mathrm{d}x_j\mathrm{d}v_j,
\end{aligned}
\end{equation} where $\omega:=\varepsilon^{-1}(x_1-x_2+t_1(v_1-v_2))$ is a unit vector. In \eqref{eq.intmini_proof_1}, once we fix $(x_2,v_2)$, then the integral in $(x_3,x_4,v_3,v_4)$ is given by a change of variables in the function $Q$, so we have
\begin{equation}\label{eq.intmini_proof_2}
\Ic_\Mb(Q)=\varepsilon^{-(d-1)}\int_\Rb\int_{\Rb^{2d}}
\dirac\big(|x_1-x_2+t_1(v_1-v_2)|-\varepsilon\big)\cdot\big[(v_1-v_2)\cdot\omega\big]_-\cdot Q(x_1,v_1,\cdots,x_4,v_4,t_1)\,\mathrm{d}t_1\mathrm{d}x_2\mathrm{d}v_2,
\end{equation} where $\omega$ is as above, and
\begin{equation}\label{eq.intmini_proof_3}
\left\{
\begin{aligned}
x_3&=x_1+t_1\big[(v_1-v_2)\cdot\omega\big]\omega,\\
v_3&=v_1-\big[(v_1-v_2)\cdot\omega\big]\omega,\\
x_4&=x_2-t_1\big[(v_1-v_2)\cdot\omega\big]\omega,\\
v_4&=v_2+[(v_1-v_2)\cdot\omega\big]\omega.
\end{aligned}
\right.
\end{equation} Now in \eqref{eq.intmini_proof_2}, we may fix $v_2$ and make the change of variables $x_2\leftrightarrow \omega$, which leads to
\begin{equation}\label{eq.intmini_proof_4}
\Ic_\Mb(Q)=\int_\Rb\int_{\Rb^{2d}}
\dirac(|\omega|-1)\cdot\big[(v_1-v_2)\cdot\omega\big]_-\cdot Q(x_1,v_1,\cdots,x_4,v_4,t_1)\,\mathrm{d}t\mathrm{d}\omega\mathrm{d}v_2,
\end{equation} where $x_2=x_1+t_1(v_1-v_2)-\varepsilon\omega$, which matches \eqref{eq.delta_supp} together with \eqref{eq.intmini_proof_3}. The integral in $\omega$ in \eqref{eq.intmini_proof_4} is done on the unit sphere $\omega\in \Sb^{d-1}$, which then gives \eqref{eq.intmini_3}. 

\textbf{Proof of \eqref{it.intmini_3}.} By integrating in the variables $(x_1,v_1)$ in \eqref{eq.intmini_3}, we also obtain \eqref{eq.intmini_6}. Then \eqref{eq.int_elem_est_deg4} immediately follows from \eqref{eq.intmini_6} and the fact that $\big[(v_1-v_2)\cdot\omega\big]_-\le |\log\varepsilon|^{C^*}$ on the support of $Q$ (which is due to our support assumptions for $v_e$).

\textbf{Proof of \eqref{it.intmini_empty_line}.} This immediately follows from (\ref{eq.associated_int_op}) and the support assumptions for $x_e$ and $v_e$.
\end{proof}

\begin{proposition}[Better upper bounds for one-atom $\Ic_\Mb$]\label{prop.int_mini_good}Consider the same setting (with $\Ic_\Mb$ for an one-atom molecule $\Mb$ etc.) as in Proposition \ref{prop.intmini} above, but with extra restrictions on the support of $Q$, which are specified below.
\begin{enumerate}
    \item\label{it.bound_one_atom_good_1} Suppose $\Mb$ is a \{3\} molecule (by symmetry we may assume $e_1$ is fixed and $(e_2,e_3,e_4)$ are free), and we further restrict the support of $Q$ to be a subset of
        \begin{equation}\label{eq.intmini_new}
        |t_1-t^*|\leq \varepsilon^{\upsilon}\quad\mathrm{or}\quad |v_i-v_j|\leq\varepsilon^{\upsilon}\quad\mathrm{or}\quad |x_j-x^*|\leq\varepsilon^{1-\upsilon}\quad\mathrm{or}\quad|x_i-x_j|\leq \varepsilon^{1-\upsilon}
    \end{equation} for some $0<\upsilon<1/2$. Here $(x^*,v^*,t^*)$ are some external parameters (it is possible that $x^*$ or $v^*$ coincide with $x_1$ or $v_1$), and the support of $Q$ may depend on these parameters. Moreover, in the second and fourth cases in \eqref{eq.intmini_new}, we require $i\neq j$, and also $\{i,j\}\neq \{1,3\},\{2,4\}$ when $\Mb$ is an O-atom. In the third case we require $j\neq 1$, and also $j\neq 3$ when $\Mb$ is an O-atom.
    
    Then we have the following estimate (uniform in all parameters):
    \begin{equation}\label{eq.int_elem_est_deg3_degen}
        \|\Ic_\Mb(Q) \|_{L^\infty}\leq \varepsilon^{\upsilon}|\log\varepsilon|^{C^*}\cdot\|Q\|_{L^\infty}.
    \end{equation}
    
    \item\label{it.bound_one_atom_good_2} 
    Suppose $\Mb$ is a \{4\} molecule, and we further restrict the support of $Q$ to be a subset of
        \begin{equation}\label{eq.intmini_6+}
    |x_i-x_j|\leq\varepsilon^{1-\upsilon}\quad\mathrm{or}\quad |x_i-x^*|\leq \varepsilon^{1-\upsilon}\quad\mathrm{or}\quad |v_i-v_j|\leq\varepsilon^{\upsilon} 
    \end{equation}
    for some $0<\upsilon<1/2$. Here $x^*$ is an external parameter as in (\ref{it.bound_one_atom_good_1}). Moreover we require $i\neq j$, and also $\{i,j\}\neq\{1,3\},\{2,4\}$ when $\Mb$ is an O-atom. Then we have the following estimate (uniform in all parameters):
    \begin{equation}\label{eq.int_elem_est_deg4_degen}
        \|\Ic_\Mb(Q) \|_{L^\infty}\leq \varepsilon^{-(d-1)+\upsilon}|\log\varepsilon|^{C^*}\cdot\|Q\|_{L^\infty}
    \end{equation}
    \item\label{it.bound_one_atom_good_empty_line} Suppose $\Mb$ is an empty end with an edge $e_1$ and associated variable $z_1 = (x_1, v_1)$, and we replace the $Q$ in (\ref{eq.associated_int_op}) by $\mathbbm{1}_*^{\varepsilon}\cdot Q$, where
    \begin{equation}\label{eq.intmini_empty_line_good}
   \mathbbm{1}_*^{\varepsilon} = \left\{\begin{aligned}
        &1 \quad &&|x_1-x^*|\leq \varepsilon^{1-\upsilon},
        \\
        &\varepsilon^\upsilon \quad &&|x_1-x^*|\geq \varepsilon^{1-\upsilon},
    \end{aligned}\right.
    \end{equation}
    for some $0<\upsilon<1/2$. Here $x^*$ is an external parameter as in (\ref{it.bound_one_atom_good_1}). Then we have the following estimate (uniform in $x^*$)
    \begin{equation}\label{eq.int_elem_est_empty_line_degen}
        \|\Ic_\Mb (\mathbbm{1}_*^{\varepsilon}\cdot Q) \|_{L^\infty}\leq \varepsilon^{-(d-1)+\upsilon}|\log\varepsilon|^{C^*}\cdot\|Q\|_{L^\infty}.
    \end{equation}
    \item\label{it.bound_one_atom_good_link} Suppose $\Mb$ is \{3\}  or \{4\} molecules, and we \emph{do not} make any restriction on the support of $Q$ as in (\ref{it.bound_one_atom_good_1})--(\ref{it.bound_one_atom_good_2}), but \emph{replace} $\Ic_\Mb(Q)$ with $\Ic_\Mb(\mathbbm{1}_*^{\varepsilon} Q)$ where     \begin{equation}
        \mathbbm{1}_* = \left\{\begin{aligned}
            &1 \quad &&|x-x'|\le \varepsilon^{1-\upsilon},
            \\
            &\varepsilon^\upsilon \quad &&|x-x'|\ge \varepsilon^{1-\upsilon}
        \end{aligned}\right.
    \end{equation} for some $0<\upsilon<1/2$. Here $(x, x')$ is either $(x_i, x^*)$ or $(x_i, x_j)$, $x^*$ is an external parameter, and the indices $i$ and $j$ also satisfy the requirements in \eqref{it.bound_one_atom_good_1} and \eqref{it.bound_one_atom_good_2} respectively. Then the estimates \eqref{eq.int_elem_est_deg3_degen} and \eqref{eq.int_elem_est_deg4_degen} hold true.
\end{enumerate}
\end{proposition}
\begin{proof} We may assume $\|Q\|_{L^\infty}=1$. We divide the proof into several parts.

\textbf{Proof of \eqref{it.bound_one_atom_good_1} and \eqref{it.bound_one_atom_good_2}.} Recall that the expression of $\Ic_\Mb$ is given by \eqref{eq.intmini_3} and \eqref{eq.intmini_6}, where the variables $(x_i, v_i)$ are given by \eqref{eq.delta_supp}. Consider the conditions \eqref{eq.intmini_new} and \eqref{eq.intmini_6+}. The case $|t_1-t^*|\leq\varepsilon^{\upsilon}$ in \eqref{eq.intmini_new} is trivial as we gain $\varepsilon^{\upsilon}$ in the integral in $t_1$ in $\Ic_\Mb$, which proves \eqref{eq.int_elem_est_deg3_degen}. Next assume $|v_i-v_j|\leq\varepsilon^{\upsilon}$ in \eqref{eq.intmini_new} and $\Mb$ is a C-atom (the O-atom case is the same as long as $\{i,j\}\neq \{1,3\},\{2,4\}$); let $v_1-v_2:=u$ and $u-(u\cdot\omega)\omega:=u^{\perp}$, then we have 
\begin{equation}\label{eq.intmini_proof_new}
    \textrm{either $|u|\leq\varepsilon^{\upsilon}$, or $|u\cdot\omega|\leq \varepsilon^{\upsilon}$, or $|u^\perp|\leq \varepsilon^{\upsilon}$, or $|u-2(u\cdot\omega)\omega|\leq \varepsilon^{\upsilon}$.}
\end{equation} 
Note that $|u-2(u\cdot\omega)\omega|=|u|$, and $|u|^2=|u\cdot\omega|^2+|u^\perp|^2$, we may reduce to one of the cases $|u\cdot\omega|\leq \varepsilon^{\upsilon}$, or $|u^\perp|\leq \varepsilon^{\upsilon}$. With fixed $\omega$ (say $\omega=(1,0,\cdots,0)$ by symmetry), in either case, at least one coordinate of $u$ is bounded by $\varepsilon^{\upsilon}$, and all the other coordinates (as well as the weight $((v_1-v_2)\cdot\omega)_-$) are bounded by $|\log\varepsilon|^{C^*}$, which implies (\ref{eq.int_elem_est_deg3_degen}) by this volume gain.

Next, assume $|x_j-x^*|\leq\varepsilon^{1-\upsilon}$ (or similarly $|x_i-x_j|\leq\varepsilon^{1-\upsilon}$) in \eqref{eq.intmini_new}, then since the distance between any two $x_j+t_1v_j$ is at most $\varepsilon$, we must have $|(x_1-x^*)+t_1(v_1-v_j)|\leq 2\varepsilon^{1-\upsilon}$; by dyadically decomposing $|t_1|$ and gaining from both the support of $t_1$ and the restriction on $v_1-v_j$ as above, we can obtain a volume gain which is $\varepsilon^{\upsilon}$, and this is enough to prove \eqref{eq.int_elem_est_deg3_degen}. The cases of \eqref{eq.intmini_6+} involving $x_i-x_j$ and $v_i-v_j$ follow similarly, and the case involving $x_i-x^*$ is trivial as we gain from integrating in $x_1$.

\textbf{Proof of \eqref{it.bound_one_atom_good_empty_line} and \eqref{it.bound_one_atom_good_link}.} The bound in \eqref{it.bound_one_atom_good_empty_line} is immediate (we gain from the $x_1$ integral in (\ref{eq.intmini_empty_line}) if $|x_1-x^*|\leq\varepsilon^{1-\upsilon}$). For \eqref{it.bound_one_atom_good_link}, the case $|x - x'| \geq \varepsilon^{1-\upsilon}$ is also immediate. For the case $|x - x'| \le \varepsilon^{1-\upsilon}$, the estimate follows by the same argument as in the $|x_j - x^*| \leq \varepsilon^{1 - \upsilon}$ case from the proofs of \eqref{it.bound_one_atom_good_1} and \eqref{it.bound_one_atom_good_2}.
\end{proof}

\begin{proposition}[Better upper bounds for two-atom $\Ic_\Mb$]\label{prop.intmini_2} Let $\Ic_\Mb$ be defined as in \eqref{eq.associated_int_op}, where $\Mb$ is a regular molecule containing only two atoms $\nf_1, \nf_2$ with 7 edges $e_1,\dots,e_7$, see {\color{blue}Figure \ref{fig.elementary}.} We also denote by $e$ the edge between $\nf_1$ and $\nf_2$, which equals some $e_j$. Let $Q$ be a function supported in $|v_j|\leq|\log\varepsilon|^{C^*}$ and $|x_j-x_j^0|\leq |\log\varepsilon|^{C^*}$, and also $\vt_\Mb\in\Dc$, as in Proposition \ref{prop.intmini}. Consider the following cases:
\begin{enumerate}
\item\label{it.intmini_4} Suppose $\Mb$ is a \{33A\} molecule with atoms $(\nf_1,\nf_2)$. Let $e$ be the bond connecting them and $(e_1,e_7)$ be the two fixed ends at $\nf_1$ and $\nf_2$ respectively. Assume moreover that either (i) $\ell[\nf_1]\neq\ell[\nf_2]$ for the two atoms $\nf_1$ and $\nf_2$, or (ii) the support of $Q$ is a subset of
\begin{equation}\label{eq.intmini_8-}
\max(|x_1-x_7|,|v_1-v_7|)\geq \varepsilon^{1-\upsilon}\quad \mathrm{or}\quad |t_1-t_2|\geq\varepsilon^{\upsilon}
\end{equation} 
for some $0<\upsilon<1/3$.
\item\label{it.intmini_5} Suppose $\Mb$ is a \{33B\} molecule with two C-atoms. Assume the support of $Q$ is a subset of 
\begin{equation}\label{eq.intmini_9-}
|t_1-t_2|\geq \varepsilon^{\upsilon} \quad\text{and}\quad |v_i-v_j|\geq \varepsilon^{\upsilon}\quad\textrm{for\ any\ edge\ }(e_i,e_j)\textrm{\ at\ the\ same\ atom},
\end{equation}
for some $0<\upsilon<1/(6d)$.
\item\label{it.intmini_6} Suppose $\Mb$ is a \{44\} molecule with atoms $(\nf_1,\nf_2)$ and $(e_1,e_7)$ being two free ends at $\nf_1$ and $\nf_2$ respectively, such that $\Mb$ becomes a \{33A\} molecule after turning these two free ends into fixed ends. Moreover assume the support if $Q$ is a subset of 
\begin{equation}\label{eq.intmini_10-}
\max(|x_1-x_7|,|v_1-v_7|)\leq \varepsilon^{1-\upsilon}\quad\textrm{and}\quad |t_1-t_2|\leq\varepsilon^{\upsilon}.
\end{equation}
for some $0<\upsilon<1/(4d)$. 
\end{enumerate}
Then, in any of the above cases (\ref{it.intmini_4})--(\ref{it.intmini_6}), we have 
\begin{equation}\label{eq.intmini_8}
    \|\Ic_\Mb(Q) \|_{L^\infty}\leq \varepsilon^{\upsilon}|\log\varepsilon|^{C^*}\cdot\|Q\|_{L^\infty}.
\end{equation} 
\end{proposition}
\begin{proof} We may assume $\|Q\|_{L^\infty}=1$. We divide the proof into several parts.

\textbf{Proof of \eqref{it.intmini_4}.} Since $\Mb$ is a \{33A\} molecule (Definition \ref{def.elementary}), by symmetry, we may assume that $\{\nf_2\}$ becomes a \{2\} molecule after cutting $\nf_1$ as free. First assume (ii) is true. To calculate $\Ic_\Mb(Q)$, we cut $\nf_1$ as free to obtain $\Ic_\Mb(Q) = \Ic_{\{\nf_1\}}\circ\Ic_{\{\nf_2\}}(Q)$. By assumption, we know that $\nf_2$ becomes a deg 2 atom with two top fixed ends. The inner integral $\Ic_{\{\nf_2\}}$ can be computed using \eqref{eq.intmini_1} in Proposition \ref{prop.intmini}.

Let the bond between $\nf_1$ and $\nf_2$ be $e$,  denote the edges at $\nf_1$ by $(e_1.\cdots,e_4)$ as in Proposition \ref{prop.intmini}, then $e=e_j$ for some $j\in\{2,3,4\}$ and $j\neq 3$ if $\nf_1$ is an O-atom (because $\Mb$ is regular, a fixed end cannot be serial with a bond at O-atom). After plugging in the above formula for the inner integral $\Ic_{\{\nf_2\}}$, we can reduce $\Ic_\Mb(Q)$ to an integral of form $\Ic_{\{\nf_1\}}$ as described in \eqref{it.intmini_2} of Proposition \ref{prop.intmini}, namely
\begin{equation}\label{eq.intmini_proof_new2}
    \Ic_\Mb(Q)=\int_{\Rb\times\Sb^{d-1}\times\Rb^d}\big[(v_1-v_2)\cdot\omega\big]_-\cdot \mathbbm{1}_{\mathrm{col}}(z_j,z_7)\cdot Q\,\mathrm{d}t_1\mathrm{d}\omega\mathrm{d}v_2,
\end{equation} 
where the input variables of $Q$ are explicit functions of $(t_1,\omega,v_2)$ and $(z_1,z_7)$ using Proposition \ref{prop.intmini} and the fact that $z_j$ satisfies \eqref{eq.delta_supp}. To prove \eqref{eq.intmini_8} for \eqref{it.intmini_4}, it suffices to control the volume of the support of $\mathbbm{1}_{\mathrm{col}}$ (which is a set depending on $z_1$ and $z_7$).

By definition \eqref{eq.intmini_1+} of $\mathbbm{1}_{\mathrm{col}}$, we know this restriction implies that \begin{equation}\label{eq.intmini_proof_5}|(x_j-x_7)\wedge(v_j-v_7)|\leq\varepsilon|\log\varepsilon|^{C^*}
\end{equation} for some $2\leq j\leq 4$, where $j\neq 3$ if $\nf$ is O-atom. Let $v_1-v_2=u$ and $u-(u\cdot\omega)\omega=u^\perp$. Plugging in the formula \eqref{eq.delta_supp}, and noting that $|x_j-(x_1+t_1u_j)|\leq\varepsilon$ where $u_j=v_1-v_j$ and $(u_2,u_3,u_4)=(u,(u\cdot\omega)\omega,u^\perp)$ if $\nf_1$ is a C-atom and $(u_2,u_4)=(u,u)$ if $\nf_1$ is an O-atom, we obtain that
\begin{equation}\label{eq.intmini_proof_6}\big|\big(x_1-x_7+t_1(v_1-v_7)\big)\wedge u_j-\boldsymbol{g}\big|\leq\varepsilon|\log\varepsilon|^{C^*}
\end{equation} for some constant 2-form $\boldsymbol{g}$ depending only on $(x_1,v_1,x_7,v_7)$. 

In \eqref{eq.intmini_proof_6}, by dyadic decomposition, we may assume $|x_1-x_7+t_1(v_1-v_7)|\sim\varepsilon^\beta$ (say with $0<\beta<1$, the endpoint cases are treated in the standard way), then with fixed $(x_1,v_1,x_7,v_7)$ and $t_1$, the vector $u_j$ must belong to a one-dimensional tube in $\Rb^d$ which has length $\leq |\log\varepsilon|^{C^*}$ and width $\varepsilon^{1-\beta}|\log\varepsilon|^{C^*}$. Note that $u_j\in\{u,(u\cdot\omega)\omega,u-(u\cdot\omega)\omega\}$, by elementary calculations, we see that this restriction on $u_j$ implies that $(u,\omega)$ must belong to a subset of $\Rb^d\times\Sb^{d-1}$ with volume $\leq \varepsilon^{1-\beta}|\log\varepsilon|^{C^*}$. This already implies \eqref{eq.intmini_8}, if $\beta\leq 1-\upsilon$.

Now suppose \[|x_1-x_7+t_1(v_1-v_7)|\sim\varepsilon^{\beta},\quad\beta>1-\upsilon,\] we shall exploit the $t_1$ integral in \eqref{eq.intmini_proof_new2}. If $|v_1-v_7|\gtrsim \varepsilon^{1-\upsilon}$, then $t_1$ belongs to an interval of length $\leq \varepsilon^{\beta-1+\upsilon}$, which then implies \eqref{eq.intmini_8} together with the above volume bound for $(u,\omega)$. If $|v_1-v_7|\ll \varepsilon^{1-\upsilon}$, then we must also have $|x_1-x_7|\ll\varepsilon^{1-\upsilon}$, while by the support of $\mathbbm{1}_{\mathrm{col}}(z_j,z_7)$ we also get
\begin{equation}\label{eq.intmini_proof_7}
|x_j-x_7+t_2(v_j-v_7)|=\varepsilon,\quad |x_j-x_1+t_1(v_j-v_1)|\in\{0,\varepsilon\},
\end{equation} which implies $|t_2-t_1|\cdot |v_j-v_1|\lesssim_* \varepsilon^{1-\upsilon}$ given that $(x_7,v_7)$ is close to $(x_1,v_1)$. Since now we must have $|t_2-t_1|\geq \varepsilon^{\upsilon}$ by \eqref{eq.intmini_8-}, we conclude that $|v_j-v_1|\lesssim_* \varepsilon^{1-2\upsilon}$ which again implies \eqref{eq.intmini_8} as in the proof of Proposition \ref{prop.intmini} \eqref{it.intmini_2}.

Finally, assume (i) is true but (ii) is not, by the above proof, we may assume that $Q$ is supported in the set where $|t_1-t_2|\leq \varepsilon^{\upsilon}$. But since $t_1$ and $t_2$ must belong to disjoint time intervals in view of $\ell[\nf_1]\neq\ell[\nf_2]$, we conclude that $t_1$ itself must belong to a fixed interval of length at most $\varepsilon^{\upsilon}$. Then \eqref{eq.intmini_8} follows from exploiting the volume of support of $t_1$. This completes the proof of \eqref{it.intmini_4}.

\textbf{Proof of \eqref{it.intmini_5}.} Now consider the case of \{33B\} molecules. Recall that in this case both atoms $\nf_j$ are $C$ atoms, the higher atom $\nf_1$ has one top fixed end, and the lower atom $\nf_2$ has a bottom fixed end. Concerning the serial relations of edges and up to symmetry, there are 3 cases for the molecule $\Mb$ (see {\color{blue}Figure \ref{fig.elementary}}), which we list below.

\textbf{Case 1.} This is the case when both fixed ends are serial with the common edge $e$. Here
\begin{equation}
\label{eq.intmini_proof_10}
\begin{aligned}
\Ic_\Mb(Q)&=\varepsilon^{-(d-1)}\int_{\Rb^{10d}\times\Rb^2}\dirac\big(x_3-x_1+t_1(v_3-v_1)\big)\cdot\dirac\big(x_4-x_2+t_1(v_4-v_2)\big)\dirac\big(v_3-v_1+[(v_1-v_2)\cdot\omega]\omega\big)\\
&\times\dirac\big(v_4-v_2-[(v_1-v_2)\cdot\omega]\omega\big)\dirac\big(|x_1-x_2+t_1(v_1-v_2)|-\varepsilon\big)\cdot\big[(v_1-v_2)\cdot\omega\big]_+\cdot \dirac\big(x_6-x_5+t_2(v_6-v_5)\big)\\
&\times\dirac\big(x_7-x_3+t_2(v_7-v_3)\big)\dirac\big(v_6-v_5+[(v_5-v_3)\cdot\eta]\eta\big)\dirac\big(v_7-v_3-[(v_5-v_3)\cdot\eta]\eta\big)\\
&\times\dirac\big(|x_5-x_3+t_2(v_5-v_3)|-\varepsilon\big)\cdot\big[(v_5-v_3)\cdot\eta\big]_+\cdot Q(x_1,v_1,\cdots,x_7,v_7,t_1,t_2)\prod_{j=2}^6\mathrm{d}x_j\mathrm{d}v_j\prod_{j=1}^2\mathrm{d}t_j,
\end{aligned}
\end{equation} where $\omega:=\varepsilon^{-1}(x_1-x_2+t_1(v_1-v_2))$ and $\eta:=\varepsilon^{-1}(x_5-x_3+t_2(v_5-v_3))$ are two unit vectors. Here we first integrate in $(x_4,v_4,x_6,v_6)$, then make the change of variables $x_5\leftrightarrow\eta$ to reduce to
\begin{equation}
\label{eq.intmini_proof_11}
\begin{aligned}
\Ic_\Mb(Q)&=\int_{\Rb^{4d}\times\Rb^2}\dirac\big(x_3-x_1+t_1(v_3-v_1)\big)\dirac\big(v_3-v_1+[(v_1-v_2)\cdot\omega]\omega\big)\\
&\times\dirac\big(|x_1-x_2+t_1(v_1-v_2)|-\varepsilon\big)\dirac\big(x_7-x_3+t_2(v_7-v_3)\big)\cdot\big[(v_1-v_2)\cdot\omega\big]_+\prod_{j=2}^3\mathrm{d}x_j\mathrm{d}v_j\prod_{j=1}^2\mathrm{d}t_j\\
&\times\int_{\Sb^{d-1}\times\Rb^d}\dirac\big(v_7-v_3-[(v_5-v_3)\cdot\eta]\eta\big)\cdot \big[(v_5-v_3)\cdot\eta\big]_+\cdot Q(\cdots)\,\mathrm{d}\eta\mathrm{d}v_5,
\end{aligned}
\end{equation} where the exact expressions of the input variables of $Q$ are omitted for simplicity. Note that in \eqref{eq.intmini_proof_11}, the $\dirac$ functions imply that $\eta=(v_7-v_3)/|v_7-v_3|$ and $(v_5-v_3)\cdot\eta=(v_7-v_3)\cdot\eta$, we may decompose $v_5-v_3$ into $[(v_5-v_3)\cdot\eta]\eta$, and a vector $w\in\Pi_\eta^\perp$ which is orthogonal to $\eta$. This leads to
\begin{equation}
\label{eq.intmini_proof_12}
\begin{aligned}
\Ic_\Mb(Q)&=\int_{\Rb^{4d}\times\Rb^2}\dirac\big(x_3-x_1+t_1(v_3-v_1)\big)\dirac\big(v_3-v_1+[(v_1-v_2)\cdot\omega]\omega\big)\cdot\dirac\big(|x_1-x_2+t_1(v_1-v_2)|-\varepsilon\big)\\&\times\big[(v_1-v_2)\cdot\omega\big]_+\cdot\dirac\big(x_7-x_3+t_2(v_7-v_3)\big)|v_7-v_3|^{-(d-2)}\prod_{j=2}^3\mathrm{d}x_j\mathrm{d}v_j\prod_{j=1}^2\mathrm{d}t_j\cdot\int_{\Pi_\eta^\perp}Q(\cdots)\,\mathrm{d}w\mathrm,
\end{aligned}
\end{equation} where $\eta=(v_7-v_3)/|v_7-v_3|$ as above, $\Pi_\eta^\perp$ is the orthogonal complement of $\eta$, and the input variables of $Q$ are again omitted. The factor $|v_7-v_3|^{-(d-2)}$ comes from the factor $\big[(v_5-v_3)\cdot\eta\big]_+$, together with the Jacobian coming from viewing $[(v_5-v_3)\cdot\eta]\eta$ as the polar coordinates of $v_7-v_3$. Now for fixed $(x_3,v_3)$, we make the change of variables $x_2\leftrightarrow \omega$, and argue similarly as above (namely decomposing $v_1-v_2$ into $[(v_1-v_2)\cdot\omega]\omega$ and a vector $\theta\in\Pi_\omega^\perp$ orthogonal to $\omega$), to get
\begin{equation}
\label{eq.intmini_proof_13}
\begin{aligned}
\Ic_\Mb(Q)&=\varepsilon^{d-1}\int_{\Rb^{2d}\times\Rb^2}\dirac\big(x_3-x_1+t_1(v_3-v_1)\big)\dirac\big(x_7-x_3+t_2(v_7-v_3)\big)\\&\times|v_7-v_3|^{-(d-2)}|v_1-v_3|^{-(d-2)}\,\mathrm{d}x_3\mathrm{d}v_3\mathrm{d}t_1\mathrm{d}t_2\cdot\int_{\Pi_\eta^\perp\times\Pi_\omega^\perp}Q(\cdots)\,\mathrm{d}w\mathrm{d}\theta,
\end{aligned}
\end{equation} where $\omega=-(v_3-v_1)/|v_3-v_1|$ and $\Pi_\omega^\perp$ is the orthogonal complement of $\omega$. Finally, by integrating the $\dirac$ functions in \eqref{eq.intmini_proof_13}, first in $x_3$ and then in $v_3$ and then in $(t_1,t_2)$, we get
\begin{equation}
\label{eq.intmini_proof_14}
\Ic_\Mb(Q)=\varepsilon^{d-1}\int_{\Rb^2}|t_1-t_2|^{-d}|v_7-v_3|^{-(d-2)}|v_1-v_3|^{-(d-2)}\,\mathrm{d}t_1\mathrm{d}t_2\int_{\Pi_\eta^\perp\times\Pi_\omega^\perp}Q(\cdots)\,\mathrm{d}w\mathrm{d}\theta,
\end{equation} 
where each of $(v_3,\eta,\omega)$ is a function of $(t_1,t_2)$ and the fixed parameters. The weight $|t_1-t_2|^{-d}|v_7-v_3|^{-(d-2)}|v_1-v_3|^{-(d-2)}$ in $\Ic_\Mb(Q)$ is bounded above by $\varepsilon^{-1/2}$ due to \eqref{eq.intmini_9-} and $0<\upsilon<1/(6d)$, and can be absorbed by $\varepsilon^{d-1}$, and hence \eqref{eq.intmini_8} follows.

\textbf{Case 2.} This is the case when the top fixed edge is serial with $e$, but the bottom fixed edge is not. Here
\begin{equation}
\label{eq.intmini_proof_15}
\begin{aligned}
\Ic_\Mb(Q)&=\varepsilon^{-(d-1)}\int_{\Rb^{10d}\times\Rb^2}\dirac\big(x_3-x_1+t_1(v_3-v_1)\big)\cdot\dirac\big(x_4-x_2+t_1(v_4-v_2)\big)\cdot \dirac\big(v_3-v_1+[(v_1-v_2)\cdot\omega]\omega\big)\\
&\times\dirac\big(v_4-v_2-[(v_1-v_2)\cdot\omega]\omega\big)\dirac\big(|x_1-x_2+t_1(v_1-v_2)|-\varepsilon\big)\cdot\big[(v_1-v_2)\cdot\omega\big]_+\cdot \dirac\big(x_7-x_5+t_2(v_7-v_5)\big)\\
&\times\dirac\big(x_6-x_3+t_2(v_6-v_3)\big)\dirac\big(v_7-v_5+[(v_5-v_3)\cdot\eta]\eta\big)\dirac\big(v_6-v_3-[(v_5-v_3)\cdot\eta]\eta\big)\\
&\times\dirac\big(|x_5-x_3+t_2(v_5-v_3)|-\varepsilon\big)\cdot\big[(v_5-v_3)\cdot\eta\big]_+\cdot Q(x_1,v_1,\cdots,x_7,v_7,t_1,t_2)\prod_{j=2}^6\mathrm{d}x_j\mathrm{d}v_j\prod_{j=1}^2\mathrm{d}t_j,
\end{aligned}
\end{equation} where $\omega$ and $\eta$ are as above. Here we first integrate in $(x_4,v_4)$, make the change of variables $x_2\leftrightarrow\omega$ and argue as in Case 1 above (i.e. decomposing $v_1-v_2$ into $[(v_1-v_2)\cdot\omega]\omega$ and a vector $\theta\in\Pi_\omega^\perp$ orthogonal to $\omega$), to reduce
\begin{equation}
\label{eq.intmini_proof_16}
\begin{aligned}
\Ic_\Mb(Q)&=\int_{\Rb^{6d}\times\Rb^2}\dirac\big(x_7-x_5+t_2(v_7-v_5)\big)\dirac\big(x_6-x_3+t_2(v_6-v_3)\big)\dirac\big(v_7-v_5+[(v_5-v_3)\cdot\eta]\eta\big)\\
&\times\dirac\big(v_6-v_3-[(v_5-v_3)\cdot\eta]\eta\big)\dirac\big(|x_5-x_3+t_2(v_5-v_3)|-\varepsilon\big)\cdot \big[(v_5-v_3)\cdot\eta\big]_+\\&\times  
\dirac\big(x_3-x_1+t_1(v_3-v_1)\big)\cdot|v_1-v_3|^{-(d-2)}\prod_{j\in\{3,5,6\}}\mathrm{d}x_j\mathrm{d}v_j\prod_{j=1}^2\mathrm{d}t_j\cdot\int_{\Pi_\omega^\perp} Q(\cdots)\,\mathrm{d}\theta,
\end{aligned}
\end{equation} where $\omega=-(v_3-v_1)/|v_3-v_1|$ and $\Pi_\omega^\perp$ is the orthogonal complement of $\omega$. Now in \eqref{eq.intmini_proof_16}, we argue as in \eqref{it.intmini_2} of Proposition \ref{prop.intmini}, namely make the change of variables $x_6\leftrightarrow\eta$, noting that
\[\varepsilon\eta=x_5-x_3+t_2(v_5-v_3)=x_7-x_6+t_2(v_7-v_6)\] as well as $(v_5-v_3)\cdot\eta=(v_6-v_7)\cdot\eta$, to reduce it to
\begin{equation}
\label{eq.intmini_proof_17}
\Ic_\Mb(Q)=\varepsilon^{d-1}\int_{\Rb^{d}\times\Sb^{d-1}\times\Rb^2}\big[(v_6-v_7)\cdot\eta\big]_+ 
\,\mathrm{d}v_6\mathrm{d}\eta\mathrm{d}t_1\mathrm{d}t_2 \cdot \dirac\big(x_3-x_1+t_1(v_3-v_1)\big)\cdot|v_1-v_3|^{-(d-2)}\cdot\int_{\Pi_\omega^\perp} Q(\cdots)\,\mathrm{d}\theta,
\end{equation}
where $(x_3,v_3)$ is determined by $(v_6,\eta,t_1,t_2)$ and the fixed parameters. Note that with \emph{fixed} $(t_1,t_2,\eta)$, we have \[x_3-x_1+t_1(v_3-v_1)=\gamma+(t_1-t_2)\sigma\] where $\gamma$ depends only on $(t_1,t_2,\eta)$ and the fixed parameters, and $\sigma$ is the projection of $v_6-v_7$ to the orthogonal complement $\Pi_\eta^\perp$ of $\eta$. Then the $\dirac$ function in \eqref{eq.intmini_proof_17} can be written as $\dirac_\Rb(\gamma \cdot \eta)\dirac_{\Pi_\eta^{\perp}}(\gamma+(t_1-t_2)  \sigma)$ where $\Pi_\eta^{\perp}$ is the orthogonal complement of $\eta$. The latter $\dirac_{\Pi_\eta^{\perp}}$ yields $|t_1-t_2|^{-(d-1)}$ upon integrating in $\sigma$, and leaves us with the remaining component $\alpha:=(v_6-v_7) \cdot \eta$ from the integral in $v_6$. Note also that 
$$
\begin{aligned}
\gamma\cdot \eta&=\eta\cdot \big[(x_3+t_1v_3)-(x_1+t_1v_1)\big]=\eta\cdot \big[(x_3+t_2v_3)-(x_1+t_1v_1)+(t_1-t_2)v_3\big]\\&=\eta\cdot \big[(x_7+t_2v_7)-\varepsilon \eta-(x_1+t_1v_1)+(t_1-t_2)v_3\big]=\eta\cdot \zeta-\varepsilon,
\end{aligned}
$$
where $\zeta:=(x_7+t_1v_7)-(x_1+t_1v_1)$, and we have used that $v_3\cdot \eta=v_7\cdot \eta$. We also have $|\zeta|\geq \frac12 \varepsilon^{2\upsilon}$ since 
\begin{multline}
(x_7+t_1v_7)-(x_1+t_1v_1)=(x_7+t_2v_7)-(x_1+t_1v_1)+(t_1-t_2)v_7\\=(x_3+t_2v_3) +\varepsilon \eta-(x_3+t_1v_3) +(t_1-t_2)v_7=(t_1-t_2)(v_7-v_3)+\varepsilon \eta
\end{multline}
due to \eqref{eq.intmini_9-}. As a result, we now have that
\begin{equation}
\label{eq.intmini_proof_18}
\Ic_\Mb(Q)=\varepsilon^{d-1}\int_{\Rb^\times\Sb^{d-1}\times\Rb^2}\max(\alpha,0)\cdot|t_1-t_2|^{-(d-1)}|v_1-v_3|^{-(d-2)} \dirac(\eta\cdot \zeta-\varepsilon)\,\mathrm{d}\alpha\mathrm{d}\eta\mathrm{d}t_1\mathrm{d}t_2\cdot\int_{\Pi_\omega^\perp} Q(\cdots)\,\mathrm{d}\theta,
\end{equation} where $v_3$ and $\omega$ are determined by $(\alpha,\eta,t_1,t_2)$ and the fixed parameters. The integral in $\eta$ can be performed by decomposing $\eta$ into a vector $\eta^\circ$ parallel to $\zeta$ and a vector $\eta^\perp$ orthogonal to $\zeta$. The part of integrand involving $\eta^\circ$ is $|\xi|^{-1}\dirac_\Rb(\eta^\circ-\frac{\varepsilon\zeta}{|\zeta|^2})$, and the integral in $\eta^\perp$ can be reparametrized into an integral on the unit sphere $\Sb^{d-2}$ in the orthogonal complement of $\zeta$ up to a multiplicative factor of $(1-\varepsilon |\zeta|^{-1})^{d-2}\sim 1$. Thus
\begin{equation}
\label{eq.intmini_proof_18+}
\Ic_\Mb(Q)=\varepsilon^{d-1}\int_{\Rb^\times\Sb^{d-2}\times\Rb^2}\max(\alpha,0)|t_1-t_2|^{-(d-1)}|v_1-v_3|^{-(d-2)}\frac{|\zeta|^{-1}}{(1-\varepsilon |\zeta|^{-1})^{d-2}}\,\mathrm{d}\alpha\mathrm{d}\eta^\perp \mathrm{d}t_1\mathrm{d}t_2\int_{\Pi_\omega^\perp} Q(\cdots)\,\mathrm{d}\theta,
\end{equation} 
where $\eta$ is determined by $(t_1, t_2, \eta^\perp)$ and the fixed parameters, and $v_3$ and $\omega$ are determined by $(\alpha,\eta,t_1,t_2)$ and the fixed parameters. The weight $|t_1-t_2|^{-(d-1)}|v_1-v_3|^{-(d-2)}|\zeta|^{-1}$ in $\Ic_\Mb(Q)$ are bounded by $\varepsilon^{-2/3}$ due to \eqref{eq.intmini_9-} and can be absorbed by $\varepsilon^{d-1}$, and hence \eqref{eq.intmini_8} follows.

\textbf{Case 3.} This is the case when neither fixed end is serial with $e$. Here, we define the unit vectors $\omega$ and $\eta$ similar to the above, but instead of substituting, we shall keep both $(\omega,\eta)$ and the original $(x_j,v_j)$ variables in the same integral, but with additional $\dirac$ factors. This leads to the integral (note that the $\varepsilon^{-(d-1)}$ is turned into $\varepsilon^{d-1}$ due to introducing the two $\omega$ and $\eta$ variables)
\begin{equation}
\label{eq.intmini_proof_19}
\begin{aligned}
\Ic_\Mb(Q)&=\varepsilon^{d-1}\int_{\Rb^{10d}\times(\Sb^{d-1})^2\times\Rb^2}\dirac\big(x_3-x_1+t_1(v_3-v_1)\big)\cdot\dirac\big(x_4-x_2+t_1(v_4-v_2)\big)\\
&\times\dirac\big(v_3-v_1+[(v_1-v_2)\cdot\omega]\omega\big)\dirac\big(v_4-v_2-[(v_1-v_2)\cdot\omega]\omega\big)\cdot \dirac\big(x_6-x_4+t_2(v_6-v_4)\big)\\
&\times\dirac\big(x_7-x_5+t_2(v_7-v_5)\big)\dirac\big(v_6-v_4+[(v_4-v_5)\cdot\eta]\eta\big)\dirac\big(v_7-v_5-[(v_4-v_5)\cdot\eta]\eta\big)\\
&\times\dirac\big(x_1-x_2+t_1(v_1-v_2)-\varepsilon\omega\big)\dirac\big(x_4-x_5+t_2(v_4-v_5)-\varepsilon\eta\big)\\
&\times \big[(v_1-v_2)\cdot\omega\big]_+\big[(v_4-v_5)\cdot\eta\big]_+\cdot Q(x_1,v_1,\cdots,x_7,v_7,t_1,t_2)\prod_{j=2}^6\mathrm{d}x_j\mathrm{d}v_j\cdot\mathrm{d}\omega\mathrm{d}\eta\cdot\prod_{j=1}^2\mathrm{d}t_j.
\end{aligned}
\end{equation} Now let $(t_1,t_2)$ be fixed, then $y^*:=x_1+t_1v_1$ and $z^*:=x_7+t_2v_7$ are also fixed, and the $\dirac$ functions occurring in \eqref{eq.intmini_proof_19} that involve $x_j$ (there are six of them) can be rearranged as
\begin{equation}\label{eq.intmini_proof_20}
\begin{aligned}
&\dirac\big(x_2+t_1v_2-(y^*-\varepsilon\omega)\big)\cdot\dirac\big(x_3+t_1v_3-y^*\big)\cdot\dirac\big(x_4+t_1v_4-(y^*-\varepsilon\omega)\big)\\
&\times\dirac\big(x_4+t_2v_4-(z^*+\varepsilon\eta)\big)\dirac\big(x_5+t_2v_5-z^*\big)\dirac\big(x_6+t_2v_6-(z^*+\varepsilon\eta)\big).
\end{aligned}
\end{equation} Therefore, once $\omega, \eta$ and $v_j$ are fixed for some $j\in\{2,3,4,5,6\}$, the corresponding $\dirac$ factor will take care of the integration in $x_j$, except that the two $\dirac$ factors involving $x_4$ in \eqref{eq.intmini_proof_20} will lead to the new $\dirac$ factor
    $$\dirac((t_2-t_1) v_4 +y^*-z^*-\varepsilon(\omega+\eta)).$$
    
Now, suppose $(\omega,\eta, v_4)$ are fixed, we shall integrate the remaining four $\dirac$ functions in \eqref{eq.intmini_proof_19} involving only $v_j$, in the variables $(v_2,v_3,v_5,v_6)$, as follows. The integration in $(v_3,v_6)$ result in substitution, and for the integration in $(v_2,v_5)$ we note that
\begin{equation}\label{eq.intmini_proof_22}
\int_{\Rb^d}\dirac\big(v_4-v_2-[(v_1-v_2)\cdot\omega]\omega\big)\cdot \widetilde{Q}(v_2)\,\mathrm{d}v_2=\dirac\big((v_4-v_1)\cdot\omega\big)\cdot\int_\Rb \widetilde{Q}(v_4+\alpha\omega)\,\mathrm{d}\alpha
\end{equation} for any function $\widetilde{Q}$, which is easily proved by assuming $\omega$ is a coordinate vector. The same holds for the integration in $v_5$, and putting them together, we can reduce
\begin{multline}\label{eq.intmini_proof_23}
\Ic_\Mb(Q)=\varepsilon^{d-1}\int_{(\Sb^{d-1})^2\times\Rb^2\times\Rb^2}\max(\alpha_1,0)\max(\alpha_2,0)\int_{\Rb^d} \dirac((t_2-t_1) v_4 +y^*-z^*-\varepsilon(\omega+\eta))\\
\times \dirac\big((v_4-v_1)\cdot\omega\big)\dirac\big((v_4-v_7)\cdot\eta\big)\cdot Q(\cdots)\,\mathrm{d}v_4\cdot \mathrm{d}\omega\mathrm{d}\eta \prod_{j=1}^2\mathrm{d}t_j\mathrm{d}\alpha_j.
\end{multline}

Next, with $v_4$ fixed, we integrate in $\eta$ and $\omega$ using the formula 
$$
\int_{\Sb^{d-1}}\dirac(X\cdot\omega) f(\omega) \mathrm{d}\omega=|X|^{-1}\int_{\Sb^{d-2}} f(R_X\omega')\,\mathrm{d}\omega', \qquad X\in \{v_4-v_1, v_4-v_7\}
$$
where $\Sb^{d-2}$ is the unit sphere in the orthogonal complement of a fixed coordinate vector $e$, and $R_X$ is the rotation such that $R(e)=\frac{X}{|X|}$. This follows from applying the spherical coordinates. Therefore, we have
\begin{multline}\label{eq.intmini_proof_23+}
\Ic_\Mb(Q)=\varepsilon^{d-1}\int_{(\Sb^{d-2})^2\times\Rb^2\times\Rb^2}|t_1-t_2|^{-d}|v_4-v_1|^{-1}|v_4-v_7|^{-1}\max(\alpha_1,0)\max(\alpha_2,0)\\\times\int_{\Rb^d} \dirac\bigg( v_4 +\frac{y^*-z^*-\varepsilon(R_{v_4-v_1}\omega'+R_{v_4-v_7}\eta')}{t_2-t_1}\bigg)
\cdot Q(\cdots)\,\mathrm{d}\omega'\mathrm{d}\eta' \mathrm{d}v_4\prod_{j=1}^2\mathrm{d}t_j\mathrm{d}\alpha_j.
\end{multline}
Now, with $(t_1, t_2,\omega',\eta')$ fixed, we evaluate the integral in $v_4$ by first changing variables into $v_4\to \xi$ where $\xi$ is the argument of the $\dirac$. Since $\varepsilon^{1/2} \ll |t_2-t_1|$ by assumption, the Jacobian of this change of variables is close to $1$, so \eqref{eq.intmini_8} follows, using the lower bounds in \eqref{eq.intmini_9-}.

\textbf{Proof of \eqref{it.intmini_6}.} This is easy and follows from \eqref{it.intmini_4}: simply note that $(x_1,v_1,x_7,v_7)$ belongs to a set of volume $\varepsilon^{2d(1-\upsilon)}|\log\varepsilon|^{C^*}$ with $0<\upsilon<1/(4d)$. Once $(x_1,v_1,x_7,v_7)$ is fixed, we can perform the the integration in the remaining variables just as in the \{33A\} molecule case in \eqref{it.intmini_4}. This leads to the desired bound \eqref{eq.intmini_8}.
\end{proof}

Now we define the normal, good, and bad components, based on Propositions \ref{prop.intmini}--\ref{prop.intmini_2}. The cutting algorithm aims to produce good components while avoiding bad ones. These properties depend not only on the molecule $\Mb$, but also on the restriction on the support of $Q$, which comes from splitting (Definition \ref{def.splitting}).

\begin{definition}[Normal, bad and good molecules]
\label{def.good_normal} Suppose $\Mb'$ is formed from a full molecule $\Mb$ by an operation sequence (Definition \ref{def.cutting_algorithm}) and picking one sub-case at each splitting, and $\Mb'$ contains only elementary molecules. The initial link structure $\Lc$ between the bottom ends of $\Mb'$ are inherited from that of $\Mb$ (as the operations do not affect any bottom end). Note also that in the splitting steps, by selecting a sub-case each time, we may have made certain extra restrictions to the support of $Q$ in (\ref{eq.associated_int_op}).

Based on all this information, we now define the notions of good, normal and bad molecules. Below, for any one-atom molecule $\{\nf_1\}$, we denote the 4 edges by $e_1,\dots,e_4$, where $e_1$ is serial with $e_3$, and $e_2$ is serial with $e_4$. For two-atom molecule $\{\nf_1,\nf_2\}$, we denote the 7 edges by $e_1,\cdots,e_7$ with $e_1$ and $e_7$ being the fixed ends in the case of \{33\} molecules, see {\color{blue} Figure \ref{fig.elementary}}. Recall also the constant $\upsilon:=3^{-d-1}$.
\begin{enumerate}
    \item The following molecules are \newterm{normal}:
    \begin{enumerate}
        \item Any \{2\} molecule.
        \item Any \{3\} molecule, except those described in \eqref{it.good_1} below.
        \item\label{it.normal_33A} Any \{33A\} molecule, except those described in \eqref{it.good_2} below.
    \end{enumerate}
    
    \item The following molecules are \newterm{bad}:
    \begin{enumerate}
        \item Any \{4\} molecule.
        \item By convention, we also view any empty end as bad.
    \end{enumerate}
    
    \item The following molecules are \newterm{good}:
    \begin{enumerate}
    \item\label{it.good_1} Any \{3\} molecule $\Xc=\{\nf\}$ with fixed end $e_1$, where either the extra restriction includes one of the following conditions:
    \begin{equation}\label{eq.good_normal_1}
    |t_\nf-t_{\nf'}|\leq \varepsilon^{\upsilon}\quad\mathrm{or}\quad|v_{e_i}-v_{e_j}|\leq\varepsilon^{\upsilon}\quad\mathrm{or}\quad|x_e-x_{e'}|\leq\varepsilon^{1-\upsilon}\quad\mathrm{or}\quad|x_{e_i}-x_{e_j}|\leq\varepsilon^{1-\upsilon}
    \end{equation}
    or $(e_i,e_j)\in\Lc$, or $(e,e')\in\Lc$ (i.e. the relevant pair of ends are initial linked as in Definition \ref{def.molecule}; see also the $\mathbbm{1}_{\Lc}^{\varepsilon}$ factor in Definition \ref{def.associated_op}).
    
     Here in \eqref{eq.good_normal_1} we assume that (i) $\nf'$ is an atom in a molecule $\Yc$ with $\Yc\prec_{\mathrm{cut}}\Xc$ (Definition \ref{def.ordered_disjoint}), and $e'$ is a free end or bond at $\nf'$; (ii) $e$ is a free end at $\nf$ and $e_i\neq e_j$ are two ends at $\nf$; (iii) if $\nf$ is O-atom, then $e'$ is not serial with $e_1$, and $e_i$ not serial with $e_j$. In this case, $t_{\nf'}$ and $x_{e'}$ are fixed variables, so the extra restrictions (or initial link) makes $\Xc$ good.
    \item\label{it.good_2} Any \{33A\} molecule $\Xc=\{\nf_1,\nf_2\}$ with fixed ends $(e_1,e_7)$, where either the extra restriction includes one of the following conditions:
    \begin{equation}\label{eq.good_normal_2}|x_{e_1}-x_{e_7}|\geq\varepsilon^{1-\upsilon}\quad \mathrm{or}\quad |v_{e_1}-v_{e_7}|\geq\varepsilon^{1-\upsilon}\quad\mathrm{or}\quad|t_{\nf_1}-t_{\nf_2}|\geq\varepsilon^{\upsilon},
    \end{equation}
    or $\ell[\nf_1] \neq \ell[\nf_2]$. In this case $\Xc$ is good as a \{33A\} molecule without strong degeneracy.
    \item Any \{33B\} molecule $\Xc=\{\nf_1,\nf_2\}$ of two C-atoms with fixed ends $(e_1,e_7)$, where the extra restriction includes both of the following conditions: 
    \begin{equation}\label{eq.good_normal_4}|v_{e_i}-v_{e_j}|\geq \varepsilon^{\upsilon}\quad\textrm{and}\quad |t_{\nf_1}-t_{\nf_2}|\geq \varepsilon^{\upsilon}.
    \end{equation} 
Here $(e_i,e_j)$ are two different ends at the same atom ($\nf_1$ or $\nf_2$). In this case $\Xc$ is good as a \{33B\} molecule without weak degeneracy.
    \item\label{it.good_44} Any \{44\} molecule $\Xc=\{\nf_1,\nf_2\}$, where the extra restriction includes both of the following conditions: 
    \begin{equation}\label{eq.good_normal_3}
        |x_{e_1}-x_{e_7}|\leq \varepsilon^{1-\upsilon}\quad\mathrm{and}\quad |v_{e_1}-v_{e_7}|\leq \varepsilon^{1-\upsilon}.
    \end{equation} 
    Here $(e_1,e_7)$ are two free ends at $\nf_1$ and $\nf_2$ respectively, such that $\Mb$ becomes a \{33A\} molecule after turning these two free ends into fixed ends. In this case the strong extra restrictions (\ref{eq.good_normal_3}) makes $\Xc$ good.
    \item\label{it.good_3}

Any \{4\} molecule $\Xc=\{\nf\}$, where either the extra restriction includes one of the following conditions:
    \begin{equation}\label{eq.good_normal_4+}
   |v_{e_i}-v_{e_j}|\leq\varepsilon^{\upsilon}\quad \mathrm{or}\quad|x_e-x_{e'}|\leq\varepsilon^{1-\upsilon}\quad\mathrm{or}\quad |x_{e_i}-x_{e_j}|\leq\varepsilon^{1-\upsilon},
    \end{equation}
    or $(e_i,e_j)\in\Lc$, or $(e,e')\in\Lc$ (i.e. the relevant pair of ends are initial linked as in Definition \ref{def.molecule}; see also the $\mathbbm{1}_{\Lc}^{\varepsilon}$ factor in Definition \ref{def.associated_op}).
    
     Here in \eqref{eq.good_normal_4+} we assume that (i) $\nf'$ is an atom in a molecule $\Yc$ with $\Yc\prec_{\mathrm{cut}}\Xc$ (Definition \ref{def.ordered_disjoint}), and $e'$ is a free end or bond at $\nf'$; (ii) $e$ is a free end at $\nf$; (iii) $(e_i,e_j)$ are two different ends at $\nf$, and are non-serial when $\nf$ is O-atom. In this case, $x_{e'}$ is a fixed variable, so the extra restrictions (or initial link) makes $\Xc$ good.
    \item Finally, by convention, we also view an empty end $e$ as good, if $(e,e')\in\Lc$ for some $e'$, i.e. $e$ is involved in some initial link. In this case the initial link condition (see the $\mathbbm{1}_{\Lc}^{\varepsilon}$ factor in Definition \ref{def.associated_op}) will ensure the gain and thus we view $e$ as good.
    \end{enumerate}
\end{enumerate}
\end{definition}
\begin{remark}
\label{rem.elementary} We explain the classification of elementary molecules in Definition \ref{def.good_normal} as follows:
\begin{enumerate}
\item Note that Definition \ref{def.good_normal} allows for some \{4\} molecules and empty ends to be both good and bad; this will not cause any confusion, because in practice the property of good/non-good and bad/non-bad are always exploited independently. The rule of thumb is that, each good molecule leads to a power $\varepsilon^{\upsilon}$ gain in the estimate of the integral $\Ic_\Mb(Q)$ in (\ref{eq.associated_int_op}), and each bad molecule leads to a power $\varepsilon^{-(d-1)}$ loss in this integral (and each normal molecule leads to no gain or loss). These gain and loss estimates easily follow from Definition \ref{def.good_normal} and Propositions \ref{prop.intmini}--\ref{prop.intmini_2}.
\item Note that the definition of good molecules depends on various extra restrictions imposed in Definition \ref{def.good_normal}. In fact, these extra restrictions (corresponding to degeneracy and non-degeneracy conditions) are what motivates the \emph{splitting} steps in the operation sequence. Basically, we may split into degenerate and non-degenerate sub-cases, and exploit different good molecules based on the extra restrictions in each different sub-case. We may also \emph{design different cutting strategies} based on different sub-cases, which will be clear in the algorithm below (cf. the strongly and weakly degeneracies in Propositions \ref{prop.case2} and \ref{prop.comb_est_case4}, where we have a different cutting strategy from the main case).
\end{enumerate}
\end{remark}
\begin{remark}\label{rem.type_elem} Most of the elementary molecules that will occur in the algorithm (Sections \ref{sec.toy}--\ref{sec.maincr}) are \{2\}, \{3\}, \{4\} and \{33A\} molecules. In fact, \{33B\} molecules will occur \emph{only} in the \textbf{MAINUD} algorithm (Definition \ref{def.alg_maincr}), and the \{44\} molecule will occur \emph{only} in \textbf{Choice 1} in the proof of Proposition \ref{prop.case2}. Our construction will guarantee that the \{33B\} and \{44\} molecules are \emph{always good} every time they occur.
\end{remark}
\subsection{Proof of Proposition \ref{prop.cumulant_est}}\label{sec.summary} In this subsection we finish the proof of Proposition \ref{prop.cumulant_est}, assuming Proposition  \ref{prop.comb_est} below, which states that there exists an operation sequence that satisfies the nice property (\ref{eq.overall_alg}). This condition is natural in view of the gains and losses associated with good and bad molecules (see Remark \ref{rem.elementary}; the choice of $\upsilon/2$ instead of $\upsilon$ is due to the need to cover certain log losses).

The construction of the operation sequence in Proposition  \ref{prop.comb_est} will be left to Sections \ref{sec.prepare}--\ref{sec.maincr}. The proof in this subsection relies on Propositions \ref{prop.cumulant_formula}, \ref{prop.layerrec3}, \ref{prop.local_int}, \ref{prop.cutting_algorithm} and \ref{prop.intmini}--\ref{prop.intmini_2}. We also rely on two estimates for $(v_e)$ volume and cross sections in (\ref{eq.associated_int_op}), namely Propositions \ref{prop.weight} and \ref{prop.volume}. The proof of these estimates are technical, and we leave it to Section \ref{sec.aux}.
\begin{proposition}
\label{prop.comb_est} Let $\Mb=(\Mc,\Ec,\Pc,\Lc)$ be a full (unlabeled) molecule as in Definition \ref{def.set_T_F}, which belongs to either $ \mathcal{F}_{\boldsymbol{\Lambda}_\ell}$ or $\mathcal{F}_{\boldsymbol{\Lambda}_\ell}^{\mathrm{err}}$. Let $r(\Mb)=H$ and $Q=Q_\Mb$ as in (\ref{eq.Q_M}). Let also $\rho$ be defined as in (\ref{eq.def_rho_old}) in Definition \ref{def.parameter_rho_old}.

Then, there exist an operation sequence (Definition \ref{def.cutting_algorithm}) which produces at most $C^{|\Mb|}\cdot |\log\varepsilon|^{C^*\rho}$ total sub-cases, such that in each sub-case the $\Mb'$ contains only good, normal and bad elementary molecules (Definition \ref{def.good_normal}). Moreover, if $\Mb\in \mathcal{F}_{\boldsymbol{\Lambda}_\ell}$, we have the following estimate (recall $\upsilon=3^{-d-1}$):
\begin{equation}\label{eq.overall_alg}(\upsilon/2)\cdot\#_{\mathrm{good}}+(d-1)(|H|-\#_{\mathrm{bad}})\geq (C_{13}^*)^{-1}\cdot \rho,
\end{equation} where $\#_{\mathrm{good}}$ and $\#_{\mathrm{bad}}$ are the number of good and bad molecules respectively. If $\Mb\in \mathcal{F}_{\boldsymbol{\Lambda}_\ell}^{\mathrm{err}}$, then (\ref{eq.overall_alg}) holds with an extra term $-C^* |H|$ on the right hand side.
\end{proposition}
\begin{proof} The proof will occupy Sections \ref{sec.prepare}--\ref{sec.maincr}, and will be finished in Section \ref{sec.finish} (after reducing to Proposition \ref{prop.case5} at the end of Section \ref{sec.layer}).
\end{proof}
\begin{proof}[Proof of Proposition \ref{prop.cumulant_est} assuming Proposition \ref{prop.comb_est}] We first make some reductions. By Proposition \ref{prop.cumulant_formula}, we can estimate $|E_H|$ and $|\mathrm{Err}^1|$ by sum of $|\Ic\Nc_\Mb|$ over $\Mb\in\mathcal{F}_{\boldsymbol{\Lambda}_\ell}$ and $\Mb\in \mathcal{F}_{\boldsymbol{\Lambda}_\ell}^{\mathrm{err}}$ respectively. We can bound the number of choices of $\Mb$ by Proposition \ref{prop.layerrec3}, and reduce the $L^1$ norms of $|\Ic\Nc_\Mb|$ to $\Ic_\Mb(Q_\Mb)$ by By Proposition \ref{prop.local_int}. We then apply the operation sequence to $\Mb$ as in Proposition \ref{prop.comb_est} (which includes the upper bound on the number of sub-cases produced), and apply Proposition \ref{prop.cutting_algorithm} to reduce to estimating $\Ic_{\Mb'}(Q')$ for each sub-case. In summary, we get
\begin{equation}\label{eq.proof_int_est_from_alg_step0_3}\|E_H\|_{L^1}\leq \sum_{m;\rho\geq |H|}C^{m}|\log\varepsilon|^{C^*\rho}\cdot\varepsilon^{(d-1)|H|}\sup_{\substack{\Mb\in\Fc_{\vLambda_\ell}:r(\Mb)=H\\|\Mb|=m,\,\rho\mathrm{\,fixed}}}\sup_{\textrm{sub-cases}}\Ic_{\Mb'}(Q'),\end{equation}
and the same for $\|\mathrm{Err^1}\|_{L^1}$ but with $\Fc_{\vLambda_\ell}$ replaced by $\Fc_{\vLambda_\ell}^{\mathrm{err}}$. Here note that $|H|\leq\rho$ because $|H|=|r(\Mb_\ell)|=s_\ell$; moreover, for $\Mb\in \Fc_{\vLambda_\ell}^{\mathrm{err}}$ we must have $|\Mb|\geq \Lambda_\ell$ (as there exists a component of $\Mb_{\ell'}$ that has exactly $\Lambda_\ell$ clusters, see Definition \ref{def.set_T_F} (\ref{it.set_F_err})).

Now, using (\ref{eq.proof_int_est_from_alg_step0_3}) and the above discussions, it is clear that Proposition \ref{prop.cumulant_est} would follow if we can prove the following estimate for $\Ic_{\Mb'}(Q')$ (note that $3^{-d-2}$ on the right hand side of (\ref{eq.cumulant_est_1}) equals $\upsilon/3$):
\begin{equation}\label{eq.proof_int_est_from_alg_step0_*}
    \Ic_{\Mb'}(Q')\leq \tau^{|\Mb|/9}\cdot\varepsilon^{- (d - 1) |H| + (2\upsilon/5)+(5C_{13}^*)^{-1}\cdot \rho}\cdot |\log\varepsilon|^{C^*\rho}
\end{equation} for any $\Mb\in \Fc_{\vLambda_\ell}$ and any sub-case in the operation sequence that leads to $\Mb'$ and $Q'$. Note that for $\Mb\in \Fc_{\vLambda_\ell}^{\mathrm{err}}$ the right hand side of (\ref{eq.proof_int_est_from_alg_step0_*}) should include an extra factor of $\varepsilon^{-C^*|H|}$, and we will leave it to Part 6 of the proof below. In practice this factor will not matter because $|H|\leq A_{\ell}$ and $\Lambda_{\ell}\geq A_\ell^4\geq |\log\varepsilon|^4$ by (\ref{eq.defLambdaseq}), so $\varepsilon^{-C^*|H|}$ is easily absorbed by $\tau^{|\Mb|/2}$ provided $|\Mb|\geq\Lambda_\ell$.

The rest part of this proof is devoted to the proof of \eqref{eq.proof_int_est_from_alg_step0_*}. We divide into 6 parts.

\textbf{Proof part 1: properties of final components.} In this part, as a preparation, we prove a simple fact for the (elementary) final components in $\Mb'$ after the operation sequence. More precisely, we will show that
    \begin{equation}\label{eq.not_3_comp_upp_3}
        \#_{\neg\mathrm{N}\{3\}}\le C(\#_{\mathrm{good}}+\rho),
    \end{equation} where $\#_{\neg\mathrm{N}\{3\}}$ is the number of final components in $\Mb'$ that are either nor normal or not \{3\} molecules, and $\#_{\mathrm{good}}$ and $\#_{\mathrm{bad}}$ are the number of good and bad final components in $\Mb'$ respectively.

To prove (\ref{eq.not_3_comp_upp_3}), we will consider the quantity $3|\Mb|-|\Ec_*|$ under the operation sequence. It is easy to see that this quantity decreases by $1$ under deleting, and does not change under cutting. The latter is obvious when there is no O-atom (since a bond is cut into a free and a fixed end, and fixed ends do not count in $|\Ec_*|$); in general it can be directly verified using Definition \ref{def.cutting} (\ref{it.cutting_2})).

Next consider the initial value of $3|\Mb|-|\Ec_*|$ for $\Mb$ (for which $\Ec=\Ec_*$ as there is no fixed end). In Part 1 of the proof of Proposition \ref{prop.local_int}, we have proved that $4|\Mb|=2|\Ec|-|\Ec_{\mathrm{end}}^-|-|\Ec_{\mathrm{end}}^+|$. Note that each edge in $\Ec$ is classified as either bond, or top end in $\Ec_{\mathrm{end}}^+$, or bottom end in $\Ec_{\mathrm{end}}^-$, and each empty end counts as both top end end bottom end, we deduce that \begin{equation}\label{eq.rec_new}3|\Mb|-|\Ec|=|\Ec|-|\Ec_{\mathrm{end}}^-|-|\Ec_{\mathrm{end}}^+|-|\Mb|=\#_{\mathrm{bond}}-|\Mb|-\#_{\mathrm{ee}},\end{equation} where $\#_{\mathrm{bond}}$ and $\#_{\mathrm{ee}}$ are the number of bonds and empty ends in $\Mb$ respectively. Note that $\#_{\mathrm{ee}}\geq 0$, and the number of bonds between different layer $\Mb_{\ell'}$ is bounded by $\sum_{\ell'}|r(\Mb_{\ell'})|$ (as each such bond corresponds to a particle line crossing two different layer which will be in some $r(\Mb_{\ell'})$), we see that the initial value
\begin{equation}\label{eq.extra_rec}
3|\Mb|-|\Ec|\leq \sum_{\ell'}|r(\Mb_{\ell'})|+\sum_{\ell'}(\#_{\mathrm{bond}(\Mb_{\ell'})}-|\Mb_{\ell'}|)\leq \sum_{\ell'}|r(\Mb_{\ell'})|+\sum_{\ell'}\rho(\Mb_{\ell'})\leq \rho
\end{equation} using also Definition \ref{def.recollision_number}.

Now we prove \eqref{eq.not_3_comp_upp_3}. By the above discussions, we see that $3|\Mb'|-|\Ec_*'|\leq\rho$ for the final result $\Mb'$ (and corresponding $\Ec_*'$). Note that this value equals 0 for \{3\} molecules, equals $-1$ for \{2\}  and \{33\}  molecules, and equals $1$ for \{4\}  and \{44\}  molecules and empty ends, so we conclude that
\begin{equation}\label{eq.proof_not_3_comp_upp_2}
    \#_{\{2\}}+\#_{\{33\}}\leq \rho+\#_{\mathrm{ee}}+\#_{\{4\}}+\#_{\{44\}}\leq 4(\#_{\mathrm{good}}+\#_{\mathrm{bad}}+\rho)
\end{equation} but now for $\Mb'$, where $\#_{\{2\}}$ and $\#_{\{33\}}$ etc. are the number of corresponding molecules in $\Mb'$. Here the last inequality is because each \{4\}  and \{44\} molecule and empty end is either good or bad by Definition \ref{def.good_normal}. Finally, since there are only elementary components in $\Mb'$, by Definition \ref{def.elementary}, we have \[\#_{\neg \mathrm{N}\{3\}} = \#_{\{2\}}+\#_{\{33\}}+\#_{\mathrm{ee}}+\#_{\{4\}}+\#_{\{44\}}+\#(\textrm{good \{3\} molecules}),\] thus $\#_{\neg \mathrm{N}\{3\}}\leq 9(\#_{\mathrm{good}}+\#_{\mathrm{bad}}+\rho)$ thanks to (\ref{eq.proof_not_3_comp_upp_2}). Since we also have $\#_{\mathrm{bad}}\leq C(\#_{\mathrm{good}}+|H|)$ by (\ref{eq.overall_alg}) and $|H|\leq\rho$, we conclude that $\#_{\neg \mathrm{N}\{3\}}\leq C(\#_{\mathrm{good}}+\rho)$. This proves \eqref{eq.not_3_comp_upp_3}.

\textbf{Proof part 2: decomposition of the support of $Q$.} Recall $Q'$ corresponds to a sub-case in the operation sequence obtained by splittings; as in Proposition \ref{prop.cutting_algorithm}, we have $Q'=Q_\Mb\cdot\mathbbm{1}_{\Pc}$ where $\Pc$ represents the restriction in this sub-case. Recall also the $\mathbbm{1}_\Dc$ and $\mathbbm{1}_\Lc^{\varepsilon}$ in (\ref{eq.Q_M}). In this part, we will decompose $Q_\Mb$ as \begin{equation}\label{eq.proof_int_est_from_alg_step1_*}
    Q_\Mb\le \sum_\alpha w_\alpha Q_\alpha,
\end{equation} where $\alpha$ is the summation index, $w_\alpha\geq 0$ is a scalar, and $Q_\alpha$ are nonnegative functions that have localized support in the $x_e$ and $v_e$ variables, specified as below:
\begin{enumerate}
    \item\label{it.proof_int_est_from_alg_step1_1} Each $Q_\alpha$ has the following form for some $(x_e^*)$ with $x_e^*\in\Zb^d$: 
    \begin{equation}\label{eq.proof_int_est_from_alg_step1_*2}
    \begin{aligned}
        Q_\alpha= \mathbbm{1}_\Dc\cdot \mathbbm{1}_\Lc^{\varepsilon}\cdot\prod_{e\in \Ec_{\mathrm{end}}^-} \big(e^{-\beta|v_e|^2/3}\cdot \oneb_{|x_e-x_e^*|\leq |\log\varepsilon|^{C^*},\,|v_e|\leq |\log\varepsilon|^{C^*}}\big).
    \end{aligned}
    \end{equation}
\item\label{it.proof_int_est_from_alg_step1_2} The $w_\alpha$ satisfies
    \begin{equation}\label{eq.proof_int_est_from_alg_step1_*1}
        \sum_\alpha w_\alpha\le C^{|\Mb|}.
    \end{equation}
    \item\label{it.proof_int_est_from_alg_step1_3} In the support of $Q_\alpha$ and $\prod_{\nf\in\Mb}\boldsymbol{\Delta}_\nf$ in $\Ic_{\Mb}$ (see \eqref{eq.associated_int_op}), for each edge $e\in \Ec$, we have $|x_{e}-x_{e'}^*|\leq|\log\varepsilon|^{C^*}$ and $|v_e|\leq |\log\varepsilon|^{C^*}$, where $e'\in \Ec_{\mathrm{end}}^-$ is the bottom end of the particle line containing $e$. 
\end{enumerate} 

Now we proceed to prove (\ref{eq.proof_int_est_from_alg_step1_*}) and (\ref{it.proof_int_est_from_alg_step1_1})--(\ref{it.proof_int_est_from_alg_step1_3}). Recall that $Q_\Mb$ is given by \eqref{eq.Q_M},
\begin{equation}\label{eq.proof_int_est_from_alg_step1_1}
    Q_\Mb = \mathbbm{1}_\Dc(\vt_\Mb)\cdot \mathbbm{1}_\Lc^{\varepsilon}(\vz_\Ec)\cdot \prod_{e\in \Ec_{\mathrm{end}}^-}\big|f^\Ac\big((\ell_1[e]-1)\tau,x_e-(\ell_1[e]-1)\tau\cdot v_e,v_e\big)\big|.
\end{equation} We begin by decomposing $f^\Ac$ and then inserting into (\ref{eq.proof_int_est_from_alg_step1_1}).

\textbf{Decomposition of $f^\Ac$.} By \eqref{eq.approxfA_1}, \eqref{eq.approxfA_2},  \eqref{eq.boltzmann_decay_4} and \eqref{eq.decayseq}, we know that $\|f^\Ac(\ell\tau)\|_{\mathrm{Bol}^{2\beta/5}}\le C$. The definition of $\|\cdot\|_{\mathrm{Bol}^{2\beta/5}}$ \eqref{eq.boltzmann_decay_2} implies that $f^\Ac$ decay exponentially in $v$, so we may restrict $|v_e|\leq |\log\varepsilon|^{C^*}$ for each bottom end $e$, otherwise the exponential decay will overwhelm the powers $\varepsilon^{-(d-1)(|\Ec_*|-2|\Mb|)}$ in \eqref{eq.associated_int_op} and the desired bound becomes trivial. 

Next, with the assumption that $f^\Ac$ is supported in $|v_e|\leq |\log\varepsilon|^{C^*}$ , we will decompose
\begin{equation}\label{eq.proof_int_est_from_alg_step1_dec}
\begin{gathered}
    |f^\Ac(x, v)| \le \sum_{x^*\in \Zb^d} w_{x^*}\cdot g_{x^*}(x, v),
    \\
    g_{x^*}(x, v)\le e^{-(2\beta/5)|v|^2}\cdot \oneb_{|x-x^*|\leq C,\,|v|\leq |\log\varepsilon|^{C^*}},\quad \sum_{x^*\in \Zb^d} w_{x^*}\le C,
\end{gathered}
\end{equation} In fact, this can be proved by choosing $C = \sqrt{d}$ and 
\[w_{x^*}:=\sup_{|x-x^*|\leq C,\, v\in \Rb^d} e^{\beta |v|^2} |f^\Ac(x, v)|,\quad g_{x^*}(x, v)=(w_{x^*})^{-1}\cdot g(x, v)\cdot \oneb_{|x-x^*|\leq C},\] and (\ref{eq.proof_int_est_from_alg_step1_dec}) follows from the definition of $\|g\|_{\mathrm{Bol}^{2\beta/5}}$.

\textbf{Inserting the decomposition of $f^\Ac$.} Now, by inserting the decomposition (\ref{eq.proof_int_est_from_alg_step1_dec}) into \eqref{eq.proof_int_est_from_alg_step1_1}, we get
\begin{equation}\label{eq.proof_int_est_from_alg_step1_3}
    Q_\Mb \le \sum_{\alpha} \underbrace{\bigg(\prod_{e\in \Ec_{\mathrm{end}}^-}w_{x_e^*}\bigg)}_{w_\alpha} \underbrace{\bigg(\mathbbm{1}_\Dc\cdot \mathbbm{1}_\Lc^{\varepsilon}\cdot \prod_{e\in \Ec_{\mathrm{end}}^-} g_{x_e^*}((\ell_1[e]-1)\tau,x_e-(\ell_1[e]-1)\tau\cdot v_e,v_e)\bigg)}_{Q_\alpha}=\sum_\alpha w_\alpha Q_\alpha,
\end{equation}
where $\alpha=(x_e^*)_{e\in \Ec_{\mathrm{end}}^-}$, and $g_{x_e^*}$ is as in (\ref{eq.proof_int_est_from_alg_step1_dec}) but with $x_e^*$ instead of $x^*$. This already establishes (\ref{eq.proof_int_est_from_alg_step1_*}), and (\ref{it.proof_int_est_from_alg_step1_1})--(\ref{it.proof_int_est_from_alg_step1_2}) also follow from (\ref{eq.proof_int_est_from_alg_step1_dec}) where we note that $|(\ell_1[e]-1)\tau\cdot v_e|\leq|\log\varepsilon|^{C^*}$. 

Finally we prove \eqref{it.proof_int_est_from_alg_step1_3}. The bound for $v_e$ follows from (\ref{eq.proof_int_est_from_alg_step1_dec}) and energy conservation implied by the support of $\Dirac_\nf$: $|v_e|^2\leq\sum_{e'\in\Ec_{\mathrm{end}}^-}|v_{e'}|^2$ (and that the size $|\Mb|_p\leq |\log\varepsilon|^{C^*}$). As for $x_e$, let $e$ be any edge and $e'$ be the bottom end of the particle line containing $e$. Then we have $|x_{e'}-x_{e'}^*|\leq |\log\varepsilon|^{C^*}$ by (\ref{eq.proof_int_est_from_alg_step1_dec}); moreover, since $e$ is equal to or connected to $e'$ by a sequence of C- and O-atoms, we can recursively apply the relation between incoming and outgoing positions in \eqref{eq.delta_supp} to show that $|x_{e}-x_{e'}|\le \sum_j \Lf\tau\cdot|v_{e_j}|\leq|\log\varepsilon|^{C^*}$ where $e_j$ are edges on the path from $e$ to $e'$. This proves \eqref{it.proof_int_est_from_alg_step1_3}.

Now that we have the decomposition (\ref{eq.proof_int_est_from_alg_step1_dec}), using that $Q' = \mathbbm{1}_\Pc\cdot Q_\Mb$, we can get the associated decomposition of $Q'$ by defining $Q_\alpha':=\mathbbm{1}_\Pc\cdot    Q_\alpha$, so that $Q'\leq \sum_\alpha w_\alpha Q_\alpha'$.

\textbf{Proof part 3: estimating the local integral.} In this part, we reduce the main goal \eqref{eq.proof_int_est_from_alg_step0_*} to estimates of the local integrals for the elementary molecules $(\Mb_1,\cdots,\Mb_k)$ forming $\Mb'$. By \eqref{eq.proof_int_est_from_alg_step1_*} and \eqref{eq.proof_int_est_from_alg_step1_*1}, we know that \eqref{eq.proof_int_est_from_alg_step0_*} would follow once we can prove (uniformly in $\alpha$)
\begin{equation}\label{eq.proof_int_est_from_alg_step2_1}
    \Ic_{\Mb'}(Q'_\alpha)\leq \tau^{|\Mb|/8}\cdot\varepsilon^{- (d - 1) |H| + (2\upsilon/5)+(5C_{13}^*)^{-1}\cdot \rho}\cdot |\log\varepsilon|^{C^*\rho}.\end{equation}

Below we will prove (\ref{eq.proof_int_est_from_alg_step2_1}). Before continuing, let us discuss some useful heuristics. By \eqref{eq.associated_op_disjoint} we have $\Ic_{\Mb'} = \Ic_{\Mb_1}\circ \Ic_{\Mb_2}\circ\cdots\circ \Ic_{\Mb_k}$, and by Propositions \ref{prop.intmini}--\ref{prop.intmini_2} we have upper bounds for each $\Ic_{\Mb_j}$. Basically $\Ic_{\Mb'_j} = O(\tau)$ for normal molecules, $\Ic_{\Mb'_j} =O(\varepsilon^{\upsilon})$ for good molecules, and $\Ic_{\Mb'_j} =O(\varepsilon^{-(d-1)})$ for bad molecules. We would like to insert these bounds into \eqref{eq.proof_int_est_from_alg_step2_1}, and derive an estimate of the form $\Ic_{\Mb'}(Q'_\alpha)\le \tau^{\mathrm{power}} \varepsilon^{\mathrm{power}}$ with exponents matching those in \eqref{eq.proof_int_est_from_alg_step2_1}.

However, there is one subtlety regarding normal \{3\} molecules (the other molecules are fine as their number is bounded by (\ref{eq.not_3_comp_upp_3}), so logarithmic loss coming from them will be negligible in view of the gain from good molecules and the gain on the right hand side of (\ref{eq.proof_int_est_from_alg_step2_1})). For a normal \{3\} molecule $\Mb_j$, by Proposition \ref{prop.intmini} (\ref{it.intmini_2}) we know $\Ic_\Mb(Q) = \int_{\Rb\times\Sb^{d-1}\times\Rb^d}\big[(v_1-v_2)\cdot\omega\big]_-\cdot Q\,\mathrm{d}t_1\mathrm{d}\omega\mathrm{d}v_2$, and $t_1\in [(\ell[\nf]-1)\tau,\ell[\nf]\tau]$ in this integration. If $Q$ is supported in $|v_e|\lesssim 1$, this leads to $\|\Ic_{\Mb}(Q)\|_{L^\infty}\lesssim\tau\|Q\|_{L^\infty} = O(\tau)$ from the integral over $t$. Nevertheless, in the support of $Q'_\alpha$, we only have the weaker bound $|v_e|\leq|\log\varepsilon|^{C^*}$, which leads to logarithmic losses due to both the volume of $v_2$ and the cross section $\big[(v_1-v_2)\cdot\omega\big]_-$. If estimated naively, we would have losses of $O \left(|\log\varepsilon|^{C^*|\Mb|}\right)$, which cannot be controlled by \eqref{eq.proof_int_est_from_alg_step2_1}.

We will address these issues in Part 4 below, using Propositions \ref{prop.weight}--\ref{prop.volume} whose proof will be postponed to Section \ref{sec.aux}. For now, let us introduce some preliminary reductions on \eqref{eq.proof_int_est_from_alg_step2_1}.

Recall $\Ic_{\Mb'}(Q'_\alpha)=\Ic_{\Mb_1}\circ\cdots\circ \Ic_{\Mb_k}(Q'_\alpha)$. For a final component $\Mb_j$ that is a normal \{3\} molecule with one atom $\nf$, we have $\Ic_{\Mb_j} (Q) = \int_{\Rb\times\Sb^{d-1}\times\Rb^d}\big[(v_{e_1(\nf)}-v_{e_2(\nf)})\cdot\omega_\nf\big]_-\cdot Q\,\mathrm{d}t_\nf\mathrm{d}\omega_\nf\mathrm{d}v_{e_2(\nf)}$ (by \eqref{eq.intmini_3}), where $e_1(\nf)$, $e_2(\nf)$ are the two bottom (or top) edges at $\nf$. We introduce the following notions:
\begin{equation}\label{eq.proof_int_est_from_alg_step2_2}
\begin{aligned}
    \int(\cdots) \,\mathrm{d}\Omega_{\Mb_j} &:= \int(\cdots) \,\mathrm{d}t_\nf\mathrm{d}\omega_\nf\mathrm{d}v_{e_2(\nf)},& \int(\cdots) \,\mathrm{d}\Omega_{\mathrm{N}\{3\}} &= \int(\cdots) \prod_{\Mb_j \in\Sf}\,\mathrm{d}\Omega_{\Mb_j},
    \\
    B_{\Mb_j} &:= \big[(v_{e_1(\nf)}-v_{e_2(\nf)})\cdot\omega_\nf\big]_-,& B_{\mathrm{N}\{3\}} &:= \prod_{\Mb_j\in\Sf} B_{\Mb_j},
    \\ 
    \Ic^{\neg\mathrm{N}\{3\}}_{\Mb'} &= \prod_{\Mb_j\not\in\Sf} \Ic_{\Mb_j},&\Sf&:=\{\textrm{normal \{3\} molecules}\}.
\end{aligned}
\end{equation}
Then $\Ic_{\Mb'}(Q'_\alpha)$ can be rewritten as 
\begin{equation}\label{eq.proof_int_est_from_alg_step2_3}
    \Ic_{\Mb'}(Q'_\alpha) = \int \Ic^{\neg\mathrm{N}\{3\}}_{\Mb'} \Big(B_{\mathrm{N}\{3\}}\cdot Q'_\alpha\Big)\ \mathrm{d}\Omega_{\mathrm{N}\{3\}}.
\end{equation}
Here, in order to obtain (\ref{eq.proof_int_est_from_alg_step2_3}), we have chosen to integrate the $\mathrm{d}\Omega_{\mathrm{N}\{3\}}$ variables (which are independent variables) in the end, and treat these variables as fixed when applying the remaining operators $\Ic_{\Mb_j}$ and integrating in the remaining variables. Schematically, we start with a composition operator $\int_{(1)} f_1\int_{(2)} f_2\cdots(Q)$ and rewrite it as $\int_{(1)} f_1\int_{(2)} f_2\cdots(Q) = \int_{(1)}\int_{(2)}\cdots(f_1 f_2 \cdots Q)$, then change the order of integration to move a part of the variables to the (left) end of the integral.

\textbf{Proof part 4: weight and volume bounds.} In this part, we derive the upper bounds of the weights $B_{\mathrm{N}\{3\}}$ and the integral over $\mathrm{d}\Omega_{\mathrm{N}\{3\}}$. We will rely on the estimates in Propositions \ref{prop.weight} and \ref{prop.volume}; for their statements and proof see Section \ref{sec.aux}.

\textbf{Application of the upper bound of weights.} We first derive an upper bound of the cross sections $B_{\mathrm{N}\{3\}}$. The goal is to prove (\ref{eq.proof_int_est_from_alg_step31_5}). Thanks to Proposition \ref{prop.weight}, we have
\begin{equation}\label{eq.proof_int_est_from_alg_step31_1}
 B_{\mathrm{N}\{3\}} = \prod_{\Mb_j\in\Sf} B_{\Mb_j} \leq  \prod_{\Mb_j\in\Sf} |v_{e_1(\nf)}-v_{e_2(\nf)}|\le\prod_{\nf}(1+|v_{e_1(\nf)}-v_{e_2(\nf)}|)\leq \Bf\cdot e^{(\beta/6)\sum_{e\in\Ec_{\mathrm{end}}^-}|v_e|^2},
\end{equation} 
where $e_1(\nf)$, $e_2(\nf)$ are the two bottom (or top) edges at $\nf$, and in the last step we have used Proposition \ref{prop.weight} with $\gamma = \beta/6$. The quantity $\Bf$ is as in Proposition \ref{prop.weight} and will be discussed only in Part 5 below.

Now, combining (\ref{eq.proof_int_est_from_alg_step31_1}) with the definition of $Q_\alpha$ in \eqref{eq.proof_int_est_from_alg_step1_*2} (with $Q_\alpha'=\mathbbm{1}_\Pc\cdot    Q_\alpha$), we obtain
\begin{equation}\label{eq.proof_int_est_from_alg_step31_3}
\begin{split}
        B_{\mathrm{N}\{3\}}\cdot Q'_\alpha&\leq \Bf\cdot\mathbbm{1}_\Pc\cdot\mathbbm{1}_\Dc\cdot \mathbbm{1}_\Lc^{\varepsilon}\cdot\prod_e\oneb_{|x_e-x_e^*|\leq C,\,|v_e|\leq |\log\varepsilon|^{C^*}}\cdot e^{-(\beta/6)\sum_{e\in \Ec_{\mathrm{end}}^-}|v_e|^2}
        \\
        &\leq \Bf\cdot\mathbbm{1}_\Pc\cdot\mathbbm{1}_\Dc\cdot \mathbbm{1}_\Lc^{\varepsilon}\cdot\prod_e\oneb_{|x_e-x_e^*|\leq C,\,|v_e|\leq |\log\varepsilon|^{C^*}}
                \\
        &
        \times\sum_{(X_e)_{e\in \Ec_{\mathrm{end}}}:\, X_e\in [1,|\log\varepsilon|^{C^*}]\cap 2^\Zb} \prod_{e\in \Ec_{\mathrm{end}}}e^{-(\beta/48)|X_e|^2}\cdot \oneb_{|v_e|\le X_e}.
\end{split}
\end{equation} Here in (\ref{eq.proof_int_est_from_alg_step31_3}), we have used energy conservation to replace $(\beta/6)\sum_{e\in \Ec_{\mathrm{end}}^-}|v_e|^2$ by $(\beta/12)\sum_{e\in \Ec_{\mathrm{end}}}|v_e|^2$ (with $\Ec_{\mathrm{end}}$ being the set of all ends of $\Mb$), and performed a dyadic decomposition in $|v_e|$ for each $e\in \Ec_{\mathrm{end}}$ (i.e. we have $e^{-(\beta/12)|v|^2}\le \sum_{X} e^{-(\beta/48)|X|^2}\oneb_{|v|\le X}$; note that we may assume $\beta\leq 1$).

Inserting \eqref{eq.proof_int_est_from_alg_step31_3} into \eqref{eq.proof_int_est_from_alg_step2_3}, we get
\begin{equation}\label{eq.proof_int_est_from_alg_step31_4}
    \begin{split}
        &\Ic_{\Mb'}(Q'_\alpha) \le \sum_{(X_e)_{e\in \Ec_{\mathrm{end}}}} e^{-(\beta/96)\sum_{e\in \Ec_{\mathrm{end}}}|X_e|^2} \int \Ic^{\neg\mathrm{N}\{3\}}_{\Mb'} \left(\Bf\cdot \mathbbm{1}_{(\heartsuit)}\cdot e^{-(\beta/96)\sum_{e\in \Ec_{\mathrm{end}}}|X_e|^2}\right)\ \mathrm{d}\Omega_{\mathrm{N}\{3\}},
    \end{split}
\end{equation}
where 
\begin{equation}\label{eq.proof_int_est_from_alg_step31_4'}
    \mathbbm{1}_{(\heartsuit)} \coloneqq \mathbbm{1}_\Pc\cdot\mathbbm{1}_\Dc\cdot \mathbbm{1}_\Lc^{\varepsilon}\cdot\oneb_{|x_e-x_e^*|\leq C,\,|v_e|\leq |\log\varepsilon|^{C^*},\, \forall e\in \Ec_{\mathrm{end}}^-}\cdot\oneb_{|v_e|\le X_e,\, \forall e\in \Ec_{\mathrm{end}}}.
\end{equation}
Clearly the sum $\sum_{(X_e)} e^{-(\beta/96)\sum_e|X_e|^2}$ only leads to a factor $C^{|\Mb|}$. Therefore we derive that
\begin{equation}\label{eq.proof_int_est_from_alg_step31_5}
    \Ic_{\Mb'}(Q'_\alpha) \le C^{|\Mb|} \sup_{(X_e)_{e\in \Ec_{\mathrm{end}}}}\, \Bf\cdot e^{-(\beta/96)\sum_{e\in \Ec_{\mathrm{end}}}|X_e|^2}\cdot\int \Ic^{\neg\, \{3\}}_{\Mb'} \left(\mathbbm{1}_{(\heartsuit)}\right)\ \mathrm{d}\Omega_{\mathrm{N}\{3\}},
\end{equation}
where $X_e\in [1,|\log\varepsilon|^{C^*}]\cap 2^\Zb$ in the supremum.

\textbf{Treating the integral in $\mathrm{d}\Omega_{\mathrm{N}\{3\}}$.} In this section, we apply Proposition \ref{prop.volume} to bound the volume of the set which the variables in $\mathrm{d}\Omega_{\mathrm{N}\{3\}}$ belong to, thereby controlling the integral in these variables by the corresponding supremum. The goal is to prove (\ref{eq.proof_int_est_from_alg_step33_1}).

Recall the set $\Sf$ of normal \{3\} molecules in $\Mb'$, which also occurs in Proposition \ref{prop.volume}. Let the associated notations be as in Proposition \ref{prop.volume}: $\nf$ is the unique atom in $\Mb_j$, $e_1(\nf)$ is the unique fixed end at $\nf$, and $e_2(\nf)$ is the other top (bottom) end at $\nf$ if $e_1(\nf)$ is a top (bottom) end. Recall that $\mathrm{d}\Omega_{\mathrm{N}\{3\}} := \prod_{\Mb_j\in\Sf}\mathrm{d}\Omega_{\Mb_j}$ and $\mathrm{d}\Omega_{\Mb_j} := \mathrm{d}t_\nf\mathrm{d}\omega_\nf\mathrm{d}v_{e_2(\nf)}$, see \eqref{eq.proof_int_est_from_alg_step2_2}. We also denote by $e_1'(\nf)$ and $e_2'(\nf)$ the two remaining (free) ends at $\nf$. Then, using the support of $\boldsymbol{\Delta}_\nf$ in $\Ic_{\Mb}$, we always have 
\begin{equation}\label{eq.proof_int_est_from_alg_step32_1}
    v_{e_1'(\nf)}=v_{e_1(\nf)}-[(v_{e_1(\nf)}-v_{e_2(\nf)})\cdot\omega_\nf]\omega_\nf,\quad v_{e_2'(\nf)}=v_{e_2(\nf)}+[(v_{e_1(\nf)}-v_{e_2(\nf)})\cdot\omega_\nf]\omega_\nf
\end{equation}  (this is for C-atoms; for O-atoms we simply have $v_{e_j'(\nf)}=v_{e_j(\nf)}$).

The idea in this part is to fix $(\omega_\nf)$ (which belongs to a compact set) and control the volume of the sets of $(t_\nf)$ and $(v_{e_2(\nf)})$ for the normal \{3\} molecule variables. For $(t_\nf)$ this is simple: let $A$ be the set of atoms $\nf$ \emph{not} in one of the normal \{3\} molecules in $\Mb'$, then $\Mb\backslash A$ is associated with the normal \{3\} molecules. Since $Q$ is supported in $\Dc$ in (\ref{eq.associated_domain}), we know that $(t_\nf)_{\nf\in \Mb\backslash A}\in \widetilde{\Dc}$, where $\widetilde{\Dc}$ is the projection of the set $\Dc$ onto the variables $\vt_{\Mb\backslash A}$. Note that the volume $|\widetilde{\Dc}|=\int_{\widetilde{\Dc}}1\,\mathrm{d}\vt_{\Mb\backslash A}$ is precisely the quantity in (\ref{eq.weight_est_2}) that controls the quantity $\Bf$ in Proposition \ref{prop.weight}.

As for $(v_{e_2(\nf)})$, we will use Proposition \ref{prop.volume}. Recall that $X_e\in [1,|\log\varepsilon|^{C^*}]\cap2^\Zb$ has been fixed for all $e\in \Ec_{\mathrm{end}}$ as in (\ref{eq.proof_int_est_from_alg_step31_5}). We may also fix $(\omega_\nf)$. Then, by Proposition \ref{prop.volume}, there exists a set $Y=Y(\omega)\subseteq (\Rb^d)^{|\Fc|}$ (where $\Fc$ is the set of all $e_2(\nf)$), which depends on $(\omega_\nf)$ and $(X_e)$, and on $(\Mb,\Mb')$ and the cutting sequence, but not on the concrete choices of $(v_{e_2(\nf)})$, such that $(v_{e_2(\nf)})$ always takes value in $Y$, and \begin{equation}\label{eq.proof_int_est_from_alg_step32_2}
   |Y|\leq \prod_{e\in \Ec_{\mathrm{end}}}X_e^d\cdot C^{|\Mb|}\cdot |\log\varepsilon|^{C^*\#_{\neg\mathrm{N}\{3\}}}.
    \end{equation}

Putting together the above arguments, we get that 
\begin{equation}\label{eq.proof_int_est_from_alg_step32_3}
\begin{split}
    \int \Ic^{\neg\mathrm{N}\{3\}}_{\Mb'} \left(\mathbbm{1}_{(\heartsuit)}\right)\,\mathrm{d}\Omega_{\mathrm{N}\{3\}} &=\int\prod_{\nf\in \Mb\backslash A}\mathrm{d}\omega_\nf\cdot\int_{\widetilde{\Dc}}\prod_{\nf\in \Mb\backslash A}\mathrm{d}t_\nf\cdot\int_Y\prod_{\nf\in \Mb\backslash A}\mathrm{d}v_{e_2(\nf)}\cdot \Ic^{\neg\mathrm{N}\{3\}}_{\Mb'} \left(\mathbbm{1}_{(\heartsuit)}\right)
    \\
    &\leq |\Sb^{d-1}|^{|\Mb\backslash A|}\cdot |\widetilde{\Dc}|\cdot|Y|\cdot\sup_{(\omega_\nf):|\omega_\nf| = 1} \sup_{(t_\nf)\in \widetilde{\Dc}} \sup_{(v_{e_2(\nf)})\in Y(\omega)}
     \Ic^{\neg\mathrm{N}\{3\}}_{\Mb'} \left(\mathbbm{1}_{(\heartsuit)}\right).
\end{split}
\end{equation}
Combining this with (\ref{eq.proof_int_est_from_alg_step31_5}) and (\ref{eq.proof_int_est_from_alg_step32_2}), we get
\begin{equation}\label{eq.proof_int_est_from_alg_step33_1}
    \Ic_{\Mb'}(Q'_\alpha)  \le C^{|\Mb|}\cdot e^{-(\beta/96)\sum_{e\in \Ec_{\mathrm{end}}}|X_e|^2}\cdot \Bf\cdot |\widetilde{\Dc}|\cdot \prod_{e\in \Ec_{\mathrm{end}}}X_e^d\cdot |\log\varepsilon|^{C^*\#_{\neg\mathrm{N}\{3\}}}\cdot \sup\left[\Ic^{\neg\mathrm{N}\{3\}}_{\Mb'} \left(\mathbbm{1}_{(\heartsuit)}\right)\right],
\end{equation}
where $(X_e)$ is fixed, $|\Sb^{d-1}|^{|\Mb\backslash A|}$ is merged into $C^{|\Mb|}$, the supremum is taken over $(\omega_\nf)$, $(t_\nf)$ and $(v_{e_2(\nf)})$, and $\mathbbm{1}_{(\heartsuit)}$ is given by 
\begin{equation}\label{eq.proof_int_est_from_alg_step33_2}
    \mathbbm{1}_{(\heartsuit)} \coloneqq \mathbbm{1}_\Pc\cdot\mathbbm{1}_\Dc\cdot \mathbbm{1}_\Lc^{\varepsilon}\cdot\oneb_{|x_e-x_e^*|\leq C,\,|v_e|\leq |\log\varepsilon|^{C^*},\,\forall e\in \Ec_{\mathrm{end}}^-}\cdot\oneb_{|v_e|\le X_e,\, \forall e\in \Ec_{\mathrm{end}}}.
\end{equation}

\textbf{Proof part 5: Finishing the proof for $E_H$.} In this step, we finish the proof of (\ref{eq.proof_int_est_from_alg_step2_1}) by applying the upper bounds of local integrals $\Ic_{\Mb_j}$. First we have
\begin{equation}\label{eq.proof_molecule_est_step4_1}
    \Ic^{\neg\mathrm{N}\{3\}}_{\Mb'} \left(\mathbbm{1}_{(\heartsuit)}\right) = \prod_{\Mb_j\not\in\Sf} \Ic_{\Mb_j}\left(\mathbbm{1}_{(\heartsuit)}\right)\le \prod_{\Mb_j\not\in\Sf} \|\Ic_{\Mb_j}\|_{L^\infty\to L^\infty}.
\end{equation} Then, for each $\Mb_j\not\in\Sf$, depending on whether it is good, normal or bad, we can apply Propositions \ref{prop.intmini}--\ref{prop.intmini_2} to bound $\|\Ic_{\Mb_j}\|_{L^\infty\to L^\infty}$. Here we note that (i) if a good molecule $\Xc$ in Definition \ref{def.good_normal} involves another molecule $\Yc\prec_{\mathrm{cut}}\Xc$, then the corresponding $z_{e'}$ (or $t_{\nf'}$) variables in $\Yc$ has been fixed when integrating in the $\Xc$ variables, so they play the role of the external variables in Proposition \ref{prop.int_mini_good}; (ii) in the initial link case corresponding to Proposition \ref{prop.int_mini_good} \eqref{it.bound_one_atom_good_empty_line}--\eqref{it.bound_one_atom_good_link}, we should extract the indicator function $\mathbbm{1}_{e_1\sim_{\mathrm{in}} e_2}$ from $\mathbbm{1}_\Lc^{\varepsilon}$ and bound $\Ic_{\Mb_j}\circ\mathbbm{1}_{e_1\sim_{\mathrm{in}} e_2}^{\varepsilon}$ instead. In the end we get

\begin{equation}\label{eq.Linftybd}
\|\Ic_{\Mb_j}\|_{L^\infty\to L^\infty}\leq|\log\varepsilon|^{C^*} \cdot\left\{
\begin{aligned}
&\varepsilon^{\upsilon}, &&\Mb_j\mathrm{\ is\ good};\\
&1, &&\mathrm{otherwise}
\end{aligned}
\right.
\times
\left\{
\begin{aligned}
&\varepsilon^{-(d-1)}, &&\Mb_j\mathrm{\ is\ bad};\\
&1, &&\mathrm{otherwise,}
\end{aligned}
\right.
\end{equation}

where recall $\upsilon=3^{-d-1}$. This leads to
\begin{equation}\label{eq.proof_molecule_est_step4_2}
\begin{split}
    \Ic^{\neg \mathrm{N}\{3\}}_{\Mb'} \left(\mathbbm{1}_{(\heartsuit)}\right) &\le \prod_{\Mb_j\not\in\Sf} \|\Ic_{\Mb_j}\|_{L^\infty\to L^\infty}\le |\log\varepsilon|^{C^*\#_{\neg\mathrm{N}\{3\}}}\cdot\varepsilon^{\upsilon\cdot\#_{\mathrm{good}}-(d-1)\#_{\mathrm{bad}}}.
\end{split}
\end{equation}
Combining with \eqref{eq.proof_int_est_from_alg_step33_1}, we get
\begin{equation}\label{eq.proof_molecule_est_step4_3}
\begin{split}
    \Ic_{\Mb'}(Q'_\alpha) &\le C^{|\Mb|}\cdot e^{-(\beta/96)\sum_{e\in \Ec_{\mathrm{end}}}|X_e|^2}\cdot \Bf\cdot |\widetilde{\Dc}|\cdot \prod_{e\in \Ec_{\mathrm{end}}}X_e^d\cdot |\log\varepsilon|^{C^*\#_{\neg\mathrm{N}\{3\}}}\cdot \sup\left[\Ic^{\neg\mathrm{N}\{3\}}_{\Mb'} \left(\mathbbm{1}_{(\heartsuit)}\right)\right]
    \\
    &\le (\Bf\cdot |\widetilde{\Dc}|)\cdot\bigg(\prod_{e\in \Ec_{\mathrm{end}}} X_e^d\cdot e^{-(\beta/96)|X_e|^2}\bigg)\cdot |\log\varepsilon|^{C^*\#_{\neg\mathrm{N}\{3\}}}\cdot\varepsilon^{\upsilon\cdot\#_{\mathrm{good}}-(d-1)\#_{\mathrm{bad}}}.
\end{split}
\end{equation}
By Proposition \ref{prop.weight} we know that $|\widetilde{\Dc}|\cdot \Bf\le (C\tau)^{|\Mb|}\cdot |\log\varepsilon|^{C^*(\rho + |A|)}$, and it is also obvious that $\prod_e (X_e^d\cdot e^{-(\beta/96)|X_e|^2})\le C^{|\Mb|}$. Recall $\#_{\neg\mathrm{N}\{3\}}$ is the number of molecules in $\Mb'$ that are not normal \{3\} molecules, and we have $|A|\le 2\cdot\#_{\neg\mathrm{N}\{3\}}$ (as an elementary component has at most two atoms). Using also (\ref{eq.not_3_comp_upp_3}), we know that $|A|+\#_{\neg\mathrm{N}\{3\}}\leq C(\#_{\mathrm{good}}+\rho)$. Inserting all the above into \eqref{eq.proof_molecule_est_step4_3}, we obtain
\begin{equation}\label{eq.proof_molecule_est_step4_4}
    \Ic_{\Mb'}(Q'_\alpha) \le |\log\varepsilon|^{C^*(\rho+\#_{\mathrm{good}})}\cdot \varepsilon^{\upsilon\cdot\#_{\mathrm{good}}-(d-1)\#_{\mathrm{bad}}}\cdot (C\tau)^{|\Mb|}.
\end{equation}

Now consider (\ref{eq.overall_alg}) in Proposition \ref{prop.comb_est}. Since $\rho>0$, we know the left hand side of (\ref{eq.overall_alg}) is at least $\upsilon/2$; by interpolation we get that this is also $\geq  (2\upsilon/5)+(5C_{13}^*)^{-1}\cdot \rho$. Using this and (\ref{eq.proof_molecule_est_step4_4}), we conclude that
\[\Ic_{\Mb'}(Q'_\alpha) \le (C\tau)^{|\Mb|}\cdot \big(|\log\varepsilon|^{C^*}\varepsilon^{\upsilon/2}\big)^{\#_{\mathrm{good}}}\cdot\varepsilon^{-(d-1)|H|}\cdot\varepsilon^{(2\upsilon/5)+(5C_{13}^*)^{-1}\cdot \rho}\cdot|\log\varepsilon|^{C^*\rho}.\] Since $C\tau\leq\tau^{1/8}$ and $|\log\varepsilon|^{C^*}\varepsilon^{\upsilon/2}\leq 1$, this proves (\ref{eq.proof_int_est_from_alg_step2_1}) and thus completes the proof of Proposition \ref{prop.cumulant_est} for $E_H$. 

\textbf{Proof part 6: the case of $\mathrm{Err}^1$.} Finally, to prove Proposition \ref{prop.cumulant_est} for $\mathrm{Err}^1$, we only need to prove (\ref{eq.proof_int_est_from_alg_step0_*}) with an extra factor of $\varepsilon^{-C^*|H|}$ on the right hand side, i.e. \begin{equation}\label{eq.proof_int_est_from_alg_step5_1}
    \Ic_{\Mb'}(Q')\leq \tau^{|\Mb|/9}\cdot \varepsilon^{-C^*|H|}\cdot \varepsilon^{- (d - 1) |H| + (2\upsilon/5)+(5C_{13}^*)^{-1}\cdot \rho}\cdot |\log\varepsilon|^{C^*\rho}.
\end{equation} However, this follows the same arguments as in Parts 1--5 above; the only difference is that we now have an extra term $-C^*|H|$ on the right hand side of \eqref{eq.overall_alg}, which leads to the loss $\varepsilon^{-C^*|H|}$. This proves Proposition \ref{prop.cumulant_est} for $\mathrm{Err}^1$, and completes the proof of Proposition \ref{prop.cumulant_est}.
\end{proof}

\subsection{Auxiliary results in the proof of Proposition \ref{prop.cumulant_est}} \label{sec.aux}In this subsection we state and prove the auxiliary results (Propositions \ref{prop.weight} and \ref{prop.volume}) that are used in the proof of Proposition \ref{prop.cumulant_est} in Section \ref{sec.summary}.

\begin{proposition}\label{prop.weight} Let $\Mb$ be a molecule as in Definition \ref{def.set_T_F}, which satisfies Definition \ref{def.set_T_F} \eqref{it.set_F_l_1} and \eqref{it.set_F_l_3}. Let $v_e\in\Rb^d$ be the associated variables defined for each edge $e$ of $\Mb$, satisfying
\begin{equation}\label{eq.energy_con}
    |v_e|\leq|\log\varepsilon|^{C^*},\quad|v_{e_1}|^2+|v_{e_2}|^2=|v_{e_1'}|^2+|v_{e_2'}|^2;\quad v_{e_j}=v_{e_j'}\ \textrm{if }\nf\ \textrm{is O-atom}\end{equation} 
for each atom $\nf$, where $(e_1,e_2)$ are the two bottom edges and $(e_1',e_2')$ are the two top edges at $\nf$, and $e_1$ is serial with $e_1'$ and $e_2$ is serial with $e_2'$ as in Definition \ref{def.associated_op} \eqref{it.associated_dist}. Fix also a constant $\gamma$ with $|\log\varepsilon|^{-1}\leq\gamma\leq 1$. Then there exists a scalar quantity $\Bf\ge 1$ depending only on $\Mb$, such that
\begin{equation}\label{eq.weight_est_1}
\prod_{\nf}(1+|v_{e_1}-v_{e_2}|)\leq \Bf\cdot \exp\bigg(\gamma\sum_{e\in\Ec_{\mathrm{end}}^-}|v_e|^2\bigg),
\end{equation} 
where the product on the left hand side of \eqref{eq.weight_est_1} is taken over all atoms $\nf$ and $(v_{e_1},v_{e_2})$ is as above. Moreover, this quantity $\Bf$ satisfies the following estimates. For any atom set $A\subseteq\Mb$, we have
\begin{equation}\label{eq.weight_est_2}
\Bf\cdot \bigg(\int_{\widetilde{\Dc}}1\,\mathrm{d}\vt_{\Mb\backslash A}\bigg)\leq \big(C\gamma^{-1/2}\tau\big)^{|\Mb|}\cdot|\log\varepsilon|^{C^*(\rho + |A|)},
\end{equation} 
where the integral is taken over all variables $\vt_{\Mb\backslash A}=(t_\nf)_{\nf\in\Mb\backslash A}$, and $\widetilde{\Dc}$ is the projection of the set $\Dc$ defined in \eqref{eq.associated_domain} onto the variables $\vt_{\Mb\backslash A}$.
\end{proposition}
\begin{proof} We divide the proof into 4 parts.

\textbf{Proof part 1: the binary tree case.} In this part, we consider the case when $\Mb=T$ has only one layer $\ell'$ and is a \emph{binary tree}, where $T$ has a unique highest atom $\mf$, and each atom in $T\backslash\{\mf\}$ has a unique parent in $T$.

\textbf{Proof of \eqref{eq.weight_est_1} for binary trees.} By considering $\gamma^{1/2}\cdot v_e$ instead of $v_e$, we may reduce to the case $\gamma=1$. For the domain $\Dc_T$ of time variables $(t_\nf)_{\nf\in T}$ defined in the same way as in \eqref{eq.associated_domain}, we know $t_\nf\in[(\ell'-1)\tau,\ell'\tau]$ because $T$ only has one layer. It is easy to prove by induction that
\begin{equation}\label{eq.weight_est_proof_1}
    \int_{\Dc_T}1\,\mathrm{d}\vt_{T}=\tau^{|T|}\cdot\Cf_T^{-1},\quad \Cf_T:=\prod_{\nf\in T}\sigma(\nf),
\end{equation} 
where $\sigma(\nf)$ is the number of descendant atoms of $\nf$ in $T$ (Definition \ref{def.molecule_order}; this includs $\nf$ itself).

For any atom $\nf\in T$, by iteratively applying \eqref{eq.energy_con} for the atom $\nf$ and its descendants in $T$, we get that
\begin{equation}\label{eq.weight_est_proof_2}
    (1+|v_{e_1}-v_{e_2}|)^2\leq 2+2(|v_{e_1}|^2+|v_{e_2}|^2)\leq 2+2\sum_{(\lf,e)}|v_e|^2,
\end{equation} 
where the sum is taken over all lowest atoms $\lf$ of $T$ that are descendants of $\nf$, and all bottom ends $e$ at $\lf$. Note that for each fixed $e$ there is at most one $\lf$, and for each fixed $\lf$ there are at most two $e$ in the summation. From \eqref{eq.weight_est_proof_2} we have
\begin{equation}\label{eq.weight_est_proof_3}
\begin{aligned}
\prod_{\nf\in T}(1+|v_{e_1}-v_{e_2}|)^2&\leq \prod_{\nf\in T}(\sigma(\nf))^{5/3}\cdot\prod_{\nf\in T}\bigg(\frac{2}{(\sigma(\nf))^{5/3}}+\frac{2}{(\sigma(\nf))^{5/3}}\sum_{(\lf,e)}|v_e|^2\bigg)\\
&\leq \Cf_T^{5/3}2^{|T|}\cdot\prod_{\nf\in T}\bigg(1+(\sigma(\nf))^{-5/3}\bigg(\sum_{(\lf,e)}|v_e|^{4/3}\bigg)^{3/2}\bigg).
\end{aligned}
\end{equation} 
where $\Cf_T$ is as in (\ref{eq.weight_est_proof_1}), and in the last line we used $\sum x_j^{3/2}\leq (\sum x_j)^{3/2}$ for $x_j\geq 0$. By elementary inequality, for any $Z>0$ we have $1+Z^{3/2}\leq C e^Z$; applying this to $Z:=(\sigma(\nf))^{-10/9}\sum_{(\lf,e)}|v_e|^{4/3}$, we get from \eqref{eq.weight_est_proof_3} that
\begin{equation}\label{eq.weight_est_proof_4}
\begin{aligned}
\prod_{\nf\in T}(1+|v_{e_1}-v_{e_2}|)^2&\leq \Cf_T^{5/3}(2C)^{|T|}\exp\bigg(\sum_{\nf\in T}(\sigma(\nf))^{-10/9}\sum_{(\lf,e)}|v_e|^{4/3}\bigg)\\
&= \Cf_T^{5/3}(2C)^{|T|}\exp\bigg(\sum_{(\lf,e)}|v_e|^{4/3}\sum_{\nf}(\sigma(\nf))^{-10/9}\bigg).
\end{aligned}
\end{equation} 
Here, for fixed $(\lf,e)$, the summation is taken over all $\nf\in T$ that are ancestors of $\lf$. Note that since $T$ is a binary tree, each atom in $T$ has at most one parent atom in $T$, so all such $\nf$ (for fixed $\lf$) form a monotonic sequence on which the value of $\sigma(\nf)$ is strictly increasing. This implies that $\sum_{\nf}(\sigma(\nf))^{-10/9}\leq\sum_{n=1}^\infty n^{-10/9}\leq C$ uniformly in $(\lf,e)$, and hence
\begin{equation}\label{eq.weight_est_proof_5}
\prod_{\nf\in T}(1+|v_{e_1}-v_{e_2}|)\leq \Cf_T^{5/6}C^{|T|}\cdot\exp\bigg(C\sum_{e\in \Ec_T^-}|v_e|^{4/3}\bigg),
\end{equation} 
where $\Ec_T^-$ is the set of all bottom ends of atoms in $T$.

\textbf{Analog of \eqref{eq.weight_est_2} for binary trees.} In the same setting as above (where $\Mb=T$ is a single layer binary tree), now let $A\subseteq T$ be a fixed atom set. Next we shall prove that
\begin{equation}\label{eq.weight_est_proof_6}
    \int_{\widetilde{\Dc}_T}1\,\mathrm{d}\vt_{T\backslash A}\leq \tau^{|T|-|A|}\cdot \Cf_T^{-1}\cdot\frac{(|T|)!}{(|T|-|A|)!},
\end{equation} 
where $\widetilde{\Dc}_T$ is projection of $\Dc_T$ (as defined above) onto the variables $\vt_{T\backslash A}$. We will prove \eqref{eq.weight_est_proof_6} by induction; the base case is trivial.

Suppose \eqref{eq.weight_est_proof_6} is true for $|T|<m$, now consider some $T$ with $|T|=m$. By time translation and dilation we may assume $\tau=1$ and $\ell'=0$ so $t_\nf\in[0,1]$ for each $\nf$. Let the unique highest atom of $T$ be $\mf$, we will assume $\mf$ has 2 children atoms; the case when $\mf$ has only one child atom is similar with slightly adjusted proof. Let the two children atoms of $\mf$ be $\mf_1$ and $\mf_2$, and let $T_j$ be the set of descendants of $\mf_j$ with $|T_j|=m_j$, then $m=m_1+m_2+1$. Let $A_j=A\cap T_j$ and $|A_j|=k_j$, we consider two cases, namely when $\mf\in A$ or $\mf\not\in A$.

Suppose $\mf\not\in A$, then $0\leq k_j\leq m_j$ and $|A|:=k=k_1+k_2$. For any $\vt_{T\backslash A}\in \widetilde{\Dc}_T$, we must have $\vt_{T_j\backslash A_j}\in \widetilde{\Dc}_{T_j}$ for $j\in\{1,2\}$, and $t_\nf\in[0,t_\mf]$ for all $\nf\in T_j\backslash A_j$. By applying induction hypothesis to $(T_j,A_j)$ and then integrating in $t_\mf$, we get that
\begin{multline}\label{eq.weight_est_proof_7}\int_{\widetilde{\Dc}_T}1\,\mathrm{d}\vt_{T\backslash A}\leq \Cf_{T_1}^{-1}\Cf_{T_2}^{-1}\frac{m_1!m_2!}{(m_1-k_1)!(m_2-k_2)!}\int_0^1t_\mf^{m_1+m_2-k_1-k_2}\,\mathrm{d}t_\mf\\\leq \Cf_{T_1}^{-1}\Cf_{T_2}^{-1}\frac{m_1!m_2!}{(m_1-k_1)!(m_2-k_2)!}\frac{1}{m_1+m_2-k_1-k_2+1}.\end{multline}
 Note also that
\begin{equation}\label{eq.weight_est_proof_8}\Cf_T^{-1}=\Cf_{T_1}^{-1}\Cf_{T_2}^{-1}\cdot\frac{1}{m_1+m_2+1},
\end{equation} we see that \eqref{eq.weight_est_proof_6} will follows from \eqref{eq.weight_est_proof_7} and \eqref{eq.weight_est_proof_8}, provided we can show that
\begin{equation}\label{eq.weight_est_proof_9}\frac{m_1!m_2!}{(m_1-k_1)!(m_2-k_2)!}\leq\frac{(m_1+m_2)!}{(m_1+m_2-k_1-k_2)!}.
\end{equation} Now to prove \eqref{eq.weight_est_proof_9}, simply note that $(x+y)!/x!$ is increasing in $x\geq 0$ for fixed $y\geq 0$; by applying it to $y=m_2-k_2$ and $x=m_1-k_1$ (and using symmetry), we see that the ratio between the left hand side and right hand side of \eqref{eq.weight_est_proof_9} is maximized when $x$ and $y$ are both maximized, i.e. $k_1=k_2=0$. In this extremal  case \eqref{eq.weight_est_proof_9} is trivially true, so it is true in general.

Suppose $\mf\in A$, then $0\leq k_j\leq m_j$ and $|A|:=k=k_1+k_2+1$. By similar argument as above (except that now we do not integrate in $t_\mf$ as $\mf\in A$), we get
\begin{equation}\label{eq.weight_est_proof_10}\int_{\widetilde{\Dc}_T}1\,\mathrm{d}\vt_{T\backslash A}\leq \Cf_{T_1}^{-1}\Cf_{T_2}^{-1}\frac{m_1!m_2!}{(m_1-k_1)!(m_2-k_2)!}.\end{equation} Using also \eqref{eq.weight_est_proof_8} and \eqref{eq.weight_est_proof_9}, this again implies \eqref{eq.weight_est_proof_6}.

\textbf{Proof part 2: the forest case.} Now consider the case when $\Mb$ has only one layer $\ell'$, and is a forest. By restricting to components, we may assume $\Mb$ is a tree (which need not be a binary tree). We shall divide the atom set of $\Mb$ into finitely many binary trees and upside down binary trees (where the notions of parent and child etc. are the opposite from binary trees), as follows, see {\color{blue}Figure \ref{fig.tree_decompose}}.

First choose any atom $\nf\in\Mb$, and let $T_{\nf}$ be the set of descendants of $\nf$. This is a binary tree, because each atom has at most one parent in $T_\nf$ (as $\Mb$ has no cycle). Moreover, all bottom ends of $T_\nf$ are also bottom ends of $\Mb$. Consider all the atoms $\nf_j\not\in T_\nf$ that are connected to $T_\nf$ by a bond. Since $\Mb$ is a tree, if we remove $T_\nf$ from $\Mb$ (and turn any bond between $\nf_j$ and $T_\nf$ to a free end at $\nf_j$), then $\Mb$ becomes the disjoint union of trees $\Mb^j$, where $\nf_j\in \Mb^j$ for each $j$. Note also that each $\nf_j$ has only one bond connected to $T_\nf$, which is turned to a bottom free end in $\Mb^j$.

Now for each $j$, let $T_{\nf_j}^*$ be the set of ancestors of $\nf_j$ in $\Mb^j$, which is now an upside down binary tree. Clearly all top ends of $T_{\nf_j}^*$ are also top ends of $\Mb$, and all the discussions in Step 1 above also work for upside down binary trees due to time reversal symmetry. Consider all the atoms $\nf_{jk}\in \Mb^j\backslash T_{\nf_j}^*$ that are connected to $T_{\nf_j}^*$ by a bond, then after removing $T_{\nf_j}^*$ from $\Mb^j$, the $\Mb^j$ will become the disjoint union of trees $\Mb^{jk}$ where $\nf_{jk}\in\Mb^{jk}$ for each $k$. Also each $\nf_{jk}$ has only one bond connected to $T_{\nf_j}^*$, which is turned to a top free end in $\Mb^{jk}$. Therefore, we can choose $T_{\nf_{jk}}$ to be the set of descendants of $\nf_{jk}$ in $\Mb^{jk}$, which is a binary tree, so that all bottom ends of $T_{\nf_{jk}}$ are also bottom ends of $\Mb$, and so on.
\begin{figure}[h!]
\includegraphics[width=0.35\linewidth]{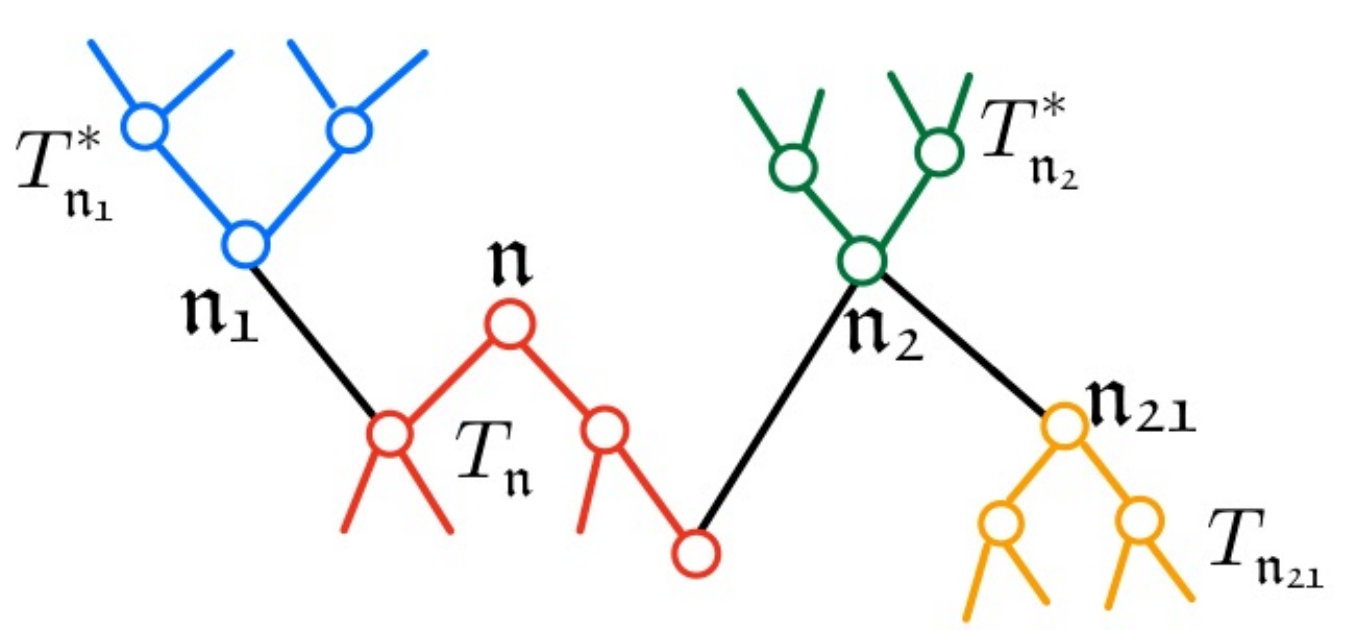}
\caption{The decomposition of a general tree $\Mb$ into binary trees and upside down binary trees, as described in Part 2 of the proof of Proposition \ref{prop.weight}. Here $\nf$, $\nf_j$, $\nf_{jk}$ etc. and $T_\nf$, $T_{\nf_j}^*$, $T_{\nf_{jk}}$ etc., are defined as in Part 2. Different colors represent different binary trees and upside down binary trees.}
\label{fig.tree_decompose}
\end{figure}

By repeating the above, we have divided the atom set of $\Mb$ into the disjoint union of finitely many $T^{(j)}$, where each $T^{(j)}$ is either a binary tree or an upside down binary tree. If we apply \eqref{eq.weight_est_proof_5} and \eqref{eq.weight_est_proof_6} to each $T^{(j)}$, we then obtain that
\begin{equation}\label{eq.weight_est_proof_11}
\prod_{\nf\in \Mb}(1+|v_{e_1}-v_{e_2}|)\leq \prod_{j}\Cf_{T^{(j)}}^{5/6}\cdot C^{|\Mb|}\cdot\exp\bigg(C\sum_{e\in \Ec_{\mathrm{end}}}|v_e|^{4/3}\bigg),
\end{equation}
\begin{equation}\label{eq.weight_est_proof_12}
\int_{\widetilde{\Dc}}1\,\mathrm{d}t_{\Mb\backslash A}\leq \tau^{|\Mb|-|A|}\prod_{j}\Cf_{T^{(j)}}^{-1}\cdot\prod_{j}\frac{(|T^{(j)}|)!}{(|T^{(j)}|-|A\cap T^{(j)}|)!}\leq  \tau^{|\Mb|}\prod_{j}\Cf_{T^{(j)}}^{-1}\cdot |\log\varepsilon|^{C^*|A|},
\end{equation} in view of the upper bound $|\Mb|\leq |\log\varepsilon|^{C^*}$ (see Part 1 of the proof of Proposition \ref{prop.layerrec3}) and the trivial bound $\tau^{-1}\leq|\log\varepsilon|$. Here $\Ec_{\mathrm{end}}$ in \eqref{eq.weight_est_proof_11} is the set of all ends of $\Mb$, and $\widetilde{\Dc}$ in \eqref{eq.weight_est_proof_12} is the same as in \eqref{eq.weight_est_2}. 

\textbf{Proof part 3: single layer case.} Now consider a molecule $\Mb=\Mb_{\ell'}$ with a single layer $\ell'$ satisfying the relevant conditions in Definition \ref{def.set_T_F} \eqref{it.set_F_l_1} and \eqref{it.set_F_l_3}.

Consider the clusters $\Mb_j$ of $\Mb$ (Definition \ref{def.cluster}). For each $j$, we may choose at most $\rho(\Mb_j)\leq\Gamma$ C-atoms in $\Mb_j$, such that if we remove these C-atoms, and turn any bond between these atoms and other atoms $\nf$ into a free end at $\nf$, then the recollision number $\rho(\Mb_j)$ will be reduced to $0$. Therefore, in the molecule $\Mb_{\ell'}$, if we remove all such C-atoms (their total number is $\leq \sum_{j}\rho(\Mb_j)=\Rf_{\ell'}$, see Definition \ref{def.parameter_rho_old}), then the resulting molecule $\Mb_{\ell'}$ will be a forest. By applying \eqref{eq.weight_est_proof_11} and \eqref{eq.weight_est_proof_12} to this forest, and using \eqref{eq.energy_con} for each of the removed atoms, we get that
\begin{equation}\label{eq.weight_est_proof_13}
\prod_{\nf\in \Mb_{\ell'}}(1+|v_{e_1}-v_{e_2}|)\leq \Cf_{\ell'}^{5/6}\cdot C^{|\Mb_{\ell'}|}|\log\varepsilon|^{C^*\Rf_{\ell'}}\cdot\exp\bigg(C\sum_{e\in \Ec_{\mathrm{end}}(\Mb_{\ell'})}|v_e|^{4/3}\bigg)\cdot\exp\bigg(C\sum_{e\in \Ec_{\ell'}^*}|v_e|^{4/3}\bigg),
\end{equation}
\begin{equation}\label{eq.weight_est_proof_14}
\int_{\widetilde{\Dc}}1\,\mathrm{d}t_{\Mb_{\ell'}\backslash A}\leq \tau^{|\Mb_{\ell'}|}\cdot \Cf_{\ell'}^{-1}\cdot |\log\varepsilon|^{C^*|A\cap \Mb_{\ell'}|}.
\end{equation} 
Here $\Cf_{\ell'}:=\prod_j\Cf_{T^{(j)}}$ as in \eqref{eq.weight_est_proof_11} and \eqref{eq.weight_est_proof_12} (where the decomposition is made for the resulting forest after removing these C-atoms, following Part 2), and the other notions are defined as above but associated with $\Mb_{\ell'}$. Moreover in \eqref{eq.weight_est_proof_13}, $\Ec_{\ell'}^*$ is the set of all bonds at all the C-atoms removed above, so that each such bond becomes a free end in the reduced molecule $\Mb_{\ell'}$ (which is a forest) after removing these C-atoms.

Consider each $e\in \Ec_{\ell'}^*$, which is at some C-atom $\nf$ removed above. Assume that this C-atom $\nf$ belongs to the cluster $\Mb_j$. By iterating \eqref{eq.energy_con} and noticing that all collisions happen between particles in $p(\Mb_j)$, we see that
\begin{equation}\label{eq.weight_est_proof_15}
\sum_e|v_e|^2\leq 2 \sum_{e\in \Ec_{\mathrm{end}}(\Mb_j)}|v_e|^2,\quad\mathrm{and\ hence}
\end{equation}
\begin{equation}\label{eq.weight_est_proof_15+}
\sum_e|v_e|^{4/3}\leq C \sum_{e\in \Ec_{\mathrm{end}}(\Mb_j)}|v_e|^{4/3},
\end{equation} 
where the left hand side of \eqref{eq.weight_est_proof_15} and \eqref{eq.weight_est_proof_15+} is the sum over (at most 4) bonds $e$ at some fixed atom $\nf$ as above. Notice that the each end in $\Ec_{\mathrm{end}}(\Mb_j)$ naturally corresponds to an end in $\Ec_{\mathrm{end}}(\Mb_{\ell'})$ (identify $e\in \Ec_{\mathrm{end}}(\Mb_j)$ with the bottom end of the particle line in $\Mb_{\ell'}$ that $e$ belongs to), and each $\Mb_j$ has at most $\Gamma\leq C$ different atoms $\nf$ removed. We obtain, upon summing over $\nf$ and then over $j$, that
\[
\sum_{e\in\Ec_{\ell'}^*}|v_e|^{4/3}\leq C\sum_{j}\sum_{e\in \Ec_{\mathrm{end}}(\Mb_j)}|v_e|^{4/3}\leq C\sum_{e\in \Ec_{\mathrm{end}}(\Mb_{\ell'})}|v_e|^{4/3},
\] 
so the last exponential factor on the right hand side of \eqref{eq.weight_est_proof_13} can be absorbed by the one before it.

\textbf{Proof part 4: general multi-layer case.} Now we put together all layers $\Mb_{\ell'}$ and consider the general multi-layer molecule $\Mb$. By applying \eqref{eq.weight_est_proof_13}--\eqref{eq.weight_est_proof_14} for each $\ell'$ we get that
\begin{equation}\label{eq.weight_est_proof_16}
\prod_{\nf\in \Mb}(1+|v_{e_1}-v_{e_2}|)\leq \Cf^{5/6}\cdot C^{|\Mb|}|\log\varepsilon|^{C^*\rho}\cdot\exp\bigg(C\sum_{\ell'}\sum_{e\in \Ec_{\mathrm{end}}(\Mb_{\ell'})}|v_e|^{4/3}\bigg),
\end{equation}
\begin{equation}\label{eq.weight_est_proof_17}
\int_{\widetilde{\Dc}}1\,\mathrm{d}t_{\Mb\backslash A}\leq \tau^{|\Mb|}\cdot \Cf^{-1}\cdot |\log\varepsilon|^{C^*|A|},
\end{equation} 
where $\Cf:=\prod_{\ell'}\Cf_{\ell'}$ and the other notions are defined above for $\Mb$.

Now for each $(e,\ell')$ such that $e\in\Ec_{\mathrm{end}}(\Mb_{\ell'})$, this $e$ is either an end in the molecule $\Mb$ (which we call case 1), or a bond connecting layer $\ell'$ and some other layer $<\ell'$ (which we call case 2). For case 1, recall $\gamma=1$, we simply apply $\exp(C|v|^{4/3})\leq C\cdot \exp(|v^2|/20)$ to get
\begin{equation}\label{eq.weight_est_proof_18}\exp\bigg(C\sum_{\textrm{case\ 1}\textrm{ (for all }\ell')}|v_e|^{4/3}\bigg)\leq C^{|\Mb|}\cdot \exp\bigg(\frac{1}{20}\sum_{e\in\Ec_{\mathrm{end}}}|v_e|^2\bigg)=C^{|\Mb|}\cdot \exp\bigg(\frac{1}{10}\sum_{e\in\Ec_{\mathrm{end}}^-}|v_e|^2\bigg),
\end{equation} 
where $\sum_{e\in\Ec_{\mathrm{end}}}|v_e|^2 = 2\cdot \sum_{e\in\Ec_{\mathrm{end}}^-}|v_e|^2$ by energy conservation. For case 2, any such bond $e$ belongs to a particle line $p\in r(\Mb_{\ell'-1})$. Note also that the number of all such particles $p$ is bounded by
\begin{equation}\label{eq.weight_est_proof_20}
    \sum_{\ell'}|r(\Mb_{\ell'})|=\sum_{\ell'}s_{\ell'}\le \rho,
\end{equation} 
following Definition \ref{def.parameter_rho_old}. Now we apply the inequality $\exp(C|v|^{4/3})\leq C^*\cdot \exp(|v|^2/20d\Lf)$ (recall our convention of $C$ and $C^*$ in Definition \ref{def.notation} (\ref{it.defC*})), and use also \eqref{eq.weight_est_proof_20}, to obtain that
\begin{equation}\label{eq.weight_est_proof_21}
    \exp\bigg(C\sum_{\textrm{case\ 2}\textrm{ (fix }\ell')}|v_e|^{4/3}\bigg)\leq (C^*)^{C^*\rho}\cdot\exp\bigg(\frac{1}{20d\Lf}\sum_{\textrm{case\ 2}\textrm{ (fix }\ell')}|v_e|^2\bigg).
\end{equation} 
By energy conservation, we have
\begin{equation}\label{eq.weight_est_proof_19}
    \exp\bigg(\frac{1}{20d\Lf}\sum_{\textrm{case\ 2}\textrm{ (fix }\ell')}|v_e|^2\bigg)\leq \exp\bigg(\frac{1}{20d\Lf}\sum_{e\in\Ec_{\mathrm{end}}^-}|v_e|^2\bigg)
\end{equation} 
Multiplying all possible $\ell'$ yields
\begin{equation}\label{eq.weight_est_proof_22}
    \exp\bigg(C\sum_{\textrm{case\ 2}\textrm{ (for all }\ell')}|v_e|^{4/3}\bigg)\leq (C^*)^{C^*\rho}\cdot\exp\bigg(\frac{1}{10d}\sum_{e\in\Ec_{\mathrm{end}}^-}|v_e|^{2}\bigg).
\end{equation} Now, we may put together \eqref{eq.weight_est_proof_16},\eqref{eq.weight_est_proof_18} and \eqref{eq.weight_est_proof_22}, and choose $\Bf$ to be the factor $\Cf^{5/6}\cdot C^{|\Mb|}|\log\varepsilon|^{C^*\rho}$ in \eqref{eq.weight_est_proof_16}, multiplied by the $C^{|\Mb|}$ and $(C^*)^{C^*\rho}$ factors in \eqref{eq.weight_est_proof_18}, \eqref{eq.weight_est_proof_21} and \eqref{eq.weight_est_proof_22}. It is then clear that $\Bf\geq \Cf^{6/5}\geq 1$, and \eqref{eq.weight_est_1} and \eqref{eq.weight_est_2} are true with these choices. This completes the proof.
\end{proof}
\begin{proposition}\label{prop.volume} Let $\Mb$ be a molecule, and fix a cutting sequence that reduces $\Mb$ to $\Mb'$, which contains only empty ends and elementary molecules. Let $\Sf$ be a subset of the set of \{3\} molecules in $\Mb'$. For each atom $\nf$ with $\{\nf\}\in\Sf$, it has a unique fixed end $e_1(\nf)$, and a unique free end $e_2(\nf)$ that is either both top or both bottom $e_1(\nf)$. Let the collection of all such $e_2(\nf)$ be $\Fc$.

Now suppose for each atom $\nf$ with $\{\nf\}\in \Sf$ we have fixed $\omega_\nf\in\Sb^{d-1}$. Assume also $v_e\in\Rb^d$ is the associated variable for each edge $e\in\Mb$, respectively. Moreover, assume $|v_e|\leq|\log\varepsilon|^{C^*}$ for each edge $e$ of $\Mb'$, and for each atom $\nf$ with $\{\nf\}\in\Sf$ and edges $(e_1,e_2,e_1',e_2')$ as in Definition \ref{def.associated_int}, assume also that
\begin{equation}\label{eq.volume_preserve}
    v_{e_1'}=v_{e_1}-[(v_{e_1}-v_{e_2})\cdot\omega_\nf]\omega_\nf,\quad v_{e_2'}=v_{e_2}+[(v_{e_1}-v_{e_2})\cdot\omega_\nf]\omega_\nf
\end{equation} 
(this is for C-atoms; for O-atoms we simply have $v_{e_j'}=v_{e_j}$). For each $e\in\Ec_{\mathrm{end}}(\Mb)$, fix 
\[
    X_e\in [1,|\log\varepsilon|^{C^*}]\cap2^\Zb,
\] 
and assume that $|v_e|\leq X_e$. Finally let $K$ be the number of molecules in $\Mb'$ that do not belong to $\Sf$ (in reality $\Sf$ will be the set of normal \{3\} molecules, so $K=\#_{\neg\mathrm{N}\{3\}}$ using the notion in (\ref{eq.not_3_comp_upp_3})).

Then, there exists a set $Y\subseteq (\Rb^d)^{|\Fc|}$, which depends on $(\omega_\nf)$ and $(X_e)$, and on $(\Mb,\Mb')$ and the cutting sequence, but not on the concrete choices of $(v_e)$, such that
\begin{equation}\label{eq.volume_1}
    (v_e)_{e\in\Fc}\in Y,\quad |Y|\leq \prod_{e\in \Ec_{\mathrm{end}}}X_e^d\cdot C^{|\Mb|}\cdot |\log\varepsilon|^{C^*K}.
\end{equation}
\end{proposition}
\begin{proof} For each atom $\nf$ with $\{\nf\}\in\Sf$, let $(e_1(\nf),e_2(\nf)$ and $(e_1'(\nf),e_2'(\nf)$ be defined as in Proposition \ref{prop.volume} and \eqref{eq.volume_preserve}. Recall $\Fc$ is the set of all $e_2(\nf)$ for $\{\nf\}\in\Sf$. Let $\Gc$ be the set of all edges of molecules in $\Mb'$ that are not in $\Sf$, plus the ends of $\Mb$ (which are preserved in cutting, so they are also free ends of $\Mb'$). Note that $|\Fc|,|\Gc|\leq 10|\Mb|$. The idea is to show that, with $(\omega_\nf)$ fixed and using \eqref{eq.volume_preserve}, each vector $v_{e_2(\nf)}$ can be written as a constant coefficient linear combination of vectors $(v_{e'})_{e'\in\Gc}$, and 
\begin{equation}\label{eq.balance_2}
    \sum_{\nf}|v_{e_2(\nf)}|^2\leq \sum_{e'\in\Gc}|v_{e'}|^2.
\end{equation}

We first prove the linear relation and (\ref{eq.balance_2}) as claimed above. Start from any $e_2(\nf)$. If $\nf$ is O-atom then let $e=e_2(\nf)$; if $\nf$ is C-atom then we apply (\ref{eq.volume_preserve}) to get a linear combination of $v_{e_1'(\nf)}$ and $v_{e_2'(\nf)}$ (and the energy inequality $|v_{e_2(\nf)}|^2\leq|v_{e_1'(\nf)}|^2+|v_{e_2'(\nf)}|^2$), and then let $e$ be either $e_1'(\nf)$ or $e_2'(\nf)$. If $e\in\Gc$ then this is already what we need; if $e\not\in\Gc$ we need to further explore from $e$ to reach $\Gc$, as follows.

Note that $e$ is not an end in $\Mb$ (as $e\not\in\Gc$), so it must occur during a cutting operation, say while operating on a maximal ov-segment $\sigma$ (with associated notions $\nf^\pm$ and $\pf_j$) as in Definition \ref{def.cutting}. It is easy to verify that $e$ does not belong to a simple pair (Definition \ref{def.reg}), so by Definition \ref{def.cutting} (\ref{it.cutting_2}), we see that the only way this cutting creates a new free end $e$ is that $\pf_1\neq\nf^+$ or $\pf_s\neq \nf^-$ (for example in {\color{blue}Figure \ref{fig.cutting2+}}, we get a new free end $e_5'$ due to $\pf_3\neq \nf^-$). Assume the former by symmetry, then $\nf^+\in\sigma\backslash A$, so the same cutting creates a fixed end $e_1$ at $\nf^+$, with $v_e=v_{e_1}$. If $\{\nf^+\}\not\in\Sf$, then we have $e_1\in\Gc$ by definition, which is what we need; thus we will assume below that $\{\nf^+\}\in\Sf$.

if $\nf^+$ is O-atom, then $e_1$ must be serial with a free end $e'$, so $v_e=v_{e'}$. If $\nf^+$ is C-atom, then $e_1=e_1(\nf^+)$, so we can add in $e_2(\nf^+)\in\Fc$: both $v_e$ and $v_{e_2(\nf^+)}$ are linear combinations of $v_{e_1'(\nf^+)}$ and $v_{e_2'(\nf^+)}$, and $|v_e|^2+|v_{e_2(\nf^+)}|^2=|v_{e_1'(\nf^+)}|^2+|v_{e_2'(\nf^+)}|^2$, then we can choose $e'$ be either $e_1'(\nf^+)$ or $e_2'(\nf^+)$. In either case, if $e'\in\Gc$ then this is what we need; if not then we can explore from $e'$ as above. In the end, every exploration will stop at some $e'\in\Gc$, which proves the desired linear relation and (\ref{eq.balance_2}).

Now, by (\ref{eq.balance_2}), we have $(v_e)_{e\in\Fc}=\Pc\big((v_{e'})_{e'\in\Gc}\big)$ for some constant matrix $\Pc$ which matrix norm $\|\Pc\|\leq 1$.
Note that $|v_{e'}|\leq X_{e'}$ is $e'$ is an end in $\Mb$, and $|v_{e'}|\leq |\log\varepsilon|^{C^*}$if $e'$ belongs to a molecule in $\Mb'$ that is not in $\Sf$, and the number of the latter case does not exceed $10K$. It follows that the vector $(v_{e'})_{e'\in\Gc}$ belongs to the union of at most
\[
A:=C^{|\Gc|}|\log\varepsilon|^{C^*K}\prod_{e\in\Ec_{\mathrm{end}}}X_{e}^d
\] 
fixed unit boxes in $\Rb^{d|\Gc|}$. However, the image of any unit box in $\Rb^{d|\Gc|}$ under $\Pc$ is contained in a ball of radius $(d|\Gc|)^{1/2}$ in $\Rb^{d|\Fc|}$, which is covered by at most \[C^{F}\cdot (F+G)^{F/2}\cdot \mathrm{Vol}\,(\textrm{unit\ ball\ in\ }\Rb^{F})=C^{F} \cdot (F+G)^{F/2}\cdot\pi^{F/2}\cdot \big(\Gamma(F/2+1)\big)^{-1}\leq C^{F+G}\] unit boxes, where $G:=d|\Gc|$ and $F:=d|\Fc|$. It then follows that $(v_e)_{e\in\Fc}$ belongs to the union of at most $C^{d(|\Fc|+|\Gc|)}\cdot A$ fixed unit boxes in $\Rb^{d|\Fc|}$, so \eqref{eq.volume_1} follows.
\end{proof}

\section{Preparations for the algorithm}\label{sec.prepare} From this section until Section \ref{sec.maincr}, we will prove Proposition \ref{prop.comb_est}, the main goal of the third part of the paper. In this section, we make some preparations for the proof by restating the properties of the molecule $\Mb$ in Proposition \ref{prop.comb_est} (originally stated in Definition \ref{def.set_T_F}) in a more convenient way for the algorithm. 

Start from the initial links $\Lc$; we shall replace it by some new initial link set $\widetilde{\Lc}$ where each particle line belongs to \emph{at most two} new initial links in $\widetilde{\Lc}$. This is because, see the proof of Proposition \ref{prop.comb_est_case4}, for any initial link, we will gain a power of $\varepsilon$ from the particle line in it that is cut \emph{after} the other; then the condition for $\widetilde{\Lc}$ ensures that the number of gains from this is comparable with the number of initial links.
\begin{proposition}\label{prop.init_link} Recall the set $\Lc$ of initial links (see Proposition \ref{prop.initial_cumulant}), let the set of particle lines involved in $\Lc$ be $H_0$. Then we have
\begin{equation}\label{eq.1_L_2}\mathbbm{1}_\Lc^\varepsilon\leq\sum_{\widetilde{\Lc}}\mathbbm{1}_{\widetilde{\Lc}}^{\varepsilon+}.
\end{equation} Here $\mathbbm{1}_\Lc^{\varepsilon}$ is defined as in (\ref{eq.1_L'}), $\widetilde{\Lc}$ is a set of pairs $(p_1,p_2)$ where $p_1,p_2\in H_0$ satisfying that each particle line $p\in H_0$ belongs to \emph{one or two} pairs in $\widetilde{\Lc}$, and $\mathbbm{1}_{\widetilde{\Lc}}^{\varepsilon+}$ is defined as 
 \begin{equation}\label{eq.1_L_3}
        \mathbbm{1}_{\widetilde{\Lc}}^{\varepsilon+} = \prod_{(p_1, p_2)\in \widetilde{\Lc}} \mathbbm{1}_{e_1\sim_{\mathrm{in}} p_2}^{\varepsilon+}\qquad\textrm{and}\qquad \mathbbm{1}_{e_1\sim_{\mathrm{in}} p_2}^{\varepsilon+} = \left\{\begin{aligned}
            &1 \quad &&|x_{p_1}-x_{p_2}|\le \varepsilon^{1-\upsilon},
            \\
            &\varepsilon^\upsilon \quad &&|x_{p_1}-x_{p_2}|\ge \varepsilon^{1-\upsilon},
        \end{aligned}\right.
    \end{equation} 
\end{proposition}
\begin{proof} For fixed values of $(x_p)_{p\in H_0}$, we only need to prove that there exists $\widetilde{\Lc}$ such that $\mathbbm{1}_\Lc^\varepsilon(x_p)\leq \mathbbm{1}_{\widetilde{\Lc}}^{\varepsilon+}(x_p)$. Define a graph $\Lc_0$ which is a subgraph of $\Lc$, such that $(p_1,p_2)\in\Lc_0$ if and only if $(p_1,p_2)\in\Lc$ and $|x_{p_1}-x_{p_2}|\leq \varepsilon$. Clearly we may restrict to a single component of the forest $\Lc$ (which is a tree), then in this component, $\Lc_0$ is obtained from $\Lc$ by deleting those edges $(p_1,p_2)$ with $|x_{p_1}-x_{p_2}|\geq\varepsilon$. Let the number of deleted edges be $q$, then $\mathbbm{1}_\Lc^\varepsilon(x_p)=\varepsilon^{\upsilon q}$ by definition. On the other hand, $\Lc_0$ has exactly $q+1$ connected components (as it is formed by deleting $q$ edges from a tree), and for any $(p_1,p_2)$ that belong to the same component, we must have $|x_{p_1}-x_{p_2}|\leq \varepsilon|H_0|\ll \varepsilon^{1-\upsilon}$ by triangle inequality. Now we simply list all the particle lines in $H_0$ by enumerating each component of $\Lc_0$ one after another, and construct $\widetilde{\Lc}$ by forming an edge between any two adjacent particle lines in the above list, then it is easy to see that $|x_{p_1}-x_{p_2}|\geq\varepsilon^{1-\upsilon}$ occurs at most $q$ times (when entering new components), which implies $\mathbbm{1}_{\widetilde{\Lc}}^{\varepsilon+}(x_p)\geq\varepsilon^{\upsilon q}=\mathbbm{1}_\Lc^\varepsilon(x_p)$, as desired.
\end{proof}
\begin{remark}\label{rem.init_link_2} In view of Proposition \ref{prop.init_link}, below we will replace $\Lc$ by $\widetilde{\Lc}$ in all instances, and refer to $\widetilde{\Lc}$ as the (new) \textbf{initial links}, and change Definition \ref{def.good_normal} etc. using these new initial links. Note that (i) the number of choices of $\widetilde{\Lc}$ (i.e. the number of terms in the summation (\ref{eq.1_L_3})) is bounded by $(C|H_0|)!$. Since $|H_0|\leq\rho\leq|\log\varepsilon|^{C^*}$, this resulting combinatorial factor is bounded by $|\log\varepsilon|^{C^*\rho}$ which can be absorbed by Proposition \ref{prop.layerrec3}. (ii) The proof of Proposition \ref{prop.cumulant_est} in Section \ref{sec.treat_integral} is not affected by these new initial links. For example, in Proposition \ref{prop.int_mini_good} (\ref{it.bound_one_atom_good_empty_line}), we are only assuming the weaker assumption $|x-x'|\leq\varepsilon^{1-\upsilon}$ consistent with $\mathbbm{1}_{\widetilde{\Lc}}^{\varepsilon+}$, so any good molecule with respect to the new initial links still carries the same power gain as in (\ref{eq.Linftybd}).
\end{remark}
\begin{definition}\label{def.layer_interval} Let $\Mb_{\ell'}$ be the layers of $\Mb$. Define the \textbf{layer intervals}
\[\Mb_{[\ell_1:\ell_2]}:=\bigcup_{\ell_1\leq\ell'\leq\ell_2}\Mb_{\ell'},\quad \Mb_{(\ell_1:\ell_2)}:=\bigcup_{\ell_1<\ell'<\ell_2}\Mb_{\ell'}\quad \Mb_{<\ell_1}=\bigcup_{\ell'<\ell_1}\Mb_{\ell'},\quad \textrm{etc.}\] After defining the refined layers $\Mb_{\zeta}^T$ (Definition \ref{def.layer_refine}) below. we will use the same notations $\Mb_{[\zeta_1:\zeta_2]}^T$ etc. for the intervals of refined layers.
\end{definition}

\begin{figure}[h!]
    \centering
    \includegraphics[width=0.25\linewidth]{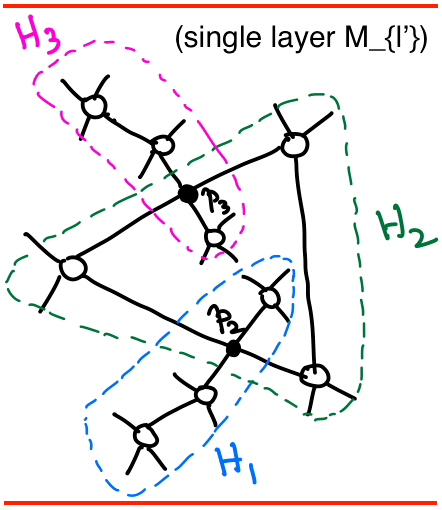}
    \caption{An example of Proposition \ref{prop.mol_axiom} (\ref{it.axiom2}); here $\Qb_j$ are clusters in a single layer $\Mb_{\ell'}$. Since the cluster graph is a forest, it follows that the only cycles in $\Mb_{\ell'}$ are formed by ov-segments between C-atoms in the same cluster (the triangle in $\Qb_2$ formed by ov-segments), so one cannot transit to another ov-segment (i.e. non-serial) at any O-atom ($\pf_2$ or $\pf_3$).}
    \label{fig.no turn at o}
\end{figure}

\begin{figure}[h!]
    \centering
    \includegraphics[width=0.58\linewidth]{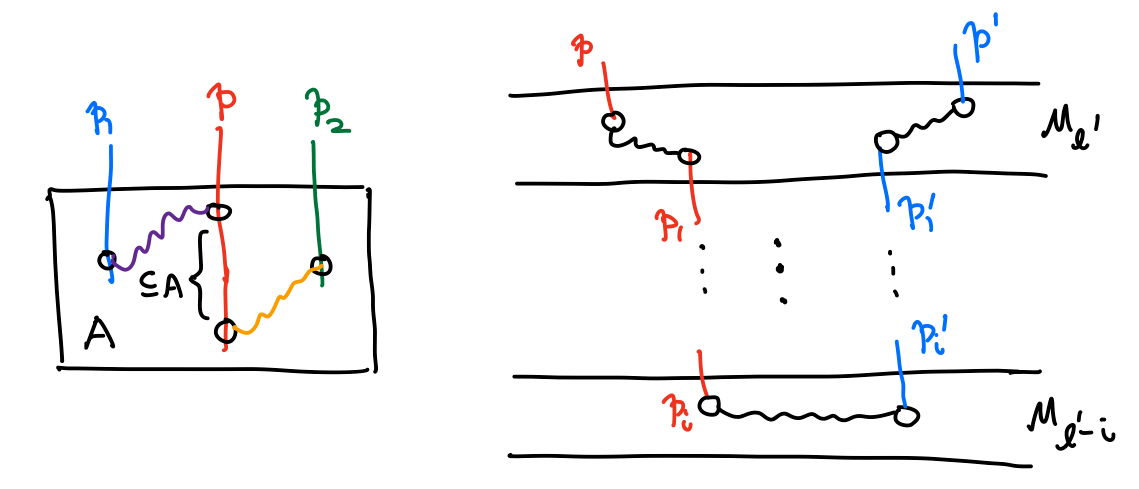}
    \caption{Left: a set $A$ as in Remark \ref{rem.layer_interval} with particle lines $\pb$, $\pb_1$ and $\pb_2$ all connected via $A$; note that the part of $\pb$ between two atoms in $A$ is again $\subseteq A$. Right: the procedure of finding $(\pb_1,\cdots,\pb_i)$ and then $(\pb_i,\cdots,\pb')$ in the proof of Proposition \ref{prop.mol_axiom} (\ref{it.axiom6}).}
    \label{fig.connectedbelow}
\end{figure}

Now we can state the main result of this section, namely Proposition \ref{prop.mol_axiom} below.
\begin{proposition}
\label{prop.mol_axiom} Let $\Mb\in\Fc_{\vLambda_\ell}$ as in Definition \ref{def.set_T_F} (\ref{it.set_F}) and $r(\Mb)=H$. Then we have the following:
\begin{enumerate}
\item\label{it.axiom1-} The total number of atoms and particle lines in $\Mb$ is bounded by $|\log\varepsilon|^{C^*}$.
\item\label{it.axiom1} No two atoms can be connected by two different ov-segments.\footnote{This property is specific to $\Rb^d$, and needs to be appropriately modified for $\Tb^d$.}
\end{enumerate}
Properties (\ref{it.axiom2})--(\ref{it.axiom3}) below correspond to the clusters (Definition \ref{def.cluster}) in each layer $\Mb_{\ell'}\,(1\leq \ell'\leq \ell)$.
\begin{enumerate}[resume]
\item\label{it.axiom2} Any cycle in a single layer $\Mb_{\ell'}$ must be formed by ov-segments connecting C-atoms in the same cluster. In particular, for any cycle in a single layer $\Mb_{\ell'}$ and any O-atom in this cycle, the two bonds of the cycle at this O-atom must be serial. see {\color{blue} Figure \ref{fig.no turn at o}}.
\item\label{it.axiom3} For any cluster $\Qb$ in each layer $\Mb_{\ell'}$, we have that $\rho(\Qb)\leq\Gamma$ (Definition \ref{def.recollision_number}). Recall that a cluster $\Qb$ or component is cyclic (or recollisional) if $\rho(\Qb)>0$, i.e. if it contains a cycle.
\end{enumerate}
Properties (\ref{it.axiom4})--(\ref{it.axiom5}) below correspond to the particle line sets $p(\Mb_{\ell'})$ and $r(\Mb_{\ell'})$ (Definition \ref{def.prsets}).
\begin{enumerate}[resume]
\item\label{it.axiom4} Each particle line $\pb$ must belong to some $p(\Mb_{\ell'})$ for some $\ell'$. If $\pb$ intersects $\Mb_{\ell'}$ (Definition \ref{def.connectedvia}), then $\pb\in p(\Mb_{\ell'})$. If $\pb$ contains an atom in both $\Mb_{\ell_1}$ and $\Mb_{\ell_2}$, then $\pb\in p(\Mb_{\ell'})$ for all $\ell'\in[\ell_1:\ell_2]$.

\item \label{it.axiom4+} We have $r(\Mb_\ell)=H$, and $r(\Mb_{\ell'})=p(\Mb_{\ell'})\cap p(\Mb_{\ell'+1})$ for $1\leq \ell'\leq\ell-1$. Define also $r(\Mb_0)=H_0$.
\item\label{it.axiom5} For $1\leq\ell'\leq\ell$, each component of $\Mb_{\ell'}$ must intersect a particle line in $r(\Mb_{\ell'})$. Each particle line in $p(\Mb_{\ell'})$ is either connected to a particle line in $r(\Mb_{\ell'})$ via $\Mb_{\ell'}$, or belongs to $r(\Mb_{\ell'})$. See {\color{blue} Figure \ref{fig.moleculeprop}}.
\end{enumerate}
Property (\ref{it.axiom7}) below corresponds to the (new) initial links, see Proposition \ref{prop.init_link}.
\begin{enumerate}[resume]
\item\label{it.axiom7} The bottom end of each particle line in $H_0$ is involved in 1 or 2 initial links. Note that initial links can be understood as either between particle lines or between bottom ends (which are in bijection).
\end{enumerate}
Property (\ref{it.axiom6}) below ensures the connectedness of particle lines in $r(\Mb_{\ell'})$. Recall Definition \ref{def.connectedvia} and Definition \ref{def.layer_interval}.
\begin{enumerate}[resume]
\item\label{it.axiom6} For each $1\leq\ell'\leq \ell$, and any particle line $\pb\in r(\Mb_{\ell'})$, this $\pb$ is either connected to another particle line $\pb'\in r(\Mb_{\ell'})$ via $\Mb_{[1:\ell']}$, or connected to a cycle within $\Mb_{[1:\ell']}$, or forms an initial link within $\Mb_{[1:\ell']}$. See {\color{blue} Figure \ref{fig.moleculeprop}}.
\end{enumerate}

\textbf{The $\Fc_{\vLambda_\ell}^{\mathrm{err}}$ case.} Finally, suppose $\Mb\in\Fc_{\vLambda_\ell}^{\mathrm{err}}$ and $r(\Mb)=H$. Then $\Mb$ satisfies all of the above properties (\ref{it.axiom1-})--(\ref{it.axiom6}), except that (\ref{it.axiom6}) is true only for $\ell'\leq\ell-1$.
\end{proposition}
\begin{proof} We first prove (\ref{it.axiom1-})--(\ref{it.axiom7}), then prove the harder (\ref{it.axiom6}). The proof in the $\Mb\in\Fc_{\vLambda_\ell}^{\mathrm{err}}$ case is the same (since Definition \ref{def.set_T_F} (\ref{it.set_F_l_1})--(\ref{it.set_F_l_4}) are the same except that (\ref{it.set_F_l_4}) now requires $\ell'\leq\ell-1$).

\textbf{Proof of (\ref{it.axiom1-})--(\ref{it.axiom7}).} For (\ref{it.axiom1-}): this is already proved in Part 1 of the proof of Proposition \ref{prop.layerrec3}.

For (\ref{it.axiom1}): this follows from the correspondence between molecules and physical trajectories in Definition \ref{def.top_reduction}, and the fact that two (physical) particles in $\Rb^d$ cannot collide/overlap and immediately collide/overlap again without other any other collision in between.

For (\ref{it.axiom2}): recall Definition \ref{def.cluster}. We know $\Mb_{\ell'}$ can be constructed, starting from the C-molecule with only the clusters, by adding the O-atoms corresponding to the overlaps joining different clusters, according to cluster graph of $\Mb_{\ell'}$. Since this cluster graph is a forest (Definition \ref{def.set_T_F} (\ref{it.set_F_l_3})), we know that the process of adding O-atoms will \emph{not} generate new cycles in $\Mb_{\ell'}$, and will only insert O-atoms into edges of existing cycles (turning these edges into ov-segments). Therefore the first statement in (\ref{it.axiom2}) is true. The second statement follows from the first, as the O-atom must be an interior atom of an ov-segment connecting two C-atoms.

For (\ref{it.axiom3}): this is just Definition \ref{def.set_T_F} (\ref{it.set_F_l_1}).

For (\ref{it.axiom4}): this follows directly from Definition \ref{def.prsets}.

For (\ref{it.axiom4+}): this follows directly from Definition \ref{def.prsets} (see also Definition \ref{def.root_particle_line}) that $r(\Mb_\ell)=r(\Mb)=H$, and $r(\Mb_{\ell'})=p(\Mb_{\ell'})\cap p(\Mb_{\ell'+1})$.

For (\ref{it.axiom5}): this is just Definition \ref{def.set_T_F} (\ref{it.set_F_l_2}).

For (\ref{it.axiom7}): this directly follows from Proposition \ref{prop.init_link}.

\textbf{Proof of (\ref{it.axiom6}).} Fix any $\pb\in r(\Mb_{\ell'})$. If $\pb$ is connected to another particle line in $r(\Mb_{\ell'})$ via $\Mb_{\ell'}$ then we are done. If not, by Definition \ref{def.set_T_F} (\ref{it.set_F_l_4}), there exists $\pb_1\in r(\Mb_{\ell'-1})$ such that $\pb$ is either equal to or connected to $\pb_1$ via $\Mb_{\ell'}$. We then replace $\pb$ by $\pb_1$ and $\ell'$ by $\ell'-1$ and repeat the above discussion: either $\pb_1$ is connected to $\pb_1\neq \pb_1'\in r(\Mb_{\ell'-1})$ via $\Mb_{\ell'-1}$, or we can find $\pb_2$, and so on.

In the end, either we reach $\pb_{\ell'}\in r(\Mb_0)=H_0$, or we get $\pb_i\in r(\Mb_{\ell'-i})$ for some $0\leq i<\ell'$, which is connected to $\pb_{i}\neq \pb_i'\in r(\Mb_{\ell'-i})$ via $\Mb_{\ell'-i}$. In the former case, either $\pb=\pb_{\ell'}$ or $\pb$ is connected to $\pb_{\ell'}$ via $\Mb_{[1:\ell']}$; by Proposition \ref{prop.mol_axiom} (\ref{it.axiom7}), this $\pb_{\ell'}$ is involved in an initial link, so $\pb$ forms an initial link within $\Mb_{[1:\ell']}$ by definition.

Now assume we get $\pb_i$ with $0\leq i<\ell'$ as said above. By (\ref{it.axiom5}) we just proved (and note that $r(\Mb_{\ell'-i})\subseteq p(\Mb_{\ell'-i+1})$), there exists $\pb_{i-1}'\in r(\Mb_{\ell'-i+1})$ such that $\pb_i'$ either equals or is connected to $\pb_{i-1}'$ via $\Mb_{\ell'-i+1}$. We then replace $\pb_i'$ by $\pb_{i-1}'$ and $\ell'-i$ by $\ell'-i+1$ and repeat this process, and subsequently find $\pb_{i-2}',\cdots,\pb_1',\pb'$, see {\color{blue} Figure \ref{fig.connectedbelow}}. Note that $\pb'\in r(\Mb_{\ell'})$.

Now, if $\pb\neq \pb'$, then we can transit from $\pb$ to $\pb_1$ and then $\pb_2,\cdots,\pb_i$ in this order, and then transit to $\pb_i$, and then transit to $\pb_{i-1}',\cdots,\pb_1',\pb'$ in this order. In each transit step, the two particle lines involved are either equal or connected via bonds within a subset of $\Mb_{[1:\ell']}$; putting together, we get that $\pb$ and $\pb\neq\pb'\in r(\Mb_{\ell'})$ are connected via $\Mb_{[1:\ell']}$.

If instead $\pb=\pb'$, note that $\pb_i\neq \pb_i'$, so we choose the smallest $1\leq n\leq i$ such that $\pb_n\neq \pb_n'$. In this case, we know that $\pb_n$ and $\pb_n'$ are both connected to $\pb_{n-1}$ via $\Mb_{\ell'-n+1}$, so $\pb_n$ is also connected to $\pb_n'$ via $\Mb_{\ell'-n+1}$ (see Remark \ref{rem.layer_interval} below). Moreover $\pb_n$ is also connected to $\pb_n'$ via $\Mb_{[1:\ell'-n]}$, by transiting from $\pb_n$ to $(\pb_{n+1},\cdots,\pb_i)$ in this order and then to $(\pb_i',\pb_{i-1}',\cdots,\pb_n')$ in this order similar to the above paragraph. Then, by putting together a chain of bonds in $\Mb_{\ell'-n+1}$ connecting $\pb_n$ to $\pb_n'$, another chain of bonds in $\Mb_{[1:\ell'-n]}$ connecting $\pb_n$ to $\pb_n'$, and possibly a chain of bonds within the particle line $\pb_n$ and another chain of bonds within the particle line $\pb_n'$, we obtain a cycle in $\Mb_{[1:\ell'-n+1]}\subseteq \Mb_{[1:\ell']}$. Since $\pb$ is connected to $\pb_n$ (and thus to this cycle) via $\Mb_{[1:\ell']}$ by transiting from $\pb$ to $\pb_1,\cdots \pb_{n}$ in this order similar to the above paragraph, we conclude that $\pb$ is connected to cycle within $\Mb_{[1:\ell']}$. This completes the proof.
\end{proof}
To end this section, we record some intuitive facts about molecules, particle lines and cuttings, which will be used frequently in the algorithm below.
\begin{remark}\label{rem.layer_interval} Let $A$ be a layer interval $\Mb_{[\ell_1:\ell_2]}$ (Definition \ref{def.layer_interval}), or any union of connected components of a layer interval, or the same for thin layers (defined in Definition \ref{def.layer_refine} below). It is easy to see that: for any atoms $\nf,\nf'\in A$ that belong to the same particle line, and any atom $\mf$ in the same particle line between $\nf$ and $\nf'$, we must have $\mf\in A$. This implies that (i) any particle line $\pb$ can intersect at most one component of $A$, and (ii) if two particle lines $\pb_1,\pb_2$ are both connected to a particle line $\pb$ via $A$ (Definition \ref{def.connectedvia}), then $\pb_1$ and $\pb_2$ are also connected via $A$. This is used in the proof of Proposition \ref{prop.mol_axiom} (\ref{it.axiom6}) above; see {\color{blue} Figure \ref{fig.connectedbelow}}.
\end{remark}
\begin{remark}\label{rem.full_cut} The following fact will be convenient to note: suppose we cut $A\subseteq\Mb$ \emph{as free}, then this may divide the remaining set $\Mb\backslash A$ into connected components; however, any \emph{full} component of $\Mb\backslash A$ must also be a full component of $\Mb$ before cutting $A$ as free (in other words, any \emph{newly created} connected component must have at least one fixed end). This is because, when any bond or ov-segment is broken in cutting $A$ as free, we get a fixed end at \emph{every atom} in this ov-segment that belongs to $\Mb\backslash A$, see Remark \ref{rem.cutting}. 
\end{remark}
\begin{remark}\label{rem.aftercut} In the algorithm below, we usually adopt the following convention: once we cut out any elementary molecule $S$, or cut out any subset $S$ (as free or as fixed) and subsequently cut it into elementary molecules, we will ignore this $S$ for the rest of the cutting sequence and consider the remaining molecule $\Mb\backslash S$ only. For example, when we say ``keep cutting deg 2 atoms as free until there is no more deg 2 atom", we always mean there is no more deg 2 atom in the remaining molecule. 
\end{remark}

\section{Main thread of the algorithm}\label{sec.toy}

We now start the proof of Proposition \ref{prop.comb_est}. Most of the proof below will focus on the case $\Mb\in\Fc_{\vLambda_\ell}$; the case $\Mb\in\Fc_{\vLambda_\ell}^{\mathrm{err}}$ is left to the end of Section \ref{sec.finish}, and will only take two paragraphs. In this section, we first discuss a simplified setting which contains most of the main ideas.

\subsection{Toy models}\label{sec.toy_intro}

To prove Proposition \ref{prop.comb_est}, we need to construct an operation sequence (with cutting, deleting and splitting). This will be done using a number of \textbf{cutting algorithms}, which will be presented in Sections \ref{sec.layer}--\ref{sec.maincr}. 

The full cutting algorithms are quite complicated; however, there is a much simpler \textbf{main thread}, and many of the complications actually come from the additional structures and conditions (``accessories") that have to be added to this main thread. In this section, we present the proof in a simplified scenario, which basically isolates this main thread. The accessories needed in full generality will be discussed in Section \ref{sec.toy_reduce}.

Throughout this section, we will make the following two simplifications (and more later):
\begin{enumerate}
\item\label{it.simp_1} $\Mb$ has no O-atoms, all \{33\} molecules are good, and $H_0=\varnothing$ (see Proposition \ref{prop.mol_axiom} (\ref{it.axiom4+}));
\item\label{it.simp_2} Each layer $\Mb_{\ell'}$ of $\Mb$ is a forest.
\end{enumerate} Here, regarding (\ref{it.simp_1}), recall that \{33\} molecules are good provided certain \emph{non-degeneracy} conditions are met (which is different for \{33A\} and \{33B\} molecules, see (\ref{eq.good_normal_2})--(\ref{eq.good_normal_4}) in Definition \ref{def.good_normal}). As such, simplification (\ref{it.simp_1}) allows us to ignore any O-atoms, degeneracies and initial (time $t=0$) cumulants, and (\ref{it.simp_2}) allows us to ignore any cycles (or ``recollisions") happening within each layer.

Given a molecule $\Mb$ satisfying Proposition \ref{prop.mol_axiom} and (\ref{it.simp_1})--(\ref{it.simp_2}), the goal of this section is to construct a cutting sequence that cuts $\Mb$ into elementary molecules, such that (recall $\upsilon=3^{-d-1}$)
\begin{equation}\label{eq.goal}
(\upsilon/2)\cdot\#_{\{33\}}-d\cdot\#_{\{4\}}\geq c\cdot\rho,
\end{equation} where $\#_{\{33\}}$ is the number of \{33\} molecules etc., $c$ is a positive constant and $\rho$ is in Definition \ref{def.parameter_rho_old} (the splitting and deletion of O-atoms is not needed in this simplified model; also note that there cannot be any empty end if $H_0=\varnothing$, which easily follows from Proposition \ref{prop.mol_axiom} (\ref{it.axiom6})).

In Sections \ref{sec.toy1}--\ref{sec.toy3}, we will consider a further simplified case, when $\Mb$ contains only \emph{two layers} (for treatment of multi-layer case, see Section \ref{sec.toy_multi}). We make the following definition:
\begin{definition}\label{def.simplified} We say $\Mb$ is a \emph{2-layer model}, if it satisfies (\ref{it.simp_1})--(\ref{it.simp_2}), and that
\begin{enumerate}[resume]
\item\label{it.simp_3} $\Mb$ can be written as $\Mb=\Mb_U\cup\Mb_D$, such that no atom in $\Mb_D$ is parent of an atom in $\Mb_U$;
\item\label{it.simp_4} Each of $\Mb_U$ and $\Mb_D$ is a tree, has no fixed end, and the number of bonds connecting $\Mb_U$ and $\Mb_D$ is at least $\rho$.
\end{enumerate}
Note that in practice, the $\rho$ in (\ref{it.simp_4}) may be replaced by some different quantity, but they will always be comparable, so this will not affect the proof.
\end{definition}
\begin{figure}[h!]
    \centering
    \includegraphics[width=1\linewidth]{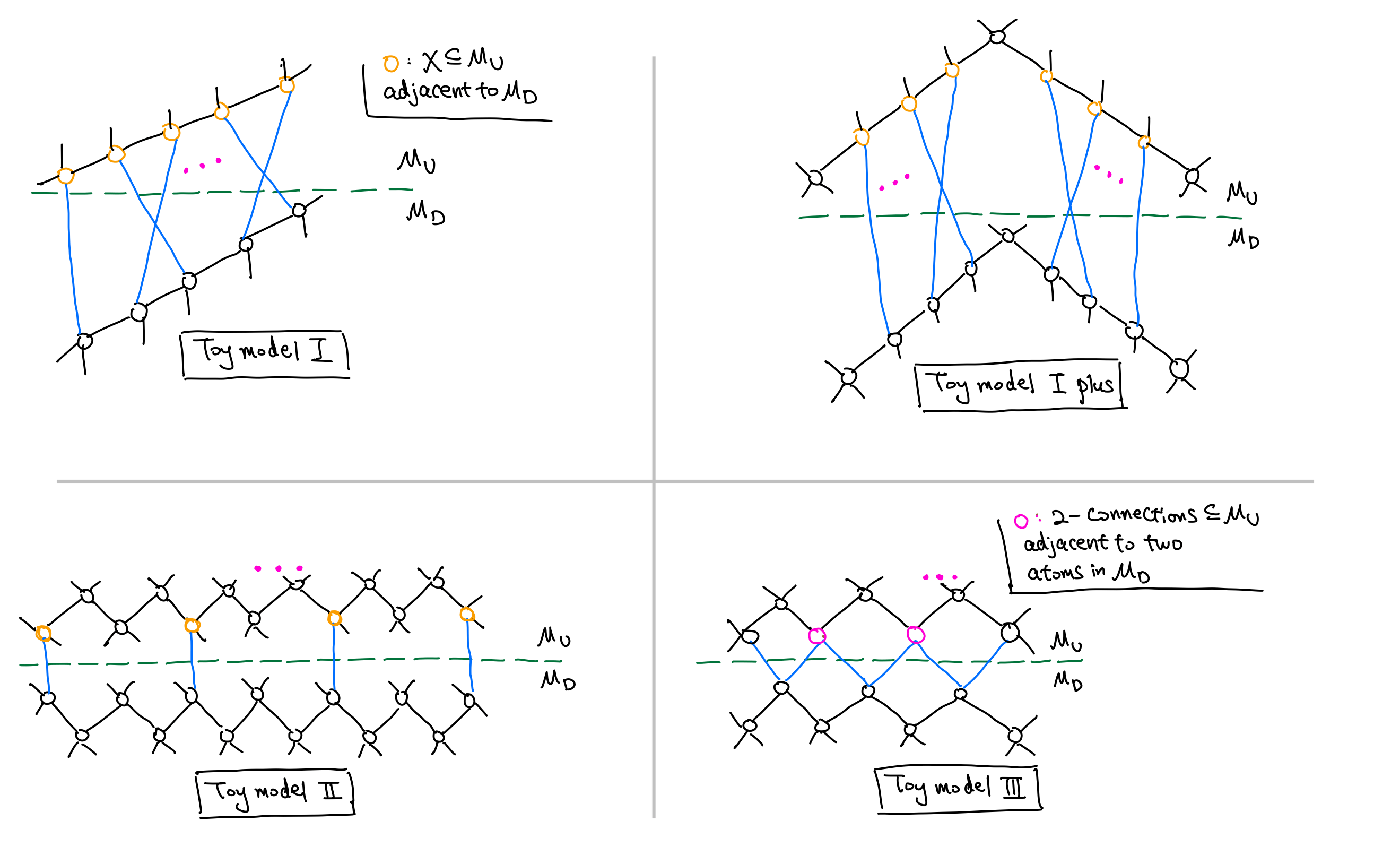}
    \caption{Examples of toy models. Upper left: toy model I ($X=\Mb_U$). Upper right: toy model I plus (where $\#_{\mathrm{comp}(X)}=2$). Lower left: toy model II (where $\#_{\mathrm{comp}(X)}=|X|$). Lower right: toy model III (where $\#_{\mathrm{2conn}}\sim |X|/2$). The blue bonds connect $\Mb_U$ and $\Mb_D$. The orange and pink atoms indicate those in $X$ and 2-connections respectively. The pink dots represent patterns that can be repeated.}
    \label{fig.toy models}
\end{figure}
In our proof, the 2-layer model in Definition \ref{def.simplified} will be further divided into 3 cases, corresponding to \textbf{toy models I--III} defined below. This classification is based on the following two considerations:
\begin{enumerate}[{(1)}]
\item Consider the set of atoms in $\Mb_U$ that are adjacent to \emph{two} atoms in $\Mb_D$ by two bonds. Is this set big or small relative to $\rho$?
\item Consider all the atoms in $\Mb_U$ that are adjacent to \emph{at least one} atom in $\Mb_D$ by at least one bond. Are they \emph{all connected} to each other, or are they \emph{separated} by other atoms in $\Mb_U$?
\end{enumerate}

The reason for considering these two questions will be clear later. For convenience, in this section we will refer to the atoms in $\Mb_U$ that are adjacent to \emph{two} atoms in $\Mb_D$ as \emph{2-connections}, and denote their number by $\#_{\mathrm{2conn}}$. We also define the set of atoms in $\Mb_U$ that are adjacent to \emph{at least one} atom in $\Mb_D$ by $X$ (note that $|X|\gtrsim\rho$), and define the number of connected components of $X$ (with $X$ viewed by itself, without passing through atoms in $\Mb_U\backslash X$) by $\#_{\mathrm{comp}(X)}$. We now define a 2-layer molecule to be:
\begin{itemize}
\item \emph{Toy model I}: if $\#_{\mathrm{2conn}}=0$ and $X=\Mb_U$. This case involves the core arguments in our algorithm in its simplest form, and will be proved in Section \ref{sec.toy1}.
\item \emph{Toy model I plus}: if $\#_{\mathrm{2conn}}=0$ and $\#_{\mathrm{comp}(X)}\ll |X|$ (no 2-connections, $X$ has few components). This is an extension of the toy model I, and can be treated in similar manners, in Section \ref{sec.toy1+}. Here we need to introduce the new idea of \textbf{UP} algorithm, which will play an important role in the later proofs. 
\item \emph{Toy model II}: if $\#_{\mathrm{2conn}}=0$ and $\#_{\mathrm{comp}(X)}\gtrsim |X|$ (no 2-connections, $X$ has many components). This requires a different argument, which involves making a clever choice in the \textbf{UP} algorithm introduced above, see Section \ref{sec.toy2}.
\item \emph{Toy model III}: if $\#_{\mathrm{2conn}}\gtrsim |X|$ (many 2-connections). This relies on a simpler argument, which again involves making suitable choices in the \textbf{UP} algorithm, see Section \ref{sec.toy3}.
\end{itemize} 
See {\color{blue} Figure \ref{fig.toy models}} for examples of toy models I--III (we distinguish toy models I and I plus only for convenience of exposition). Also the case when $0<\#_{\mathrm{2conn}}\ll|X|$, which is missing in the toy models, can be treated similarly to toy models I--II but with an extra step of removing 2-connections, see Section \ref{sec.reduce5}.

Once all the toy models (and hence the 2-layer model) have been addressed in Sections \ref{sec.toy1}--\ref{sec.toy3}, we can then treat the multi-layer case in Section \ref{sec.toy_multi}. This relies on the process of \textbf{layer selection}, which effectively reduces any multi-layer $\Mb$ to (basically) a 2-layer model, for which the proofs in Sections \ref{sec.toy1}--\ref{sec.toy3} apply.

Finally, in Section \ref{sec.toy_reduce} we discuss the simplifications made in (\ref{it.simp_1})--(\ref{it.simp_4}) above compared to the full generality case, as well as the corresponding accessories (additions, extra conditions and proofs) that are needed to treat the full generality setting.

Before proceeding, we make one remark about the presentation in this section and in Sections \ref{sec.layer}--\ref{sec.maincr}: in this section, the cases are arranged for \emph{exposition purposes}, for example we present the 2-layer model before layer selection, and present toy models I and I plus before the other toy models, in order to showcase the core ideas of our proof as early as possible. However, in Sections \ref{sec.layer}--\ref{sec.maincr} we need to present the materials in \emph{logical order}, so we need to introduce layer selection (Section \ref{sec.layer_select}) before the toy models (Section \ref{sec.maincr}), and toy model III (Definition \ref{def.alg_2connup}--Proposition \ref{prop.alg_2connup}) before toy model I plus (Definition \ref{def.alg_maincr}--Proposition \ref{prop.alg_maincr}), etc.

To make the proof easier to read, each time we discuss any argument in Sections \ref{sec.layer}--\ref{sec.maincr} that has a counterpart in this section, we will explicitly point out which argument in this section does it correspond to. See {\color{blue}Table \ref{tab.toy}} for a summary of the toy models and layer selection, and where each of them occurs in Sections \ref{sec.layer}--\ref{sec.maincr}. Also, in {\color{blue}Table \ref{tab.reduce}} we list the simplifications made in Definition \ref{def.simplified}, the corresponding accessories needed in the general setting, and where each of them occurs in Sections \ref{sec.layer}--\ref{sec.maincr}.
\begin{table}[H]
\centering
\begin{tabular}{|c|c|c|c|c|}
\hline
 &Case & Property & Location &Algorithm name\\
 \hline
 (1)& Toy model I & \makecell{No 2-connections \\and $X=\Mb_U$}& (Special case) &(Special case)\\
 \hline
 (2) & Toy model I plus & \makecell{No 2-connections \\$X$ has few components} & \makecell{Def. \ref{def.alg_maincr}, Prop. \ref{prop.alg_maincr}\\\textbf{Stage 6}} &\textbf{MAINUD}\\
 \hline
(3) & Toy model II & \makecell{No 2-connections \\$X$ has many components} & \makecell{Def. \ref{def.3comp_alg}, Prop. \ref{prop.3comp}\\\textbf{Stage 5}}&\textbf{3COMPUP}\\
\hline
(4) & Toy model III & Many 2-connections & \makecell{Def. \ref{def.alg_2connup}, Prop. \ref{prop.alg_2connup}\\ Prop. \ref{prop.comb_est_case6}}&\textbf{2CONNUP/DN}\\
\hline
(5) & Layer selection & Multi-layer & Def. \ref{def.layer_select}--Prop. \ref{prop.layer_cutting}&(Layer selection)\\
\hline
\end{tabular}
\caption{The different toy models and their definitions, the places where they occur and the corresponding names of algorithms.}
\label{tab.toy}
\end{table}
\begin{table}[H]
\centering
\begin{tabular}{|c|c|c|c|}
\hline
 &Simplification & Accessories & Location\\
 \hline
 (1)& No O-atoms & General cutting, ov-segments etc.& (Various)\\
 \hline
 (2) & All \{33\} are good & Strong and weak degeneracies & \makecell{Def. \ref{def.strdeg}, Prop. \ref{prop.case2}\\ Def. \ref{def.weadeg}, Prop. \ref{prop.comb_est_case4}}\\
 \hline
  (3) & $H_0=\varnothing$ & Treating initial cumulants & Prop. \ref{prop.comb_est_case4}\\
 \hline
(4) & $\Mb_{\ell'}$ is a forest & Layer refining and large $\Rf$ case & \makecell{Def. \ref{def.layer_refine}, Prop. \ref{prop.layer_refine_2}\\Prop. \ref{prop.alg_up_recl}--\ref{prop.comb_est_case3}}\\
\hline
(5) & \makecell{Additional \\(see Section \ref{sec.reduce5})} & Pre-processing and others & \makecell{Def. \ref{def.func_select}, Prop. \ref{prop.func_select}\\\textbf{Stage 1--4}}\\
\hline
\end{tabular}
\caption{Simplifications in Definition \ref{def.simplified} and corresponding accessories which will be needed in real problem. See Section \ref{sec.toy_reduce} for more details.}
\label{tab.reduce}
\end{table}
\subsection{Toy model I}\label{sec.toy1} 
\begin{figure}[h!]
    \centering
    \includegraphics[width=0.55\linewidth]{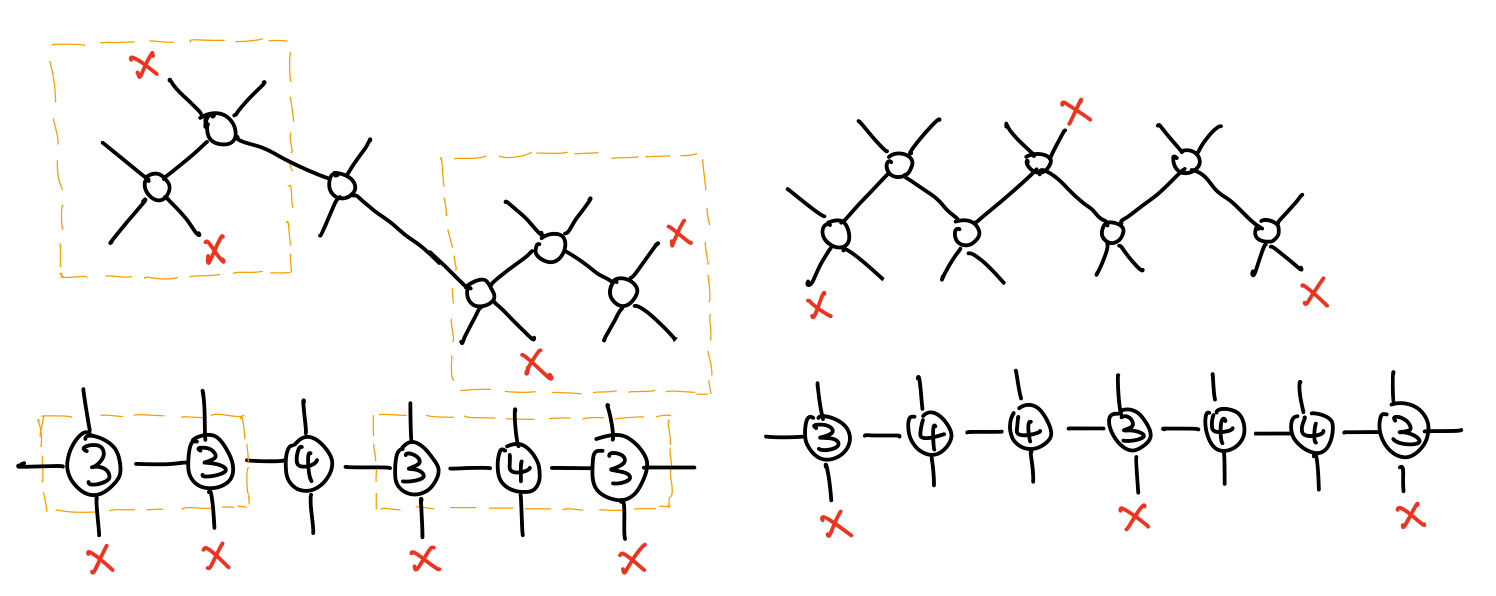}
    \caption{Examples of \{33\} and \{343\} molecules (left), and proper molecules (right, see Definition \ref{def.toy_proper}). Note these notions are applied to $\Mb_D$, in which the top/bottom distinction is less important, hence we draw an equivalent simplified picture below, with no distinction of top/bottom, and with the deg (i.e. degree) of each atom indicated. Similar simplified illustrations will occur a few times below.}
    \label{fig.proper}
    \end{figure}

\begin{figure}[h!]
    \centering
    \includegraphics[width=0.9\linewidth]{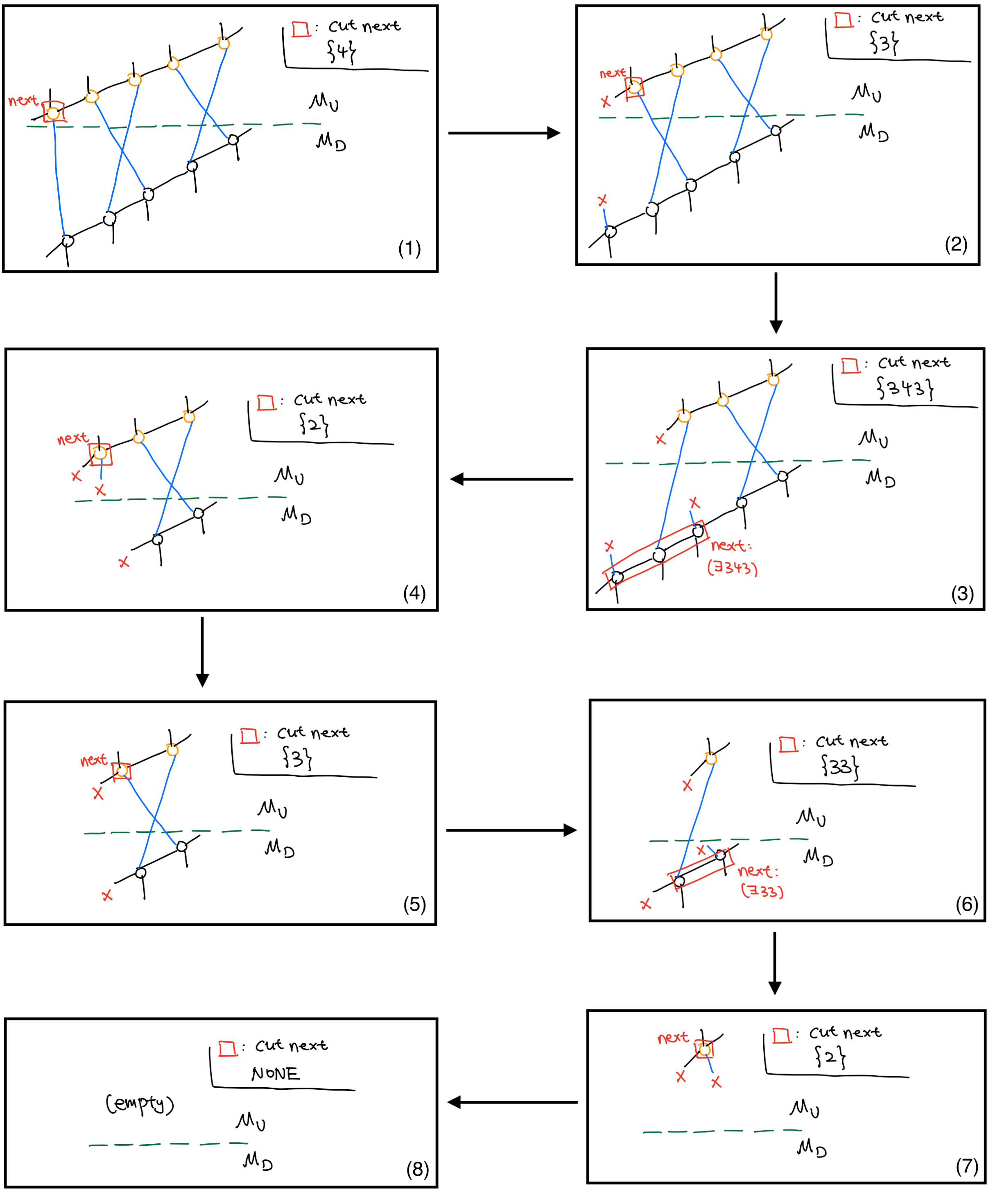}
    \caption{Algorithm for toy model I (Definition \ref{def.toy_alg}), applied to the example in {\color{blue} Figure \ref{fig.toy models}}. Each time we cut the atom(s) in the red box, creating fixed ends. The type of the elementary molecules to be cut is listed at the upper right corner.}
    \label{fig.toy model 1 alg}
\end{figure}

\begin{figure}[h!]
    \centering
    \includegraphics[width=0.3\linewidth]{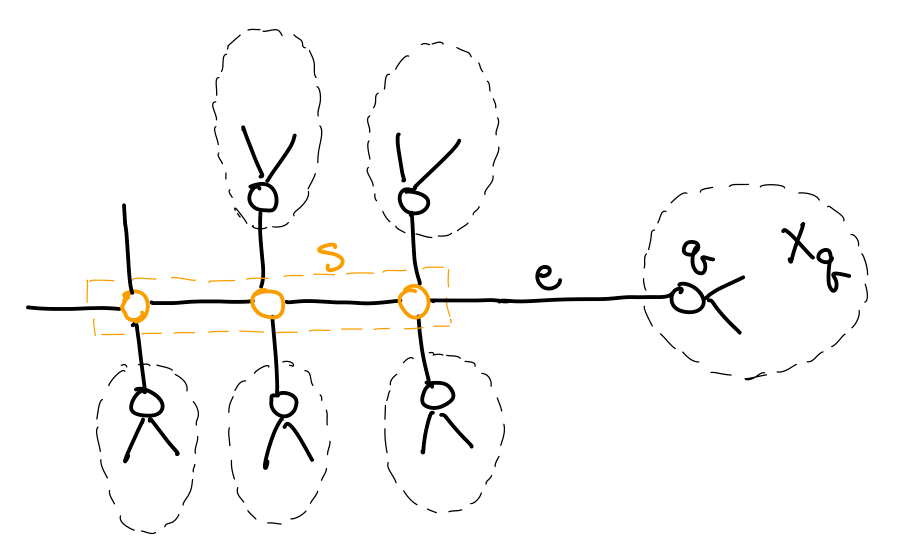}
    \caption{An illustration of Lemma \ref{lem.cutconnected0}. Here $S$ is a connected subset, and there are bijections between $\qf$ and $X_\qf$, where $\qf\in \Mb_D\backslash S$ is connected to $S$ by a bond $e$, and $X_\qf$ is the component containing $\qf$ after cutting $S$ as free. When viewed in $X_q$, the operation of cutting $S$ as free creates exactly one fixed end in $X_\qf$. All these crucially depend on the fact that $\Mb_U$ and $\Mb_D$ have no cycles.}
    \label{fig.heart}
\end{figure}

In this subsection we consider the following (first) toy model:
\begin{definition}\label{def.toy} We say $\Mb$ is a \emph{toy model I}, if it is a 2-layer model as in Definition \ref{def.simplified}, and that each atom in $\Mb_U$ has \emph{exactly} one bond connecting to an atom in $\Mb_D$ (i.e. $\#_{\mathrm{2conn}}=0$ and $X=\Mb_U$).
\end{definition}
\begin{definition}
\label{def.toy_proper} Let $\Mb_D$ be part of the 2-layer model as in Definition \ref{def.simplified} (or any resulting molecule after an operation sequence). We define $\Mb_D$ to be \emph{proper}, if it does not contain deg 1 or 2 atoms, does not contain two adjacent deg 3 atoms, and does not contain two deg 3 atoms which are both adjacent to the same deg 4 atom, see {\color{blue} Figure \ref{fig.proper}}.
\end{definition}
If $\Mb$ is a toy model I, we define the following algorithm applied to it:
\begin{definition}[algorithm for toy model I]\label{def.toy_alg} Let $\Mb$ be a toy model I. We define the following cutting sequence. For an example see {\color{blue} Figure \ref{fig.toy model 1 alg}}:
\begin{enumerate}[{(1)}]
\item\label{it.toy_alg_1} If not all atom in $\Mb_U$ have been cut, then choose a lowest atom $\nf\in\Mb_U$ that has not been cut (i.e. either $\nf$ has not child in $\Mb_U$ or all its children in $\Mb_U$ have been cut).
\item\label{it.toy_alg_2} If $\nf$ is adjacent to an atom $\mf\in\Mb_D$ which has deg 3, then cut $\{\nf,\mf\}$ as free and cut $\mf$ as free from $\{\nf,\mf\}$ if $\nf$ has deg 4; otherwise just cut $\nf$ as free.
\item\label{it.toy_alg_3} If $\Mb_D$ contains two adjacent deg 3 atoms (say $\rf$ and $\rf'$), then cut $\{\rf,\rf'\}$ as free; repeat until $\Mb_D$ has no adjacent deg 3 atoms.
\item\label{it.toy_alg_4} If $\Mb_D$ contains two degree 3 atoms (say $\rf$ and $\rf''$) both adjacent to the same deg 4 atom (say $\rf'$), then cut $\{\rf,\rf',\rf''\}$ as free and cut $\rf$ as free from $\{\rf,\rf',\rf''\}$; then go to \ref{it.toy_alg_3}--\ref{it.toy_alg_4} and repeat, until $\Mb_D$ becomes proper (Definition \ref{def.toy_proper}).
\item \label{it.toy_alg_5} Go to \ref{it.toy_alg_1} and choose the next lowest atom $\nf\in\Mb_U$ that has not been cut, and so on.
\item\label{it.toy_alg_6} Suppose all atoms in $\Mb_U$ have been cut. If $\Mb_D$ is proper, then cut any deg 3 atom; if not, then go to \ref{it.toy_alg_3}--\ref{it.toy_alg_4} and repeat to remove them, and next cut any deg 3 atom, and so on.
\end{enumerate}
\end{definition}
\begin{remark}\label{rem.toy_1}In Definition \ref{def.toy_alg}, we are essentially applying a ``naive" algorithm in $\Mb_U$ where we simply cut a lowest atom each time, and a ``greedy" algorithm in $\Mb_D$ where we exploit a \{33\} molecule (or what we call a \{343\} molecule\footnote{The reason to include this \{343\} molecule is technical: removing only \{33\} cannot exclude all forbidden scenarios, i.e. Lemma \ref{lem.toy_1} would not be true if we did this.} where two deg 3 atoms are adjacent to the same deg 4 atom). As mentioned above, this is the core idea of our proof, in its bare simplest form. It is a simple special case of the \textbf{MAINUD} algorithm, which corresponds to the toy model I plus in Section \ref{sec.toy1+}.
\end{remark}
Before stating the main result for the algorithm in Definition \ref{def.toy}, we first prove a useful lemma that will be applied frequently below, see {\color{blue} Figure \ref{fig.heart}}.
\begin{lemma}\label{lem.cutconnected0} Suppose $\Mb$ is a forest and we cut a \emph{connected subset} $S\subseteq\Mb$ as free (for example, when $S$ contains a single atom, or 2 atoms connected by a bond, or 3 atoms connected by 2 bonds, etc.). Then, for each atom $\qf\in\Mb\backslash S$ that is connected to an atom in $S$ by a (pre-cutting) bond $e$, the post-cutting molecule $\Mb_D\backslash S$ has exactly one component $X_\qf\ni\qf$, which is in bijection with $\qf$ (and $\qf$ is the unique atom in $X_\qf$ that is adjacent to $S$ pre-cutting, via the bond $e$). When viewed in $X_q$, the bond $e$ is viewed as a free end at $\qf$ before cutting, and is turned into a fixed end after cutting.
\end{lemma}
\begin{proof} This follows from the property of forests. In fact, any $\qf_1$ and $\qf_2$ (with corresponding bonds $e_1$ and $e_2$) must belong to different components of $\Mb_D$ after cutting $S$; if not, then (pre-cutting) they can be connected by a path involving atoms in $\Mb_D\backslash S$, but they can also be connected by a path with all intermediate atoms being in $S$ (using the bonds $e_1$ and $e_2$ and the connectivity of $S$ by itself), which produces a cycle and contradicts the forest assumption. Thus, for each $\qf$ we get a unique post-cutting component $X_\qf$, such that $e$ is the only pre-cutting bond between $X_\qf$ and $S$. Clearly, when viewed in $X_\qf$, the effect of cutting $S$ as free is equivalent to turning one and only one free end at $\qf$ (the one corresponding to $e$) into a fixed end.
\end{proof}
\begin{proposition}\label{prop.toy} Let $\Mb$ be a toy model I, and we apply the algorithm defined in Definition \ref{def.toy_alg}. Then in this process we only obtain elementary molecules, and we have
\begin{equation}\label{eq.toy_1}\#_{\{4\}}=1,\qquad \#_{\{33\}}\geq (|\Mb_U|-1)/5.\end{equation} Since $|\Mb_U|=|X|\gtrsim\rho\gg 1$, this clearly implies (\ref{eq.goal}).
\end{proposition}
\begin{proof} Recall Definition \ref{def.toy_alg} \ref{it.toy_alg_1}--\ref{it.toy_alg_6}. We divide the proof into 4 parts.

\textbf{Proof part 1.} We first prove that each atom in $\Mb_U$ is contained in an elementary molecule (after cutting). By the choices we make in \ref{it.toy_alg_1} and \ref{it.toy_alg_5} in Definition \ref{def.toy_alg}, we know that for each $\nf\in\Mb_U$, at the time it is cut, it must be a \emph{lowest} atom that still remains in $\Mb_U$; thus no \emph{child} of it still remains in $\Mb_U$, so cutting $\nf$ will not introduce any \emph{top fixed end} in $\Mb_U$ (and neither does cutting anything in $\Mb_D$). This means that $\Mb_U$ will not contain any top fixed end throughout the process. Since each $\nf\in\Mb_U$ is either cut by itself (when we cut $\{\nf\}$ as free in Definition \ref{def.toy_alg} \ref{it.toy_alg_2}) or belongs to a \{33\} molecule (when we cut $\{\nf,\mf\}$ as free in Definition \ref{def.toy_alg} \ref{it.toy_alg_2}), it is clear by definition of elementary molecules that $\nf$ always belongs to an elementary molecule.

\textbf{Proof part 2.} We next prove that each atom in $\Mb_D$ is contained in an elementary molecule. Note that $\Mb_D$ initially is a tree, and will remain a forest in the whole cutting process (cutting an atom will not create any cycle, but may cause a tree to become a disconnected forest). Now consider each step in Definition \ref{def.toy_alg} \ref{it.toy_alg_2} involving $\nf\in\Mb_U$, at which time $\Mb_D$ must be proper by the construction of our algorithm. Consider the possible effect of this step on $\Mb_D$ (we may assume that $\nf$ is connected to one atom $\mf\in\Mb_D$ by a bond, otherwise cutting $\nf$ will have no effect on $\Mb_D$), see {\color{blue} Figure \ref{fig.operab}}:
\begin{enumerate}[{(a)}]
\item\label{it.stepa} If we cut $\nf$ as free (which means $\mf$ has deg 4), then when viewed in $\Mb_D$, this step simply turns one and only one free end at $\mf\in\Mb_D$ into a fixed end, turning $\mf$ from deg 4 to deg 3.
\item\label{it.stepb} If we cut $\{\nf,\mf\}$ as free (which means $\mf$ has deg 3), then when viewed in $\Mb_D$, this step simply cuts $\mf$, which has deg 3, as free from $\Mb_D$.
\end{enumerate} Here in \ref{it.stepa}, we note that $\mf$ must have deg 4 before cutting $\nf$ (otherwise $\mf$ must have deg 3 because $\Mb_D$ is proper, in which case we will cut $\{\nf,\mf\}$ by Definition \ref{def.toy_alg} \ref{it.toy_alg_2}); also when viewed in $\Mb_D$, the bond between $\nf$ and $\mf$ is initially viewed as a free end at $\mf$, and is then turned into a fixed end after cutting $\nf$ as free, turning $\mf$ from deg 4 to deg 3. Likewise $\mf$ must have deg 3 in \ref{it.stepb}, and when viewed in $\Mb_D$, cutting $\{\nf,\mf\}$ as free is equivalent to cutting $\{\mf\}$ as free because $\nf$ has no other bonds connecting to other atoms in $\Mb_D$.
\begin{figure}[h!]
    \centering
    \includegraphics[width=0.6\linewidth]{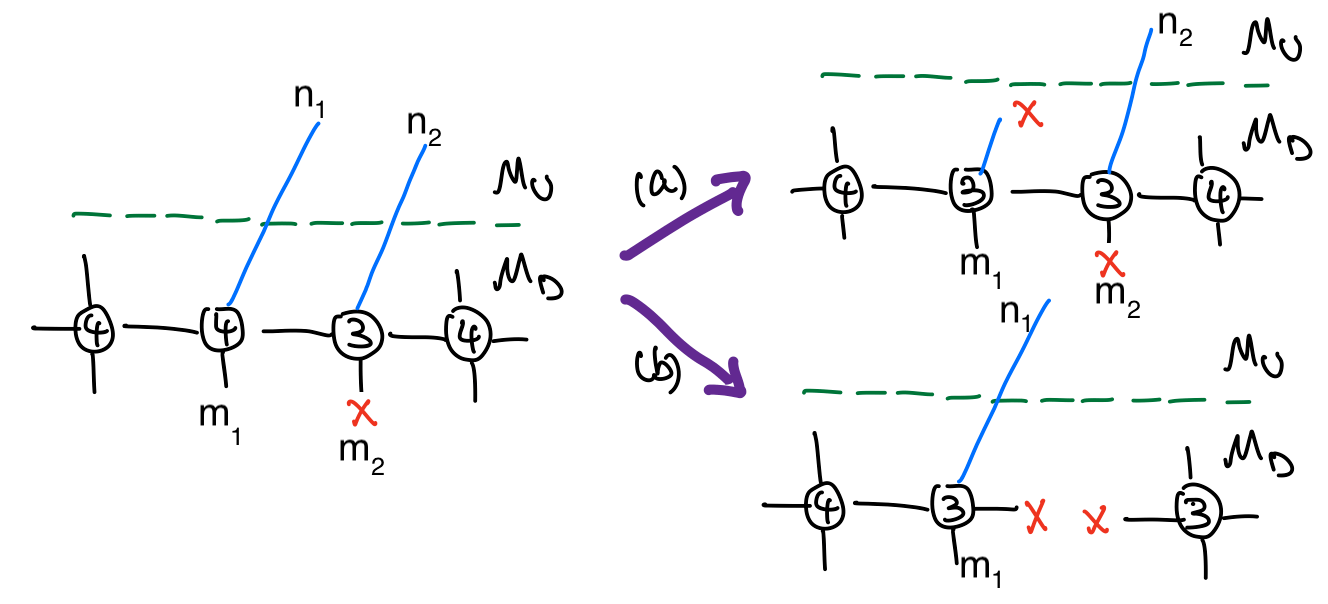}
    \caption{Operations (a) and (b) defined in the proof of Proposition \ref{prop.toy}. In either case an atom $\mf_j\in\Mb_D$ is connected to an atom $\nf_j\in\Mb_U$ by a blue bond. In case (a) $\mf_1$ has deg 4, and we cut $\{\nf_1\}$ as free, creating a fixed end at $\mf_1\in\Mb_D$ and turning it into deg 3; in case (b) $\mf_2$ has deg 3, and we cut $\{\nf_2,\mf_2\}$ as free, which when viewed in $\Mb_D$ just corresponds to cutting $\{\mf_2\}$ as free, creating fixed ends at its adjacent atoms.}
    \label{fig.operab}
\end{figure}

Now, regarding the effects described in \ref{it.stepa}--\ref{it.stepb}, we have the following lemma:
\begin{figure}[h!]
    \centering
    \includegraphics[width=0.5\linewidth]{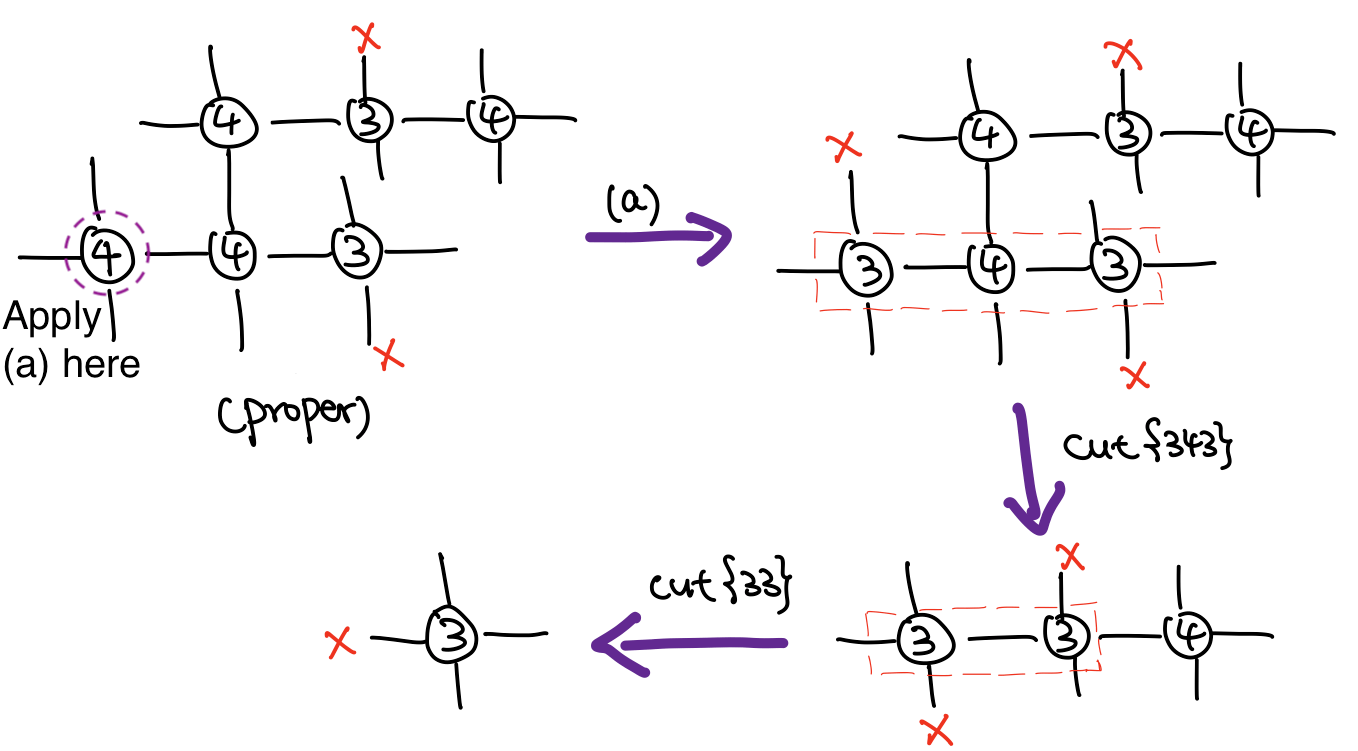}
    \caption{The greedy algorithm in Definition \ref{def.toy_alg}. Starting from upper left graph (which is proper), we apply operation (a), which creates a \{343\} molecule (in red box). We cut it, which then creates a \{33\} molecule. We cut it, and the graph becomes proper again.}
    \label{fig.greedy}
\end{figure}
\begin{lemma}\label{lem.toy_1} Suppose $\Mb_D$ is proper at some point of the algorithm, and becomes non-proper after either \ref{it.stepa} or \ref{it.stepb}. Then, $\Mb_D$ will once more become proper after applying Definition \ref{def.toy_alg} \ref{it.toy_alg_3}--\ref{it.toy_alg_4}, via cuttings of \{33\} and \{343\} molecules only (with each \{343\} molecule further cut into a \{3\} and a \{33\} molecule), and \emph{without seeing any deg 2 atom}, see {\color{blue} Figure \ref{fig.greedy}} for an illustration.
\end{lemma}
\begin{proof} Note that $\Mb_D$ is always a forest, so we can apply Lemma \ref{lem.cutconnected0}. Note that right before performing \ref{it.stepb} at a deg 3 atom $\mf$, any atom in $\Mb_D$ adjacent to $\mf$ must have deg 4 (because $\Mb_D$ is proper). Then by Lemma \ref{lem.cutconnected0}, $\Mb_D$ will break into  at most 2 components after performing \ref{it.stepb}, and the effect of \ref{it.stepb} \emph{when viewed in each component} is equivalent to \ref{it.stepa}. Therefore, we only need to prove Lemma \ref{lem.toy_1} for \ref{it.stepa}.

Assume now $\Mb_D$ becomes non-proper, say a \{33\} molecule forms, after \ref{it.stepa}. Let this molecule be $\{\rf,\rf'\}$, where $\rf'$ has deg 3 before \ref{it.stepa}, while $\rf$ has deg 4 before \ref{it.stepa} and is turned into deg 3 after \ref{it.stepa}. Then, \emph{before \ref{it.stepa}}, any atom in $\Mb_D$ adjacent to $\rf$ or $\rf'$ \emph{must have deg 4} (see {\color{blue} Figure \ref{fig.prooftoy1}}) because otherwise we either have two adjacent deg 3 atoms before \ref{it.stepa}, or have two deg 3 atoms adjacent to the same deg 4 atom before \ref{it.stepa}, contradicting properness. Therefore, by Lemma \ref{lem.cutconnected0}, we see that cutting the \{33\} molecule breaks $\Mb_D$ into at most 4 components, and the effect of this cutting \emph{when viewed in each component} is again equivalent to \ref{it.stepa}. We then induct on the size of $\Mb_D$ and apply induction hypothesis to each of these components, to conclude the proof.

Finally, assume $\Mb_D$ becomes non-proper, but with no \{33\} molecule forming, after \ref{it.stepa}. In this case there must be a \{343\} molecule forming, say $\{\rf,\rf',\rf''\}$, where $\rf'$ and $\rf''$ respectively have deg 4 and 3 before \ref{it.stepa}, and $\rf$ has deg 4 before \ref{it.stepa} and is turned into deg 3 after \ref{it.stepa}. Then, \emph{before \ref{it.stepa}}, any atom adjacent to $\rf'$ or $\rf''$ \emph{must have deg 4} (see {\color{blue} Figure \ref{fig.prooftoy1}}) for the same reason as above thanks to properness. Moreover, before \ref{it.stepa}, any atom adjacent to $\rf$ must also have deg 4, because otherwise we would get a \{33\} molecule after \ref{it.stepa}, contradicting our assumption. Therefore, by Lemma \ref{lem.cutconnected0}, we see that cutting the \{343\} molecule breaks $\Mb_D$ into at most 6 components, and the effect of this cutting \emph{when viewed in each component} is again equivalent to \ref{it.stepa}. We then conclude the proof by the same induction as above. This proves Lemma \ref{lem.toy_1}.
\end{proof}
\begin{figure}[h!]
    \centering
    \includegraphics[width=0.6\linewidth]{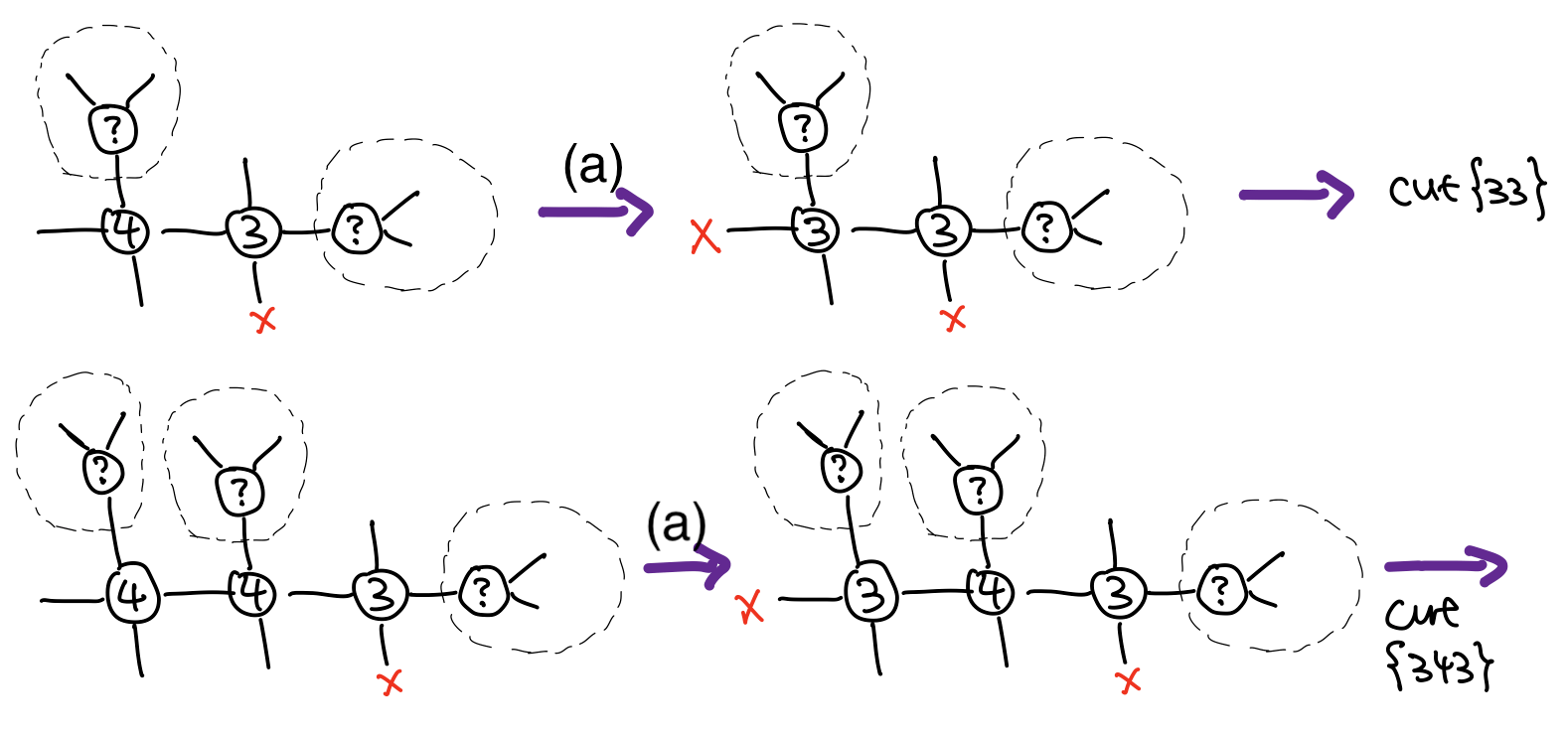}
    \caption{Proof of Lemma \ref{lem.toy_1}. After applying (a) to a proper $\Mb_D$, a deg 4 atom is turned into deg 3, thus creating a \{33\} (top) or \{343\} molecule (bottom). Crucially, all the atoms with ``?" must be deg 4, because $\Mb_D$ is proper before applying (a), and because applying (a) does not create \{33\} molecule in the bottom case. Then, after cutting the \{33\} or \{343\} molecule, when viewed in each resulting component (which contains one ``?" atom), we create only one fixed end and turn the ``?" atom from deg 4 to deg 3, and so on.}
    \label{fig.prooftoy1}
\end{figure}
Now, with Lemma \ref{lem.toy_1}, we can finish Part 2 of the proof, i.e. each atom in $\Mb_D$ is contained in an elementary molecule. In fact, those cut in Definition \ref{def.toy_alg} \ref{it.toy_alg_2} clearly belong to \{33\} molecule. Those cut in Definition \ref{def.toy_alg} \ref{it.toy_alg_3}--\ref{it.toy_alg_4} also belong to \{33\} or \{343\} (and subsequently \{3\} or \{33\}) molecules by Lemma \ref{lem.toy_1}. Finally, those cut in Definition \ref{def.toy_alg} \ref{it.toy_alg_6} either belong to \{3\} molecule, or are cut in the process of applying Definition \ref{def.toy_alg} \ref{it.toy_alg_3}--\ref{it.toy_alg_4} and thus belong to \{3\} or \{33\} molecules by Lemma \ref{lem.toy_1}.

\textbf{Proof part 3.} With Parts 1 and 2 of the proof, we now know that every resulting molecule is elementary; it suffices to calculate $\#_{\{4\}}$ and $\#_{\{33\}}$. We first show that $\#_{\{4\}}=1$. Since the first atom cut in Definition \ref{def.toy_alg} \ref{it.toy_alg_2} must have deg 4, it suffices to show that $\#_{\{4\}}\leq 1$. In fact, for each $\nf\in\Mb_U$, at the time $\nf$ is cut, it must be the \emph{lowest remaining atom} in $\Mb_U$, which means any possible child of $\nf$ in $\Mb_U$ must have been cut before, creating a fixed end at $\nf$. Therefore, \emph{if $\nf$ has any child in $\Mb_U$ at all} before any cutting, then it cannot have deg 4 when it is cut.

On the other hand, before any cutting, $\Mb_U$ \emph{can have at most one atom without child in $\Mb_U$} (see {\color{blue} Figure \ref{fig.binary tree}}); if not, say there exist two such atoms $\nf$ and $\nf'$, then consider a path $(\nf=\nf_0,\cdots,\nf_j=\nf')$ between $\nf$ and $\nf'$ within $\Mb_U$ (since $\Mb_U$ is a tree). Then $\nf_1$ must be parent of $\nf$ as $\nf$ has no child in $\Mb_U$, and then $\nf_2$ must be parent of $\nf_1$ as $\nf_1$ has a child $\nf\in\Mb_U$ and another child in $\Mb_D$ by Definition \ref{def.toy}, and so on. In the end $\nf'=\nf_j$ must be parent of $\nf_{j-1}\in\Mb_U$, which contradicts the assumption on $\nf'$. This proves $\#_{\{4\}}=1$.
\begin{figure}[h!]
    \centering
    \includegraphics[width=0.35\linewidth]{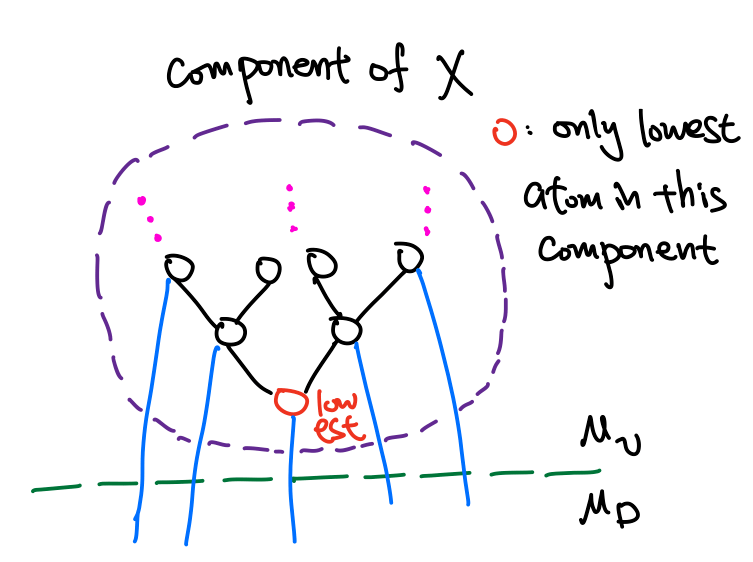}
    \caption{Illustration of each component of $X$ (or $\Mb_U$ if $X=\Mb_U$ in toy model I), in the proof of Propositions \ref{prop.toy} and \ref{prop.toy1+_alg}. It has a unique lowest atom (in red). The other atoms can only be constructed from the red one by taking \emph{parents}, as they must have a bond connecting to $\Mb_D$ (in blue) and cannot have a third child.}
    \label{fig.binary tree}
\end{figure}

\textbf{Proof part 4.} Finally we show that $\#_{\{33\}}\geq (|\Mb_U|-1)/5$. Consider the bahevior of the quantity
\[\lambda:=\#_{\mathrm{deg} 3}-\#_{\mathrm{comp}},\] where $\#_{\mathrm{deg} 3}$ and $\#_{\mathrm{comp}}$ represent the number of (remaining) deg 3 atoms and components in $\Mb_D$ respectively, so that initially $\lambda=-1$ and finally $\lambda=0$. Note that when viewed from $\Mb_D$, each step in Definition \ref{def.toy_alg} \ref{it.toy_alg_1}--\ref{it.toy_alg_6} is one of the followings:
\begin{itemize}
\item Operation \ref{it.stepa} or \ref{it.stepb} defined above (which occur in Definition \ref{def.toy_alg} \ref{it.toy_alg_2});
\item Cutting \{33\} or \{343\} molecule as free (which occur in \ref{it.toy_alg_3}--\ref{it.toy_alg_4} and \ref{it.toy_alg_6} in Definition \ref{def.toy_alg});
\item Cutting \{3\} molecule as free (which occur in Definition \ref{def.toy_alg} \ref{it.toy_alg_6}).
\end{itemize} Define the number of these operations to be $\#_{(a)}$, $\#_{(b)}$ and $\#_{\mathrm{cut}\{33\}}$ etc., then we have
\begin{equation}\label{eq.toy_2}\#_{\{33\}}=\#_{(b)}+\#_{\mathrm{cut}\{33\}}+\#_{\mathrm{cut}\{343\}}.\end{equation}

Moreover, using Lemma \ref{lem.cutconnected0} in the proof of Lemma \ref{lem.toy_1} (and the fact that $\Mb_D$ is proper right before each \ref{it.stepa} or \ref{it.stepb} or cutting of \{3\}), it is easy to see that
\begin{itemize}
\item Operation \ref{it.stepa} increased $\lambda$ by 1, and operation \ref{it.stepb} does not change $\lambda$;
\item Cutting \{33\} or \{343\} decreases $\lambda$ by $1$, and cutting \{3\} does not change $\lambda$.
\end{itemize} Using the initial and final values of $\lambda$ we deduce that \begin{equation}\label{eq.toy_3}\#_{(a)}-(\#_{\mathrm{cut}\{33\}}+\#_{\mathrm{cut}\{343\}})=1.\end{equation}Finally, consider the bonds between $\Mb_U$ and $\Mb_D$ (there are $|\Mb_U|$ of them). Each such bond is broken either in Definition \ref{def.toy_alg} \ref{it.toy_alg_2} (which corresponds to exactly one operation \ref{it.stepa} or \ref{it.stepb}), or in cutting one of \{33\} or \{343\} molecules; moreover each \ref{it.stepa}, each \ref{it.stepb}, each cutting of \{33\} and each cutting of \{343\} can break at most 1, 2, 3 and 4 such bonds respectively. This implies that $\#_{(a)}+2\#_{(b)}+3\#_{\mathrm{cut}\{33\}}+4\#_{\mathrm{cut}\{343\}}\geq|\Mb_U|$. Combining this with (\ref{eq.toy_2}) and (\ref{eq.toy_3}), we deduce that $\#_{\{33\}}\geq (|\Mb_U|-1)/5$, as desired.
\end{proof}
\subsection{Toy model I plus}\label{sec.toy1+} In this subsection we extend the toy model I to allow more general connections between $\Mb_U$ and $\Mb_D$.
\begin{definition}\label{def.toy1+} We say $\Mb$ is a \emph{toy model I plus}, if it is a 2-layer model as in Definition \ref{def.simplified}, and each atom $\nf\in\Mb_U$ has \emph{at most} one bond connecting to an atom in $\Mb_D$ (i.e. $\#_{\mathrm{2conn}}=0$), and the number of connected components of $X$ satisfies \begin{equation}\label{eq.toy1+_1}
\#_{\mathrm{comp}(X)}\ll |X|,
\end{equation} where $X$ is defined at the beginning of this note. Here, if two atoms in $X$ can only be connected via atoms in $\Mb_U\backslash X$, then they belong to different components of $X$.
\end{definition}
Let us explain the intuition behind Definition \ref{def.toy1+}, in particular the role of (\ref{eq.toy1+_1}). In the case of toy model I, we have $X=\Mb_U$ which is connected, hence $\#_{\mathrm{comp}(X)}=1$. Moreover, as in Part 3 of the proof of Proposition \ref{prop.toy}, we see that $\Mb_U$ has at most one lowest atom (i.e. atom without child in $\Mb_U$; in fact it is not hard to see that $\Mb_U$ must have an ``upside down binary tree" shape). This latter fact allows us to apply a ``naive" strategy in $\Mb_U$ where we \emph{always} cut the lowest remaining atom (and avoid any top fixed ends), \emph{without producing too many \{4\} atoms}. 

In the current case, $\Mb_U$ may have many lowest atoms, so we cannot always choose and cut a lowest remaining atom. Instead, we improve the strategy by choosing a \emph{lowest remaining deg 3} atom $\nf$ (here we do not consider deg 2 atoms, as they will always be lowest in $\Mb_U$ and always have to be cut at top priority). Now cutting $\nf$ may create top fixed ends, but only at descendants of $\nf$ (which all have deg 4 when $\nf$ is chosen). As such, we need to subsequently cut all descendants of $\nf$ after cutting $\nf$. This leads to the idea of the \textbf{UP} algorithm, see Definition \ref{def.up_toy} below, which is another basic ingredient in our proof.

Note, however, that this \textbf{UP} algorithm cannot be simply combined with the ``greedy" algorithm in $\Mb_D$ as in Definition \ref{def.toy_alg}, because $\nf$ may have (deg 4) descendants in $\Mb_U$, and this naive attempt may break the key monotonicity property (see Part 1 of the proof of Proposition \ref{prop.toy1+_alg}) and lead to non-elementary molecules. As such, we should only apply \textbf{UP} \emph{at places in $\Mb_U$ which will not affect $\Mb_D$}, i.e. in $\Mb_U\backslash X$. Within $X$ we still apply the naive strategy, which may lead to \{4\} molecules due to $X$.

We can then see the role of (\ref{eq.toy1+_1}) and the \emph{dichotomy} concerning $\#_{\mathrm{comp}(X)}$; in fact, same as in Part 3 of the proof of Proposition \ref{prop.toy}, \emph{each component} of $X$ can have \emph{at most one lowest atom}. In the naive strategy, this leads to \{4\} molecules whose number is comparable to $\#_{\mathrm{comp}(X)}$. Under the assumption (\ref{eq.toy1+_1}), this is negeligible compared to the number of good \{33\} molecules, which will be comparable to $|X|$ (i.e. the number of bonds between $\Mb_U$ and $\Mb_D$).

With the above intuition, we start by defining the \textbf{UP} algorithm in the toy model case. For the full version in the real problem, see Definition \ref{def.alg_up}.
\begin{definition}[The \textbf{UP} algorithm, toy model case]\label{def.up_toy} Let $\Mb$ be any molecule which \emph{does not have top fixed end}. Define the following cutting sequence. For an example see {\color{blue} Figure \ref{fig.up}}:
\begin{figure}[h!]
    \centering
    \includegraphics[width=0.65\linewidth]{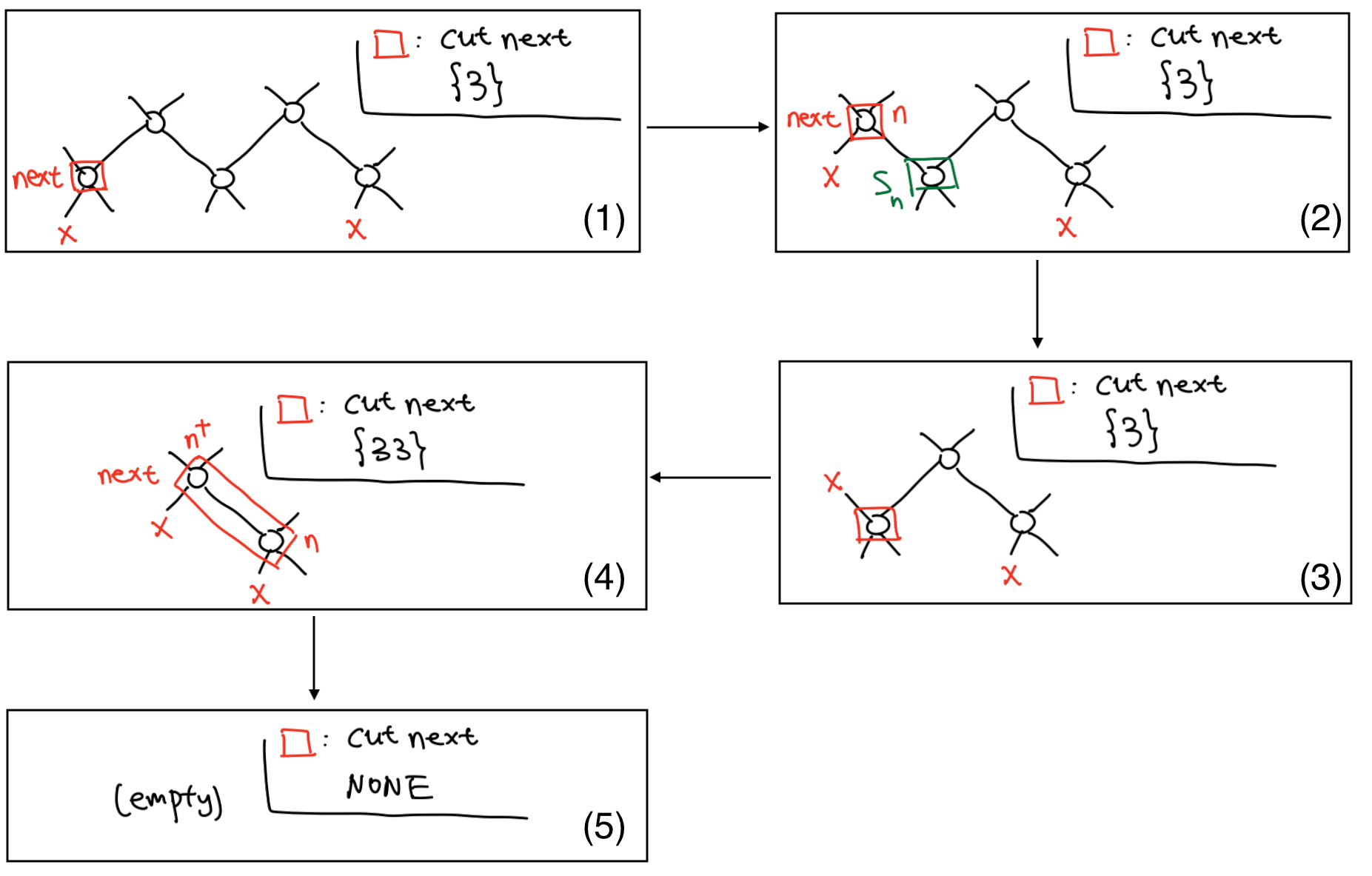}
    \caption{The \textbf{UP} algorithm (Definition \ref{def.up_toy}). As in the above figures, the red box indicates the atom(s) to be cut next, the green box indicates the atom(s) in $S_\nf$ (see Definition \ref{def.up_toy} \ref{it.up_toy_2}), which will be cut after cutting $\nf$. Note also when $\nf$ and its parent $\nf^+$ (or child $\nf^-$) form a \{33\} molecule, we should cut them together; see Definition \ref{def.up_toy} \ref{it.up_toy_3}.}
    \label{fig.up}
\end{figure}
\begin{enumerate}[{(1)}]
\item\label{it.up_toy_1} If $\Mb$ has a deg 2 atom, then cut it as free. Repeat until $\Mb$ does not have deg 2 atoms.
\item\label{it.up_toy_2} Choose a lowest deg 3 atom $\nf$ that has not bee cut (or a lowest atom if $\Mb$ has no deg 3 atoms). Let $S_\nf$ be the set of descendants of $\nf$. This $\nf$ and $S_\nf$ is fixed until all atoms in $S_\nf$ have been cut.
\item\label{it.up_toy_3} Starting from $\nf$, each time choose a highest atom $\mf$ in $S_\nf$ that has not been cut. If $\mf$ has deg 3 and has a parent $\mf^+$ or child $\mf^-$ that also has deg 3, then cut $\{\mf,\mf^\pm\}$ as free; otherwise cut $\mf$ as free. Repeat until all atoms in $S_\nf$ have been cut, then go to \ref{it.up_toy_1}--\ref{it.up_toy_2} and choose the next $\nf$, and so on.
\end{enumerate}
\end{definition}
\begin{remark}
\label{rem.up_toy} The \textbf{UP} algorithm is a combination of a ``global low-to-high" (where each time we select the \emph{lowest} remaining deg 3 atom in Definition \ref{def.up_toy} \ref{it.up_toy_2}) and ``local high-to-low" procedure (where each time we cut the \emph{highest} remaining atom in $S_\nf$ in Definition \ref{def.up_toy} \ref{it.up_toy_3}). This leads to a key \emph{monotonicity} property, see Proposition 
\ref{prop.up_toy}.

The step \ref{it.up_toy_3} in Definition \ref{def.up_toy} may seem complicated, we explain it as follows: let $\mf$ be the highest remaining atom in $S_\nf$, the simple option is to just cut $\mf$ as free. However, in the case described in Definition \ref{def.up_toy} \ref{it.up_toy_3}, the $\mf^\pm$ will become deg 2 once $\mf$ is cut as free; moreover $\mf^-\in S_\nf$ will be a highest remaining atom in $S_\nf$, and $\mf^+\not\in S_\nf$ will occur as a deg 2 atom in Definition \ref{def.up_toy} \ref{it.up_toy_1} once all atoms in $S_\nf$ have been cut. As such, this deg 2 atom $\mf^\pm$ have to be cut at top priority, so instead of cutting $\{\mf\}$ we can combine this with the cutting of $\mf^\pm$, which also forms a good \{33\} molecule $\{\mf,\mf^\pm\}$. This then leads to the step \ref{it.up_toy_3} in Definition \ref{def.up_toy}.
\end{remark}
We next prove a key monotonicity property of the \textbf{UP} algorithm and a sufficient condition for it to produce at least one \{33\} molecule (in fact \{33A\} molecule). See Proposition \ref{prop.alg_up} for the full version.
\begin{proposition}\label{prop.up_toy} Let $\Mb$ be as in Definition \ref{def.up_toy} and connected (otherwise consider each component of $\Mb$). Then after applying algorithm \textbf{UP} to $\Mb$ (and same for \textbf{DOWN}):
\begin{enumerate}
\item\label{it.up_toypf_1} We have the following \emph{monotonicity} property, which we refer to as (MONO-toy): at any time, any atom in $(\Mb\backslash S_\nf)\cup\{\nf\}$ has no \emph{top fixed end}, and any atom in $S_\nf\backslash\{\nf\}$ has no \emph{bottom fixed end}. In particular, applying \textbf{UP} only produces elementary molecules.
\item\label{it.up_toypf_2} We have $\#_{\{4\}}\leq 1$, and $\#_{\{4\}}=1$ if and only if $\Mb$ is full. 
\end{enumerate}
In (\ref{it.up_toypf_3}) below we assume $\Mb$ has no deg 2 atoms.
\begin{enumerate}[resume]
\item\label{it.up_toypf_3} If $\Mb$ contains a cycle or contains at least two deg $3$ atoms, then $\#_{\{33\}}\geq 1$. If not (i.e. $\Mb$ is a tree and contains at most one deg 3 atom), then all the elementary molecules are \{3\} molecules with at most one exception of \{4\} molecule.
\end{enumerate}
\end{proposition}
\begin{proof} \textbf{Proof of (\ref{it.up_toypf_1}).} We first prove (MONO-toy), which implies all components are elementary (since each atom is either cut by itself or belongs to a \{33\} molecule, due to \ref{it.up_toy_1}--\ref{it.up_toy_3} in Definition \ref{def.up_toy}). Initially there is no $S_\nf$ and $\Mb$ has no top fixed end by assumption, so it suffices to prove that (MONO-toy) is perserved. First, any deg 2 atom $\nf$ cut in Definition \ref{def.up_toy} \ref{it.up_toy_1} must have two \emph{bottom} fixed ends by (MONO-toy), so it has no child remaining when it is cut, thus cutting it will not create any \emph{top fixed end}.

Now let $\nf$ and $S_\nf$ be as in Definition \ref{def.up_toy} \ref{it.up_toy_2}, and $\mf$ be a highest remaining atom in $S_\nf$ as in Definition \ref{def.up_toy} \ref{it.up_toy_3}. When $\nf$ is chosen in Definition \ref{def.up_toy} \ref{it.up_toy_2}, all atoms in $S_\nf\backslash\{\nf\}$ \emph{must have deg 4} by construction, and thus have no bottom fixed end. Consider the operation of cutting $\mf$ as free in Definition \ref{def.up_toy} \ref{it.up_toy_3}. We then know that that this cutting does not create any \emph{bottom fixed end} in $S_\nf\backslash\{\nf\}$ (as $\mf$ is highest), and does not create any \emph{top fixed end in $(\Mb\backslash S_\nf)\cup\{\nf\}$} (as any child of $\mf$ must belong to $S_\nf\backslash\{\nf\}$, see {\color{blue} Figure \ref{fig.monotonicity}}). The same holds for the operation of cutting $\{\mf,\mf^\pm\}$, because this is equivalent to cutting $\{\mf\}$ and then cutting a deg 2 atom (either $\mf^+\not\in S_\nf$ with two bottom fixed ends or $\mf^-\in S_\nf\backslash\{\nf\}$ with two top fixed ends). This proves that (MONO-toy) is perserved, as desired.
\begin{figure}[h!]
    \centering
    \includegraphics[width=0.65\linewidth]{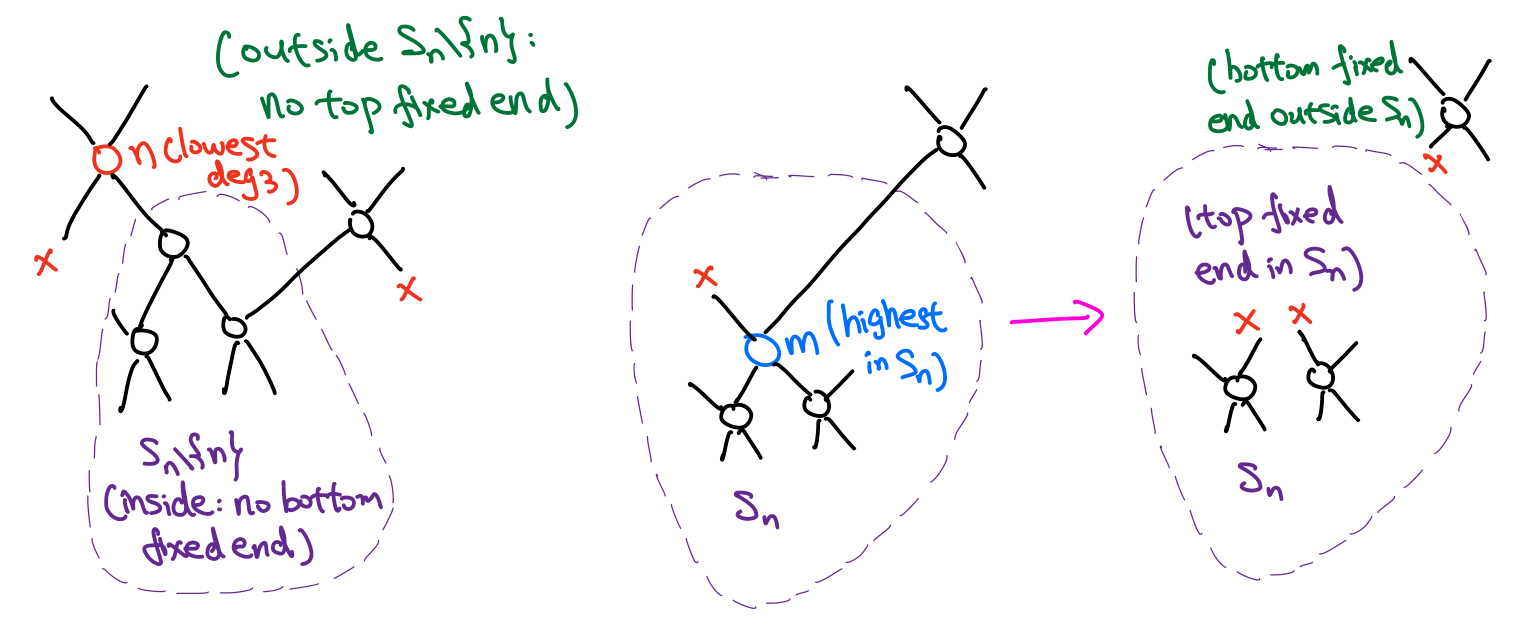}
    \caption{Illustration of the (MONO-toy) property of the \textbf{UP} algorithm (Proposition \ref{prop.up_toy}). Left: when we choose $\nf$ in Definition \ref{def.up_toy} \ref{it.up_toy_2}, there is no top fixed end outside $S_\nf\backslash\{\nf\}$ and no fixed end at all in $S_\nf\backslash\{\nf\}$ (i.e. all deg 4). Right: when we cut a highest remaining atom $\mf$ in $S_\nf$, we only get top fixed end in $S_\nf$ and bottom fixed end outside $S_\nf$.}
    \label{fig.monotonicity}
\end{figure}

\textbf{Proof of (\ref{it.up_toypf_2}).} Note that after the first cutting, $\Mb$ will \emph{always} contain an atom that is not deg 4 (cf. Remark \ref{rem.full_cut}). By construction, this means that the $\nf$ chosen in Definition \ref{def.up_toy} \ref{it.up_toy_2} will always be deg 3; moreover, any $\mf\in S_\nf\backslash \{\nf\}$ that is cut in Definition \ref{def.up_toy} \ref{it.up_toy_3} also cannot be deg 4, because $\mf$ has a parent in $S_\nf$ which must have been cut before $\mf$ thanks to the highest assumption of $\mf$. Thus \emph{only} the first atom can produce a possible \{4\} molecule, so $\#_{\{4\}}\leq 1$; moreover $\#_{\{4\}}=1$ if and only if the first atom cut has deg 4, i.e. when $\Mb$ is full.

\textbf{Proof of (\ref{it.up_toypf_3}).} Suppose $\Mb$ has no deg 2 atom, we first prove that either $\#_{\{33\}}\geq 1$, or all resulting molecules are \{4\} and \{3\} (with $\#_{\{4\}}\leq 1$). If not, then we must have a \{2\} molecule. Consider the first atom $\mf$ that has deg 2 when we cut it, and consider the cutting operation that turns $\mf$ into deg 2, which must be cutting some deg 3 atom $\pf$ as free (this $\pf$ cannot have deg 4, because a degree 4 atom is cut only when no degree 3 atoms are present, and in this case cutting a degree 4 atom cannot produce a degree 2 atom).  This cutting creates a fixed end at $\mf$, so this $\pf$ must be a parent or child of $\mf$; moreover this $\mf$ and this $\pf$ must both be deg 3 at the time of cutting $\pf$, but this violates Definition \ref{def.up_toy} \ref{it.up_toy_3}.

Finally we prove $\#_{\{33\}}=0$ if and only if $\Mb$ is a tree and has at most one deg 3 atom. For this, consider the increment of the quantity $\sigma:=\#_{\mathrm{bo/fr}}-3\#_{\mathrm{atom}}$, where $\#_{\mathrm{bo/fr}}$ is the number of bonds plus the number of free ends (not counting fixed ends) and $\#_{\mathrm{atom}}$ is the number of remaining atoms. It is easy to see that cutting each \{2\}, \{3\}, \{33\} and \{4\} molecule increases $\sigma$ by 1, 0, 1 and $-1$ respectively. Since finally $\sigma=0$, we see that $\#_{\{4\}}-(\#_{\{2\}}+\#_{\{33\}})=\sigma_{\mathrm{init}}$ (initial value of $\sigma$).

Note that initially $\#_{\mathrm{bo/fr}}=\#_{\mathrm{deg}}-\#_{\mathrm{bond}}$ where $\#_{\mathrm{deg}}$ is the total degree of all atoms (not counting fixed ends) and $\#_{\mathrm{bond}}$ is the number of bonds (because each bond contributes twice to the degrees), so $\sigma_{\mathrm{init}}=(|\Mb|-\#_{\mathrm{bond}})-|\Mb_3|$ where $|\Mb_3|$ is the number of (initially) deg 3 atoms in $\Mb$. Moreover by connectedness we have $\#_{\mathrm{bond}}\geq |\Mb|-1$ with equality if and only if $\Mb$ is a tree, and note also $\#_{\{4\}}\in\{0,1\}$ (and $\#_{\{4\}}=1$ if and only if $|\Mb_3|=0$ by the same proof in Part 2), it is then easy to see that $\#_{\{2\}}+\#_{\{33\}}=0$ if and only if $\Mb$ is a tree and $|\Mb_3|\leq 1$. This completes the proof.
\end{proof}
Now we can define the algorithm for toy model I plus. As described in the above heuristics, this features a naive algorithm (always cutting the lowest) in $X$, an \textbf{UP} algorithm in $\Mb_U\backslash X$, and a greedy algorithm in $\Mb_D$. Up to some simplifications, this is just the \textbf{MAINUD} algorithm, which is the main algorithm in Section \ref{sec.maincr}, see Definition \ref{def.alg_maincr}. Note also that the part of Definition \ref{def.toy1+_alg} involving $\Mb_D$ is exactly the same as Definition \ref{def.toy_alg}; the only difference occurs in $\Mb_U$.
\begin{definition}[Algorithm for toy model I plus]\label{def.toy1+_alg} Let $\Mb$ be a toy model I plus with $X$ be defined in the beginning of this note (this $X$ is fixed throughout the algorithm). We define the following cutting sequence. For an example see {\color{blue}Figures \ref{fig.toy model 1+ alg1}--\ref{fig.toy model 1+ alg2}}:
\begin{figure}[h!]
    \centering
    \includegraphics[width=0.85\linewidth]{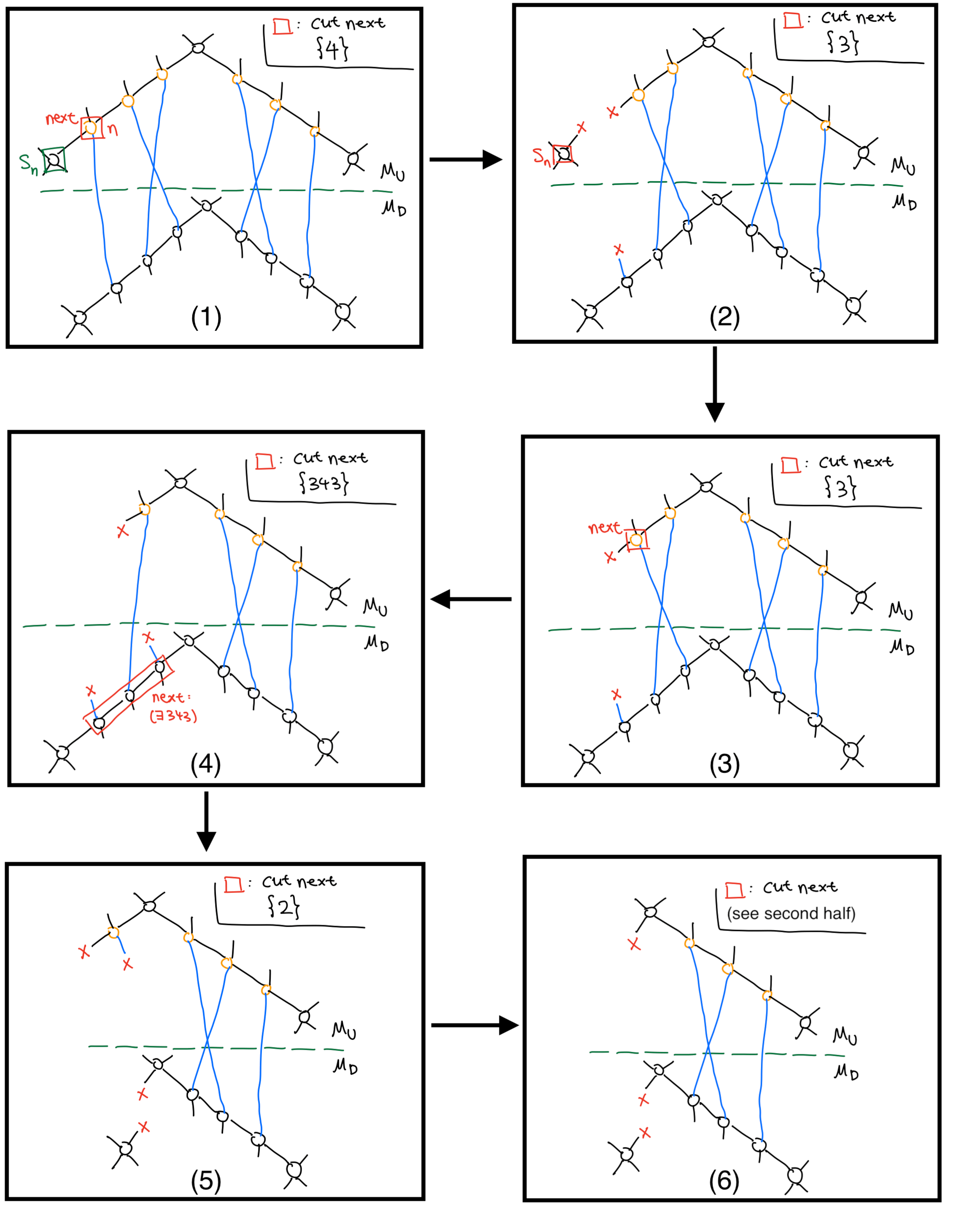}
    \caption{Algorithm for toy model I plus (first half) in Definition \ref{def.toy1+_alg}, applied to the example in {\color{blue} Figure \ref{fig.toy models}}. In addition to {\color{blue} Figure \ref{fig.toy model 1 alg}}, we use green box to indicate atom(s) in $S_\nf$ (see Definition \ref{def.toy1+_alg} \ref{it.toy1+_alg_2}) which will be cut after cutting $\nf$.}
    \label{fig.toy model 1+ alg1}
\end{figure}
\begin{enumerate}[{(1)}]
\item\label{it.toy1+_alg_1} If $\Mb_U$ has a deg 2 atom, then cut it as free. Repeat until $\Mb_U$ does not have deg 2 atoms.
\item\label{it.toy1+_alg_2} Consider the set of all (remaining) atoms that either belong to $X$, or has deg 3; choose $\nf$ to be a lowest atom in this set (if this set is empty then choose $\nf$ to be any lowest atom). Let $S_\nf$ be the set of descendants of $\nf$ in $\Mb_U$, and fix $\nf$ and $S_\nf$ until the end of (\ref{it.alg_maincr_7}).
\item\label{it.toy1+_alg_3} Starting from $\nf$, each time choose a highest atom $\pf$ in $S_\nf$ that has not been cut. If $\pf$ is adjacent to an atom $\qf\in\Mb_D$ which has deg 3, then cut $\{\pf,\qf\}$ as free and cut $\qf$ as free from $\{\pf,\qf\}$ if $\pf$ has deg 4; otherwise just cut $\pf$ as free.
\item\label{it.toy1+_alg_4} If $\Mb_D$ becomes non-proper, then perform the same steps \ref{it.toy_alg_3}--\ref{it.toy_alg_4} in Definition \ref{def.toy_alg} until it becomes proper again. Then go to \ref{it.toy1+_alg_3} and repeat, until all atoms in $S_\nf$ have been cut.
\item\label{it.toy1+_alg_5} When all atoms in $S_\nf$ have been cut, go to \ref{it.toy1+_alg_1}--\ref{it.toy1+_alg_2} and choose the next $\nf$, and so on.
\item\label{it.toy1+_alg_7} When all atoms in $\Mb_U$ have been cut, perform the same step \ref{it.toy_alg_6} in Definition \ref{def.toy_alg} to finish off $\Mb_D$.
\end{enumerate}
\end{definition}
\begin{figure}[h!]
    \centering
    \includegraphics[width=1\linewidth]{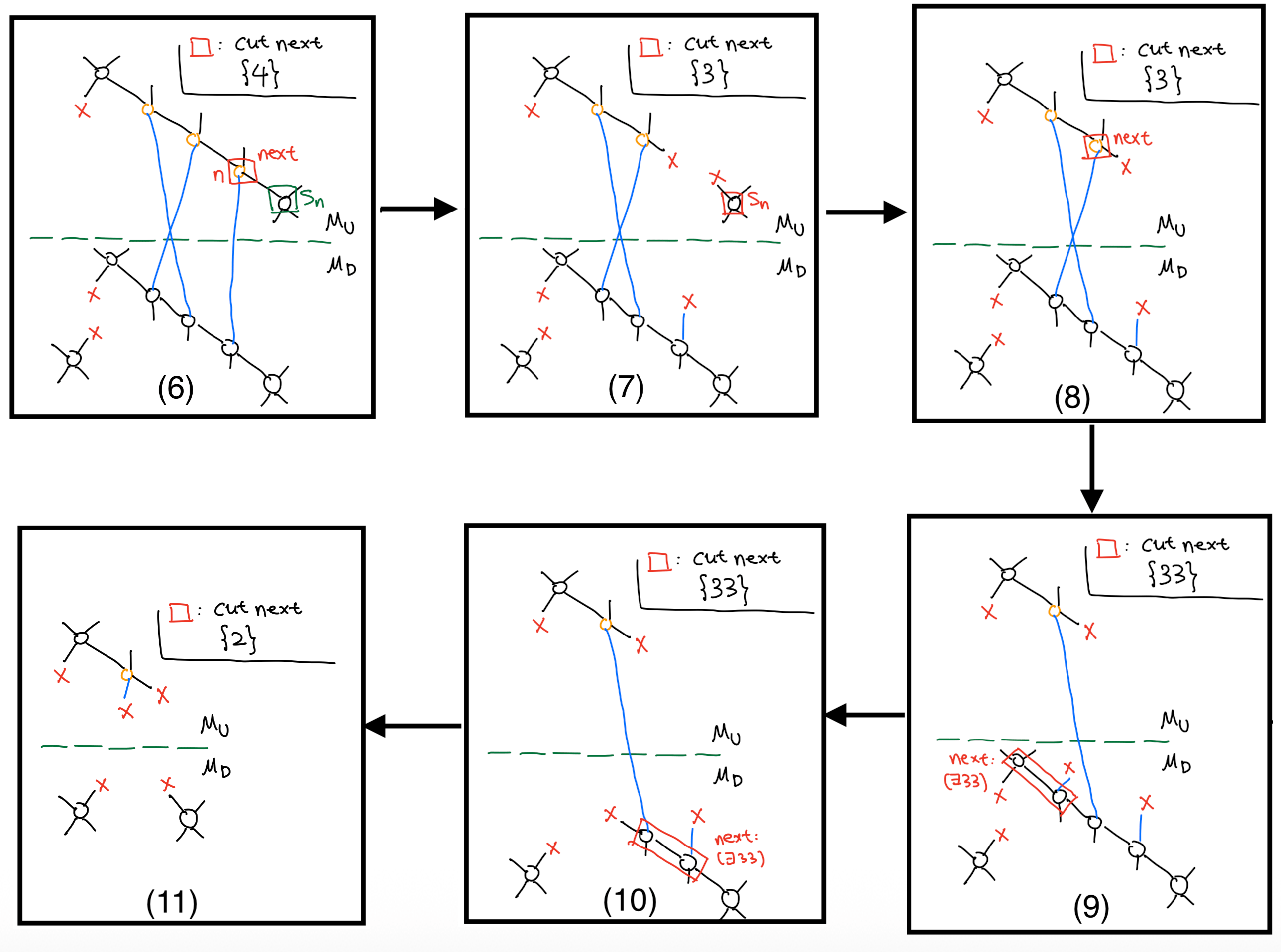}
    \caption{Algorithm for toy model I plus (second half). Here, the next atom $\nf$ we choose is not the unique deg 3 atom in $\Mb_U$, but rather the \emph{atom in $X$} below it, see Definition \ref{def.toy1+_alg} \ref{it.toy1+_alg_2}. This leads to another \{4\} loss, but the total loss is bounded by $\#_{\mathrm{comp}(X)}=2$. In the end, we subsequently cut two deg 2 atoms in $\Mb_U$, and two deg 3 atoms in $\Mb_D$, to conclude.}
    \label{fig.toy model 1+ alg2}
\end{figure}
\begin{proposition}\label{prop.toy1+_alg} Let $\Mb$ be a toy model I plus, and we apply the algorithm defined in Definition \ref{def.toy1+_alg}. Then in this process we only obtain elementary molecules, and we have
\begin{equation}\label{eq.toy1+_2}\#_{\{4\}}\leq \#_{\mathrm{comp}(X)},\qquad \#_{\{33\}}\geq (|X|-1)/5-\#_{\mathrm{comp}(X)}.\end{equation}
Since $|X|\gtrsim\rho\gg 1$ and using (\ref{eq.toy1+_1}), this clearly implies (\ref{eq.goal}).
\end{proposition}
\begin{proof} We divide the proof into 4 parts just as Proposition \ref{prop.toy}.

\textbf{Proof part 1.} We show that each atom in $\Mb_U$ is contained in an elementary molecule. As in Part 1 of the proof of Proposition \ref{prop.up_toy}, this follows from the same (MONO-toy) property but \emph{for $\Mb_U$ only}: any atom in $(\Mb_U\backslash S_\nf)\cup\{\nf\}$ has no \emph{top fixed end}, and any atom in $S_\nf\backslash\{\nf\}$ has no \emph{bottom fixed end}.

To prove the monotonicity property, we first note that, for any $\nf$ chosen in Definition \ref{def.toy1+_alg} \ref{it.toy1+_alg_2} and at the time it is chosen, each atom in $S_\nf\backslash\{\nf\}$ has deg 4 and thus has no bottom fixed end. Then consider the process of cutting in Definition \ref{def.toy1+_alg} \ref{it.toy1+_alg_3} and below, see {\color{blue} Figure \ref{fig.mono2layer}}. Following exactly Part 1 of the proof of Proposition \ref{prop.up_toy}, we see that cutting any atom in $\Mb_U$ in this process does not create any top fixed end in $(\Mb_U\backslash S_\nf)\cup\{\nf\}$ or any bottom fixed end in $S_\nf\backslash\{\nf\}$; moreover, for the cutting of any atom in $\Mb_D$ in this process (such as the one in Definition \ref{def.toy1+_alg} \ref{it.toy1+_alg_4}), it \emph{only creates bottom fixed ends} in $\Mb_U$, and \emph{only at atoms in $(\Mb_U\backslash S_\nf)\cup\{\nf\}$}. This is because $S_\nf\backslash\{\nf\}$ contains no atom in $X$ by construction, and thus has no bond connecting to $\Mb_D$, hence no fixed end can be created at any atom in this set when cutting anything in $\Mb_D$. This proves (MONO-toy) for $\Mb_U$.
\begin{figure}[h!]
    \centering
    \includegraphics[width=0.5\linewidth]{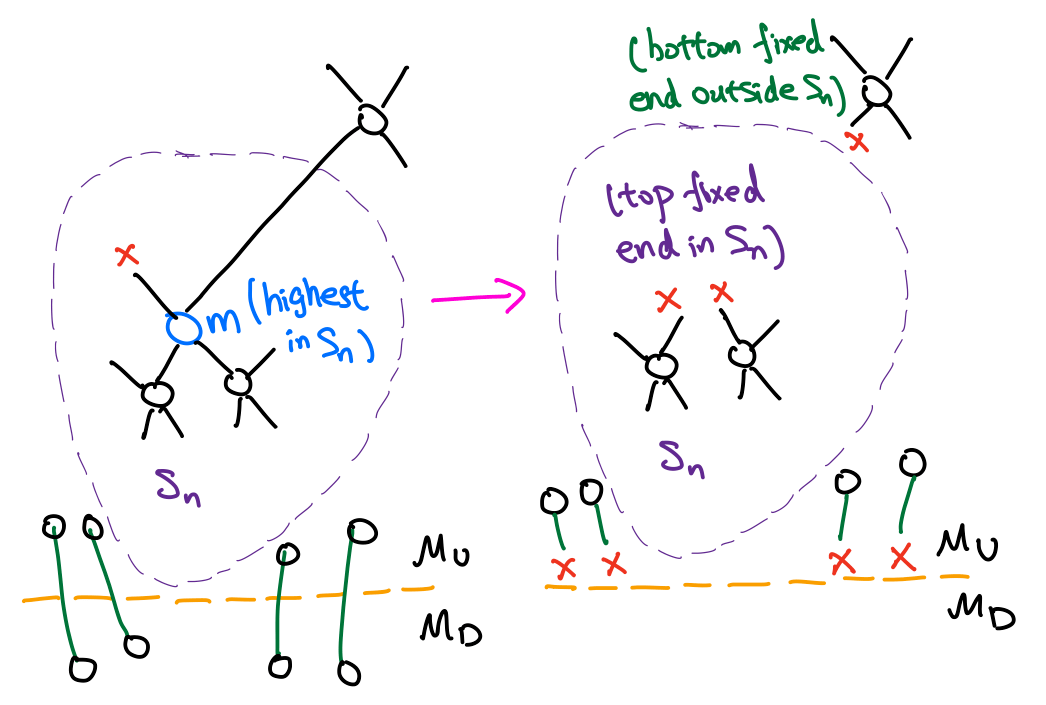}
    \caption{Part 1 of the proof of Proposition \ref{prop.toy1+_alg}: the preservation of the monotonicity property (just as for \textbf{UP}) in $\Mb_U$. This is essentially the same as the right half of {\color{blue} Figure \ref{fig.monotonicity}}, except the green bonds connecting $\Mb_D$ to atoms in $\Mb_U$ (\textbf{those atoms must be in $X$}). By our choice, $S_\nf\backslash\{\nf\}$ is \textbf{disjoint} with $X$, so none of these bonds involve any atom in $S_\nf$ (once $\nf$ has been cut). Thus, if we do any cutting in $\Mb_D$ (as in Definition \ref{def.toy1+_alg} \ref{it.toy1+_alg_4}), this only creates bottom ends outside $S_\nf\backslash\{\nf\}$, hence preserving the monotonicity property in $\Mb_U$.}
    \label{fig.mono2layer}
\end{figure}

\textbf{Proof part 2.} We show that each atom in $\Mb_D$ is contained in an elementary molecule. Since the part of Definition \ref{def.toy1+_alg} involving $\Mb_D$ is exactly the same as Definition \ref{def.toy_alg}, this part is also exactly the same as Part 2 of the proof of Proposition \ref{prop.toy}.

\textbf{Proof part 3.} We show that $\#_{\{4\}}\leq \#_{\mathrm{comp}(X)}$. Same as Part 2 of the proof of Proposition \ref{prop.up_toy}, we see that any $\pf\in S_\nf\backslash \{\nf\}$ in Definition \ref{def.toy1+_alg} \ref{it.toy1+_alg_3} cannot have deg 4 when it is cut. Therefore, the only possible \{4\} atoms must come from those $\nf$ chosen in Definition \ref{def.toy1+_alg} \ref{it.toy1+_alg_2}. By our choice in Definition \ref{def.toy1+_alg} \ref{it.toy1+_alg_2}, such $\nf$ must belong to $X$ and must be \emph{lowest atoms} in $X$, i.e. it belongs to $X$ but does not have any child in $X$ before any cutting.

Now we claim that any two lowest atoms in $X$ \emph{must belong to different components of $X$} (see {\color{blue} Figure \ref{fig.binary tree}}). If not, say $\nf$ and $\nf'$ are connected by a path $(\nf=\nf_0,\cdots,\nf_j=\nf')$ in $X$, then arguing as Part 3 of the proof of Proposition \ref{prop.toy} and noticing that each atom in $X$ has a child in $\Mb_D$, we deduce that each $\nf_i$ must be parent of $\nf_{i-1}$ which then leads to contradiction with properties of $\nf_j=\nf'$. We then conclude that the number of lowest atoms in $X$ does not exceed $\#_{\mathrm{comp}(X)}$, hence $\#_{\{4\}}\leq \#_{\mathrm{comp}(X)}$.

\textbf{Proof part 4.} Finally we show that $\#_{\{33\}}\geq (|X|-1)/5-\#_{\mathrm{comp}(X)}$. Since the part of Definition \ref{def.toy1+_alg} involving $\Mb_D$ is exactly the same as Definition \ref{def.toy_alg}, this part is essentially the same as Part 4 of the proof of Proposition \ref{prop.toy}. There is only one difference: in the current case, each step of cutting $\{\pf,\qf\}$ as free in Definition \ref{def.toy1+_alg} \ref{it.toy1+_alg_3} (which corresponds to the operation \ref{it.stepb} in Part 4 of the proof of Proposition \ref{prop.toy} when viewed in $\Mb_D$) may produce two \{3\} molecules instead of one \{33\} molecule, but this only happens when\footnote{This cannot happen for toy model I, because there $\pf$ can have deg 4 only when it is the first atom cut, in which case $\qf$ cannot have deg 3.} $\pf$ has deg 4, and the number of such exceptions is bounded by $\#_{\mathrm{comp}(X)}$ following Part 3 above. Therefore, instead of (\ref{eq.toy_2}) we have
\begin{equation}\label{eq.toy1+_3}\#_{\{33\}}\geq\#_{(b)}+\#_{\mathrm{cut}\{33\}}+\#_{\mathrm{cut}\{343\}}-\#_{\mathrm{comp}(X)},\end{equation} and the same proof then yields that, instead of $\#_{\{33\}}\geq (|X|-1)/5$ as in (\ref{eq.toy_1}), we now have $\#_{\{33\}}\geq (|X|-1)/5-\#_{\mathrm{comp}(X)}$. This completes the proof.
\end{proof}
\subsection{Toy model II}\label{sec.toy2} In this subsection we introduce the toy model II, which is opposite to toy model I plus via the dichotomy, where we assume $\#_{\mathrm{comp}(X)}\gtrsim |X|$ instead of $\#_{\mathrm{comp}(X)}\ll |X|$ in Definition \ref{def.toy1+}.
\begin{definition}\label{def.toy2} We say $\Mb$ is a toy model II, if it satisfies all the assumptions in Definition \ref{def.toy1+} for toy model I plus with the same notations, except that (\ref{eq.toy1+_1}) is replaced by $|X|\geq 2$ and
\begin{equation}\label{eq.toy2_1}
\#_{\mathrm{comp}(X)}\gtrsim |X|.
\end{equation}
\end{definition}
Let us explain the strategy for toy model II, which will be completely different from toy models I and I plus. Here we will not gain from the greedy algorithm in $\Mb_D$; instead we first cut all of $\Mb_D$ as free without exploiting any gain, and then gain from $\Mb_U$ under the assumption (\ref{eq.toy2_1}). Note that, after cutting $\Mb_D$ as free, there will not be any deg 2 atoms in $\Mb_U$, and the deg 3 atoms in $\Mb_U$ are \emph{exactly those in $X$}.

The intuition is that, due to the lower bound of the number of components if $X$ given by (\ref{eq.toy2_1}), the deg 3 atoms in $\Mb_U$ are now \emph{separated by deg 4 atoms}. Essentially, we expect that for such molecules, we can design a cutting algorithm with minimal loss caused \{4\} molecules, such that \emph{each connected component of the set of deg 3 atoms} provides a good \{33\} molecule. In fact, this can be achieved by a clever choice in the \textbf{UP} algorithm.

\begin{figure}[h!]
    \centering
    \includegraphics[width=0.45\linewidth]{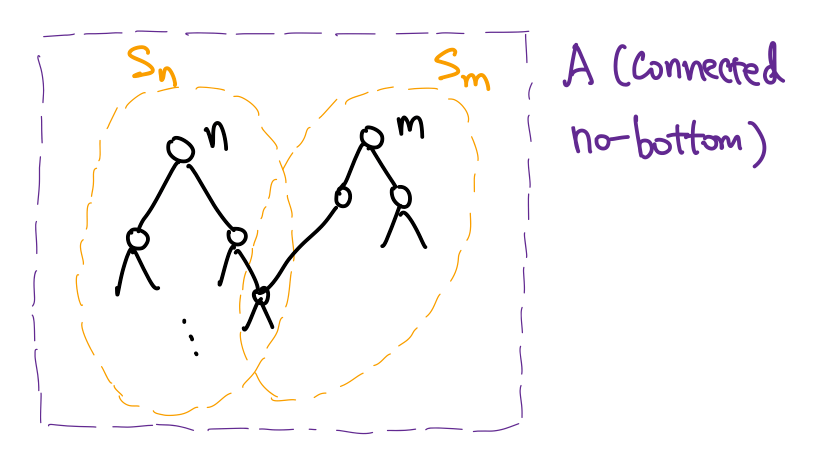}
    \caption{An example of connected, no-bottom sets in Lemma \ref{lem.toy2_aux}; this one is a union of two intersecting descendant sets $S_\mf$ and $S_\nf$. In general each such set is a union of descendant sets $S_\nf$ with $\nf$ forming a connected set; see the proof of Lemma \ref{lem.toy2_aux}.}
    \label{fig.nobottom}
\end{figure}

We next introduce the main algorithm in this case (together with an auxiliary lemma needed in this algorithm). Up to simplifications, this is just the \textbf{3COMPUP} algorithm in Definition \ref{def.3comp_alg}.
\begin{lemma}\label{lem.toy2_aux} Suppose $\Mb$ is a molecule that is a tree, contains only deg 3 and 4 atoms and has no top fixed end. Let the set of deg 3 atoms in $\Mb$ be $\Mb_3$ with $|\Mb_3|\geq 2$. Define a set $A\subseteq\Mb$ to be \emph{no-bottom} if any child of any atom in $A$ must also be in $A$ (see {\color{blue}Figure \ref{fig.nobottom}}). Then, there always exists a connected, no-bottom subset $A\subseteq \Mb$ with $|A\cap \Mb_3|=2$.
\end{lemma}
\begin{proof} For any $\nf\in\Mb$, let $S_\nf$ be the set of descendants of $\nf$ (including $\nf$ itself), then it is always connected and no-bottom. If there exists $\nf$ such that $|S_\nf\cap \Mb_3|\geq 2$, then choose a lowest such $\nf$. Since $\nf$ has no top fixed end, we know that either (i) $\nf$ has one child $\nf_1$, or (ii) $\nf$ has deg 4 and has two children $\nf_1,\nf_2$. In either case, by assumption we have $|S_{\nf_j}\cap \Mb_3|\leq 1$ for each $j$, which then implies $|S_\nf\cap \Mb_3|\leq 2$ because $S_\nf=(\cup_j S_{\nf_j})\cup\{\nf\}$. This implies $|S_\nf\cap \Mb_3|= 2$, as desired.

Now assume $|S_\nf\cap \Mb_3|\leq 1$ for each $\nf\in\Mb$. If $|S_\nf\cap \Mb_3|=1$, there is a unique $\pf=\pf(\nf)\in S_\nf\cap \Mb_3$; now consider all atom pairs $(\nf,\nf')$ such that $\pf(\nf)$ and $\pf(\nf')$ both exist and \emph{are different} (which exist because $|\Mb_3|\geq 2$), and one such pair $(\nf,\nf')$ that are connected by \emph{a shortest path} in $\Mb$. Let this path be $(\nf_0,\nf,\cdots,\nf_j=\nf')$, then for each $0<i<j$ we must have either $S_{\nf_i}\cap \Mb_3=\varnothing$ or $\pf(\nf_i)=\pf(\nf)$ (otherwise the pair $(\nf,\nf_i)$ corresponds to a shorter path, contradiction). Now let $A=\cup_{i=0}^j S_{\nf_i}$, it is clear that $A$ is connected, no-bottom and $A\cap \Mb_3=\{\pf(\nf),\pf(\nf')\}$, as desired. 
\end{proof}
\begin{definition}[The \textbf{3COMPUP} algorithm, toy model case]\label{def.toy2_alg} Let $\Mb$ be a molecule that is a tree, contains only deg 3 and 4 atoms and has no top fixed end. Let the set of deg 3 atoms in $\Mb$ be $\Mb_3$ with $|\Mb_3|\geq 2$. Define the following cutting sequence. For an example see {\color{blue}Figure \ref{fig.toy model 2 alg}}:
\begin{figure}[h!]
    \centering
    \includegraphics[width=0.55\linewidth]{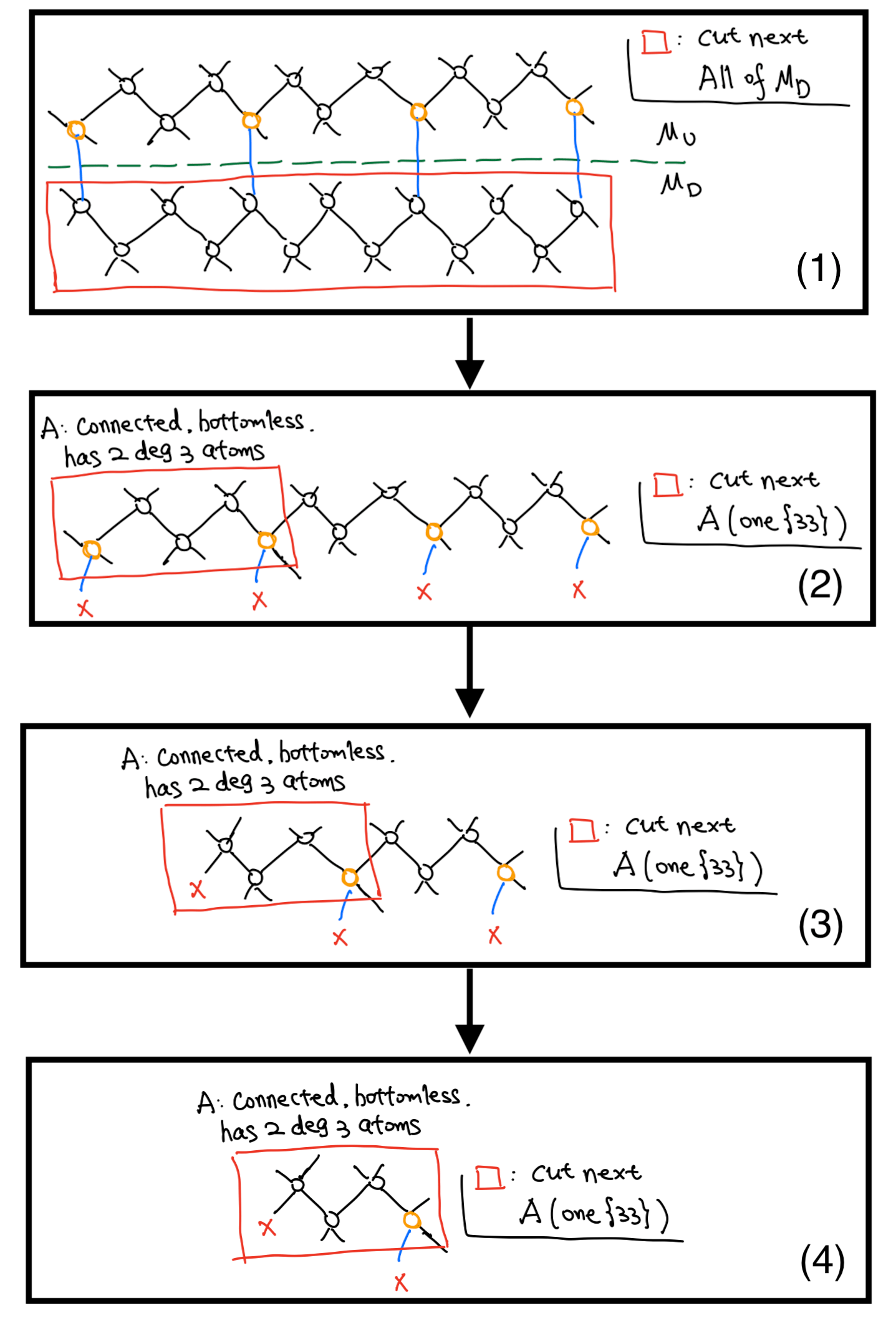}
    \caption{Algorithm for toy model II, applied to the example in {\color{blue}Figure \ref{fig.toy models}}. Here we first cut all of $\Mb_D$ as free. Then, each time we choose a connected, no-bottom subset $A$ that contains exactly two deg 3 atoms, see Definition \ref{def.toy2_alg} (\ref{it.toy2_1}). We then apply \textbf{UP} to cut all atoms in $A$ (cf. {\color{blue}Figure \ref{fig.up}} for the application of \textbf{UP}) and repeat. Note that each application of \textbf{UP} guarantees at least one \{33\} molecule.}
    \label{fig.toy model 2 alg}
\end{figure}
\begin{enumerate}[{(1)}]
\item \label{it.toy2_1} Apply Lemma \ref{lem.toy2_aux} to $\Mb$ and select a connected no-bottom subset $A$. This set is fixed until the end of \ref{it.toy2_5}.
\item \label{it.toy2_2} Choose a lowest deg 3 atom $\nf$ in $A$ that has not been cut. Let $S_\nf$ be the set of descendants of $\nf$. This $\nf$ and $S_\nf$ are fixed until the end of \ref{it.toy2_3}.
\item \label{it.toy2_3} Starting from $\nf$, each time choose a highest atom $\mf$ in $S_\nf$ that has not been cut. If $\mf$ has deg 3 and has a parent $\mf^+$ or child $\mf^-$ that also has deg 3, then cut $\{\mf,\mf^\pm\}$ as free; otherwise cut $\mf$ as free. Repeat until all atoms in $S_\nf$ have been cut.
\item \label{it.toy2_4} If $\Mb$ has deg 2 atom, then cut it; repeat until $\Mb$ has no deg 2 atoms. Then go to \ref{it.toy2_2} and choose the next $\nf$, and so on.
\item \label{it.toy2_5} When all atoms in $A$ have been cut, we subsequently cut all deg 2 atoms as in \ref{it.toy2_4}, and then cut all \emph{components of $\Mb$ with only one deg 3 atom} using \textbf{UP}. 
\item\label{it.toy2_6} Go to \ref{it.toy2_1} and select the next $A$ from one of the remaining components of $\Mb$, and so on.
\end{enumerate}
\end{definition}
\begin{remark}\label{rem.toy2} Let us explain the ideas in Lemma \ref{lem.toy2_aux} and Definition \ref{def.toy2_alg}. Note that steps \ref{it.toy2_2}--\ref{it.toy2_4} in Definition \ref{def.toy2_alg} are exactly the same as \textbf{UP} but applied to a particular set $A$, so we can view \textbf{3COMPUP} as a clever way of carrying out \textbf{UP}. The role of the set $A$ is explained as follows: first it is connected (otherwise we can cut each of its components separately); it should also be no-bottom, because when cutting any atom we should also cut its descendants to guarantee the monotonicity property. Finally we need it to contain two deg 3 atoms so that this sequence of cutting will produce at least one \{33\} molecule, and we need \emph{exactly two} deg 3 atoms so that we do not waste too many components of the deg-3-atom set (which are the key objects leading to the gain).
\end{remark}
\begin{proposition}\label{prop.toy2} Let $\Mb$ be a toy model II. We first cut $\Mb_D$ as free and cut it into elementary molecules by \textbf{DOWN} (i.e. the dual of \textbf{UP} with all directions reversed), then apply \textbf{3COMPUP} (Definition \ref{def.toy2_alg}) to $\Mb_U$. Then all resulting components are elementary, and we have
\begin{equation}\label{eq.toy2_1+}\#_{\{4\}}=1,\qquad \#_{\{33\}}\geq \#_{\mathrm{comp}(X)}/10.
\end{equation}
Since $|X|\gtrsim\rho\gg 1$ and using (\ref{eq.toy2_1}), this clearly implies (\ref{eq.goal}).
\end{proposition}
\begin{proof} Note that steps \ref{it.toy2_2}--\ref{it.toy2_4} in Definition \ref{def.toy2_alg} are exactly the same as \textbf{UP}, so the proof that all molecules are elementary (as well as $\#_{\{4\}}=1$ which comes only from $\Mb_D$) follows in the same way as Proposition \ref{prop.up_toy}. We now focus on the proof that $\#_{\{33\}}\geq  \#_{\mathrm{comp}(X)}/10$.

For any molecule $\Mb$, define $\#_{\mathrm{comp}(\Mb)}$ to be the number of components of $\Mb$, and $\#_{3\mathrm{comp}(\Mb)}$ to be the number of components of the set of deg 3 atoms in $\Mb$, and \begin{equation}\label{eq.toy2_2-}\lambda(\Mb):=\frac{1}{5}\#_{\mathrm{comp}(\Mb)}+ \frac{1}{10}\#_{3\mathrm{comp}(\Mb)}.\end{equation} In Definition \ref{def.toy2_alg} (applied to $\Mb_U$), the whole algorithm consists of finitely many loops where one loop is formed by steps \ref{it.toy2_1}--\ref{it.toy2_5} in Definition \ref{def.toy2_alg}. Consider each such loop which reduces $\Mb$ to $\Mb'$, we shall prove the following:
\begin{itemize}
\item[($\clubsuit$)] If $\Mb$ is a tree, has no top fixed end and contains only deg 3 and 4 atoms with at least two deg 3 atoms, then each component of $\Mb'$ satisfies the same property. Moreover for this loop we have
\begin{equation}\label{eq.toy2_2}\#_{\{33\}}\geq\lambda(\Mb)-\lambda(\Mb').
\end{equation}
\end{itemize} Note that adding up (\ref{eq.toy2_2}) leads to the desired lower bound of $\#_{\{33\}}$ (initially $\Mb_U$ has no deg 2 atoms nor top fixed end, and $\#_{3\mathrm{comp}(\Mb_U)}=\#_{\mathrm{comp}(X)}$), thus it suffices to prove ($\clubsuit$) for a given loop. We divide the proof into 4 parts.

\textbf{Proof part 1.} We prove that each component of $\Mb'$ satisfies the desired properties. This is obvious, as the lack of top fixed ends follows from the same proof as Proposition \ref{prop.up_toy}, and the lack of deg 2 atoms and components with only one deg 3 atom follows from the cuttings done in Definition \ref{def.toy2_alg} \ref{it.toy2_5}.

\textbf{Proof part 2.} Now we need to prove (\ref{eq.toy2_2}). We start by making some preparations. Note that $A$ is a connected subset of $\Mb$; by Lemma \ref{lem.cutconnected0} in the proof of Lemma \ref{lem.toy_1}, each atom $\mf\in \Mb\backslash A$ that is connected to an atom $\nf=\nf(\mf)$ in $A$ by a pre-loop bond $e$ corresponds to a unique component $X_\mf$ after cutting the atoms in $A$ as free, see {\color{blue}Figure \ref{fig.setab}} (the loop will in fact cut as free more atoms). Note also that $\mf$ must be a parent of $\nf$ due to the no-bottom property of $A$.

Let $\Mb_3$ be the set of deg 3 atoms in $\Mb$ (before executing the loop), and consider all components of $\Mb_3$ that intersect $A$ or are connected to $A$ by a bond. We know $A$ has exactly two atoms in $\Mb_3$; apart from the one or two components that contain these two atoms, each component under consideration must be a subset of a unique $X_\mf$ where $\mf$ has deg 3 and $\nf=\nf(\mf)$ has deg 4 (call it $Y_\mf$). Let $B$ be the union of these components with $A$ ({\color{blue}Figure \ref{fig.setab}}). We claim that after steps \ref{it.toy2_1}--\ref{it.toy2_4} in Definition \ref{def.toy2_alg} and cutting the deg 2 atoms in Definition \ref{def.toy2_alg} \ref{it.toy2_5}, exactly those atoms in $B$ have been cut as free.
\begin{figure}[h!]
    \centering
    \includegraphics[width=0.5\linewidth]{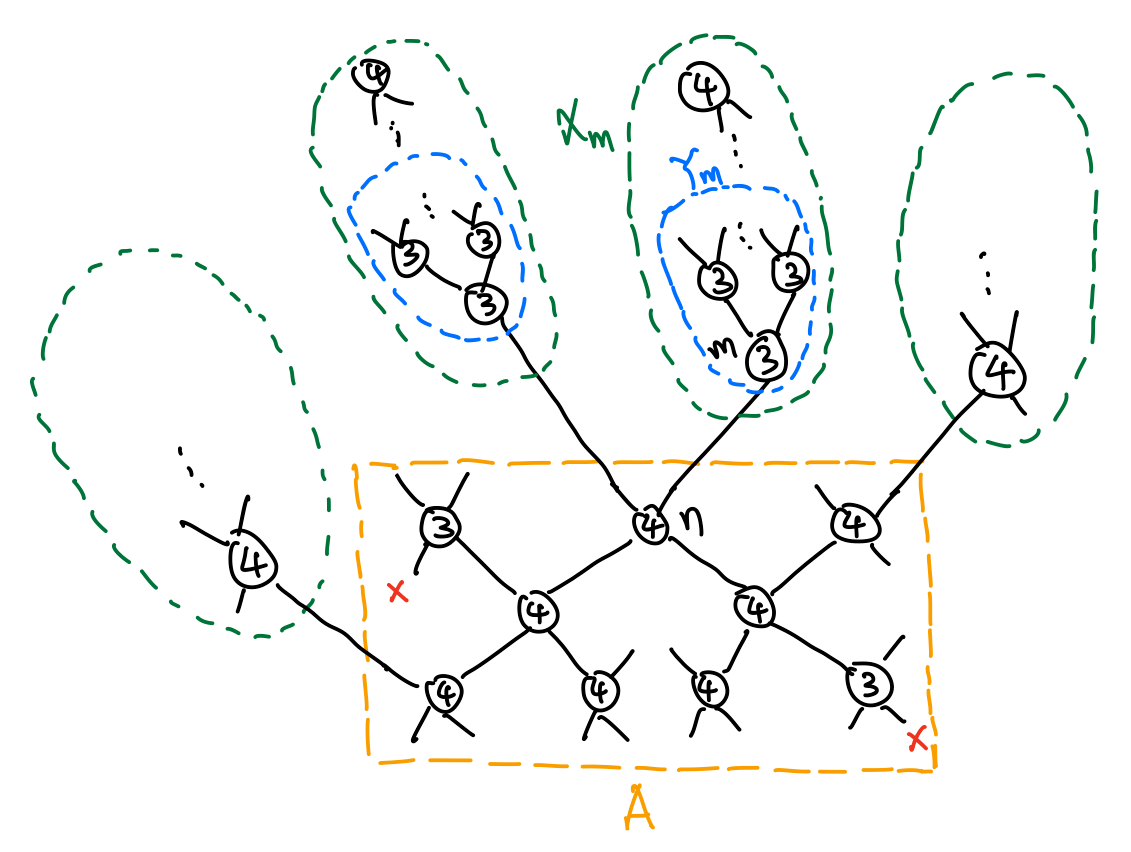}
    \caption{Part 2 of the proof of Proposition \ref{prop.toy2}. Here $A$ is the connected, no-bottom set as in Lemma \ref{lem.toy2_aux}; various atoms $\mf$ (together with the corresponding component $X_\mf$ as in Lemma \ref{lem.cutconnected0} are drawn, each of which connected to its child $\nf\in A$ by a bond. Apart from 2 exceptions, this $\nf$ will have deg 4. If $\mf$ has deg 4 we do nothing; if $\mf$ has deg 3, we take all the cubsequently connected deg 3 atoms in $X_\mf$ and call them $Y_\mf$ (i.e. the blue circle within the green circle $X_\mf$). \textbf{The $B$ is defined as the union of $A$ and all blue circles.} Within one loop \ref{it.toy2_1}--\ref{it.toy2_5} in Definition \ref{def.toy2_alg}, we cut everything in $A$, and subsequently everything in the blue circles as deg 2 (or as part of \{33\} molecule). That is, we cut everything in $B$ and nothing else (as cutting $B$ only turns deg 4 atoms into deg 3 but does not produce further deg 2 atom).}
    \label{fig.setab}
\end{figure}

In fact, the atoms $A$ must have been cut by definition. For any atom in the $\Mb_3$ components forming $B$, by definition it is connected to an atom in $A$ by a chain of deg 3 atoms; once the atom in $A$ has been cut, the atoms on this chain will successively become deg 2, so they must have been cut (either in Definition \ref{def.toy2_alg} \ref{it.toy2_4}--\ref{it.toy2_5}, or as part of \{33\} molecule in Definition \ref{def.toy2_alg} \ref{it.toy2_3}). On the other hand, each atom in $\Mb\backslash B$ has at most one bond connecting to $B$ (as $B$ is connected and $\Mb$ is a tree), and must be deg 4 if it has one such bond (otherwise it would belong to one of the $\Mb_3$ components forming $B$), so each atom in $\Mb\backslash B$ will have deg 3 or 4 after cutting $B$ as free, thus it will not be cut in Definition \ref{def.toy2_alg} \ref{it.toy2_1}--\ref{it.toy2_4} or in cutting the deg 2 atoms in Definition \ref{def.toy2_alg} \ref{it.toy2_5}. 

\textbf{Proof part 3.} Next we obtain a lower bound on $\#_{\{33\}}$. Note that $A$ contains two deg 3 atoms, we can show that $\#_{\{33\}}\geq 1$ in the same way as in Part 3 of the proof of Proposition \ref{prop.up_toy}. Moreover, since $A$ contains \emph{exactly} two deg 3 atoms, by arguing as in Part 4 of the proof of Proposition \ref{prop.up_toy} and calculating the quantity $\sigma$ defined there but within $A$, we can show that, apart from at most two exceptions, each deg 4 atom in $A$ will have deg 3 when it is cut.

Now let $\lambda$ be the number of $\Mb_3$ components forming $B$. With at most two exceptions, each such component corresponds to a deg 3 atom $\mf\in X_{\mf}$ together with its deg 4 child $\nf=\nf(\mf)\in A$ as in Part 2 above. With at most two more exceptions, this $\nf$ will have deg 3 when it is cut, at which time $\mf$ also has deg 3 (see {\color{blue}Figure \ref{fig.cuta33}}); then by Definition \ref{def.toy2_alg} \ref{it.toy2_3}, we know that $\nf$ must belong to a \{33\} molecule (which may or may not be $\{\nf,\mf\}$). Note that each \{33\} molecule can be so obtained by at most one $\nf$, and each $\nf$ can be obtained by at most two $\mf$ (hence at most two components), we conclude that
\begin{equation}\label{eq.toy2_3}\#_{\{33\}}\geq\max\bigg(1,\frac{\lambda-2}{2}-2\bigg).\end{equation}
\begin{figure}[h!]
    \centering
    \includegraphics[width=0.5\linewidth]{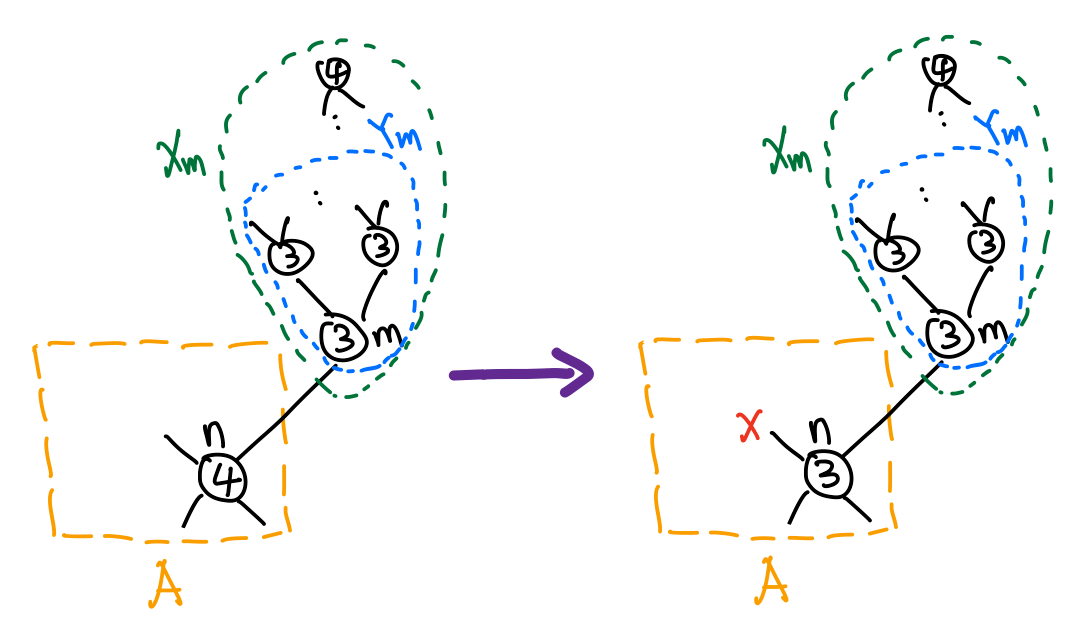}
    \caption{Part 3 of the proof of Proposition \ref{prop.toy2}: here $\nf\in A$ is connected to $\mf\in X_\mf$ which has deg 3, as in {\color{blue}Figure \ref{fig.setab}}. With at most 2 more exceptions, this $\nf$ will also have deg 3 when it is cut (as we are essentially doing \textbf{UP} within $A$, see Proposition \ref{prop.up_toy}), which leads to a \{33\} molecule formed by $\nf$ and $\mf$.}
    \label{fig.cuta33}
\end{figure}

\begin{figure}[h!]
    \centering
    \includegraphics[width=0.5\linewidth]{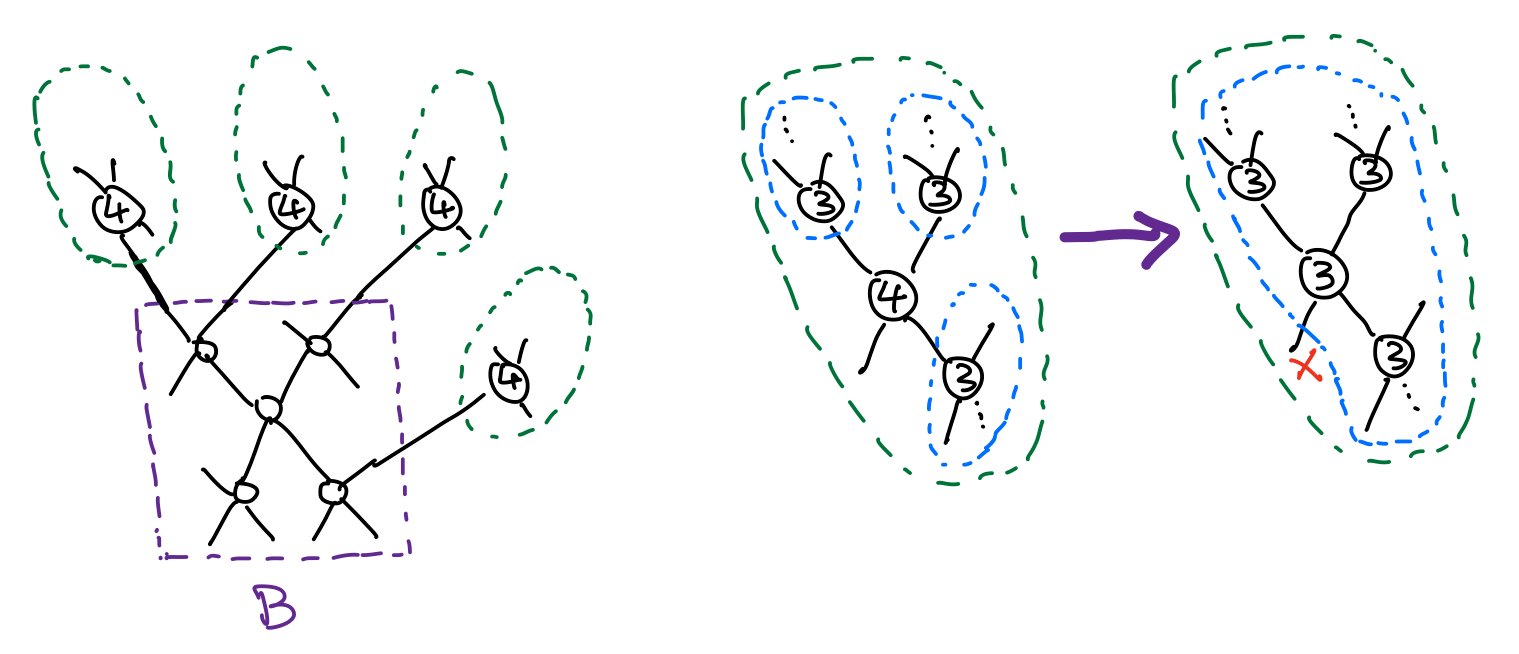}
    \caption{Part 4 of the proof of Proposition \ref{prop.toy2}. Left: suppose we cut $B$ and are left with all the components as in Lemma \ref{lem.cutconnected0}. Within each such component, a fixed end is created, turning a deg 4 atom into deg 3. Right: for each component, turning one deg 4 atom into deg 3 may decrease the number of \textbf{components of the set of deg 3 atoms} by 2 (by joining 3 such components into 1 with the new deg 3 atom), but \textbf{no more}. Note that the blue circles here are unrelated to the blue circles in {\color{blue}Figures \ref{fig.setab}--\ref{fig.cuta33}} and represent different things.}
    \label{fig.cutb}
\end{figure}

\textbf{Proof part 4.} Now we consider the components of $\Mb'$ and $\Mb_3'$ (set of deg 3 atoms in $\Mb'$). Note that $\Mb'$ is formed from $\Mb$ by cutting $B$ as free and then removing all components with only one deg 3 atom. Since $B$ is connected, by applying Lemma \ref{lem.cutconnected0} in the proof of Lemma \ref{lem.toy_1} again, we see that each component $X_\qf$ of $\Mb'$ corresponds to a unique atom $\qf\in \Mb\backslash B$ that is adjacent to an atom in $B$, see {\color{blue}Figure \ref{fig.cutb}}. Let $\mu$ be the number of such bonds (i.e. corresponding to those new components with at least two deg 3 atoms after cutting $B$), then we have \begin{equation}\label{eq.toy2_4}\#_{\mathrm{comp}(\Mb)}-\#_{\mathrm{comp}(\Mb')}=1-\mu.\end{equation}

Finally, consider the set of deg 3 atoms, i.e. $\Mb_3$ and $\Mb_3'$, and their components. The components of $\Mb_3$ includes those $\lambda$ components forming $B$ as in Part 2 (which have been cut when forming $\Mb'$), and the other components which are contained in some $X_\qf$ defined above. Note that for each new component $Y$ formed after cutting $B$ that contains only one deg 3 atom, it will be removed and will not occur in $\Mb'$; however, all atoms in $Y$ originally have deg 4 before any cutting, so initially $Y$ also has no contribution to $\#_{3\mathrm{comp}(\Mb)}$, and thus we can ignore such $Y$.

By Lemma \ref{lem.cutconnected0}, when viewed in each $X_\qf$, the effect of cutting $B$ exactly turns one free end into fixed end (turning one deg 4 atom into deg 3). Therefore the part of $\Mb_3'$ within $X_\qf$ is formed by \emph{adding one deg 3 atom} into the existing part of $\Mb_3$ within $X_\qf$. This new deg 3 atom may \emph{decrease} $\#_{3\mathrm{comp}(X_\qf)}$ (i.e. the number of components of the set of deg 3 atoms in $X_\qf$) by joining existing components together, but this decrease is \emph{at most 2} because this new deg 3 atom can join at most 3 components into 1. From this, and accounting for the $\lambda$ components removed when cutting $B$, we deduce that
\begin{equation}\label{eq.toy2_5}\#_{3\mathrm{comp}(\Mb)}-\#_{3\mathrm{comp}(\Mb')}\leq \lambda+2\mu.\end{equation}
Putting together (\ref{eq.toy2_2-}) and (\ref{eq.toy2_3})--(\ref{eq.toy2_5}), we deduce that
\begin{equation}\label{eq.toy2_6}
\begin{aligned}\lambda(\Mb)-\lambda(\Mb')&=\frac{1}{5}\big(\#_{\mathrm{comp}(\Mb)}-\#_{\mathrm{comp}(\Mb')}\big)+\frac{1}{10}\big(\#_{3\mathrm{comp}(\Mb)}-\#_{3\mathrm{comp}(\Mb')}\big)
\\&\leq \frac{1}{5}+\frac{\lambda}{10}\leq\max\bigg(1,\frac{\lambda-6}{2}\bigg)\leq \#_{\{33\}}.
\end{aligned}\end{equation} This proves (\ref{eq.toy2_2}) and completes the proof of Proposition \ref{prop.toy2}.
\end{proof}

\subsection{Toy model III}\label{sec.toy3} 
\begin{figure}[h!]
    \centering
    \includegraphics[width=0.7\linewidth]{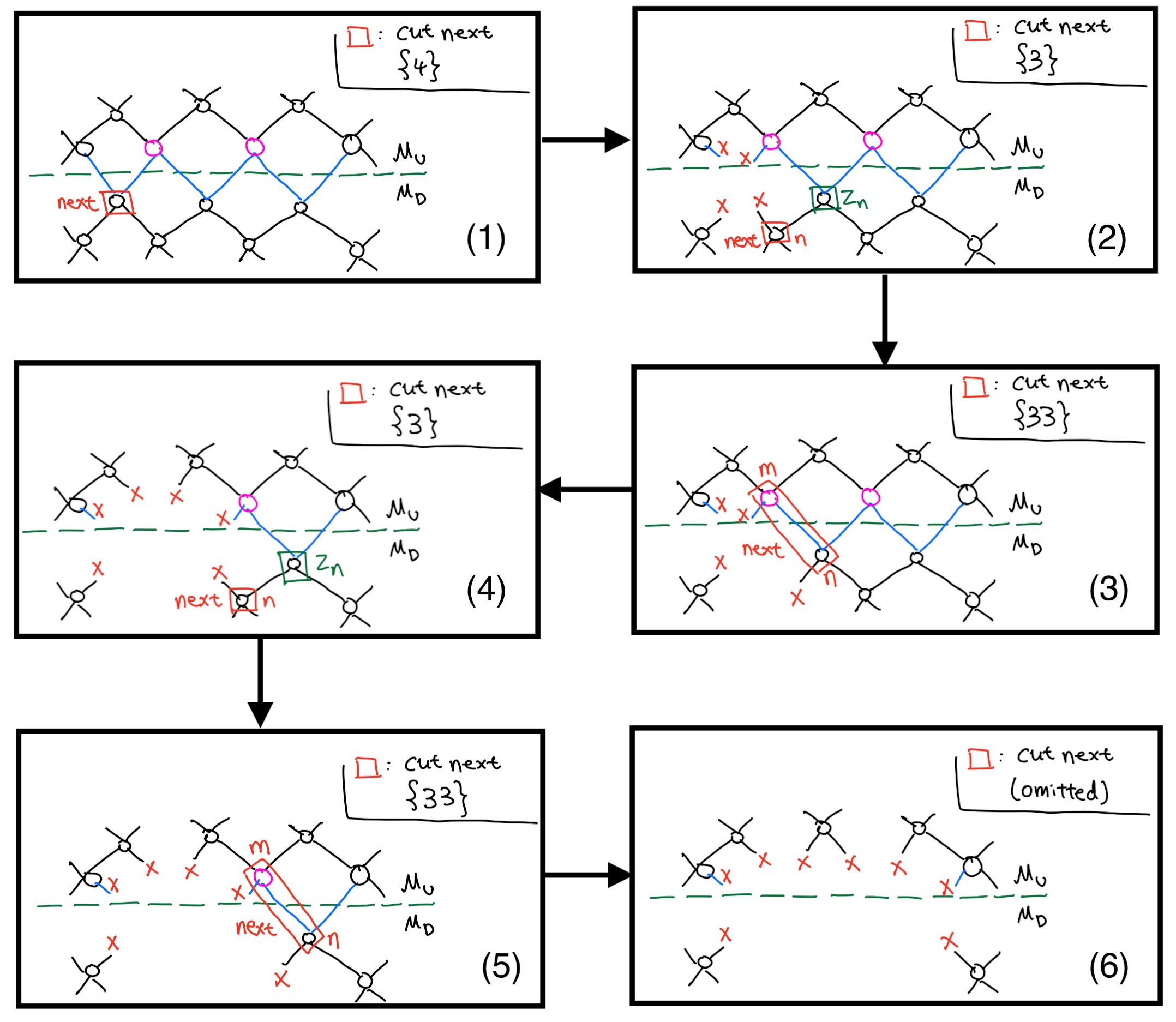}
    \caption{Algorithm for toy model III, applied to the example in {\color{blue}Figure \ref{fig.toy models}}. Here we perform \textbf{DOWN} which is the dual of \textbf{UP} in $\Mb_D$, so in Definition \ref{def.up_toy} \ref{it.up_toy_2} in \textbf{UP} algorithm we choose the \emph{highest} deg 3 atom instead of the lowest, and in in Definition \ref{def.up_toy} \ref{it.up_toy_3} in \textbf{UP} algorithm we choose the \emph{lowest} remaining atom instead of the highest, etc. Here green box indicates atoms in $Z_\nf$ (the set of \emph{ancestors} of $\nf$), and when an atom $\nf\in\Mb_D$ forms a \{33\} molecule with an atom $\mf\in\Mb_U$ we always cut it. In the end we see that each 2-connection (pink atom) belongs to one \{33\} molecule, as desired. As the gain is already enough, the rest can be cut into elementary molecules using \textbf{UP} in $\Mb_U$ and \textbf{DOWN} in $\Mb_D$.}
    \label{fig.toy model 3 alg}
\end{figure}

\begin{figure}[h!]
    \centering
    \includegraphics[width=0.5\linewidth]{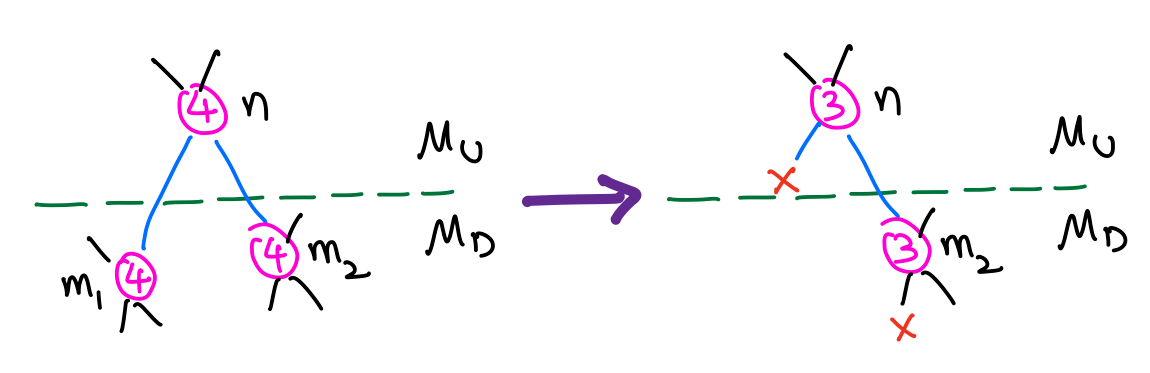}
    \caption{The proof of Proposition \ref{prop.toy3}. Here $\nf\in\Mb_U$ is a 2-connection connected to $\mf_{1,2}\in \Mb_D$; assume $\mf_1$ is cut before $\mf_2$ when applying \textbf{DOWN} in $\Mb_D$. Then at the time $\mf_2$ is cut, this $\mf_2$ has deg 3 due to properties of \textbf{DOWN} (Proposition \ref{prop.up_toy}), and $\nf$ also has deg 3 because $\mf_1$ has been cut.}
    \label{fig.2conn}
\end{figure}

In this subsection we introduce the toy model III. Note that in both toy models I plus and II, we are assuming the absence of 2-connections, i.e. $\#_{\mathrm{2conn}}=0$. Now we define the toy models III, which are (essentially) opposite to toy models I plus and II via the \emph{dichotomy}, where we assume\footnote{There is also the case $0<\#_{\mathrm{2conn}}\ll |X|$; in this case we will first cut all such 2-connections, which generates an error which is eventually negligible but still requires some extra care; see Section \ref{sec.reduce5} for more discussions.} $\#_{\mathrm{2conn}}\gtrsim |X|$ instead of $\#_{\mathrm{2conn}}=0$.
\begin{definition}\label{def.toy3} We say $\Mb$ is a toy model III, if it is a simplified 2-layer model in Definition \ref{def.simplified}, and
\begin{equation}\label{eq.toy3_1}\#_{\mathrm{2conn}}\gtrsim |X|,
\end{equation} where $\#_{\mathrm{2conn}}$ and $X$ are defined at the beginning of this note.
\end{definition}
The strategy for toy models III is simple: we work in $\Mb_D$ as if $\Mb_U$ were absent, and apply the \textbf{DOWN} algorithm (the dual of \textbf{UP}) to it. Now consider any 2-connection $\nf\in \Mb_U$, which will become deg 2 when all atoms in $\Mb_D$ have been cut. Assume $\nf$ is adjacent to two atoms $\mf_1,\mf_2\in\Mb_D$, and assume $\mf_1$ is cut before $\mf_2$ in the \textbf{DOWN} algorithm applied to $\Mb_D$. Then, right before $\mf_2$ is cut, this $\mf_2$ will have deg 3 (say) by Proposition \ref{prop.up_toy}; also since $\mf_1$ has been cut at this time, now $\nf$ must also have deg 3, which leads to a \{33\} molecule, so we simply cut it and enjoy the gain. This way we guarantee that each 2-connection belongs to a \{33\} molecule, which is sufficient under the assumption (\ref{eq.toy3_1}). We shall refer to this as the \textbf{2CONNDN} algorithm (and the dual version \textbf{2CONNUP}), see Definition \ref{def.alg_2connup} for the full version (which is more complicated due to the various necessary additions including the O-atoms).
\begin{definition}[Algorithm for toy model III]\label{def.toy3_alg} Let $\Mb$ be toy model III. We define the following cutting sequence. For an example see {\color{blue}Figure \ref{fig.toy model 3 alg}}:
\begin{enumerate}[{(1)}]
\item\label{it.toy3_1} If $\Mb_D$ has a deg 2 atom, then cut it as free. Repeat until $\Mb_D$ does not have deg 2 atoms.
\item\label{it.toy3_2} Choose a highest deg 3 atom $\nf$ in $\Mb_D$ that has not bee cut (or a highest atom if $\Mb_D$ has no deg 3 atoms). Let $Z_\nf$ be the set of ancestors of $\nf$ in $\Mb_D$. This $\nf$ and $Z_\nf$ is fixed until all atoms in $Z_\nf$ have been cut.
\item\label{it.toy3_3} Starting from $\nf$, each time choose a lowest atom $\mf$ in $Z_\nf$ that has not been cut. If $\mf$ has deg 3 and has a parent $\pf\in \Mb_U$ that also has deg 3, then cut $\{\mf,\pf\}$ as free; otherwise cut $\mf$ as free. Repeat until all atoms in $Z_\nf$ have been cut, then go to \ref{it.toy3_1}--\ref{it.toy3_2} and choose the next $\nf$, and so on.
\item \label{it.toy3_4} When all atoms in $\Mb_D$ have been cut, finish off $\Mb_U$ by applying the \textbf{UP} algorithm.
\end{enumerate}
\end{definition}
Note that steps \ref{it.toy3_1}--\ref{it.toy3_3} in Definition \ref{def.toy3_alg} are essentially just \textbf{DOWN} (i.e. the dual of \textbf{UP}) applied to $\Mb_D$, except we propritize the \{33\} molecules formed by an atom in $\Mb_D$ with an atom in $\Mb_U$.
\begin{proposition}\label{prop.toy3} Let $\Mb$ be toy model III. Then, after applying the cutting sequence defined in Definition \ref{def.toy3_alg}, all molecules we get are elementary, and we have
\begin{equation}\label{eq.toy3_2}\#_{\{4\}}=1,\qquad \#_{\{33\}}\geq \frac{1}{2}\cdot\#_{\mathrm{2conn}}.
\end{equation}
\end{proposition}
\begin{proof} Since steps \ref{it.toy3_1}--\ref{it.toy3_3} in Definition \ref{def.toy3_alg} are essentially the same as in \textbf{DOWN}, it is easy to prove all molecules we get are elementary and $\#_{\{4\}}=1$, in the same way as the proof of Proposition \ref{prop.up_toy}. Now we only need to prove that $\#_{\{33\}}\geq \#_{\mathrm{2conn}}/2$.

Let $\nf\in\Mb_U$ be any 2-connection, which is connected to two atoms $\mf_1,\mf_2\in\Mb_D$ by two bonds. Assume $\mf_1$ is cut before $\mf_2$ in Definition \ref{def.toy3_alg} (they cannot be cut at the same time, as Definition \ref{def.toy3_alg} \ref{it.toy3_3} does not provide an option to cut a \{33\} molecule with both atoms in $\Mb_D$). Arguing in the same way as in Part 2 of the proof of Proposition \ref{prop.up_toy}, we know that $\mf_2$ cannot have deg 4 when it is cut. On the other hand, note that $\Mb_D$ is a tree initially without deg 3 atoms; by arguing in the same way as in Part 4 of the proof of Proposition \ref{prop.up_toy} and calculating the quantity $\sigma$ defined there but within $\Mb_D$, we conclude that $\mf_2$ also cannot have deg 2 when it is cut.

Therefore, $\mf_2$ must have deg 3 when it is cut; at this time $\nf$ must also have deg 3 because $\mf_1$ has been cut, $\mf_2$ has not, and $\nf$ also has two top edges (not fixed ends) that are not affected, see {\color{blue}Figure \ref{fig.2conn}}. By Definition \ref{def.toy3_alg} \ref{it.toy3_3} we conclude that $\mf_2$ must belong to a \{33\} molecule, which means that $\nf$ must be adjacent with an $\Mb_D$ atom in \emph{some} \{33\} molecule. Since this is true for each 2-connection, and each \{33\} molecule may correspond to at most two 2-connections in this way (i.e. given the \{33\} molecule there can be at most two $\nf$ that are adjacent with the $\Mb_D$ atom in this \{33\} molecule), we conclude that $\#_{\{33\}}\geq \#_{\mathrm{2conn}}/2$.
\end{proof}
\subsection{Reduction to 2-layers: layer selection} \label{sec.toy_multi} 
\begin{figure}[h!]
    \centering
    \includegraphics[width=0.6\linewidth]{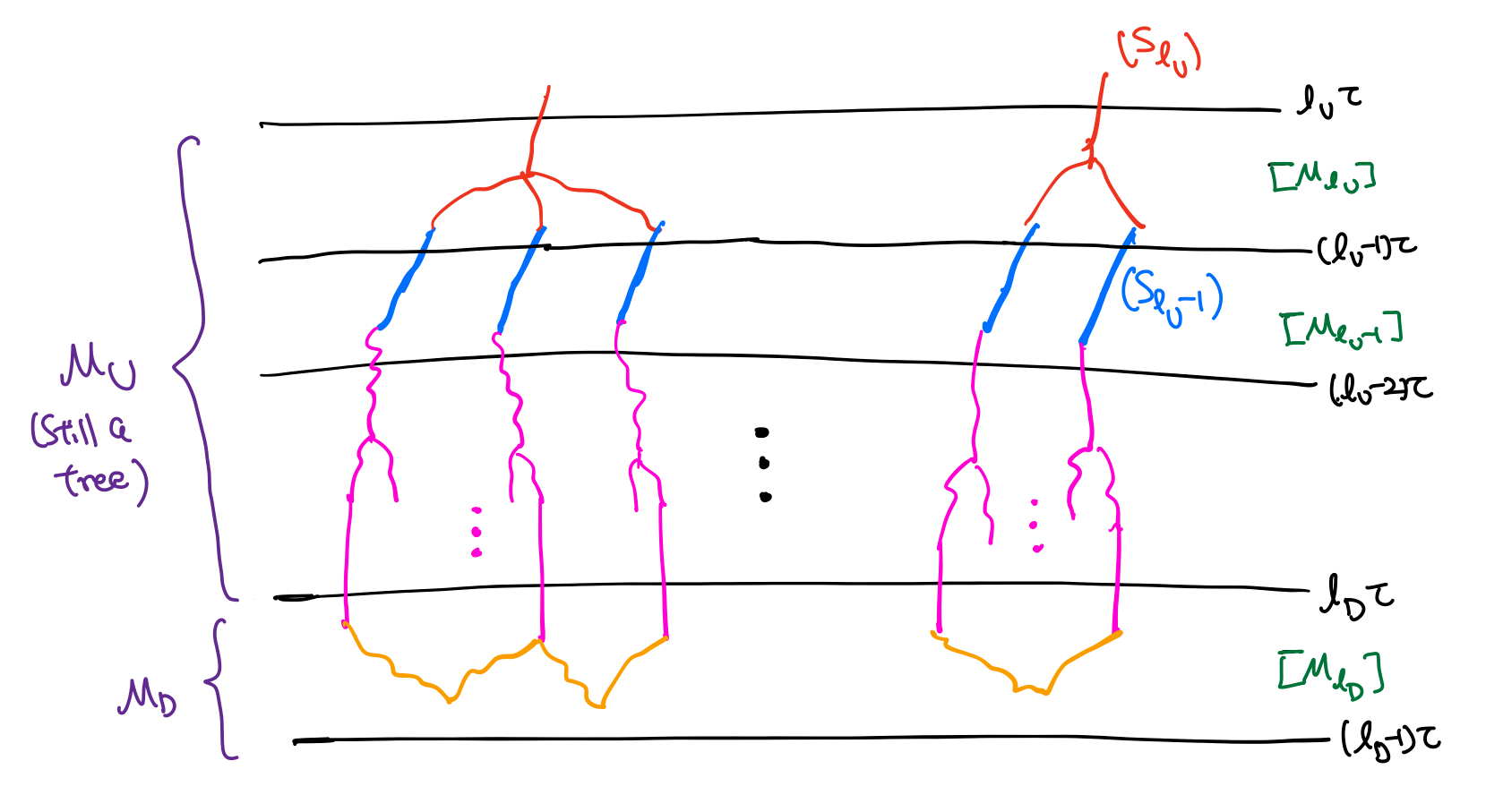}
    \caption{Layer selection (Definition \ref{def.layer_select_toy}). The blue edges in layer $\Mb_{\ell_U-1}$ belong to particle lines in $r(\Mb_{\ell_U-1})$, whose size satisfies  $s_{{\ell_U}-1}\gtrsim\rho$. Each of them is connected to a red edge crossing $\ell_U\tau$ (they belong to the particle lines in $r(\Mb_{\ell_U})$, and their number $s_{\ell_U}\ll s_{\ell_U-1}$, i.e. these particle lines are essentially connected to each other (basically, this condition is used to guarantee that all the bonds between $\Mb_U$ and $\Mb_D$ do lead to cycles). in layer $\Mb_{\ell_U}$. 
    \newline
\-\hspace{8pt}    On the other hand, these particle lines are connected to pink edges below, but are \textbf{not connected to each other} until reaching layer $\Mb_{\ell_D}$ (by yellow edges in $\Mb_{\ell_D}$). This corresponds to the case $\lambda_{{\ell_U};\ell_D+1}=0$ and $\lambda_{\ell_U;\ell_D}\gtrsim s_{\ell_U}$ (the general case $0<\lambda_{{\ell_U};\ell_D+1}\ll s_{\ell_U}$ requires a little more work). We then set $\Mb_U$ to contain layers $\Mb_{[\ell_D+1:\ell_U]}$ (which is a forest in our case; in general case we need to cut some exceptional components to make it a forest), and $\Mb_D$ to contain layer $\Mb_{\ell_D}$, to reduce to the 2-layer case, where each of $\Mb_U$ and $\Mb_D$ contains no cycle, and the number of bonds between $\Mb_U$ and $\Mb_D$ being $\gtrsim\rho$.}
    \label{fig.layerselect}
\end{figure}

In Sections \ref{sec.toy1}--\ref{sec.toy3} we have considered the 2-layer models in Definition \ref{def.simplified}. In fact, this two-layer case already captures most of the essential difficulties in the multi-layer case, and a simple process allows us to reduce the general $\Lf$-layer case (where $\Lf$ is any large constant) to the 2-layer case. This is refered to as \emph{layer selection}, which we define as follows; for the precise definition see Definition \ref{def.layer_select}.
\begin{definition}\label{def.layer_select_toy} Let $\Mb$ satisfy Proposition \ref{prop.mol_axiom} and simplifications (\ref{it.simp_1})--(\ref{it.simp_2}) in Section \ref{sec.toy_intro}. Recall the set $r(\Mb_\ell)$ of root particles at layer $\ell$ (Definition \ref{def.prsets}), and $s_\ell=|r(\Mb_\ell)|$; they are the particle lines that \emph{cross} into layer $\ell$. Assume also $|H|=s_\ell\ll\rho$ (the case when $|H|\gtrsim\rho$ is easy due to the positive term $|H|$ in (\ref{eq.overall_alg}); see Proposition \ref{prop.case1}). We define the \textbf{layer selection} process as follows, see {\color{blue}Figure \ref{fig.layerselect}}:
\begin{enumerate}[{(1)}]
\item \label{it.select_toy1}Choose $\ell_U$ such that \begin{equation}\label{eq.select_toy}\rho':=s_{\ell_U-1}\sim\max_{\ell'}s_{\ell'}\sim\rho,\qquad \rho'\gg\max_{\ell'\geq\ell_U}s_{\ell'};\end{equation} here the first condition in (\ref{eq.select_toy}) states that the majority of crossings happen between layers $\ell_U$ and $\ell_U-1$, and the second condition states that the crossings in layers above $\ell_U$ are negligible.
\item\label{it.select_toy2} Consider each of the $s_{\ell_U-1}$ particle lines in $r(\Mb_{\ell_U-1})$. For $\ell_D<\ell_U$, define $v_{\ell_D}^*$ to be the number of particle lines in $r(\Mb_{\ell_U-1})$ that are connected to either a particle line in in $r(\Mb_{\ell_U-1})$ or to a cycle \emph{within $\Mb_{[\ell_D:\ell_U)}$}; in other words, we start from layer $\ell_U-1$ and move down, keeping track of the number of crossings that are connected to other crossings or to a cycle \emph{at each layer}. Note that $v_1^*=s_{\ell_U-1}$ by Proposition \ref{prop.mol_axiom} (\ref{it.axiom6}) (and $H_0=\varnothing$); now as in \ref{it.select_toy1}, we choose $\ell_D$ such that
\begin{equation}\label{eq.select_toy2}\rho'':=v_{\ell_D}^*\sim s_{\ell_U-1}=\rho',\qquad \rho''\gg\max_{\ell''>\ell_D}v_{\ell''}^*.\end{equation}
\item\label{it.select_toy3} Now we focus on the layers $\Mb_{[\ell_D:\ell_U]}:=\cup_{\ell_D\leq\ell'\leq\ell_U}\Mb_{\ell'}$, and ignore all other layers. Consider each component $Z$ of the set $\Mb_{(\ell_D:\ell_U)}$. If $Z$ \emph{intersects at least two particle lines} in $r(\Mb_{\ell_U-1})$ or \emph{contains a cycle}, then by the choice of $\ell_D$ in \ref{it.select_toy2}, the number of \emph{such exceptional components} is $\leq v_{\ell_D+1}^*\ll\rho''$ which is negligible. We will assume this number is zero (otherwise we need an additional cutting and \emph{pre-processing} step, see Section \ref{sec.reduce5} below).  
\item \label{it.select_toy4} Define the result of layer selection, to be $\Mb_{UD}:=\Mb_U\cup\Mb_D$, where $\Mb_U:=\Mb_{(\ell_D:\ell_U]}$ and $\Mb_D:=\Mb_{\ell_D}$.
\end{enumerate}
\end{definition}
\begin{proposition}\label{prop.layer_select_toy} Assume $\Mb$ satisfies Proposition \ref{prop.mol_axiom} and simplifications (\ref{it.simp_1})--(\ref{it.simp_2}) in Section \ref{sec.toy_intro}, and $|H|\ll\rho$, and assume that the \emph{number of exceptional components is zero} in Definition \ref{def.layer_select_toy} \ref{it.select_toy3}. Then the process in Definition \ref{def.layer_select_toy} is well-defined, and for $\Mb_{UD}=\Mb_U\cup\Mb_D$ in Definition \ref{def.layer_select_toy} \ref{it.select_toy4}, we have that
\begin{itemize}
\item Each of $\Mb_U$ and $\Mb_D$ is a \emph{forest} with no fixed end, the number of components of $\Mb_U$ being $\ll\rho''$, and no atom in $\Mb_D$ is a parent of an atom in $\Mb_U$;
\item There exist at least $\rho''$ bonds between $\Mb_U$ and $\Mb_D$, and each such bond is connected to another bond between $\Mb_U$ and $\Mb_D$, by a path within $\Mb_D$.
\end{itemize}
\end{proposition}
\begin{proof} To prove that the procedure in Definition \ref{def.layer_select_toy} is well-defined, we just need to prove the existence of $\ell_U$ in Definition \ref{def.layer_select_toy} \ref{it.select_toy1} and $\ell_D$ in Definition \ref{def.layer_select_toy} \ref{it.select_toy2}. In fact, we can choose two large constants (see Definition \ref{def.layer_select} and Proposition \ref{prop.layer_select} for more details) $C\gg (C')^\ell$ such that
\begin{equation}\label{eq.select_toy0}\ell_U-1:=\max\big\{\ell': s_{\ell'}\geq C^{-\ell'}\cdot\rho\big\},\qquad\ell_D:=\max\big\{\ell''<\ell_U:v_{\ell''}^*\geq (C')^{-\ell''}\cdot\rho'\big\}.\end{equation} Note that $\rho\sim\max_{\ell'}s_{\ell'}$ by (\ref{eq.def_rho_old}) (as $\Rf_{\ell'}=0$ by assumption) and $s_{\ell}=|H|\ll C^{-\ell}\cdot\rho$ by assumption, so $\ell_U$ exists; similarly $\ell_D$ exists because $v_1^*=\rho'$ by Proposition \ref{prop.mol_axiom} (\ref{it.axiom6}).

Now, we have assumed $v_{\ell_D+1}^*=0$ (i.e. no exceptional component) in Definition \ref{def.layer_select_toy} \ref{it.select_toy3}. It is obvious that $\Mb_D$ is a forest (by simplification (\ref{it.simp_2}) in Section \ref{sec.toy_intro}), $\Mb_U$ and $\Mb_D$ has no fixed end (because nothing is cut in $\Mb_{UD}=\Mb_U\cup\Mb_D$), and no atom in $\Mb_D$ is parent of an atom in $\Mb_U$ (because the layers forming $\Mb_U$ are above that forming $\Mb_D$).

We next show that $\Mb_U$ is a forest. In fact, each component of $\Mb_{(\ell_D:\ell_U)}$ must intersect a particle line in $r(\Mb_{\ell_U-1})$ (it must intersect a particle line in $p(\Mb_{\ell'})$ for some $\ell'<\ell_U$ by Proposition \ref{prop.mol_axiom} (\ref{it.axiom4}); by applying Proposition \ref{prop.mol_axiom} (\ref{it.axiom5}) we know that this component must also intersect a particle line in $r(\Mb_{\ell'})\subseteq p(\Mb_{\ell'+1})$, by repeating this and using Remark \ref{rem.layer_interval}, we eventually get to a particle line in $r(\Mb_{\ell_U-1})$). As there is no exceptional component in Definition \ref{def.layer_select_toy} \ref{it.select_toy3}, we know that each component of $\Mb_{(\ell_D:\ell_U)}$ is a tree by itself, and intersects a unique particle line in $r(\Mb_{\ell_U})$. Since $\Mb_{\ell_U}$ is a forest by itself, we see that $\Mb_U=\Mb_{(\ell_D:\ell_U]}$ is also a forest (namely, if we start with a forest, and attach a number of trees to it with each tree connected to the original forest by a unique edge, then the result is still a forest).

We next show that the number of components of $\Mb_D$ does not exceed $s_{\ell_U}$, which is $\ll\rho''$ because
\[s_{\ell_U}\leq C^{-1}\cdot s_{\ell_U-1}= C^{-1}\cdot\rho'\leq C^{-1}\cdot (C')^{\ell}\cdot\rho''\ll\rho''\] by (\ref{eq.select_toy0}) and $C\gg (C')^\ell$. In fact, by the same arguments for $\Mb_{(\ell_D:\ell_U)}$ as above, we see that each component of $\Mb_{(\ell_D:\ell_U]}$ must intersect a particle line in $r(\Mb_{\ell_U})$, so the number of these components is bounded by $|r(\Mb_{\ell_U})|=s_{\ell_U}$ by Remark \ref{rem.layer_interval}.

Finally, consider each of the $\rho''=v_{\ell_D}^*$ particle lines that is connected to either a particle line in in $r(\Mb_{\ell_U-1})$ or to a cycle within $\Mb_{[\ell_D:\ell_U)}$, say $\pb$. By definition there exists a path going from an atom on $\pb$ to an atom on either a cycle or another particle line in $r(\Mb_{\ell_U-1})$. This path \emph{must enter} $\Mb_{\ell_D}$, because the above-stated property of $\pb$ will not hold without layer $\ell_D$; moreover it also \emph{must exist} $\Mb_{\ell_D}$ (including exiting $\Mb_D$ when getting onto another particle line in $r(\Mb_{\ell_U-1})$), because $\Mb_D$ by itself does not contain a cycle. By considering the \emph{first entry} and \emph{first exit} at $\Mb_D$, we obtain two bonds between $\Mb_U$ and $\Mb_D$ that are connected within $\Mb_D$; moreover the first entry bond must be connected to $\pb$ via $\Mb_{(\ell_D:\ell_U)}$, and is thus in bijection with $\pb$. This proves the desired result, namely the existence of at least $\rho''$ bonds between $\Mb_U$ and $\Mb_D$, each of which being connected to another by a path within $\Mb_D$.
\end{proof}
\begin{remark}\label{rem.select_toy1} We make two remarks concerning Definition \ref{def.layer_select_toy} and Proposition \ref{prop.layer_select_toy}.
\begin{enumerate}
\item In Definition \ref{def.layer_select_toy}, we focus on $\Mb_{UD}=\Mb_U\cup\Mb_D=\Mb_{[\ell_D:\ell_U]}$ and ignore other layers. In practice, we need to \emph{cut $\Mb_{[\ell_D:\ell_U]}$ as free} from $\Mb$ and cut the other layers into elementary components (first $\Mb_{>\ell_U}$ using \textbf{UP}, then $\Mb_{<\ell_D}$ using \textbf{DOWN}). This may create some \{4\} molecules, but their number is easily bounded using Proposition \ref{prop.up_toy} and the upper bound on the number of components of $\Mb_{>\ell_U}$ and $\Mb$ (see Definition \ref{def.layer_cutting} and Proposition \ref{prop.layer_cutting}).

\item Compared to Definition \ref{def.simplified}, the $\Mb_U$ and $\Mb_D$ in Proposition \ref{prop.layer_select_toy} are \emph{forests} instead of \emph{trees}. However, the number of components of $\Mb_U$ is $\ll\rho''$ by Proposition \ref{prop.layer_select_toy} and is thus negligible; on the other hand, $\Mb_D$ may have many components, but this is compensated by the fact that the $\rho''$ bonds between $\Mb_U$ and $\Mb_D$ in Proposition \ref{prop.layer_select_toy} are connected within $\Mb_D$. Indeed, in the case when $\Mb_D$ has many components, this (in fact a related property of $\Mb_D$, see Proposition \ref{prop.layer_cutting} (3b)) allows us to first cut $\Mb_U$ as free and gain from subsequently cutting $\Mb_D$, see Proposition \ref{prop.comb_est_case6}.
\end{enumerate}
\end{remark} 
\begin{remark}\label{rem.select_toy2} In practice, there are two complications to layer selection: first, if some single layer $\Mb_{\ell'}$ contains cycles (without simplification (\ref{it.simp_2}) in Section \ref{sec.toy_intro}), then we need an extra \emph{layer refining} process (see Definition \ref{def.layer_select}), which is discussed in Section \ref{sec.reduce3}. Second, if $v_{\ell_D+1}^*\neq 0$ (i.e. there are exceptional components in Definition \ref{def.layer_select_toy} \ref{it.select_toy3}, though their number is still $\ll\rho''$ and negligible), we need to first cut them as free before defining $\Mb_U$ and $\Mb_D$ (Definition \ref{def.layer_cutting}). This will create some \{4\} molecules but their number is bounded (see Propositions \ref{prop.layer_select} and \ref{prop.layer_cutting}), and will also introduce fixed ends to $\Mb_U$ and $\Mb_D$ which necessitates the \emph{pre-processing} step, see discussions in Section \ref{sec.reduce5}.
\end{remark} 
\subsection{Additions in the general case}\label{sec.toy_reduce} In this subsection, we discuss the additional ``accessories" needed to address the general case, without the simplifications (\ref{it.simp_1})--(\ref{it.simp_2}) in Section \ref{sec.toy_intro}. These should be viewed as branches that are attached to the main thread of proof, which are less essential but necessary. These accessories are listed in {\color{blue}Table \ref{tab.reduce}}, and will be discussed in Sections \ref{sec.reduce1}--\ref{sec.reduce5} below.
\subsubsection{O-atoms and ov-segments}\label{sec.reduce1} First, in the real problem, the molecule $\Mb$ may contain O-atoms, and the cutting operation in this case is more complicated (see Definition \ref{def.cutting} (\ref{it.cutting_2}) compared to Definition \ref{def.cutting} (\ref{it.cutting_1})). 

As we already see in Definition \ref{def.cutting}, the notion of ov-segments is important in the analysis of O-atoms, and should be viewed as a generalization of edges/bonds. In Definition \ref{def.ov_connect}, we have introduced the notions of being ov-adjacent, ov-parent and ov-child, and ov-components etc., which are direct generalizations of adjacent, parent and child, and connected components. Fortunately, in almost all cases, this generalization \emph{does not cause any essential change} to the proof; we simply translate every sentence in the proofs in Sections \ref{sec.toy1}--\ref{sec.toy3} into the language of ov-segments, with ``adjacent" replaced by ``ov-adjacent" and ``components" replaced by ``ov-components", etc. As an example, by replacing ``connected" with ``ov-connected" and ``bond" with ``ov-segment" in Lemma \ref{lem.cutconnected0}, we obtain the following statement, which is just Lemma \ref{lem.cutconnected}:
\begin{itemize}
\item Suppose we cut $S\subseteq\Mb$ as free, where $S$ is an \emph{ov-connected} subset of $\Mb$ and $\Mb$ is a tree. Then, for each atom $\qf\in \Mb\backslash S$ that is connected to an atom in $S$ by a (pre-cutting) \emph{ov-segment $e$}, the post-cutting molecule $\Mb_D$ has exactly one component $X_\qf$ corresponding to $\qf$. Moreover, when viewed in $X_\qf$, the operation of cutting $S$ turns exactly one free end at $\qf$ into a fixed end.
\end{itemize}
\subsubsection{Reducing each layer to a tree: layer refining}\label{sec.reduce3} In both the toy models in Sections \ref{sec.toy1}--\ref{sec.toy3} and layer selection in Section \ref{sec.toy_multi}, it is crucial that each layer $\Mb_{\ell'}$ of $\Mb$ is a forest. In reality, to deal with the case where some $\Mb_{\ell'}$ contains cycles, we need to define the \textbf{layer refining} process as follows.

Recall (Proposition \ref{prop.mol_axiom} (\ref{it.axiom2})--(\ref{it.axiom3})) that each $\Mb_{\ell'}$ is formed by joining finitely many C-molecules by O-atoms, such that each C-molecule has circuit rank $\rho(W)\leq\Gamma$, and joining these C-molecules by O-atoms does not add any cycle. We start with the following simple intuitive lemma:
\begin{itemize}
\item For each molecule $\Mb$ with at most $\Gamma$ independent cycles, we can always divide it into at most $\Gamma+1$ subsets $\Mb_j\,(1\leq j\leq \Gamma+1)$, such that each subset is a tree, and no atom in $\Mb_i$ can be parent of atom in $\Mb_j$ for $i<j$.
\end{itemize} 

For the proof of the above lemma, see Lemma \ref{lem.layer_refine}. Now, to complete the layer refining process for each single layer $\Mb_{\ell'}$, we need to delete some O-atoms to ensure each component of $\Mb_{\ell'}$ has $\rho$ value at most $\Gamma$, and then apply this lemma to divide it into at most $\Gamma+1$ \textbf{thin layers}, such that each thin layer is now a forest. For details see Definition \ref{def.layer_refine} and Proposition \ref{prop.layer_refine_2}. This then allows us to repeat the layer selection process in Section \ref{sec.toy_multi}; the use of thin layers instead of layers leads to one extra step in layer selection, but this is easily adjusted, see Definition \ref{def.layer_select}.

We remark that the above layer refining process does require the total number of independent cycles in these single-layer C-molecules to be negligible (compared to $\rho$). When this total number of cycles is large, we can perform an alternative cutting sequence that exploits these cycles, see Propositions \ref{prop.alg_up_recl} and \ref{prop.comb_est_case3}. This requires layer refining of a subset of layers and also leads to a minor auxiliary notion of \emph{primitivity}, see Definition \ref{def.strdeg} and Lemma \ref{lem.alg_up_ex}.
\subsubsection{Strong and weak degeneracies}\label{sec.reduce2} In Definition \ref{def.simplified} we have assumed that all \{33\} molecules are good. In reality there may be exceptions to this when the $(x_e,v_e)$ for \emph{fixed ends} $e$ at this \{33\} molecule takes some very specific values, as stated in the definition of good molecules (Definition \ref{def.good_normal}). For example, a \{33A\} molecule will not be good, if \begin{equation}\label{eq.degen_toy}|x_{e_1}-x_{e_7}|\leq \varepsilon^{1-\upsilon} \quad\textrm{and}\quad |v_{e_1}-v_{e_7}|\leq \varepsilon^{1-\upsilon},\end{equation} where $e_1$ and $e_7$ are the two fixed ends of the \{33\} molecule.

As such, we define a \emph{strong degeneracy} to be an instance where (\ref{eq.degen_toy}) happens at two edges that \emph{would become the two fixed ends if we cut the two atoms as a \{33A\} molecule}; see Definition \ref{def.strdeg}. Similarly we can define the notion of \emph{weak degeneracy} associated with the \{33B\} molecules, which involves one or two atoms and some conditions opposite to the ones stated in Definition \ref{def.good_normal}; see Definition \ref{def.weadeg}. Note that these definitions \emph{only involve} certain edges and atoms in $\Mb$, and is \emph{not} a priori associated with any concrete cutting sequence. This provides the flexibility, so that we can classify the degeneracies and non-degeneracies \emph{before} performing any cutting, and perform different cutting operations according to different cases, which allows us to minimize the need of backtracking.

Our strategy is that, when there are many degeneracies, we shall perform an \emph{alternative cutting sequence}, where we do not gain from \{33\} molecules but instead gain from the extra degeneracy assumptions such as (\ref{eq.degen_toy}). Compared to the sophisticated algorithms exploiting \{33\} molecules, this alternative cutting sequence essentially follows the simple procedure where we cut every atom by itself, so that most of them are cut as \{3\} molecules. Usually integrating a \{3\} molecule does not gain any power, but if this \{3\} molecule contains an edge or atom that is involved in some degeneracy condition like (\ref{eq.degen_toy}), then it is clear that \emph{this} \{3\} molecule integration should gain some power, which makes it a good molecule, see Definition \ref{def.good_normal}.

In practice, there will be two places in the overall algorithm where we need to consider degeneracies: the first time for strong degeneracies, see Proposition \ref{prop.case2} and the second time for weak degeneracies, see Proposition \ref{prop.comb_est_case4}. If there are many (strong or weak) degeneracies, then we can apply the alternative cutting sequence to conclude; if there are few (strong and weak) degeneracies, then their effect will be negligible and we can get rid of them and proceed with the main thread of our algorithm.
\subsubsection{Additional ingredients and pre-processing}\label{sec.reduce5} Finally, we discuss a few minor ingredients that are needed in the full case, that are not discussed in Sections \ref{sec.reduce1}--\ref{sec.reduce2} above.

\emph{(a) Initial cumulants.} In this section we have assumed $r(\Mb_0)=H_0=\varnothing$ in Proposition \ref{prop.mol_axiom}; in reality this $H_0$ may exist, which corresponds to the initial (time $t=0$) cumulants. This may affect the layer selection process, as it is now possible in Proposition \ref{prop.mol_axiom} (\ref{it.axiom6}) that most particle lines in $r(\Mb_{\ell'})$ form initial links instead of being connected to other particle lines or cycles. In this case we simply exploit the initial link condition to obtain good \{3\} molecules (see Definition \ref{def.good_normal}); the will be quite easy, see Proposition \ref{prop.comb_est_case4}.

\emph{(b) Cutting exceptional components.} Recall the notion of exceptional components in Definition \ref{def.layer_select_toy} \ref{it.select_toy3}; in Section \ref{sec.toy_multi} we have assumed there is no such exceptional components, so $\Mb_U=\Mb_{\ell_U}\cup\Mb_{(\ell_D:\ell_U)}$ is a forest (Proposition \ref{prop.layer_select_toy}). In general, when such exceptional components exist (but with a negligible number), we need to first cut them as free, before defining $\Mb_U$ and $\Mb_D$ which can be forests. This cutting process is described in Definition \ref{def.layer_cutting}; it may create some \{4\} molecules, but their number can be bounded in Propositions \ref{prop.layer_select} and \ref{prop.layer_cutting}. It may also introduce fixed ends to $\Mb_D$, which then requires the \emph{pre-processing} step to treat, see (d) below.

\emph{(c) Cutting 2-connections.} Recall that in the 2-layer models, there is a case when the number of 2-connections satisfies $0<\#_{2\mathrm{conn}}\ll |X|$ (see Section \ref{sec.toy3}), which is not covered in any toy model. In this case we first cut all these 2-connections as free, so no more 2-connections remain in the resulting molecule; this corresponds to \textbf{Stage 1} in the proof of Proposition \ref{prop.case5} in Section \ref{sec.finish}. After this, there will be fixed ends in $\Mb_D$ similar to (b) above, so we should apply the extra pre-processing step (see (d) below).

\emph{(d) Pre-processing.} In Definition \ref{def.simplified}, we have assumed that the 2-layer molecule $\Mb$ has no fixed ends; however, in reality, such fixed ends may exist as shown in (b)--(c) above, which also leads to fixed ends and deg 3 atoms after we subsequently cut all deg 2 atoms. In fact those fixed ends in $\Mb_U$ will not matter as there will not be any top fixed ends; however the fixed ends in $\Mb_D$ may affect the properness of $\Mb_D$ (which is needed in the \textbf{MAINUD} algorithm) and needs to be addressed.

In this case, depending on the number of components of deg 3 atoms in $\Mb_D$ (similar to the dichotomy in Sections \ref{sec.toy1+}--\ref{sec.toy2}), we either apply \textbf{3COMPDN} (i.e. the dual of \textbf{3COMPUP}) in $\Mb_D$ similar to Section \ref{sec.toy2} with enough gain (this corresponds to \textbf{Stage 2} in Section \ref{sec.finish}, or apply an extra \emph{pre-processing} step to turn $\Mb_D$ into a proper molecule (this corresponds to \textbf{Stage 3}) . The pre-processing depends on a specific process \textbf{SELECT}, see Definition \ref{def.func_select} and Proposition \ref{prop.func_select}.

After pre-processing, we are then ready to apply the dichotomy in Sections \ref{sec.toy1+}--\ref{sec.toy2}, to reduce to either toy model II (corresponding to \textbf{Stage 5}), or toy model I plus (corresponding to \textbf{Stage 6}). Finally, note that the pre-processing step may also cause some complications, as the set $X$ involved in the dichotomy needs to be replaced by some
natural modification $X_1$, which will be explained in details in the \textbf{MAINUD}, see the definition of $X_0$ and $X_1$ in Definition \ref{def.alg_maincr}, and \textbf{Stage 4}.
\subsubsection{Strategy flowchart} 
\begin{figure}[h!]
    \centering
    \includegraphics[width=0.9\linewidth]{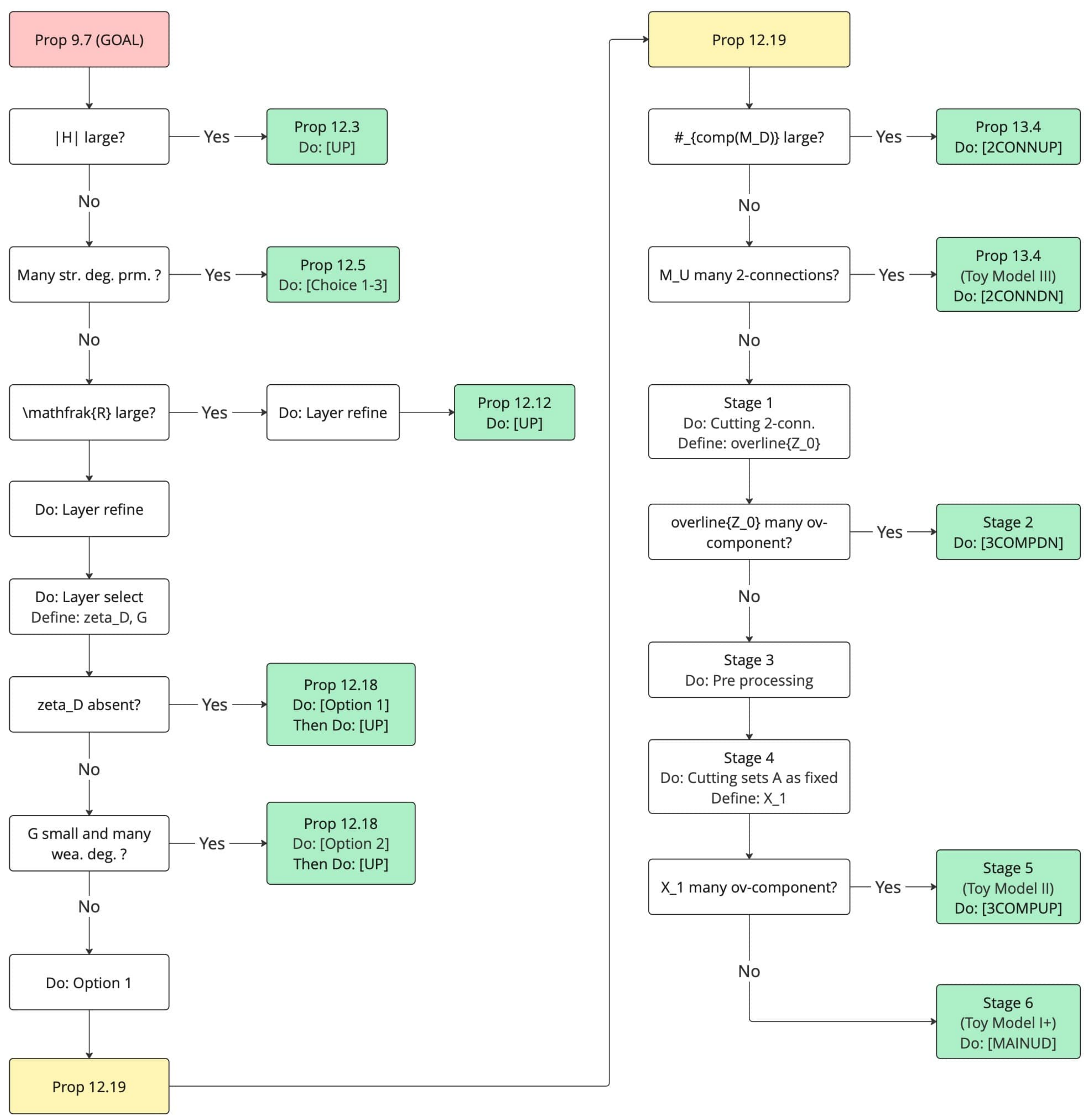}
    \caption{Yes/No flow chart for the proof of Proposition \ref{prop.comb_est}. Here each green box indicates end of a proof path; each ``Do" describes the algorithm, process or steps we take, and ``Define" indicates the new notions defined. The notions of strongly degenerate, primitive and weakly degenerate (abbreviated str. deg., prm. and wea. deg.) pairs are defined in Definitions \ref{def.strdeg} and \ref{def.weadeg}; the other notions are defined at the place indicated. The ``many/large/small" etc. are all relative to $\rho$; for the precise threshold of these statements, see the corresponding propositions.}
    \label{fig.flowchart}
\end{figure}

We present a flowchart summarizing the proof of Proposition \ref{prop.comb_est}, see {\color{blue}Figure \ref{fig.flowchart}}.
\section{Cutting algorithm I: reduction to UD molecules}\label{sec.layer} In this section and Section \ref{sec.maincr} we will construct the cutting algorithm in Proposition \ref{prop.comb_est}. As already mentioned in the beginning of Section \ref{sec.toy}, these algorithms and arguments will be presented in a different order from Section \ref{sec.toy}, to be consistent with the logical flow of the proof. The goal of this section is to prove several cases of Proposition \ref{prop.comb_est}, namely Propositions \ref{prop.case1}, \ref{prop.case2}, \ref{prop.comb_est_case3} and \ref{prop.comb_est_case4}, and reduce Proposition \ref{prop.comb_est} to Proposition \ref{prop.case5}, by getting rid of degeneracies (Propositions \ref{prop.case2} and \ref{prop.comb_est_case4}), introducing layer refinement (Section \ref{sec.layer_refine}) and layer selection (Section \ref{sec.layer_select}).

\subsection{The \textbf{UP} algorithm and the large $|H|$ case} In this subsection, we prove the first case of Proposition \ref{prop.comb_est}, namely Proposition \ref{prop.case1}, where we assume $|H|$ is large. We start by introducing the full version of algorithm \textbf{UP}, which is a direct extension of its toy model version in Definition \ref{def.up_toy}, by replacing parent/child with ov-parent/ov-child (cf. Section \ref{sec.reduce1}). This is a general purpose or ``default" algorithm that minimizes the number of \{4\} molecules, which also guarantees sufficiently many good molecules \emph{only in some cases}. However, this algorithm alone is enough to treat the simple scenarios studied in this section. Moreover, even in later cases where other algorithms are needed to achieve gain in a key subset of the molecule, it is usually helpful to use \textbf{UP} (or its dual \textbf{DOWN}) to finish the decomposition of the remaining parts of the molecule without loss.
\begin{definition}[The algorithm \textbf{UP}]\label{def.alg_up} Define the following cutting algorithm \textbf{UP}. It takes as input any (regular) molecule $\Mb$, which has \emph{no top fixed end at any C-atom} (we henceforth abbreviate it as \textbf{C-top fixed ends}; same for \textbf{C-bottom fixed ends}). For any such $\Mb$, we define the cutting sequence as follows:
\begin{enumerate}
\item\label{it.alg_up_1} If $\Mb$ contains any deg $2$ atom $\nf$, then cut it as free, and repeat until there is no deg $2$ atom left.
\item\label{it.alg_up_2} Consider the set of all remaining deg 3 atoms in $\Mb$ (i.e. those deg 3 that have not been cut), and choose a lowest atom $\nf$ in this set. If $\Mb$ has no deg 3 atoms, then choose a lowest atom $\nf$ in $\Mb$. Let $S_\nf$ be the set of descendants of $\nf$ (including $\nf$); this $\nf$ and $S_\nf$ will be fixed until the end of (\ref{it.alg_up_3}).
\item \label{it.alg_up_3} Starting from $\nf$, each time choose a \emph{highest} remaining atom $\mf$ in $S_\nf$. If $\mf$ has deg 3, and either an ov-parent $\mf^+$ or an ov-child $\mf^-$ of $\mf$ also has deg 3, then cut $\{\mf,\mf^+\}$ or $\{\mf,\mf^-\}$ as free; otherwise cut $\mf$ as free. Repeat until all atoms in $S_\nf$ have been cut. Then go to (\ref{it.alg_up_1}).
\end{enumerate}

We also define the dual algorithm \textbf{DOWN}, by reversing the notions of parent/child, lowest/highest, top/bottom etc. It applies to any molecule that has no C-bottom fixed end.
\end{definition}

The next proposition is the full version of Proposition \ref{prop.up_toy}.
\begin{proposition}
\label{prop.alg_up} Let $\Mb$ be any molecule as in Definition \ref{def.alg_up}. We may assume $\Mb$ is connected (otherwise consider each component of $\Mb$). Then after applying algorithm \textbf{UP} to $\Mb$ (and same for \textbf{DOWN}):
\begin{enumerate}
\item\label{it.up_proof_1} We have the following monotonicity property (MONO): throughout the process of \textbf{UP}, the set $S_\nf\backslash\{\nf\}$ will not have any \emph{C-bottom} fixed end, and $(\Mb\backslash S_\nf)\cup\{\nf\}$ will not have any \emph{C-top} fixed end (in particular, when all atoms in $S_\nf$ have been cut and before the next $\nf$ is chosen in Definition \ref{def.alg_up} (\ref{it.alg_up_2}), there is no C-top fixed end in whole $\Mb$). In particular, the molecule $\Mb$ is cut into elementary molecules, and $\#_{\{33B\}}=\#_{\{44\}}=0$.
\item\label{it.up_proof_2} We have $\#_{\{4\}}\leq 1$, and $\#_{\{4\}}=1$ if and only if $\Mb$ is full. 
\end{enumerate}
In (\ref{it.up_proof_3})--(\ref{it.up_proof_4}) below we assume $\Mb$ has no deg 2 atoms.
\begin{enumerate}[resume]
\item\label{it.up_proof_3} If $\Mb$ contains a cycle or contains at least two deg $3$ atoms, then $\#_{\{33A\}}\geq 1$. If not (i.e. $\Mb$ is a tree and contains at most one deg 3 atom), then all the elementary molecules are \{3\} molecules with at most one exception which might be a \{4\} molecule.
\item\label{it.up_proof_4} If $\Mb$ contains at most one deg 3 atom \emph{and} contains a cycle, then we can split into $2$ sub-cases such that for each sub-case we have $\#_{\mathrm{good}}\geq 1$.
\end{enumerate}
\end{proposition}
\begin{proof} We first make some comments on the result and the proof. Here, compared to (MONO-toy) in Proposition \ref{prop.up_toy}, the (MONO) property changes the absence of fixed ends to the absence of C-fixed ends, because of the presence of O-atoms (we do not need to worry about fixed ends at O-atoms because they cannot affect any other atom; see Remark \ref{rem.reg}).

The main new ingredient compared to Proposition \ref{prop.up_toy} is (\ref{it.up_proof_4}), which guarantees a \emph{good} molecule instead of just a \{33A\} molecule, and is thus useful in strongly degenerate cases (Definition \ref{def.strdeg}, Proposition \ref{prop.case2}) where \{33A\} molecules might not be good. The proof of (\ref{it.up_proof_4}) also provides an example how the splitting into sub-cases in Proposition \ref{prop.comb_est} works, by decomposing 1 into suitable indicator functions with prescribed supports. Apart from these, the proof is essentially the same as Proposition \ref{prop.up_toy}.

\textbf{Proof of (\ref{it.up_proof_1}).} We first prove (MONO). Note that this statement is true initially, as $\Mb$ has no C-top fixed end. Moreover at each time of choosing $\nf$ in Definition \ref{def.alg_up} (\ref{it.alg_up_2}), there is also no C-bottom fixed end in $S_\nf\backslash\{\nf\}$, because all atoms in $S_\nf\backslash\{\nf\}$ must have deg 4 (since $\nf$ is the lowest deg 3 atom). Therefore, it suffices to show that (MONO) is preserved. For this, first note that if $\nf$ has deg 2 as in Definition \ref{def.alg_up} (\ref{it.alg_up_1}), and has no C-top fixed end, then cutting $\nf$ as free will only generate bottom fixed ends at its ov-parents (cf. Remark \ref{rem.reg}), so the absence of C-top fixed end in the whole molecule is preserved.

Next, suppose $\nf$ is chosen as in Definition \ref{def.alg_up} (\ref{it.alg_up_2}), and $\mf$ is a highest remaining atom in $S_\nf$ as in Definition \ref{def.alg_up} (\ref{it.alg_up_3}). If we cut $\mf$ as free, this does not create any C-bottom fixed end in $S_\nf\backslash\{\nf\}$ (because $\mf$ is a highest remaining atom, so any ov-parent of $\mf$ in $S_\nf\backslash\{\nf\}$ has already been cut), and does not create any C-top fixed end in $(\Mb\backslash S_\nf)\cup\{\nf\}$ (because any ov-child atom of $\mf$ must again be in $S_\nf\backslash\{\nf\}$), so (MONO) is preserved.

The same can be proved for the operation of cutting $\{\mf,\mf^+\}$ or $\{\mf,\mf^-\}$; let us consider $\{\mf,\mf^+\}$ as the other case is similar. Here, cutting $\{\mf,\mf^+\}$ may create bottom fixed ends at ov-parents of $\mf$ and $\mf^+$, but none of them belong to $S_\nf\backslash\{\nf\}$ (so (MONO) is not affected), because $\mf$ is a highest remaining atom (so no ov-parent of $\mf$ still remains in $S_\nf$, in particular $\mf^+\not\in S_\nf$ and no ov-parent of $\mf^+$ can be in $S_\nf$). Moreover, cutting $\{\mf,\mf^+\}$ may create top fixed ends at ov-children of $\mf$ and $\mf^+$, but they all have to be in $S_\nf\backslash\{\nf\}$ (so again (MONO) is not affected). This last statement is obvious for $\mf$ as $\mf\in S_\nf$. As for $\mf^+$, since it has deg 3, has no C-top fixed end, and has one ov-child being $\mf$, we know it must have a bottom fixed end (or is O-atom, cf. Remark \ref{rem.reg}). Since $\mf$ is ov-child of $\mf^+$, we know that any other ov-child of $\mf^+$ must belong to the maximal ov-segment containing $\mf^+$ and $\mf$. Using the structure of ov-segments, we know that any other ov-child of $\mf^+$ that is C-atom must occur at the bottom of this ov-segment, thus it must be an ov-child of $\mf$ and belong to $S_\nf\backslash\{\nf\}$. This shows (MONO) is preserved in the whole process, thus it is true for all time.

Now, note that by construction, any atom in $\Mb$ is either cut by itself, or belongs to a \{33\} molecule. By (MONO), all the single-atom molecules we cut must have deg at least 2, and those of deg 2 must have two fixed ends either both being top or both being bottom (or non-serial for O-atoms, due to regularity). Similarly, all the \{33\} molecules $\{\mf,\mf^+\}$ and $\{\mf,\mf^-\}$ cut in Definition \ref{def.alg_up} (\ref{it.alg_up_3}) must be \{33A\}, because at the time it is cut, (MONO) implies that the $\mf^+\not\in S_\nf$ must be an ov-parent of $\mf$ with no C-top fixed end, and $\mf^-\in S_\nf\backslash\{\nf\}$ must be an ov-child of $\mf$ with no C-bottom fixed end (in particular, if we cut $\mf$ as free from $\{\mf,\mf^\pm\}$, then $\{\mf^\pm\}$ will become an elementary \{2\} molecule). This proves that all molecules we cut are elementary, and $\#_{\{33B\}}=\#_{\{44\}}=0$.

\textbf{Proof of (\ref{it.up_proof_2}).} If $\Mb$ has only deg 4 atoms, then the first cutting in $\Mb$ must create a \{4\} molecule; at any time thereafter, $\Mb$ will contain at least one non-deg-4 atom (cf. Remark \ref{rem.full_cut}), and by Definition \ref{def.alg_up} (\ref{it.alg_up_2})--(\ref{it.alg_up_3}), any subsequent cutting will \emph{not} create any \{4\} molecule (the $\nf$ chosen in Definition \ref{def.alg_up} (\ref{it.alg_up_2}) must have deg 3; also any $\mf\in S_\nf\backslash\{\nf\}$ will have deg at most 3 when it is cut in Definition \ref{def.alg_up} (\ref{it.alg_up_3}), because $\mf$ has an ov-parent in $S_\nf$ that has already been cut before). The same holds true for any $\Mb$ that contains non-deg-4 atoms initially, in which case we do not get any \{4\} molecule. This proves (\ref{it.up_proof_2}).

\textbf{Proof of (\ref{it.up_proof_3}).} Here and in the proof of (\ref{it.up_proof_4}) below we assume $\Mb$ has no deg 2 atoms. We claim that if $\#_{\{33A\}}=0$, then necessarily also $\#_{\{2\}}=0$. Assume the contrary, consider the first atom $\mf$ that has deg 2 when it is cut, and consider the cutting operation that turns $\mf$ into deg 2, which must be cutting a single atom $\pf$ (otherwise $\pf$ would belong to a \{33A\} molecule). Note that $\mf$ must have deg 3 at the time $\pf$ is cut (because cutting $\pf$ only generates one fixed end at $\mf$, see Proposition \ref{prop.mol_axiom} (\ref{it.axiom1})), and $\pf$ also has deg 3 or 4 when it is cut (it cannot have deg 2 because $\mf$ is the first). If $\pf$ has deg 4 when it is cut (so $\{\pf\}$ is cut as a \{4\} molecule), then by the proof in (\ref{it.up_proof_2}) we know that all atoms must have deg 4 initially, and $\pf$ must the first atom cut; but that would imply that $\mf$ has deg 4 at the time $\pf$ is cut, contradiction.

We now know that $\{\pf\}$ is cut as a \{3\} molecule, and both $\pf$ and $\mf$ have deg 3 at the time $\pf$ is cut. As $\pf$ must be an ov-parent or ov-child of $\mf$ at this time (because cutting $\pf$ creates a fixed end at $\mf$), by Definition \ref{def.alg_up} (\ref{it.alg_up_3}), we know that we should have cut $\pf$ together with another atom instead of cutting it alone as a \{3\} molecule, contradiction.

Next, define $\sigma:=\#_{\mathrm{bo/fr}}-3|\Mb|$ and consider its increment under cutting operations, where $\#_{\mathrm{bo/fr}}$ is the number of bonds plus the number of free ends (not counting fixed ends), and $|\Mb|$ is the number of (remaining) atoms in $\Mb$. It is straightforward to check that cutting a \{2\}, \{33\}, \{3\} and \{4\} molecule increases the value of $\sigma$ by $1$, $1$, $0$ and $-1$. For example, cutting a \{4\} molecule turns $4$ bonds/free ends into fixed ends when there is no O-atom (the case of O-atoms is similar and the result is the same) which decreases $\#_{\mathrm{bo/fr}}$ by $4$, and also decreases $|\Mb|$ by $1$, leading to $\Delta\sigma=-1$. In view of the initial and final values of $\sigma$, we then get
\begin{equation}\label{eq.incr}\#_{\{4\}}-(\#_{\{2\}}+\#_{\{33A\}})=\sigma_{\mathrm{init}}=|\Mb|-|\Mb_3|-\#_{\mathrm{bond}},\end{equation} where $\sigma_{\mathrm{init}}$ is the initial value of $\sigma$ (with final value $\sigma_{\mathrm{fin}}=0$), $\#_{\mathrm{bond}}$ is the number of bonds (initially) in $\Mb$ and $|\Mb_3|$ is the number of deg 3 atoms; here the last equality is because each atom has deg $3$ or $4$, and each bond is exactly counted twice, so $\#_{\mathrm{bo/fr}}=(\textrm{total\ deg\ of\ all\ atoms})-\#_{\mathrm{bond}}=4|\Mb|-|\Mb_3|-\#_{\mathrm{bond}}$.

Since $\Mb$ is connected, we know $|\Mb|-\#_{\mathrm{bond}}\leq 1$, with equality if and only if $\Mb$ is a tree. Now, if $\#_{\{33A\}}=0$ (thus $\#_{\{2\}}=0$), then we have $\#_{\{4\}}+|\Mb_3|=|\Mb|-\#_{\mathrm{bond}}$ by (\ref{eq.incr}), which is only possible when $|\Mb|-\#_{\mathrm{bond}}=1$ (so $\Mb$ is a tree) and $|\Mb_3|\leq 1$; conversely, if $\Mb$ is a tree and has at most one deg 3 atom, then $|\Mb|-\#_{\mathrm{bond}}=1$ and $|\Mb_3|+\#_{\{4\}}=1$ by (\ref{it.up_proof_2}), so we get $\#_{\{2\}}=\#_{\{33A\}}=0$. This proves (\ref{it.up_proof_3}).

\textbf{Proof of (\ref{it.up_proof_4}).} Finally we prove (\ref{it.up_proof_4}). Following \textbf{UP}, we first cut $\nf\in\Mb$ as free, which is either the only deg 3 atom, or some deg 4 atom in $\Mb$; note that this does not create any deg 2 atom (because it creates at most one fixed end at any atom, see Proposition \ref{prop.mol_axiom} (\ref{it.axiom1})). In the subsequent cuttings, consider the first \{33A\} molecule (which exists due to (\ref{it.up_proof_2})) and the two fixed ends $(e_1,e_2)$ at it, which are created when cutting previous atoms $\pf_1$ and $\pf_2$ as free, and neither $\pf_1$ nor $\pf_2$ has deg 2 when it is cut, by the same proof in (\ref{it.up_proof_3}). If $\pf_1\neq\pf_2$, say $\pf_1$ is cut before $\pf_2$, then there exist two free ends $e_1'$ and $e_2'$ in the elementary molecules $\{\pf_1\}$ and $\{\pf_2\}$ respectively (i.e. with $e_1'$ belonging to the same ov-segment as $e_1$ in the original molecule, and same for $e_2'$ and $e_2$), such that ${z}_{e_1'}={z}_{e_1}$ and ${z}_{e_2'}={z}_{e_2}$ (cf. Remark \ref{rem.cut_order}; see also {\color{blue}Figure \ref{fig.up_good}}).
\begin{figure}[h!]
    \centering
    \includegraphics[width=0.65\linewidth]{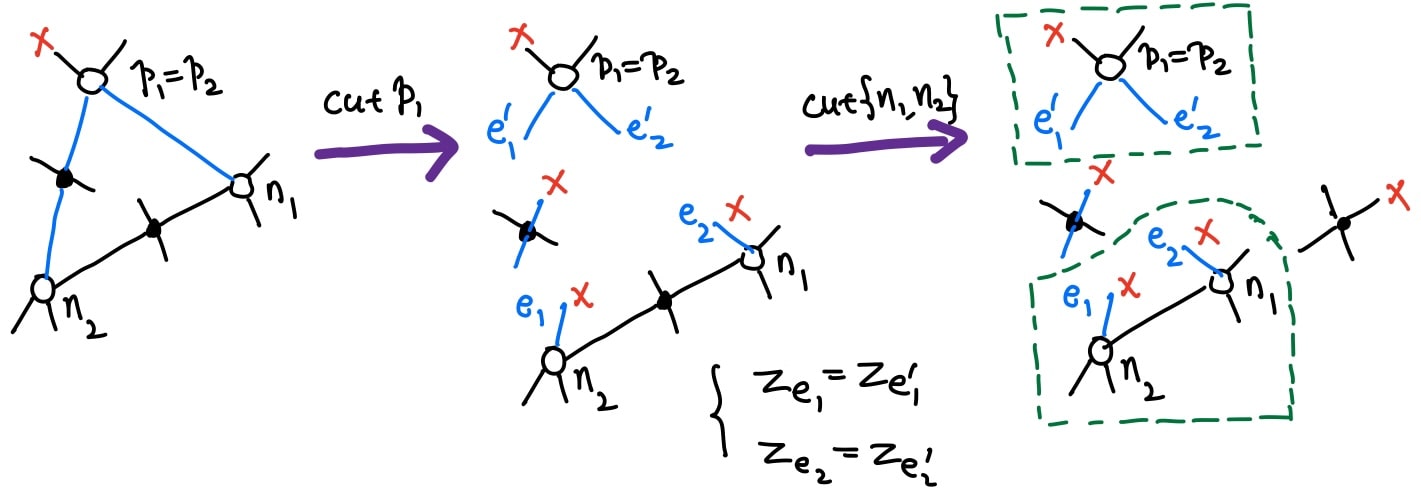}
    \caption{The proof of Proposition \ref{prop.alg_up} (\ref{it.up_proof_4}). Here both $\pf_1=\pf_2$ is cut as free before the \{33A\} molecule $\{\nf_1,\nf_2\}$ (and thus $\prec_{\mathrm{cut}}\{\nf_1,\nf_2\}$ as in Definition \ref{def.ordered_disjoint}), so there exist free ends $e_1'$ and $e_2'$ in $\{\pf_1\}$ such that we have $z_{e_1}=z_{e_1'}$ and $z_{e_2}=z_{e_2'}$ for the fixed ends $e_1$ and $e_2$ in $\{\nf_1,\nf_2\}$ (cf. Remark \ref{rem.cut_order}), which leads to the splitting.}
    \label{fig.up_good}
\end{figure}

We perform the \textbf{splitting} by considering the dichotomy: \[\mathrm{either}\quad|{x}_{e_1'}-{x}_{e_2'}|\leq\varepsilon^{1-\upsilon} \quad \mathrm{or}\quad |{x}_{e_1'}-{x}_{e_2'}|>\varepsilon^{1-\upsilon}.\] In the former case, we may insert the corresponding indicator function, and this extra restriction implies that $\{\pf_2\}$ is cut as a \emph{good} \{3\} molecule due to Definition \ref{def.good_normal} (note that $\pf_1$ is cut before $\pf_2$). In the latter case we have $|{x}_{e_1'}-{x}_{e_2'}|=|{x}_{e_1}-{x}_{e_2}|>\varepsilon^{1-\upsilon}$, and this extra restriction implies that the \{33A\} molecule is good due to Definition \ref{def.good_normal}.

Similarly, if $\pf_1=\pf_2$, then both $e_1'$ and $e_2'$ are free ends at $\pf_1$. Moreover, if $\pf_1$ is an O-atom, then $e_1'$ and $e_2'$ cannot be serial (otherwise the two atoms in the \{33A\} molecule are connected by an ov-segment containing both $(e_1,e_1')$ and $(e_2,e_2')$ and goes through $\pf_1$, so they are connected by two different ov-segments, violating Proposition \ref{prop.mol_axiom} (\ref{it.axiom1})). We then perform the same splitting; if $|{x}_{e_1'}-{x}_{e_2'}|\leq \varepsilon^{1-\upsilon}$ , then $\{\pf_1\}$ is cut as a good \{3\} or \{4\} molecule due to Definition \ref{def.good_normal}; if $|{x}_{e_1'}-{x}_{e_2'}|=|{x}_{e_1}-{x}_{e_2}|> \varepsilon^{1-\upsilon}$, then again the \{33A\} molecule is good due to Definition \ref{def.good_normal}. This completes the proof.
\end{proof}
Now, with the \textbf{UP} algorithm, we can already prove the first simple case of Proposition \ref{prop.comb_est}, namely when $|H|\geq (C_{12}^*)^{-1}\cdot\rho$, see Proposition \ref{prop.case1}. In this case the proof is easy due to the positive term $|H|$ in (\ref{eq.overall_alg}).
\begin{proposition}\label{prop.case1} Proposition \ref{prop.comb_est} is true in the case when $|H|\geq (C_{12}^*)^{-1}\cdot\rho$.
\end{proposition}
\begin{proof} Under the assumption $|H|\geq  (C_{12}^*)^{-1}\rho$, we simply cut the whole molecule $\Mb$ using \textbf{UP}. By Definition \ref{def.good_normal} we know that $\#_{\mathrm{bad}}=\#_{\{4\}}+\#_{\mathrm{ee}}$, where $\#_{\mathrm{ee}}$ is the number of empty ends. Recall (\ref{eq.overall_alg}) in Proposition \ref{prop.comb_est}. We shall prove that $\nu:=|H|-\#_{\{4\}}-\#_{\mathrm{ee}}\geq 0$, and that (note that $d-1\geq 1$ and $C_{13}^*\gg C_{12}^*$, so (\ref{eq.overall_alg_1_rep}) is stronger than (\ref{eq.overall_alg}))
\begin{equation}\label{eq.overall_alg_1_rep}
(\upsilon/2)\cdot\#_{\mathrm{good}}+\nu\geq (\upsilon/8)\cdot(C_{12}^*)^{-1}\cdot \rho.
\end{equation}

In the proof below, in Part 1 we will show $\nu\geq 0$ by using Proposition \ref{prop.alg_up} and bounding the number of \emph{connected components} of $\Mb$ by $|H|$ (the empty ends will be trivial), which is done by repeatedly iterating Proposition \ref{prop.mol_axiom} (\ref{it.axiom4})--(\ref{it.axiom5}) and showing that every component of $\Mb$ must intersect a particle line in $H$, a process we call the (PL) argument below. In Parts 2--3, by using Proposition \ref{prop.mol_axiom} (\ref{it.axiom6}), we will show that either $\nu\gtrsim|H|$ (i.e. the number of components of $\Mb$ is strictly less than $|H|$), or a constant percentage of these components contain cycles or initial links. In the former case we already gain enough powers, while in the latter case we can obtain sufficiently many good molecules either by Proposition \ref{prop.alg_up} (\ref{it.up_proof_4}) or by exploiting the initial links.

\textbf{Proof part 1.} We first prove $\nu\geq 0$. By Proposition \ref{prop.alg_up} (\ref{it.up_proof_2}) we have $\#_{\{4\}}\leq \#_{\mathrm{comp}}$, where $\#_{\mathrm{comp}}$ is the number of components of $\Mb$, so it suffices to show $\#_{\mathrm{comp}}+\#_{\mathrm{ee}}\leq |H|$. Here we present an argument which will be used many times later in different variants, which we refer to as the \textbf{(PL) argument} below (see {\color{blue}Figure \ref{fig.plargument}}). Start from any component $W$ of $\Mb$; by Proposition \ref{prop.mol_axiom} (\ref{it.axiom4}), this $W$ must intersect a particle line $\pb\in p(\Mb_{\ell'})$ for some $\ell'$. By Proposition \ref{prop.mol_axiom} (\ref{it.axiom5}), this $\pb$ is equal to or connected to some $\pb'\in r(\Mb_{\ell'+1})$ via $\Mb_{\ell'+1}$; in either case $\pb'$ also intersects $W$ (note that each particle line can intersect at most one component of $\Mb$, see Remark \ref{rem.layer_interval}). We then replace $\pb$ by $\pb'$ and repeat the same argument to get $\pb''\in r(\Mb_{\ell'+2})$ and so on, in the end we see that $W$ must intersect a particle line $\boldsymbol{q}\in r(\Mb_{\ell})=H$ (in particular this $\boldsymbol{q}\in H$ is not an empty end). Since each $\boldsymbol{q}$ intersects at most one component of $\Mb$, this shows that $\#_{\mathrm{comp}}$ is bounded by the number of particle lines in $H$ that are not empty ends.
\begin{figure}[h!]
    \centering
    \includegraphics[width=0.28\linewidth]{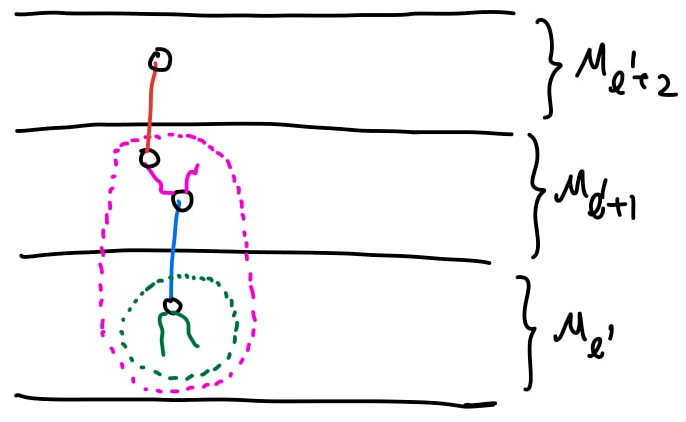}
    \caption{The (PL) argument in the proof of Proposition \ref{prop.case1}: start from the blue particle line $\pb$ intersecting a component of $\Mb_{\ell'}$ (green dotted circle, green edges) that is $\subseteq W$. This $\pb$ is connected to the red particle line $\pb'$ via $\Mb_{\ell'+1}$, leading to a larger component of $\Mb_{[\ell':\ell'+1]}$ (pink dotted circle, pink edges) that is also $\subseteq W$ and now intersects $\pb'$; then connect $\pb'$ to another particle line $\pb''$, and so on. This argument will be used in a few variants below (e.g. Proposition \ref{prop.layer_select} relying on Proposition \ref{prop.layer_refine_2}).}
    \label{fig.plargument}
\end{figure}

On the other hand, if $\pb$ is any particle line which is an empty end, then $\pb\in p(\Mb_{\ell'})$ for some $\ell'$, and contains no atom; by repeating Proposition \ref{prop.mol_axiom} (\ref{it.axiom5}) we know that $\pb$ must also belong to $H$. This means that $\#_{\mathrm{ee}}$ is equal to the number of particle lines in $H$ that are empty ends, so putting together we obtain that $\#_{\mathrm{comp}}+\#_{\mathrm{ee}}\leq |H|$, as desired.

\textbf{Proof part 2.} Now we have $\nu\geq |H|-\#_{\mathrm{comp}}-\#_{\mathrm{ee}}\geq 0$, it suffices to prove (\ref{eq.overall_alg_1_rep}). Define
\begin{itemize}
\item $X$ to be the number of particle lines in $H$ that are connected to another particle line in $H$ via $\Mb$,
\item $Y$ to be the number of components of $\Mb$ that contains a cycle,
\item and $Z$ to be the number of initial links involving two non-empty-end particle lines, such that at least one of them belongs to a component that does not contain a cycle.
\end{itemize} We claim that
\begin{equation}\label{eq.case1_1}X+Y+2Z\geq |H|-4\cdot\#_{\mathrm{ee}}.
\end{equation} In fact, by Proposition \ref{prop.mol_axiom} (\ref{it.axiom6}), each particle line $\pb$ in $H$ either (i) belongs to the $X$ particle lines listed above, or (ii) is not connected to other particle lines in $H$ but is connected to a cycle, or (iii) does not belong to (i) or (ii) but forms an initial link. Here in (ii) the $\pb$ must intersect one of the $Y$ components with cycle, and each component is counted at most once (if two particle lines correspond to the same component then they would be connected and thus would be in (i)), leading to at most $Y$ particle lines. In (iii), suppose $\pb$ forms an initial link that involves $\pb'$ as in Definition \ref{def.connectedvia}, then each $\pb'$ corresponds to a unique $\pb$ (otherwise they are again in (i)), and the number of $\pb'$ is at most $2Z+2\cdot\#_{\mathrm{ee}}$ (those $Z$ initial links not involving empty ends contribute at most $2Z$ choices of $\pb'$, and those initial links involving empty ends contribute at most $4\cdot\#_{\mathrm{ee}}$ choices of $\pb'$ because each empty end belongs to at most two initial links). This proves (\ref{eq.case1_1}).

Consider all components of $\Mb$ (indexed by some $\alpha$), and the number of particle lines in $H$ intersecting each component (denoted by $m_\alpha$); by the above proof we have $m_\alpha\geq 1$. Note that the $X$ defined above satisfies $X=\sum_{m_\alpha\geq 2}m_\alpha$, and the total number of non-empty-end particle lines is equal to $\sum_\alpha m_\alpha$, we conclude that\begin{equation}\label{eq.case1_2}\nu\geq |H|-\#_{\mathrm{ee}}-\#_{\mathrm{comp}}=\sum_{\alpha}(m_\alpha-1)=\sum_{m_\alpha\geq 2}(m_\alpha-1)\geq X/2.\end{equation}

\textbf{Proof part 3.} Next we prove that
\begin{equation}\label{eq.case1_3}
\#_{\mathrm{good}}\geq Y+Z/2+\#_{\mathrm{ee}}.
\end{equation} In fact, each empty end must belong to a link by Proposition \ref{prop.mol_axiom} (\ref{it.axiom6}), so it must be good. In the process of \textbf{UP} we may choose to first cut the components with cycles; by Proposition \ref{prop.alg_up} (\ref{it.up_proof_4}), each of these components contributes at least one good molecule under \textbf{UP} (possibly after splitting into 2 sub-cases, but the total number of sub-cases caused by this is bounded by $2^{|\Mb|}$ and is acceptable in Proposition \ref{prop.comb_est}), so we get at least $Y$ good molecules.

After this, we are left with a forest without fixed end, so each of the subsequent elementary molecules we cut will be \{3\} or \{4\} molecules by Proposition \ref{prop.alg_up} (\ref{it.up_proof_3}). Moreover, for each of the $Z$ initial links, at least one of the involved particle lines (i.e. the one belonging to a no-cyclic component) has still not been cut; we call the corresponding bottom free end $e_1$ and call the one initial linked to it $e_2$, note that $e_2$ may or may not have been cut.

In the subsequent cuttings of the remaining forest, let $\sigma_j$ be the maximal ov-segment containing $e_1$, and assume $\sigma_2$ is broken (i.e. has any atom on it cut) before $\sigma_1$. Let $\pf$ be the first atom in $\sigma_1$ that is cut, then $\pf$ must be cut as a \{3\} or \{4\} molecule and has no fixed end along $\sigma_1$ (as it is the first atom cut in $\sigma_1$), so the initial link condition guarantees that $\pf$ is a good (\{3\} or \{4\}) molecule by Definition \ref{def.good_normal}, see {\color{blue}Figure \ref{fig.linkproof}}. Similarly, if $\sigma_1$ and $\sigma_2$ are broken at the same time, say when cutting an atom $\pf$, then the initial link condition again guarantees that $\{\pf\}$ is a good \{3\} or \{4\} molecule.
\begin{figure}[h!]
    \centering
    \includegraphics[width=0.25\linewidth]{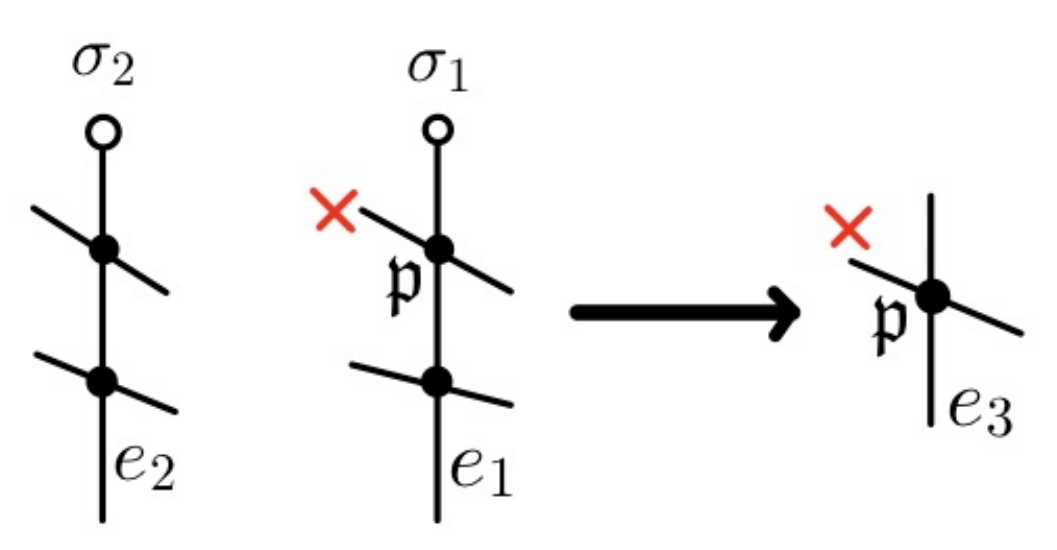}
    \caption{Part 3 of the proof of Proposition \ref{prop.case1}. Here $\sigma_2$ is broken before $\sigma_1$, and $\pf$ is the first atom cut in $\sigma_1$ (which is a \{3\} molecule in the picture). In the picture, $z_{e_2}$ has been fixed, and an initial link exists between $e_1$ and $e_2$, and $z_{e_3}=z_{e_1}$ (where $e_3$ is the edge at $\pf$), so $\{\pf\}$ is a good \{3\} molecule, thanks to the initial link condition and Definition \ref{def.good_normal}.}
    \label{fig.linkproof}
\end{figure}

In summary, each of the $Z$ initial links provides one good molecule in the process of cutting forests, leading to $Z/2$ good molecules (each good molecule is counted at most twice as each particle line is involved in at most 2 initial links). This produces $Z/2$ extra good molecules and thus proves (\ref{eq.case1_3}). Finally, putting together (\ref{eq.case1_1})--(\ref{eq.case1_3}), we get
\begin{equation}\label{eq.case1_4}\#_{\mathrm{good}}+\nu\geq\frac{1}{2}(X+2Y+Z+2\cdot\#_{\mathrm{ee}})\geq\frac{1}{4}|H|\geq (4C_{12}^*)^{-1}\cdot\rho.\end{equation} As the left hand side of (\ref{eq.case1_4}) is strictly positive and is trivially $\geq 1$, this proves (\ref{eq.overall_alg_1_rep}).
\end{proof}
\subsection{Strong degeneracies}\label{sec.strdeg} From now on we only need to consider the case $|H|<(C_{12}^*)^{-1}\cdot\rho$ in Proposition \ref{prop.comb_est} (we will make this assumption until the end of Section \ref{sec.maincr}), and the goal is to prove (\ref{eq.overall_alg}). Note also that $\#_{\mathrm{bad}}=\#_{\{4\}}+\#_{\mathrm{ee}}$ and $|H|\geq \#_{\mathrm{ee}}$ by the same proof as in Proposition \ref{prop.case1}, so it suffices to prove
\begin{equation}\label{eq.case2_1}
(\upsilon/2)\cdot \#_{\mathrm{good}}-d\cdot\#_{\{4\}}\geq (C_{13}^*)^{-1}\cdot\rho.
\end{equation} In this subsection we treat the case of sufficiently many strong degeneracies (Definition \ref{def.strdeg} below), see Proposition \ref{prop.case2} (also cf. Section \ref{sec.reduce2}). For technical reasons (explained in the proof of Proposition \ref{prop.case2} below) we also need the extra notion of \textbf{primitivity} in Definition  \ref{def.strdeg} below. 
\begin{definition}\label{def.strdeg}
Let $\Mb$ be a molecule as in Proposition \ref{prop.mol_axiom}. Define a pair of atoms $\{\nf,\nf'\}$ in $\Mb$ to be \textbf{strongly degenerate} if \begin{enumerate}[{(i)}]
\item They are ov-adjacent and in the same layer ($\ell[\nf]=\ell[\nf']$).
\item There exist edges $(e,e')$ at $\nf$ and $\nf'$ respectively, such that cutting $\{\nf,\nf'\}$ as free in $\Mb$ and turning $(e,e')$ into fixed ends results in a \{33A\} molecule, as in Definition \ref{def.elementary};
\item We have the restriction that $|x_{e}-x_{e'}|\leq\varepsilon^{1-\upsilon}$ and $|v_{e}-v_{e'}|\leq\varepsilon^{1-\upsilon}$.
\end{enumerate}

We also define an atom pair $\{\nf,\nf'\}$ in $\Mb$ to be \textbf{primitive}, if they are in the same layer ($\ell[\nf]=\ell[\nf']$), $\nf$ is ov-parent of $\nf'$, and there \emph{do not exist} C-atoms $\qf_j\,(1\leq j\leq r)$ with $r\geq 1$, such that $\qf_{j+1}$ is ov-child of $\qf_j$ for each $j$ with $\qf_0=\nf$ and $\qf_{r+1}=\nf'$.
\end{definition} 
\begin{proposition}\label{prop.case2} Recall Definition \ref{def.strdeg}. Suppose we make the further restriction that $\Mb$ contains at least $(C_{12}^*)^{-1}\rho$ disjoint strongly degenerate primitive pairs (by inserting the corresponding indicator function $\mathbbm{1}_{\mathrm{str.deg.prm}}$), then Proposition \ref{prop.comb_est} is true.
\end{proposition}
\begin{proof} Recall $|H|<(C_{12}^*)^{-1}\cdot\rho$ and the goal is to prove (\ref{eq.case2_1}). The idea is to restrict sufficiently many strong degeneracies to a single layer $\Mb_{\ell'}$ (by definition, the two atoms in any strongly degenerate pair must be in the same layer), and cut most of the atoms involved in these degeneracies as \{3\} molecules which will be good due to the extra restrictions given by strong degeneracy, see \textbf{Choice 3} below. If this is not possible, then there must exist sufficiently many cycles in $\Mb_{\ell'}$ (i.e. $Y\gtrsim\rho$ in our notions below), so we can try to gain from the cyclic clusters using Proposition \ref{prop.alg_up} (\ref{it.up_proof_4}), see \textbf{Choice 2}. If this is also not possible, then $\Mb_{\ell'}$ must have many components (i.e. $X\gtrsim\rho$ in our notions below); in this case we can cut each pair of strongly degenerate atoms as a good \{44\} molecule, see \textbf{Choice 1}. However, we need to make sure that cutting this \{44\} molecule does not violate the (MONO) property needed to execute the \textbf{UP} algorithm, which then naturally leads to the notion of primitive atom pairs in Definition \ref{def.strdeg}.

First assume there exists at least $\Af:=(C_{12}^*)^{-1}\rho$ disjoint strongly degenerate  primitive pairs. By splitting into at most $|\log\varepsilon|^{C^*\Af}$ sub-cases (i.e. decomposing $1$ into the corresponding indicator functions), we may fix one choice of these pairs. By using pigeonhole principle and replacing $\Af$ with $\Lf^{-1}\cdot\Af$, we may restrict all these atom pairs $\{\nf_j,\nf_j'\}$ to a single layer, say $\Mb_{\ell'}$. Define $\Qb$ to be the union of all connected omponents of $\Mb_{\ell'}$ that contain at least one strongly degenerate primitive pair $\{\nf_j,\nf_j'\}$. We then
\begin{enumerate}[{(i)}]
\item Cut $\Qb$ as free, then cut $\Mb_{>\ell'}$ as free from $\Mb\backslash\Qb$ and cut it into elementary molecules using \textbf{UP};
\item Then cut $\Mb_{<\ell'}\cup (\Mb_{\ell'}\backslash \Qb)$ into elementary molecules using \textbf{DOWN}.
\end{enumerate} It is easy to see that this can always be done, and the contribution of these steps (i.e. coming from atoms outside $\Qb$) to $\#_{\{4\}}$ is at most $2|H|$. This because the number of components of $\Mb_{>\ell'}$ is $\leq|H|$ by repeating the (PL) argument in the proof of Proposition \ref{prop.case1} (i.e. starting with a particle line intersecting one component and iteratively applying Proposition \ref{prop.mol_axiom} (\ref{it.axiom5})), and moreover each \emph{full} component of $\Mb_{<\ell'}\cup (\Mb_{\ell'}\backslash \Qb)$ \emph{after the rest has been cut} must also be a full component in the original molecule $\Mb$ (see Remark \ref{rem.full_cut}) and there are at most $|H|$ of them, so the upper bound on $\#_{\{4\}}$ follows from Proposition \ref{prop.alg_up} (\ref{it.up_proof_2}).

Now we only need to cut $\Qb$ into elementary molecules. For each connected component $\Qb'$ of $\Qb$, consider the clusters in $\Qb'$ as defined in Proposition \ref{prop.mol_axiom} (\ref{it.axiom3}) (i.e. ov-components of the set of C-atoms in $\Qb'$); denote these clusters by $\Qb_j$, then $\Qb'$ is formed by joining these $\Qb_j$ with O-atoms. Denote by $X$ the number of components of $\Qb$, and by $Y$ the total number of $\Qb_j$ that are cyclic (i.e. $\rho(\Qb_j)>0$ or $\Qb_j$ contains a cycle). Depending on the sizes of $X$ and $Y$, we may then choose to perform one of the following three cutting sequences.

\textbf{Choice 1}: in this choice we will gain from $X$. For each component $\Qb'$ of $\Qb$, choose one pair $\{\nf_j,\nf_j'\}$ in it, then cut $\{\nf_j,\nf_{j}'\}$ as free ({\color{blue}Figure \ref{fig.choice_1}}), which is a good \{44\} molecule due to the restrictions $|x_{e_j}-x_{e_j'}|\leq\varepsilon^{1-\upsilon}$ and $|v_{e_j}-v_{e_j'}|\leq\varepsilon^{1-\upsilon}$. Recall that $S_{\nf_j}$ is the set of descendants of $\nf_j$. We next prove that $\nf_j'$ has no ov-parent $\qf\neq\nf_j$ that is a C-atom in $S_{\nf_j}$. In fact, if such a $\qf$ exists, then $\qf$ is a descendant of $\nf_j$, so there exist a path from $\nf_j$ to $\qf$ by iteratively taking children. This path forms a cycle together with the ov-segments from $\qf$ to $\nf_j'$ and from $\nf_j'$ to $\nf_j$. As the cycle belongs to a single layer $\Mb_{\ell'}$, by Proposition \ref{prop.mol_axiom} (\ref{it.axiom2}) it must be formed by ov-segments connecting C-atoms; in particular, by considering the path between $\nf_j$ and $\qf$, we see that the $\qf_j$ must exist as in Definition \ref{def.strdeg}, so $\{\nf_j,\nf_j'\}$ is not primitive, contradicting our assumption that $\{\nf_j,\nf_j'\}$ is a strongly degenerate primitive pair.

Now, we have proved that the above $\qf$ does not exist; this implies that, after cutting $\{\nf_j,\nf_j'\}$ as free, there is no C-top fixed end in $\Qb'\backslash S_{\nf_j}$ (recall that $\nf_j'\in S_{\nf_j}$), and no C-bottom fixed end in $S_{\nf_j}$, in particular the monotonicity property (MONO) in Proposition \ref{prop.alg_up} is true. Since this property is preserved under \textbf{UP}, we can then repeat the steps in Definition \ref{def.alg_up} (\ref{it.alg_up_3}) and subsequently apply \textbf{UP} to the rest of $\Qb'$ to cut it into elementary molecules. By applying this to each component $\Qb'$ of $\Qb$, we obtain $\#_{\mathrm{good}}\geq X$ and $\#_{\{4\}}=0$.
\begin{figure}[h!]
    \centering
    \includegraphics[width=0.5\linewidth]{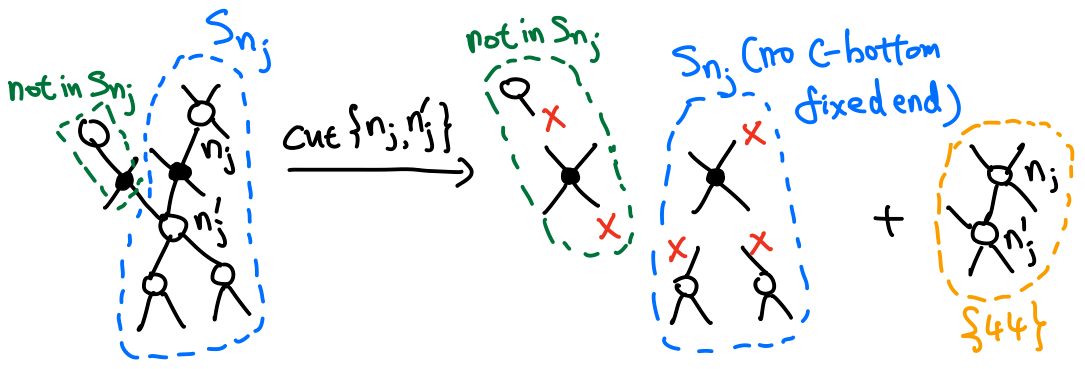}
    \caption{\textbf{Choice 1} in the proof of Proposition \ref{prop.case2}, cutting out a \{44\} molecule $\{\nf_j,\nf_j'\}$ which is good due to the strong degeneracy assumption. Note that cutting this \{44\} molecule does not violate the (MONO) property (i.e. creates no C-bottom fixed end in $S_\nf$) if and only if this pair is primitive.}
    \label{fig.choice_1}
\end{figure}

\textbf{Choice 2}: in this choice we will gain from $Y$. For each component $\Qb'$ of $\Qb$, let $\Qb_j$ be the clusters in $\Qb'$ defined as above, which are joined by O-atoms to form the full component $\Qb'$. For convenience, we allow some $\Qb_j$ to be empty ends (where joining an edge with an empty end simply means creating an O-atom in the interior of this edge, and drawing two serial free ends at this O-atom). We may reorder them as follows: first choose $\Qb_1$ arbitrarily, and each time choose one remaining cluster $\Qb_j$ that is joined with some $\Qb_{j'}\,(j'<j)$ by some O-atom (such $\Qb_j$ must exist because $\Qb'$ is connected). 

In this case, suppose (for the sake of contradiction) that $\Qb_j$ is joined with the clusters $\Qb_{j'}\,(j'<j)$ by \emph{two} different O-atoms $\pf_1$ and $\pf_2$, then we get a cycle by collecting (i) the two non-serial ov-segments connecting $\pf_1$ to C-atoms $\nf_1\in\Qb_j$ and $\mf_1\in \Qb_{j_1}\,(j_1<j)$, and the same for $\pf_2$ (with $\nf_2\in \Qb_j$ and $\mf_2\in \Qb_{j_2}\,(j_2<j)$), (ii) the ov-segments connecting $\nf_1$ to $\nf_2$ via O-atoms and C-atoms in $\Qb_j$, and (iii) the ov-segments $\mf_1$ to $\mf_2$ via O-atoms and C-atoms in $\Qb_{j'}\,(j'<j)$ (note that $\Qb_{j_1}$ is indirectly joined to $\Qb_{j_2}$ by a sequence of O-atoms due to our construction). In this way we get a cycle in $\Qb'$ with an O-atom $\pf_1$ such that the two bonds of this cycle at $\pf_1$ are \emph{not serial}, contradicting Proposition \ref{prop.mol_axiom} (\ref{it.axiom2}).

By the above arguments, we have reordered the $\Qb_j$ such that each $\Qb_j$ is joined with \emph{exactly one} $\Qb_{j'}\,(j'<j)$ by \emph{exactly one} O-atom, say $\pf_j$. We then cut $\Qb_1$ as free, which creates the \emph{unique} deg 3 atom $\pf_2$ in an ov-segment connecting two C-atoms in $\Qb_2$ (as $\pf_2$ is the only atom joining $\Qb_2$ with $\Qb_1$); then we cut $\Qb_2\cup\{\pf_2\}$ as free, after which we have the unique deg 3 atom $\pf_3$ in an ov-segment connecting two C-atoms in $\Qb_3$ (as $\pf_3$ is the only atom joining $\Qb_3$ with $\Qb_1\cup \Qb_2$), and so on, see {\color{blue}Figure \ref{fig.choice_2}}. Note that when we cut each $\Qb_j\cup\{\pf_j\}$, it always contains \emph{at most one} deg 3 atom (i.e. $\pf_j$ if $j\geq 2$ and none if $j=1$), so by Proposition \ref{prop.alg_up} (\ref{it.up_proof_4}), we can subsequently cut it into elementary molecules using \textbf{UP} or \textbf{DOWN} depending on whether the unique fixed end is bottom or top, and get a good molecule for each $\Qb_j$ that is \emph{cyclic}. 

By applying the above process to each component $\Qb'$ of $\Qb$, we obtain $\#_{\mathrm{good}}\geq Y$ (number of cyclic $\Qb_j$) and $\#_{\{4\}}\leq X$ (total number of components). We note that strongly degenerate and primitive is not used in this part.
\begin{figure}[h!]
    \centering
    \includegraphics[width=0.65\linewidth]{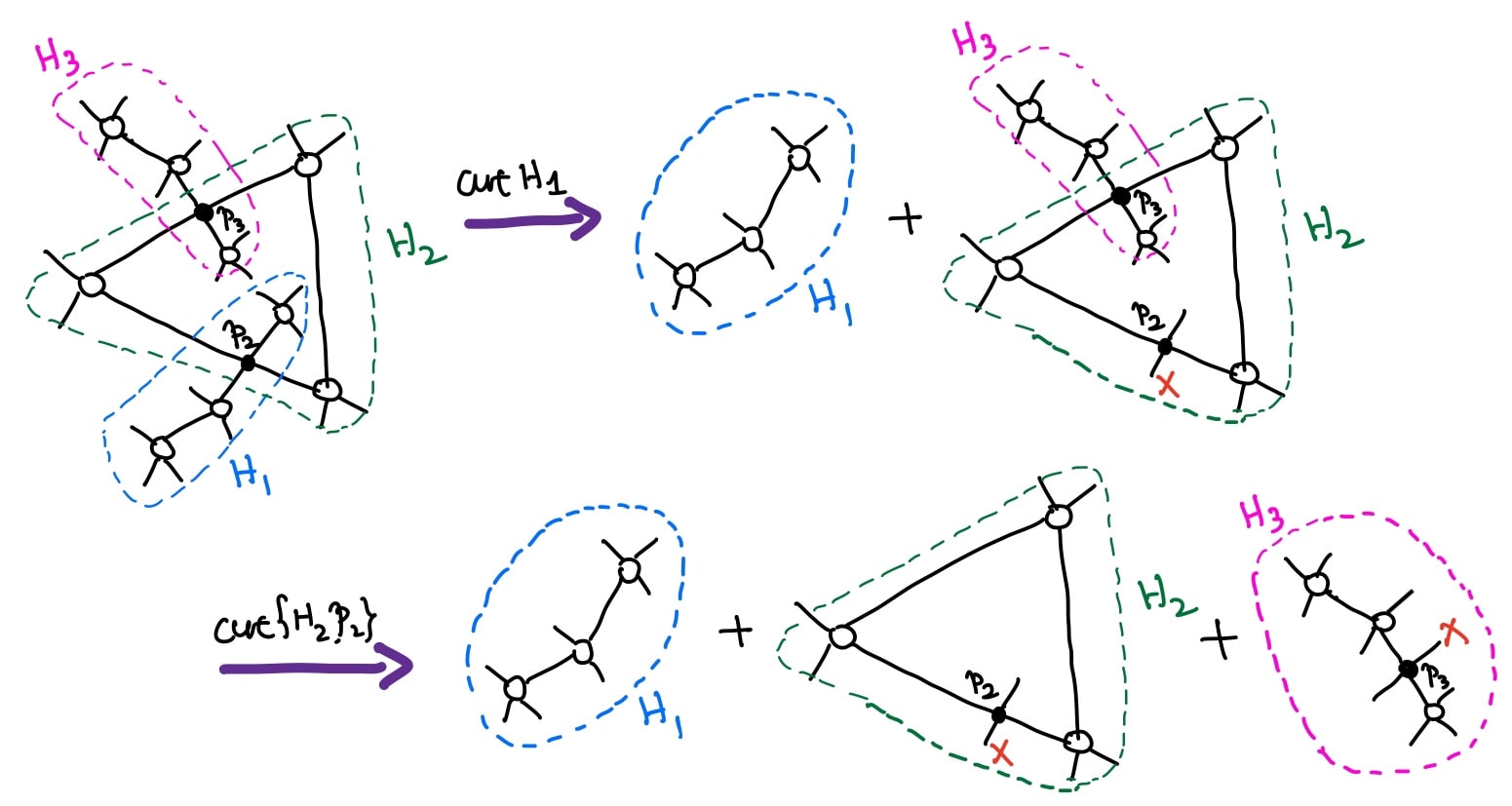}
    \caption{\textbf{Choice 2} in the proof of Proposition \ref{prop.case2}, each time cutting $\Qb_j$ together with $\pf_j$; note that each $\Qb_j\cup\{\pf_j\}$ contains a unique deg 3 atom which is $\pf_j$. This is also an example of the cutting order $\prec_{\mathrm{cut}}$ in Definition \ref{def.ordered_disjoint}, i.e. we have $\Qb_1\prec_{\mathrm{cut}}(\Qb_2\cup\{\pf_2\})\prec_{\mathrm{cut}}(\Qb_3\cup\{\pf_3\})$.}
    \label{fig.choice_2}
\end{figure}

\textbf{Choice 3}: in this choice we will gain from $\Af$. For each component $\Qb'$ of $\Qb$, we cut it as free using \textbf{UP}. To analyze this choice, for each $\Qb'$, we define the quantity $\sigma=\#_{\mathrm{bo/fr}}-3|\Qb'|$ as in the proof of Proposition \ref{prop.alg_up} (\ref{it.up_proof_3}). By the same arguments there, we get
\begin{equation}\label{eq.case2_2}\#_{\{4\}}-(\#_{\{2\}}+\#_{\{33A\}})=\sigma_{\mathrm{init}}=|\Qb'|-\#_{\mathrm{bond}}
\end{equation} (note that initially $\Qb'$ has no deg 3 atom), where $\#_{\mathrm{bond}}$ is the number of bonds in $\Qb'$. To calculate $\sigma_{\mathrm{init}}$, assume $\Qb'$ is formed by joining $q$ C-molecules with $q-1$ O-atoms as discussed in \textbf{Choice 2}, then we have 
\begin{equation}\label{eq.calcsigma}\sigma_{\mathrm{init}}=|\Qb'|-\#_{\mathrm{bond}}=\sum_{j=1}^q(1-\rho(\Qb_j))-(q-1)=1-\sum_{j=1}^q\rho(\Qb_j),\end{equation} where $\rho(\Qb_j)$ is the circuit rank of the C-molecule $\Qb_j$ (to see this, just note that joining two C-molecules by an O-atom decreases the value of $|\Qb'|-\#_{\mathrm{bond}}$ by 1 as it increases $|\Qb'|$ by 1 and increases $-\#_{\mathrm{bond}}$ by 2, while value of $|\Qb'|-\#_{\mathrm{bond}}$ before introducing any O-atom equals $\sum_{j=1}^q(1-\rho(\Qb_j))$ by definition).

Now by adding up (\ref{eq.case2_2}) and (\ref{eq.calcsigma}) for each component $\Qb'$, and note that $0\leq\rho(\Qb_j)\leq\Gamma$ with $\rho(\Qb_j)>0$ if and only if $\Qb_j$ is cyclic, we conclude that
\[\#_{\{4\}}-(\#_{\{2\}}+\#_{\{33A\}})\geq X-\Gamma\cdot Y\] (recall $X$ is the number of components $\Qb'$, and $Y$ is the number of cyclic $\Qb_j$); since also $\#_{\{4\}}\leq X$, we then have
\[(\textrm{number of molecules that are not \{3\}})=\#_{\{4\}}+(\#_{\{2\}}+\#_{\{33A\}})\leq X+\Gamma\cdot Y.\] In particular, by excluding this many exceptions, we see that the number of pairs $\{\nf_j,\nf_j'\}$ such that \emph{both $\nf_j$ and $\nf_j'$ are cut as \{3\} molecules} is at least $\Lf^{-1}\cdot\Af-2(X+\Gamma Y)$.

For each such pair $\{\nf_j,\nf_j'\}$, if $\nf_j$ is cut before $\nf_j'$, since $\{\nf_j'\}$ is cut as a \{3\} molecule, we know that $\{\nf_j'\}$ must be good due to the restriction $|x_{e_j}-x_{e'_j}|\leq\varepsilon^{1-\upsilon}$ and Definition \ref{def.good_normal}. Here when $\nf_j'$ is O-atom, the edge $e_j'$ cannot be serial with the bonds between $\nf_j$ and $\nf_j'$, otherwise $\{\nf_j,\nf_j'\}$ will not be regular (and thus not a \{33A\} molecule) after turning $e_j'$ into a fixed end. In this way we get $\#_{\mathrm{good}}\geq \Lf^{-1}\cdot\Af-2(X+\Gamma Y)$ and $\#_{\{4\}}\leq X$.

Summing up, by choosing one of the \textbf{Choices 1--3}, we get
\begin{equation}\label{eq.case2_5} (\upsilon/2)\cdot (\#_{\mathrm{good}})-d\cdot\#_{\{4\}}\geq\max\big((\upsilon/2)X,(\upsilon/2)Y-d\cdot X,(\Lf^{-1}\upsilon/2)\Af-(d+1)(X+\Gamma Y)\big)\geq (C_{2}^*)^{-1}\Af,
\end{equation} which proves (\ref{eq.case2_1}) as $\Af=(C_{12}^*)^{-1}\rho$.
\end{proof}
Note that Proposition \ref{prop.case2} deals with the case where there exist many strongly degenerate primitive pairs. As for the non-primitive pairs, it turns out that they always lead to good molecules if they are cut as \{33A\}, as shown in the following lemma.
\begin{lemma}\label{lem.alg_up_ex}
For each \{33A\} molecule $\{\nf,\nf'\}$ obtained \emph{in any cutting sequence after any deletion of O-atoms}, if the atom pair $\{\nf,\nf'\}$ is \emph{not primitive in $\Mb$ before any operation}, then by splitting into 2 sub-cases, we can either make this \{33A\} molecule good, or cut it into a good \{3\} and a normal \{2\} molecule.
\end{lemma}
\begin{proof} Consider any \{33A\} molecule $\{\nf,\nf'\}$ that is not primitive before any operation, and let $\qf_j$ be as in Definition \ref{def.strdeg}. Note that $\nf$ and $\nf'$ must belong to the same (original) layer $\Mb_{\ell'}$ (otherwise $\{\nf,\nf'\}$ is good by Definition \ref{def.good_normal}), and so does every $\qf_j$. Consider the cycle formed by the ov-segments from $\qf_j$ to $\qf_{j+1}$ (with $\nf=\qf_0$ and $\nf'=\qf_{r+1}$) and from $\nf'$ to $\nf$, which is contained in the single layer $\Mb_{\ell'}$. By Proposition \ref{prop.mol_axiom} (\ref{it.axiom2}), we see that both $\nf$ and $\nf'$ must be C-atom. Moreover the C-atoms $\qf_j$, and the ov-adjacency relations between them, are not affected by deletion of O-atoms. Now, consider the C-atom $\qf_1$ which is ov-child of $\nf$, and $\qf_r$ which is ov-parent of $\nf'$ (they could be the same atom).

Consider any cutting sequence. If $\qf_1$, or any (O-) atom $\rf$ in the ov-segment between $\qf_1$ and $\nf$ is cut before $\{\nf,\nf'\}$ (i.e. belongs to a molecule after cutting that is $\prec_{\mathrm{cut}}\{\nf,\nf'\}$), then we can make the \textbf{splitting} depending on whether (i) $|t_\nf-t_{\nf'}|\geq \varepsilon^{\upsilon}$ or (ii) $|t_\nf-t_{\nf'}|\leq \varepsilon^{\upsilon}$. In case (i) the $\{\nf,\nf'\}$ is already good by Definition \ref{def.good_normal}; in case (ii), due to the ov-child relations we have \[t_{\nf}\geq t_{\rf}\geq t_{\qf_1}\geq\cdots \geq t_{\qf_{r+1}}=t_{\nf'},\] in particular $|t_\nf-t_{\nf'}|\leq \varepsilon^{\upsilon}$ implies that \[|t_{\nf}-t_{\rf}|\leq \varepsilon^{\upsilon},\quad |t_{\nf'}-t_{\rf}|\leq \varepsilon^{\upsilon}.\] Now we may cut either $\nf$ or $\nf'$ from $\{\nf,\nf'\}$ as a \{3\} molecule such that the remaining atom becomes a \{2\} molecule (one of the two choices must be possible due to the definition of \{33A\} molecules), then this \{3\} molecule will be good by Definition \ref{def.good_normal}, so the lemma is proved.

In the same way, if $\qf_r$, or any (O-) atom $\rf$ in the ov-segment between $\qf_r$ and $\nf'$ is cut before $\{\nf,\nf'\}$, then the lemma is also proved by the same arguments as above.

Finally, if none of $\qf_1$, $\qf_r$ and none of the above-mentioned O-atoms is cut before $\{\nf,\nf'\}$, then at the time $\{\nf,\nf'\}$ is cut, the atom $\nf$ cannot have any bottom fixed end (because it is ov-parent of $\nf'$, and the ov-segment between $\nf$ and $\qf_1$ does not lead to any fixed end at $\nf$ as the other atoms on this ov-segment are all cut after $\nf$); in the same way $\nf'$ has no top fixed end at the time $\{\nf,\nf'\}$ is cut. However $\{\nf,\nf'\}$ is a \{33\} molecule, which means that $\nf$ must have a top fixed end and $\nf'$ has a bottom fixed end, and thus $\{\nf,\nf'\}$ will be a \{33B\} molecule instead of \{33A\}, contradiction. This completes the proof. 
\end{proof}

\begin{reduct}\label{red.1} With Propositions \ref{prop.case1} and \ref{prop.case2} and Lemma \ref{lem.alg_up_ex} proved, from now on, we may always assume in Proposition \ref{prop.comb_est} that $|H|\leq (C_{12}^*)^{-1}\rho$ and the number of disjoint strongly degenerate primitive pairs is at most $(C_{12}^*)^{-1}\rho$.

Using also Lemma \ref{lem.alg_up_ex}, we see that each \{33A\} molecule must be good with at most $(C_{12}^*)^{-1}\rho$ exceptions (these exceptions are when the corresponding atom pair is both strongly degenerate and primitive; if it is not strongly degenerate then we get a good molecule by Definition \ref{def.good_normal}, if it is not primitive we get a good molecule by Lemma \ref{lem.alg_up_ex}).

As such, we see that Proposition \ref{prop.comb_est} will follow if we can prove that
\begin{equation}\label{eq.overall_alg_3} (\upsilon/2)\cdot(\#_{\mathrm{good}})-d\cdot(\#_{\{4\}})\geq (C_{11}^*)^{-1}\cdot\rho
\end{equation} under the assumption that \textbf{all \{33A\} molecules are good}, which we will assume throughout the rest of the proof.
\end{reduct}

\subsection{Layer refinement and the large $\Rf$ case} \label{sec.layer_refine} 
\begin{figure}[h!]
    \centering
    \includegraphics[width=0.6\linewidth]{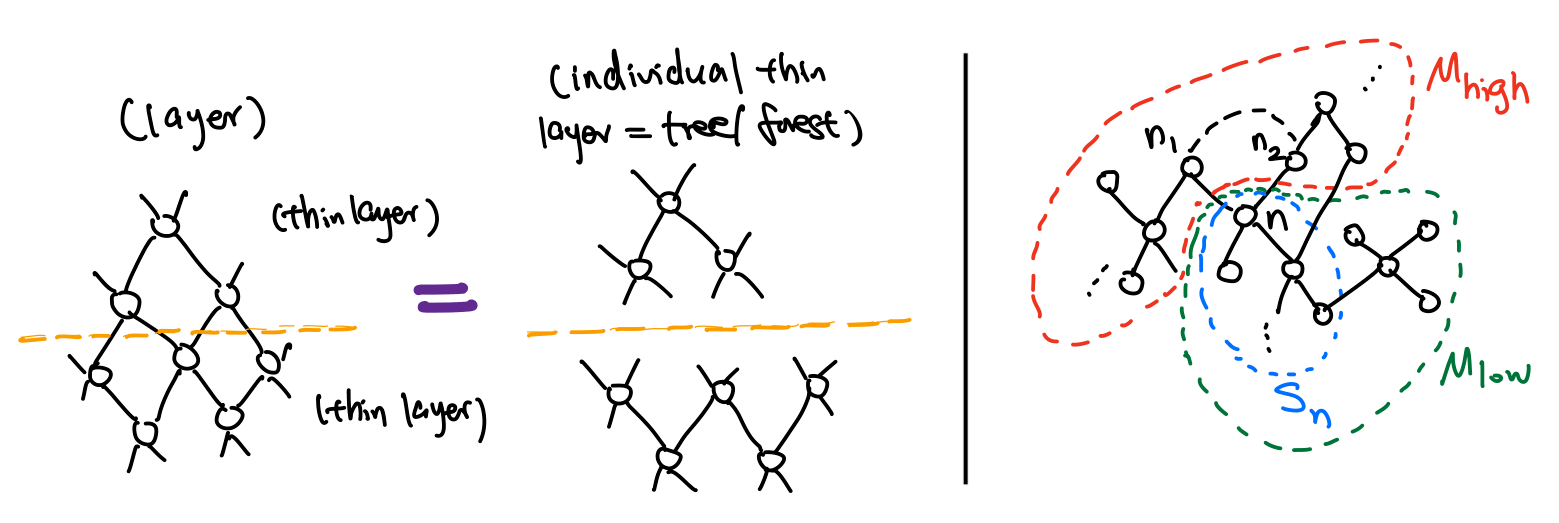}
    \caption{Left: example of layer refining (Lemma \ref{lem.layer_refine}, each thin layer is a forest). Right: the proof of Lemma \ref{lem.layer_refine}. Note that $\nf_1$ and $\nf_2$ are connected by a path above $\nf$ (forming a cycle); also $\Mb_{\mathrm{high}}$ consists of those atoms that can be connected to $\nf_1$ without entering $S_\nf$.}
    \label{fig.layerrefine}
\end{figure}

The goal of this subsection is to introduce the layer refining process (discussed in Section \ref{sec.reduce3}) and prove the next case of Proposition \ref{prop.comb_est} when $\Rf$ is large, namely Proposition \ref{prop.comb_est_case3}. The following lemma is already stated in Section \ref{sec.reduce3}:
\begin{lemma}\label{lem.layer_refine} Let $\Mb$ be a connected molecule with $\rho(\Mb)\leq\Gamma$, see Definition \ref{def.recollision_number}. Then we can divide $\Mb$ into at most $\Gamma+1$ disjoint subsets $\Mb_{(k)}$ for $1\leq k\leq\gamma\leq\Gamma+1$, such that (i) each $\Mb_{(k)}$ is a tree, and (ii) if $k<k'$ then no atom in $\Mb_{(k)}$ can be parent of any atom in $\Mb_{(k')}$. See {\color{blue}Figure \ref{fig.layerrefine}}. Note that, unlike the original layer defined by time intervals, there is \emph{no time inequality} between atoms in different thin layers within the same original layer, other than those coming from parent-child relations.
\end{lemma}
\begin{proof} We induct on $\rho(\Mb)$. If $\rho(\Mb)=0$ then $\Mb$ is a tree and there is nothing to prove. If $\rho(\Mb)>0$, then $\Mb$ must contain a cycle of atoms. Let a \emph{lowest} atom in this cycle be $\nf$, which is adjacent to two atoms $\nf_1$ and $\nf_2$ in this cycle, then both $\nf_j$ must be parents of $\nf$. Let $S_\nf$ be the set of descendants of $\nf$ (including $\nf$ itself), we then define
\begin{equation}\label{eq.layer_refine_1}\Mb_{\mathrm{high}}:=\big\{\mf\in\Mb:\mf\textrm{\ is\ connected\ to\ }\nf_1\textrm{\ by\ a\ path\ not\ containing\ atoms\ in\ }S_\nf\big\},
\end{equation} and $\Mb_{\mathrm{low}}:=\Mb\backslash\Mb_{\mathrm{high}}$, see {\color{blue}Figure \ref{fig.layerrefine}}. Note in particular $S_\nf\subseteq\Mb_{\mathrm{low}}$ and $\nf_1,\nf_2\in\Mb_{\mathrm{high}}$ ($\nf_2\in\Mb_{\mathrm{high}}$ because any atom in the cycle other than $\nf$ cannot belong to $S_\nf$). Clearly $\Mb_{\mathrm{high}}$ is connected, because for any $\mf\in\Mb_{\mathrm{high}}$ and any $\mf'$ on the path connecting $\mf$ to $\nf_1$ defined by (\ref{eq.layer_refine_1}), we must also have $\mf'\in\nf_1$ by definition, so each atom in $\Mb_{\mathrm{high}}$ is connected to $\nf_1\in \Mb_{\mathrm{high}}$ by a path completely within $\Mb_{\mathrm{high}}$.

To prove $\Mb_{\mathrm{low}}$ is connected, choose any atom $\mf\in\Mb_{\mathrm{low}}\backslash S_\nf$. Since $\Mb$ is connected, there is a path in $\Mb$ connecting $\mf$ to $\nf$. Let the first atom in this path (starting from $\mf$) that belongs to $S_\nf$ be $\pf$, and consider any atom $\qf\neq \pf$ between $\mf$ and $\pf$ on this path, then $\mf$ is connected to $\qf$ by a path not containing atoms in $S_\nf$, which means $\qf\in \Mb_{\mathrm{low}}$ (as $\qf\in\Mb_{\mathrm{high}}$ would imply $\mf\in\Mb_{\mathrm{high}}$). Thus we know each $\mf\in \Mb_{\mathrm{low}}\backslash S_\nf$ is connected to an atom $\pf\in S_\nf$ by a path completely within $\Mb_{\mathrm{low}}$, but $S_\nf$ is connected (every atom in $S_\nf$ is connected to $\nf$ by a path completely within $S_\nf$), so $\Mb_{\mathrm{low}}$ is also connected.

Moreover, if an atom $\mf_1\in \Mb_{\mathrm{high}}$ is connected to an atom $\mf_2\in\Mb_{\mathrm{low}}$ by a bond, then we must have $\mf_2\in S_\nf$ (otherwise $\mf_1\in\Mb_{\mathrm{high}}$ would imply $\mf_2\in\Mb_{\mathrm{high}}$). This means $\mf_1$ must be a parent of $\mf_2$, since otherwise $\mf_1$ is a child of $\mf_2$, so $\mf_2\in S_\nf$ implies $\mf_1\in S_\nf\subseteq \Mb_{\mathrm{low}}$, which is impossible. Finally, since there are at least two bonds between $\Mb_{\mathrm{high}}$ and $\Mb_{\mathrm{low}}$ (those connecting $\nf$ to $\nf_1$ and $\nf_2$), we see that
\[|\Mb|=|\Mb_{\mathrm{high}}|+|\Mb_{\mathrm{low}}|,\,\,\,|\Bc|\geq |\Bc_{\mathrm{high}}|+|\Bc_{\mathrm{low}}|+2\,\,\Rightarrow\,\,\rho(\Mb_{\mathrm{high}})+\rho(\Mb_{\mathrm{low}})\leq \rho(\Mb)-1\leq\Gamma-1,\] where $\Bc$ (resp. $\Bc_{\mathrm{high}}$ and $\Bc_{\mathrm{low}}$) is the set of bonds in $\Mb$ (resp. $\Mb_{\mathrm{high}}$ and $\Mb_{\mathrm{low}}$). We now apply induction hypothesis to $\Mb_{\mathrm{high}}$ and $\Mb_{\mathrm{low}}$, and note that no atom in $\Mb_{\mathrm{low}}$ is parent of an atom in $\Mb_{\mathrm{high}}$, to divide $\Mb$ into at most $(\rho(\Mb_{\mathrm{high}})+1)+(\rho(\Mb_{\mathrm{low}})+1)\leq\Gamma+1$ subsets that satisfy requirements (i) and (ii), so the proof is complete.
\end{proof}
With Lemma \ref{lem.layer_refine} we can then define the full layer refining process. The idea is to \emph{delete} those O-atoms joining two cyclic clusters (or cyclic components), so that each component contains at most one cyclic cluster after the deletions; in particular each component will have $\rho$ value at most $\Gamma$, so we can apply Lemma \ref{lem.layer_refine} to each component.
\begin{definition}[Layer refining]
\label{def.layer_refine} Let $\Mb$ be a molecule formed by layers $\Mb_{\ell'}$, as in Proposition \ref{prop.mol_axiom}. We define the following operation of \textbf{layer refining} to turn $\Mb_{\ell'}$ into a new layered molecule $\widetilde{\Mb}_{\ell'}$. Here the atoms and edges of $\widetilde{\Mb}_{\ell'}$ are the same as $\Mb_{\ell'}$ except some O-atoms being deleted as in Definition \ref{def.delete} (for convenience, we will not distinguish between $\Mb_{\ell'}$ and $\widetilde{\Mb}_{\ell'}$ below), but the layering of the atoms are different: the atoms in $\widetilde{\Mb}_{\ell'}$ are divided into at most $\Gamma+2$ disjoint subsets, each of them indicated by a refined layer number. We will call these \textbf{thin layers}. The layer refining process is defined as follows:
\begin{enumerate}
\item\label{it.thin_1} Recall all the clusters (i.e. C-molecules) in $\Mb_{\ell'}$. We repeatedly delete O-atoms in $\Mb_{\ell'}$, as long as each deletion disconnects two \emph{cyclic} clusters $\Qb_j$ (i.e. these two clusters belong to different components after deletion), until no more O-atoms can be deleted. For each connected component of $\Mb_{\ell'}$ that does not contain a cyclic cluster, we assign thin layer $(\ell',0)$ to all its atoms; note that with this definition, the thin layer $(\ell',0)$ plays a special role different from those thin layers $(\ell',k)$ for $k\geq 1$.
\item\label{it.thin_2} Apply Lemma \ref{lem.layer_refine} to each component $\Qb$ of $\Mb_{\ell'}$ that contains a cyclic cluster (we are guaranteed to have $\rho(\Qb)\leq\Gamma$, which will be proved in Proposition \ref{prop.layer_refine_2} below), to divide it into $\Qb_k$, where $1\leq k\leq \Gamma+1$. Then, for each atom in $\Mb_{\ell'}$, assign it thin layer $(\ell',k)$ if it belongs to $\Qb_k$ for some component $\Qb$.
\end{enumerate}

The thin layers defined above are linearly ordered by lexicographic ordering (where we first compare $\ell'$ and then compare $k$); if layer $\ell'$ is not refined, we identify it with $(\ell',0)$. Note that the original layer $\ell'$ precisely contains all the thin layers $(\ell',k)$. For any $\zeta$ and $\xi$, we define $\Mb_\zeta^T$ to be the set of atoms with thin layer $\zeta$; the associated notions $\Mb_{(\zeta:\xi)}^T$, $\Mb_{[\zeta:\eta]}^T$ etc. are defined similarly.
\end{definition}
In Proposition \ref{prop.layer_refine_2} below we state and prove some properties of the molecule $\Mb$ after layer refining. The important points are that (a) each thin layer after refinement is a forest, see Proposition \ref{prop.layer_refine_2} (\ref{it.layer_refine_22})--(\ref{it.layer_refine_23}); (b) the number of O-atoms deleted is bounded by $\Rf_{\ell'}:=\rho(\Mb_{\ell'})$ for the refined layers $\Mb_{\ell'}$, see Proposition \ref{prop.layer_refine_2} (\ref{it.layer_refine_21}); and (c) the original properties of $\Mb$ (notably Proposition \ref{prop.mol_axiom} (\ref{it.axiom4})--(\ref{it.axiom5}) and (\ref{it.axiom6})) still hold after layer refining, with a number of exceptions that is bounded by the number of deleted O-atoms, see Proposition \ref{prop.layer_refine_2} (\ref{it.layer_refine_24})--(\ref{it.layer_refine_27}).

\begin{figure}[h!]
    \centering
    \includegraphics[width=0.3\linewidth]{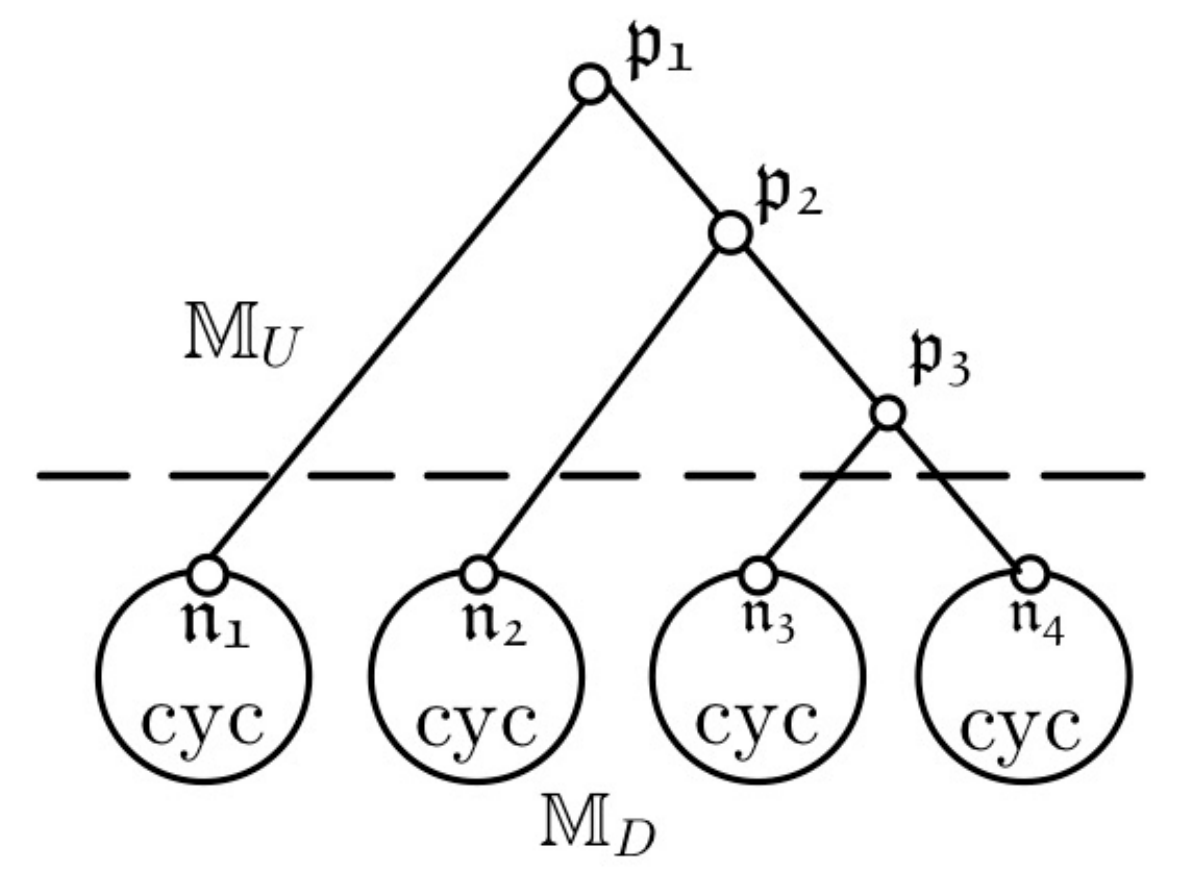}
    \caption{Proposition \ref{prop.alg_up_recl}. Here the cyclic clusters in $\Mb_D$ are shown, with $\nf_j$ being the first atom in each component that has deg 2 when cut, and $\pf_j\in\Mb_U$ is its ov-parent; see the proof of Proposition \ref{prop.alg_up_recl} below.}
    \label{fig.cyc_cluster}
\end{figure}

\begin{figure}[h!]
    \centering
    \includegraphics[width=0.5\linewidth]{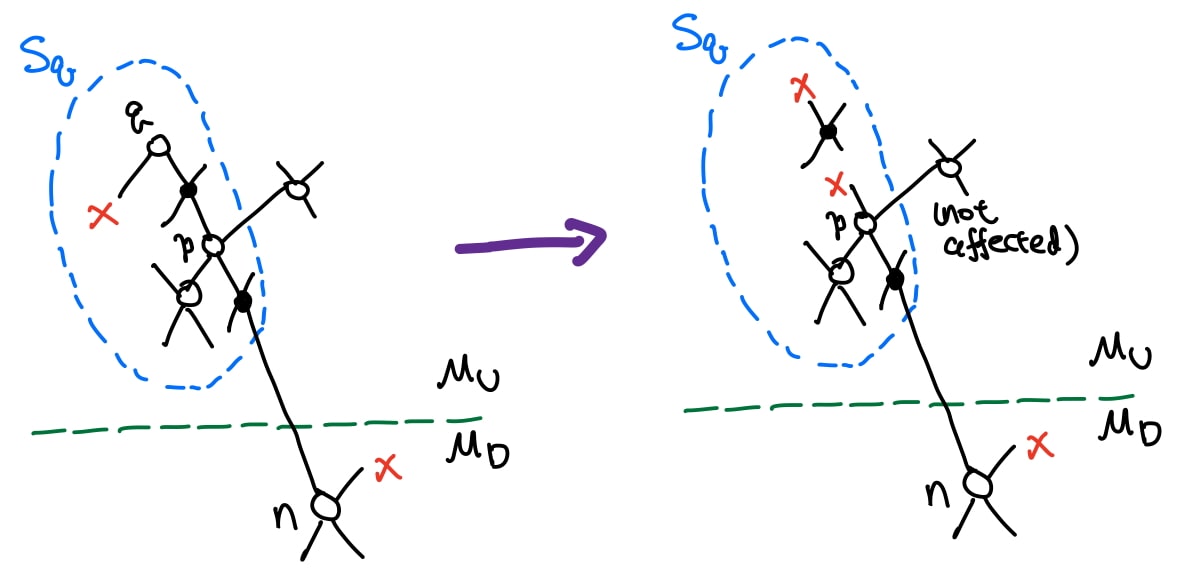}
    \caption{The proof of Proposition \ref{prop.alg_up_recl} showing $\pf$ has deg 3 when it is cut; this is because $\pf$ has two bottom free ends/bonds and one top free end/bond coming from a parent that does not belong to $S_\qf$ (and is not affected after cutting atoms in $S_\qf$). This fact will be used a few times below (such as in the proof of Proposition \ref{prop.comb_est_case4}).}
    \label{fig.up_aux}
\end{figure}

The proof of (a)--(b) above relies on the construction in Definition \ref{def.layer_refine}, in particular the fact that we delete all O-atoms joining two cyclic clusters, so each component contains at most one cyclic cluster after layer refining; see Parts 1--2 in the proof of Proposition \ref{prop.layer_refine_2} below. The proof of (c) is straightforward, basically because deleting each O-atom creates at most two exceptions to the relevant properties in Proposition \ref{prop.layer_refine_2} (\ref{it.layer_refine_24})--(\ref{it.layer_refine_27}), see Part 3 in the proof of Proposition \ref{prop.layer_refine_2} below.
\begin{proposition}
\label{prop.layer_refine_2} If we perform the layer refining for $\ell'$ as in Definition \ref{def.layer_refine}, then
\begin{enumerate}
\item\label{it.layer_refine_21} The number of O-atoms deleted in the process is at most $\Rf_{\ell'}:=\rho(\Mb_{\ell'})$.
\item\label{it.layer_refine_21+} After deleting the O-atoms, each component $\Qb$ of $\Mb_{\ell'}$ that contains a cyclic cluster satisfies $\rho(\Qb)\leq\Gamma$, in particular Lemma \ref{lem.layer_refine} is applicable in Definition \ref{def.layer_refine}.
\item\label{it.layer_refine_22} Each $\Mb_{(\ell',k)}^T$ (where $1\leq k\leq \Gamma+1$) is a forest, and has most $\Rf_{\ell'}$ connected components.
\item\label{it.layer_refine_23} $\Mb_{(\ell',0)}^T$ is a forest, and the number of components of $\Mb_{(\ell',0)}^T$ is not more than the number of components of $\Mb_{\ell'}$ before refining.
\end{enumerate}

Now consider the result of refining any set of layers in $\Mb$, which we still denote by $\Mb$ for convenience. Then we have the followings
\begin{enumerate}[resume]
\item\label{it.layer_refine_24} For $\zeta<\xi$, no atom in $\Mb_\zeta^T$ can be parent of any atom in $\Mb_\xi^T$.
\item\label{it.layer_refine_25} Each connected component of $\Mb_{(\ell',0)}^T$ (where $\ell'\leq\ell$) must intersect some particle line $\pb\in r(\Mb_{\ell'})$.
\item\label{it.layer_refine_26} With at most $2\Rf_{\ell'}$ exceptions, for each particle line $\pb\in p(\Mb_{\ell'})$, either $\pb$ intersects a component of $\Mb_{(\ell',0)}^T$, or it intersects a component of $\Mb_{(\ell',k)}^T$ for some $1\leq k\leq\Gamma+1$, or $\pb\in r(\Mb_{\ell'})$.
\item\label{it.layer_refine_27} Recall $\Rf=\sum_{\ell'}\Rf_{\ell'}$. For each $\ell'<\ell$, with at most $2\Rf$ exceptions, each particle line $\pb\in r(\Mb_{\ell'})$ either forms an initial link within $\Mb_{<(\ell'+1,0)}^T$, or is connected to a cycle or another particle line $\pb'\in r(\Mb_{\ell'})$ via the subset $\Mb_{<(\ell'+1,0)}^T$.
\end{enumerate}
\end{proposition}
\begin{proof} \textbf{Proof part 1.} First note the following facts about the refinement in Definition \ref{def.layer_refine}:
\begin{enumerate}[{(i)}]
\item Deleting any O-atom $\of$ in $\Mb_{\ell'}$ always breaks a component of $\Mb_{\ell'}$ into two;
\item Each component of $\Mb_{\ell'}$ that does not contain a cyclic cluster remain the same after refining. Each of the other components after refining contains a unique cyclic cluster.
\end{enumerate}

In fact, (i) is true because if $\nf_1$ and $\nf_2$ are the two atoms adjacent to $\of$ along two ov-segments in two different particle lines, then they cannot be connected after deleting $\of$, otherwise we get a cycle (before deleting $\of$) in view of the two bonds from $\of$ to $\nf_j\,(j\in\{1,2\})$, such that the two bonds at $\of$ are not serial in the cycle, contradicting Proposition \ref{prop.mol_axiom} (\ref{it.axiom2}). Also (ii) is obvious by our construction (each component not containing cyclic clusters is not affected, and each component after refining cannot contain two cyclic clusters, otherwise we would be able to delete one more O-atom).

\textbf{Proof part 2.} Next we prove (\ref{it.layer_refine_21})--(\ref{it.layer_refine_23}). Note first that the number of cyclic clusters in $\Mb_{\ell'}$ does not exceed $\Rf_{\ell'}$ (because each cyclic cluster contributes one the to the value of $\rho(\Mb_{\ell'})$). Now (\ref{it.layer_refine_21}) is easy because each deletion increases the number of cyclic components by 1 (due to (i) above), and this number cannot exceed $\Rf_{\ell'}$ which is the total number of cyclic clusters (imagine we delete all O-atoms, then the number of cyclic components would be equal to $\Rf_{\ell'}$), so the number of O-atoms deleted cannot exceed $\Rf_{\ell'}$. Moreover, (\ref{it.layer_refine_21+}) follows from the fact that each component contains only one cyclic cluster (the other clusters are trees and do not add to the $\rho$ value, and joining the clusters by O-atoms also does not increase the $\rho$ value).

Also (\ref{it.layer_refine_23}) is obvious because $\Mb_{(\ell',0)}^T$ has not been affected by the deletion of O-atoms, and is a forest because it does not contain any cyclic cluster (cf. Proposition \ref{prop.mol_axiom} (\ref{it.axiom2})). As for (\ref{it.layer_refine_22}), the number of components of $\Mb_{(\ell',k)}^T\,(k\geq 1)$ is at most $\Rf_{\ell'}$ because this is an upper bound for the number of components of $\Mb_{\ell'}$ containing one cyclic cluster; also each $\Mb_{(\ell',k)}^T\,(k\geq 1)$ is a forest which directly follows from Lemma \ref{lem.layer_refine}.

\textbf{Proof part 3}. Now we prove (\ref{it.layer_refine_24})--(\ref{it.layer_refine_27}), where we perform layer refining to any set of layers $\ell'$, and denote the result still by $\Mb$. To prove (\ref{it.layer_refine_24}), due to the lexicographic ordering of thin layers, we only need to prove it for a single layer $\ell'$, and we only need to consider those components containing one cyclic cluster $\Qb_j$ (otherwise all thin layers are $(\ell',0)$). But by definition, this follows from Lemma \ref{lem.layer_refine} (for each component of $\Mb_{\ell'}$) and the fact that no bond exists between different components of $\Mb_{\ell'}$.

Next we prove (\ref{it.layer_refine_25})--(\ref{it.layer_refine_26}). Note that particle lines can be defined in the same way before and after the refinement, since deleting an O-atom  corresponds to removing one atom from a particle line $\pb$ and merging two consecutive edges. From Proposition \ref{prop.mol_axiom} (\ref{it.axiom4})--(\ref{it.axiom5}) we already know that \emph{before layer refining}, any connected component in $\Mb_{\ell'}$ must intersect some particle line $\pb$ with $\pb\in r(\Mb_{\ell'})$. This then proves (\ref{it.layer_refine_25}) because $\Mb_{(\ell',0)}^T$ is not affected by layer refining: there is no bond connecting it to any $\Mb_{(\ell',k)}^T$ for $k>0$, and it also does not contain any of the O-atoms deleted in the refinement process. Similarly (\ref{it.layer_refine_26}) follows from Proposition \ref{prop.mol_axiom} (\ref{it.axiom4})--(\ref{it.axiom5}): any particle line $\pb\in p(\Mb_{\ell'})$ must either belong to $r(\Mb_{\ell'})$, or intersects $\Mb_{\ell'}$ (which is then divided into $\Mb_{(\ell',k)}^T$ for $0\leq k\leq \Gamma+1$). The only exceptions are those containing one of the deleted O-atoms, but the number of these particle lines is at most $2\Rf_{\ell'}$ as each O-atom belongs to at most two particle lines.

Finally we prove (\ref{it.layer_refine_27}). Note that before layer refining, this just corresponds to Proposition \ref{prop.mol_axiom} (\ref{it.axiom6}); however, after the layer refinement, we have deleted at most $\Rf$ O-atoms in $\Mb_{[1:\ell']}$ (note that by our conventions, $\Mb_{[1:\ell']}$ is just $\Mb_{<(\ell'+1,0)}^T$ after the refinement). If we add back these O-atoms one by one, then each time we create at most one cycle in one component, or merge at most two components of $\Mb_{[1:\ell']}$, and does not affect any initial link. Define a particle line $\pb$ to be \emph{single} if it does not satisfy the requirement of (\ref{it.layer_refine_27}), then adding back one O-atom at most reduces the number of single particle lines by $2$ (for example, if adding back an O-atom merges two components, then it may turn two single particle lines into non-single as they are now connected to each other; the other particle lines are unaffected, as they cannot intersect the relevant components by the single assumption). By Proposition \ref{prop.mol_axiom} (\ref{it.axiom6}), it then follows that the number of single particle lines is at most $2\Rf$, as desired.
\end{proof}
Before stating and proving Proposition \ref{prop.comb_est_case3} below (i.e. the case $\Rf\gtrsim\rho$), we need one more auxiliary result concerning a special case of \textbf{UP}. The intuition will be explained in the proof of Proposition \ref{prop.comb_est_case3}.
\begin{proposition}
\label{prop.alg_up_recl} Suppose $\Mb$ is a full molecule divided into two subsets $\Mb_U$ and $\Mb_D$, such that (i) $\Mb_U$ is a forest, and (ii) no atom in $\Mb_D$ is parent of atom in $\Mb_U$. Define $\#_{\mathrm{comp}(\Mb)}$ to be the number of components of $\Mb$, and define $\#_{\mathrm{cyc.comp}(\Mb_D)}$ to be the number of components of $\Mb_D$ that contains a cycle within this component. Then, see {\color{blue}Figure \ref{fig.cyc_cluster}}, after applying the algorithm \textbf{UP}, we have that $\#_{\{33B\}}=\#_{\{44\}}=0$, and
\begin{equation}\label{eq.alg_up_recl}\#_{\{4\}}\leq \#_{\mathrm{comp}(\Mb)},\quad\mathrm{and}\quad \#_{\{33A\}}\geq \#_{\mathrm{cyc.comp}(\Mb_D)}/4.
\end{equation}
\end{proposition}
\begin{proof} The statements about $\#_{\{33B\}}$, $\#_{\{44\}}$ and $\#_{\{4\}}$ already follow from Proposition \ref{prop.alg_up}. Now we prove the statement about $\#_{\{33A\}}$. First we prove the following claim: if $A$ is a component of $\Mb_D$ that contains a cycle, and no atom in $A$ belongs to a \{33A\} molecule, then at least one atom in $A$ must be deg 2 when it is cut. To prove this claim, choose a cycle $\theta$ within $A$, we may think of $\theta$ as formed by ov-segments, and consider the endpoints of these ov-segments. Let $\nf$ be the one among these endpoints that is cut \emph{after all the other endpoints}, then $\nf$ is ov-adjacent to other endpoints which are cut before $\nf$, which creates two fixed ends at $\nf$, so $\nf$ must be deg 2 when it is cut.

Now, if $A$ is as above, and suppose $\nf\in A$ is the first atom in $A$ that has deg 2 when it is cut, and consider the cutting operation that turns $\nf$ into deg 2. Then this cutting operation must involve an atom $\pf$ that is ov-adjacent to $\nf$; note that this $\pf$ belongs to either $A$ or $\Mb_U$. If $\pf\in A$, then $\pf$ cannot be deg 2 when it is cut (as $\nf$ is the first); this $\pf$ also cannot be deg 4 when it is cut (because we cut deg 4 only when every atom has deg 4, and cutting a deg 4 atom $\pf$ cannot turn the deg 4 atom $\nf$ into deg 2, cf. Proposition \ref{prop.mol_axiom} (\ref{it.axiom1})), so it must be deg 3 when it is cut, but $\nf$ also must has deg 3 when $\pf$ is cut, which violates Definition \ref{def.alg_up} (\ref{it.alg_up_3}).

Therefore we must have $\pf\in\Mb_U$. Since $\pf$ is cut before $\nf$, it cannot be cut in Definition \ref{def.alg_up} (\ref{it.alg_up_1}) (as cutting any deg 2 C-atom in Definition \ref{def.alg_up} (\ref{it.alg_up_1}) requires that this atom has two bottom fixed ends, see also Remark \ref{rem.reg}), nor can it have deg 4 when it is cut (for the same reason as in the last paragraph). Suppose $\pf$ belongs to $S_\qf$ for some $\qf$ as in  Definition \ref{def.alg_up} (\ref{it.alg_up_3}), since $\Mb_U$ is a forest, we know that $\pf$ has at most one parent in $S_\qf$ (otherwise we have two different paths going from $\pf$ to $\qf$ that form a cycle in $\Mb_U$), so it must have deg 3 when it is cut (the top edge at $\pf$ corresponding to the non $S_\qf$ parent, as well as the two bottom edges at $\pf$, are not fixed when $\pf$ is cut; see {\color{blue}Figure \ref{fig.up_aux}}), while at the same time $\nf$ also has deg 3. By Definition \ref{def.alg_up} (\ref{it.alg_up_3}), we then know that $\pf$ belongs to a \{33A\} molecule.

As a conclusion, we know that each component $A$ of $\Mb_D$ that contains a cycle, must either contain one atom in a \{33A\} molecule, or intersect a particle line of some atom in some \{33A\} molecule. Clearly each \{33A\} molecule can be so obtained by at most 4 components (as each particle line can intersect at most one component of $\Mb_D$, see Remark \ref{rem.layer_interval}), so we get $\#_{\{33A\}}\geq \#_{\mathrm{cyc.comp}(\Mb_D)}/4$.
\end{proof}
We are now ready to prove Proposition \ref{prop.comb_est}, in the case when $\Rf\geq (C_{10}^*)^{-1}\rho$ (discussed in Section \ref{sec.reduce3}).
\begin{proposition}
\label{prop.comb_est_case3}
Proposition \ref{prop.comb_est} is true, if the total number of cyclic clusters $\Rf\geq (C_{10}^*)^{-1}\rho$.
\end{proposition}
\begin{proof} In this proof we may assume that $|H|\leq (C_{12}^*)^{-1}\rho$, and restrict that $\Mb$ contains at most $(C_{12}^*)^{-1}\rho$ disjoint strongly degenerate primitive pairs (by inserting the indicator function $1-\mathbbm{1}_{\mathrm{str.deg.prm}}$ opposite to the one in Proposition \ref{prop.case2}). In fact, the case when  $|H|\geq (C_{12}^*)^{-1}\rho$ follows from Proposition \ref{prop.case1}, and the case of many strongly degenerate primitive pairs follows from Proposition \ref{prop.case2}. Recall our goal is to prove (\ref{eq.overall_alg_3}) (see Reduction \ref{red.1}), which we repeat here as
\begin{equation}\label{eq.overall_alg_2_copy} (\upsilon/2)\cdot \#_{\mathrm{good}}-d\cdot\#_{\{4\}}\geq (C_{11}^*)^{-1}\cdot \rho,
\end{equation} where we also treat all \{33A\} molecules as good.

In the proof, we fix a layer $\ell_1$ such that $\Rf_{\ell_1}\gtrsim \Rf\gtrsim\rho$ and $\Rf_{\ell_1}\gg\Rf_{\ell'}$ for each $\ell'>\ell_1$ (see (\ref{eq.comb_est_case31}) below), then the number of cyclic cluster in $\Mb_{\ell_1}$ is at least $\geq \Gamma^{-1}\cdot\Rf_{\ell_1}\gtrsim\rho$ as each cyclic cluster contributes at most $\Gamma$ to $\rho(\Mb_{\ell'})$. We could exploit the $\gtrsim\rho$ cyclic clusters in $\Mb_{\ell_1}$ using Proposition \ref{prop.alg_up} (\ref{it.up_proof_4}), \emph{except} when these cyclic clusters belong to many different components in $\Mb_{\ell_1}$. In the latter case, we start from $\Mb_{\ell_1}$ and keep going up to include more layers, until these cyclic clusters are almost connected to each other (in the sense that the number of components containing them becomes negligible). Here we need to perform layer refining for layers $>\ell_1$ the number of deleted atoms being bounded by $\Rf_{\ell'}$ for $\ell'>\ell_1$ and thus negligible), so that the new thin layer we include at each time is a forest. Then, at the layer of transition, the relevant molecule will be a union $\Mb_U\cup\Mb_D$, where $\Mb_D$ is the union of thin layers before the transition, and $\Mb_U$ is the newly included thin layer (which is a forest). We are then in the situation of Proposition \ref{prop.alg_up_recl}, so applying it leads to an acceptable upper bound on $\#_{\{4\}}$ and lower bound on $\#_{\{33\}}$.

To start, we hoose a maximal $\ell_1\in[1:\ell]$ (which always exists because $\Rf\geq (C_{10}^*)^{-1}\rho$), such that
\begin{equation}\label{eq.comb_est_case31}\Rf_{\ell_1}>(C_{6}^*)^{-\ell_1}\cdot (C_{10}^*)^{-1}\rho;\qquad \Rf_{\ell'}\leq (C_{6}^*)^{-\ell'}\cdot (C_{10}^*)^{-1}\rho< (C_{6}^*)^{-1}\cdot \Rf_{\ell_1},\quad\forall \ell_1<\ell'\leq\ell.
\end{equation} Note that, since each cyclic cluster $W\subseteq\Mb_{\ell_1}$ satisfies $\rho(W)\leq\Gamma$ and thus contribute at most $\Gamma$ to $\rho(\Mb_{\ell_1})$, and joining the clusters by O-atoms does not create new independent cycles (and thus does not increase $\rho(\Mb_{\ell_1})$), we know that the number of cyclic clusters in $\Mb_{\ell_1}$ is at least $\Gamma^{-1}\cdot\Rf_{\ell_1}$.

Consider all the components of $\Mb_{\ell_1}$ (i.e. layer $\ell_1$) that contain at least one cyclic cluster, let them be $\Xb_j\,(1\leq j\leq N)$ where $N$ is the total number of such components, and let their union be $\Xb$. If $N< (C_{5}^*)^{-1}\cdot\Rf_{\ell_1}$ (which includes when $\ell_1=\ell$, because in this case $N\leq |H|$ by repeating the (PL) argument in the proof of Proposition \ref{prop.case1}, and note also that $C_{12}^*\gg C_{10}^*$), then we 
\begin{enumerate}[{(i)}]
\item Cut $\Xb$ as free, then cut $\Mb_{>\ell_1}$ as free from $\Mb\backslash\Xb$ and cut it into elementary molecules using \textbf{UP};
\item Then cut $\Mb_{<\ell_1}\cup (\Mb_{\ell_1}\backslash \Xb)$ into elementary molecules using \textbf{DOWN}.
\end{enumerate} It is easy to see that this can always be done, and the contribution of these steps to $\#_{\{4\}}$ is at most $2|H|$ by repeating the (PL) argument as above (this allows us to bound the number of components of $\Mb_{>\ell_1}$; moreover, the number of \emph{full} components of $\Mb_{<\ell_1}\cup (\Mb_{\ell_1}\backslash \Xb)$ after cutting $X$ and $\Mb_{>\ell_1}$ as free does not exceed the number of components of $\Mb$ before cutting, cf. the proof of Proposition \ref{prop.case2}). Then, to cut $\Xb$ into elementary molecules, we argue as \textbf{Choice 2} in the proof of Proposition \ref{prop.case2}; however instead of the $\Qb_j$ in the proof of Proposition \ref{prop.case2}, we consider all the clusters forming each component of $\Xb$. The same proof of Proposition \ref{prop.case2} implies that each cyclic cluster provides one good molecule, so for the contribution of the cuttings involving $\Xb$ we have $\#_{\mathrm{good}}\geq \Gamma^{-1}\cdot \Rf_{\ell_1}$ and $\#_{\{4\}}\leq N$, which easily implies (\ref{eq.overall_alg_2_copy}), in view of $N< (C_{5}^*)^{-1}\cdot\Rf_{\ell_1}$, and (\ref{eq.comb_est_case31}) and that $|H|\leq (C_{12}^*)^{-1}\cdot\rho$.

Now assume $N\geq (C_{5}^*)^{-1}\Rf_{\ell_1}$ (in particular we must have $\ell_1<\ell$). We then perform layer refining (Definition \ref{def.layer_refine}) for each layer $\ell_1<\ell'\leq\ell$. After that, for each thin layer $\zeta\geq(\ell_1,0)$ (where we identify the unrefined layer $\ell_1$ with $(\ell_1,0)$), define
\begin{equation}\label{eq.comb_est_case32} N_\zeta:=\textrm{number of components of $\Mb_{[(\ell_1,0):\zeta]}^T$ that contains some $\Xb_j\,(1\leq j\leq N)$}.
\end{equation} Note that for $\zeta=(\ell_1,0)$ we have $N_\zeta=N$, and for $\zeta$ being the highest thin layer of $\Mb$, we have \begin{equation}\label{eq.comb_est_case32.5}N_\zeta\leq |H|+(2\Gamma+4)\sum_{\ell'>\ell_1}\Rf_{\ell'}<C_{2}^*(C_{6}^*)^{-1}\cdot \Rf_{\ell_1};\end{equation}this is proved by repeating the (PL) argument in the proof of Proposition \ref{prop.case1} but using Proposition \ref{prop.layer_refine_2} (\ref{it.layer_refine_22}), (\ref{it.layer_refine_25}) and (\ref{it.layer_refine_26}). In fact, start from any component $Z$ in $\Mb_{[(\ell_1,0):\zeta]}^T$; by dividing each layer $\ell'>\ell_1$ into thin layers and using Proposition \ref{prop.layer_refine_2} (\ref{it.layer_refine_22}), with at most $(\Gamma+1)\sum_{\ell'>\ell_1}\Rf_{\ell'}$ exceptions, we may assume $Z$ contains a component of $\Mb_{(\ell',0)}^T$ for some $\ell'\geq\ell_1$ (note that the layer $\ell_1$ or $(\ell_1,0)$ is unrefined). Then by Proposition \ref{prop.layer_refine_2} (\ref{it.layer_refine_25}) and Proposition \ref{prop.mol_axiom} (\ref{it.axiom4}), this $Z$ must intersect a particle line $\pb\in r(\Mb_{\ell'})\subseteq p(\Mb_{\ell'+1})$. By Proposition \ref{prop.layer_refine_2} (\ref{it.layer_refine_26}), with at most $2\Rf_{\ell'+1}$ exceptions, this $\pb$ must either intersect a component of $\Mb_{(\ell'+1),0}^T$ (in which case this component must be contained in $Z$, as $\pb$ can intersect at most one component of $\Mb_{[(\ell_1,0):\zeta]}^T$, see Remark \ref{rem.layer_interval}), or intersect a component of $\Mb_{(\ell'+1,k)}^T\,(k\geq 1)$ (in which case this component must be contained in $Z$, leading to at most $(\Gamma+1)\Rf_{\ell'+1}$ choices for $Z$ due to Proposition \ref{prop.layer_refine_2} (\ref{it.layer_refine_22})), or belong to $r(\Mb_{\ell'+1})$. In any case, by excluding the $2\Rf_{\ell'+1}+(\Gamma+1)\Rf_{\ell'+1}$ exceptions, this $Z$ must intersect a particle line $\pb'\in r(\Mb_{\ell'+1})$ (if $Z$ contains a component of $\Mb_{(\ell'+1),0}^T$, then apply Proposition \ref{prop.layer_refine_2} (\ref{it.layer_refine_25}) again). Repeating this process until reaching $r(\Mb_\ell)=H$ (at which time $Z$ will intersect one ofthe $|H|$ particle lines, leading to at most $|H|$ choices) will then lead to at most
\[|H|+(\Gamma+1)\sum_{\ell'>\ell_1}\Rf_{\ell'}+\sum_{\ell'\geq\ell_1}\big(2\Rf_{\ell'+1}+(\Gamma+1)\Rf_{\ell'+1}\big)\] choices for $Z$, which proves (\ref{eq.comb_est_case32.5}).

Now, using (\ref{eq.comb_est_case32.5}), we see that there exists a (maximal) thin layer $\zeta\geq(\ell_1,0)$, which is not the highest thin layer of $\Mb$, such that
\begin{equation}\label{eq.comb_est_case33} N_\zeta\geq (C_{5}^*)^{n((\ell_1,0))-n(\zeta)-1}\cdot \Rf_{\ell_1};\qquad N_{\zeta'}<(C_{5}^*)^{n((\ell_1,0))-n(\zeta')-1}\cdot \Rf_{\ell_1}\leq (C_{5}^*)^{-1}\cdot N_\zeta,\quad\forall \zeta'>\zeta,
\end{equation} where $n(\zeta)$ is the number of thin layers below or equal to $\zeta$; see {\color{blue}Figure \ref{fig.many_cyclic}}. Let $\zeta^+$ be the lowest thin layer above $\zeta$, then define $\Yb$ to be the union of all components of $\Mb_{[(\ell_1,0):\zeta^+]}^T$ that contains some $\Xb_j\,(1\leq j\leq N)$, and write it as $\Yb=\Yb_U\cup\Yb_D$, where $\Yb_U=\Yb\cap \Mb_{\zeta^+}^T$ and $\Yb_D=\Yb\cap\Mb_{[(\ell_1,0):\zeta]}^T$. It is easy to verify that the assumptions of Proposition \ref{prop.alg_up_recl} are satisfied, and that
\begin{equation}\label{eq.comb_est_case34} \#_{\mathrm{comp}(\Yb)}\leq N_{\zeta^+},\qquad \#_{\mathrm{cyc.comp}(\Yb_D)}\geq N_\zeta
\end{equation} (note that each $\Xb_j$ already contains a cycle, and so does any component that contains some $\Xb_j$; therefore each of the $N_\zeta$ components of $\Mb_{[(\ell_1,0):\zeta]}^T$ containing some $\Xb_j$ will contribute to $\#_{\mathrm{cyc.comp}(\Yb_D)}$).
\begin{figure}[h!]
    \centering
    \includegraphics[width=0.63\linewidth]{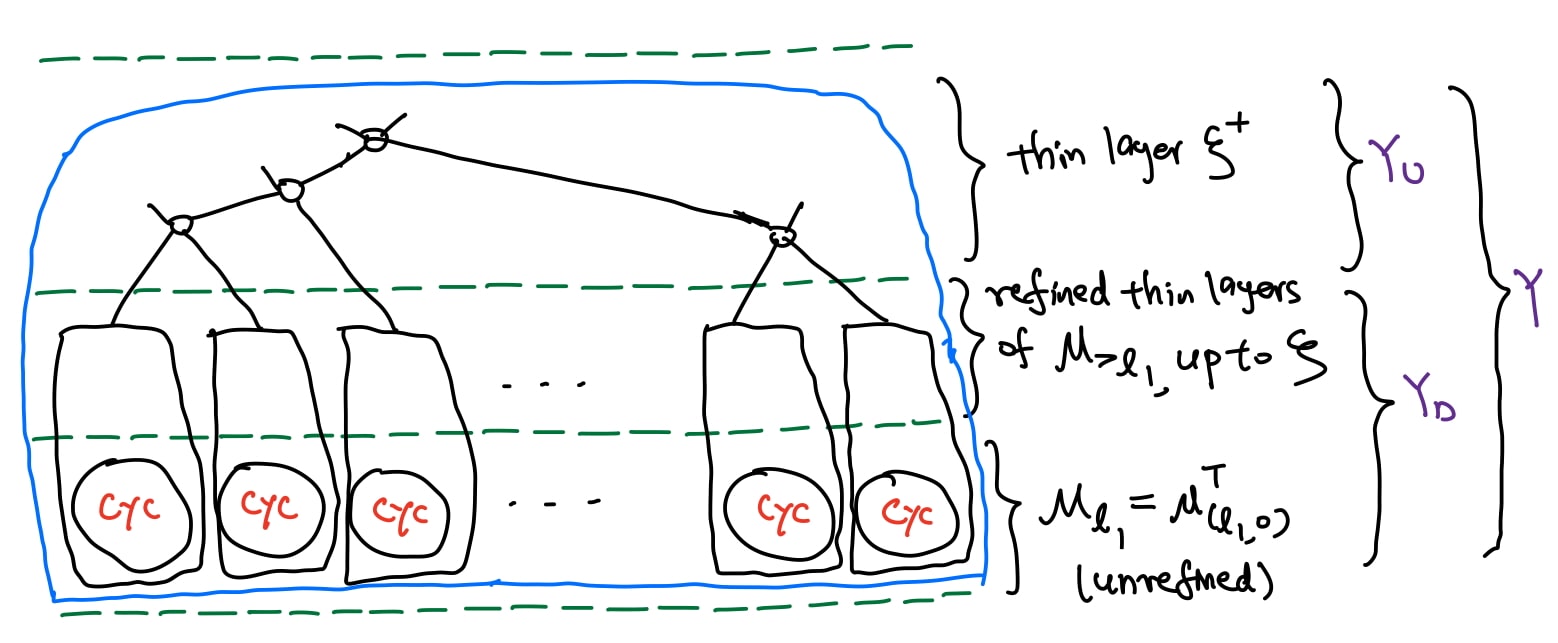}
    \caption{The proof of Proposition \ref{prop.comb_est_case3}, where we start from the cyclic components of the unrefined layer $\Mb_{\ell_1}$, and keep going up (adding one thin layer each time) until the first thin layer after which these components are connected together (with a negligible number of resulting components if we account for all thin layers). Define $\Yb=\Yb_U\cup\Yb_D$ as in the figure, then we can apply Proposition \ref{prop.alg_up_recl} to get many \{33\} molecules.}
    \label{fig.many_cyclic}
\end{figure}

Then, as above, we 
\begin{enumerate}[{(i)}]
\item Cut $\Yb$ as free, then cut $\Mb_{>\zeta^+}^T$ as free from $\Mb\backslash\Yb$  and subsequently cut it using \textbf{UP};
\item Then cut $\Mb_{<(\ell_1,0)}^T\cup (\Mb_{[(\ell_1,0):\zeta^+]}^T\backslash \Yb)$ into elementary molecules using \textbf{DOWN}.
\end{enumerate} It is easy to see that this can always be done (as there is no bond between $\Yb$ and $\Mb_{[(\ell_1,0):\zeta^+]}^T\backslash \Yb$), and the contribution of these steps to $\#_{\{4\}}$ is at most 
\begin{equation}\label{eq.comb_est_case35}
2|H|+8(\Gamma+1)\sum_{\ell'>\ell_1}\Rf_{\ell'}<C_{2}^*(C_{6}^*)^{-1}\cdot \Rf_{\ell_1};
\end{equation}
this follows from the upper bounds on the number of components of $\Mb_{>\zeta^+}^T$ and $\Mb$, which is proved by repeating the above (PL) arguments (as in the proof of (\ref{eq.comb_est_case32.5})) which uses Proposition \ref{prop.layer_refine_2} (\ref{it.layer_refine_22}), (\ref{it.layer_refine_25}) and (\ref{it.layer_refine_26}).
Finally we cut $\Yb$ into elementary molecules using \textbf{UP}, and by Proposition \ref{prop.alg_up_recl}, for the contribution of the cuttings involving $\Yb$, we have \begin{equation}\label{eq.case35.5}\#_{\{33A\}}\geq N_\zeta/4,\qquad\#_{\{4\}}\leq N_{\zeta^+}.
\end{equation} By putting together (\ref{eq.comb_est_case31}), (\ref{eq.comb_est_case33}), (\ref{eq.comb_est_case35}) and (\ref{eq.case35.5}), and recall that all \{33A\} molecules are treated as good (Reduction \ref{red.1}), we get
\begin{equation}\label{eq.comb_est_case37}
(\upsilon/2)\cdot(\#_{\mathrm{good}})-d\cdot(\#_{\{4\}})\geq (\upsilon/8)\cdot N_\zeta-d\cdot (N_{\zeta^+}+C_2^*(C_6^*)^{-1}\cdot\Rf_{\ell_1})
\geq (\upsilon/16)\cdot N_\zeta\geq (C_{10}^*)^{-2}\cdot\rho,
\end{equation} which proves (\ref{eq.overall_alg_2_copy}). The proof is now complete.
\end{proof}
\subsection{Layer selection and weak degeneracies}\label{sec.layer_select} From now on we will focus on proving (\ref{eq.overall_alg_3}) under the assumptions in Reduction \ref{red.1}. The goal of this subsection is to reduce this result to Proposition \ref{prop.case5} (and prove Proposition \ref{prop.comb_est_case4} which contains the simple cases of initial cumulant and weak degeneracies), by applying the procedure of \textbf{layer selection} (discussed in Section \ref{sec.toy_multi}).

The full layer selection process, defined in Definition \ref{def.layer_select} below, is essentially the same as the toy version in Definition \ref{def.layer_select_toy} in Section \ref{sec.toy}. The only difference is that we need to perform layer refining for layers $\Mb_{\ell'}\,(\ell'\leq\ell_1)$, which leads to an extra step of picking a thin layer $\zeta_U$ within layer $\Mb_{\ell_1}$, see Definition \ref{def.layer_select} (\ref{it.layer_select_2}). The rest is the same as in Definition \ref{def.layer_select_toy}, except that $v_\zeta^*$ in Definition \ref{def.layer_select} (\ref{it.layer_select_2}) is defined for thin layers instead of layers, and the $\zeta_D$ we pick is also a thin layer (following the same rule as in Definition \ref{def.layer_select_toy}).
\begin{definition}[Layer selection]
\label{def.layer_select} Let $\Mb$ be a molecule as in Proposition \ref{prop.mol_axiom}. Assume that $|H|\leq (C_{12}^*)^{-1}\rho$ and $\Rf\leq (C_{10}^*)^{-1}\rho$ as in Reduction \ref{red.1}. Then we will perform the following steps to determine the values of two \textbf{thin layers} $\zeta_U$ and $\zeta_D$ (in one case it is possible that $\zeta_D$ is absent and we only determine $\zeta_U$), see {\color{blue}Figure \ref{fig.layerselect_full}}. Recall that $r(\Mb_{\ell})=H$ and $r(\Mb_0)=H_0$, and that we denote $s_{\ell}=|r(\Mb_{\ell'})|$.
\begin{figure}[h!]
    \centering
    \includegraphics[width=0.55\linewidth]{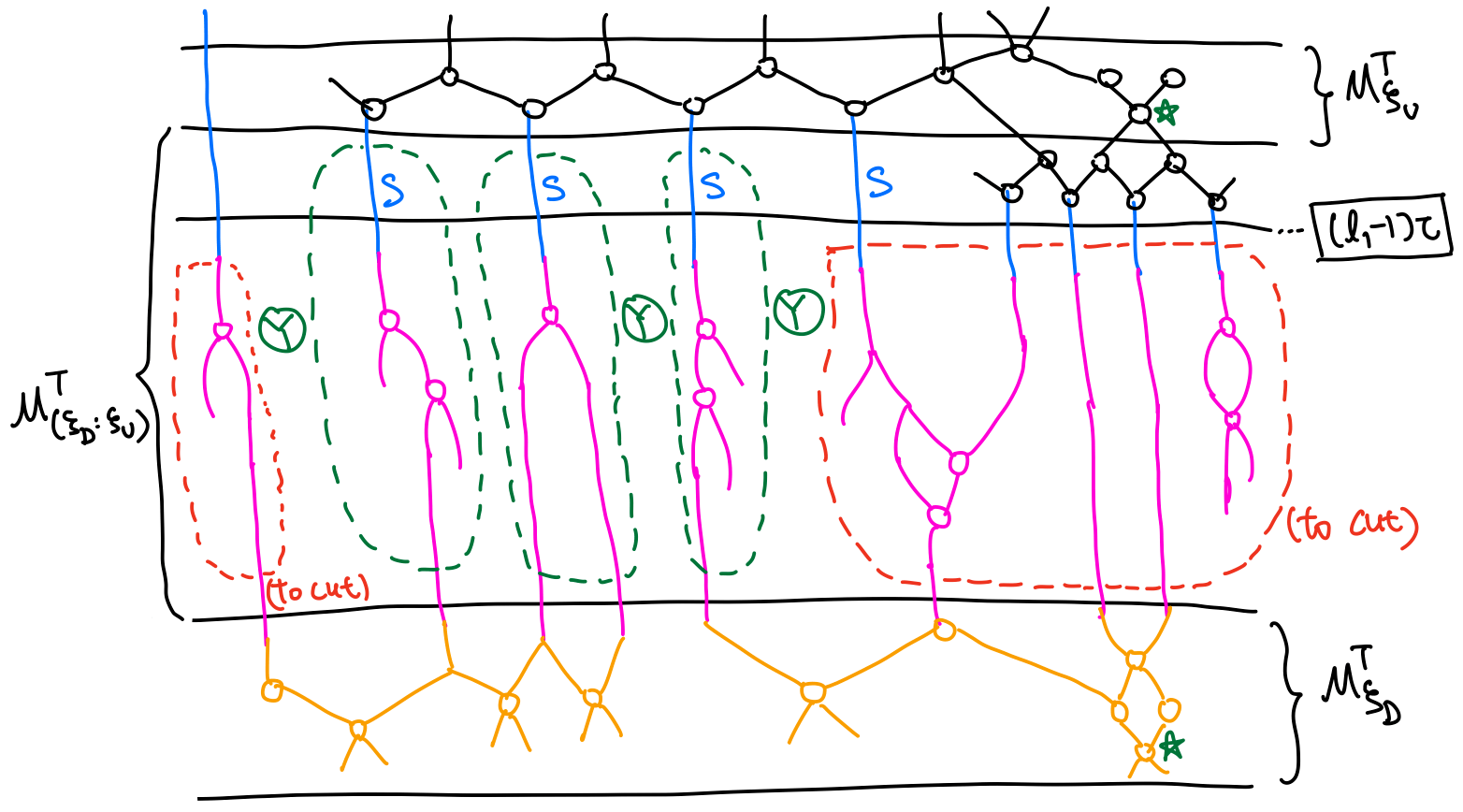}
    \caption{The full layer selection algorithm (Definition \ref{def.layer_select}). Compared to the toy case (Definition \ref{def.layer_select_toy}, {\color{blue}Figure \ref{fig.layerselect}}), we have an extra step of layer refining and selecting $\zeta_U$. Here the blue particle lines are in $r(\Mb_{\ell_1-1})$, but we only care for those in $S$ (i.e. those that have atoms in the thin layer $\Mb_{\zeta_U}^T$ but not in lower thin layers within $\Mb_{\ell_1}$).
    \newline
        \-\hspace{8pt} Another difference is that there are exceptional components in $\Mb_{(\zeta_D:\zeta_U)}^T$ (the red ones, compared to the green ones in $Y$) that contain cycles, or intersect a particle line not in $S$, or intersect two particle lines in $S$. These will be cut in Definition \ref{def.layer_cutting}, \textbf{Option 1} below.}
    \label{fig.layerselect_full}
\end{figure}
\begin{enumerate}
\item\label{it.layer_select_1} Choose $1\leq \ell_1\leq \ell$ such that
\begin{equation}\label{eq.layer_select_1}s_{\ell_1-1}\geq (C_{7}^*)^{-(\ell_1-1)}\cdot\rho;\qquad s_{\ell'}<(C_{7}^*)^{-\ell'}\cdot\rho\leq (C_{7}^*)^{-1}
\cdot s_{\ell_1-1},\quad\forall \ell_1\leq\ell'\leq\ell,\end{equation} and refine layer $\Mb_{\ell_1}$ as in Definition \ref{def.layer_refine}. Define also $s_{\ell_1-1}:=\rho'$.
\item\label{it.layer_select_2} After refining, consider all particle lines $\pb\in r(\Mb_{\ell_1-1})$ (there are $\rho'$ of them), for any $0\leq\gamma\leq \Gamma+1$, define $u_\gamma^*$ to be the number of these $\rho'$ particle lines that \emph{do not} contain any atom in $\Mb_{(\ell_1,k)}^T$ for any $k\leq\gamma$ (denote $u_{-1}^*:=s_{\ell_1-1}$). Then choose $0\leq\gamma\leq\Gamma+1$ such that
\begin{equation}\label{eq.layer_select_2}u_{\gamma-1}^*\geq (C_{6}^*)^{-\gamma}\cdot \rho',\quad u_\gamma^*<(C_{6}^*)^{-\gamma-1}\cdot \rho'\leq (C_{6}^*)^{-1}\cdot u_{\gamma-1}^*.\end{equation} Fix \textbf{the thin layer $\zeta_U:=(\ell_1,\gamma)$}. Define also $u_{\gamma-1}^*-u_\gamma^*:=\rho''$. If $\ell_1=1$, then we stop the process with $\zeta_D$ absent.
\item\label{it.layer_select_3} If $\ell_1>1$, refine the layers $\Mb_{\ell'}$ for all $\ell'\leq\ell_1-1$, and consider all the particle lines $\pb\in r(\Mb_{\ell_1-1})$ that does not contain any atom in $\Mb_{(\ell_1,k)}^T$ for any $k\leq\gamma-1$ but does contain an atom in $\Mb_{(\ell_1,\gamma)}^T$; there are $\rho''$ of them, and define the set of them to be $S$. Let $\zeta_0$ be the highest thin layer below $\zeta_U$. For each thin layer $\zeta\leq\zeta_0$, define $v_\zeta^*$ to be the number of particle lines in $S$ that are either \emph{connected to a cycle} via $\Mb_{[\zeta:\zeta_0]}^T$, or \emph{are connected to another particle line} $\pb'\in r(\Mb_{\ell_1-1})$ via $\Mb_{[\zeta:\zeta_0]}^T$ (this $\pb'$ need not be in $S$). Now, if for the lowest thin layer $(1,0)$ we have $v_{(1,0)}^*\leq (C_{5}^*)^{-1}\cdot\rho''$, then stop the process with $\zeta_D$ absent; otherwise \textbf{choose a thin layer} $\zeta=\zeta_D$ such that
\begin{equation}\label{eq.layer_select_3}v_{\zeta_D}^*\geq (C_{5}^*)^{-n(\zeta_D)}\cdot\rho'',\qquad v_{\zeta}^*<(C_{5}^*)^{-n(\zeta)}\cdot \rho''\leq(C_{5}^*)^{-1}\cdot v_{\zeta_D}^*,\quad \forall \zeta_D<\zeta\leq\zeta_0,\end{equation} where $n(\zeta)$ is the number of thin layers below or equal to $\zeta$. Define also $v_{\zeta_D}^*:=\rho'''$. In the case when $\zeta_D$ is absent, we also define $\rho''':=\rho''$.
\end{enumerate}
\end{definition}
In Proposition \ref{prop.layer_select} below we state and prove several properties of layer selection. That this process is well-defined is proved in the same way as in Proposition \ref{prop.layer_select_toy} (except for the extra step in Definition \ref{def.layer_select} (\ref{it.layer_select_2})). Additionally, this proposition contains the upper bound on the number of components of various subsets of $\Mb$, which will be used to bound $\#_{\{4\}}$ for the cutting procedures defined in Definition \ref{def.layer_cutting}, see Proposition \ref{prop.layer_cutting} (\ref{it.layer_cutting_1}); most important among them is the set $\Mb_{(\zeta_D:\zeta_U)}^T\backslash Y$, which precisely contains the \emph{exceptional components} as in Definition \ref{def.layer_select_toy}. The proof of these upper bounds follows from repeating the post-refining (PL) argument using Proposition \ref{prop.layer_refine_2}, similar to the proof of (\ref{eq.comb_est_case32.5}) in Proposition \ref{prop.comb_est_case3}.
\begin{proposition}
\label{prop.layer_select} The construction process in Definition \ref{def.layer_select} is well-defined. Moreover, in Definition \ref{def.layer_select}, consider all the connected components of $\Mb_{(\zeta_D:\zeta_U)}^T$ (or $\Mb_{<\zeta_U}^T$ if $\zeta_D$ is absent) that (i) intersect exactly one particle line in $S$ (as defined in Definition \ref{def.layer_select} (\ref{it.layer_select_3})), and (ii) does not intersect any other particle line in $r(\Mb_{\ell_1-1})$, and (iii) does not contain a cycle. Denote the union of these components by $Y$. Then, each of the following quantities is bounded by $10(C_{5}^*)^{-1}\cdot\rho'''$:
\begin{enumerate}
\item\label{it.layer_select_new1} The number of components of $\Mb$, $\Mb_{>\zeta_U}^T$, $\Mb_{>\zeta_D}^T$, $\Mb_{\zeta_U}^T$, and $\Mb_{[\zeta_D:\zeta_U]}^T$ and $\Mb_{(\zeta_D:\zeta_U]}^T$;
\item\label{it.layer_select_new2} The number of components of $\Mb_{\zeta_D}^T$ that are not connected to $\Mb_{\zeta_U}^T\cup\Mb_{(\zeta_D:\zeta_U)}^T$ by a bond.
\item\label{it.layer_select_new3} The number of components of $\Mb_{(\zeta_D:\zeta_U)}^T\backslash Y$;
\end{enumerate}
\end{proposition}
\begin{proof} Recall from (\ref{eq.def_rho_old}) that $\rho=\Rf+\sum_{\ell'\geq 0}s_{\ell'}$; since $\Rf\leq (C_{10}^*)^{-1}\rho$, there exists at least one $\ell'$ (including $0$) such that $s_{\ell'}\geq (C_{7}^*)^{-\ell'}\cdot\rho$. Using also that $s_\ell=|H|\leq (C_{12}^*)^{-1}\rho$, we must have $\ell'<\ell$, and may then define the largest such $\ell'$ to be $\ell_1-1$ in in Definition \ref{def.layer_select} (\ref{it.layer_select_1}), so (\ref{eq.layer_select_1}) holds. 

Next, consider $u_\gamma^*$ defined in Definition \ref{def.layer_select} (\ref{it.layer_select_2}); by Proposition \ref{prop.layer_refine_2} (\ref{it.layer_refine_26}) we know the number of particle lines in $r(\Mb_{\ell_1-1})$ that do not contain any atom in $\Mb_{(\ell_1,k)}^T$ for any $k$, is bounded by \[u_{\Gamma+1}^*\leq 2\Rf_{\ell_1}+s_{\ell_1}<\big(2 (C_{10}^*)^{-1}(C_{7}^*)^{\Lf}+(C_{7}^*)^{-1}\big)\cdot\rho'<(C_{6}^*)^{-10\Gamma}\cdot\rho'\] using (\ref{eq.layer_select_1}). Since also $u_{-1}^*=\rho'$, we know that there exists a largest value $\gamma-1$ satisfying $u_{\gamma-1}^*\geq (C_{6}^*)^{-\gamma}\cdot\rho'$ in Definition \ref{def.layer_select} (\ref{it.layer_select_2}), so (\ref{eq.layer_select_2}) holds.

Next, if $\ell_1>1$, consider $v_\zeta^*$ defined in Definition \ref{def.layer_select} (\ref{it.layer_select_3}). Let $\zeta_1$ be the highest thin layer within $\Mb_{\ell_1-1}$. For the lowest thin layer above $\zeta_1$ (which is $(\ell_1,0)$) we obviously have $v_{(\ell_1,0)}^*=0$ (by definition of $S$, these particle lines do not intersect any thin layer between $(\ell_1,0)$ and $\zeta_0$); now if $v_{(1,0)}^*\geq  (C_{5}^*)^{-1}\cdot\rho''$, we know that there exists a largest value $\zeta_D$ satisfying $v_{\zeta_D}^*\geq (C_{5}^*)^{-n(\zeta_D)}\cdot\rho''$ in Definition \ref{def.layer_select} (\ref{it.layer_select_3}), so (\ref{eq.layer_select_3}) holds. This shows that the thin layers $\zeta_U$ and $\zeta_D$ are well defined and $\zeta_D\leq\zeta_1$.

Now we bound the number of components as stated in (\ref{it.layer_select_new1})--(\ref{it.layer_select_new3}). We will consider the case where $\zeta_D$ exists; the case where $\zeta_D$ is absent only requires trivial modifications. In fact, by repeating the (PL) argument using Proposition \ref{prop.layer_refine_2} (\ref{it.layer_refine_22}), (\ref{it.layer_refine_25}) and (\ref{it.layer_refine_26}) as in the proof of (\ref{eq.comb_est_case32.5}) in Proposition \ref{prop.comb_est_case3}, we can bound the number of components of $\Mb$, $\Mb_{>\zeta_U}^T$ and $\Mb_{>\zeta_D}^T$by $|H|+(2\Gamma+4)\Rf\leq 2 (C_{10}^*)^{-1}\rho$. The number of components of $\Mb_{\zeta_U}^T$ (or of any union of thin layers within layer $\ell_1$) is also bounded by $s_{\ell_1}+(\Gamma+1)\Rf\leq 2(C_{7}^*)^{-1}\cdot\rho'$ by Proposition \ref{prop.layer_refine_2} (\ref{it.layer_refine_22})--(\ref{it.layer_refine_23}) and Proposition \ref{prop.mol_axiom} (\ref{it.axiom5}). This establishes the first four estimates in (\ref{it.layer_select_new1}).

In addition, consider each component $X$ of $\Mb_{(\zeta_D:\zeta_U)}^T$ or $\Mb_{[\zeta_D:\zeta_U]}^T$ or $\Mb_{(\zeta_D:\zeta_U]}^T$. If it intersects some thin layer $\Mb_{(\ell_1,k)}^T$, then the number of such components is bounded by $s_{\ell_1}+(\Gamma+1)\Rf\leq 2(C_{7}^*)^{-1}\cdot\rho'$ as above. If not, it must intersect some thin layer $\Mb_{(\ell',k)}^T$ for $\ell'\leq \ell_1-1$, which allows us to repeat the (PL) argument using Proposition \ref{prop.layer_refine_2} (\ref{it.layer_refine_22}), (\ref{it.layer_refine_25}) and (\ref{it.layer_refine_26}). In summary we see that, with at most $4(C_{7}^*)^{-1}\cdot\rho'$ exceptions, each component $X$ must intersect a particle line $\pb\in r(\Mb_{\ell_1-1})$. Next, with at most $u_\gamma^*<2(C_{6}^*)^{-1}\cdot\rho''$ exceptions (and the total number of exceptions $\ll (C_{5}^*)^{-2}\cdot\rho'''$), we may assume this $\pb$ intersects $\Mb_{(\ell_1,k)}^T$ for some $k\leq\gamma$. Now, if we are considering $\Mb_{[\zeta_D:\zeta_U]}^T$ or $\Mb_{(\zeta_D:\zeta_U]}^T$ which includes all thin layers $\Mb_{(\ell_1,k)}^T\,(k\leq\gamma)$, then we know that the component $X$ must intersect some $\Mb_{(\ell_1,k)}^T\,(k\leq\gamma)$, which leads to $\leq2(C_{7}^*)^{-1}\cdot\rho'$ choices for $X$ (and $\ll (C_5^*)^{-2}\cdot\rho'''$ accounting for all the above exceptions, recall again Remark \ref{rem.layer_interval} that each $\rho$ intersects at most one $X$). This establishes the last two estimates in (\ref{it.layer_select_new1}). Moreover, (\ref{it.layer_select_new2}) follows because any component of $\Mb_{\zeta_D}^T$ that is not connected to $\Mb_{\zeta_U}^T\cup\Mb_{(\zeta_D:\zeta_U)}^T$ by a bond is also a component of $\Mb_{[\zeta_D:\zeta_U]}^T$.

Finally consider (\ref{it.layer_select_new3}). As above, with $\ll (C_{5}^*)^{-2}\cdot\rho'''$ exceptions, for each component $Z$ of $\Mb_{(\zeta_D:\zeta_U)}^T$, it must intersect a particle line $\pb\in r(\Mb_{\ell_1-1})$ that intersects $\Mb_{(\ell_1,k)}^T$ for some $k\leq\gamma$. We may assume $\pb$ does not intersect $\Mb_{(\ell_1,k)}^T$ for $k<\gamma$ (otherwise $Z$ must intersect $\Mb_{(\ell_1,k)}^T$, leading to $\leq2(C_{7}^*)^{-1}\cdot\rho'$ possibilities for $Z$), so in particular $\pb\in S$, which means that $Z$ intersects some particle line $\pb\in S$. Let $\eta$ be the lowest thin layer above $\zeta_D$. By excluding at most $v_{\eta}^*\leq(C_{5}^*)^{-1}\cdot \rho'''$ exceptions, we may assume that $\pb$ is not connected to a cycle within $\Mb_{[\eta:\zeta_0]}^T$, and is not connected to another particle line $\pb'\in r(\Mb_{\ell_1-1})$ via $\Mb_{[\eta:\zeta_0]}^T$, but this means that $Z\subseteq Y$ by definition of $Y$. Putting all these together, we see that the number of components of $\Mb_{(\zeta_D:\zeta_U)}^T$ that are not contained in $Y$ is at most $2(C_{5}^*)^{-1}\cdot \rho'''$, as desired.
\end{proof}
Now, once the thin layers $\zeta_U$ and $\zeta_D$ are chosen by Definition \ref{def.layer_select}, we need to cut as free all the exceptional components in $\Mb_{(\zeta_D:\zeta_U)}^T\backslash Y$ in Definition \ref{def.layer_select_toy} (discussed Section \ref{sec.reduce5} (b)). Together with a few auxiliary operations treating other layers (and cutting components of $\Mb_{\zeta_D}^T$ that have no connections to $\Mb_{\zeta_U}^T\cup Y$ which is purely for convenience), this forms the \textbf{Option 1} in Definition \ref{def.layer_cutting} below, which is the main option to be executed in the rest of the proof. On the other hand, the \textbf{Option 2} in Definition \ref{def.layer_cutting} involves cutting $\Mb_{\zeta_D}^T$ first, which guarantees that all relevant atoms in $\Mb_D$ belong to \{3\} or \{4\} molecules, and will be used in the \emph{weakly degenerate} case (Definition \ref{def.weadeg}, Proposition \ref{prop.comb_est_case4}).
\begin{definition}\label{def.layer_cutting} Recall the layer selection process in Definition \ref{def.layer_select} and the definition of $Y$ in Proposition \ref{prop.layer_select}. Define also the set $\Vb$ to be the union of all components of $\Mb_{\zeta_D}^T$ that are connected to $\Mb_{\zeta_U}^T\cup Y$ by a bond, and let $G$ be the number of those components of $\Vb$ that are connected to $\Mb_{\zeta_U}^T\cup \Mb_{(\zeta_D:\zeta_U)}^T$ by at least two bonds. Below we define two options of cutting sequences to be applied to $\Mb$.
\begin{itemize}
\item \textbf{Option 1}. We perform the following cutting operations (note that if $\zeta_D$ is absent, then the steps here involving $\Mb_{\zeta_D}^T$ and $\Mb_{<\zeta_D}^T$ will be skipped); see {\color{blue}Figure \ref{fig.layer_cutting}}. 
\begin{figure}[h!]
    \centering
    \includegraphics[width=0.55\linewidth]{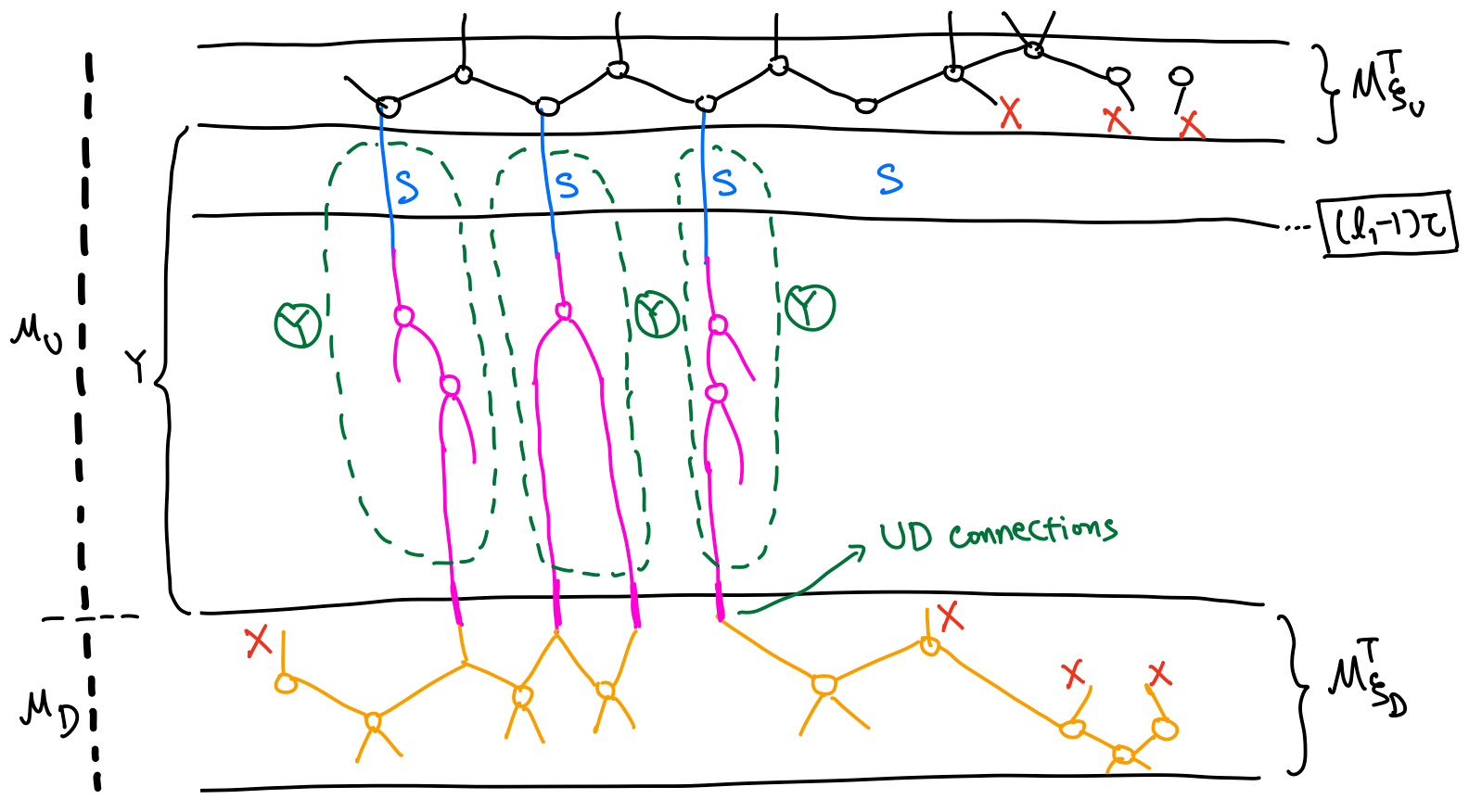}
    \caption{The result of \textbf{Option 1} in Definition \ref{def.layer_cutting}, starting from {\color{blue}Figure \ref{fig.layerselect_full}}. We cut as free all exceptional components in $\Mb_{(\zeta_D:\zeta_U)}^T\backslash Y$ (red components in {\color{blue}Figure \ref{fig.layerselect_full}}), so we are only left with $Y$ (green components). We also subsequently cut any atom that becomes deg 2 in $\Mb_{\zeta_U}^T$ or $\Mb_{\zeta_D}^T$ (the ones marked by green star in {\color{blue}Figure \ref{fig.layerselect_full}}).
    \newline
    \-\hspace{8pt} In the end $\Mb_U$ becomes a forest, while $\Mb_D$ is still a forest. The number of \{4\} molecules created is bounded using Proposition \ref{prop.layer_select}; in addition, this will create fixed ends in both $\Mb_U$ and $\Mb_D$, but only \emph{bottom} fixed ends in the former and only \emph{top} fixed ends in the latter, which is all we need.
}
    \label{fig.layer_cutting}
\end{figure}
\begin{enumerate}[{(i)}]
\item We cut $\Mb_{[\zeta_D:\zeta_U]}^T$ as free, then cut $\Mb_{>\zeta_U}^T$ as free from $\Mb_{>\zeta_U}^T\cup\Mb_{<\zeta_D}^T$ and cut it into elementary molecules using \textbf{UP}, then cut $\Mb_{<\zeta_D}^T$ into elementary molecules using \textbf{DOWN}.
\item We cut as free all the components of $\Mb_{(\zeta_D:\zeta_U)}^T$ (or $\Mb_{<\zeta_U}^T$ if $\zeta_D$ is absent) not in $Y$, and cut each of them into elementary molecules using \textbf{UP}. If any atom in $\Mb_{\zeta_U}^T$ or $\Mb_{\zeta_D}^T$ becomes deg 2 then we cut it as free, and repeat until there is no deg 2 atom in $\Mb_{\zeta_U}^T$ and $\Mb_{\zeta_D}^T$.
\item We cut as free all components of $\Mb_{\zeta_D}^T$ that are not connected to $\Mb_{\zeta_U}^T\cup Y$ by a bond, and cut each of them into elementary molecules using \textbf{DOWN}.
\end{enumerate}
After the cutting sequence in \textbf{Option 1}, define the remaining molecule to be $\Mb_{UD}$ (or $\Mb_U$ if $\zeta_D$ is absent); it contains parts of $\Mb_{\zeta_U}^T$, $\Mb_{\zeta_D}^T$ and $Y$ and will be studied in Section \ref{sec.maincr}.

\item \textbf{Option 2}. This cutting sequence applies only to the case where $\zeta_D$ exists, and cuts $\Mb$ completely into elementary molecules. We perform the following cutting operations, see {\color{blue}Figure \ref{fig.option2}}.
\begin{figure}[h!]
    \centering
    \includegraphics[width=0.55\linewidth]{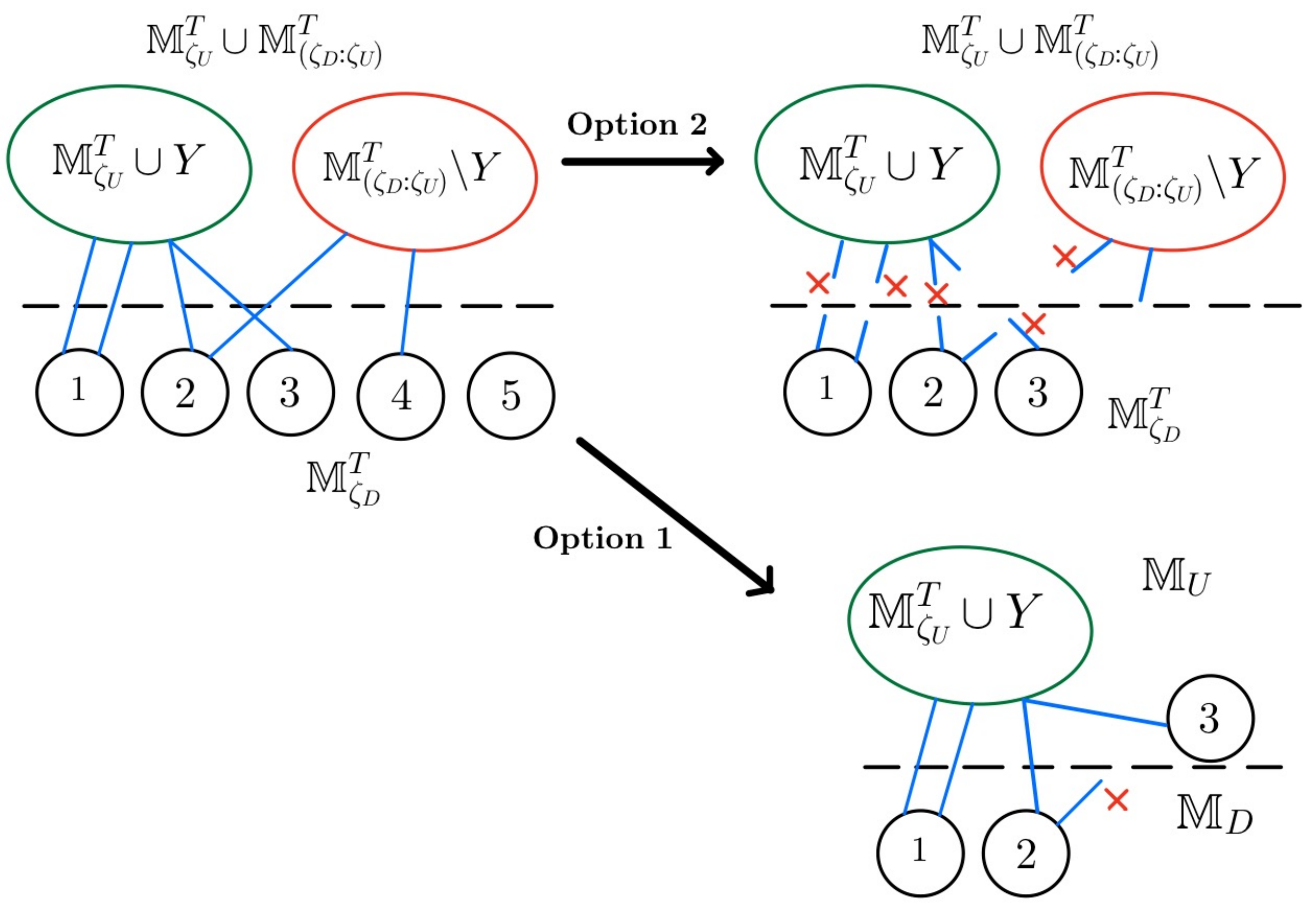}
    \caption{Illustration of \textbf{Option 2} (and \textbf{Option 1}) in Definition \ref{def.layer_cutting}. Among the connected components of $\Mb_{\zeta_D}^T$, components $4\sim 5$ \emph{do not belong to $\Vb$} and will be ignored in either option; components $1\sim 3$ belong to $\Vb$, but only $1\sim 2$ count in the $G$ components (i.e. $G=2$). By Definition \ref{def.layer_cutting}, in \textbf{Option 2}, components $1\sim 2$ are cut as free, and $3$ is cut as fixed. In \textbf{Option 1}, component $3$ will be included in $\Mb_U$, and $1\sim2$ remain in $\Mb_D$.
}
    \label{fig.option2}
\end{figure}
\begin{enumerate}[{(i)}]
\item With $\Vb$ defined as above, we cut $\Mb_{\zeta_U}^T\cup\Mb_{(\zeta_D:\zeta_U)}^T\cup \Vb$ as free, then cut $\Mb_{>\zeta_U}^T$ as free from $\Mb_{<\zeta_D}^T\cup \Mb_{>\zeta_U}^T\cup (\Mb_{\zeta_D}^T\backslash\Vb)$ and cut it into elementary molecules using \textbf{UP}, then cut $\Mb_{<\zeta_D}^T\cup (\Mb_{\zeta_D}^T\backslash\Vb)$ into elementary molecules using \textbf{DOWN}.
\item If a component of $\Vb$ has only one bond connecting to $\Mb_{\zeta_U}^T\cup \Mb_{(\zeta_D:\zeta_U)}^T$, then cut it as fixed, otherwise cut it as free. Then cut each component into elementary molecules using \textbf{DOWN}, and then cut $\Mb_{\zeta_U}^T\cup \Mb_{(\zeta_D:\zeta_U)}^T$ into elementary molecules using \textbf{UP}.
\end{enumerate}
\end{itemize}
\end{definition}
In Proposition \ref{prop.layer_cutting} we summarize a few properties of the molecule after the cutting process in Definition \ref{def.layer_cutting}. Parts (\ref{it.layer_cutting_1}) and (\ref{it.layer_cutting_2}) of Proposition \ref{prop.layer_cutting} correspond to the initial cumulant/initial link case (discussed in Section \ref{sec.reduce5} (a)) and weakly degenerate case (where we apply \textbf{Option 2}), which are easy to deal with. The main part of Proposition \ref{prop.layer_cutting} is (\ref{it.layer_cutting_3}), in which case \textbf{Option 1} reduces $\Mb$ to the 2-layer model $\Mb_{UD}:=\Mb_U\cup\Mb_D$ as in Proposition \ref{prop.layer_select_toy}, but with possible fixed ends as discussed in Section \ref{sec.reduce5} (b). For precise statements for properties of $\Mb_{UD}$, see Proposition \ref{prop.layer_cutting} (\ref{it.layer_cutting_3}); the proof is essentially the same as Proposition \ref{prop.layer_select_toy}, except the complications due to the set $S$ etc. coming from the extra step of picking $\zeta_U$.
\begin{proposition}\label{prop.layer_cutting} Recall $G$ defined in Definition \ref{def.layer_cutting}. For \textbf{Options 1 and 2} in Definition \ref{def.layer_cutting}, we have the followings.

\begin{enumerate}
\item\label{it.layer_cutting_1} The contribution of cuttings in \textbf{Option 1} to $\#_{\{4\}}$ is at most $100 (C_{5}^*)^{-1}\cdot \rho'''$. If $\zeta_D$ is absent, then see {\color{blue} Figure \ref{fig.no_zeta_D}}, after \textbf{Option 1}, we can choose at least $(\rho'''/2)$ particle lines $\pb \in S$ such that 
\begin{enumerate}
\item\label{it.layer_cutting_1i} there exists a subset $T(\pb)$ of $\Mb_U\backslash \Mb_{\zeta_U}^T=Y$ which is either empty or a component of $Y$; 
\item\label{it.layer_cutting_1ii} if $T(\pb)\neq\varnothing$ then it is a tree and has only one bond connecting to $\Mb_{\zeta_U}^T$, which is along $\pb$; 
\item\label{it.layer_cutting_1iii} consider the set $\Ec_0$ of all bottom free ends $e$ at atoms in $T(\pb)$ (or along $\pb$ if $T(\pb)=\varnothing$) for all the chosen $\pb$, then for at least  $\rho'''/4$ chosen particle lines $\pb$, there exist $e\in \Ec_0$ at an atom in $T(\pb)$ or along $\pb$, such that an initial link exists between $e$ and $e'$, where $e'\neq e$ either belongs to $\Ec_0$ or has an atom in some component already cut in \textbf{Option 1}.
\end{enumerate}
\begin{figure}[h!]
    \centering
    \includegraphics[width=0.3\linewidth]{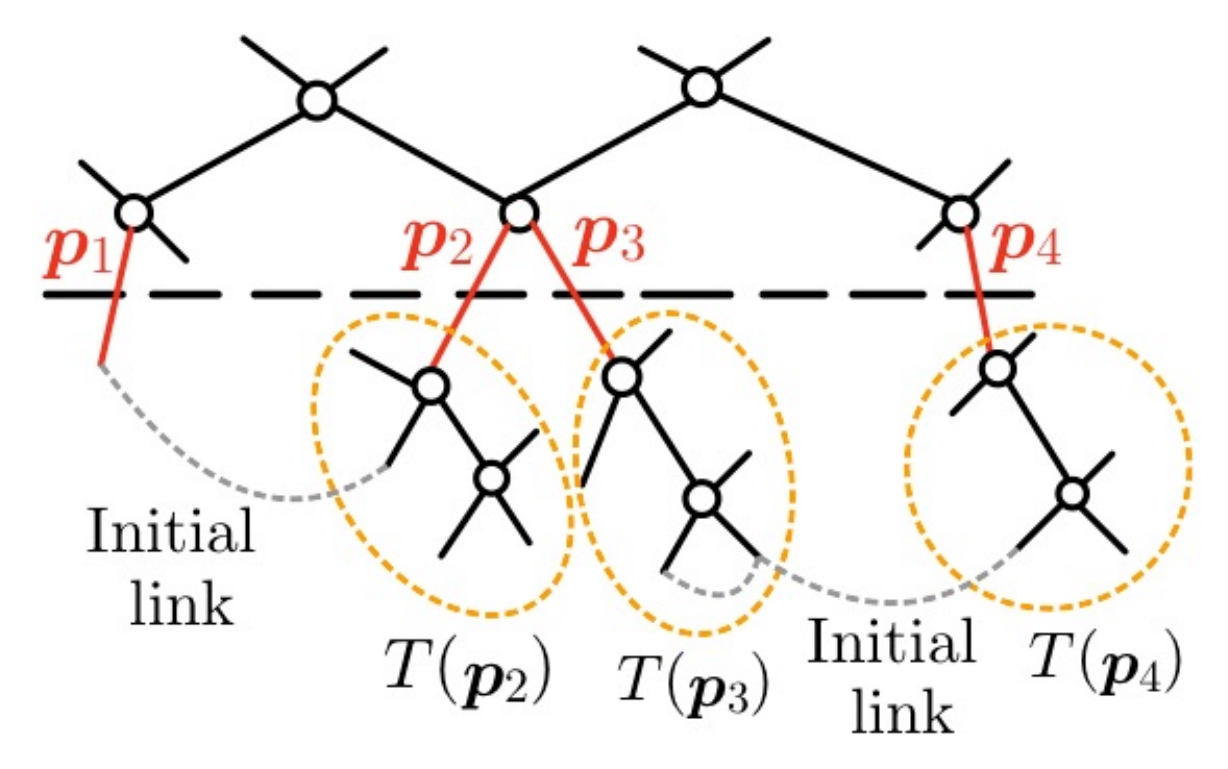}
    \caption{Proposition \ref{prop.layer_cutting} (\ref{it.layer_cutting_1}). Here $\pb_j\in S$ are the particle lines in Proposition \ref{prop.layer_cutting} (\ref{it.layer_cutting_1}) colored in red, with the corresponding set $T(\pb_j)$ and initial links.
}
    \label{fig.no_zeta_D}
\end{figure}

\item\label{it.layer_cutting_2} If $\zeta_D$ exists, then after \textbf{Option 2}, $\Mb$ is cut into elementary molecules with $\#_{\{33B\}}=\#_{\{44\}}=0$. We also have $\#_{\{4\}}\leq 100(C_{5}^*)^{-1}\cdot\rho'''+G$ (note here $G$ equals the number of components of $\Vb$ that are cut as free in \textbf{Option 2} (ii)). Moreover, each atom in $\Vb$ belongs to a \{3\} or \{4\} molecule.

\item\label{it.layer_cutting_3} If $\zeta_D$ exists, then after \textbf{Option 1}, the molecule $\Mb_{UD}$ consists of parts of $\Mb_{\zeta_U}^T$, parts of $\Mb_{\zeta_D}^T$, and $Y$. Define $\Mb_U$ to be the union of $\Mb_{\zeta_U}^T$ and $Y$, and all the components in $\Mb_{\zeta_D}^T$ that are full (after \textbf{Option 1}) and have exactly one bond connecting to $\Mb_{\zeta_U}^T\cup Y$. Let also $\Mb_D:=\Mb_{UD}\backslash \Mb_U$ (which is a subset of $\Mb_{\zeta_D}^T$). Then we have the following properties for $\Mb_U$ and $\Mb_D$:
\begin{enumerate}
\item\label{it.layer_cutting_3i} Each of $\Mb_U$ and $\Mb_D$ is a forest and has no deg 2 atoms. No atom in $\Mb_D$ is a parent of any atom in $\Mb_U$, and $\Mb_U$ (resp. $\Mb_D$) has no C-top (resp. C-bottom) fixed end.
\item\label{it.layer_cutting_3ii} The number of full (i.e. all deg 4) components of $\Mb_U$ is at most $10(C_{5}^*)^{-1}\cdot\rho'''$. The number of components of $\Mb_D$ is $\geq G$. Each component in $\Mb_D$ either has at $\geq 2$ bonds connecting to $\Mb_U$, or has one bond connecting to $\Mb_U$ and one deg 3 atom.
\item\label{it.layer_cutting_3iii} Define the set of \textbf{UD connections} to be the set of bonds connecting an atom in $\Mb_U$ and an atom in $\Mb_D$. We say two UD connections are \textbf{connected via} $\Mb_D$, if the $\Mb_D$ endpoints of them belong to the same connected component of $\Mb_D$. Then, there exist at least $(\rho'''/4)$ many UD connections, such that each UD connection is either connected to another UD connection via $\Mb_D$, or has its $\Mb_D$ endpoint belonging to component of $\Mb_D$ that contains a deg 3 atom.
\end{enumerate}
\end{enumerate}
\end{proposition}
\begin{proof} We first prove the easier statements, and leave the proof of the hardest (\ref{it.layer_cutting_3iii}) to the end.

\textbf{Proof of (\ref{it.layer_cutting_1}).} The statement about $\#_{\{4\}}$ follow from Proposition \ref{prop.layer_select} (\ref{it.layer_select_new1})--(\ref{it.layer_select_new3}) and Proposition \ref{prop.alg_up} (\ref{it.up_proof_2}) (and Remark \ref{rem.full_cut} about full components of $\Mb_{<\zeta_D}^T$ after cutting everything else as free). Now suppose $\zeta_D$ is absent (so $v_{(1,0)}^*\leq (C_{5}^*)^{-1}\cdot\rho''$ in Definition \ref{def.layer_select} (\ref{it.layer_select_3}); the case $\ell_1=1$ is similar, in which case each $T(\pb)$ will be empty). Note that $|S|=\rho''=\rho'''$ by Definition \ref{def.layer_select}. Choose all the particle lines $\pb\in S$ such that $\pb$ either intersects a component of $Y$ or does not intersect $\Mb_{<\zeta_U}^T$ at all, then all particle lines $\pb\in S$ are chosen with at most $(C_{5}^*)^{-1}\cdot\rho'''$ exceptions, by Definition \ref{def.layer_select} (\ref{it.layer_select_3}) and $v_{(1,0)}^*\leq (C_{5}^*)^{-1}\cdot\rho''$.

For each chosen $\pb$, set $T(\pb)$ to be the component of $Y$ that intersects $\pb$ (or $\varnothing$ if $\pb$ does not intersect $\Mb_{<\zeta_U}^T$). Next, by Proposition \ref{prop.layer_refine_2} (\ref{it.layer_refine_27}), we know that with at most $(C_{5}^*)^{-1}\cdot\rho'''$ exceptions, each chosen $\pb$ must form an initial link within $\Mb_{<\zeta_U}$. For each such $\pb$, define the corresponding bottom free end $e$ by the initial link condition; we then claim that these choices satisfy the requirements (\ref{it.layer_cutting_1i})--(\ref{it.layer_cutting_1iii}). 

To prove this, note that (\ref{it.layer_cutting_1i}) and (\ref{it.layer_cutting_1ii}) follow from properties of $Y$. As for (\ref{it.layer_cutting_1iii}), if $\pb$ forms an initial link with a particle line $\pb'$ (this $\pb'$ must belong to $H_0=r(\Mb_0)$), such that the corresponding free end $e'$ does not belong to $\Ec_0$ and does not have an atom in a component cut in \textbf{Option 1}, then $\pb'$ must not intersect $\Mb_{<\zeta_U}^T\cup\Mb_{\zeta_U}^T$ after layer refining (indeed, if $\pb'$ intersects a component of $\Mb_{<\zeta_U}^T\backslash Y$ then it contains an atom in a component that is cut; if $\pb'$ intersects a component of $Y$, or does not intersect $\Mb_{<\zeta_U}^T$ and intersects $\Mb_{\zeta_U}^T$, then $\pb'$ must belong to $\Ec_0$). This means that this $\pb'$ either contains a deleted O-atom or belong to $r(\Mb_{\ell_1-1})$ (by applying Proposition \ref{prop.mol_axiom} (\ref{it.axiom5})), which leads to at most $(C_{5}^*)^{-1}\cdot\rho'''$ choices for $(\pb,\pb')$ (in the former case using the upper bound of the number of deleted atoms, and in the latter case using Definition \ref{def.layer_select} (\ref{it.layer_select_2})), and is thus negligible. In summary, at least $\rho'''/2$ particle lines $\pb\in S$ must be chosen, and at least $\rho'''/4$ chosen $\pb$ must satisfy the given restrictions.
 
\textbf{Proof of (\ref{it.layer_cutting_2}).} The bound on $\#_{\{4\}}$ follows from the above proof, with one extra \{4\} molecule for each of the $G$ components of $\Vb$ cut as free; note that each component cut as fixed will not be full, so they do not contribute \{4\} molecules. Each atom in $\Vb$ belongs to a \{3\} or \{4\} molecule, simply because $\Vb\subseteq\Mb_{\zeta_D}^T$ which is a forest, and each component of $\Vb$ cut as fixed will have exactly one deg 3 atom, so there will not be any \{2\} or \{33\} molecule by Proposition \ref{prop.alg_up} (\ref{it.up_proof_1}).
 
\textbf{Proof of (\ref{it.layer_cutting_3i}) and (\ref{it.layer_cutting_3ii}).} First, each of $\Mb_{\zeta_U}^T$ and $\Mb_{\zeta_D}^T$ is a forest and does not have deg 2 atoms, by layer refining and the operations of cutting deg 2 atoms in \textbf{Option 1} (ii). Moreover $\Mb_{\zeta_U}^T$ (resp. $\Mb_{\zeta_D}^T$) does not have any C-top (resp. C-bottom) fixed end, because cutting as free any deg 2 atom in $\Mb_{\zeta_U}^T$ (resp. $\Mb_{\zeta_D}^T$) does not generate any C-top (resp. C-bottom) fixed end in $\Mb_{\zeta_U}^T$ (resp. $\Mb_{\zeta_D}^T$) if there is none initially (of course, cutting any atom in $\Mb_{(\zeta_D:\zeta_U)}^T$ also does not create any C-top fixed end in $\Mb_{\zeta_U}^T$ or C-bottom fixed end in $\Mb_{\zeta_D}^T$). The fact that no atom in $\Mb_{\zeta_D}^T$ can be parent of any atom in $\Mb_{\zeta_U}^T$ is also clear from properties of the thin layers.
 
 All the above properties are also preserved after adding $Y$ to $\Mb_{\zeta_U}^T$, and adding to $\Mb_{\zeta_U}^T\cup Y$ a full component in $\Mb_{\zeta_D}^T$ that has only one bond connecting to $\Mb_{\zeta_U}^T\cup Y$. For example, the forest property is preserved because each component of $Y$ is a tree, and each component of $\Mb_{\zeta_D}^T$ is also a tree, and each component of $Y$ has exactly one bond\footnote{This is because, each component of $Y$ contains an atom in a particle line in $S$ (this atom must belong to some layer $\Mb_{\ell'}\,(\ell'\leq\ell_1-1)$ by definition of $S$). If it is connected to $\Mb_{\zeta_U}^T\subseteq \Mb_{\ell_1}$ by a bond, then there must exist such a bond connecting an atom in $\Mb_{\ell'}\,(\ell'\leq\ell_1-1)$ to an atom in $\Mb_{\ell_1}$; by Proposition \ref{prop.mol_axiom} (\ref{it.axiom4}) this particle line must belong to $p(\Mb_{\ell_1})\cap p(\Mb_{\ell_1-1})=r(\Mb_{\ell_1-1})$, so by definition of $Y$ in Proposition \ref{prop.layer_select}, this particle line (and this bond) must be unique.} connecting to $\Mb_{\zeta_U}^T$; the other properties are easy to prove, where we also note that $Y$ has no fixed end, because cutting any atom in $\Mb_{(\zeta_D:\zeta_U)}^T\backslash Y$, or any subsequent deg 2 atom in $\Mb_U$ or $\Mb_D$, does not create a fixed end in $Y$.
 
Next, the number of full components of $\Mb_{\zeta_U}^T$ is at most $10(C_{5}^*)^{-1}\cdot\rho'''$ by Proposition \ref{prop.layer_select} (\ref{it.layer_select_new1}), and this cannot increase after cutting $\Mb_{(\zeta_U:\zeta_D)}^T\backslash Y$ as free. Moreover adding $Y$ and the components in $\Mb_{\zeta_D}^T$ will not increase this number, as these added sets are connected to $\Mb_{\zeta_U}^T$ by a bond. Next, each component in $\Mb_D$ must have at least one bond connecting to $\Mb_U$ (otherwise it would have been cut in \textbf{Option 1} (iii)), and it must either have 2 bonds connecting to $\Mb_U$, or contain a deg 3 atom (otherwise it would have been included in $\Mb_U$ by the definition of $\Mb_U$).

Finally we prove the number of components of $\Mb_D$ is at least $G$. In fact, for each of the $G$ components $W\subseteq\Vb$ that has at least two bonds connecting to $\Mb_{\zeta_U}^T\cup\Mb_{(\zeta_U:\zeta_D)}^T$, we claim that after \textbf{Option 1} in Definition \ref{def.layer_cutting}, it must still contain a component that has at least one UD connection and either another UD connection or a deg 3 atom. This claim then implies that this component within $W$ must stay in $\Mb_D$ instead of being included in $\Mb_U$ (see {\color{blue}Figure \ref{fig.option2}}); as such, we know that the number of components of $\Mb_D$ is still $\geq G$ after \textbf{Option 1} in Definition \ref{def.layer_cutting}.

To prove the above claim, note that $W$ has a bond $e$ connecting to $\Mb_{\zeta_U}^T\cup Y$ which cannot be broken by \textbf{Option 1} and remains as a UD connection (see the proof of (\ref{it.layer_cutting_3iii}) below); the other bond $e'$ connecting to $\Mb_{\zeta_U}^T\cup\Mb_{(\zeta_U:\zeta_D)}^T$ is either the same (and leads to another UD connection) or connects to $\Mb_{(\zeta_U:\zeta_D)}^T\backslash Y$, which creates a fixed end in $\Mb_D$ (or a fixed end in each newly created component of $W$ in the presence of O-atoms) after $\Mb_{(\zeta_U:\zeta_D)}^T\backslash Y$ being cut in \textbf{Option 1}. In either case, see the proof of (\ref{it.layer_cutting_3iii}) below, we conclude that the new component of $W$ containing the $\Mb_D$ atom of $e$ must either contain another UD connection or a deg 3 atom, as desired.

\textbf{Proof of (\ref{it.layer_cutting_3iii})}. For $\pb,\pb'\in r(\Mb_{\ell_1-1})$ and $\zeta_D\leq\zeta\leq\zeta_0$ (recall $\zeta_0$ is the highest thin layer below $\zeta_U$), define $\pb\sim_{\zeta}\pb'$ if the particle lines $\pb$ and $\pb'$ are connected via $\Mb_{[\zeta:\zeta_0]}^T$. It is easy to see that $\sim_{\zeta}$ is an equivalence relation, which again follows from Remark \ref{rem.layer_interval}.

Recall the definition of $S$ in Definition \ref{def.layer_select} (\ref{it.layer_select_3}); define $S_1$ to be the set of $\pb\in r(\Mb_{\ell_1-1})$ that do not contain any atom in $\Mb_{(\ell_1,k)}^T$ for $k\leq \gamma$. Define $S_2$ to be the set of $\pb\in S$ that is either connected to a cycle within $\Mb_{[\zeta_D:\zeta_0]}^T$ or satisfies $\pb\sim_{\zeta_D}\pb'$ for some other $\pb'\in r(\Mb_{\ell_1-1})$, but does not have the same property if $\zeta_D$ is replaced by any $\zeta>\zeta_D$. By (\ref{eq.layer_select_2}) and (\ref{eq.layer_select_3}), we have $|S_1|< (10C_{5}^*)^{-1}\cdot\rho'''$, $|S_2|\geq (0.9)\rho'''$ and $S_1\cap S_2=\varnothing$. Then define $S_3$ same as $S_2$, but with the extra restriction $\pb'\not\in S_1$ in the relation $\pb\sim_{\zeta_D}\pb'$ above; since $\sim_{\zeta_D}$ is an equivalence relation, we know that $|S_3|\geq |S_2|-|S_1|\geq (0.8)\rho'''$ (indeed, if $\pb\sim_{\zeta_D}\pb'$ for some $\pb'\in S_1$, then each $\pb'$ would lead to a unique $\pb$, unless two $\pb$ are $\sim_{\zeta_D}$ with the same $\pb'$, in which case these to $\pb$ are $\sim_{\zeta_D}$ with each other and they will belong to $S_3$).

Now we claim that after \textbf{Option 1}, for each $\pb\in S_3$, there exists a UD connection $e$ in $\Mb_{UD}$ such that (i) it is either connected to some other UD connection $e'$ via $\Mb_D$, or has its $\Mb_D$ endpoint in a component of $\Mb_D$ that contains at least one deg 3 atom, and (ii) it either belongs to $\pb$, or its $\Mb_U$ endpoint belongs to a connected component of $Y$ intersecting $\pb$. Once this is proved, is follows from (ii) that the UD connection obtained for different $\pb$ must be different (because each $Y$ intersects a unique $\pb$), so the desired result (\ref{it.layer_cutting_3iii}) follows from (i).

To prove the above claim, first assume there exists $\pb'\neq\pb$ and $\pb'\not\in S_1$ such that $\pb\sim_{\zeta_D}\pb'$. Since $\pb\in S_3\subseteq S_2\subseteq S$ and $\pb'\not\in S_1$, we know that $\pb$ intersects $\Mb_{\zeta_U}^T$, and either does not intersect $\Mb_{(\zeta_D:\zeta_U)}^T$, or intersects a unique component $Z$ of $\Mb_{(\zeta_D:\zeta_U)}^T$ contained in $Y$; moreover $\pb'\not\in S_1$ must either intersect $\Mb_{\zeta_U}^T$, or intersect a unique component $W\neq Z$ of $\Mb_{(\zeta_D:\zeta_U)}^T$. Since $\pb\sim_{\zeta_D}\pb'$, we consider a path going from an atom on $\pb$ to an atom on $\pb'$ within $\Mb_{[\zeta_D:\zeta_0]}^T$. This path must enter $\zeta_D$ because $\pb\not\sim_{\zeta}\pb'$ for the layer $\zeta>\zeta_D$; moreover, it must then exit $\zeta_D$ because in the end it reaches $\pb'$ which contains an atom in either $\Mb_{(\zeta_D:\zeta_U)}^T$ or $\Mb_{\zeta_U}^T$.

Consider the first entry and first exit to $\Mb_{\zeta_D}^T$, then we can find two bonds $e$ and $e'$ such that (i) each has only one endpoint in $\Mb_{\zeta_D}^T$ and these two endpoints (say $\nf$ and $\nf'$) belong to the same component $X$ of $\Mb_{\zeta_D}^T$ \emph{before} \textbf{Option 1}, (ii) either $e$ belongs to $\pb$ or the other endpoint of $e$ belongs to $Z$, and (iii) either $e'$ belongs to $\pb'$ or the other endpoint of $e'$ belongs to a component $Z_1$ of $\Mb_{(\zeta_D:\zeta_U)}^T$ other than $Z$ (which may or may not be $W$). See {\color{blue}Figure \ref{fig.udconn}}.
\begin{figure}[h!]
    \centering
    \includegraphics[width=0.24\linewidth]{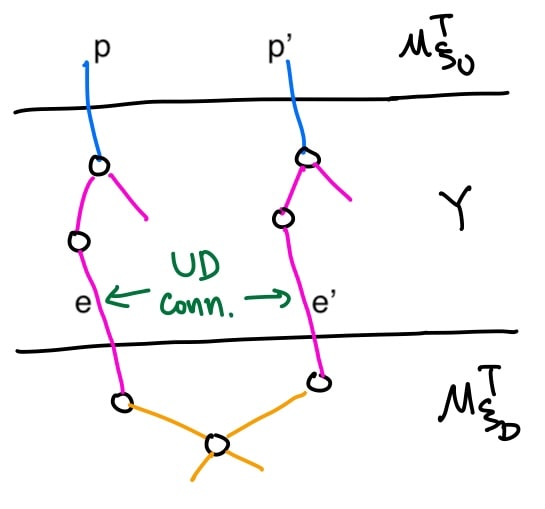}
    \caption{Proof of Proposition \ref{prop.layer_cutting} (3c). Here the two particle lines $\pb$ and $\pb'$ are connected via $\Mb_{\zeta_D}^T\cup Y$ but are \emph{not} connected if without $\Mb_{\zeta_D}^T$. Therefore, the path connecting them must first enter $\Mb_{\zeta_D}^T$ and then exit $\Mb_{\zeta_D}^T$, which leads to two UD connections $e$ and $e'$ that are connected via $\Mb_D$ (so by definition both of them will be good). If $e$ is not broken in any subsequent cutting, it is easy to see that $e$ will remain good (either with connection to $e'$, or with a deg 3 atom).}
    \label{fig.udconn}
\end{figure}

Now we consider the behavior of $(e,e')$ under the cutting operations in \textbf{Option 1}. Clearly each cutting in \textbf{Option 1} (ii) cannot break the maximal ov-segment containing $e$ (i.e. cannot cut any atom in this ov-segment as in Definition \ref{def.cutting}), because we only cut components not contained in $Y$, and any deg 2 atoms that are subsequently cut in \textbf{Option 1} (ii) must belong to $\Mb_{\zeta_D}^T$ (resp. $\Mb_{\zeta_U}^T$) and must be O-atom or have two top (resp. bottom) fixed ends. Since the component $X$ has at least one bond connecting to $\Mb_{\zeta_U}^T\cup Y$ and is thus not cut as free in \textbf{Option 1} (iii), we know that $e$ remains a UD connection after \textbf{Option 1}. Moreover, if the component $X$ is not affected by the cuttings in \textbf{Option 1} (i.e. no atom in $X$ is cut, and no fixed end is created at atoms in $X$), then we must have $Z_1\subseteq Y$ so $e'$ is also a UD connection that is still connected to $e$ via $\Mb_D$ even \emph{after} \textbf{Option 1}; if $X$ is affected, then the (possibly new) component of $X$ containing $\nf$ will contain at least one deg 3 atom \emph{after} \textbf{Option 1} (cf. Remark \ref{rem.full_cut}). In either case the desired property for $e$ follows.

Finally, if the above $\pb'$ does not exist, then $\pb$ must form a cycle. In this case, using this cycle within $\Mb_{[\zeta_D:\zeta_0]}^T$ and note that $\Mb_{\zeta_D}^T$ itself is a forest, we can still find the $e$ and $e'$ satisfying (a)--(c) above with $e\neq e'$, except that it is now possible to have $Z_1=Z$; but this does not affect the proof, and the same arguments as above then proves (3c).
\end{proof}
Now we can prove the next case of Proposition \ref{prop.comb_est}, namely Proposition \ref{prop.comb_est_case4} below, which contains the initial link case (when $\zeta_D$ is absent) and the weak degenerate case (where there are many weak degeneracies in $\Mb_D$ as defined in Definition \ref{def.weadeg}) below. In the weakly degenerate case we choose \textbf{Option 2} in Definition \ref{def.layer_cutting}, which allows to gain from the weak degeneracy conditions at the \{3\} and \{4\} molecules in $\Mb_D$.
\begin{definition}\label{def.weadeg}
Define the atom pair $\{\nf,\nf'\}$ to be \textbf{weakly degenerate} if they are ov-adjacent C-atoms, and we have the restriction (by indicator functions) $|t_{\nf}-t_{\nf'}|\leq \varepsilon^{\upsilon}$. Define also one atom $\nf$ to be weakly degenerate if it is C-atom and we have the restriction $|v_e-v_{e'}|\leq\varepsilon^{\upsilon}$ for two distinct edges $(e,e')$ at $\nf$.
\end{definition}
\begin{proposition} 
\label{prop.comb_est_case4} Suppose $\Mb$ satisfies the assumptions in Reduction \ref{red.1}, and we apply the layer selection process in Definition \ref{def.layer_select}. Then, Proposition \ref{prop.comb_est} is true if $\zeta_D$ is absent. 

Moreover, recall the definition of weak degeneracy in Definition \ref{def.weadeg}. If $G\leq (C_{3}^*)^{-1}\cdot\rho'''$ (with $G$ defined in Proposition \ref{prop.layer_cutting}) and we make the restriction that there exist at least $\Bf:=(C_{2}^*)^{-1}\cdot\rho'''$ weakly degenerate atoms and atom pairs in $\Vb$ (by inserting an indicator function $\mathbbm{1}_{\mathrm{wea.deg}}$),  then Proposition \ref{prop.comb_est} is also true.
\end{proposition}
\begin{proof} Recall our goal is to prove (\ref{eq.overall_alg_3}), see Reduction \ref{red.1}. 

\textbf{Proof part 1.} If $\zeta_D$ is absent, let $\Mb_U$ be the remaining molecule after applying \textbf{Option 1} in Definition \ref{def.layer_cutting}. Note that $\Mb_U$ is a forest by the proof of Proposition \ref{prop.layer_cutting} (\ref{it.layer_cutting_1}). We simply cut it into elementary molecules using \textbf{UP}. Clearly $\#_{\{33B\}}=\#_{\{44\}}=0$, and $\#_{\{4\}}\leq  (C_{5}^*)^{-1}\cdot\rho'''$ by Proposition \ref{prop.layer_cutting} (the number of full components of $\Mb_U$ equals that of $\Mb_{\zeta_U}^T$ after \textbf{Option 1} in Definition \ref{def.layer_cutting}, which is at most $(C_{5}^*)^{-1}\cdot\rho'''$ by Proposition \ref{prop.layer_select}).

Next, for each $\pb$ as in Proposition \ref{prop.layer_cutting} (\ref{it.layer_cutting_1}) and the corresponding bottom end $e$, consider the maximal ov-segment containing $e$ in $\Mb_U$, and the first atom $\mf$ that is cut in this ov-segment. This $\mf$ may belong to a \{4\} molecule but in at most $ (C_{5}^*)^{-1}\cdot\rho'''$ cases. Moreover it cannot belong to a \{2\} molecule; in fact, it cannot be cut as \{2\} molecule in Definition \ref{def.alg_up} (\ref{it.alg_up_1}) (as it has a bottom non-fixed edge when it is cut), and cannot be cut as \{2\} molecule in \ref{def.alg_up} (\ref{it.alg_up_3}) (because $\mf$ belongs to $S_\nf$ for some atom $\nf$ and has only one parent in $S_\nf$ thanks to the fact that $\Mb_U$ is a forest; see {\color{blue}Figure \ref{fig.up_aux}}).

As such, we know that with at most $ (C_{5}^*)^{-1}\cdot\rho'''$ exceptions, each such atom $\mf$ defined above must belong to either a \{33A\} molecule (which can be treated as good, see Reduction \ref{red.1}) or a \{3\} molecule $\Xc=\{\mf\}$. Moreover, if $e$ is initial linked to some $e'$, which belongs to an elementary molecule that is either $\Xc$ or in $\Mb_U$ and cut before $\Xc$ or cut in \textbf{Option 1} of Definition \ref{def.layer_cutting}, then this $\Xc$ must be good molecule due to Definition \ref{def.good_normal}. In summary, from the $\geq \rho'''/4$ particle lines $\pb$ in Proposition \ref{prop.layer_cutting} (\ref{it.layer_cutting_1}), we obtain at least $\rho'''/8$ good molecules (each bottom end may be involved in two initial links). This then leads to
\[\#_{\mathrm{good}}\geq \rho'''/8,\qquad \#_{\{4\}}\leq 200 (C_{5}^*)^{-1}\cdot\rho''';\qquad \rho'''\geq (C_{7}^*)^{-10\Lf}\cdot\rho\] by the above arguments and Proposition \ref{prop.layer_cutting} (\ref{it.layer_cutting_1}) and Definition \ref{def.layer_select}. This clearly implies (\ref{eq.overall_alg_3}).

\textbf{Proof part 2.} Now assume $\zeta_D$ exists, and $G\leq (C_{3}^*)^{-1}\cdot\rho'''$, and there exist at least $\Bf$ weakly degenerate atoms and atom pairs in $\Vb$. By splitting into $|\log\varepsilon|^{C^*\Bf}$ sub-cases (by decomposing 1 into the corresponding indicator functions), we may fix one choice of these atoms and pairs. Then we apply \textbf{Option 2} of Definition \ref{def.layer_cutting} and cut $\Mb$ into elementary molecules. Clearly $\#_{\{33B\}}=\#_{\{44\}}=0$, and \[\#_{\{4\}}\leq 100(C_{5}^*)^{-1}\cdot\rho'''+G\leq 2(C_{3}^*)^{-1}\cdot\rho'''\] by Proposition \ref{prop.layer_cutting} (\ref{it.layer_cutting_2}) and assumption on $G$. Moreover, if we exclude the (at most $2(C_{3}^*)^{-1}\cdot\rho'''$) exceptional atoms in $\Vb$ that belong to \{4\} molecules, then each atom in $\Vb$ belongs to a \{3\} molecule. In addition, if $\nf$ is a weakly degenerate atom then $\{\nf\}$ is a good molecule by Definition \ref{def.good_normal}; if $\{\nf,\nf'\}$ is a weakly degenerate pair, and $\nf$ is cut before $\nf'$, then $\{\nf'\}$ is also a good molecule by Definition \ref{def.good_normal}. In either case, the $\Bf$ weakly degenerate atoms and pairs lead to at least $\Bf/4$ good molecules, thus \[\#_{\mathrm{good}}\geq (4C_2^*)^{-1}\cdot\rho''',\qquad \#_{\{4\}}\leq 2(C_{3}^*)^{-1}\cdot\rho''';\qquad\rho'''\geq (C_{8}^*)^{-1}\cdot\rho,\] which implies (\ref{eq.overall_alg_3}).
\end{proof}
Now, summarizing Propositions \ref{prop.case1}, \ref{prop.case2}, \ref{prop.comb_est_case3} and \ref{prop.comb_est_case4}, we see that, for the rest of the proof of Proposition \ref{prop.comb_est}, we may assume that the assumptions in Reduction \ref{red.1} are satisfied, that $\zeta_D$ in Definition \ref{def.layer_select} exists, and that either $G\geq (C_{3}^*)^{-1}\cdot\rho'''$ with $G$ defined in Proposition \ref{prop.layer_cutting} or we can restrict the number of weakly degenerate atoms and pairs in $\Vb$ to be less than $(C_{2}^*)^{-1}\cdot\rho'''$ (by inserting the indicator function $1-\mathbbm{1}_{\mathrm{wea.deg}}$ opposite to the one defined in Proposition \ref{prop.comb_est_case4}).

We then apply \textbf{Option 1} in Definition \ref{def.layer_cutting} and reduce to the molecule $\Mb_{UD}$ that satisfies (3a)--(3c) of Proposition \ref{prop.layer_cutting}. In particular the number of components of $\Mb_D$ satisfies $\#_{\mathrm{comp}(\Mb_D)}\geq G$ by (3b) of Proposition \ref{prop.layer_cutting}. Using also the upper bound on the contribution to $\#_{\{4\}}$ of \textbf{Option 1} in Definition \ref{def.layer_cutting} (proved in Proposition \ref{prop.layer_cutting}), we see that Proposition \ref{prop.comb_est} (i.e. (\ref{eq.overall_alg_3})) now follows from the following
\begin{proposition}
\label{prop.case5} Let $\Mb_{UD}$ be a molecule satisfying the properties (3a)--(3c) in Proposition \ref{prop.layer_cutting} with some $\rho'''\geq (C_{8}^*)^{-1}\cdot\rho$. Moreover assume that either the number of components of $\Mb_D$ satisfies $\#_{\mathrm{comp}(\Mb_D)}\geq (C_{3}^*)^{-1}\cdot\rho'''$, or (we can restrict by some indicator function such that) the number of weakly degenerate atoms and pairs in $\Mb_D$ is less than $(C_{2}^*)^{-1}\cdot\rho'''$. Then in either case, we can cut $\Mb_{UD}$ into elementary molecules such that $\#_{\{44\}}=0$ and each \{33B\} molecule is good, and moreover
\begin{equation}\label{eq.comb_est_case51} (\upsilon/2)\cdot(\#_{\mathrm{good}})-d\cdot(\#_{\{4\}})\geq (C_{4}^*)^{-1}\cdot \rho''',
\end{equation} where we again understand that all \{33A\} molecules are good as in Reduction \ref{red.1}.
\end{proposition}

\section{Cutting algorithm II: analyzing UD molecules}\label{sec.maincr} The goal of this section is to prove Proposition \ref{prop.case5} for the two-layer molecule $\Mb_{UD}$, by introducing the various algorithms corresponding to the toy models in Section \ref{sec.toy}. To simplify notation in this section, if a molecule $\Mb_{UD}=\Mb_U\cup\Mb_D$ satisfies that (a) each of $\Mb_U$ and $\Mb_D$ is a forest and has no deg 2 atoms, and (b) no atom in $\Mb_D$ is parent of any atom in $\Mb_U$, then we say it is a \textbf{UD molecule}. If in addition $\Mb_U$ has no C-top fixed end and $\Mb_D$ has no C-bottom fixed end, then we say $\Mb_{UD}$ is \textbf{canonical}.
\begin{reduct} Our proof below relies on the properties (3a)--(3c) in Proposition \ref{prop.layer_cutting}. For convenience we record them here:
\label{red.2}
\begin{enumerate}[{(i)}]
\item Each of $\Mb_U$ and $\Mb_D$ is a forest and has no deg 2 atoms. No atom in $\Mb_D$ is a parent of any atom in $\Mb_U$, and $\Mb_U$ (resp. $\Mb_D$) has no C-top (resp. C-bottom) fixed end.
\item The number of full  components of $\Mb_U$ is at most $10(C_{5}^*)^{-1}\cdot\rho'''$. Each component in $\Mb_D$ either has at least two bonds connecting to $\Mb_U$, or has one bond connecting to $\Mb_U$ and one deg 3 atom.
\item There exist at least $(\rho'''/4)$ many UD connections, such that each UD connection is either connected to another UD connection via $\Mb_D$, or has its $\Mb_D$ endpoint belonging to component of $\Mb_D$ that contains a deg 3 atom.
\end{enumerate}
\end{reduct}
\subsection{The \textbf{2CONNUP} algorithm and toy model III}\label{sec.2connup}
We start by introducing the full \textbf{2CONNUP} algorithm corresponding to toy model III. Compared to the toy version (Definition \ref{def.toy3_alg}), the main extra ingredient is Definition \ref{def.alg_2connup} (\ref{it.alg_2connup_1}); see the proof of Proposition \ref{prop.alg_2connup} for more explanations.
\begin{definition}[The algorithm \textbf{2CONNUP}]
\label{def.alg_2connup}
Define the following cutting algorithm \textbf{2CONNUP}. It takes as input any canonical UD molecule $\Mb_{UD}$, such that $\Mb_D$ has no full component that has no bond connecting to $\Mb_U$. Define $\Mb_{\mathrm{2conn}}^1$ to be the set of atoms $\nf\in\Mb_D$ that have deg 3 and are ov-adjacent to one atom $\mf\in\Mb_U$, such that no other $\Mb_D$ atom on the ov-segment between $\nf$ and $\mf$ has deg 3. Define $\Mb_{\mathrm{2conn}}^2$ to be the set of  atoms $\nf\in\Mb_D$ that are ov-adjacent to two atoms in $\Mb_U$ along two different particle lines, and $\Mb_{\mathrm{2conn}}=\Mb_{\mathrm{2conn}}^1\cup\Mb_{\mathrm{2conn}}^2$. For any such $\Mb=\Mb_{UD}$, we define the cutting sequence as follows:
\begin{enumerate}
\item\label{it.alg_2connup_1} If any (maximal) ov-segment contains at least $A:=32\upsilon^{-2}$ atoms in $\Mb_{\mathrm{2conn}}^2$, then we choose $\nf_0$ to be any one of these atoms, and cut it as free. Repeat until no such ov-segment exists.
\item\label{it.alg_2connup_2} If $\Mb$ contains any deg $2$ atom $\nf$, then cut it as free, and repeat until there is no deg $2$ atom left.
\item\label{it.alg_2connup_3} Choose a lowest deg 3 atom $\nf\in\Mb_U$ among all remaining deg 3 atoms in $\Mb_U$ (or a lowest atom in $\Mb_U$ if $\Mb_U$ has no deg 3 atom). Let $S_\nf$ be the set of descendants of $\nf$ in $\Mb_U$; this $\nf$ and $S_\nf$ will be fixed until the end of (\ref{it.alg_2connup_4}).
\item\label{it.alg_2connup_4} Now, starting from $\nf$, each time choose a highest atom $\mf$ in $S_\nf$. If $\mf$ has deg 3, and $\mf$ is ov-adjacent to an atom $\pf\in \Mb_{\mathrm{2conn}}$ that also has deg 3, then cut $\{\mf,\pf\}$ as free. Otherwise just cut $\mf$ as free. Repeat until all atoms in $S_\nf$ has been cut. Then go to (\ref{it.alg_2connup_2}).
\item\label{it.alg_2connup_5} Repeat (\ref{it.alg_2connup_2})--(\ref{it.alg_2connup_4}) until all atoms in $\Mb_U$ has been cut. Then apply algorithm \textbf{DOWN} to cut the rest of $\Mb_D$ into elementary molecules.
\end{enumerate}
We also define the dual algorithm \textbf{2CONNDN}, by reversing the roles of $\Mb_U$ and $\Mb_D$ and the notions of parent/child etc., just as how we define \textbf{DOWN} from \textbf{UP} in Definition \ref{def.alg_up}.
\end{definition}

The key estimate for the algorithm \textbf{2CONNUP} is stated in the following proposition, which is the full version of Proposition \ref{prop.toy3} in Section \ref{sec.toy}.

\begin{proposition}
\label{prop.alg_2connup} For the algorithm \textbf{2CONNUP}, we have $\#_{\{33B\}}=\#_{\{44\}}=0$, and
\begin{equation}\label{eq.alg_2connup_1}
(\upsilon/2)\cdot\#_{\{33A\}}-d\cdot \#_{\{4\}}\geq C^{-1}(|\Mb_{\mathrm{2conn}}^1|+|\Mb_{\mathrm{2conn}}^2|)-C\cdot \#_{\mathrm{fucp}(\Mb_U)},
\end{equation} where $\#_{\mathrm{fucp}(\Mb_U)}$ is the number of \emph{full components} of $\Mb_U$ (note the difference with $\#_{\mathrm{comp}(\Mb_U)}$, i.e. the number of \emph{all components}). Recall also the notion of constant $C$ in Definition \ref{def.notation} (\ref{it.defC*}). Similar results hold for the dual algorithm \textbf{2CONNDN}, if we reverse the roles of $\Mb_U$ and $\Mb_D$.
\end{proposition}
\begin{proof} We first explain the relevance of Definition \ref{def.alg_2connup} (\ref{it.alg_2connup_1}), see {\color{blue}Figure \ref{fig.2conn_ov}}. Note that, in Definition \ref{def.toy3_alg} where there is no O-atom, each atom in $\Mb_{\mathrm{2conn}}^2$ is \emph{adjacent} to two atoms in $\Mb_U$, so there cannot be a \emph{bond} containing two atoms in $\Mb_{\mathrm{2conn}}^2$; however, now an atom in $\Mb_{\mathrm{2conn}}^2$ is only \emph{ov-adjacent} to two atoms in $\Mb_U$, so there might be other atoms in $\Mb_{\mathrm{2conn}}^2$ within the ov-segment between it and atoms in $\Mb_U$.
\begin{figure}[h!]
    \centering
    \includegraphics[width=0.6\linewidth]{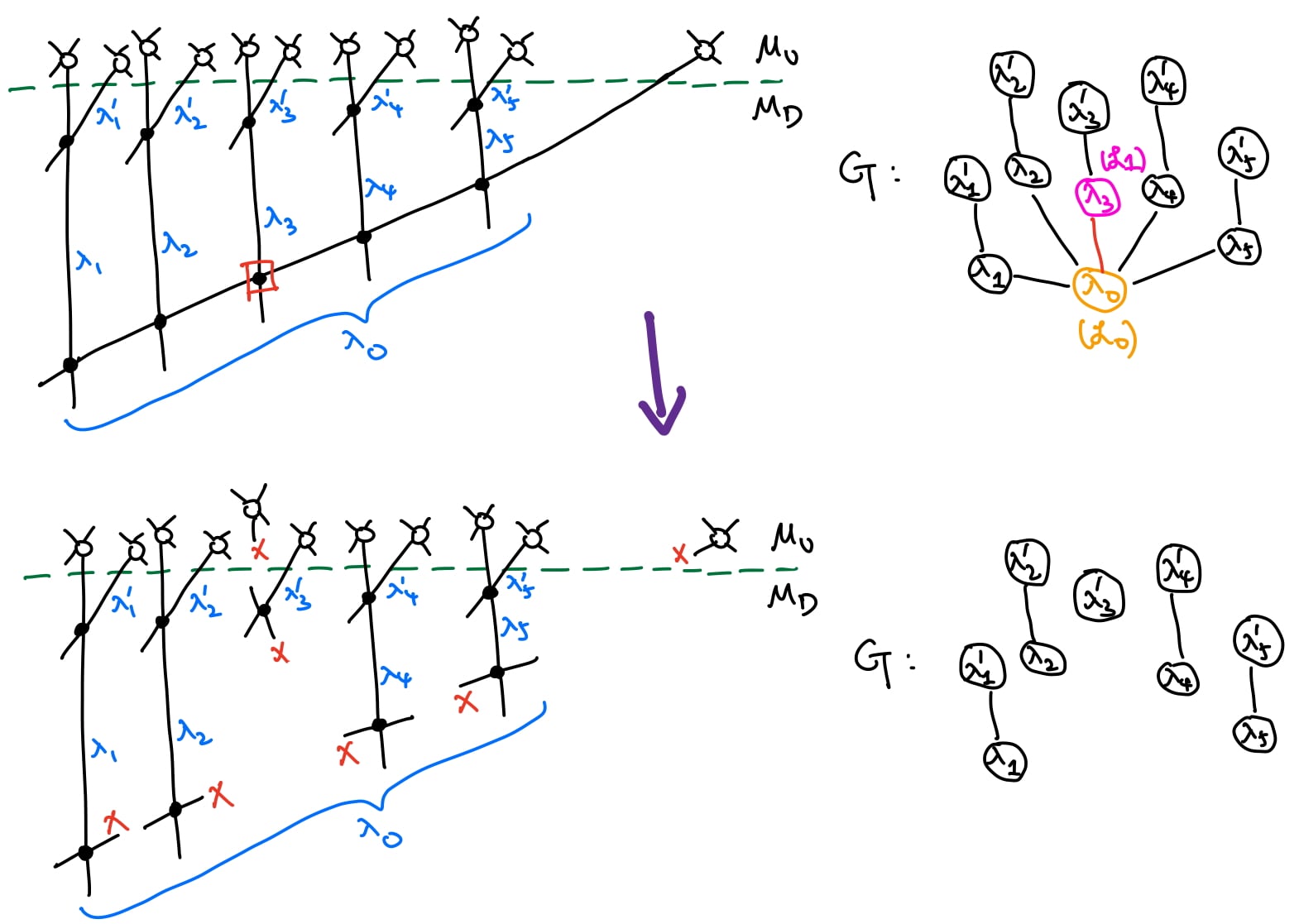}
    \caption{Illustration of Definition \ref{def.alg_2connup} (\ref{it.alg_2connup_1}). Here $(\lambda_0,\lambda_j,\lambda_j')$ are the ov-segments in $\Mb_D$ connecting an atom in $\Mb_{\mathrm{2conn}}^2$ to $\Mb_U$. They intersect at different $\Mb_{\mathrm{2conn}}^2$ atoms; see also the auxiliary graph $G$ on the right (G-edges correspond to intersections of these ov-segments in $\Mb_{\mathrm{2conn}}^2$). Now $\lambda_0$ contains many atoms in $\Mb_{\mathrm{2conn}}^2$ (i.e. has large degree in $G$), so we choose one atom (the one with $\lambda_3$ here) and cut it, breaking both $\lambda_0$ and $\lambda_3$ but not the others. The point is that (i) the number of such operations is negligible, (ii) each operation breaks two ov-segments, so the total number of ov-segments broken is also negligible, and (iii) after Definition \ref{def.alg_2connup} (\ref{it.alg_2connup_1}), each ov-segment contains at most $O(1)$ atoms in $\Mb_{\mathrm{2conn}}^2$, so the same algorithm \textbf{2CONNUP} will suffice. See Part 1 of the proof below.}
    \label{fig.2conn_ov}
\end{figure}

If we still argue naively as in Definition \ref{def.toy3_alg}, it might happen that we break an ov-segment containing many atoms in $\Mb_{\mathrm{2conn}}^2$, turning them into deg 2 simultaneously. This is undesirable, as we expect each atom in $\Mb_{\mathrm{2conn}}^2$ to contribute a \{33A\} molecule in the proof of Proposition \ref{prop.toy3}. To resolve this, we would like that each ov-segment contains at most $O(1)$ (more precisely $A=32\upsilon^{-2}$) atoms in $\Mb_{\mathrm{2conn}}^2$; this leads to the process in Definition \ref{def.alg_2connup} (\ref{it.alg_2connup_1}), where we cut and break any ov-segments not satisfying the above assumptions, so that the same arguments in Definition \ref{def.toy3_alg} will be sufficient after such cuttings.

On the other hand, since $\Mb_D$ has no cycle, it is easy to see that we can define an auxiliary graph $G$ in Part 1 below (with G-nodes being ov-segments and G-edges being atoms in $\Mb_{\mathrm{2conn}}^2$ that are intersections of these ov-segments) which is a forest, so the \emph{average} number of $\Mb_{\mathrm{2conn}}^2$ atoms on all ov-segments is \emph{less than $2$}. As such, the number of ov-segments broken in Definition \ref{def.alg_2connup} (\ref{it.alg_2connup_1}) will be negligible and does not affect the subsequent proof, see the calculations in Part 2 below.

\textbf{Proof part 1.} We claim that throughout the algorithm \textbf{2CONNUP}, the (MONO) property in Proposition \ref{prop.alg_up} is preserved in $\Mb_U$, and $\Mb_D$ has no C-bottom fixed end. In fact, the (MONO) property is proved in the same way as in Proposition \ref{prop.alg_up} (which is also not affected by cutting the $\Mb_D$ atoms in Definition \ref{def.alg_2connup} (\ref{it.alg_2connup_1}), as they do not create C-top fixed end in $\Mb_U$). Regarding $\Mb_D$, note that cutting atoms in $\Mb_U$, or atoms in $\Mb_{\mathrm{2conn}}^2$ which has two ov-segments connecting to $\Mb_U$, will not introduce C-bottom fixed end in $\Mb_D$; as for cutting $\{\mf,\pf\}$ as free in Definition \ref{def.alg_2connup} (\ref{it.alg_2connup_4}), it can be viewed as first cutting $\mf$ and then cutting $\pf$, and $\pf$ must be O-atom or have two top fixed ends (cf. Remark \ref{rem.reg}) when it is cut, so cutting $\pf$ also does not introduce C-bottom fixed end in $\Mb_D$. It then follows that we only have \{2\}, \{3\}, \{4\} and \{33A\} molecules, thus $\#_{\{33B\}}=\#_{\{44\}}=0$. Now we only need to prove (\ref{eq.alg_2connup_1}).

 Consider the set $\Lc$ of all the maximal ov-segments $\lambda$ that contain at least one atom in $\Mb_U$ and at least one atom in $\Mb_{\mathrm{2conn}}$. Clearly each atom in $\Mb_{\mathrm{2conn}}^j\,(j\in\{1,2\})$ belongs to exactly $j$ ov-segments in $\Lc$; we define $\lambda\sim\lambda'$ if they intersects at an atom in $\Mb_{\mathrm{2conn}}^2$ (note that they cannot intersect at two atoms in $\Mb_{\mathrm{2conn}}^2$ due to Proposition \ref{prop.mol_axiom} (\ref{it.axiom1})). There cannot exist $\lambda_1,\cdots,\lambda_q\in\Lc$ such that $\lambda_j\sim\lambda_{j+1}$ with $\lambda_{q+1}=\lambda_1$ (otherwise, say $\lambda_j$ intersects $\lambda_{j+1}$ at atom $\nf_j\in \Mb_{\mathrm{2conn}}^2$, then the ov-segments between $\nf_j$ and $\nf_{j+1}$ along $\lambda_{j+1}$ forms a cycle in $\Mb_D$, contradicting that $\Mb_D$ is a tree). In other words, if we define the auxiliary graph $G$ with G-nodes $\lambda\in\Lc$ and G-edges $(\lambda,\lambda')$ where $\lambda\sim\lambda'$, then this $G$ is a forest, see {\color{blue}Figure \ref{fig.aux_graph}}.
 \begin{figure}[h!]
    \centering
    \includegraphics[width=0.36\linewidth]{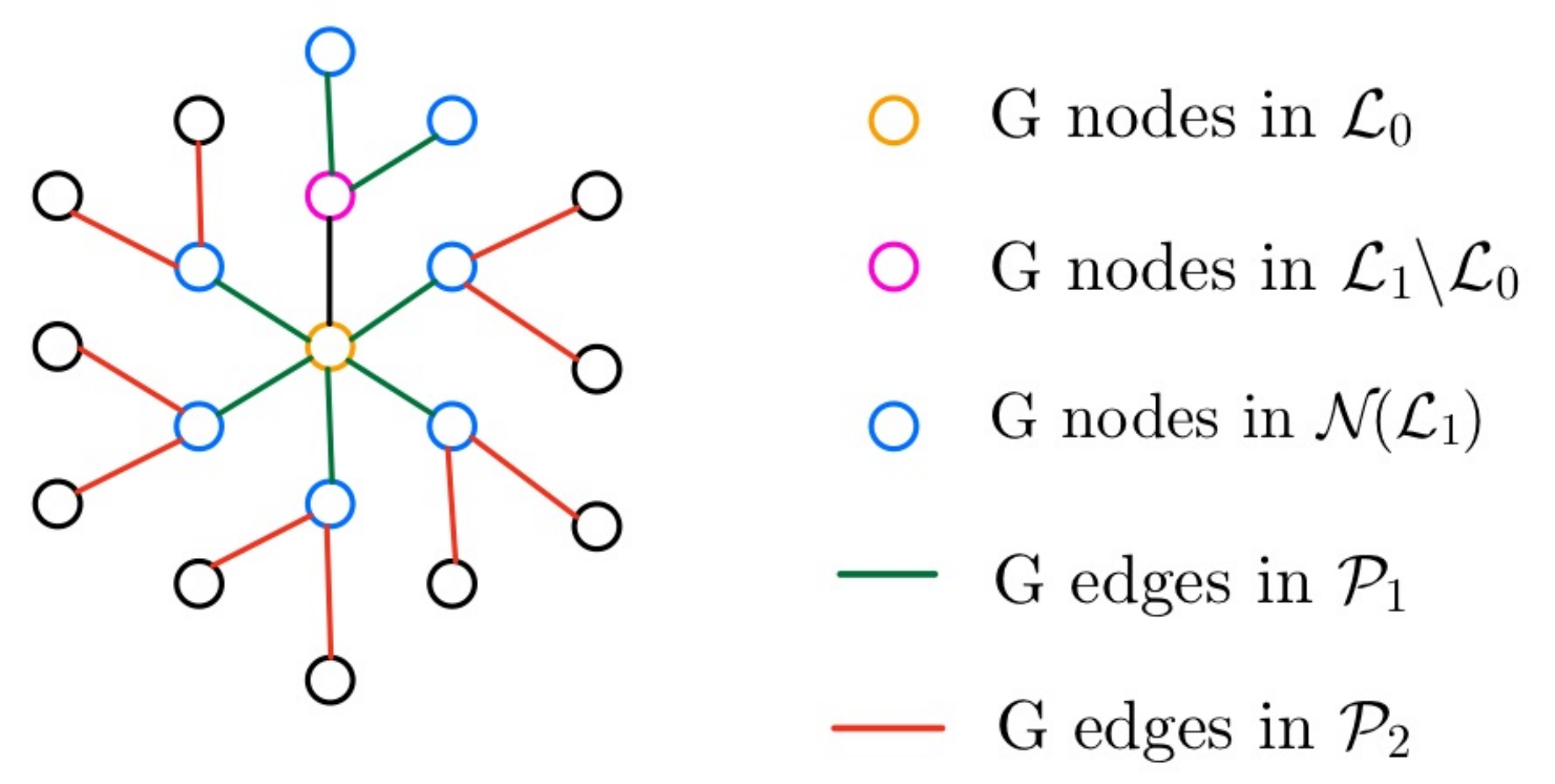}
    \caption{The auxiliary graph $G$. The sets $\Lc_0$, $\Lc_1$, $\Nc(\Lc_1)$, and $\Pc_1$, $\Pc_2$ are shown. Here $|\Lc_0|=1$, $|\Lc_1|=2$, $|\Nc(\Lc_1)|=7$, $|\Pc_1|=6$ and $|\Pc_2|=10$; they satisfy (\ref{eq.alg_2connup_2}) with $A=6$.}
    \label{fig.aux_graph}
\end{figure}
 
Recall by Definition \ref{def.cutting} that, if any atom (or any two atoms) in $\lambda\in\Lc$ is cut as free, then we either get a \{33A\} molecule, or get a fixed end at each remaining atom in $\lambda$. In this case we say $\lambda$ is \emph{broken}. Note that each cutting in Definition \ref{def.alg_2connup} (\ref{it.alg_2connup_1})--(\ref{it.alg_2connup_4}) breaks at most two $\lambda$, and each cutting in Definition \ref{def.alg_2connup} (\ref{it.alg_2connup_2}) does not break any $\lambda$. Moreover, if a cutting in Definition \ref{def.alg_2connup} (\ref{it.alg_2connup_2})--(\ref{it.alg_2connup_4}) breaks some $\lambda$, then this breaking must be due to cutting an atom in $\Mb_U$ (this is because, in Definition \ref{def.alg_2connup} (\ref{it.alg_2connup_4}), if we cut $\{\mf,\pf\}$ as free, then this is equivalent to first cutting $\mf$ and then cut $\pf$, and $\pf$ will have deg 2 when it is cut, so cutting $\pf$ will not break any $\lambda\in\Lc$). In particular, if a cutting in Definition \ref{def.alg_2connup} (\ref{it.alg_2connup_2})--(\ref{it.alg_2connup_4}) breaks both $\lambda$ and $\lambda'$, then they must intersect at an atom in $\Mb_U$, so we must have $\lambda\not\sim\lambda'$ by Proposition \ref{prop.mol_axiom} (\ref{it.axiom1}).

Now suppose $\lambda\sim\lambda'$, and $\lambda$ is broken before $\lambda'$, and $\lambda'$ is broken in a cutting operation in Definition \ref{def.alg_2connup} (\ref{it.alg_2connup_4}), then this cutting operation must involve some $\mf\in\Mb_U$ that belongs to $\lambda'$. Suppose $\mf\in S_\nf$ is as in Definition \ref{def.alg_2connup} (\ref{it.alg_2connup_4}) \emph{and does not have deg 4}, since $\Mb_U$ is a forest, we know that $\mf$ has at most one parent in $S_\nf$, so $\mf$ must have deg 3 when it is cut in Definition \ref{def.alg_2connup} (\ref{it.alg_2connup_4}) (because $\mf$ has two bottom non-fixed edges and one top non-fixed edge corresponding to its parent not in $S_\nf$, just as before). Let $\pf=\lambda\cap\lambda'\in \Mb_{\mathrm{2conn}}^2$, then the previously broken $\lambda$ provides a fixed end at $\pf$, so $\pf$ also has deg 3 at the time when $\mf$ is cut. By construction in Definition \ref{def.alg_2connup} (\ref{it.alg_2connup_4}), we get a \{33A\} molecule with higher atom $\mf$, such that $\mf\in\lambda'$. In the same way, if $\lambda$ contains an atom in $\Mb_{\mathrm{conn}}^1$ and is broken in a cutting operation in Definition \ref{def.alg_2connup} (\ref{it.alg_2connup_4}) with $\mf$ not having deg 4, then we also get a \{33A\} molecule with higher atom $\mf$, such that $\mf\in\lambda$.

\textbf{Proof part 2.} Let $\Lc_0\subseteq \Lc$ be the set of ov-segments selected in Definition \ref{def.alg_2connup} (\ref{it.alg_2connup_1}), and $\Lc_1\supseteq \Lc_0$ be the set of all ov-segments $\lambda$ that are broken in Definition \ref{def.alg_2connup} (\ref{it.alg_2connup_1}), then we have $|\Lc_1|\leq 2|\Lc_0|$. Recall the auxiliary graph $G$ defined above; let $\Nc(\Lc_1)$ be the set of G-nodes in $\Lc\backslash \Lc_1$ that are adjacent to a G-node in $\Lc_1$, and let $\Pc_1$ be the set of G-edges connecting a G-node in $\Lc_1$ to a G-node in $\Lc\backslash\Lc_1$. We claim that
\begin{equation}\label{eq.alg_2connup_2}|\Pc_1|\leq |\Nc(\Lc_1)|+|\Lc_1|,\quad |\Pc_1|\geq A|\Lc_0|-2|\Lc_1|\geq (A/2-2)|\Lc_1|.\end{equation} Here, the first inequality in (\ref{eq.alg_2connup_2}) follows because the subgraph $\Nc(\Lc_1)\cup\Lc_1\subseteq G$ of the auxiliary graph contains no cycle (and is a forest); the second inequality is because each G-node $\lambda\in\Lc_0$ has at least $A$ different G-edges, and that each G-edge in $\Lc_1$ is counted twice when adding up the G-degrees of G-nodes in $\Lc_1$, and their number is bounded by $|\Lc_1|$. Moreover, let $(\lambda_1,\cdots,\lambda_s)$ be the permutation of $\lambda\in\Lc\backslash\Lc_1$ ordered by time of breaking in Definition \ref{def.alg_2connup} (\ref{it.alg_2connup_2})--(\ref{it.alg_2connup_4}), then by construction in Definition \ref{def.alg_2connup} (\ref{it.alg_2connup_1}), each $\lambda_j$ is adjacent to at most $A$ different G-nodes $\lambda_k\,(k\neq j)$ in $G$. Let $\Pc_2$ be the set of G-edges among these $\lambda_j$, then $|\Pc_1|+|\Pc_2|\geq |\Mb_{\mathrm{2conn}}^2|-|\Lc_1|$ (where $|\Mb_{\mathrm{2conn}}^2|$ is the total number of G-edges, as each G-edge corresponds to an intersection of two $\lambda_j$).

Now, by the discussion in Part 1 above, for each $\lambda\in\Nc(\Lc_1)$, and each $\lambda$ containing an atom in $\Mb_{\mathrm{2conn}}^1$, and each $\lambda=\lambda_j$ that is adjacent to some $\lambda_k\,(k<j)$ in $G$, if $\lambda$ is broken in a cutting operation in Definition \ref{def.alg_2connup} (\ref{it.alg_2connup_4}) with $\mf$ not having deg 4, we get a \{33A\} molecule corresponding to $\lambda$ when it is broken, and each \{33A\} molecule can be so obtained from at most two such $\lambda$ (as shown in Part 1, this $\lambda$ must contain the higher atom $\mf$ in this \{33A\} molecule, so it has at most 2 choices). In the case $\lambda=\lambda_j$, each $j$ also corresponds to at most $A$ different G-edges $(\lambda_j,\lambda_k)$. This implies (using also the first inequality in (\ref{eq.alg_2connup_2})) that
\begin{equation}\label{eq.alg_2connup_3}\#_{\{33A\}}\geq \frac{1}{2}\max\bigg(|\Pc_1|-|\Lc_1|,\,|\Mb_{\mathrm{2conn}}^1|,\,A^{-1}|\Pc_2|\bigg)-\#_{\mathrm{exc}}.
\end{equation} Here we use $|\Pc_1|-|\Lc_1|\leq |\Nc(\Lc_1)|$, and that the number of $\lambda=\lambda_j$ adjacent to some $\lambda_k\,(k<j)$ in $G$ is at least $A^{-1}|\Pc_2|$. Moreover $\#_{\mathrm{exc}}$ is the number of atoms cut in Definition \ref{def.alg_2connup} (\ref{it.alg_2connup_1}) plus the number of deg 4 atoms cut in Definition \ref{def.alg_2connup} (\ref{it.alg_2connup_4}) (each of them may break two $\lambda$ which count as exceptions; also cutting any deg 2 atom in Definition \ref{def.alg_2connup} (\ref{it.alg_2connup_2}) will not break any $\lambda$).

On the other hand, the number of deg 4 atoms cut in Definition \ref{def.alg_2connup} (\ref{it.alg_2connup_4}) does not exceed the number of full components $\#_{\mathrm{fucp}(\Mb_U)}$ of $\Mb_U$, and also \[\max(\#_{\mathrm{exc}},\#_{\{4\}})\leq \#_{\mathrm{fucp}(\Mb_U)}+|\Lc_1|,\] because (i) any new component of $\Mb_U$ generated in the cutting sequence cannot be full (cf. Remark \ref{rem.full_cut}, and (ii) any component of $\Mb_D$ also cannot be full when all atoms in $\Mb_U$ are cut. By combining this with the second inequality in (\ref{eq.alg_2connup_2}) and (\ref{eq.alg_2connup_3}), and the inequality $|\Pc_1|+|\Pc_2|\geq |\Mb_{\mathrm{2conn}}^2|-|\Lc_1|$ and using the definition $A:=32\upsilon^{-2}$, we conclude that (recall that all \{33A\} molecules are viewed as good)
\begin{equation*}
\begin{aligned}
(\upsilon/2)\cdot(\#_{\{33A\}})-d\cdot(\#_{\{4\}})&\geq (\upsilon/4)\max\big(|\Pc_1|-|\Lc_1|,|\Mb_{\mathrm{2conn}}^1|,A^{-1}|\Pc_2|\big)-(d+1)(\#_{\mathrm{fucp}(\Mb_U)}+|\Lc_1|)\\
&\geq (\upsilon/16)(|\Pc_1|-8\upsilon^{-2}|\Lc_1|)+(\upsilon/16)(|\Mb_{\mathrm{2conn}}^1|+A^{-1}|\Pc_2|)-(d+1)\#_{\mathrm{fucp}(\Mb_U)}\\
&\geq(\upsilon/16)A^{-1}(|\Pc_1|+|\Lc_1|)+(\upsilon/16)(|\Mb_{\mathrm{2conn}}^1|+A^{-1}|\Pc_2|)-(d+1)\#_{\mathrm{fucp}(\Mb_U)}\\
&\geq (\upsilon/16)A^{-1}(|\Mb_{\mathrm{2conn}}^1|+|\Mb_{\mathrm{2conn}}^2|)-(d+1)\#_{\mathrm{fucp}(\Mb_U)},
\end{aligned}
\end{equation*} which proves (\ref{eq.alg_2connup_1}).
\end{proof}
Now we can prove the first case of Proposition \ref{prop.case5}, namely Proposition \ref{prop.comb_est_case6}. This corresponds to toy model III in Section \ref{sec.toy3}, where there exist (many components of $\Mb_D$ or) many 2-connections in $\Mb_U$.
\begin{proposition}
\label{prop.comb_est_case6} Proposition \ref{prop.case5} is true if either $\#_{\mathrm{comp}(\Mb_D)}\geq (C_{3}^*)^{-1}\cdot\rho'''$, or if $\#_{\mathrm{comp}(\Mb_D)}\leq (C_{3}^*)^{-1}\cdot\rho'''$ and $|\overline{\Mb}_{\mathrm{2conn}}^2|\geq (C_{2}^*)^{-1}\cdot\rho'''$, where $\overline{\Mb}_{\mathrm{2conn}}^2$ is defined as in Definition \ref{def.alg_2connup}, but with the roles of $\Mb_U$ and $\Mb_D$ reversed.
\end{proposition}
\begin{proof} Our goal is to prove that
\begin{equation}\label{eq.case51_copy} (\upsilon/2)\cdot(\#_{\mathrm{good}})-d\cdot(\#_{\{4\}})\geq (C_{4}^*)^{-1}\cdot \rho'''.
\end{equation}

\textbf{Proof part 1.} Suppose $\#_{\mathrm{comp}(\Mb_D)}\geq (C_{3}^*)^{-1}\cdot\rho'''$. Recall the definition of $\Mb_{\mathrm{2conn}}^j$ for $j\in\{1,2\}$ in Definition \ref{def.alg_2connup} (where we do not reverse the roles of $\Mb_U$ and $\Mb_D$, which is different from $\overline{\Mb}_{\mathrm{2conn}}^2$ defined above). Consider all components of $\Mb_D$, if at least half of them contain an atom in $\Mb_{\mathrm{2conn}}:=\Mb_{\mathrm{2conn}}^1\cup \Mb_{\mathrm{2conn}}^2$, then we have $|\Mb_{\mathrm{2conn}}|\geq \#_{\mathrm{comp}(\Mb_D)}/2$. In this case we apply \textbf{2CONNUP} to $\Mb_{UD}$, by Proposition \ref{prop.alg_2connup} we have
\begin{equation}\label{eq.comb_est_case61}
(\upsilon/2)\cdot(\#_{\mathrm{good}})-d\cdot(\#_{\{4\}})\geq C^{-1}|\Mb_{\mathrm{2conn}}|-C\cdot \#_{\mathrm{fucp}(\Mb_U)}\geq\big(C^{-1}(C_{3}^*)^{-1}-C(C_{5}^*)^{-1}\big)\cdot\rho'''
\end{equation} (where the upper bound of $\#_{\mathrm{fucp}(\Mb_U)}$ follows from Reduction \ref{red.2} (ii)), which implies (\ref{eq.case51_copy}).

If at least half of the components of $\Mb_D$ do not contain an atom in $\Mb_{\mathrm{2conn}}$, then we shall cut $\Mb_U$ as free and cut it into elementary molecules using \textbf{UP}. Then, for each of these at least half components $X\subseteq \Mb_D$, it will not contain any deg 2 atom (because this component does not contain any atom in either $\Mb_{\mathrm{2conn}}^1$ or $\Mb_{\mathrm{2conn}}^2$); moreover, this $X$ may be divided into new components due to the broken ov-segments, but each new component will contain at least one deg 3 atom (cf. Remark \ref{rem.full_cut}), and at least one new component will contain at least two deg 3 atoms. This last fact is because at least two atoms in $X$ will become deg 3 after cutting $\Mb_U$ (i.e. has deg 3 or has a bond with $\Mb_U$) by Reduction \ref{red.2} (ii); by considering a path in $X$ between these two atoms before cutting and whether or where this path breaks after cutting, we can obtain a component of $X$ after cutting that contains two deg 3 atoms (as the atom where the path breaks must have a fixed end). 

We then cut each resulting component of $\Mb_D$ into elementary molecules using \textbf{DOWN}. By Proposition \ref{prop.alg_up}, and note that all \{33A\} molecules are treated as good, we get that $\#_{\mathrm{good}}\geq \#_{\mathrm{comp}(\Mb_D)}/2$ and $\#_{\{4\}}\leq\#_{\mathrm{fucp}(\Mb_U)}$ (this is because $\Mb_D$ has no full component after cutting $\Mb_U$, due to Remark \ref{rem.full_cut} and Reduction \ref{red.2} (ii)), so (\ref{eq.case51_copy}) follows in the same way as (\ref{eq.comb_est_case61}).

\textbf{Proof part 2.} Assume $\#_{\mathrm{comp}(\Mb_D)}\leq (C_{3}^*)^{-1}\cdot\rho'''$ and $|\overline{\Mb}_{\mathrm{2conn}}^2|\geq (C_{2}^*)^{-1}\cdot\rho'''$. We first cut as free each full component in $\Mb_U$ that has no bond connecting to $\Mb_D$ (and then cut it into elementary molecules using \textbf{UP}), to make $\Mb_{UD}$ satisfy the assumptions of Definition \ref{def.alg_2connup} after reversing the roles of $\Mb_U$ and $\Mb_D$. Then we apply \textbf{2CONNDN} to $\Mb_{UD}$. Proposition \ref{prop.alg_2connup} then gives that (note $\#_{\mathrm{fucp}(\Mb_D)}\leq\#_{\mathrm{comp}(\Mb_D)}$)
\begin{multline}\label{eq.comb_est_case62}
(\upsilon/2)\cdot(\#_{\mathrm{good}})-d\cdot(\#_{\{4\}})\geq C^{-1}\cdot|\overline{\Mb}_{\mathrm{2conn}}^2|-C\cdot(\#_{\mathrm{comp}(\Mb_D)}+\#_{\mathrm{fucp}(\Mb_U)})\\
\geq\big(C^{-1}(C_{2}^*)^{-1}-C(C_{3}^*)^{-1}\big)\cdot\rho''',
\end{multline} which implies (\ref{eq.case51_copy}). Here, we used the bound on $\#_{\mathrm{fucp}(\Mb_U)}$ provided by Reduction \ref{red.2} (ii).
\end{proof}
\begin{reduct}
\label{red.3} With Proposition \ref{prop.comb_est_case6}, in the remaining proof of Proposition \ref{prop.case5}, we may now assume for $\Mb_{UD}$, in addition to Reduction \ref{red.2} (i)--(iii):
\begin{itemize}
\item That $\#_{\mathrm{comp}(\Mb_D)}\leq (C_{3}^*)^{-1}\cdot\rho'''$, and the number of weakly degenerate atoms and pairs in $\Mb_D$ is less than $(C_{2}^*)^{-1}\cdot\rho'''$.
\item We have $|\overline{\Mb}_{\mathrm{2conn}}^2|\leq (C_{2}^*)^{-1}\cdot\rho'''$, with the set $\overline{\Mb}_{\mathrm{2conn}}^2$ defined in Proposition \ref{prop.comb_est_case6}.
\end{itemize}
\end{reduct}
\subsection{The \textbf{3COMPUP} algorithm}\label{sec.3compup}
\begin{figure}[h!]
    \centering
    \includegraphics[width=0.5\linewidth]{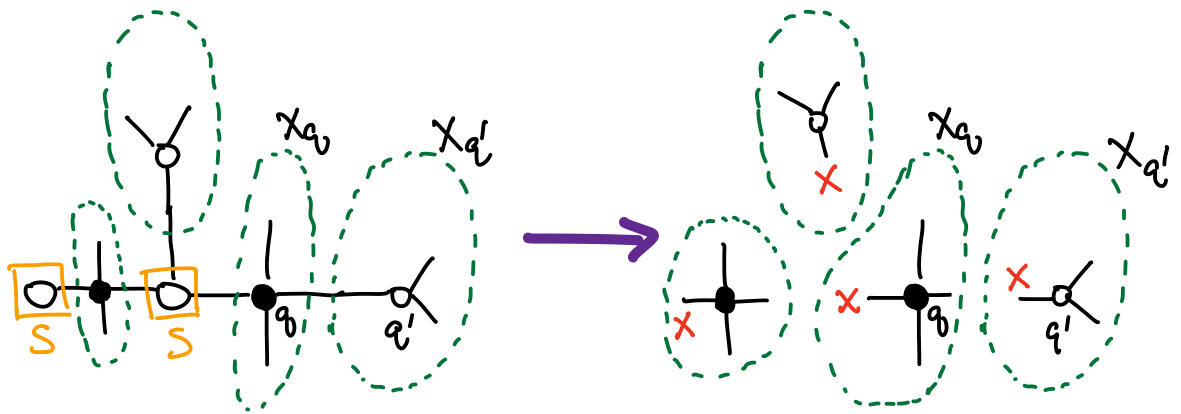}
    \caption{Illustration of Lemma \ref{lem.cutconnected}. The main differences with Lemma \ref{lem.cutconnected0} are that (i) now $S$ is allowed to be ov-connected instead of connected (with possible O-atoms in between which also belong to different components after cutting $S$ as free), and (ii) each ov-adjacent atom $\qf$ (not necessarily adjacent) now leads to its own component. For example, the $\qf$ and $\qf'$ will belong to different components after cutting $S$ as free, but again each component has exactly one new fixed end.}
    \label{fig.heart2}
\end{figure}

The rest of this section is devoted to the proof of Proposition \ref{prop.case5} under the assumption of Reductions \ref{red.2} and \ref{red.3}. This proof will be presented,  following \textbf{Stages 1--6}, in Section \ref{sec.finish}. Before this, we first introduce the remaining ingredients needed in the proof; in this subsection we define the full \textbf{3COMPUP} algorithm (see Definition \ref{def.3comp_alg}) corresponding to toy model II ((see \textbf{Stage 5} in Section \ref{sec.finish}), which is a direct extension of the toy version in Definition \ref{def.toy2_alg} by including O-atoms and translating into the ov-segment language (cf. Section \ref{sec.reduce1}). Apart from this translation, Definition \ref{def.3comp_alg} and Proposition \ref{prop.3comp} are exactly the same as Definition \ref{def.toy2_alg} and Proposition \ref{prop.toy2} in Section \ref{sec.toy}.
Before proceeding, we state and prove the full version of Lemma \ref{lem.cutconnected0} in Section \ref{sec.toy} in the presence of O-atoms (by translating adjacency to ov-adjacency etc.).
\begin{lemma}\label{lem.cutconnected} Let $\Mb$ be a molecule that is a tree, and $S\subseteq\Mb$ is an \emph{ov-connected} set. Suppose we cut $S$ as free from $\Mb$. Then, for each atom $\qf\in\Mb\backslash S$ that is \emph{ov-adjacent} to an atom in $S$ before cutting, there is a unique component $X_\qf\subseteq \Mb\backslash S$ after cutting that contains $\qf$. This $X_\qf$ is in \emph{bijection} with $\qf$, and when viewed in $X_\qf$, the operation of cutting $S$ creates \emph{only one} fixed end in $X_\qf$ at $\qf$, see {\color{blue}Figure \ref{fig.heart2}}.
\end{lemma}
\begin{proof} Consider all the ov-segments connecting an atom in $S$ to another atom (in $S$ or not in $S$), and all atoms (including the intermediate O-atoms) on these ov-segments. Let the set of all such ov-segments (and the bonds contained in these ov-segments) be $\Lambda$, and the set of all these atoms be $Z\supseteq S$, then the atoms in $Z\backslash S$ are precisely all the atoms $\qf$ as stated. Note that all the bonds in $\Lambda$ with an endpoint in $\Mb\backslash S$ are broken after cutting $S$ as free; we claim that at this time, any two different atoms $\qf_1,\qf_2\in Z\backslash S$ must be in two different components of $\Mb\backslash S$. In fact, if not, then they are connected in $\Mb\backslash S$ by a path of bonds \emph{not containing any bond in $\Lambda$}. However, $\qf_1,\qf_2\in Z$ are also connected by another path \emph{containing only bonds in $\Lambda$}, using that $S$ is ov-connected. This then produces a cycle in $\Mb$, contradiction.

Now we have proved that $\qf\mapsto X_{\qf}$ is a bijection, and for each $X_\qf$, the only atom in it that is ov-adjacent to an atom in $S$ is $\qf$. As such, cutting $S$ as free only creates fixed end at $\qf$ and not at any other atom in $X_\qf$. If we get \emph{more than one} fixed end at $\qf$, then $\qf$ must be ov-adjacent to two atoms in $\nf_1,\nf_2\in S$ along \emph{two different} ov-segments (note: breaking two bonds at $\qf$ that belong to the same ov-segment creates only one fixed end at $\qf$, see Definition \ref{def.cutting}). Note that $\nf_1$ is connected to $\nf_2$ by a path of bonds in $\Lambda$; since $\qf\not\in S$, we know that $\qf$ either does not occur in this path, or occurs in this path but as one of the intermediate O-atoms (so the two bonds at $\qf$ in this path \emph{must be serial}). However, the two ov-segments connecting $\qf$ to $\nf_j$ also form a path between $\nf_1$ and $\nf_2$, in which the two bonds at $\qf$ \emph{are not serial}. Thus these two paths cannot coincide, which again produces a cycle in $\Mb$. This contradiction then proves that cutting $S$ as free produces only one fixed end at $\qf$ (and in $X_\qf$), as desired.
\end{proof}
\begin{lemma}\label{lem.3comp_aux} Suppose $\Mb$ is a molecule that is a tree, contains only deg 3 and 4 atoms and has no C-top fixed end. Let the set of deg 3 atoms in $\Mb$ be $\Mb_{(3)}$ and assume $|\Mb_{(3)}|\geq 2$. Define a set $A\subseteq\Mb$ to be \textbf{no-bottom} if any child of any atom in $A$ must also be in $A$. Then, there always exists a connected, no-bottom subset $A\subseteq \Mb$ with $|A\cap \Mb_{(3)}|=2$.
\end{lemma}
\begin{proof} For any $\nf\in\Mb$, let $S_\nf$ be the set of descendants of $\nf$, then it is always connected and no-bottom. If there exists $\nf$ such that $|S_\nf\cap \Mb_{(3)}|\geq 2$, then choose a lowest such $\nf$. Since $\nf$ is either O-atom or has no top fixed end, we know that either (i) $\nf$ has one child $\nf_1$, or (ii) $\nf$ has deg 4 and has two children $\nf_1,\nf_2$. In either case, by assumption we have $|S_{\nf_j}\cap \Mb_{(3)}|\leq 1$ for each $j$, which then implies $|S_\nf\cap \Mb_{(3)}|\leq 2$ because $S_\nf=(\cup_j S_{\nf_j})\cup\{\nf\}$. This implies $|S_\nf\cap \Mb_{(3)}|= 2$, as desired.

Now assume $|S_\nf\cap \Mb_{(3)}|\leq 1$ for each $\nf\in\Mb$. If $|S_\nf\cap \Mb_{(3)}|=1$, there is a unique $\pf=\pf(\nf)\in S_\nf\cap \Mb_{(3)}$; now consider all atom pairs $(\nf,\nf')$ such that $\pf(\nf)$ and $\pf(\nf')$ both exist and \emph{are different} (which exist because $|\Mb_{(3)}|\geq 2$), and one such pair $(\nf,\nf')$ that are connected by \emph{a shortest path} in $\Mb$. Let this path be $(\nf_0,\nf,\cdots,\nf_j=\nf')$, then for each $0<i<j$ we must have either $S_{\nf_i}\cap \Mb_{(3)}=\varnothing$ or $\pf(\nf_i)=\pf(\nf)$ (otherwise the pair $(\nf,\nf_i)$ corresponds to a shorter path, contradiction). Now let $A=\cup_{i=0}^j S_{\nf_i}$, it is clear that $A$ is connected, no-bottom and $A\cap \Mb_{(3)}=\{\pf(\nf),\pf(\nf')\}$, as desired. 
\end{proof}
\begin{definition}[The \textbf{3COMPUP} algorithm]\label{def.3comp_alg} Let $\Mb$ be a molecule which is a forest, contains only deg 3 and 4 atoms and has no C-top fixed end, and assume that each component of $\Mb$ contains at least two deg 3 atoms. Let the set of deg 3 atoms in $\Mb$ be $\Mb_{(3)}$. Define the following cutting sequence:
\begin{enumerate}[{(1)}]
\item \label{it.3comp_1} Apply Lemma \ref{lem.3comp_aux} to one component of $\Mb$ and select a connected no-bottom subset $A$. This set is fixed until the end of \ref{it.3comp_5}.
\item \label{it.3comp_3} Choose a lowest deg 3 atom $\nf$ in $A$ that has not been cut. Let $S_\nf$ be the set of descendants of $\nf$ (note that $S_\nf\subseteq A$). This $\nf$ and $S_\nf$ are fixed until the end of \ref{it.3comp_4}.
\item \label{it.3comp_4} Starting from $\nf$, each time choose a highest atom $\mf$ in $S_\nf$ that has not been cut. If $\mf$ has deg 3 and has an ov-parent $\mf^+$ or ov-child $\mf^-$ that also has deg 3, then cut $\{\mf,\mf^\pm\}$ as free; otherwise cut $\mf$ as free. Repeat until all atoms in $S_\nf$ have been cut. 
\item \label{it.3comp_4+} If $\Mb$ has deg 2 atom, then cut it; repeat until $\Mb$ has no deg 2 atoms. Then go to \ref{it.3comp_3}.
\item \label{it.3comp_5} When all atoms in $A$ have been cut, we subsequently cut all remaining deg 2 atoms in $\Mb$ as in \ref{it.3comp_4+}, and then cut all \emph{components of $\Mb$ with only one deg 3 atom} using \textbf{UP}. 
\item\label{it.t3comp_6} Go to \ref{it.3comp_1} and select the next $A$ from one of the remaining components of $\Mb$, and so on.
\end{enumerate}
\end{definition}
\begin{proposition}\label{prop.3comp} Let $\Mb$ be as in Definition \ref{def.3comp_alg}. Define $\#_{\mathrm{3ovcp}(\Mb)}$ to be the number of ov-components of $\Mb_{(3)}$ (the det of deg 3 atoms in $\Mb$), and $\#_{\mathrm{comp}(\Mb)}$ to be the number of components of $\Mb$. Then $\Mb$ is cut into elementary molecules after \textbf{3COMPUP} (Definition \ref{def.3comp_alg}), and we have
\begin{equation}\label{eq.alg_3ovcpup_1}\#_{\{33B\}}=\#_{\{44\}}=\#_{\{4\}}=0;\qquad
\#_{\{33A\}}\geq\frac{1}{10}\cdot\#_{\mathrm{3ovcp}(\Mb)}+\frac{1}{5}\cdot\#_{\mathrm{comp}(\Mb)}.
\end{equation}
\end{proposition}
\begin{proof} Note that steps \ref{it.3comp_3}--\ref{it.3comp_4} in Definition \ref{def.3comp_alg} are exactly the same as those \textbf{UP} (see Definition \ref{def.alg_up}  (\ref{it.alg_up_2})--(\ref{it.alg_up_3}), note also the step of cutting deg 2 atoms is done in Definition \ref{def.3comp_alg} \ref{it.3comp_5}), but applied to $A$ only. Since $A$ is no-bottom which means $S_\nf\subseteq A$ whenever $\nf\in A$, it is easy to see (following the same proof as Proposition \ref{prop.alg_up}) that the (MONO) property is preserved, and thus $\Mb$ is cut into elementary molecules only, and $\#_{\{33B\}}=\#_{\{44\}}=\#_{\{4\}}=0$. We now focus on the proof of (\ref{eq.alg_3ovcpup_1}).

Note that the algorithm \textbf{3COMPUP} consists of finitely many loops where one loop is formed by steps \ref{it.3comp_1}--\ref{it.3comp_5} in Definition \ref{def.3comp_alg}, and (say) turns $\Mb$ into $\Mb'$. Define $\lambda(\Mb)$ to be the right hand side of the last inequality in (\ref{eq.alg_3ovcpup_1}), then we only need to consider each single loop \ref{it.3comp_1}--\ref{it.3comp_5} in Definition \ref{def.3comp_alg}, and prove that:
\begin{itemize}
\item[($\clubsuit$)] If $\Mb$ is a forest, has no top fixed end and contains only deg 3 and 4 atoms with each component containing at least two deg 3 atoms, then the result $\Mb'$ after executing the loop \ref{it.3comp_1}--\ref{it.3comp_5} in Definition \ref{def.3comp_alg} satisfies the same property. Moreover for this loop we have
\begin{equation}\label{eq.3comp_2}\#_{\{33\}}\geq\lambda(\Mb)-\lambda(\Mb').
\end{equation}
\end{itemize} Note that adding up (\ref{eq.3comp_2}) leads to the desired lower bound of $\#_{\{33\}}$, thus it suffices to prove ($\clubsuit$) for a given loop. We divide the proof into 4 parts, which exactly correspond to the 4 parts of the proof in toy model II (Proposition \ref{prop.toy2}).

\textbf{Proof part 1.} We prove that each component of $\Mb'$ satisfies the desired properties. This is obvious, as the lack of top fixed ends follows from the same proof as Proposition \ref{prop.alg_up}, and the lack of deg 2 atoms and components with only one deg 3 atom follows from the cuttings done in Definition \ref{def.3comp_alg} \ref{it.3comp_5}.

\textbf{Proof part 2.} Now we need to prove (\ref{eq.3comp_2}). We start by making some preparations. Note that $A$ is an ov-connected subset of $\Mb$; by Lemma \ref{lem.cutconnected}, each atom $\mf\in \Mb\backslash A$ that is ov-adjacent to an atom $\nf\in A$ by a pre-loop ov-segment corresponds to a unique component $X_\mf$ after cutting the atoms in $A$ as free (the loop will in fact cut as free more atoms). Note also that $\mf$ must be an ov-parent of $\nf$ as $A$ has no bottom.

Let $\Mb_{(3)}$ be the set of deg 3 atoms in $\Mb$ (before executing the loop), and consider all ov-components of $\Mb_{(3)}$ that intersect $A$ or are ov-adjacent to $A$ by an ov-segment (we call them \textbf{good}). We know $A$ has exactly two atoms in $\Mb_{(3)}$; apart from the one or two ov-components that contain these two atoms, for any atom in any other good ov-component, it must belong to $X_\mf$ for a unique $\mf$, and this $\mf$ must have deg 3 (otherwise the ov-component will not be ov-adjacent to $A$), with the corresponding $\nf$ having deg 4. This is similar to the scenario in toy model II, see {\color{blue}Figure \ref{fig.3comp_ov_1}}. Let $B$ be the union of all good ov-components and $A$. We claim that after steps \ref{it.3comp_1}--\ref{it.3comp_4+} in Definition \ref{def.3comp_alg} and cutting the deg 2 atoms in Definition \ref{def.3comp_alg} \ref{it.3comp_5}, exactly those atoms in $B$ have been cut as free.
\begin{figure}[h!]
    \centering
    \includegraphics[width=0.38\linewidth]{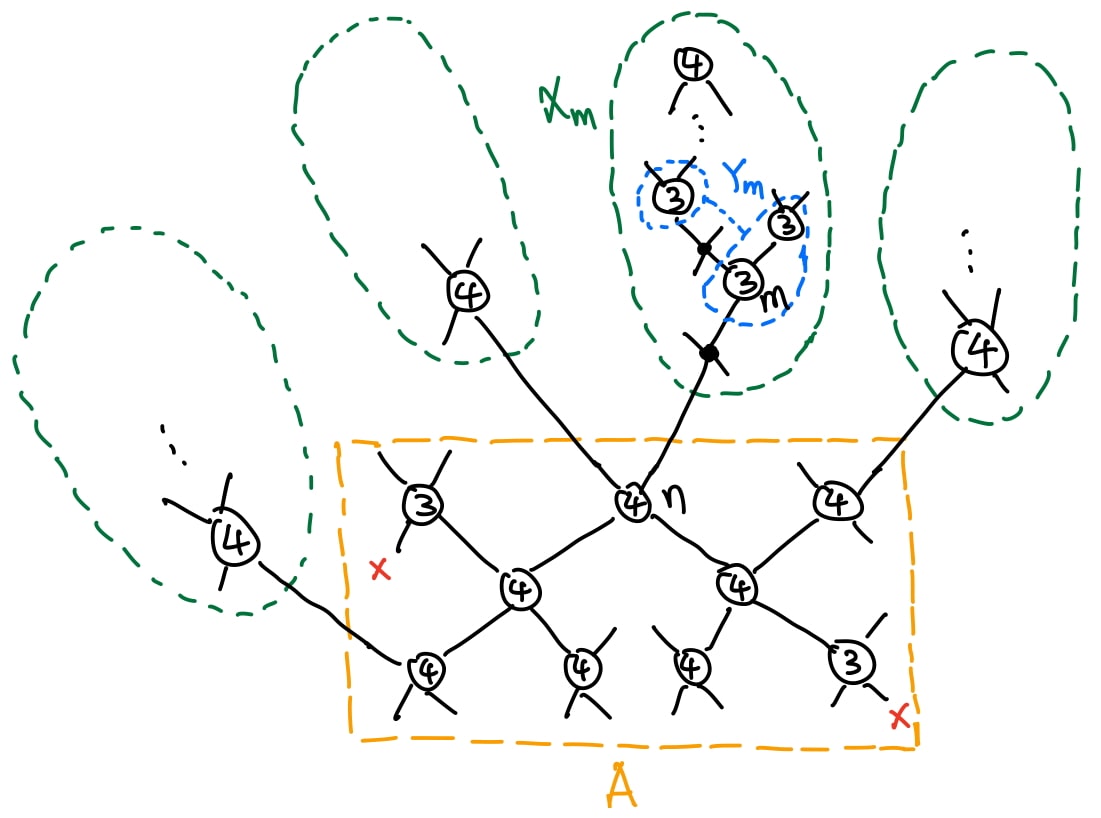}
    \caption{The set $A$ and ov-components of deg 3 atoms forming $B$, in Part 2 of the proof of Proposition \ref{prop.3comp}. This is basically the same as in toy model II ({\color{blue}Figure \ref{fig.setab}}), except now the notion of adjacency and components are replaced by ov-adjacency and ov-components.}
    \label{fig.3comp_ov_1}
\end{figure}

In fact, the atoms $A$ must have been cut by definition. For any atom in a good ov-component of $\Mb_{(3)}$, by definition it is connected to an atom in $A$ by a sequence of deg 3 atoms with each consecutive ones being ov-adjacent; once the atom in $A$ has been cut, the atoms on this chain will successively become deg 2, so they must have been cut (either in Definition \ref{def.3comp_alg} \ref{it.3comp_4+}--\ref{it.3comp_5}, or as part of \{33\} molecule in Definition \ref{def.3comp_alg} \ref{it.3comp_4}). On the other hand, each atom in $\Mb\backslash B$ has at most one ov-segment connecting to $B$ (as $B$ is ov-connected and $\Mb$ is a tree, cf. the proof of Lemma \ref{lem.cutconnected}), and must be deg 4 if it has one such ov-segment (otherwise it would belong to one of the good ov-components of $\Mb_{(3)}$), so each atom in $\Mb\backslash B$ will have deg 3 or 4 after cutting $B$ as free, thus it will not be cut in steps \ref{it.3comp_1}--\ref{it.3comp_4+} in Definition \ref{def.3comp_alg} or in cutting the deg 2 atoms in Definition \ref{def.3comp_alg} \ref{it.3comp_5}. 

\textbf{Proof part 3.} Next we obtain a lower bound on $\#_{\{33\}}$. As $A$ contains two deg 3 atoms, we can show that $\#_{\{33\}}\geq 1$ in the same way as Proposition \ref{prop.alg_up} (\ref{it.up_proof_3}). Moreover, since $A$ contains \emph{exactly} two deg 3 atoms, by arguing as in the proof of Proposition \ref{prop.alg_up} (\ref{it.up_proof_3}) and calculating the quantity $\sigma$ defined there but within $A$, we can show that, apart from at most two exceptions, each deg 4 atom in $A$ will have deg 3 when it is cut.

Now let $\lambda$ be the number of good ov-components of $\Mb_{(3)}$ (i.t. those forming $B$). By the discussions in Part 2 above, with at most two exceptions, each such ov-component must contain some $\mf$ which has deg 3 and is an ov-parent of an atom $\nf\in A$ of deg 4; consider then the maximal ov-segment $\sigma$ containing $\mf$ and $\nf$, note that each $\sigma$ corresponds to at most one good ov-component (if $\sigma$ is given, then any deg 3 atoms on $\sigma$ must be ov-adjacent to each other and thus belong to the same ov-component), so the number of such ov-segments $\sigma$ is at least $\lambda$. We claim that with at most two more exceptions, the first atom $\pf\in\sigma$ that is cut, must belong to a \{33A\} molecule. In fact, we must have $\pf\in A$ by Definition \ref{def.3comp_alg} \ref{it.3comp_4}; with at most two exceptions we may assume $\pf$ has deg 3 when it is cut. At this time the $\mf$ (which is ov-parent of $\pf$) also has deg 3 (because the only way $\mf$ is ov-connected to $A$ is by the ov-segment between $\mf$ and $\nf$, cf. Lemma \ref{lem.cutconnected}), so by Definition \ref{def.3comp_alg} \ref{it.3comp_4}, we know that $\pf$ must belong to a \{33\} molecule (which may or may not be $\{\pf,\mf\}$). Since each \{33\} molecule can be so obtained from at most two $\sigma$, we thus obtain that
\begin{equation}\label{eq.3comp_3}\#_{\{33\}}\geq\max\bigg(1,\frac{\lambda-2}{2}-2\bigg).\end{equation}

\textbf{Proof part 4.} Now consider the components of $\Mb'$ and ov-components of $\Mb_{(3)}'$ (set of deg 3 atoms in $\Mb'$). Note that $\Mb'$ is formed from $\Mb$ by cutting $B$ as free and then removing all components with only one deg 3 atom. Since $B$ is ov-connected, by applying Lemma \ref{lem.cutconnected} again, we see that each component $X_\qf$ of $\Mb'$ corresponds to a unique atom $\qf\in \Mb\backslash B$ that is ov-adjacent to an atom in $B$. Let $\mu$ be the number of such atoms (i.e. corresponding to those new components with at least two deg 3 atoms after cutting $B$), then \begin{equation}\label{eq.3comp_4}\#_{\mathrm{comp}(\Mb)}-\#_{\mathrm{comp}(\Mb')}=1-\mu.\end{equation}

Finally, consider the set of deg 3 atoms, i.e. $\Mb_{(3)}$ and $\Mb_{(3)}'$, and their ov-components. The ov-components of $\Mb_{(3)}$ includes those $\lambda$ ov-components forming $B$ as in Part 3 (which have been cut when forming $\Mb'$), and the other ov-components which are contained in some $X_\qf$ defined above. Note that for each new component $Y$ formed after cutting $B$ that contains only one deg 3 atom, it will be removed and will not occur in $\Mb'$; however, all atoms in $Y$ originally have deg 4 before any cutting, so initially $Y$ also has no contribution to $\#_{3\mathrm{ovcp}(\Mb)}$, and thus we can ignore such $Y$.

\begin{figure}[h!]
    \centering
    \includegraphics[width=0.43\linewidth]{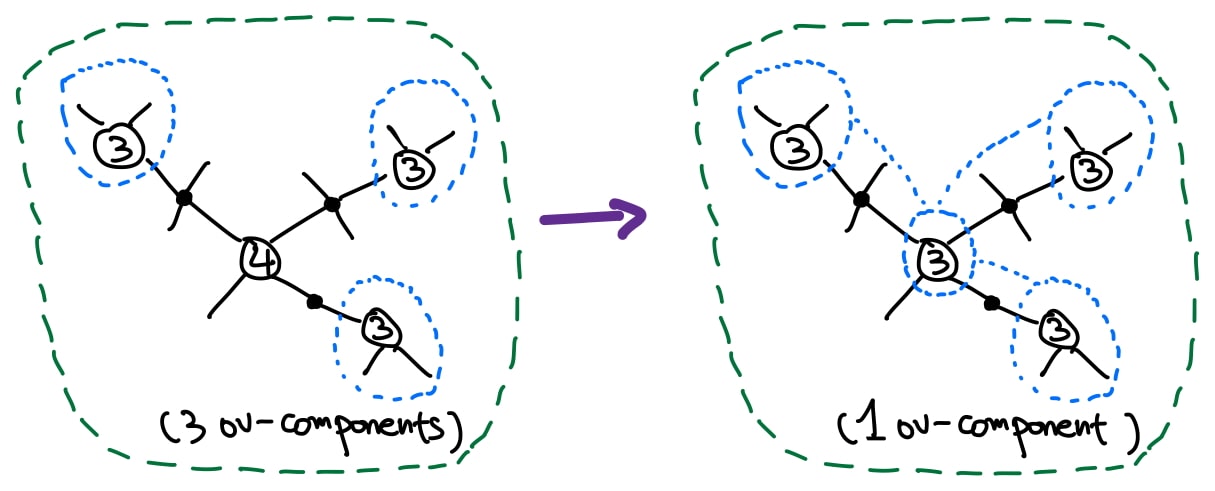}
    \caption{Part 4 of the proof of  Proposition \ref{prop.3comp}, where turning a deg 4 atom into deg 3 may decrease the number of ov-components of the deg 3 set by at most 2. This is exactly the same as the toy model II counterpart ({\color{blue}Figure \ref{fig.cutb}}), as long as we change ``components" to ``ov-components" etc.}
    \label{fig.3comp_ov_2}
\end{figure}

By Lemma \ref{lem.cutconnected}, when viewed in each $X_\qf$, the effect of cutting $B$ exactly turns one free end into fixed end (turning one deg 4 atom into deg 3). Therefore the part of $\Mb_{(3)}'$ within $X_\qf$ is formed by \emph{adding one deg 3 atom} into the existing part of $\Mb_{(3)}$ within $X_\qf$. This new deg 3 atom may \emph{decrease} $\#_{3\mathrm{ovcp}(X_\qf)}$ (i.e. the number of ov-components of the set of deg 3 atoms in $X_\qf$) by joining existing ov-components together, but this decrease is \emph{at most 2} because this new deg 3 atom can join at most 3 ov-components into 1, see {\color{blue} Figure \ref{fig.3comp_ov_2}}. From this, and accounting for the $\lambda$ ov-components removed when cutting $B$, we deduce that
\begin{equation}\label{eq.3comp_5}\#_{3\mathrm{ovcp}(\Mb)}-\#_{3\mathrm{ovcp}(\Mb')}\leq \lambda+2\mu.\end{equation}
Putting together (\ref{eq.3comp_3})--(\ref{eq.3comp_5}), and using the definition of $\lambda(\Mb)$, we deduce that
\begin{equation}\label{eq.3comp_6}
\begin{aligned}\lambda(\Mb)-\lambda(\Mb')&=\frac{1}{5}\big(\#_{\mathrm{comp}(\Mb)}-\#_{\mathrm{comp}(\Mb')}\big)+\frac{1}{10}\big(\#_{3\mathrm{ovcp}(\Mb)}-\#_{3\mathrm{ovcp}(\Mb')}\big)
\\&\leq \frac{1}{5}+\frac{\lambda}{10}\leq\max\bigg(1,\frac{\lambda-6}{2}\bigg)\leq \#_{\{33\}}.
\end{aligned}\end{equation} This proves (\ref{eq.3comp_2}) and completes the proof of Proposition \ref{prop.3comp}.
\end{proof}
\subsection{Pre-processing}\label{sec.preprocess}
In this subsection we introduce the \textbf{pre-processing} step, discussed in Section \ref{sec.reduce5} (d). We start with the notion of \emph{ov-proper}, which is the extension of proper in Definition \ref{def.toy_proper}.
\begin{definition}\label{def.ov_proper} Let $\Mb$ be a molecule, we say it is \textbf{ov-proper}, if it is a forest, contains only deg 3 and 4 atoms, and the ov-distance (Definition \ref{def.ov_connect}) between any two deg 3 atoms is at least 3.
\end{definition}
Note that $\Mb_D$ may have fixed ends as the result of the cuttings in Definition \ref{def.layer_cutting} (and further cuttings of 2-connections, see \textbf{Stage 1} in Section \ref{sec.finish}), and might not be ov-proper. The idea of pre-processing is to turn $\Mb_D$ into ov-proper again by cutting as free a suitable set $S$. Namely, given a subset of atoms in $Z\cup Y\subseteq\Mb_D$ (which in practice will include the deg 3 atoms and weakly degenerate atoms (Definition \ref{def.weadeg})), we define a function \textbf{SELECT} (Definition \ref{def.func_select}) which outputs a set $S\supseteq Z\cup Y$, such that cutting $S$ as free turns $\Mb_D$ into ov-proper, see Proposition \ref{prop.func_select} (1). Moreover, since the \emph{UD connections} are important in the algorithm \textbf{MAINUD} below, we also prove a result on how cutting $S$ as free may affect the UD connections, see Proposition \ref{prop.func_select} (\ref{it.func_select_new2}).

The idea of the \textbf{SELECT} function is a simple greedy algorithm. Start with the set $Z\cup Y$ and its ov-components. If cutting $Z\cup Y$ as free would create (deg 2 atom or) two deg 3 atoms with ov-distance $\leq 2$, then there must exist two ov-components $A,B\subseteq Z\cup Y$ with ov-distance $\leq 4$. We then add to $Z\cup Y$ the intermediate atoms in the shortest path of ov-segments between $A$ and $B$ to bring them together, and repeat this process until it is no longer possible. In the end we get $S=\textbf{SELECT}(\Mb_D,Z,Y)$ such that cutting $S$ as free makes $\Mb_D$ ov-proper; see Definition \ref{def.func_select} and Proposition \ref{prop.func_select} (1). As for Proposition \ref{prop.func_select} (2), the proof relies on Lemma \ref{lem.3_comp_property} below, which is in the same spirit as the proof of uniqueness of lowest atom in Part 3 of the proof of Proposition \ref{prop.toy} in Section \ref{sec.toy}.
\begin{lemma}
\label{lem.3_comp_property} Let $\Mb_{UD}$ be a UD molecule such that $\Mb_D=\varnothing$, and $\Mb_U$ does not have C-top fixed end. Suppose $A\subseteq \Mb_U$ is ov-connected within $\Mb_U$ (cf. Definition \ref{def.ov_connect}) such that each atom in $A$ has deg 3. Then $A$ contains a \emph{unique} lowest atom $\nf$, such that (i) each atom $\mf\in A$ is obtained from $\nf$ by iteratively taking ov-parents in $A$, and (ii) each ov-child of $\mf\in A$ is either one of the terms preceding $\mf$ in the ov-parent sequence defined in (i), or an O-atom between two consecutive terms preceding $\mf$, or an ov-child of $\nf$.

Similarly, if we assume $\Mb_D\neq \varnothing$, $A$ is ov-connected within $\Mb_U$, each atom of $A$ is C-atom, and each deg 4 atom of $A$ has at most one ov-child in $A$, then conclusion (i) is still true.
\end{lemma}
\begin{proof} Choose a lowest atom $\nf\in A$. We will prove (i), which also establishes the uniqueness of $\nf$, and prove (ii) in the process. For any $\mf\in A$, since $A$ is ov-connected, there exists a sequence $\pf_j\in A\,(0\leq j\leq r)$ such that $(\pf_0,\pf_r)=(\nf,\mf)$ and $\pf_{j+1}$ is ov-adjacent to $\pf_j$; we may assume this sequence is shortest. Since $\nf$ is lowest, $\pf_1$ must be ov-parent of $\pf_0=\nf$; since $\pf_1$ has deg 3 (thus it is O-atom or has one bottom fixed end, cf. Remark \ref{rem.reg}, we know that $\pf_2$ must be ov-parent of $\pf_1$ (otherwise $\pf_2$ is ov-child of $\pf_1$ along the same maximal ov-segment containing $\pf_1$ and $\pf_0$, which means $\pf_2$ is ov-adjacent to $\pf_0$ and we can skip $\pf_1$ to get a shorter sequence), and for the same reason $\pf_3$ must be ov-parent of $\pf_2$, and so on. This proves (i). Then (ii) follows from the same proof (any ov-child of $\pf_j$ either belongs to the ov-segment between $\pf_j$ and $\pf_{j-1}$, or is an ov-child of $\pf_{j-1}$ so we can keep iterating). The proof of (i) in the case $\Mb_D\neq\varnothing$ is similar, where we note that $\pf_j$ has at most one ov-child in $A$ in the case when $\pf_j$ has deg 4, so $\pf_{j+1}$ still cannot be ov-child of $\pf_j$ if $\pf_{j+1}\in A$ and is C-atom.
\end{proof}
\begin{definition}[The function \textbf{SELECT}]
\label{def.func_select} Let $\Mb$ be any UD molecule, such that $\Mb_D$ has no deg 2 atoms, and has no C-bottom fixed end. Let $Z$ be the set of deg 3 atoms in $\Mb_D$, and $Y$ be a subset of the set of deg 4 atoms in $\Mb_D$. Given $(\Mb_D,Z,Y)$, we define a set $S:=\textbf{SELECT}(\Mb_D,Z,Y)$ as follows, see {\color{blue}Figure \ref{fig.select}}:
\begin{figure}[h!]
    \centering
    \includegraphics[width=0.55\linewidth]{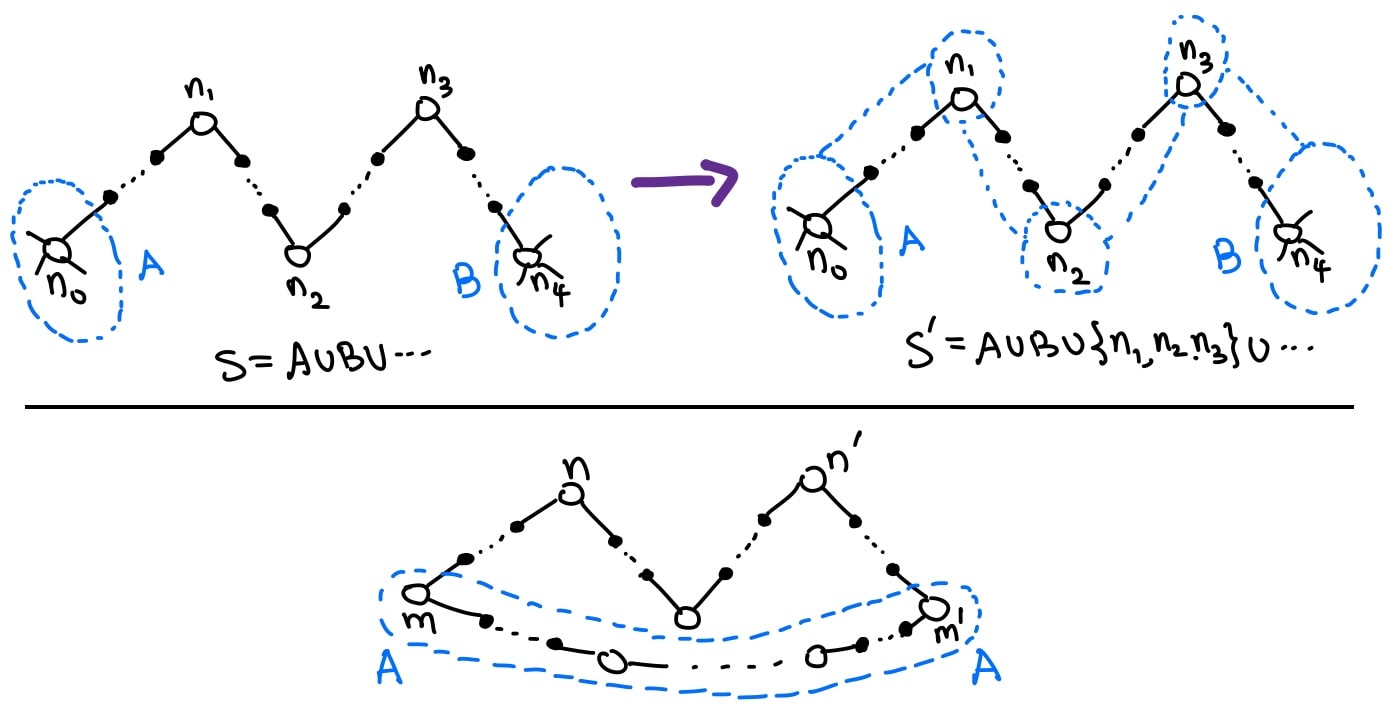}
    \caption{Top: the process of constructing the \textbf{SELECT} function in Definition \ref{def.func_select} (for two ov-components $A$ and $B$ with ov-distance $\leq 4$, add  the intermediate atoms in the ov-segments connecting them, until this cannot be done anymore). Bottom: proof of Proposition \ref{prop.func_select} (\ref{it.func_select_new1}); if $\mf$ and $\mf'$ belong to the same component $A$, then we get a cycle from $\mf\to\nf\to\cdots\to\nf'\to\mf'\to\cdots\to\mf$ (the path from $\mf'$ to $\mf$ belongs to $A$).}
    \label{fig.select}
\end{figure}
\begin{enumerate}
\item\label{it.func_select_1} Consider all the ov-components of $Z\cup Y$ in $\Mb_D$ (this means we only consider the O-atoms in $\Mb_D$ when defining ov-connectedness in Definition \ref{def.ov_connect}), which is a finite collection of disjoint subsets of $\Mb_D$.
\item\label{it.func_select_2} If any two of these subsets, say $A$ and $B$, have the shortest ov-distance (within $\Mb_D$) which is at most 4, then we choose the shortest path between an atom in $A$ and an atom in $B$, which can be divided into at most 4 ov-segments. Write this path as $(\nf_0,\cdots,\nf_q)$, where $q\leq 4$ and $\nf_0\in A$, $\nf_q\in B$, and each $(\nf_j,\nf_{j+1})$ is an ov-segment. Then replace the two sets $A$ and $B$ by a single set $A\cup B\cup\{\nf_1,\cdots,\nf_{q-1}\}$.
\item\label{it.func_select_4} Repeat (\ref{it.func_select_2}) until this can no longer be done (i.e. when only one set is left, or any two sets have of distance at least 5). Then output $S:=\textbf{SELECT}(\Mb_D,Z,Y)$ as the union of all the current sets.
\end{enumerate}
\end{definition}
\begin{proposition}
\label{prop.func_select} Let the function \textbf{SELECT} be as in Definition \ref{def.func_select}. Then we have the followings:
\begin{enumerate}
\item\label{it.func_select_new1} After cutting $S:=\textbf{SELECT}(\Mb_D,Z,Y)$ from $\Mb$ as free, then the set $\Mb_D\backslash S$ in the resulting molecule $\Mb\backslash S$ is ov-proper (Definition \ref{def.ov_proper}). Note that $\Mb_D\backslash S$ may have both C-top and C-bottom fixed ends, however this is compensated by the good property of being ov-proper.

\item\label{it.func_select_new2} Consider all the maximal ov-segments in $\Mb$ that contains at least one atom in $S$ and one atom in $\Mb_U$. Then the number of these maximal ov-segments does not exceed $10(\#_{\mathrm{ovcp}(Z)}+|Y|)$, where $\#_{\mathrm{ovcp}(Z)}$ is the number of ov-components of $Z$ (i.e. the set of all deg 3 atoms in $\Mb_D$), where the ov-components are determined by considering $\Mb_D$ only.
\end{enumerate}
\end{proposition}
\begin{proof} \textbf{Proof of (\ref{it.func_select_new1}).} At the end of the construction in Definition \ref{def.func_select} (\ref{it.func_select_4}), we get a number of disjoint subsets, such that each one is ov-connected and the ov-distance between each two of them is at least 5, and $S$ is the union of all these subsets. Note also that the sub-molecule $\Mb_D\backslash S$ of the molecule formed by cutting $S$ as free from $\Mb$, is equal to the molecule formed by cutting $S$ as free from $\Mb_D$.

Now suppose after cutting $S$ as free, there exist two deg 3 atoms $\nf,\nf'\in\Mb_D\backslash S$ of ov-distance at most 2; the case of a deg 2 atom in $\Mb_D\backslash S$ can be treated similarly. Since $\nf,\nf'\not\in S$, they must have deg 4 initially in $\Mb_D$. As they become deg 3 after cutting $S$ as free, we know that for each of $\nf$ and $\nf'$, it must be ov-adjacent to an atom in $S$. Assume $\nf$ is ov-adjacent to $\mf\in S$, and $\nf'$ is ov-adjacent to $\mf'\in S$. If $\mf$ and $\mf'$ belong to two different subsets $A$ and $B$ defined above, then the ov-distance between $\mf$ and $\mf'$ (and hence the ov-distance between $A$ and $B$) is at most $1+2+1=4$, contradicting the assumption about these subsets.

Therefore, we know that $\mf$ and $\mf'$ belong to the same set $A$ defined above. But then there exists a path $\mu_1:\mf\to\nf\to\nf'\to\mf'$ (composed of ov-segments), and another path $\mu_2$ between $\mf$ and $\mf'$, which is composed of ov-segments connecting atoms in $S$ (because $\mf$ and $\mf'$ belong to the same ov-component $A$ of $S$). The path $\mu_1$ cannot be expressed in this way (because the atoms $\nf$, $\nf'$ separating the different ov-segments, as well as the possible atom ov-adjacent to both $\nf$ and $\nf'$, all belong to $\Mb\backslash S$), unless this whole path belongs to a single ov-segment connecting $\mf$ and $\mf'$ (meaning the atoms $\nf$, $\nf'$ etc., all belong to this single ov-segment connecting $\mf$ and $\mf'$), but in this latter case the path $\nf\to\nf'$ would no longer exist after cutting $S$ as free (see Definition \ref{def.cutting}). Consequently, we know that the paths $\mu_1$ and $\mu_2$ \emph{cannot coincide}, but then they will form a cycle, which contradicts the fact that $\Mb_D$ is a forest.

\textbf{Proof of (\ref{it.func_select_new2}).} In each step in Definition \ref{def.func_select} (\ref{it.func_select_2}) we add at most 3 atoms to the union of the current subsets. Moreover, the number of such steps that can be carried out cannot exceed the initial number of subsets (as each step reduces the number of subsets by at least 1), which is the number of ov-components of $Z\cup Y$, and is at most $\#_{\mathrm{ovcp}(Z)}+|Y|$. This means that the number of newly added atoms in the construction of $S$ (starting from $Z\cup Y$) is at most $3(\#_{\mathrm{ovcp}(Z)}+|Y|)$.

Therefore, we only need to consider the required maximal ov-segments that contains an atom in $Z$. In each such maximal ov-segment, we have an atom $\mf$ belonging to an ov-component $W$ of $Z$, which has an ov-parent in $\Mb_U$. By (the dual version of) Lemma \ref{lem.3_comp_property} (ii), this ov-parent of $\mf$, which does not belong to $\Mb_D$, must also be an ov-parent of the unique highest atom $\nf$ of $W$; therefore this maximal ov-segment must contain $\nf$. As each ov-component $W$ contains a unique highest atom $\nf$, and each atom belongs to at most two maximal ov-segments, we conclude that the number of such maximal ov-segments does not exceed $2\cdot\#_{\mathrm{ovcp}(Z)}$ (as $\#_{\mathrm{ovcp}(Z)}$ is the number of ov-components of $Z$ in $\Mb_D$), as desired.
\end{proof}
\subsection{The \textbf{MAINUD} algorithm}\label{sec.mainud}
In this section we introduce the last ingredient, which is also the central algorithm in our proof, namely the \textbf{MAINUD} algorithm corresponding to toy model I plus (see \textbf{Stage 6} in Section \ref{sec.finish}). In Definition \ref{def.alg_maincr} below, we present the full version of this algorithm. Compared to the toy version (Definition \ref{def.toy1+_alg}), the main difference (apart from standard translations into ov-segment language) is that, the set $X$ in Definition \ref{def.toy1+_alg} is replaced by the set $X_1$ (see Definition \ref{def.alg_maincr} below). The reason has to do with the main dichotomy between toy models I plus and II in Section \ref{sec.toy} (basically $X_1$ is the set of deg 3 atoms after cutting $\Mb_D$ as free, which plays the same role as the $X$ in Section \ref{sec.toy}); see the discussion before \textbf{Stages 1--6} in Section \ref{sec.finish}. Apart from this, Definition \ref{def.alg_maincr}, Proposition \ref{prop.alg_maincr} and its proof are essentially the same as their toy versions (Definition \ref{def.toy1+_alg}, Proposition \ref{prop.toy1+_alg}) in Section \ref{sec.toy}.
\begin{definition}[The algorithm \textbf{MAINUD}] \label{def.alg_maincr}
Define the following cutting algorithm \textbf{MAINUD}. It takes as input any UD molecule $\Mb_{UD}$ that satisfies the following assumptions:
\begin{enumerate}[{(i)}]
\item $\Mb_U$ contains no C-top fixed end, and $\overline{\Mb}_{\mathrm{2conn}}^2=\varnothing$, where $\overline{\Mb}_{\mathrm{2conn}}^2$ is defined as in Proposition \ref{prop.comb_est_case6}.
\item $\Mb_D$ is ov-proper (Definition \ref{def.ov_proper}), and does not have any full component that is not connected to $\Mb_U$ by a bond.
\end{enumerate}

 For any such $\Mb$, we define the sets $X_0$ and $X_1$ of C-atoms in $\Mb_U$ as follows (see {\color{blue}Figure \ref{fig.x0_x1}}):
 \begin{itemize}
 \item Define $\nf\in X_0^0$ if $\nf$ has deg 3 and has one ov-child in $\Mb_D$.
 \item Define $\nf\in X_0^q$, if $\nf\in X_0^{q-1}$, or $\nf$ has deg 3 and has one ov-child in $X_0^{q-1}\cup\Mb_D$, or if $\nf$ has two ov-children in $X_0^{q-1}\cup\Mb_D$ along two ov-segments in two different particle lines. Then define $X_0=\cup_{q\geq 0}X_0^q$.
 \item With $X_0$ fixed, define $\nf\in X_1$ if and only if $\nf\not\in X_0$, and $\nf$ either has deg 3 or has one ov-child in $X_0\cup\Mb_D$.
 \end{itemize}
 \begin{figure}[h!]
    \centering
    \includegraphics[width=0.65\linewidth]{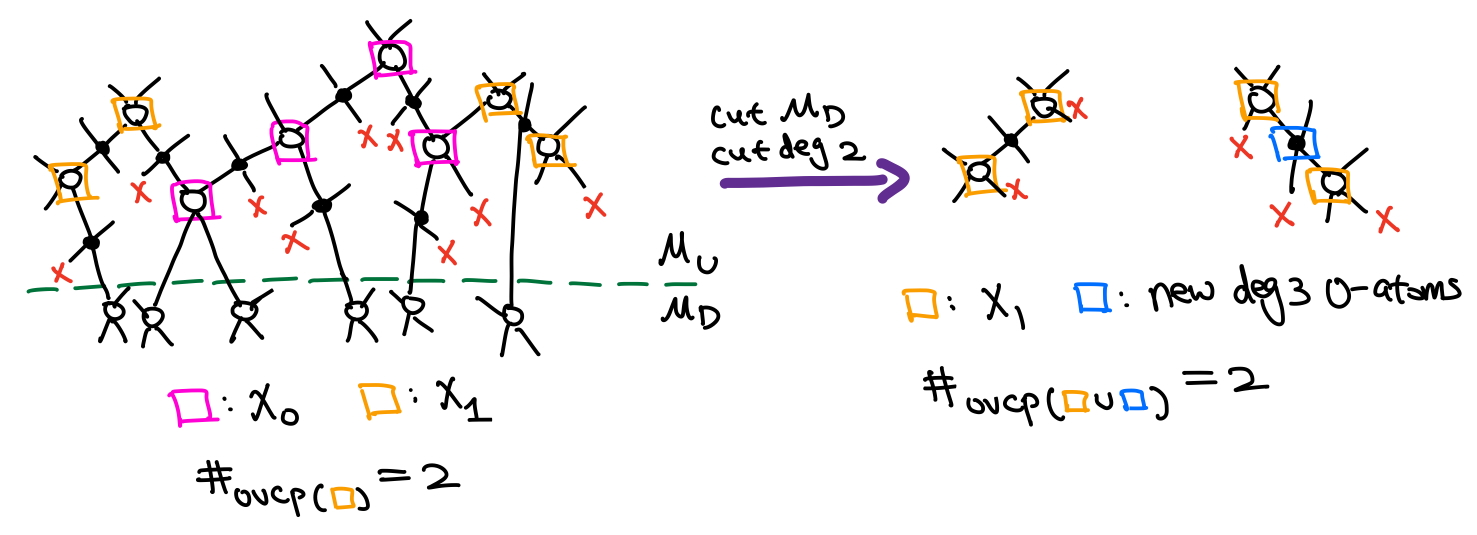}
    \caption{The sets $X_0$ and $X_1$ in Definition \ref{def.alg_maincr}, constructed by definition. The main point (which will be used in \textbf{Stage 5} in Section \ref{sec.finish}) is that, after cutting $\Mb_D$ as free, (i) all atoms in $X_0$ (pink) will be subsequently cut as deg 2, (ii) $X_1$ (yellow) is precisely the set of deg 3 C-atoms, (iii) the number of ov-components does not decrease if we add the new deg 3 O-atoms (blue) to $X_1$.}
    \label{fig.x0_x1}
\end{figure}
 
Note these sets are fixed throughout the cutting sequence below. Define the following cutting sequence:
\begin{enumerate}
\item\label{it.alg_maincr_2} If $\Mb_U$ contains any deg $2$ atom $\nf$, then cut it as free, and repeat until there is no deg $2$ atom left. 
\item\label{it.alg_maincr_3} Consider the set of all (remaining) atoms that either belong to $X_0\cup X_1$, or has deg 3; choose $\nf$ to be a lowest atom in this set (if this set is empty then choose $\nf$ to be any lowest atom). Let $S_\nf$ be the set of descendants of $\nf$ in $\Mb_U$, and fix $\nf$ and $S_\nf$ until the end of (\ref{it.alg_maincr_7}).
\item\label{it.alg_maincr_4} Starting from $\nf$, choose a highest atom $\pf$ in $S_\nf$; if $\pf$ is ov-adjacent to a deg 3 atom $\qf\in\Mb_D$, then cut $\{\pf,\qf\}$ as free, and cut $\qf$ as free from $\{\pf,\qf\}$ is $\pf$ has deg 4. Otherwise cut $\pf$ as free.
\item\label{it.alg_maincr_5} If any two deg 3 atoms $(\rf,\rf')$ in $\Mb_D$ becomes ov-adjacent, then choose $\{\rf,\rf'\}$ such that this ov-segment contains the fewest bond, and cut $\{\rf,\rf'\}$ as free. Repeat until no such instances exist.
\item\label{it.alg_maincr_6} If the ov-distance of any two deg 3 atoms $(\rf,\rf')$ in $\Mb_D$ becomes 2 (within $\Mb_D$), say $\rf$ and $\rf'$ are both ov-adjacent to some $\rf''\in\Mb_D$, then cut $\{\rf,\rf',\rf''\}$ as free and then cut $\rf$ as free from $\{\rf,\rf',\rf''\}$. Go to (\ref{it.alg_maincr_5}). Repeat until $\Mb_D$ becomes ov-proper.
\item\label{it.alg_maincr_7} When $\Mb_D$ becomes proper, go to (\ref{it.alg_maincr_4}). Repeat until all atoms in $S_\nf$ have been cut.
\item\label{it.alg_maincr_8} When all atoms in $S_\nf$ have been cut, go to (\ref{it.alg_maincr_2}). Repeat until all atoms in $\Mb_U$ have been cut.
\item\label{it.alg_maincr_9} When all atoms in $\Mb_U$ have been cut, choose any deg 3 atom in $\Mb_D$ and cut it as free. If the ov-distance between any two deg 3 atoms becomes at most 2, go to (\ref{it.alg_maincr_5}) and repeat (\ref{it.alg_maincr_5})--(\ref{it.alg_maincr_6}) until $\Mb_D$ becomes proper again. Then choose and cut the next deg 3 atom in $\Mb_D$, and so on.
\end{enumerate}
\end{definition}
\begin{proposition}
\label{prop.alg_maincr} Let $\Mb$ be any UD molecule that satisfies (i)--(ii) in Definition \ref{def.alg_maincr}. Consider the number of UD connections $e$ that are either (i) connected to some other UD connection via $\Mb_D$, or (ii) connected to a component of $\Mb_D$ that contains deg 3 atoms. We call such UD connections \textbf{good} and let the number of them be $\#_{\mathrm{UD}}$.

Let $\#_{\mathrm{fucp}(\Mb_U)}$ be the number of full components of $\Mb_U$, and let $\#_{\mathrm{ovcp}(X_1)}$ be the number of ov-components of the set $X_1$ defined in Definition \ref{def.alg_maincr}. Then, for the algorithm \textbf{MAINUD}, we have $\#_{\{44\}}=0$, and each \{33B\} molecule must have both its atoms in $\Mb_D$, and that
\begin{equation}\label{eq.alg_maincr}
\#_{\{33A\}}+\#_{\{33B\}}\geq\frac{1}{5}\cdot(\#_{\mathrm{UD}}-\#_{\mathrm{fucp}(\Mb_U)}-\#_{\mathrm{ovcp}(X_1)}),\quad\mathrm{and}\quad\#_{\{4\}}\leq \#_{\mathrm{fucp}(\Mb_U)}+\#_{\mathrm{ovcp}(X_1)}.
\end{equation}
\end{proposition}
\begin{proof} We divide the proof into 4 parts, which exactly correspond to the 4 parts of the proof in toy model I plus (Proposition \ref{prop.toy1+_alg}); the presentation of Part 4 is slightly different, but essentially it is the same counting argument (involving the same quantity $\#_{\mathrm{deg3}}-\#_{\mathrm{comp}}$).

\textbf{Proof part 1.} We prove the (MONO) property as in Proposition \ref{prop.alg_up}: at any time in the cutting sequence, there is no C-top fixed end in $(\Mb_U\backslash S_\nf)\cup\{\nf\}$ and no C-bottom fixed end in $S_\nf\backslash\{\nf\}$. In fact, when we choose $\nf$ as in Definition \ref{def.alg_maincr} (\ref{it.alg_maincr_3}), every atom $\pf\in S_\nf\backslash\{\nf\}$ must have deg 4. In the process of Definition \ref{def.alg_maincr} (\ref{it.alg_maincr_4})--(\ref{it.alg_maincr_7}), if we only consider the cutting within $\Mb_U$, then the desired property follows from the same proof as in Proposition \ref{prop.alg_up}. Now suppose an arbitrary set of atoms in $\Mb_D$ are also cut in the process, then this will only create C-bottom fixed ends, and only at C-atoms $\tf$ that are ov-adjacent to some atom in $\Mb_D$. By definition we must have $\tf\in X_0\cup X_1$, and thus $\tf\not\in S_\nf\backslash\{\nf\}$ by our choice, so it does not affect the desired (MONO) property. Similarly, for any \{33\} molecule $\{\pf,\qf\}$ cut in (\ref{it.alg_maincr_4}), the $\pf$ must either be O-atom or $\pf=\nf$, so this \{33\} molecule must be \{33A\} (as $\pf$ has not C-top fixed end). 

\textbf{Proof part 2.} We prove that each atom in $\Mb_D$ belongs to a \{3\} or \{33\} molecule (which can be either \{33A\} or \{33B\}). In fact, note first that initially $\Mb_D$ is ov-proper (thus has no deg 2 atom), and does not have full components not connected to $\Mb_U$ by a bond (thus will not have any full component after cutting $\Mb_U$ as free), so the cuttings in Definition \ref{def.alg_maincr} (\ref{it.alg_maincr_9}) will not create any \{4\} molecule. 

Next, note that $\Mb_D$ initially is a forest, and will remain a forest in the whole cutting process (cutting an atom will not create any cycle, but may cause a tree to become a disconnected forest). Now consider each cutting step in Definition \ref{def.alg_maincr} (\ref{it.alg_maincr_4}) involving $\pf\in\Mb_U$, at which time $\Mb_D$ must be proper by the construction in Definition \ref{def.alg_maincr} (\ref{it.alg_maincr_5})--(\ref{it.alg_maincr_7}). This $\pf$ has \emph{at most one} ov-segment connecting to $\Mb_D$, because $\overline{\Mb}_{\mathrm{2conn}}^2=\varnothing$ by Definition \ref{def.alg_maincr} (i). If no such ov-segment exists, then cutting $\pf$ does not affect $\Mb_D$. If \emph{one} such ov-segment exists, then the effect of the cutting step in Definition \ref{def.alg_maincr} (\ref{it.alg_maincr_4}) \emph{on $\Mb_D$}, has the following two possibilities:
\begin{enumerate}[{(a)}]
\item\label{it.stepa1} If we cut $\pf$ as free (which means all atoms in $\Mb_D$ ov-adjacent to $\pf$ have deg 4), then this step breaks $\Mb_D$ into connected components $X_\qf$ which are in bijection with atoms $\qf\in\Mb_D$ ov-adjacent to $\pf$, and this cutting creates only one fixed end at $\qf$ when viewed within $X_\qf$, turning it from deg 4 to deg 3.
\item\label{it.stepb1} If we cut $\{\pf,\qf\}$ as free (which means $\qf\in\Mb_D$ is ov-adjacent to $\pf$ and has deg 3), then when viewed in $\Mb_D$, this step equivalent to cutting $\qf$ (which has deg 3) as free in $\Mb_D$.
\end{enumerate} Here \ref{it.stepa1} follows from the same proof in Lemma \ref{lem.cutconnected} (i.e. $\qf$ in bijection with $X_\qf$, in which only one fixed end is created at $\qf$) as $\Mb_D$ is a forest. In \ref{it.stepb1}, the effect of cutting $\{\pf,\qf\}$ as free to $\Mb_D$ is equivalent to cutting $\{\qf\}$ as free, because $\pf$ has no other ov-segment connecting to other atoms in $\Mb_D$ (apart from those on the ov-segment containing $\qf$).

Now, regarding the effects described in \ref{it.stepa1}--\ref{it.stepb1}, we have the following lemma, which is extension of Lemma \ref{lem.toy_1} in the proof of Proposition \ref{prop.toy}:
\begin{lemma}\label{lem.mainudlem} Suppose $\Mb_D$ is proper at some point of the algorithm, and becomes non-proper after either \ref{it.stepa1} or \ref{it.stepb1}. Then, $\Mb_D$ will once more become proper after applying Definition \ref{def.alg_maincr} (\ref{it.alg_maincr_5})--(\ref{it.alg_maincr_6}), via cuttings of \{3\} and \{33\} molecules only, and \emph{without seeing any deg 2 atom}.
\end{lemma}
\begin{proof} Note that $\Mb_D$ is always a forest, and our proof will crucially rely on Lemma \ref{lem.cutconnected}. Consider first \ref{it.stepa1}: when viewed within each resulting component $X_\qf$, its effect is simply \emph{switching one free end at a deg 4 atom into a fixed end, turning it into deg 3}. We call the latter step $(\ast)$. As for \ref{it.stepb1}, right before performing it at a deg 3 atom $\qf$, any atom in $\Mb_D$ that is ov-adjacent to $\qf$ must have deg 4 (because $\Mb_D$ is ov-proper). By Lemma \ref{lem.cutconnected}, $\Mb_D$ will break into finitely many components after performing \ref{it.stepb1}, and the effect of \ref{it.stepb1} \emph{when viewed within each component} is again step $(\ast)$. Therefore we only need to prove Lemma \ref{lem.mainudlem} for step $(\ast)$.

Assume now $\Mb_D$ becomes non-proper, with two ov-adjacent deg 3 atoms (say $\rf$ and $\rf'$), after step $(\ast)$. We may assume $\rf'$ has deg 3 before $(\ast)$, while $\rf$ has deg 4 before $(\ast)$ and is turned into deg 3 after $(\ast)$. Then, \emph{before $(\ast)$}, any atom in $\Mb_D$ ov-adjacent to $\rf$ or $\rf'$ \emph{must have deg 4}, because otherwise we either have two ov-adjacent deg 3 atoms before $(\ast)$, or have two deg 3 atoms ov-adjacent to the same deg 4 atom before \ref{it.stepa}, contradicting properness. Therefore, by Lemma \ref{lem.cutconnected}, we see that cutting the \{33\} molecule breaks $\Mb_D$ into finitely many components, and the effect of this cutting \emph{when viewed in each component} is again equivalent to step $(\ast)$. We then induct on the size of $\Mb_D$ and apply induction hypothesis to each of these components, to conclude the proof of Lemma \ref{lem.mainudlem}.

Now suppose $\Mb_D$ becomes non-proper, without two ov-adjacent deg 3 atoms, after $(\ast)$. In this case there must be two deg 3 atoms, say $\rf$ and $\rf''$, that are both ov-adjacent to a deg 4 atom $\rf'$. We may assume $\rf'$ and $\rf''$ respectively have deg 4 and 3 before $(\ast)$, and $\rf$ has deg 4 before $(\ast)$ and is turned into deg 3 after $(\ast)$. Then,  \emph{before $(\ast)$}, any atom ov-adjacent to $\rf'$ or $\rf''$ \emph{must have deg 4} for the same reason as above thanks to properness. Moreover, before $(\ast)$, any atom adjacent to $\rf$ must also have deg 4, because otherwise we would get two ov-adjacent deg 3 atoms after $(\ast)$, contradicting our assumption. Therefore, by Lemma \ref{lem.cutconnected}, we see that cutting $\{\rf,\rf',\rf''\}$ breaks $\Mb_D$ into finitely many components, and the effect of this cutting \emph{when viewed in each component} is again equivalent to $(\ast)$. We then conclude the proof by the same induction as above. This completes the proof of Lemma \ref{lem.mainudlem}.
\end{proof}
Now, with Lemma \ref{lem.mainudlem}, we can finish Part 2 of the proof, i.e. each atom in $\Mb_D$ is contained in a \{3\} or \{33\} molecule. In fact, those cut in Definition \ref{def.alg_maincr} (\ref{it.alg_maincr_4}) clearly belong to \{33\} molecule. Those cut in Definition \ref{def.alg_maincr} (\ref{it.alg_maincr_5})--(\ref{it.alg_maincr_6}) also belong to \{3\} or \{33\} molecules by Lemma \ref{lem.mainudlem}. Finally, those cut in Definition \ref{def.alg_maincr} (\ref{it.alg_maincr_9}) either belong to \{3\} molecule, or are cut in the process of applying Definition \ref{def.alg_maincr} (\ref{it.alg_maincr_5})--(\ref{it.alg_maincr_6}) and thus belong to \{3\} or \{33\} molecules by Lemma \ref{lem.mainudlem}.

\textbf{Proof part 3.} We prove that the number of deg 4 atoms $\nf\in\Mb_U$ cut in Definition \ref{def.alg_maincr} (\ref{it.alg_maincr_4}) is at most  $\#_{\mathrm{fucp}(\Mb_U)}+\#_{\mathrm{ovcp}(X_1)}$. In fact, such $\nf$ must be the atom chosen in Definition \ref{def.alg_maincr} (\ref{it.alg_maincr_3}) (the $\pf\in S_\nf\backslash\{\nf\}$ in (\ref{it.alg_maincr_4}) will not be cut as deg 4 because it has an ov-parent in $S_\nf$ which must be cut before $\pf$). Moreover, there are two possibilities when this $\nf$ can has deg 4: either the component of $\Mb_U$ containing $\nf$ contains only deg 4 atoms, or $\nf$ is the lowest remaining atom in $X_0\cup X_1$ and has deg 4. The former case leads to at most $\#_{\mathrm{fucp}(\Mb_U)}$ choices for $\nf$.

In the latter case, $\nf$ must have deg 4 and must not have any ov-child in $X_0\cup X_1$ (otherwise this ov-child must have been cut before $\nf$ as $\nf$ is the lowest remaining atom in $X_0\cup X_1$, so $\nf$ would not have deg 4), which means that $\nf\not\in X_0$ by definition of $X_0$ (and that $\overline{\Mb}_{\mathrm{2conn}}^2=\varnothing$). Therefore $\nf\in X_1$, and it must belong to one of the initial ov-components of $X_1$, say $A$. Since $\nf$ has no ov-child in $A$, by Lemma \ref{lem.3_comp_property} (i), it must be the unique lowest atom in $A$ (here $A$ satisfies the assumption of Lemma \ref{lem.3_comp_property}, because each deg 4 atom in $X_1$ must have an ov-child in $X_0\cup\Mb_D$, thus has at most one ov-child in $X_1$), which leads to at most $\#_{\mathrm{ovcp}(X_1)}$ choices for $\nf$ due to the same number of choices of $A$. Putting together, this finishes Part 3 of the proof.

\textbf{Proof part 4.} We now complete the proof of (\ref{eq.alg_maincr}). The bound for $\#_{\{4\}}$ follows directly from Part 3, as no deg 4 atom can be cut in Definition \ref{def.alg_maincr} (\ref{it.alg_maincr_5})--(\ref{it.alg_maincr_6}), or Definition \ref{def.alg_maincr} (\ref{it.alg_maincr_9}) (as shown in Part 2). To prove the lower bound for $\#_{\{33\}}$ (which is $\#_{\{33A\}}+\#_{\{33B\}}$), note that with at most $\mu:=\#_{\mathrm{fucp}(\Mb_U)}+\#_{\mathrm{ovcp}(X_2)}$ exceptions, each \{3\} molecule in $\Mb_D$ created in Definition \ref{def.alg_maincr} (\ref{it.alg_maincr_5})--(\ref{it.alg_maincr_6}) is always paired with a \{33\} molecule (the exceptions correspond to those atoms $\pf$ in Definition \ref{def.alg_maincr} (\ref{it.alg_maincr_4}) with deg 4, which may cause $\qf$ to be cut as a \{3\} molecule without a pairing \{33\} molecule), so the number of atoms in $\Mb_D$ cut in these steps is at most $3\cdot\#_{\{33\}}+\mu$.

Now consider the $\#_{\mathrm{UD}}$ good UD connections and the corresponding atoms $\rf\in \Mb_D$. If any such $\rf$ is not cut in Definition \ref{def.alg_maincr} (\ref{it.alg_maincr_5})--(\ref{it.alg_maincr_6}), then in the molecule $\Mb_D$ after Definition \ref{def.alg_maincr} (\ref{it.alg_maincr_8}), this $\rf$ must be a deg 3 atom (it cannot be deg 4 as it belong to a UD connection, and cannot be deg 2 because $\Mb_D$ is ov-proper by Part 2), and must belong to a component of $\Mb_D$ that contains at least two deg 3 atoms (by considering the pre-algorithm path connecting $\rf$ to a deg 3 atom or another UD connection, and the possible place where a bond of this path is broken; note that a broken bond creates a fixed end and thus deg 3 atom). We collect all these atoms $\rf$ (which are not cut in Definition \ref{def.alg_maincr} (\ref{it.alg_maincr_5})--(\ref{it.alg_maincr_6})), with the number of them being $\geq \#_{\mathrm{UD}}-(3\cdot\#_{\{33\}}+\mu)$ by the above proof, and consider all the components of $\Mb_D$ (after Definition \ref{def.alg_maincr} (\ref{it.alg_maincr_8})) containing them; let the number of deg 3 atoms in each of these components be $m_j\geq 2$, then after Definition \ref{def.alg_maincr} (\ref{it.alg_maincr_8}) we have that (note that each component in $\Mb_D$ will have at least one deg 3 atom after cutting $\Mb_U$, see Remark \ref{rem.full_cut}):
\begin{equation}\label{eq.alg_maincr_2}\#_{\mathrm{deg3}} - \#_{\mathrm{comp}}\geq\sum_j(m_j-1)
\geq\frac{1}{2}\sum_j m_j\geq\,(\#_{\mathrm{UD}}-3\cdot\#_{\{33\}}-\mu)/2,
\end{equation} where $\#_{\mathrm{deg3}}$ and $\#_{\mathrm{comp}}$ are the number of deg 3 atoms and number of components respectively. Then, consider the increment of the left hand side of (\ref{eq.alg_maincr_2}) in each cutting in Definition \ref{def.alg_maincr} (\ref{it.alg_maincr_9}). Using Definition \ref{def.cutting}, and noting that each such cut can only create one fixed end at one deg 4 atom (adding only one deg 3 atom) in each resulting (extra) component, we can verify that this is invariant after cutting one deg 3 atom, and decreases by 1 after each cutting in Definition \ref{def.alg_maincr} (\ref{it.alg_maincr_5})--(\ref{it.alg_maincr_6}). Since this value should be come $0$ after all atoms of $\Mb_D$ have been cut, we conclude that \[\#_{\{33\}}\geq \frac{1}{2}(\#_{\mathrm{UD}}-3\cdot\#_{\{33\}}-\mu),\] and therefore
\[\#_{\{33\}}\geq\frac{1}{5}\cdot(\#_{\mathrm{UD}}-\mu).\]This completes the proof.
\end{proof}
\subsection{Toy model I plus/II dichotomy, and finishing the proof}\label{sec.finish}
We now put everything together and finish the proof of Proposition \ref{prop.case5}. Now we may assume (i)--(iii) in Reduction \ref{red.2}, as well as Reduction \ref{red.3} (otherwise Proposition \ref{prop.case5} already follows from Proposition \ref{prop.comb_est_case6}); namely, we assume that
\begin{itemize}
\item Each of $\Mb_U$ and $\Mb_D$ is a forest and has no deg 2 atoms. No atom in $\Mb_D$ is a parent of any atom in $\Mb_U$, and $\Mb_U$ (resp. $\Mb_D$) has no C-top (resp. C-bottom) fixed end.
\item The number of full  components of $\Mb_U$ is $\#_{\mathrm{fucp}(\Mb_U)}\leq 10(C_{5}^*)^{-1}\cdot\rho'''$. Each component in $\Mb_D$ has either at least two bonds connecting to $\Mb_U$, or one bond connecting to $\Mb_U$ and one deg 3 atom.
\item The number of good UD connections (i.e. UD connections that are connected to either another UD connection or a deg 3 atom in $\Mb_D$, see Proposition \ref{prop.alg_maincr}) satisfy $\#_{\mathrm{UD}}\geq \rho'''/4$.
\item The number of components of $\Mb_D$ is $\#_{\mathrm{comp}(\Mb_D)}\leq (C_{3}^*)^{-1}\cdot\rho'''$, and the number of weakly degenerate atoms and pairs in $\Mb_D$ is less than $(C_{2}^*)^{-1}\cdot\rho'''$.
\item We have $|\overline{\Mb}_{\mathrm{2conn}}^2|\leq (C_{2}^*)^{-1}\cdot\rho'''$, with the set $\overline{\Mb}_{\mathrm{2conn}}^2$ defined in Proposition \ref{prop.comb_est_case6}.
\end{itemize}
\begin{proof}[Proof of Proposition \ref{prop.case5}]
Before proceeding with \textbf{Stages 1--6} below, we first briefly summarize our strategy, how the different toy models (I plus and II) in Section \ref{sec.toy} occur in the proof, and how the different algorithms (\textbf{3COMPUP} and \textbf{MAINUD}) are used. First note that, from the above assumptions, each of $\Mb_U$ and $\Mb_D$ is essentially a tree (with negligible full components), with negligible number of 2-connections in $\Mb_U$, and sufficiently many good UD connections (all relative to $\rho'''$; same below).

In \textbf{Stage 1}, we first cut as free all 2-connections (i.e. atoms in $\overline{\Mb}_{\mathrm{2conn}}^2$) in $\Mb_U$ and subsequently cut any deg 2 atoms. This allows to completely get rid of 2-connections in $\Mb_U$ at a loss of negligible \{4\} molecules. The quantity $\#_{\mathrm{UD}}$ (i.e. the number of \emph{good UD connections}) may decrease, but the decrease is again negligible, so we still have $\#_{\mathrm{UD}}\gtrsim\rho'''$ after this stage.

Now, let $\overline{Z}$ be the set of deg 3 atoms in $\Mb_D$ (in reality this needs slight modification to avoid components with only one deg 3 atom, which we omit here). Consider $\#_{\mathrm{ovcp}(\overline{Z})}$ which is the number of ov-components of $Z$; if $\#_{\mathrm{ovcp}(\overline{Z})}\gtrsim\rho'''$, then in \textbf{Stage 2}, we can apply \textbf{3COMPDN} (the dual of \textbf{3COMPUP}) to $\Mb_D$, which is already sufficient by Proposition \ref{prop.3comp}.

Therefore, we may assume $\#_{\mathrm{ovcp}(\overline{Z})}\ll\rho'''$. Now in \textbf{Stage 3}, we perform pre-processing to get rid of $\overline{Z}$ as well as the set $\overline{Y}$ of weakly degenerate atoms in $\Mb_D$, by cutting as free the set $S=\textbf{SELECT}(\Mb_D,\overline{Z},\overline{Y})$ defined in Definition \ref{def.func_select}. By assumption we have $\#_{\mathrm{ovcp}(\overline{Z})}\ll\rho'''$ and $|\overline{Y}|\ll\rho'''$, so by Proposition \ref{prop.func_select}, $\Mb_D$ becomes ov-proper after cutting $S$ as free, and the decrease of $\#_{\mathrm{UD}}$ is still negligible, so we still have $\#_{\mathrm{UD}}\gtrsim\rho'''$ after this stage.

At this point we can see the dichotomy that distinguishes the cases of toy model I plus (where we apply \textbf{MAINUD}) and toy model II (where we apply \textbf{3COMPUP}). Let $X_0$ and $X_1$ be defined as in Definition \ref{def.alg_maincr}; by definition, it is easy to see that if we cut $\Mb_D$ as free, then the C-atoms in $X_0$ are subsequently cut as deg 2, and the deg 3 C-atoms in $\Mb_U$ are \emph{exactly those in $X_1$} (the O-atoms in $\Mb_U$ will not matter in this discussion). As such, this $X_1$ plays the role of the set $X$ in toy model II (see the discussion below Definition \ref{def.toy2}), so our dichotomy will be based on the value of $\#_{\mathrm{ovcp}(X_1)}$, the number of ov-components of $X_1$ (which is the analog of the number of components of $X$ in Section \ref{sec.toy}).

Before proceeding, however, we need one more step: in order to apply \textbf{3COMPUP}, we need that each component of $\Mb_U$ after cutting $\Mb_D$ as free (and subsequently cutting deg 2 atoms) to contain \emph{at least two} deg 3 atoms. If a component $A$ violates this, it is easy to see that \emph{before cutting $\Mb_D$}, the set $A$ has only one bond connecting to $\Mb\backslash A$. This leads to \textbf{Stage 4}, where we get rid of these sets $A$ by cutting them as fixed, which allows us to apply \textbf{3COMPUP} when needed. It is easy to see that we still have $\#_{\mathrm{UD}}\gtrsim\rho'''$.

Now we can consider the dichotomy discussed above: if $\#_{\mathrm{ovcp}(X_1)}\gtrsim\rho'''$, then in \textbf{Stage 5} we cut $\Mb_D$ as free and apply \textbf{3COMPUP} to $\Mb_U$ (as in toy model II), which is sufficient by Proposition \ref{prop.3comp}. Finally, if $\#_{\mathrm{ovcp}(X_1)}\ll\rho'''$, then in \textbf{Stage 6} we apply \textbf{MAINUD} to $M_{UD}$ (as in toy model I plus). Note that $\Mb_D$ is ov-proper and has no weak degeneracy atom thanks to \textbf{Stage 3}, so any \{33\} molecule will be good. We then apply Proposition \ref{prop.alg_maincr}, and use that $\#_{\mathrm{UD}}\gtrsim\rho'''$ and $\#_{\mathrm{ovcp}(X_1)}\ll\rho'''$, to conclude the proof.

\textbf{Stage 1.} We cut each atom in $\overline{\Mb}_{\mathrm{2conn}}^2$ as free (see {\color{blue}Figure \ref{fig.stage1}}). If any atom in $\Mb_{UD}$ becomes deg 2, we also cut it as free, until $\Mb_{UD}$ contains no deg 2 atoms. Note that the contribution of these operations to $\#_{\{4\}}$ is at most $|\overline{\Mb}_{\mathrm{2conn}}^2|\leq (C_{2}^*)^{-1}\cdot\rho'''$.
\begin{figure}[h!]
    \centering
    \includegraphics[width=0.65\linewidth]{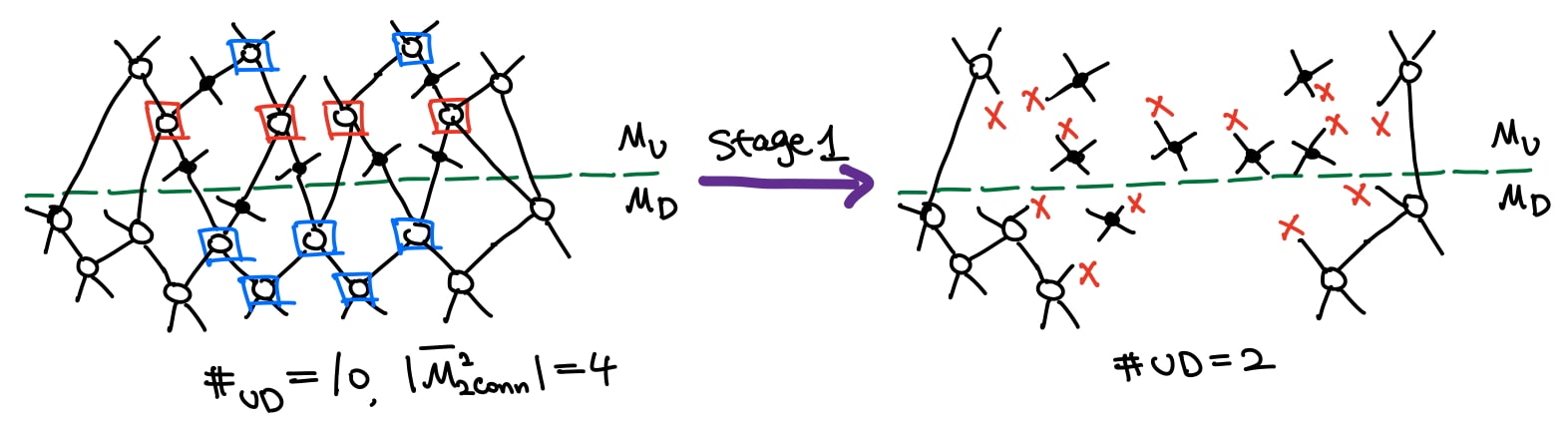}
    \caption{Cutting atoms in $\overline{\Mb}_{\mathrm{2conn}}^2$ (red) and subsequent deg 2 atoms (blue) in \textbf{Stage 1}. Note that cutting each red breaks at most 2 UD connections (cutting each blue breaks none), and any good UD connection that is not broken shall remain good.}
    \label{fig.stage1}
\end{figure}

By definition of $\overline{\Mb}_{\mathrm{2conn}}^2$, cutting any atom in it as free does not create any C-top fixed atom in $\Mb_U$ or any C-bottom fixed atom in $\Mb_D$, and neither does any subsequent operations cutting deg 2 atoms (because, when we subsequently cut any deg 2 atom in $\Mb_U$ (resp. $\Mb_D$), it must be O-atom or has two bottom (resp. top) fixed ends, so cutting it will not create any top (resp. bottom) fixed end). Moreover, consider any UD connection $e$ and the maximal ov-segment containing $e$; cutting each atom in $\overline{\Mb}_{\mathrm{2conn}}^2$ breaks at most two such ov-segments (where breaking means cutting an atom in this ov-segment as free), and any subsequent cutting of deg 2 atoms cannot break any such ov-segment (for the same reason as above).

If the ov-segment containing $e$ is not broken, then $e$ will remain a good UD connection if it is initially so. In fact, consider the component $X$ of $\Mb_D$ containing the $\Mb_D$ endpoint of $e$; if $X$ is not affected in the cutting (i.e. no atom in $X$ is cut and no fixed end is created in $X$) then $e$ will remain good as it is; if $X$ is affected then the new component containing the $\Mb_D$ endpoint of $e$ must contain a deg 3 atom.

Therefore, after \textbf{Stage 1}, we still have \begin{equation}\label{eq.comb_est_case71}\#_{\mathrm{UD}}\geq \rho'''/4-2|\overline{\Mb}_{\mathrm{2conn}}^2|\geq \rho'''/5.
\end{equation}
Moreover, after \textbf{Stage 1}, $\Mb_{UD}$ is still a canonical UD molecule, and now $\overline{\Mb}_{\mathrm{2conn}}^2=\varnothing$. The number of full components of $\Mb_U$ and $\Mb_D$ still satisfy that $\#_{\mathrm{fucp}(\Mb_U)}+\#_{\mathrm{fucp}(\Mb_D)}\leq 4(C_{3}^*)^{-1}\cdot \rho'''$.

\textbf{Stage 2.} Now consider all the components of $\Mb_D$ that contain at least two deg 3 atoms, and all the deg 3 atoms contained in one such component. Let the set of all such deg 3 atoms be $\overline{Z}_0$, and let $\overline{Y}$ be the union of all weakly degenerate atoms and atom pairs in $\Mb_D$ (see Proposition \ref{prop.comb_est_case4}). Let the number of ov-components of $\overline{Z}_0$ in $\Mb_D$ be $\#_{\mathrm{ovcp}(\overline{Z}_0)}$ as in Propositions \ref{prop.3comp} and \ref{prop.func_select}. 

Now, if $\#_{\mathrm{ovcp}(\overline{Z}_0)}\geq (10/\upsilon)^4(C_{2}^*)^{-1}\cdot \rho'''$, then we cut $\Mb_D$ as free. For each component of $\Mb_D$ containing at least two deg 3 atoms, we cut it into elementary molecules using \textbf{3COMPDN} (the dual of \textbf{3COMPUP} in Definition \ref{def.3comp_alg}); for other components we cut them into elementary molecules using \textbf{DOWN}. Then cut $\Mb_U$ into elementary molecules using \textbf{UP}.

In this case, by Proposition \ref{prop.3comp} we get $\#_{\{33B\}}=\#_{\{44\}}=0$, and
\begin{equation}\label{eq.comb_est_case72}\#_{\{33A\}}\geq (1/10)\cdot\#_{\mathrm{ovcp}(\overline{Z}_0)}\geq100\upsilon^{-4}\cdot (C_{2}^*)^{-1}\cdot\rho''',
\end{equation}
\begin{equation}\label{eq.comb_est_case73}
\#_{\{4\}}\leq \#_{\mathrm{fucp}(\Mb_U)}+\#_{\mathrm{fucp}(\Mb_D)}+(C_{2}^*)^{-1}\cdot\rho'''\leq 2\cdot(C_{2}^*)^{-1}\cdot \rho''',
\end{equation} which already proves (\ref{eq.comb_est_case51}) (recall that all \{33A\} molecules are treated as good).

\textbf{Stage 3.} Now suppose $\#_{\mathrm{ovcp}(\overline{Z}_0)}\leq (10/\upsilon)^4(C_{2}^*)^{-1}\cdot \rho'''$. Define $\overline{Z}$ to be the union of $\overline{Z_0}$, together with the deg 3 atoms in those components of $\Mb_D$ that contain only one deg 3 atom and at least one atom in $\overline{Y}$. Then for the set $\overline{Z}$ we have $\#_{\mathrm{3ovcp}(\overline{Z})}\leq \#_{\mathrm{ovcp}(\overline{Z}_0)}+|\overline{Y}|$. Note also $|\overline{Y}|\leq 10(C_{2}^*)^{-1}\cdot\rho'''$ by Reduction \ref{red.3}.

We then define $S:=\textbf{SELECT}(\widetilde{\Mb}_D,\overline{Z},\overline{Y})$ as in Definition \ref{def.func_select}, where $\widetilde{\Mb}_D$ contains only those components of $\Mb_D$ that contain an atom in $\overline{Z}\cup\overline{Y}$ (but this does not affect Definition \ref{def.func_select}). We then cut $S$ as free and subsequently cut it into elementary molecules using \textbf{DOWN}. If any atom in $\Mb_{U}$ becomes deg 2, we also cut it as free, until $\Mb_{U}$ contains no deg 2 atoms.

By Proposition \ref{prop.func_select} (\ref{it.func_select_new1}), we know that $\Mb_{UD}$ is still a UD molecule with $\overline{\Mb}_{\mathrm{2conn}}^2=\varnothing$, and $\Mb_D$ becomes ov-proper after cutting $S$ (those components of $\Mb_D$ that are not in $\widetilde{\Mb}_D$ will contain at most one deg 3 atom by definition of $\Mb_D$, and are thus automatically ov-proper). Now $\Mb_D$ may contain C-bottom fixed ends, but $\Mb_U$ still does not contain C-top fixed end. Moreover, by Proposition \ref{prop.func_select} (\ref{it.func_select_new2}), we know that the number of ov-segments corresponding to UD connections $e$, that are broken in the process of cutting $S$ as free (such an ov-segment must contain an atom in $S$; again subsequent operations cutting deg 2 atoms will not break any such ov-segments), is at most
\begin{equation}\label{eq.comb_est_case74}10\cdot(3\cdot\#_{\mathrm{ovcp}(\overline{Z})}+|\overline{Y}|)\leq C\cdot (C_{2}^*)^{-1}\cdot\rho''',
\end{equation} therefore we still have
\begin{equation}\label{eq.comb_est_case75}\#_{\mathrm{UD}}\geq \rho'''/5-C\cdot (C_{2}^*)^{-1}\cdot\rho'''\geq \rho'''/6,
\end{equation} after \textbf{Stage 3} (note that the good UD connections that are not broken will remain good, as proved in \textbf{Stage 1} above). The contribution to $\#_{\{4\}}$ of these cutting operations is also bounded by (\ref{eq.comb_est_case74}). Note that any component of $\Mb_D$ that is not
connected to $\Mb_{U}$ after cutting $S$ cannot be full.

\textbf{Stage 4.} Now the molecule $\Mb_{UD}$ satisfies the assumptions (i)--(ii) in Definition \ref{def.alg_maincr}, and has $\#_{\mathrm{UD}}\geq \rho'''/6$. Moreover, since $\overline{Y}$ has been cut, we know that $\Mb_D$ is ov-proper and does not contain any weakly degenerate atom or atom pair. Next we perform the following operations, and repeat until this can no longer be done:
\begin{itemize}
\item If a subset $A\subseteq\Mb_U$ is not connected to $\Mb_{UD}\backslash A$ by any bond, then we cut it into elementary molecules using \textbf{UP}.
\item If a subset $A\subseteq \Mb_U$ has only deg 4 atoms, and is connected to $\Mb_{UD}\backslash A$ by only one bond (this bond can be connected to either $\Mb_U\backslash A$ or $\Mb_D$), then we cut $A$ \emph{as fixed} and subsequently cut it into elementary molecules using \textbf{UP} or \textbf{DOWN} (depending on the direction of this only bond).
\end{itemize}

In the above operation, there are three possibilities: (i) when there is no bond connecting $A$ to $\Mb_{UD}\backslash A$, (ii) when there is only one such bond connecting $A$ to an atom in $\Mb_U\backslash A$, or (iii) when the only bond connects $A$ to an atom in $\Mb_D$. Note that cases (i) and (ii) will not affect any good UD connections; in case (iii), $A$ is a connected component of $\Mb_U$ by itself, so cutting $A$ as fixed just removes one full component of $\Mb_U$; in particular case (iii) can occur at most $\#_{\mathrm{ovcp}(\Mb_U)}\leq 10(C_5^*)^{-1}\cdot\rho'''$ times, and each time reduces the number of good UD connections by at most 2 (it may break one UD connection and $e$, and turn at most one more good UD connection into not-good; note that if $e'$ and $e''$ are both connected to $e$ via $\Mb_D$ before cutting $A$, then they are also connected to each other via $\Mb_D$ after cutting $A$, and thus remain good).

As such, after \textbf{Stage 4}, we have that no subset $A\subseteq\Mb_U$ as described above, and we still have
\begin{equation}\label{eq.comb_est_case755}\#_{\mathrm{UD}}\geq \rho'''/6-20 (C_5^*)^{-1}\cdot\rho'''\geq \rho'''/7;\end{equation} The assumptions (i)--(ii) in Definition \ref{def.alg_maincr} are still preserved (if needed, we can further remove any full component of $\Mb_D$ that is not connected to $\Mb_U$ by a bond, and cut it into elementary molecules \emph{following Definition \ref{def.alg_maincr} (\ref{it.alg_maincr_9})}; the number of \{4\} molecules produced is at most $\#_{\mathrm{fucp}(\Mb_D)}\leq (C_3^*)^{-1}\cdot\rho'''$ and is negligible). Now define the sets $X_0$ and $X_1$ as in Definition \ref{def.alg_maincr}, and consider the value of $\#_{\mathrm{ovcp}(X_1)}$ as in Proposition \ref{prop.alg_maincr}.

\textbf{Stage 5.} If $\#_{\mathrm{ovcp}(X_1)}\geq 10^{-6}\upsilon^2\cdot\rho'''$, then we cut $\Mb_D$ as free and cut it into elementary molecules \emph{following Definition \ref{def.alg_maincr} (\ref{it.alg_maincr_9})} (except that if a component of $\Mb_D$ is full, then we first cut any deg 4 atom as free). If any atom in $\Mb_U$ becomes deg 2 then we cut it as free, and repeat until $\Mb_U$ contains no deg 2 atoms. This produces \{2\}, \{3\}, \{4\}, \{33A\} and \{33B\} molecules (with \{33B\} molecule only in $\Mb_D$ and \{2\} molecule only in $\Mb_U$), with the number of \{4\} molecules bounded by $\#_{\mathrm{fucp}(\Mb_D)}$; moreover all \{33B\} molecules are good due to absence of weakly degenerate atoms and atom pairs.

By definitions of $X_j$ in Definition \ref{def.alg_maincr}, it is clear that all the C-atoms in $X_0$ have been cut, and all the C-atoms in $X_1$ now have deg 3, and all the other C-atoms now have deg 4, see {\color{blue}Figure \ref{fig.x0_x1}} (all atoms in $X_0$ will subsequently become deg 2 and be cut when its ov-children have been cut, and cutting all atoms in $X_0$ leaves all atoms in $X_1$ as deg 3; on the other hand, if an atom does not belong to $X_0\cup X_1$ then it is not affected by cutting atoms in $\Mb_D\cup X_0$, or cutting any subsequent deg 2 O-atoms, cf. Remark \ref{rem.reg}). 

In addition, we claim that after the above cuttings, each atom $\nf\in X_1$ must belong to a component of $\Mb_U$ with at least two deg 3 atoms. In fact, let the component of $\Mb_U$ containing $\nf$ be $A$, then before the cuttings, $A$ must have one bond connecting to $\Mb_{UD}\backslash A$, and either has another such bond or a deg 3 atom (due to the operations in \textbf{Stage 4}). In either case, after the cutting operations, these two bonds (or one bond and one deg 3 atom) will lead to two fixed ends in $A$ (or a fixed end at an O-atom), hence the result.

Now consider all the components of $\Mb_U$ after the above cutting operations, that contain an C-atom in $X_1$ (i.e. ignore the other components, which can always be cut into elementary molecules using \textbf{UP}, in view of absence of C-top fixed ends in $\Mb_U$); denote the union of these components by $\widetilde{\Mb}_U$. In $\Mb_{UD}$ before cutting, the number of ov-components of $X_1$ is $\mu:=\#_{\mathrm{ovcp}(X_1)}\geq 10^{-6}\upsilon^2\cdot\rho'''$. After the cutting operations, the set of deg 3 atoms in $\widetilde{\Mb}_U$ equals $X_1\cup X_O$, where $X_O$ is the set of deg 3 atoms in $\widetilde{\Mb}_U$ that are O-atoms (see {\color{blue}Figure \ref{fig.x0_x1}}). Now note that, the number of ov-components of $X_1\cup X_O$ is no less than the number of ov-components of $X_1$ (which is obviously no less than $\mu$, because these cuttings can only break ov-segments but not create them). This is because for each deg 3 O-atom, it must have two edges (free ends or bonds) that are serial, and one free end that is serial with one fixed end; in other words one cannot change to a different ov-segment by going through this O-atom. Therefore, if two atoms in $X_1$ are ov-connected within $X_1\cup X_O$, they must also be ov-connected within $X_1$.

Then, we can apply \textbf{3COMPUP} to $\widetilde{\Mb}_U$ (as in toy model II), and use Proposition \ref{prop.3comp} to get
\begin{equation}\label{eq.comb_est_case75-}
\#_{\{33A\}}\geq (1/10)\cdot\mu\geq C^{-1}\cdot\rho''',
\end{equation}
\begin{equation}\label{eq.comb_est_case75+}\#_{\{4\}}\leq \#_{\mathrm{fucp}(\Mb_U)}+\#_{\mathrm{fucp}(\Mb_D)}+C(C_{2}^*)^{-1}\cdot\rho'\leq C(C_{2}^*)^{-1}\cdot \rho''',
\end{equation} which already proves (\ref{eq.comb_est_case51}) (with all \{33A\} molecules treated as good).

\textbf{Stage 6.} Finally, if $\#_{\mathrm{ovcp}(X_1)}\leq 10^{-6}\upsilon^2\cdot\rho'''$, then we apply \textbf{MAINUD} to cut $\Mb_{UD}$ into elementary molecules (as in toy model I plus). By Proposition \ref{prop.alg_maincr}, and note that all \{33B\} molecules are good (as they must have both atoms in $\Mb_D$, and $\Mb_D$ does not contain weakly degenerate atoms or atom pairs; note also that all the \{33B\} molecules must have both atoms being C-atoms, because the fixed end at any O-atom cannot be serial with the bond), we get that
\begin{equation}\label{eq.comb_est_case76}
\#_{\mathrm{good}}\geq \frac{1}{5}\cdot\big(\#_{\mathrm{UD}}-\#_{\mathrm{fucp}(\Mb_U)}-\#_{\mathrm{ovcp}(X_1)})\geq\frac{1}{35}\cdot\rho'''
\end{equation}
\begin{equation}\label{eq.comb_est_case76+}\qquad \#_{\{4\}}\leq \#_{\mathrm{fucp}(\Mb_U)}+\#_{\mathrm{ovcp}(X_1)}+ C(C_{2}^*)^{-1}\cdot \rho'''\leq 10^{-5}\upsilon^2\cdot \rho''',
\end{equation} which proves (\ref{eq.comb_est_case51}). This completes the proof of Proposition \ref{prop.case5}, and thus proves Proposition \ref{prop.comb_est}, in the case when $\Mb\in\Fc_{\Lambda_\ell}$.

\textbf{The finale: the case when $\Mb\in\Fc_{\Lambda_\ell}^{\mathrm{err}}$.} Finally we prove Proposition \ref{prop.comb_est} when $\Mb\in\Fc_{\Lambda_\ell}^{\mathrm{err}}$. First, note that we can cut $\Mb$ using \textbf{UP} which leads to $\#_{\mathrm{good}}\geq 0$ and $\#_{\mathrm{bad}}\leq|H|$ (see Part 1 of the proof of Proposition \ref{prop.case1}). In view of the extra term $C^*|H|$ on the right hand side of (\ref{eq.overall_alg}), we see that (\ref{eq.overall_alg}) is true if $\rho\lesssim C^*|H|$.

Now, suppose $\rho\gg C^*|H|(\gg C_{12}^*|H|)$, then we can repeat exactly the same proof in Sections \ref{sec.toy}--\ref{sec.maincr}. The only place where Proposition \ref{prop.mol_axiom} (\ref{it.axiom6}) is used (beyond Proposition \ref{prop.case1} which is already taken cared of), is in layer selection (Definition \ref{def.layer_select}). However, in this case we always have $\ell_1-1\leq\ell-1$, so Proposition \ref{prop.mol_axiom} (\ref{it.axiom6}) remains true for the set $r(\Mb_{\ell_1-1})$, therefore the proof goes exactly the same way as in the $\Mb\in\Fc_{\Lambda_\ell}$ case.
\end{proof}

\section{Proof of the $f^\Ac$ estimates}\label{sec.fa} In this section, we prove Proposition \ref{prop.est_fa}. By induction, we only need to prove (\ref{eq.approxfA_1})--(\ref{eq.approxfA_2}) for $\ell\tau$, assuming that they are already true for $(\ell-1)\tau$. For convenience, below we will denote $f^\Ac((\ell-1)\tau):=g$ and $(f^\Ac((\ell-1)\tau))_{\mathrm{rem}}:=g_{\mathrm{rem}}$, so $g=f((\ell-1)\tau)+g_{\mathrm{rem}}$.

Start from the expansion (\ref{eq.fAterm}), where $\Mb\in\Tc_{\Lambda_\ell}$ (or $\Tc_{\Lambda_\ell}^{\mathrm{err}}$) is an unlabeled molecule with one fixed root particle line, and $\Ic\Nc_\Mb$ and $|\Ic\Nc_\Mb|$ are defined in (\ref{eq.associated_integral_molecule})--(\ref{eq.associated_integral_molecule_abs}). Since we need to estimate these quantities with \emph{fixed} value $z$, in this section we will understand that, the unique top end $e$ for the root particle line $\pb$ will be a \emph{fixed end}\footnote{Note that this $\Mb$ might not be regular (Definition \ref{def.reg}), but this does not affect our proof until we get to the cutting algorithm in Section \ref{sec.fa_subleading}, at which time we will make $\Mb$ regular again by a special operation, see \textbf{Point 4} in Section \ref{sec.fa_subleading}.}, and all other ends of $\Mb$ will be free.

For each molecule $\Mb$ with single layer $\ell$, define the operator $\Ic_\Mb$ as in (\ref{eq.associated_int_op}), and define \begin{equation}\label{eq.Q_M_2}
Q_\Mb:=\mathbbm{1}_\Dc(\vt_\Mb)\cdot\prod_{e\in\Ec_{\mathrm{end}}^-}g(x_e+(\ell-1)\tau\cdot v_e,v_e),
\end{equation} where relevant notations are in Definition \ref{def.associated_op}. Then we have the following equalities (note also that $\Ic_\Mb(Q)$ is a function of the single variable $z$ which corresponds to the unique fixed end $e_{\mathrm{tr}}$):
\begin{equation}\label{eq.local_int_2}\Ic\Nc_\Mb(x,v)=(-1)^{|\Mb|_{\mathrm{O}}}\Ic_\Mb(Q_\Mb)(x-\tau v,v),\quad|\Ic\Nc_\Mb|(x,v)=\Ic_\Mb|Q_\Mb|(x-\tau v,v),
\end{equation} where $|\Mb|_{\mathrm{O}}$ is the number of O-atoms in $\Mb$.

To prove (\ref{eq.local_int_2}), note that $\Mb$ has single layer $\ell$, so $\ell_1[\pb]=\ell$ for all particle lines $\pb$. We may make the change of variables $t_\nf\mapsto t_\nf-(\ell-1)\tau$ and $(x_e,v_e)\mapsto (x_e+(\ell-1)\tau\cdot v_e,v_e)$ in the integral (\ref{eq.associated_int_op}), and reduce to the case $\ell=1$. Then (\ref{eq.local_int_2}) follows from the same proof as in Proposition \ref{prop.local_int}, by replacing $Q_\Mb$ with $Q_\Mb\cdot \dirac(z_{e_{\mathrm{tr}}}-(y,v))$, where $(y,v)=(x-\tau v,v)$ represents the final state vector of the root particle transported backwards to time $0$.
 
 In the rest of this section, we will perform the above change of variables and thereby assume $\ell=1$. Now, using (\ref{eq.local_int_2}), we can rewrite the sum (\ref{eq.fAterm}) as
\begin{equation}\label{eq.fAterm2}f^\Ac(\ell\tau)=f^{\Ac,\mathrm{led}}+f^{\Ac,\mathrm{sub}},
\end{equation} where (recall $|O_1(A)|\leq |A|$)
\begin{equation}\label{eq.fAterm3}
\begin{aligned}f^{\Ac,\mathrm{led}}(x,v)&=\sum_{\substack{\Mb\in\Tc_{\Lambda_\ell}\\\rho(\Mb)=0,\,|\Mb|_p<\Lambda_\ell}}(-1)^{|\Mb|_{\mathrm{O}}}\Ic_\Mb(Q_\Mb)(x-\tau v,v),\\
f^{\Ac,\mathrm{sub}}(x,v)&=\sum_{\substack{\Mb\in\Tc_{\Lambda_\ell},\,\rho(\Mb)>0\\\mathrm{or\ }|\Mb|_p\geq\Lambda_\ell,\mathrm{\ or\ }\Mb\in\Tc_{\Lambda_\ell}^{\mathrm{err}}}}O_1\big(\Ic_\Mb|Q_\Mb|\big)(x-\tau v,v)
\end{aligned}
\end{equation} are the leading and sub-leading terms in $f^\Ac$ respectively.

We will study the leading term $f^{\Ac,\mathrm{led}}$ in Section \ref{sec.fa_leading}, and establish its asymptotics (see (\ref{eq.fa_est_3})) in Section \ref{sec.approx}. In Section \ref{sec.fa_subleading}, we will study the sub-leading term $f^{\Ac,\mathrm{sub}}$ and prove it vanishes in the limit. Putting these together, we can then complete the proof of Proposition \ref{prop.est_fa} at the end of Section \ref{sec.fa_subleading}.
\subsection{The leading terms}\label{sec.fa_leading} Start with the term $f^{\Ac,\mathrm{led}}$ in (\ref{eq.fAterm3}). Here $\Mb$ is a tree as $\rho(\Mb)=0$ and $\Mb$ is connected (by Definition \ref{def.set_T_F} (\ref{it.set_F_l_2}), note $\Mb$ has only one root particle line); moreover, it is straightforward to verify that $\Mb$ is a tree and $|\Mb|_p<\Lambda_\ell$ actually implies $\Mb\in\Tc_{\Lambda_\ell}$. Therefore, the summation in $f^{\Ac,\mathrm{led}}$ in (\ref{eq.fAterm3}) is taken with \emph{only} the assumptions that $\Mb$ is a tree and $|\Mb|_p<\Lambda_\ell$. Note also that $|\Mb|_p=|\Mb|+1$ (i.e. number of collisions equals number of particles minus $1$, in the absence of recollisions).

Now we introduce the following definition.
\begin{definition}\label{def.fa_notion} Let $\Mb$ be a molecule with single layer $\ell$ and one root particle line, which we denote by $\pb_{\mathrm{rt}}$. Define the highest atom in $\pb_{\mathrm{rt}}$ to be $\nf_{\mathrm{rt}}$, and the top end in $\pb_{\mathrm{rt}}$ to be $e_{\mathrm{rt}}$, which is the unique fixed end of $\Mb$.

Define $\Mf_n$ to be the set of trees $\Mb$ with $\rho(\Mb)=0$ and $|\Mb|=n$ (with $|\Mb|_p=n+1<\Lambda$). Define also $\Mf_n^*$ to be the subset of $\Mf_n$ formed by those $\Mb$ such that every atom in $\Mb$ is a descendant of $\nf_{\mathrm{rt}}$.
\end{definition}

It is known in \cite{BGSS22_2} that the contribution of the molecules $\Mb\in\Mf_n\backslash \Mf_n^*$ to $f^{\Ac,\mathrm{led}}$ in (\ref{eq.fAterm3}) is zero, due to a simple cancellation between the terms. We now present this argument, in a similar manner to \cite{BGSS22_2}.
\begin{proposition}\label{prop.cancel} There is an involution $\Mb\mapsto\Mb^{\mathrm{twi}}$ from $\Mb\in\Mf_n\backslash \Mf_n^*$ to itself, such that $|\Mb'|_{\mathrm{O}}=|\Mb|_{\mathrm{O}}\pm 1$, and $\Ic_{\Mb'}(Q_{\Mb'})=\Ic_{\Mb}(Q_\Mb)$ for $\Mb'=\Mb^{\mathrm{twi}}$. In particular, this implies that
\begin{equation}\label{eq.cancel}\sum_{\Mb\in\Mf_n\backslash\Mf_n^*}(-1)^{|\Mb|_{\mathrm{O}}}\Ic_{\Mb}(Q_\Mb)=0.
\end{equation}
\end{proposition}
\begin{proof} If $\Mb\in\Mf_n\backslash \Mf_n^*$, then $\Mb\neq S_{\nf_{\mathrm{tr}}}$ (where $S_\nf$ is the set of descendants of $\nf$, see Definition \ref{def.molecule_order}). Choose a highest atom $\mf\in\Mb\backslash S_{\mathrm{tr}}$. If there is more than one $\mf$, we can select one in an arbitrary way (for example by introducing the artificial total ordering in Part 2 of the proof of Proposition \ref{prop.local_int} and select the highest $\mf$). Then define $\Mb'=\Mb^{\mathrm{twi}}$ to be formed from $\Mb$ by turning $\mf$ into a C/O atom $\mf'$ if $\mf$ is originally an O/C atom.

Clearly $|\Mb'|_{\mathrm{O}}=|\Mb|_{\mathrm{O}}\pm 1$. To prove $\Ic_{\Mb'}(Q_{\Mb'})=\Ic_{\Mb}(Q_\Mb)$, note that $\mf\in\Mb$ (and $\mf'\in\Mb'$) has no parent, thus it has two top free ends $e_1'$ and $e_2'$ (the fixed end $e_{\mathrm{rt}}\not\in\{e_1',e_2'\}$). Now we compare the integrals $\Ic_{\Mb'}(Q_{\Mb'})$ and $\Ic_{\Mb}(Q_\Mb)$ as in (\ref{eq.associated_int_op}). Note that $Q_\Mb=Q_{\Mb'}$ does not depend on $z_{e_1'}$ and $z_{e_2'}$, thus the only difference in the two integrals occurs in the $\Dirac_\mf$ and $\Dirac_{\mf'}$ factors. However, by (\ref{eq.singlecol}), if we integrate out the variables $(z_{e_1'},z_{e_2'},t_\mf)$, then the result would be the same for C and O-atoms, provided that the $Q$ in (\ref{eq.singlecol}) does not depend on $z_{e_1'}$ and $z_{e_2'}$ (which is true in our case). This proves $\Ic_{\Mb'}(Q_{\Mb'})=\Ic_{\Mb}(Q_\Mb)$, then (\ref{eq.cancel}) follows consequently.
\end{proof}
Thanks to (\ref{eq.cancel}), we can rewrite the sum in $f^{\Ac,\mathrm{led}}$ in (\ref{eq.fAterm3}) as 
 \begin{equation}\label{eq.fB_3Moleq_3}
f^{\Ac,\mathrm{led}}(x,v)=g(x-\tau v,v)+\sum_{1\leq n< \Lambda_\ell-1}\sum_{\Mb\in \Mf_n^*}(-1)^{|\Mb|_{\mathrm{O}}} \cdot\Ic_\Mb(Q_\Mb)(x-\tau v,v).
\end{equation} Note that, for $\Mb\in\Mf_n^*$, we have $\Mb=S_{\nf_{\mathrm{tr}}}$ is a binary tree in the sense that each atom other than $\nf_{\mathrm{tr}}$ has a unique parent (and one top fixed end), and two bottom free ends/children in total. As such, it will be convenient to rename the variables at the node $\nf$ as
\begin{equation}\label{eq.renamingvar}
(z_{e_1},z_{e_2},z_{e_1'},z_{e_2'})\leftrightarrow (z_1^\nf, z_2^\nf, z_3^\nf, z_4^\nf)
\end{equation}
where $z_2^\nf$ is the variable corresponding to the top free end at $\nf$, $z_1^\nf$ is the one corresponding to the bond connecting $\nf$ to its parent (or the free end $e_{\mathrm{tr}}$ when $\nf=\nf_{\mathrm{tr}}$), $z_4^\nf$ is the one corresponding to the bottom edge at $\nf$ that is serial to the top free end, and $z_3^\nf$ is the remaining variable. Given the symmetries of $\boldsymbol{\Delta}_\nf$, we always have $\boldsymbol{\Delta}_\nf(z_1^\nf, z_2^\nf, z_3^\nf, z_4^\nf)=\boldsymbol{\Delta}_\nf(z_{e_1},z_{e_2},z_{e_1'},z_{e_2'})$.
 
Next, we will split $g=f((\ell-1)\tau)+g_{\mathrm{rem}}$ for each input factor in $\Ic_\Mb(Q_\Mb)$ in (\ref{eq.fB_3Moleq_3}). We expect those terms containing at least one $g_{\mathrm{rem}}$ should vanish in the limit, in view of (\ref{eq.approxfA_2}) satisfied by $g_{\mathrm{rem}}$; to prove this, we need the next proposition.
 \begin{proposition}\label{prop.fa_est_1} Let $\Mb\in\Mf_n^*$, and let $\widetilde{Q}$ be defined as in (\ref{eq.Q_M_2}) but with the input factors $g$ replaced by $h_e$ which depends on $e$. Recall also the $\beta_\ell$ as in (\ref{eq.decayseq}), and define $\widetilde{\beta}_\ell:=(\beta_\ell+\beta_{\ell-1})/2$. Then we have
 \begin{equation}\label{eq.fB4lossy}
\|\Ic_\Mb(\widetilde{Q})(x-\tau v,v)\|_{\mathrm{Bol}^{\widetilde{\beta}_\ell}}\leq |\log\varepsilon|^{C^*}\cdot C^n  \tau^{n/2} \prod_{e\in \Ec^-_{\mathrm{end}}} \|h_e\|_{\mathrm{Bol}^{\beta_{\ell-1}}}.
\end{equation}
 \end{proposition}
\begin{proof} Recall that $\widetilde{\beta_\ell}=\beta_{\ell-1}-(\beta/40\Lf)$ by definition. We divide the proof into 3 parts.

\textbf{Proof part 1: adjusting Maxwellian weights.} Without loss of generality, we may assume $h_e\geq 0$ and $\|h_e\|_{\mathrm{Bol}^{\beta_{\ell-1}}}=1$ for all $e\in \Ec^-_{\mathrm{end}}$. Recall $\Ec_{\mathrm{end}}^\pm$ is the set of all top/bottom ends of $\Mb$, and $e_{\mathrm{tr}}\in\Ec_{\mathrm{end}}^+$ is the unique fixed end; by energy conservation we have $\sum_{e\in \Ec^-_{\mathrm{end}}}|v_e|^2=\sum_{e\in \Ec^+_{\mathrm{end}}}$. Using also the notion $z_j^\nf$ in (\ref{eq.renamingvar}), this implies that
\begin{equation}\label{eq.consenergys9+}\sum_{e\in \Ec_{\mathrm{end}}^-}
\beta_{\ell-1}|v_e|^2=\widetilde{\beta}_\ell|v|^2+\sum_{e\in \Ec_{\mathrm{end}}^-}
\frac{\beta}{40\Lf}|v_e|^2+\sum_{\nf}\widetilde{\beta}_{\ell}|v_2^\nf|^2,
\end{equation} where $v=v_{e_{\mathrm{tr}}}$ is the output variable in the function $\Ic_\Mb(\widetilde{Q})(x-\tau v,v)$. Since $\|h\|_{\mathrm{Bol}^\beta}=\|e^{\beta|v|^2}h\|_{\mathrm{Bol}^0}$ by (\ref{eq.boltzmann_decay_2}), we can replace $h_e$ by $e^{\beta_{\ell-1}|v_e|^2}h_e$ and apply (\ref{eq.consenergys9+}), to reduce (\ref{eq.fB4lossy}) to the following estimate:
\begin{equation}\label{eq.fB4lossy2}
\|\Ic_\Mb'(\widetilde{Q})(x-\tau v,v)\|_{\mathrm{Bol}^{0}}\leq |\log\varepsilon|^{C^*}\cdot C^n  \tau^{n/2},\quad \mathrm{given\ }\|h_e\|_{\mathrm{Bol}^0}=1,
\end{equation}
where $\Ic_\Mb'(\widetilde{Q})=\Ic_\Mb\left(\widetilde{Q}\cdot \prod_{\nf} e^{-\widetilde{\beta}_\ell|v_2^\nf|^2}\cdot\prod_{e\in \Ec_{\mathrm{end}}^-}e^{-(\beta/40\Lf)|v_e|^2}\right)$.

\textbf{Proof part 2: decomposition of $\widetilde{Q}$.} By the exponential decay in the formula of $\Ic_\Mb'(\widetilde{Q})$, it is easy to see that we can restrict the support of $Q$ to $|v_e|\leq |\log\varepsilon|^{C^*}$ for each $e\in\Ec_{\mathrm{end}}^-$, and therefore for all $e\in\Ec$ by energy conservation, in the same way as in Part 2 of the proof of Proposition \ref{prop.cumulant_est}. Moreover, by repeating the same proof there (in Section \ref{sec.summary}), we can also decompose
\begin{equation}\label{eq.extra_decomp}
\widetilde{Q}=\sum_{\alpha} w_{\alpha}\widetilde{Q}_{\alpha},\quad \sum_\alpha w_{\alpha}\leq C^n,
\end{equation} where each $\widetilde{Q}_{\alpha}$ is uniformly bounded and is supported in the set $|x_e-x_e^*|\leq |\log\varepsilon|^{C^*}$ for each $e\in\Ec_{\mathrm{end}}^-$, with some choice of $(x_e^*)$ depending on $\alpha$.

Since $|v_e|\leq |\log\varepsilon|^{C^*}$, this support condition implies that $|x_{e_{\mathrm{tr}}}-x_{e_0}^*|\leq |\log\varepsilon|^{C^*}$ where $e_0$ is the bottom end of the root particle line ($e_{\mathrm{tr}}$ is its top end). As we have $x=x_{e_{\mathrm{tr}}}+\tau \cdot v_{e_{\mathrm{tr}}}$ for the variable $x$ in (\ref{eq.fB4lossy}), we know that for each $\alpha$, this $x$ belongs to a fixed ball of radius $|\log\varepsilon|^{C^*}$. This ball can be divided into $\leq |\log\varepsilon|^{C^*}$ unit cubes, so under this support condition we have $\|f\|_{\mathrm{Bol}^0}\leq|\log\varepsilon|^{C^*}\cdot\|f\|_{L^\infty}$ due to (\ref{eq.boltzmann_decay_2}).

Summing up, using (\ref{eq.extra_decomp}) and the above discussions, we see that (\ref{eq.fB4lossy2}) would follow if we can prove the following estimate (uniformly in $\alpha$):
\begin{equation}\label{eq.fB4lossy3}
\|\Ic_\Mb'(\widetilde{Q}_\alpha)\|_{L^\infty}\leq (C\tau^{1/2})^n,\quad \mathrm{given\ }\|\widetilde{Q}_\alpha\|_{L^\infty}=1.
\end{equation}

\textbf{Proof part 3: estimating the integral.} To prove \eqref{eq.fB4lossy3}, we first invoke \eqref{eq.weight_est_proof_5} in Proposition \ref{prop.weight} with $v_{e}$ replaced by $\Lf^{-1/2}v_e$, which allows us to bound the product of the cross-sections in the $\Dirac(\cdots)$ factors by
$$\prod_{\nf}(1+|v_{e_1}-v_{e_2}|)\leq \Lf^{n/2}\prod_{\nf}(1+|\Lf^{-1/2}v_{e_1}-\Lf^{-1/2}v_{e_2}|)\leq C^n \Lf^{n/2}  \cdot \mathfrak C_\Mb^{5/6}\cdot \exp\bigg(\frac{\beta}{40\Lf}\sum_{e\in\Ec_{\mathrm{end}}^-}|v_e|^2\bigg),
$$
where the products are taken over all atoms $\nf$ with $(v_{e_1},v_{e_2})$ being the two bottom edges at $\nf$, and $\mathfrak C_\Mb$ is defined in \eqref{eq.weight_est_proof_1} with $T$ replaced by $\Mb$ (which is also a binary tree as in Part 1 of the proof of Proposition \ref{prop.weight}). Here we have used the elementary inequality $\exp(C|v|^{4/3})\leq C \exp((\beta/40)|v|^2)$. This allows to bound 
\begin{equation*}
\begin{aligned}
\|\Ic_\Mb'(\widetilde{Q}_\alpha)\|_{L^\infty}&\leq C^n \Lf^{n/2}  \cdot \mathfrak C_\Mb^{5/6}\cdot \bigg(\int_{\Dc} \,\mathrm{d}\vt_\Mb\bigg)\cdot \sup_{\vt_\Mb\in \Dc} \|V_n(z, \vt_\Mb)\|_{L^\infty},\\
V_n(z, \vt_\Mb)&=\varepsilon^{-(d-1)n}\int_{\mathbb{R}^{2d|\Ec_*|}} \prod_{\nf\in\Mb}e^{-\widetilde{\beta}_\ell |v_2^\nf|^2} \boldsymbol{\Delta}^\dagger(z^\nf_{1},z_{2}^\nf,z_{3}^\nf,z_{4}^\nf,t_\nf)\,\mathrm{d}z_3^\nf \mathrm{d}z_4^\nf\mathrm{d}z_2^\nf,
\end{aligned}
\end{equation*}
where $\Dc$ is as in (\ref{eq.associated_domain}) but with $\ell'$ replaced by $1$, and 
$\boldsymbol{\Delta}^\dagger$ is as in (\ref{eq.associated_dist_C})--(\ref{eq.associated_dist_O}) but without the $[(v_{e_1}-v_{e_2})\cdot \omega]_-$ factor, and we have relabeled the variables at each atom $\nf \in \Mb$ using the variables $z_j^\nf$ as in \eqref{eq.renamingvar}. Note also that $|\Ec_*|-2|\Mb|=n$ in (\ref{eq.associated_int_op}) using Part 1 of the proof of Proposition \ref{prop.local_int} (and $|\Ec|-|\Ec_*|=1$ and $|\Mb|_p=n+1$).

Since $\Cf_\Mb\geq 1$ and $\int_\Dc \mathrm{d}\vt_\Mb=\tau^n \cdot\Cf_\Mb^{-1}$ by \eqref{eq.weight_est_proof_1}, and recall that $\tau \Lf=t_{\mathrm{fin}}\leq C$, we see that \eqref{eq.fB4lossy3} will follow once we show that $\sup_{\vt_\Mb\in \Dc} \|V_n(z, \vt_\Mb)\|_{L^\infty}\leq C^n$. The latter is straightforward as one calculates the integrals from top to bottom in the binary tree $\Mb$ as follows: for each $\nf \in \Mb$, we fix $z_1^\nf$ and integrate in the variables $(z_2^\nf, z_3^\nf, z_4^\nf)$, noting that for any node $\nf \in \Mb$ we have 
\begin{equation*}
\varepsilon^{-(d-1)}\int_{(\Rb^{2d})^3}e^{-\widetilde{\beta}_\ell |v_2^\nf|^2} \cdot\boldsymbol{\Delta}^\dagger(z^\nf_{1},z_{2}^\nf,z_{3}^\nf,z_{4}^\nf,t_\nf)\,\mathrm{d}z_3^\nf \mathrm{d}z_4^\nf\mathrm{d}z_2^\nf=\int_{\Sb^{d-1}\times\Rb^d} e^{-\widetilde{\beta}_\ell|v_2^\nf|^2}\,\mathrm{d}\omega\mathrm{d}v_2^\nf\leq C
\end{equation*}
by a simple variant of (\ref{eq.intmini_3}). This finishes the proof.
\end{proof}

Now, using Proposition \ref{prop.fa_est_1}, we may split $g=f((\ell-1)\tau)+g_{\mathrm{rem}}((\ell-1)\tau)$ for each input factor in in $\Ic_\Mb(Q_\Mb)$ in (\ref{eq.fB_3Moleq_3}). Define $f_1^{\Ac,\mathrm{led}}$ to be the same as in (\ref{eq.fB_3Moleq_3}) but with all input factors $g=f^\Ac((\ell-1)\tau)$ replaced by $f((\ell-1)\tau)$. By Proposition \ref{prop.fa_est_1} and (\ref{eq.approxfA_2}) and (\ref{eq.boltzmann_decay_4}), it is easy to see that
\begin{equation}\label{eq.fa_est_2}\big\|f^{\Ac,\mathrm{led}}-f_1^{\Ac,\mathrm{led}}\big\|_{\mathrm{Bol}^{\beta_\ell}}\leq|\log\varepsilon|^{C^*}\sum_{n\geq 0}(C\tau^{1/2})^n\cdot\varepsilon^{\theta_{\ell-1}}\leq \varepsilon^{\theta_\ell}/10.
\end{equation} Here we have used the fact that the number of molecules $\Mb\in\Mf_n^*$ is bounded by $C^n$, which follows from the well-known upper bound for the number of binary trees.

\subsection{Approximation by Boltzmann iterates}\label{sec.approx} In view of (\ref{eq.fa_est_2}), now we only need to estimate $f_1^{\Ac,\mathrm{led}}$. In this subsection, we will approximate it by the iterates, up to order $\Lambda-2$, of the Boltzmann equation. To do this, we first perform a change of variables in the integral appearing in the definition of $f_1^{\Ac,\mathrm{led}}$. We adopt the same relabeling of the variables $z_j^\nf$ at each atom $\nf\in\Mb$, then there are two cases of $\Dirac=\Dirac_\nf$ in (\ref{eq.fB_3Moleq_3}):
\begin{itemize}
\item When $\nf$ is C-atom, we have $\Dirac_+$ defined by
\begin{equation}
\label{eq.newboldDelta}
\begin{aligned}
\boldsymbol{\Delta}_+(z^{\nf}_{1},z^{\nf}_{2},z^{\nf}_{3},z^{\nf}_{4},t_\nf)&:=\boldsymbol{\delta}\big(x^{\nf}_{1}-x^{\nf}_{3}+t_\nf(v^{\nf}_{1}-v^{\nf}_{3})\big)\cdot \boldsymbol{\delta}\big(x^{\nf}_{2}-x^{\nf}_{4}+t_\nf(v^{\nf}_{2}-v^{\nf}_{4})\big)\\&\times\boldsymbol{\delta}\big(v^\nf_{3}-v^\nf_{1}+[(v^\nf_{1}-v^\nf_{2})\cdot\omega]\omega\big)\cdot\boldsymbol{\delta}\big(v^\nf_4-v^\nf_{2}-[(v^\nf_{1}-v^\nf_{2})\cdot\omega]\omega\big)\\&\times \boldsymbol{\delta}\big(|x^\nf_{2}-x^\nf_{1}+t_\nf(v^\nf_{2}-v^\nf_{1})|-\varepsilon\big)\cdot \big[(v^\nf_{1}-v^\nf_{2})\cdot \omega\big]_+;
\end{aligned}
\end{equation}
\item When $\nf$ is O-atom, we have $\Dirac_-$ defined by
\begin{equation}\label{eq.newboldDelta_2}
\begin{aligned}
\boldsymbol{\Delta}_-(z^{\nf}_{1},z^{\nf}_{2},z^{\nf}_{3},z^{\nf}_{4},t_\nf)&:=\boldsymbol{\delta}(x^{\nf}_{1}-x^{\nf}_{3})\cdot \boldsymbol{\delta}(x^{\nf}_{2}-x^{\nf}_{4})\cdot \boldsymbol{\delta}\big(v^\nf_{3}-v^\nf_{1}\big)\cdot\boldsymbol{\delta}\big(v^\nf_4-v^\nf_{2}\big)\\
&\times\boldsymbol{\delta}\big(|x^\nf_{2}-x^\nf_{1}+t_\nf(v^\nf_{2}-v^\nf_{1})|-\varepsilon\big)\cdot \big[(v^\nf_{1}-v^\nf_{2})\cdot \omega\big]_+.
\end{aligned}
\end{equation}
\end{itemize}
Here we have used that  $\big[(v_{e_1}-v_{e_2})\cdot \omega\big]_- =\big[(v^\nf_{3}-v^\nf_{4})\cdot \omega\big]_-=\big[(v^\nf_{1}-v^\nf_{2})\cdot \omega\big]_+$ for the given relabeling in \eqref{eq.renamingvar}; also recall that $\omega:=\varepsilon^{-1}(x_{e_1}-x_{e_2}+t_\nf(v_{e_1}-v_{e_2}))=\varepsilon^{-1}(x^\nf_{2}-x^\nf_{1}+t_\nf(v^\nf_{2}-v^\nf_{1}))$. Now, recall that the integration of the $\boldsymbol{\Delta}_{\pm}$ factor at each atom is performed by first integrating in $(z_3^\nf, z_4^\nf)$ (with $(z_1^\nf, z_2^\nf)$ fixed), followed by integrating in $z^\nf_2$ followed by that in $z^\nf_1$ and $t_\nf$. Since the integral over $x^\nf_2$ involves a $\boldsymbol{\delta}$ measure supported on a sphere centered at $x_1^\nf-t_\nf(v_2^\nf -v_1^\nf)$, we can make the change variables $x_2^\nf=x_1^\nf-t_\nf(v_2^\nf -v_1^\nf)+\varepsilon \omega_\nf$ as in the proof of Proposition \ref{prop.intmini}, and write
\[\dirac\big(|x^\nf_{2}-x^\nf_{1}+t_\nf(v^\nf_{2}-v^\nf_{1})|-\varepsilon\big)=\varepsilon^{d-1}\int_{\Sb^{d-1}}\dirac\big(x_1^\nf-x_2^\nf-t_\nf(v_2^\nf -v_1^\nf)+\varepsilon \omega_\nf\big)\,\mathrm{d}\omega_\nf.\] As a result, we can now rewrite the integral defining $f_1^{\Ac,\mathrm{led}}$ as follows: 
\begin{equation}\label{eq.fB_5}
 \begin{aligned}
f_1^{\Ac,\mathrm{led}}(x,v)&=f\big((\ell-1)\tau, x-\tau v, v\big)+\sum_{1\leq n<\Lambda_\ell-1}\sum_{\Mb\in \mathfrak M_{n}^*} (-1)^{|\Mb|_{\mathrm{O}}}\cdot\Ic_\Mb(\widetilde{Q}_\Mb)(x-\tau v, v),\\
\Ic_\Mb(\widetilde{Q}_\Mb)&:=\int_{\Dc}\int_{(\Sb^{d-1})^{|\Mb|}}\int_{\mathbb{R}^{2d|\Ec_*|}} \prod_{\nf\in\Mb}\widetilde {\boldsymbol{\Delta}}_{\pm}(z^{\nf}_{1},z^{\nf}_{2},z^{\nf}_{3},z^{\nf}_{4},\omega_\nf, t_\nf)\prod_{e\in \Ec_{\mathrm{end}}^-}f(({\ell-1})\tau, z_e)\,\mathrm{d}\vz_{\Ec_*}\mathrm{d}\boldsymbol{\omega}_\Mb\mathrm{d}\vt_{\Mb},
\end{aligned}
\end{equation} 
where $\mathrm{d}\boldsymbol{\omega}_\Mc:=\prod_{\nf \in \Mc}\mathrm{d}\omega_\nf$ and $\widetilde {\boldsymbol{\Delta}}_{\pm}$ is the same as ${\boldsymbol{\Delta}_{\pm}}$ in \eqref{eq.newboldDelta}--\eqref{eq.newboldDelta_2} except that $\omega$ is replaced by $\omega_\nf$ and the last $\boldsymbol \delta$ in both expressions is replaced by $\boldsymbol{\delta}(x_1^\nf-x_2^\nf-t_\nf(v_2^\nf -v_1^\nf)+\varepsilon \omega_\nf)$.

Next, we will explain the relation of the above expansion in \eqref{eq.fB_5} with the Duhamel expansion of the solution $f(t, z)$ to the Boltzmann equation on the interval $[(\ell-1) \tau, \ell \tau]$. 
\begin{proposition}\label{prop.BoltzExpansion}
The solution $f(s, z)$ of the Botlzmann equation on the interval $[(\ell-1) \tau, \ell \tau]$
can be written as $f(s, z):=F(s-(\ell-1) \tau, x-sv+(\ell-1)\tau v, v)$ where $F(t, z)$ is defined for $t\in [0, \tau]$ as
\begin{equation}\label{eq.tildenexp}
 \begin{aligned}
F(t, z)&:=F^{\mathrm{(app)}}+ F^{\mathrm{(rem)}},\\
F^{\mathrm{(app)}}&:=f((\ell-1)\tau,z)+\sum_{1\leq k<\Lambda_\ell-1} F^{(k)}(t,z), \qquad F^{\mathrm{(rem)}}:=\sum_{k\geq \Lambda_\ell-1}F^{(k)}(t, z),
\end{aligned}
\end{equation} and $F^{(k)}(t,z)$ is defined by
\begin{equation}\label{eq.tildenexp_2}
\begin{aligned}
F^{(k)}(t, z)&:= \sum_{\Mb\in \mathfrak M_{k}^*}(-1)^{|\Mb|_{\mathrm{O}}}\cdot\widetilde{\Ic}_\Mb(\widetilde{Q}_\Mb)(t,z),\\
\widetilde{\Ic}_\Mb(\widetilde{Q}_\Mb)&:=\int_{\Dc^t}\int_{(\Sb^{d-1})^{|\Mb|}}\int_{\mathbb{R}^{2d|\Ec_*|}} \prod_{\nf\in\Mb}\widetilde {\boldsymbol{\Delta}}_{\pm}^{(0)}(z^{\nf}_{1},z^{\nf}_{2},z^{\nf}_{3},z^{\nf}_{4}, \omega_\nf,t_\nf)\prod_{e\in \Ec_{\mathrm{end}}^-}f(({\ell-1})\tau, z_e)\,\mathrm{d}\vz_{\Ec_*}\mathrm{d}\boldsymbol{\omega}_\Mb\mathrm{d}\vt_{\Mb}.
\end{aligned}
\end{equation} 
Here $\Dc^t$ is the same as $\Dc$ in (\ref{eq.associated_domain}) but with $\ell'$ replaced by $1$ and with the extra restriction that $t_\nf<t$ for all $\nf\in\Mb$, and $\widetilde {\boldsymbol{\Delta}}^{(0)}_{\pm}$ is the same as ${\boldsymbol{\Delta}_{\pm}}$ in \eqref{eq.newboldDelta} but with $\omega$ replaced by $\omega_\nf$ and the last $\boldsymbol \delta$ in both expressions replaced by $\boldsymbol{\delta}(x_1^\nf-x_2^\nf-t_\nf(v_2^\nf -v_1^\nf))$. Note that this makes $\widetilde{\Ic}_\Mb(\widetilde{Q}_\Mb)$ independent of $\varepsilon$. Moreover, the following estimate holds for any $\Mb\in \mathfrak M_{k}^*$:
\begin{equation}\label{eq.INBoltzest}
\big\|\widetilde{\Ic}_\Mb(\widetilde{Q}_\Mb)(x-\tau v,v)\big\|_{\mathrm{Bol}^{\widetilde{\beta}_\ell}}\leq C^k \tau^{k/2} \prod_{e\in \Ec^-_{\mathrm{end}}} \|f((\ell-1)\tau)\|_{\mathrm{Bol}^{\beta}}.
\end{equation}

\end{proposition}
\begin{proof} First notice that $F$ satisfies the following integral equation on the interval $[0, \tau]$:
\begin{equation}\label{eq.BoltIntEq}
\begin{aligned}
F(t,z)&=f((\ell-1)\tau, z)+\Cc(F,F) (t,z),\\
\Cc(f, g)(t,z)&=\int_0^t\int_{\Rb^d}\int_{\Sb^{d-1}}\big[(v-v_1)\cdot\omega\big]_+\cdot\big(f(s, x+vs-v's, v')g(s, x+vs-v_1's, v_1')\\
&\phantom{\int_0^t\int_{\Rb^d}\int_{\Sb^{d-1}}\big[(v-v_1)\cdot\omega\big]_+}-f(s, x, v)g(s, x+vs-v_1s, v_1)\big)\,\mathrm{d}\omega\mathrm{d}v_1\mathrm{d}s\\
&=\sum_{(\pm)}\pm\int_0^t\int_{\Sb^{d-1}}\int_{\Rb^{6d}}\widetilde{\Dirac}_{\pm}^{(0)}(z,z_2,z_3,z_4,s)\cdot f(s,z_3)\cdot g(s,z_4)\prod_{i=2}^4\mathrm{d}z_i\mathrm{d}\omega\mathrm{d}s,
\end{aligned}
\end{equation} 
where $v', v_1'$ are defined from $(v,v_1)$ as in \eqref{eq.boltzmann2}.

Next, notice that $F^{(k)}$ defined in \eqref{eq.tildenexp_2} satisfies the following rrecurrence relation:
$$
F^{(k)}=\sum_{k_1+k_2=k-1}\Cc(F^{(k_1)}, F^{(k_2)}),\qquad F^{(0)}(t,z):=f((\ell-1)\tau,z),
$$ which can be checked readily by induction using the last equation in \eqref{eq.BoltIntEq}. Here we note that each molecule $\Mb\in\Mf_k^*$ consists of one highest atom $\nf_{\mathrm{tr}}$ which is either C-atom (corresponding to $\widetilde{\Dirac}_+^{(0)}$ in \eqref{eq.BoltIntEq}) or O-atom (corresponding to $\widetilde{\Dirac}_-^{(0)}$ in \eqref{eq.BoltIntEq}), and two sub-molecules $\Mb^i\in\Mf_{k_i}^*$ which are the sets of descendants of the two children $\nf_i$ of $\nf^*$, with $k_1+k_2=k-1$. The cases when some of these two children are absent are modified trivially (with $k_1=0$ or $k_2=0$).

The convergence of the infinite series $\sum_{k=0}^\infty F^{(k)}$ then follows, once we establish \eqref{eq.INBoltzest}, in view that $\mathfrak M_k$ contains $O(C^k)$ molecules (by counting binary trees), and our choice of $\tau$. The fact that this sum solves \eqref{eq.BoltIntEq} follows by writing down the equation satisfied by $F^{\leq k}=\sum_{j \leq k} F^{(i)}$, which reads
\begin{equation}\label{eq.T_Bol_err}
F^{(0)}+\Cc(F^{\leq k}, F^{\leq k})-F^{\leq k}=\sum_{k+1\leq i \leq 2k}\sum_{i_1+i_2=i-1}\Cc(F^{(i_1)}, F^{(i_2)}),
\end{equation}
and noticing that $\Cc(F^{(i_1)}, F^{(i_2)})$ is an algebraic sum of at most $C^i$ terms $\widetilde{\Ic}_\Mb(\widetilde{Q}_\Mb)(t,z)$ where $\Mb\in \mathfrak M_{i}^*$. This gives that the right hand side of (\ref{eq.T_Bol_err}) can be bounded by $C^k \tau^{k/2}$ using \eqref{eq.INBoltzest} again, which converges geometrically to zero. This gives that $F=\sum_{k=0}^\infty F^{(k)}$. 

As such, it remains to justify \eqref{eq.INBoltzest}. But, this follows from exactly the same proof of Proposition \ref{prop.fa_est_1}, only with $\Ic_\Mb$ replaced by $\widetilde{\Ic}_{\Mb}$. 
\end{proof}
Now we can compare $f_1^{\Ac,\mathrm{led}}$ with $f(\ell\tau)$. The key step is contained in the next proposition.
\begin{proposition}\label{prop.fa_est_2} With the expressions defined in (\ref{eq.fB_5}) and (\ref{eq.tildenexp_2}), we have that
\begin{equation}\label{eq.ININ0aim}
\|\big[\Ic_\Mb(\widetilde{Q}_\Mb)-\widetilde{\Ic}_\Mb(\widetilde{Q}_\Mb)(\tau)\big](x-\tau v,v)\|_{\mathrm{Bol}^{\widetilde{\beta}_\ell}}\leq\varepsilon^{1/2} C^k \tau^{k/2},
\end{equation} for any $\Mb \in \mathfrak M^*_k$ with $k < \Lambda_{\ell}-1$.
\end{proposition}
\begin{proof}
First define $\widetilde {\boldsymbol{\Delta}}_{\pm}^{(1)}:=\widetilde {\boldsymbol{\Delta}}_{\pm}-\widetilde {\boldsymbol{\Delta}}_{\pm}^{(0)}$. Then we can write 
\begin{equation}\label{eq.diffININ0}
\Ic_\Mb(\widetilde{Q}_\Mb)(z)-\widetilde{\Ic}_\Mb(\widetilde{Q}_\Mb)(\tau,z)=\sum_{\nf_0 \in \Mc}\Ic_{\Mb,\nf_0}(\widetilde{Q}_\Mb)(z),
\end{equation}
where $\Ic_{\Mb,\nf_0}(\widetilde{Q}_\Mb)$ is an expression of the following form: 
\begin{equation}\label{eq.diff_z0}
\begin{aligned}
\Ic_{\Mb,\nf_0}(\widetilde{Q}_\Mb)(z)&:=\int_\Dc\int_{(\Sb^{d-1})^{|\Mb|}}\int_{\mathbb{R}^{2d|\Ec_*|}}\prod_{\nf\neq\nf_0}\widetilde{\Dirac}_{\pm}^{(*)}(z^{\nf}_{1},z^{\nf}_{2},z^{\nf}_{3},z^{\nf}_{4},t_\nf)\\
&\phantom{\int_\Dc\int_{(\Sb^{d-1})^{|\Mb|}}\int_{\mathbb{R}^{2d|\Ec_*|}} }\times\widetilde{\Dirac}_{\pm}^{(1)}(z^{\nf_0}_{1},z^{\nf_0}_{2},z^{\nf_0}_{3},z^{\nf_0}_{4},t_{\nf_0}) \prod_{e\in \Ec_{\mathrm{end}}^-}f(({\ell-1})\tau, z_e)\,\mathrm{d}\vz_{\Ec_*}\mathrm{d}\boldsymbol{\omega}_\Mb\mathrm{d}\vt_{\Mb},
\end{aligned}
\end{equation}
where $\widetilde {\boldsymbol{\Delta}}_{\pm}^{(*)}\in \{\widetilde {\boldsymbol{\Delta}}_{\pm}, \widetilde {\boldsymbol{\Delta}}_{\pm}^{(0)}\}$.

Let $\Mb_0=S_{\nf_0}=\{\nf_0\}\cup\Mb_1\cup\Mb_2$, where $\Mb_j=S_{\nf_j}$ and $\nf_j$ are the children of $\nf_0$ (with trivial modifications if any $\nf_i$ is absent), and $\Mb_3=\Mb\backslash \Mb_0$. Note in particular that the parent $\nf_0^{\mathrm{pr}}$ of $\nf_0$ belongs to $\Mb_3$, and the bond $e_{\mathrm{pr}}$ between $\nf_0$ and $\nf_0^{\mathrm{pr}}$ in $\Mb$ becomes a bottom free end at $\nf_0^{\mathrm{pr}}$ in $\Mb_3$, and a top fixed end at $\nf_0$ in $\Mb_0$. We then have $\Mb_i\in\Mf_{k_i}^*$, where $k=k_0+k_3$ and $k_0=k_1+k_2+1$. Define also the same notions $\Dc$ and $(\Ec,\Ec_*)$ etc. for the sub-molecules $\Mb_i$, but with subscript $i$ (the fixed ends will be specified below).

By first integrating in the $\Mb_1$ and $\Mb_2$ variables, then in $\nf_0$ and finally in $\Mb_3$ variables, we get
\begin{equation}\label{eq.diff_z0_1}
\begin{aligned}
\Ic_{\Mb,\nf_0}(\widetilde{Q}_\Mb)(z)&=\int_{\Dc_3}\int_{(\Sb^{d-1})^{|\Mb_3|}}\int_{\mathbb{R}^{2d|\Ec_{*3}|}} \prod_{\nf\in\Mb_3}\widetilde {\boldsymbol{\Delta}}_{\pm}^{(*)}(z^{\nf}_{1},z^{\nf}_{2},z^{\nf}_{3},z^{\nf}_{4},t_\nf)\\
&\phantom{\int_{\mathbb{R}^{2d|\Ec_{*3}|}\times (\Sb^{d-1})^{|\Mb_3|}\times  \Dc_3} }\times \Qc_{\nf_0}(t_{\nf_0^{\mathrm{pr}}},z_{e_{\mathrm{pr}}})\prod_{e\in \Ec_{\mathrm{end}}^{3,-}\backslash \{e_{\mathrm{pr}}\}}f(({\ell-1})\tau, z_e)\,\mathrm{d}\vz_{\Ec_{*3}}\mathrm{d}\boldsymbol{\omega}_{\Mb_3}\mathrm{d}\vt_{\Mb_3},
\end{aligned}
\end{equation} where $\Qc_{\nf_0}$ is defined as (with the fixed top end in $\Mb_0$ being at $\nf_0$)
\begin{equation}\label{eq.diff_z0_2}
\begin{aligned}
\Qc_{\nf_0}(t,z):=&\int_{\Dc^{t}}\int_{(\Sb^{d-1})^{|\Mb_0|}}\int_{\mathbb{R}^{2d|\Ec_{*0}|}} \prod_{\nf\in\Mb_0\backslash\{\nf_0\}}\widetilde {\boldsymbol{\Delta}}_{\pm}^{(*)}(z^{\nf}_{1},z^{\nf}_{2},z^{\nf}_{3},z^{\nf}_{4},t_\nf)\\&\phantom{\int_{\mathbb{R}^{2d|\Ec^0|}\times (\Sb^{d-1})^{|\Mc^0|}\times  \Dc^{0,t}}}\times\widetilde{\boldsymbol{\Delta}}_{\pm}^{(1)}(z^{\nf_0}_{1},z^{\nf_0}_{2},z^{\nf_0}_{3},z^{\nf_0}_{4},t_{\nf_0})\prod_{e\in \Ec_{\mathrm{end}}^{0,-}}f(({\ell-1})\tau, z_e)\,\mathrm{d}\vz_{\Ec_{*0}}\mathrm{d}\boldsymbol{\omega}_{\Mb_0}\mathrm{d}\vt_{\Mb_0}\\
=&\int_0^t\int_{\Sb^{d-1}}\int_{\Rb^{6d}}\widetilde{\boldsymbol{\Delta}}_{\pm}^{(1)}(z,z_2,z_{3},z_{4},s)\cdot
\Ic_{\Mb_1}^*(\widetilde{Q}_{\Mb_1})(s,z_3)\cdot\Ic_{\Mb_2}^*(\widetilde{Q}_{\Mb_2})(s,z_4)
\prod_{i=2}^4\mathrm{d}z_i\mathrm{d}\omega\mathrm{d}s.
\end{aligned}
\end{equation} Here in (\ref{eq.diff_z0_2}), the $\Ic_{\Mb}^*(\widetilde{Q}_{\Mb})$ is defined in the same way as $\widetilde{\Ic}_{\Mb}(\widetilde{Q}_{\Mb})$ in (\ref{eq.tildenexp_2}), but allowing arbitrary choices between $\widetilde{\Dirac}_{\pm}$ and $\widetilde{\Dirac}_{\pm}^{(0)}$ at each factor. If we define $f^{(j)}(z):=\Ic_{\Mb_j}^*(\widetilde{Q}_{\Mb_j})(s,z)$, then we obtain, after integrating in the $x$ variables, that
$$
\begin{aligned}
\Qc_{\nf_0}(t,x,v)=\int_0^t\int_{\Sb^{d-1}}\int_{\Rb^d}\big[f^{(2)}(x+s(v-v_4)+\varepsilon \omega, v_4)-f^{(2)}(x+s(v-v_4), v_4)\big]\\
\times f^{(1)}(x+s(v -v_3), v_3)\cdot\big[(v-v_2)\cdot \omega\big]_+\,\mathrm{d}v_2\mathrm{d}\omega \mathrm{d}s,
\end{aligned}
$$ where $(v_3,v_4)$ equals $(v,v_2)$ if $\nf_0$ is O-atom, and equals $(v',v_2')$ defined by (\ref{eq.boltzmann2}) from $(v,v_2)$ if $\nf$ is C-atom. By mean value theorem, we can rewrite this expression of $\Qc_{\nf_0}(t,x,v)$ as
$$
\begin{aligned}
&\phantom{=}\varepsilon\int_0^1\int_0^t\int_{\Sb^{d-1}}\int_{\Rb^d}f^{(1)}(x+s(v -v_3), v_3) \cdot(\omega\cdot\nabla_x)f^{(2)}(x+s(v-v_4))+\theta\varepsilon\omega,v_4)\cdot\big[(v-v_2)\cdot \omega\big]_+\,\mathrm{d}v_2\mathrm{d}\omega \mathrm{d}s\mathrm{d}\theta\\
&=\varepsilon \int_0^1\int_0^t\int_{\Sb^{d-1}}\int_{\Rb^{6d}}\widetilde{\Dirac}_{\pm}^{(\theta)}(z,z_2,z_3,z_4,s)\cdot\Ic_{\Mb_1}^*(\widetilde{Q}_{\Mb_1})(s,z_3)\cdot(\omega\cdot\nabla_x)\Ic_{\Mb_2}^*(\widetilde{Q}_{\Mb_2})(s,z_4)
\prod_{i=2}^4\mathrm{d}z_i\mathrm{d}\omega\mathrm{d}s\mathrm{d}\theta,
\end{aligned}
$$
where $\widetilde{\Dirac}_{\pm}^{(\theta)}$ is the same as ${\boldsymbol{\Delta}_{\pm}}$ in \eqref{eq.newboldDelta}--\eqref{eq.newboldDelta_2} except that the last $\boldsymbol \delta$ in both expressions is replaced by $\boldsymbol{\delta}(x_1^\nf-x_2^\nf-t_\nf(v_2^\nf -v_1^\nf)+\theta\varepsilon \omega_\nf)$.

Now, given the definition of $\Ic_{\Mb}^*(\widetilde{Q}_{\Mb})(s,z_4)$ similar to $\Ic_{\Mb}(\widetilde{Q}_{\Mb})$ and $\widetilde{\Ic}_{\Mb}(\widetilde{Q}_{\Mb})$, it is clear that its $\nabla_x$ gradient can be written as a sum of at most $|\Ec_{\mathrm{end}}^{2,-}|$ terms, each of which having the same form as $\Ic_{\Mb}^*(\widetilde{Q}_{\Mb})$ itself, but with one of the input factors $f((\ell-1)\tau, z_e)$ replaced by $\nabla_x f((\ell-1)\tau, z_e)$. As a result, we can bound \eqref{eq.diffININ0} by at most $\Lambda_\ell^2\leq |\log\varepsilon|^{C^*}$ expressions of the following form: for some fixed $\nf_0$ and some bottom free end $e_0$ at a descendant atom of $\nf_0$, the expression is given by
\begin{equation*}
\begin{aligned}
&\varepsilon\int_0^1\int_{\Dc}\int_{(\Sb^{d-1})^{|\Mb|}}\int_{\mathbb{R}^{2d|\Ec_*|}}\prod_{\nf\neq\nf_0}\widetilde {\boldsymbol{\Delta}}_{\pm}^{(*)}(z^{\nf}_{1},z^{\nf}_{2},z^{\nf}_{3},z^{\nf}_{4},t_\nf)\cdot\widetilde{\Dirac}_{\pm}^{(\theta)}(z^{\nf_0}_{1},z^{\nf_0}_{2},z^{\nf_0}_{3},z^{\nf_0}_{4},t_{\nf_0})\\
&\phantom{\varepsilon\int_0^1\int_{\Dc}\int_{(\Sb^{d-1})^{|\Mc|}}\int_{\mathbb{R}^{2d|\Ec|}}}\times\prod_{e\in \Ec_{\mathrm{end}}^-\backslash \{e_0\}}|f(({\ell-1})\tau, z_e)|\cdot|\nabla_x f(({\ell-1})\tau, z_{e_0})|\,\mathrm{d}\vz_{\Ec_*}\mathrm{d}\boldsymbol{\omega}_\Mb\mathrm{d}\vt_{\Mb}\mathrm{d}\theta.
\end{aligned}
\end{equation*}
Now the integral in $(\vz_{\Ec_*},\boldsymbol{\omega}_\Mb,\vt_\Mb)$ is uniformly bounded in $\theta\in[0,1]$, following exactly the same proof of Proposition \ref{prop.fa_est_1}, using the $\mathrm{Bol}^\beta$ bound for both $n$ and $\nabla_xn$ in (\ref{eq.boltzmann_decay_4}). This proves Proposition \ref{prop.fa_est_2}.
\end{proof}

Now, by combining (\ref{eq.fa_est_2}), (\ref{eq.fB_5}), (\ref{eq.tildenexp})--(\ref{eq.tildenexp_2}), we conclude that
\begin{equation}\label{eq.fa_est_3}\big\|f^{\Ac,\mathrm{led}}-f(\ell\tau)\big\|_{\mathrm{Bol}^{\beta_{\ell}}}\leq \varepsilon^{\theta_\ell}/5.
\end{equation} Here the difference between $f_1^{\Ac,\mathrm{led}}$ and $F^{(\mathrm{app})}$ is bounded using (\ref{eq.ININ0aim}) and counting binary trees, and the contribution of $F^{(\mathrm{rem})}$ is bounded using (\ref{eq.INBoltzest}) and our choice of $\Lambda_{\ell}$.
\subsection{The sub-leading terms}\label{sec.fa_subleading} We now turn to $f^{\Ac,\mathrm{sub}}$ in (\ref{eq.fAterm2})--(\ref{eq.fAterm3}). This is a sum over molecules $\Mb$. Let $|\Mb|=n$ and $\rho(\Mb)=\rho$, note that either $\rho>0$ or $|\Mb|_p\geq\Lambda_\ell$ by definition of $f^{\Ac,\mathrm{sub}}$ (if $\Mb\in\Tc_{\Lambda_\ell}^{\mathrm{err}}$ then the extra O-atom creates a cycle). The same proof of Proposition \ref{prop.layerrec3} allows us to bound the number of molecules $\Mb$ by $C^n|\log\varepsilon|^{C^*\rho}$, so it suffices to bound each individual term $\Ic_\Mb|Q_\Mb|$. Just as in the proof of Proposition \ref{prop.cumulant_est} (see Section \ref{sec.summary}), it suffices to construct an operation sequence (after turning $\Mb$ into regular, see \textbf{Point 4} below) that produces at most $C^n|\log\varepsilon|^{C^*\rho}$ sub-cases, and prove that each sub-case satisfies (cf. (\ref{eq.proof_int_est_from_alg_step0_*}))
 \begin{equation}\label{eq.s9Merror}
\|\Ic_{\Mb'}|Q'|(x-\tau v,v)\|_{\mathrm{Bol}^{\widetilde{\beta}_\ell}}\leq \tau^{n/9}\cdot\varepsilon^{(\upsilon/3)+(C_{14}^*)^{-1}\cdot\rho}.
\end{equation}

The proof of (\ref{eq.s9Merror}) will follow the same arguments as in Sections \ref{sec.treat_integral}--\ref{sec.maincr}. We will only elaborate on the few places where the proof needs to be adapted to the current setting.

\textbf{Point 1: adjusting Maxwellian weights.} Start with the input factors $g=f^\Ac((\ell-1)\tau)$ in $\Ic_\Mb|Q_\Mb|$. We decompose $g=f((\ell-1)\tau)+g_{\mathrm{rem}}$ as before, where $f((\ell-1)\tau)$ and $g_{\mathrm{rem}}$ satisfy (\ref{eq.boltzmann_decay_4}) and (\ref{eq.approxfA_2}) respectively. Let $\Vc\subseteq\Ec_{\mathrm{end}}^-$ be such that we have $g_{\mathrm{rem}}$ for $e\in \Vc$ and $f((\ell-1)\tau)$ for $e\not\in\Vc$, then we already gain a factor $\varepsilon^{\theta_{\ell-1}\cdot|\Vc|}$ from the inputs in $\Vc$.

Next, we use $\|h\|_{\mathrm{Bol}^\beta}=\|e^{\beta |v|^2}h\|_{\mathrm{Bol}^0}$ to reduce the $\mathrm{Bol}^\beta$ norms to $\mathrm{Bol}^0$ norms. By conservation of energy $\sum_{e\in \Ec^-_{\mathrm{end}}}|v_e|^2=\sum_{e\in \Ec^+_{\mathrm{end}}}$ we have
\[\sum_{e\in \Vc}\beta_{\ell-1}|v_e|^2+\sum_{e\in\Ec_{\mathrm{end}}^-\backslash \Vc}\beta|v_e|^2\geq \widetilde{\beta}_\ell|v|^2+\sum_{e\in \Vc}(\beta/40\Lf)|v_e|^2+\sum_{e\in\Ec_{\mathrm{end}}^-\backslash \Vc}(\beta/2)|v_e|^2+(\beta/40\Lf)|v|^2+\sum_{e\in\Ec_{\mathrm{end}}^+\backslash\{e_{\mathrm{tr}}\}}\widetilde{\beta}_\ell|v_e|^2,\] so arguing as in Part 1 of the proof of Proposition \ref{prop.fa_est_1}, we can reduce $Q_\Mb$ to $\widetilde{Q}$ where the input factors are bounded in $\mathrm{Bol}^0$ norm, and we replace $\Ic_\Mb$ by $\widetilde{\Ic}_\Mb$ with 
\begin{equation}\label{eq.adjustweight}\widetilde{\Ic}_\Mb|\widetilde{Q}|:=\Ic_\Mb\bigg(|\widetilde{Q}|\cdot\varepsilon^{\theta_{\ell-1}\cdot|\Vc|}\cdot e^{-(\beta/40\Lf)|v|^2}\cdot\prod_{e\in\Vc}e^{-(\beta/40\Lf)|v_e|^2}\prod_{e\in\Ec_{\mathrm{end}}^-\backslash \Vc}e^{-(\beta/2)|v_e|^2}\prod_{e\in\Ec_{\mathrm{end}}^+\backslash\{e_{\mathrm{tr}}\}}e^{-\widetilde{\beta}_\ell|v_e|^2}\bigg).\end{equation} Define the post-operation quantity $\widetilde{\Ic}_{\Mb'}|\widetilde{Q}'|$ in the same way.

\textbf{Point 2: decomposition of $\widetilde{Q}'$.} This is the same as in Part 2 of the proof of Proposition \ref{prop.fa_est_1}. We decompose $\widetilde{Q}$ (and consequently $\widetilde{Q}'$ by decomposing the input functions as in Part 2 of the proof of Proposition \ref{prop.cumulant_est}, to reduce to estimating a single $\widetilde{Q}_\alpha'$ under the support condition that $|v_e|\leq|\log\varepsilon|^{C^*}$ and each $x_e$ belonging to a fixed ball of radius $|\log\varepsilon|^{C^*}$. By the same arguments in Part 2 of the proof of Proposition \ref{prop.fa_est_1}, we can also bound the $\mathrm{Bol}^0$ norm of $\widetilde{\Ic}_{\Mb'}|\widetilde{Q}_\alpha'|$ by its $L^\infty$ norm up to a factor $|\log\varepsilon|^{C^*}$.

\textbf{Point 3: applying Propositions \ref{prop.weight} and \ref{prop.volume}.} In view of the weights in (\ref{eq.adjustweight}), we will choose $\gamma=\beta/(10^5\Lf)$ in Proposition \ref{prop.weight}, then the $\exp(\cdots)$ factor on the right hand side of (\ref{eq.weight_est_1}) is acceptable and can be absorbed by the decay factor in (\ref{eq.adjustweight}), so we can apply Proposition \ref{prop.weight} in the same way as in Part 3 of the proof of Proposition \ref{prop.fa_est_1}, except that we now have an additional factor $\gamma^{-n/2}\leq (C\Lf)^{n/2}$ on the right hand side of (\ref{eq.weight_est_1}). Note that $\tau\Lf\leq C$, this additional factor can then be absorbed by the time integral in (\ref{eq.weight_est_2}), which provides a power $\tau^{n}$.

As for Proposition \ref{prop.volume}, let $e$ be any (top or bottom) free end of $\Mb$. Thanks to the decay factor in (\ref{eq.adjustweight}), if we restrict to $|v_e|\sim X_e$, then we will gain a factor $e^{-\beta X_e^2/4}$ \emph{unless} $e\in\Vc\cup\{e_{\mathrm{tr}}\}$, which we treat as exceptional cases. In non-exceptional cases, the gain $e^{-\beta X_e^2/4}$ absorbs the power $X_e^d$ in (\ref{eq.volume_1}) leading to a constant $\leq C$; in exceptional cases we simply apply the rough bound $|v_e|\leq |\log\varepsilon|^{C^*}$. This leads to the loss of $|\log\varepsilon|^{C^*(|\Vc|+1)}$, but can be absorbed by the gain $\varepsilon^{\theta_{\ell-1}\cdot|\Vc|}$ in (\ref{eq.adjustweight}) (if $|\Vc|=0$ then we only lose $|\log\varepsilon|^{C^*}$ which is acceptable in view of the gain in (\ref{eq.s9Merror})).

\begin{figure}[h!]
    \centering
    \includegraphics[width=0.7\linewidth]{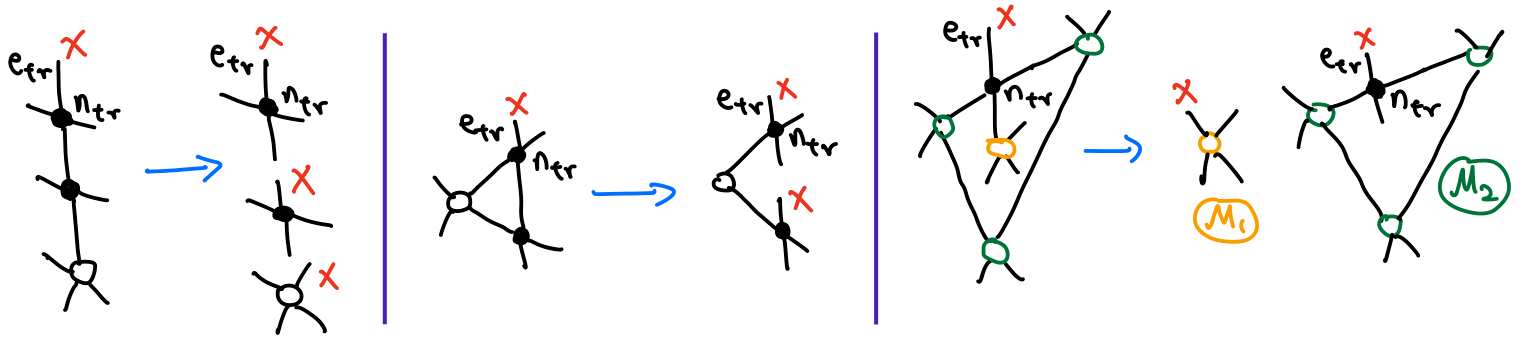}
    \caption{Left: operation in \textbf{Point 4}. Middle: this operation may break a cycle and create two deg 3 atoms. Right: this may also disconnect different clusters, creating one deg 3 atom in each component.}
    \label{fig.finaloper}
\end{figure}

\textbf{Point 4: turning $\Mb$ to regular.} Note that $\Mb$ has a fixed end $e_{\mathrm{tr}}$ at an atom $\nf_{\mathrm{tr}}$. By Definition \ref{def.reg}, this might not be regular if $\nf_{\mathrm{tr}}$is O-atom and $e_{\mathrm{tr}}$ is serial with a bond. In this case, we then turn it into a regular molecule as follows: take the maximal ov-segment containing $e_{\mathrm{tr}}$, then remove every bond in $\sigma$, draw a bottom free end at $\nf_{\mathrm{tr}}$, draw a simple pair (Definition \ref{def.reg}) at each intermediate O-atom in $\sigma$, and draw a top fixed end at the lowest atom of $\sigma$ ($\nf_{\mathrm{tr}}$ is the highest atom), see {\color{blue}Figure \ref{fig.finaloper}}.

This operation is akin to the cutting operation in Definition \ref{def.cutting}, and it is easy to prove (similar to Proposition \ref{prop.cutting}) that we can match the edges before and after the operation (similar to Definition \ref{def.associated_vars_op}), such that the value of $\Ic_\Mb|Q_\Mb|$ is unchanged after the operation. For convenience, we still call the molecule $\Mb$ after the operation, which is now regular.

\textbf{Point 5: the cutting algorithm.} By all the above preparations, and by repeating the same proof as in Sections \ref{sec.treat_integral}--\ref{sec.maincr}, it now suffices to construct an operation sequence that produces at most $C^n|\log\varepsilon|^{C^*\rho}$ sub-cases, and for each sub-case $\Mb'$ contains only elementary molecules, and that \begin{equation}\label{eq.T_final}
\#_{\{33B\}}=\#_{\{4\}}=\#_{\{44\}}=0,\quad\textrm{and }\#_{\mathrm{good}}\geq (10\Gamma)^{-1}\cdot\rho.
\end{equation} Here we may assume $\rho>0$ (otherwise $|\Mb|_p\geq\Lambda_\ell$ and thus $n=|\Mb|\geq \Lambda_\ell-1$ using also $\rho=0$, in which case the power $\tau^{n}$ coming from time integration (\ref{eq.weight_est_2})  provides sufficient gain, in view of our choice of $\Lambda_{\ell}$).

First assume $\Mb\in\Tc_{\Lambda_\ell}^{\mathrm{err}}$, then it is obtained from some $\Mb'\in\Tc_{\Lambda_\ell}$ by creating one more O-atom. If $\rho(\Mb')>0$, we can delete this O-atom and reduce to $\Mb'$ which will be discussed below; if not, then $\rho(\Mb)=1$.
It is easy to see that, after the operation in \textbf{Point 4}, $\Mb$ has no full component, no deg 2 atom and no C-bottom fixed end, and at least one component must contain a cycle or two deg 3 atoms (the two deg 3 atoms case may occur if the operation in \textbf{Point 4} breaks a cycle). We then apply \textbf{DOWN} to $\Mb$, so we get $\#_{\{33B\}}=\#_{\{4\}}=\#_{\{44\}}=0$, and $\#_{\{33A\}}\geq 1$ by Proposition \ref{prop.alg_up} (\ref{it.up_proof_3}). Moreover, the two fixed ends of this \{33A\} molecule cannot both belong to the ov-segment containing $e_{\mathrm{tr}}$ (this follows from Proposition \ref{prop.mol_axiom} (\ref{it.axiom1})), so one of them must come from cutting as free another deg 3 atom. This then leads to at least one \emph{good} molecule, by repeating the proof of Proposition \ref{prop.alg_up} (\ref{it.up_proof_4}).

Now let us assume $\Mb\in\Tc_{\Lambda_\ell}$. Let $\Mb_j$ be the clusters of $\Mb$, with $\Mb_1$ containing the unique fixed end $e_{\mathrm{tr}}$ (before the operation in \textbf{Point 4}). Then, after the operation in \textbf{Point 4}, there is only one deg 3 atom in the component containing $\Mb_1$, and any other component of $\Mb$ also contains exactly one deg 3 atom $\of_j$ on an ov-segment within some $\Mb_j$. For convenience we call such $\Mb_j$ \textbf{lovely}.

We then follow the strategy in \textbf{Choice 2} in the proof of Proposition \ref{prop.case2} (and formally allow clusters to be empty ends as in that proof). Namely, we cut $\Mb_1$ as free and cut it into elementary components using \textbf{DOWN}; for any lovely $\Mb_j$, we also cut $\Mb_j\cup\{\of_j\}$ as free and cut it into elementary components using \textbf{DOWN}. This then creates more lovely $\Mb_j$ and then we repeat this process until all atoms have been cut. Here each $\Mb_j$ (or $\Mb_j\cup\{\of_j\}$) contains exactly one deg 3 atom when it is cut, which is proved using Proposition \ref{prop.mol_axiom} (\ref{it.axiom2}), as in \textbf{Choice 2} in the proof of Proposition \ref{prop.case2}.

As such, by Proposition \ref{prop.alg_up} (\ref{it.up_proof_3}), we get a good component for each $\Mb_j$ with $\rho(\Mb_j)>0$ (i.e. with a cycle). Since $\rho=\sum_j\rho(\Mb_j)$ and $\rho(\Mb_j)\leq\Gamma$ by Proposition \ref{prop.mol_axiom} (\ref{it.axiom3}), we know that the number of $j$ with $\rho(\Mb_j)>0$, and hence the number of good components, is $\geq \rho/\Gamma$. This proves (\ref{eq.T_final}), which then completes the proof of Proposition \ref{prop.est_fa}.

\section{Proof of the error term estimates}\label{sec.error}In this section we prove Proposition \ref{prop.cumulant_error}. The key ingredient here is a cutting algorithm specifically designed for the layer $\Mb_\ell\in\Fc_{\Lambda_\ell}^{\mathrm{trc.err}}$ for the molecule $\Mb\in\Fc_{\vLambda_\ell}^{\mathrm{trc.err}}$. This is stated in Proposition \ref{prop.comb_est_extra} below.
\begin{proposition}\label{prop.comb_est_extra} Let $\Mb$ be a connected molecule with only C-atoms, with $|\Mb|\leq |\log\varepsilon|^{C^*}$ and $\Gamma<\rho(\Mb)\leq 2\Gamma$. Then there exists an operation sequence (formed by cutting and splitting) which produces at most $C^{|\Mb|}|\log\varepsilon|^{C^*}$ sub-cases, such that in each sub-case, $\Mb'$ consists of good, normal, and bad elementary molecules, and the numbers $\#_{\mathrm{good}}$ and $\#_{\mathrm{bad}}$ of good and bad molecules satisfy that
\begin{equation}
\label{eq.extra_case_1} (\upsilon/2)\cdot(\#_{\mathrm{good}})-d\cdot (\#_{\mathrm{bad}})\geq \upsilon/2.
\end{equation}
\end{proposition}
\begin{proof} The proof will occupy Section \ref{sec.last_proof}.
\end{proof}
We now prove Proposition \ref{prop.cumulant_error} assuming Proposition \ref{prop.comb_est_extra}.
\begin{proof}[Proof of Proposition \ref{prop.cumulant_error} assuming Proposition \ref{prop.comb_est_extra}] Start with (\ref{eq.cumulant_formula_err_trunc}). The term $\mathrm{Err}((\ell-1)\tau)$ is already controlled using (\ref{eq.cumulant_formula_err2}) and (\ref{eq.cumulant_est_2}), so we only need to estimate the summation over $\Mb$. Note that only the root particles of $\Mb$ are labeled.

Let $\Mb\in\Fc_{\vLambda_\ell}^{\mathrm{trc.err}}$ and $r(\Mb)=[s]$ be fixed. We will assume $\rho(\Mb_\ell)>\Gamma$ in (\ref{eq.trc_err}); otherwise $|\Mb_\ell|_p>\Lambda_\ell$, in which case we can argue similarly as below, and the powers of $\tau^{\Lambda_\ell/10}$ coming from time integration will provide enough gain, similar to the proof of (\ref{eq.cumulant_est_2}).

Define $\rho$ as in (\ref{eq.def_rho_old}), but only summing $s_{\ell'}$ and $\Rf_{\ell'}$ for $\ell'\leq\ell-1$. Define also $\kappa$ to be the number of root particle lines that have an atom in $\Mb_{\ell'}$ (the other root particle lines occur as empty ends in $\Mb_{\ell}$). By Proposition \ref{prop.layerrec3}, the number of choices for layers $\Mb_{\ell'}\,(\ell'\leq\ell-1)$, initial links and merging top and bottom ends, is bounded by $C^{|\Mb|}|\log\varepsilon|^{C^*\rho}$.

As for $\Mb_{\ell'}$, by Definition \ref{def.set_T_F} (\ref{it.set_F_err_2}), we know it equals $s-\kappa$ empty ends (with fixed labels) plus a graph, which equals a tree plus $\rho(\Mb_\ell)\leq2\Gamma$ edges. Note also that the $\kappa$ root particles have fixed labels and must be identified in the graph. By the same arguments as in Part 2 of the proof of Proposition \ref{prop.layerrec3}, the number of choices of $\Mb_{\ell}$ is bounded by $C^{|\Mb|}\cdot|\log\varepsilon|^{C^*(1+\kappa)}$ (the factor $|\log\varepsilon|^{C^*\kappa}$ comes from identifying the $\kappa$ root particles in the graph).

With the number of choices for $\Mb$ under control, now fix one $\Mb$. By the same proof of Proposition \ref{prop.cumulant_est} (in Section \ref{sec.summary}), we know that Proposition \ref{prop.cumulant_error} holds, provided that there exists an operation sequence on $\Mb$ that produces at most $C^{|\Mb|}|\log\varepsilon|^{C^*\rho}$ sub-cases, and for each sub-case we have only good, normal, and bad elementary molecules in $\Mb'$, and the following analog of (\ref{eq.overall_alg}) in Proposition \ref{prop.comb_est} holds:
\begin{equation}\label{eq.cutting_err_2}(\upsilon/2)\cdot\#_{\mathrm{good}}+(d-1)(s-\#_{\mathrm{bad}})\geq (100d^2C_{13}^*)^{-1}\cdot(1+\rho+\kappa).
\end{equation}

Now we focus on the proof of (\ref{eq.cutting_err_2}). If a particle line $\pb\in r(\Mb_\ell)=[s]$ is an empty end in $\Mb$ and $\pb\not\in r(\Mb_{\ell-1})$, then we remove it from $\Mb$, which does not change the value of (\ref{eq.cutting_err_2}) (an empty end counts as a bad molecule by Definition \ref{def.good_normal}), so below we will assume no such particle line exists. Consider two cases.

\textbf{Case 1}: assume the sole component of $\Mb_\ell$ has no particle line in $r(\Mb_{\ell-1})$. In this case, $\Mb_\ell$ and $\Mb_{<\ell}$ are two disconnected components of $\Mb$, which we can cut separately. Note that $|r(\Mb_{<\ell})|=s-\kappa$. By applying Proposition \ref{prop.comb_est} to $\Mb_{<\ell}$ and Proposition \ref{prop.comb_est_extra} to $\Mb_{\ell}$, we get
\begin{equation}\label{eq.cutting_err_3}(\upsilon/2)\cdot\#_{\mathrm{good}}+(d-1)(s-\#_{\mathrm{bad}})\geq (d-1)\kappa+(\upsilon/2)+(C_{13}^*)^{-1}\cdot\rho,\end{equation} which implies (\ref{eq.cutting_err_2}).

\textbf{Case 2}: assume the sole component of $\Mb_\ell$ has at least one particle line in $r(\Mb_{\ell-1})$. In this case, we can verify that Definition \ref{def.set_T_F} (\ref{it.set_F_l_4}) actually holds for $\ell'=\ell$. Moreover Definition \ref{def.set_T_F} (\ref{it.set_F_l_1}) and (\ref{it.set_F_l_3}) also hold for $\ell'=\ell$ with slightly worse constants (which do not matter in the proof in Sections \ref{sec.treat_integral}--\ref{sec.maincr}), as does Definition \ref{def.set_T_F} (\ref{it.set_F_l_2}) if we artificially add a particle line in the sole component of $\Mb_\ell$ to $r(\Mb)$ (this increases $s$ by $1$, so we should replace $s=|H|$ in (\ref{eq.overall_alg}) by $s+1$). Therefore, by Proposition \ref{prop.comb_est}, we have
\begin{equation}\label{eq.cutting_err_4}(\upsilon/2)\cdot\#_{\mathrm{good}}+(d-1)(s-\#_{\mathrm{bad}})\geq (C_{13}^*)^{-1}\cdot\rho-(d-1),
\end{equation} where the $-(d-1)$ term is due to the change $s\mapsto s+1$.

Alternatively, we can also perform the following cutting sequence: first cut $\Mb_\ell$ as free, then cut it using Proposition \ref{prop.comb_est_extra}, then cut $\Mb_{<\ell}$ using \textbf{DOWN}. In this case, note that after cutting $\Mb_\ell$ as free, the number of empty ends plus full components of $\Mb_{<\ell}$ does not exceed $s-\kappa$ (this is because, by the (PL) argument in Part 1 of the proof of Proposition \ref{prop.case1}, each such component must intersect a particle line in $r(\Mb)=[s]$, and this particle line must not have an atom $\Mb_{\ell'}$ if this component is to be full; same for empty ends). By Proposition \ref{prop.alg_up}, we know that the number of bad molecules produced in cutting $\Mb_{<\ell}$ using \textbf{DOWN}, is at most $s-\kappa$. Therefore for this alternative we have 
\begin{equation}\label{eq.cutting_err_5}(\upsilon/2)\cdot\#_{\mathrm{good}}+(d-1)(s-\#_{\mathrm{bad}})\geq (\upsilon/2)+(d-1)\kappa.
\end{equation} Finally, by optimizing between (\ref{eq.cutting_err_4}) and (\ref{eq.cutting_err_5}) we easily get (\ref{eq.cutting_err_2}). This completes the proof of Proposition \ref{prop.cumulant_error}, assuming Proposition \ref{prop.comb_est_extra}.
\end{proof}
\subsection{The proof of Proposition \ref{prop.comb_est_extra}} \label{sec.last_proof} Now we prove Proposition \ref{prop.comb_est_extra}. First we recall the theorem of Burago-Ferleger-Kononenko \cite{BFK98}, which establishes an upper bound for the number of collisions with fixed number of particles in $\Rb^d$; we restate it below in a more convenient form.
\begin{proposition}\label{prop.upper_bound_col} Suppose $\Mb$ is a full molecule of C-atoms and $A\subseteq \Mb$ is an atom set, such that for any particle line $\pb$ and any atoms $\nf_1\prec \nf_2\prec \nf_3$ on the particle line $\pb$ (in the ordering in Definition \ref{def.molecule_order}) with $\nf_1,\nf_3\in A$, we must have $\nf_2\in A$. Let $q$ be the total number of bottom edges at atoms of $A$ that are not bonds between two atoms of $A$ (i.e. these edges can be bottom ends of $\Mb$ or bonds connecting an atom in $A$ to an atom not in $A$), then we must have $|A|\leq G(q)$, where $G(q):=(2q+1)\cdot\binom{q}{2}\cdot(32q^{3/2})^{q^2}$.
\end{proposition}
\begin{proof} We may identify particles with the corresponding particle lines, see {\color{blue}Table \ref{tab.trans}}. Let $B$ be the set of particle lines $\pb$ that contain at least one atom (collision) in $A$. Each such particle line must contain a bottom edge at an atom of $A$ that is not a bond between two atoms of $A$, which implies that $|B|\leq q$.

For each particle $\pb\in B$, consider the collision $\nf_\pb^-$ immediately before the first collision of the particle $\pb$ in $A$, and the collision $\nf_\pb^+$ immediately after the last collision of the particle $\pb$ in $A$. Note that one or both of $\nf_\pb^\pm$ may be absent, but the total number of such collisions is at most $2q$. We may fix a total time ordering of these collisions, and consider each of the at most $2q+1$ time intervals between two adjacent collisions time $t_{\nf_\pb^\pm}$. For any such time interval $I$, define the set \begin{multline}B'=\big\{\pb\in B:\textrm{$t_{\nf_\pb^-}$ is either absent, or is on the left of $I$, or is the left endpoint of $I$,}\\
\textrm{$t_{\nf_\pb^+}$ is either absent, or is on the right of $I$, or is the right endpoint of $I$}\big\}.\end{multline} We also shrink $I$ a little bit to $\widetilde{I}$ such that $I\backslash \widetilde{I}$ contains no collision (other than the ones at the endpoint of $I$, which do not belong to $A$). It is then easy to see that (i) any collision in $A$ that happens within time interval $\widetilde{I}$ must be between two particles in $B'$, and (ii) any collision that happens within time interval $\widetilde{I}$ and involves an atom in $B'$ must belong to $A$. In particular, within this time interval $\widetilde{I}$, no particles in $B'$ collides with any other particle, so the restriction of the full dynamics to these particles is equal to the partial duynamics with only atoms in $B'$. Let $\lambda_I$ be the number of collisions in $A$ that happen within time interval $\widetilde{I}$, then the collision dynamics for the $|B'|$ particles in $B'$ within interval $\widetilde{I}$ produces at least $\lambda_I$ collisions. By the BFK theorem, we know that\footnote{The extra factor $\binom{q}{2}$ is due to some extra considerations needed for extended dynamics; see Appendix \ref{app.aux}.} $\lambda_I\leq G_0(q)$ for $G_0(q)=\binom{q}{2}\cdot(32q^{3/2})^{q^2}$. Summing in $I$ we then get $|A|\leq (2q+1)G_0(q):=G(q)$, as desired.
\end{proof}
Before getting to the proof of Proposition \ref{prop.comb_est_extra}, we first briefly explain the strategy.

\textbf{Cutting deg 3 and deg 2 atoms.} A naive idea to exploit the recollisions is to keep cutting deg 3 atoms (the first being deg 4), until hitting a recollision and getting a \{33\} molecule. The main difficulty here, also the main difficulty in general, is that we may be forced to cut long sequences of deg 2 atoms without any gain. Let the set of deg 3 atoms being cut be $A$, then the set of subsequent deg 2 atoms being cut is obtained by starting from $A$ and iteratively taking atoms with 2 bonds with previous atoms, which is just $X(A)$ in the notion of Definition \ref{def.trans}.

We note that each atom in $A$ is cut as deg 3, so it has only one bond connecting to previously cut atoms; in particular keep adding atoms in $A$ will not create any new cycle, and we can think of these sets $A$ as close to being trees (i.e. with $\rho(A)$ controlled). On the other hand, adding each atom in $X(A)$ will create one new cycle and increase the $\rho$ value by $1$. Finally, in choosing the set $A$ we can optimize this process by each time cutting the highest (or lowest) deg 3 atom as in \textbf{UP}, which will ensure this set $A$ is \emph{transversal}\footnote{The intuition is as follows. Suppose we apply \textbf{UP} until getting the first deg 2 atom (or \{33\} molecule), and let the set of deg 3 atoms cut be $A$. Then $A$ basically consists of the union of sets of $S_\nf$ for some $\nf$ (which has no bottom in the sense of Lemma \ref{lem.toy2_aux}), and a subset of the next $S_{\nf'}$ which will have no top (as we keep cutting the highest atom in $S_{\nf'}$). Note that no-bottom (resp. no-top) subsets are transversal sets where $A^-=\varnothing$ (resp. $A^+=\varnothing$), so the general notion of transversality then arises naturally as combination of both.} as in Definition \ref{def.trans}.

\textbf{The sets $A_j$ and reduction to toy models.} Now, suppose we keep cutting sequences of deg 3 and deg 2 atoms as above, then we get sets $A_j$ such that $A_{j+1}\supseteq A_j\cup X(A_j)$. Note that forst each sequence of deg 3 and the next sequence of deg 2 atoms, we get a \{33\} molecule at the transition of the two sequences; thus, we may assume there are only $O(1)$ sequences. Note also that $\rho(A_{j+1}=\rho(A_j\cup X(A_j)=\rho(A_j)+|X(A_j)|$ (cf. Proposition \ref{prop.trans} (\ref{it.trans4})), so in order for $A_J$ to reach $\Mb$ and include many ($>\Gamma$) recollisions for some $J=O(1)$, it must occur at some step that for some $A=A_j$, the size of $X(A)$ is significantly larger than $\rho(A)$. Now, by applying \textbf{the BFK theorem}, we know that the number of bonds between $X^+(A)$ and $A$ must also be sufficiently larger than $\rho(A)$ (Proposition \ref{prop.upper_bound_col}).

Effectively, we can now think $A$ as almost a tree, compare it to $\Mb_D$ in the 2-layer model, and compare $X(A)$ (or more precisely $X^+(A)$ where we only take parents in the construction of $X(A)$) to the $\Mb_U$ in the 2-layer model (this set is not clse to a tree, but this will not be used in the algorithms below). The role of UD connections is now played by the bonds between $X^+(A)$ and $A$. We then repeat the same proof for (part of) Proposition \ref{prop.case5}: if there are many 2-connections (corresponding to a subset $X_0^+(A)\subseteq X^+(A)$) then apply \textbf{2CONNDN} as in toy model III; if not, then we cut these 2-connections as in \textbf{Stage 1}, then perform pre-processing as in \textbf{Stege 3}, then get rid of weak degeneracies and finally apply \textbf{MAINUD} as in toy model I plus.

\begin{figure}[h!]
    \centering
    \includegraphics[width=0.34\linewidth]{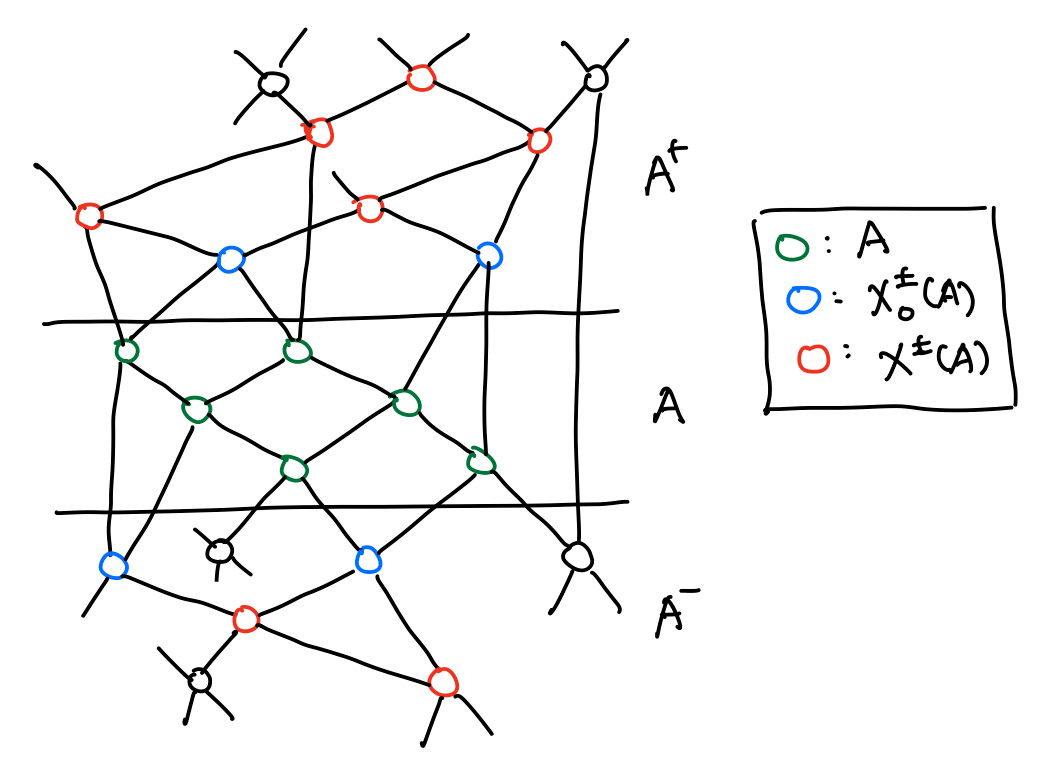}
    \caption{Example of sets $X(A)$ and $X_0(A)$ for a specific connected transversal set $A$. Note that if we cut $A$ as free, then those atoms in $X_0(A)$ (blue) will immediately become deg 2; if we keep cutting deg 2 atoms, then those atoms in $X(A)$ (red) will eventually become deg 2 and be cut.}
    \label{fig.xa set}
\end{figure}

We now define the \textbf{transversal} subsets $A$ of $\Mb$, and prove some properties.
\begin{definition}
\label{def.trans} Let $\Mb$ be a full molecule of C-atoms, and $A\subseteq \Mb$ is an atom set. We say $A$ is transversal, if we can decompose $\Mb\backslash A$ into two disjoint subsets $A^+$ and $A^-$, such that no atom in $A^-$ is parent of any atom in $A\cup A^+$, and no atom in $A$ is parent of any atom in $A^+$.

For any transversal set $A$, define the set $X_0(A)$ such that, an atom $\nf\in X_0(A)$ if and only if $\nf\not\in A$ and $\nf$ has two bonds connected to two atoms in  $A$. Inductively define $X_q(A)$ such that, $\nf\in X_q(A)$ if and only if $\nf\not\in A$, and either $\nf\in X_{q-1}(A)$, or $\nf$ has two bonds connected to two atoms in $X_{q-1}(A)\cup A$. Define $X(A)=\cup_{q\geq 0}X_q(A)$, see {\color{blue}Figure \ref{fig.xa set}}.
\end{definition}
\begin{proposition}
\label{prop.trans}
Let $A$ be a connected transversal subset and $\Mb$ is a connected molecule.
\begin{enumerate}
\item\label{it.trans1} We can always choose $(A^+,A^-)$ in Definition \ref{def.trans} such that each component of $A^+$ has at least one bond connected to $A$.
\item\label{it.trans2} Define the sets $X_q^\pm (A)$ in the same way as $X_q(A)$, but with the sentence ``$\nf$ has two bonds connected to two atoms in $X_{q-1}(A)\cup A$" replaced by ``$\nf$ has two \emph{children} (or \emph{two parents}) in $X_{q-1}^\pm (A)\cup A$''. Then we have $X_q^\pm(A)=X_q(A)\cap A^\pm$. Define $X^\pm(A)=\cup_{q\geq 0}X_q^\pm(A)$, then $X^\pm(A)=X(A)\cap A^\pm$.
\item\label{it.trans3} The set $X(A)\cup A$ is also connected and transversal. Moreover any atom in $\Mb\backslash (X(A)\cup A)$ has at most one bond with atoms in $X(A)\cup A$, so $X_0(X(A)\cup A)=\varnothing$.
\item\label{it.trans4} Recall $\rho(A)$ defined in Definition \ref{def.recollision_number}. Then there exists a connected transversal set $B\supseteq A$, such that $\rho(B)=\rho(A)$, and either $B=\Mb$ or $X_0(B)\neq\varnothing$.
\end{enumerate}
\end{proposition}
\begin{proof} \textbf{Proof of (\ref{it.trans1}).} Suppose $A^+$ has a component $U$ that is not connected to $A$ by a bond, then we replace $A^+$ by $A^+\backslash U$ and $A^-$ by $A^-\cup U$, which still satisfies the requirements for $(A^+,A^-)$ in Definition \ref{def.trans}, but the number of components of $A^+$ decreases by $1$. Repeat until each component of $A^+$ is connected to $A$ by a bond (or $A^+$ becomes empty).

\textbf{Proof of (\ref{it.trans2}).} First $X_q^\pm(A)\subseteq X_q(A)\cap A^\pm$, which easily follows from definition and induction. To prove the opposite direction, first if an atom $\nf\in X_0(A)$ belongs to $A^\pm$, then $\nf$ must have two children (or two parents) in $A$ so $\nf\in X_0^\pm(A)$. We then proceed by induction: suppose $X_{q-1}^\pm(A)=X_{q-1}(A)\cap A^\pm$, now consider $\nf\in (X_q(A)\cap A^+)\backslash X_{q-1}(A)$ which has two adjacent atoms $\nf_1$ and $\nf_2$ in $X_{q-1}(A)\cup A$. If $\nf_1$ belongs to $A^-$, then $\nf$ must be child of $\nf_1$ (it is easy to see that both parents of $\nf_1\in A^-$ must be in $A\cup A^-$) and thus $\nf\in A^-$ which is impossible. This means $\nf_1\in A^+\cup A$ and same for $\nf_2$. Since $\nf_{1,2}\in X_{q-1}(A)\cup A$, we know by induction hypothesis that $\nf_{1,2}\in X_{q-1}^+(A)\cup A$, thus $\nf\in X_q(A)$ by definition. Same is true for $A^-$.

\textbf{Proof of (\ref{it.trans3}).} Connectedness is obvious by definition. To prove transversality, we simply decompose $\Mb\backslash (A\cup X(A))=B^+\cup B^-$, where $B=A^+\backslash X^+(A)$ and $B^-=A^-\backslash X^-(A)$. This then satisfies the requirements (using the fact that any child (or parent) of any atom $\nf\in X^\pm(A)$ must belong to $X^\pm(A)\cup A$). The second statement follows from the definition of $X(A)$ in Definition \ref{def.trans}.

\textbf{Proof of (\ref{it.trans4}).} Suppose $A$ is connected and transversal; we may assume $A\neq\Mb$ and $X_0(A)=\varnothing$ (otherwise choose $B=A$). Since $\Mb$ is connected, there must exist an atom $\nf\in \Mb\backslash A$ that is either a parent or child of an atom in $A$. We may assume it is a parent, and then choose a lowest atom $\nf\in \Mb\backslash A$ among these parents. Since $X_0(A)=\varnothing$, we know that $\nf$ has only one bond with atoms in $A$; let $A_1=A\cup\{\nf\}$, then $A_1$ is connected and $\rho(A_1)=\rho(A)$, see {\color{blue}Figure \ref{fig.trans}}.
\begin{figure}[h!]
    \centering
    \includegraphics[width=0.45\linewidth]{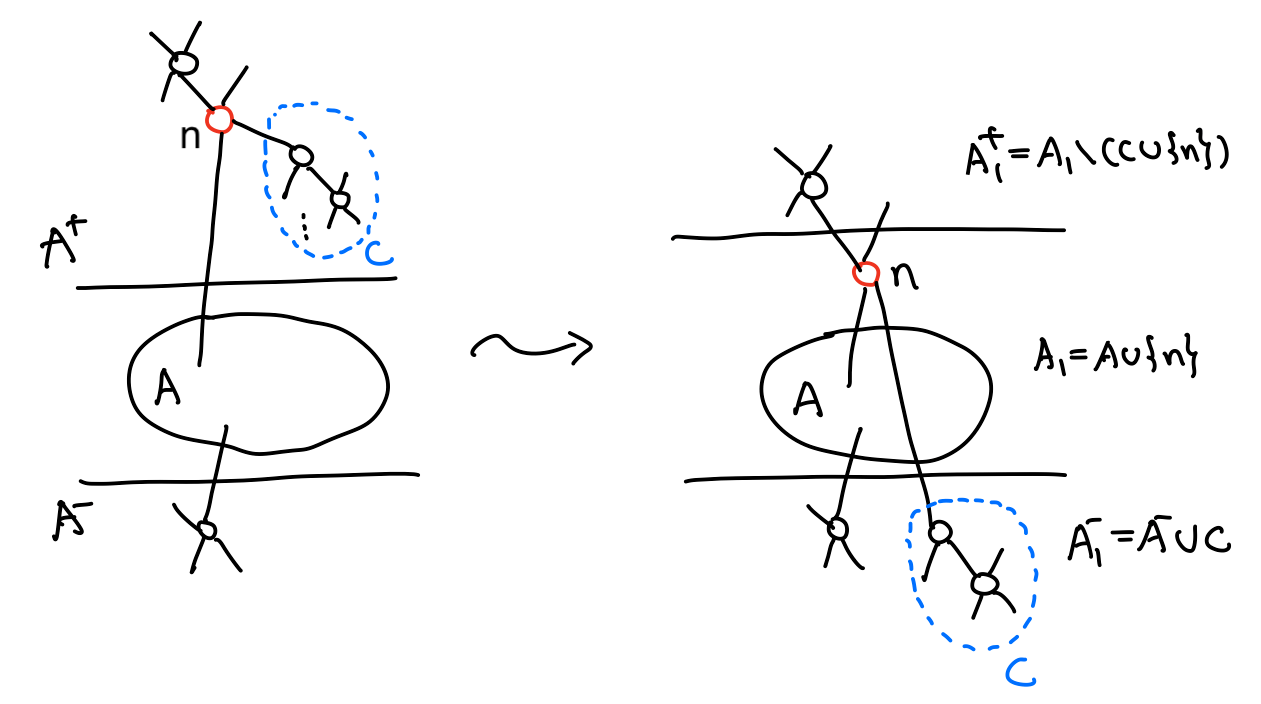}
    \caption{Proof of Proposition \ref{prop.trans} (\ref{it.trans4}): $\nf$ is a lowest parent of atom in $A$, which is natural in terms of the \textbf{UP} or \textbf{DOWN} algorithm; $C$ is the set of descendants of $\nf$ in $A^+$.}
    \label{fig.trans}
\end{figure}

Now we claim that $A_1$ is transversal. To see this, define $C$ to be the set of atoms $\mf\in A^+\backslash\{\nf\}$ that are descendants of $\nf$, and decompose $\Mb\backslash A_1=A_1^+\cup A_1^-$, where $A_1^+=A^+\backslash (C\cup\{\nf\})$ and $A_1^-=A^-\cup C$. It is clear by checking all cases that no atom in $A_1$ can be parent of any atom in $A_1^+$, and no atom in $A^-$ can be parent of any atom in $A_1\cup A_1^+$. If an atom $\mf\in C$ is parent of some atom $\pf\in A_1^+$, then $\pf\in A^+\backslash C$, but $\mf\in C$ and $\pf$ is child of $\mf$, so we also have $\pf\in C$, contradiction. Finally if $\mf\in C$ is parent of $\pf\in A_1=A\cup \{\nf\}$, clearly $\pf\in A$, but then $\mf$ is also parent of an atom in $A$, contradicting the lowest assumption of $\nf$.

Now we know that $A_1$ is also connected transversal; by replacing $A$ by $A_1$ and repeating the above discussion, we eventually will reach some $B$ such that either $B=\Mb$ or $X_0(B)\neq\varnothing$, as desired.
\end{proof}

\begin{figure}[h!]
    \centering
    \includegraphics[width=0.55\linewidth]{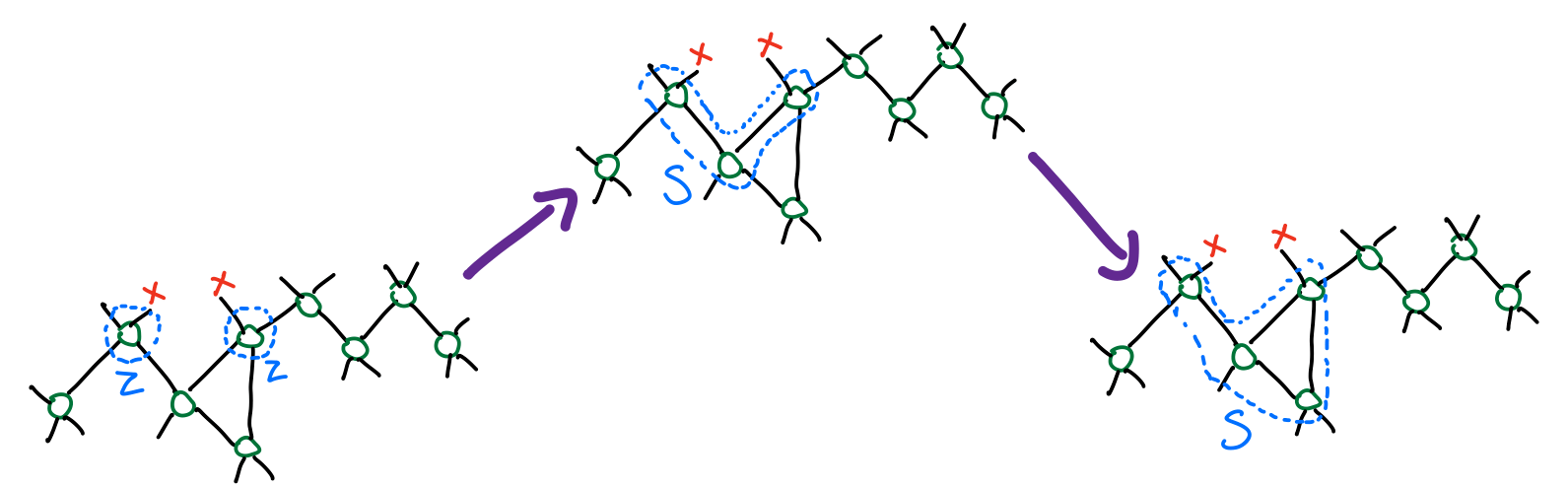}
    \caption{The \textbf{SELECT2} function, which is basically the same as the \textbf{SELECT} function (Definition \ref{def.func_select}, {\color{blue}Figure \ref{fig.select}}). The left arrow: Definition \ref{def.func_select2} (\ref{it.func_select2_2} (same as Definition \ref{def.func_select}), where we add new atoms to connect two components. The right arrow: Definition \ref{def.func_select2} (\ref{it.func_select2_4}), where we add the new atoms that form cycles.}
    \label{fig.select2}
\end{figure}

Now we define the algorithm \textbf{TRANSDN}, which is analogous to the \textbf{2CONNDN} algorithm and is applicable in the case when $|X_0^+(A)|$ (corresponding to the 2-connections in toy model III) is large.
\begin{definition}[Algorithm \textbf{TRANSDN}]
\label{def.alg_transup}
Suppose $\Mb$ is connected full molecule and $A$ is a connected transversal subset. We may choose $(A^+,A^-)$ as in Proposition \ref{prop.trans} (\ref{it.trans1}). Define the following algorithm:
\begin{enumerate}
\item\label{it.alg_transup_1} If $A$ contains any deg $2$ atom $\nf$, then cut it as free, and repeat until there is no deg $2$ atom left.
\item\label{it.alg_transup_2} Choose a \emph{highest} atom $\nf$ in the set of all deg 3 atoms in $A$ (or a highest atom in $A$, if $A$ only contains deg 4 atoms). Let $Z_\nf$ be the set of ancestors of $\nf$ in $A$.
\item \label{it.alg_transup_3} Starting from $\nf$, choose a \emph{lowest} atom $\mf$ in $Z_\nf$ that has not been cut. If $\mf$ has deg 3, and is adjacent to an atom $\pf\in X_0^+(A)$ that also has deg 3, then cut $\{\mf,\pf\}$ as free; otherwise cut $\mf$ as free. Repeat until all atoms in $Z_\nf$ have been cut. Then go to (\ref{it.alg_up_1}).
\item \label{it.alg_transup_4} Repeat (\ref{it.alg_transup_1})--(\ref{it.alg_transup_3}) until all atoms in $A$ have been cut. Then cut the remaining part of $A^+$ as free and cut it into elementary components using \textbf{UP}, then cut $A^-$ into elementary components using \textbf{DOWN}.
\end{enumerate}
We also define the dual algorithm \textbf{TRANSUP} by reversing the notions of parent/child etc. Note that we replace $A^+$ and $X_0^+(A)$ by $A^-$ and $X_0^-(A)$, and also replace $A^+$ by $A^-$ in Proposition \ref{prop.trans} (\ref{it.trans1}).
\end{definition}
\begin{proposition}
\label{prop.alg_transup} For the algorithm \textbf{TRANSDN} (and same for \textbf{TRANSUP}), we have $\#_{\{33B\}}=\#_{\{44\}}=0$ and $\#_{\{4\}}=1$, and $\#_{\{33A\}}\geq |X_0^+(A)|/2-\rho(A)-1$.
\end{proposition}
\begin{proof} It is clear, in the same way as in the proof of Proposition \ref{prop.alg_up}, that during the process, there is no bottom fixed end in $(A\backslash Z_\nf)\cup\{\nf\}$ or in $A^-$, and there is no top fixed end in $Z_\nf\backslash\{\nf\}$ or in $A^+$, so $\#_{\{33B\}}=\#_{\{44\}}=0$. Moreover $\#_{\{4\}}=1$ because no component of the remaining part of $A^+$ after Definition \ref{def.alg_transup} (\ref{it.alg_transup_3}) can be full when it is first cut in Definition \ref{def.alg_transup} (\ref{it.alg_transup_4}) (thanks to Proposition \ref{prop.trans} (\ref{it.trans1})), and no component of $A^-$ can be full when it is first cut in Definition \ref{def.alg_transup} (\ref{it.alg_transup_4}) (thanks to $\Mb$ being connected).

To prove the lower bound for $\#_{\{33A\}}$, note that only one atom in $A$ belongs to a \{4\} molecule (as $A$ is connected). Let the number of atoms in $A$ that belongs to a \{2\} molecule be $q$, then by considering the increment of the same quantity $\sigma:=\#_{\mathrm{bo/fr}}-3|\Mb|$ and applying the same arguments in the proof of Proposition \ref{prop.alg_up} (but restrected to $A$), and noticing that initially all atoms in $A$ have deg 4, we deduce that
\[\#_{\{2\}}+\#_{\{33\}}-\#_{\{4\}}=-\sigma_{\mathrm{init}}=\#_{\mathrm{bond}}-|A|=\rho(A)-1,\] and $\#_{\{4\}}=1$, hence $q\leq\rho(A)$. Then, for each $\nf\in X_0^+(A)$, which is adjacent to two atoms $\nf_1,\nf_2\in A$ by two bonds, assume say $\nf_2$ is cut after $\nf_1$. The total number of such $\nf_2$ is at least $|X_0^+(A)|/2$ as each $\nf_2$ can be obtained from at most two $\nf$. For each such $\nf_2$, the corresponding $\nf$ must have deg 3 when it is cut, so by Definition \ref{def.alg_transup} (\ref{it.alg_transup_3}), it must belong to either \{4\}-, or \{2\}-, or \{33A\} molecule. Using the upper bound for the number of \{4\} and \{2\} molecules, this completes the proof.
\end{proof}
Next we define the \textbf{SELECT2} function, which is an extension of the \textbf{SELECT} function in the case when $A$ may contain cycles (but with $\rho(A)$ bounded). The proof follows similar ideas, except for the extra step (\ref{it.func_select2_4}), which is needed because of the possible cycle in $A$.
\begin{definition}[The function \textbf{SELECT2}]
\label{def.func_select2} Let $A$ be a connected molecule with only C-atoms, no bottom fixed end, and no deg 2 atoms. Let $Z$ be the set of deg 3 atoms in $A$, and let $Y$ be a subset of atoms in $A$ such that $A$ becomes a forest after removing the atoms in $Y$. We define the function $\textbf{SELECT2}=\textbf{SELECT2}(A,Z,Y)$ as follows, see {\color{blue}Figure \ref{fig.select2}}. 
\begin{enumerate}
\item\label{it.func_select2_1} Consider all the components of $Z\cup Y$ in $A$, which is a finite collection of disjoint subsets of $A$.
\item\label{it.func_select2_2} Since $A$ only has C-atoms, the shortest ov-distance between sets (Definition \ref{def.ov_connect}) coincides with the standard definition of distance on a graph. If any two of the subsets in (\ref{it.func_select2_1}), say $U$ and $V$, have the shortest distance (within $A$) which is at most 4, then we choose one shortest path between an atom in $U$ and an atom in $V$, and let the atoms on this path be $\nf_j$. Then replace the two sets $U$ and $V$ by a single set which is $U\cup V\cup \{\nf_j\}$.
\item\label{it.func_select2_4} Repeat (\ref{it.func_select2_2}) until this can no longer be done. Next, if a single subset $U$ contains two atoms (which may be the same) that are connected by a path of length at most 4 with none of the intermediate atoms belonging to $U$ or any other subset, then add the atoms on this path to $U$. If this causes the shortest distance between two subsets to be $\leq 4$, then proceed to (\ref{it.func_select2_2}) and repeat it as above.
\item\label{it.func_select2_5} When no scenario in (\ref{it.func_select2_2}) or (\ref{it.func_select2_4}) occurs, we output $S:=\textbf{SELECT2}(A,Z,Y)$ as the union of all the current sets.
\end{enumerate}
\end{definition}
\begin{proposition}\label{prop.func_select2} The molecule $A$ becomes a proper forest after cutting $S$ as free, and we also have $|S|\leq 10(|Y|+|Z|+\rho(A))$.
\end{proposition}
\begin{proof} The proof is essentially the same as Proposition \ref{prop.func_select}. First $A$ becomes a forest after cutting $S$ as free, because $Y\subseteq S$. Also $A$ becomes proper after cutting $S$ as free, due to the same reason in the proof of Proposition \ref{prop.func_select}; in the language of that proof, the absence of the scenario in Definition \ref{def.func_select2} (\ref{it.func_select2_2}) prevents the case when $\mf$ and $\mf'$ belong to two different components $U$ and $V$, and the absence of the scenario in Definition \ref{def.func_select2} (\ref{it.func_select2_4}) prevents the case when $\mf$ and $\mf'$ belong to the same components $U$.

As for the upper bound for $|S|$, simply note that each single step in Definition \ref{def.func_select2} (\ref{it.func_select2_2})--(\ref{it.func_select2_4}) adds at most $3$ atoms to $S$. Moreover each single step in Definition \ref{def.func_select2} (\ref{it.func_select2_2}) decreases the number of subsets by 1 (as it connectes $U$ and $V$ together) and does not decrease $\rho(S)$ (as $\rho(S)$ is the number of independent cycles, so adding more atoms cannot decrease $\rho(S)$). Each single step in Definition \ref{def.func_select2} (\ref{it.func_select2_4}) does not change the number of subsets, but \emph{adds one more cycle to $S$}; in particular adding this cycle increases the value of $\rho(S)$ by at least 1. Putting together, we know that the value of
\[\gamma:=(\textrm{number of sets in }S) -\rho(S)\] decreases by at least 1 in each step Definition \ref{def.func_select2} (\ref{it.func_select2_2})--(\ref{it.func_select2_4}). Since we always have $\rho(S)\leq \rho(A)$ because $S\subseteq A$ (again it is easy to see that adding one atom or arbitrarily many atoms does not decrease the $\rho$ value), we know that the initial value of $\gamma$ is $\leq |Y|+|Z|$ and the final value of $\gamma$ is $\geq -\rho(A)$, which means that the total number of steps in Definition \ref{def.func_select2} (\ref{it.func_select2_4}) is bounded by $|Y|+|Z|+\rho(A)$, hence the result.
\end{proof}
Now we can define the algorithm \textbf{MAINTRUP}, which is the main algorithm in this section. It has two options, \textbf{Option 1}, which is applied only in the (simple) weakly degenerate case, and essentially trivially cuts most atoms in $A$ as \{3\} molecules, so one can gain from extra restrictions from weak degeneracies. In \textbf{Option 2}, which is the main part, we basically combine \textbf{Stage 1} (cutting 2-connections in $X_0^+(A)$) and \textbf{Stage 3} (pre-processing) in the proof of Proposition \ref{prop.case5} at the end of Section \ref{sec.maincr}, and the \textbf{MAINUD} algorithm in toy model I plus, but wih the roles of $\Mb_U$ and $\Mb_D$ played by $X^+(A)$ and $X$ respectively.
\begin{definition}[Algorithm \textbf{MAINTRUP}] 
\label{def.alg_maintrup}
Suppose $\Mb$ is a connected full molecule and $A$ is a connected transversal subset. We may choose $(A^+,A^-)$ as in Proposition \ref{prop.trans} (\ref{it.trans1}). Let the number of bonds connecting an atom in $X^+(A)$ to an atom in $A$ be $\#_{\mathrm{conn}}^+$. It is easy to prove that there exists a set $Y_0\subseteq A$ of at most $\rho(A)$ atoms, such that $A$ becomes a forest after cutting atoms in $Y_0$ as free. Moreover, by splitting into at most $C^{|A|+\rho(A)}$ sub-cases, we can identify a set of weakly degenerate atoms and atom pairs in $A$ in the sense of Proposition \ref{prop.comb_est_case4}; let $Y_1\subseteq A$ be the atoms involved in these weak degeneracies. We define the following algorithm, which contains two Options that we can choose at the beginning.

In \textbf{Option 1} we do the followings:
\begin{enumerate}
\item Cut $A$ as free, then cut it into elementary components using \textbf{UP}.
\item Then cut $A^+$ as free and cut it into elementary components using \textbf{UP}, then cut $A^-$ into elementary components using \textbf{DOWN}.
\end{enumerate}

In \textbf{Option 2} we do the followings:
\begin{enumerate}
\item\label{it.alg_maintrup_1}  Cut all atoms in $X_0^+(A)$ as free. If any atom in $A$ becomes deg 2, also cut it as free until $A$ has no deg 2 atom.
\item\label{it.alg_maintrup_2} If $A$ remains connected after the above step, let $Z$ be the set of deg 3 atoms in $A$, and $Y$ be those atoms in $Y_0\cup Y_1$ that have not been cut. Define $S=\textbf{SELECT2}(A,Z,Y)$ as in Definition \ref{def.func_select2}, then cut $S$ as free and cut it into elementary components using \textbf{DOWN}. If $A$ is not connected, apply this step to each connected component of $A$.
\item\label{it.alg_maintrup_3} If not all atoms in $X^+(A)$ have been cut, then choose a lowest atom $\nf$ in $X^+(A)$ that has not been cut. If $\nf$ is adjacent to an atom $\pf\in A$ that has deg 3, then cut $\{\nf,\pf\}$ as free; otherwise cut $\nf$ as free. If $A$ becomes non-proper (Definition \ref{def.toy_proper}) then repeat the steps in Definition \ref{def.alg_maincr} (\ref{it.alg_maincr_5})--(\ref{it.alg_maincr_6}) until $A$ becomes proper again.
\item Repeat (\ref{it.alg_maintrup_3}) until all atoms in $X^+(A)$ have been cut. Then repeat the step in Definition \ref{def.alg_maincr} (\ref{it.alg_maincr_9}) until all atoms in $A$ have been cut.
\item\label{it.alg_maintrup_5} Finally, cut (the remaining parts of) $A^+$ as free and cut it into elementary components using \textbf{UP}, and then cut $A^-$ into elementary components using \textbf{DOWN}.
\end{enumerate}

We define the dual algorithm \textbf{MAINTRDN} in the same way (so $\#_{\mathrm{conn}}^+$ is replaced by $\#_{\mathrm{conn}}^-$ etc.).
\end{definition}
\begin{proposition}
\label{prop.alg_maintrup} In \textbf{Option 1} of \textbf{MAINTRUP} in Definition \ref{def.alg_maintrup} (same for \textbf{MAINTRDN}), we have
\begin{equation}\label{eq.alg_maintrup_1}
\#_{\{33B\}}=\#_{\{44\}}=0,\quad \#_{\{4\}}=1;\qquad \#_{\mathrm{good}}\geq\frac{1}{10}\cdot |Y_1|-\rho(A)-1.
\end{equation}
 In \textbf{Option 2} of of \textbf{MAINTRUP} in Definition \ref{def.alg_maintrup}, we have $\#_{\{44\}}=0$ and all \{33B\} molecules are good, and moreover
\begin{equation}\label{eq.alg_maintrup_2}
\begin{aligned}
\#_{\{33A\}}+\#_{\{33B\}}&\geq\frac{1}{10}\cdot\big(\#_{\mathrm{conn}}^+-10^5(|Y_1|+\rho(A)+|X_0^+(A)|)\big),\\
\#_{\{4\}}&\leq |Y_1|+\rho(A)+|X_0^+(A)|.
\end{aligned}
\end{equation}
\end{proposition}
\begin{proof} The case of \textbf{Option 1} is easy; $\#_{\{33B\}}=\#_{\{44\}}=0$ is obvious by definition, and $\#_{\{4\}}=1$ because $A$ is connected and $A^+$ satisfies the assumption in Proposition \ref{prop.trans} (\ref{it.trans1}), in the same way as in the proof of Proposition \ref{prop.alg_transup}. Moreover, each atom in $A$ belongs to a \{4\} or\{3\} or \{2\} or \{33A\} molecule, and the total number of \{2\} and \{33A\} molecules is equal to $\rho(A)$ (see the proof of Proposition \ref{prop.alg_transup}). For any atom $\nf$ that belongs to a  \{3\} molecule, if it is a weakly degenerate atom, or if it is part of a weakly degenerate pair and is cut after the other atom $\nf'$ of the pair, then $\{\nf\}$ must be a good component by Definition \ref{def.good_normal} (\ref{it.good_1}). This gives at least $|Y_1|/10-\rho(A)-1$ good components, as desired.

In the case of \textbf{Option 2}, note that in the whole process, there is no top fixed end in $A^+$ and no bottom fixed end in $A^-$; this is because for the lowest atom $\nf$ in $X^+(A)$ chosen in Definition \ref{def.alg_maintrup} (\ref{it.alg_maintrup_3}), any child of $\nf$ must either belong to $A$ or belong to $X^+(A)$ (and thus will have already been cut). Moreover, any $\nf$ chosen in this step either has deg 2 and no bond connecting to $A$ or has deg 3 and exactly one bond connecting to $A$ (because all atoms in $X_0^+(A)$ have been cut in Definition \ref{def.alg_maintrup} (\ref{it.alg_maintrup_1})). Also any \{33\} molecule $\{\nf,\pf\}$ cut in this way must be \{33A\} molecule, and no full component can be cut (hence no \{4\} molecule created) in Definition \ref{def.alg_maintrup} (\ref{it.alg_maintrup_5}). Therefore, all \{33B\} molecule must have both atoms in $A$, and the only \{4\} molecules created are those created in Definition \ref{def.alg_maintrup} (\ref{it.alg_maintrup_1}) (contributing at most $|X_0^+(A)|$ many \{4\} molecules) and Definition \ref{def.alg_maintrup} (\ref{it.alg_maintrup_2}) (contributing at most $|Y_1|+\rho(A)$ many \{4\} molecules), hence the upper bound $\#_{\{4\}}\leq|Y_1|+\rho(A)+|X_0^+(A)|$.

Now we prove the lower bound on $\#_{\{33\}}$. Note that after Definition \ref{def.alg_maintrup} (\ref{it.alg_maintrup_1}), there is no bottom fixed end in $A$ (and hence none in $S$), so we can cut $S$ into elementary components using \textbf{DOWN} as in Definition \ref{def.alg_maintrup} (\ref{it.alg_maintrup_2}). By definition of $Y_0$ and $Y_1$, and by Proposition \ref{prop.func_select2}, we know that after Definition \ref{def.alg_maintrup} (\ref{it.alg_maintrup_2}) is finished, $A$ will become a forest which is proper, and contains no weakly degenerate atoms or atom pairs. Next we prove an upper bound on $|S|$; let $Z_0$ be the set of atoms in $A$ cut in Definition \ref{def.alg_maintrup} (\ref{it.alg_maintrup_1}), then $|Z|\leq 4(|Z_0|+|X_0^+(A)|)$ (because the deg 3 atoms in $Z$ must be adjacent to one of the atoms in $X_0^+(A)$ or $Z_0$ that has been cut) and $|S|\leq 10(4|Z_0|+4|X_0^+(A)|+|Y_1|+2\rho(A))$ by Proposition \ref{prop.func_select2}. By Definition \ref{def.alg_maintrup} (\ref{it.alg_maintrup_1}), we know that each atom in $Z_0$ must have two parents with each of them in either  $Z_0$ or $X_0^+(A)$. This then leads to $\rho(Z_0)\geq |Z_0|-2|X_0^+(A)|$ (by counting the total number of top bonds at all atoms of $Z_0$, note that each bond is counted only once, and the number of such bonds that do not connect two atoms in $Z_0$ is at most $2|X_0^+(A)|$). Since also $\rho(Z_0)\leq\rho(A)$ (because $Z_0\subseteq A$), we know $|Z_0|\leq \rho(A)+2|X_0^+(A)|$, thus
\[|S|\leq 2\cdot 10^2(\rho(A)+|Y_1|+|X_0^+(A)|).\]

Now, after cutting $X_0^+(A)\cup Z_0\cup S$ as free, the number of bonds connecting an atom in $A$ to an atom in $X^+(A)$ is still at least \[(\#_{\mathrm{conn}}^+)'=\#_{\mathrm{conn}}^+-2(|X_0^+(A)|+|Z_0|+|S|)\geq\#_{\mathrm{conn}}^+-10^4(\rho(A)+|Y_1|+|X_0^+(A)|).\] At this point we can apply the same arguments as in the proof of Proposition \ref{prop.alg_maincr} (especially using Part 2 of the proof and (\ref{eq.alg_maincr_2})) to show that
\begin{equation}\label{eq.alg_maintrup_3}\#_{\{33A\}}+\#_{\{33B\}}\geq\frac{1}{10}\cdot \big(\#_{\mathrm{conn}}^+-10^5(\rho(A)+|Y_1|+|X_0^+(A)|)\big).\end{equation} Here the role of UD connections is played by the $(\#_{\mathrm{conn}}^+)'$ bonds connecting an atom in $A$ to an atom in $X^+(A)$. Note that after cutting $S$ as free, the number of components of $A$ is at most
\[\#_{\mathrm{comp}(A)}\leq 1+4(|X_0^+(A)|+|S|+|Z_0|)\leq 10^5(\rho(A)+|Y_1|+|X_0^+(A)|),\] so among these $(\#_{\mathrm{conn}}^+)'$ bonds (we denote this set by $Q$), there exists a subset $Q'\subseteq Q$ of at least $(\#_{\mathrm{conn}}^+)'-\#_{\mathrm{comp}(A)}$ bonds, such that each bond in $Q'$ is connected to some other bond in $Q$ via $A$ (which matches the definition of UD connections). The same proof for Proposition \ref{prop.alg_maincr} then applies, which leads to (\ref{eq.alg_maintrup_3}). This completes the proof.
\end{proof}
Now, we can finally prove Proposition \ref{prop.comb_est_extra}.
\begin{proof}[Proof of Proposition \ref{prop.comb_est_extra}] First assume $\Mb$ contains at least one strongly degenerate primitive pair $(\mf,\nf)$, then we can cut $\{\mf,\nf\}$ as free and cut the rest of $\Mb$ into elementary components using either \textbf{UP} or \textbf{DOWN} (cf. proof of Proposition \ref{prop.case2}, \textbf{Choice 1}). This creates one good \{44\} molecule and no \{4\} molecule (as $\Mb$ is connected), so $\#_{\mathrm{good}}\geq 1$ and $\#_{\mathrm{bad}}=0$, which implies (\ref{eq.extra_case_1}).

Now assume $\Mb$ contains no strongly degenerate primitive pair, then we can show (by splitting into two sub-cases) that each \{33A\} molecule is either good or can be cut into one \{2\} and one good \{3\}-atom. The proof is similar to Lemma \ref{lem.alg_up_ex}; suppose this \{33A\} molecule is $\{\nf,\nf'\}$ where $\nf$ is parent of $\nf'$, then we may assume this is strongly degenerate and $|t_\nf-t_{\nf'}|\leq \varepsilon^{\upsilon}$ otherwise this is already good by Definition \ref{def.good_normal}. Now our assumption implies that $\{\nf,\nf'\}$ must be non-primitive; in the absence of C-atoms, this implies that there exists another chain of atoms going from $\nf$ to $\nf'$ by iteratively taking children. If any atom $\qf$ in this chain is cut before $\{\nf,\nf'\}$, then the inequality $t_\nf\geq t_{\qf}\geq t_\nf'$ allows us to argue as in the proof of Lemma \ref{lem.alg_up_ex} to get a good molecule; if not, this means that at the time $\{\nf,\nf'\}$ is cut, $\nf$ still has two children (i.e. $\nf'$ and a child on the chain), and $\nf'$ still has two parents (i.e. $\nf$ and a parent on the chain), which means that $\nf$ must have a top fixed end and $\nf'$ must have a bottom fixed end, but this contradicts that $\{\nf,\nf'\}$ is \{33A\} molecule.

Therefore, from now on, we will treat \{33A\} molecules as good. We start by choosing a connected transversal subset $A_1$ of $\Mb$ which is a tree (i.e. $\rho(A_1)=0$) and either $A_1=\Mb$ or $X_0(A_1)\neq\varnothing$; the existence of such $A_1$ follows in the same way as in Proposition \ref{prop.trans} (\ref{it.trans4}) starting from a single atom with no children. If $A_1\neq\Mb$, let $X_1=X(A_1)$, then $A_1\cup X_1$ is connected transversal by Proposition \ref{prop.trans} (\ref{it.trans3}), so we can find $A_2\supseteq A_1\cup X_1$ starting from $A_1\cup X_1$, using Proposition \ref{prop.trans} (\ref{it.trans4}). If $A_2\neq \Mb$, then let $X_2=X(A_2)$ and find $A_3$ by Proposition \ref{prop.trans} (\ref{it.trans4}), then define $X_3=X(A_3)$ and so on, see {\color{blue}Figure \ref{fig.aj set}}.
\begin{figure}[h!]
    \centering
    \includegraphics[width=0.2\linewidth]{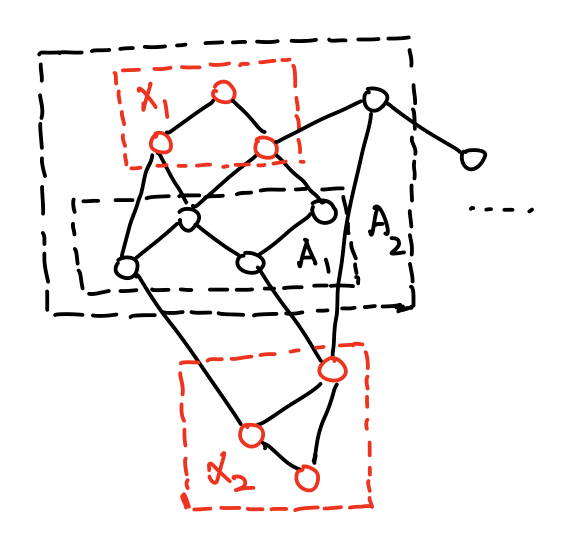}
    \caption{Sets $A_j$ and $X_j$ in the proof of Proposition \ref{prop.comb_est_extra} where $X_j=X(A_j)$ and $A_{j+1}\supseteq A_j\cup X_j$. Note that if we cut $A_1$, then $X_1$, then $A_2\backslash (A_1\cup X_1)$, then $X_2$ etc., then we get a \{33\} molecule when cutting each $X_j$.}
    \label{fig.aj set}
\end{figure}

Recall the function $G(q)=(2q+1)\cdot \binom{q}{2}\cdot(32q^{3/2})^{q^2}$ as in Proposition \ref{prop.upper_bound_col}. Define $D=K_1:=(60d)^{60d}$ and $K_{j+1}=(60d\cdot G((60dK_{j})^{60d}))^{60d}$ for $j\geq 1$. We will assume $\Gamma>(60d)^{60d} K_D$, and consider several cases.

(1) Suppose $A_j$ exists and $X_0(A_j)\neq\varnothing$ for all $j\leq D$. By construction $\rho(A_j)=\rho(A_{j-1}\cup X_{j-1})$, so we can list the elements of $B_j:=A_j\backslash (A_{j-1}\cup X_{j-1})$ as $\nf_1^j,\cdots ,\nf_{|B_j|}^j$ such that $\nf_1^j$ has only one bond with $A_{j-1}\cup X_{j-1}$ and $\nf_i^j$ has only one bond with $A_{j-1}\cup X_{j-1}\cup\{\nf_1^j,\cdots,\nf_{i-1}^j\}$. Similarly we can list the atoms in $B_1:=A_1$. Choose also an element $\mf_j\in X_0(A_j)$ for each $j$, then this $\mf_j$ is adjacent to two atoms in $A_j$, with at least one being in $B_j$. We then perform the following cutting sequence:

Starting from $j=1$, for each $j$, cut the atoms $(\nf_i^j:i\leq |B_j|)$ as free in the increasing order of $i$. However, if $\nf_i^j$ is adjacent to $\mf_{j}$ and $\mf_{j}$ has deg 3 when we cut $\nf_i^j$, then we cut $\{\nf_i^j,\mf_{j}\}$ instead of $\nf_i^j$. After all atoms in $B_j$ have been cut, the atom $\mf_{j}$ will also have been cut; we then cut the remaining atoms in $X_j$ as free (starting from the lowest ones; each will have deg 2 when we cut it), and proceed with $B_{j+1}=A_{j+1}\backslash (A_j\cup X_j)$ and so on.

In this way, it is clear that each $\mf_j$ must belong to a \{33\} molecule, which also has to be \{33A\} molecule because $\mf_j$ will have two top (or two bottom) free ends when it is cut, assuming $\mf_j\in X_0^+(A_j)$ (or $\mf_j\in X_0^-(A_j)$). This produces at least $D$ many \{33A\} molecules (which are treated as good), while $\#_{\{4\}}=1$, which implies (\ref{eq.extra_case_1}).

(2) Suppose $A_j=\Mb$ for some $j<D$, then $\rho(A_j)\geq\Gamma$ for this $j$. Therefore we may choose the smallest $j$ such that $\rho(A_{j+1})\geq K_{j+1}$. Let $A=A_j$, then we have
\[\rho(A)<K_j,\quad K_{j+1}\leq \rho(A_{j+1})=\rho(A\cup X(A))\leq \rho(A)+10|X(A)|,\] so $|X(A)|\geq (20)^{-1}\cdot K_{j+1}$. We may assume $|X^+(A)|\geq (40)^{-1}\cdot K_{j+1}$; consider the set $X^+(A)$, which obviously satisfies the assumptions in Proposition \ref{prop.upper_bound_col}, so by Proposition \ref{prop.upper_bound_col} we conclude that $\#_{\mathrm{conn}}^+>(60d)^{60d}\cdot K_j$ in Definition \ref{def.alg_maintrup}. We then
\begin{itemize} 
\item Run algorithm \textbf{TRANSUP} and apply Proposition \ref{prop.alg_transup} if $|X_0^+(A)|\geq (30d)^{30d}\cdot K_j$;
\item Run algorithm \textbf{MAINTRUP}, \textbf{Option 1} and apply Proposition \ref{prop.alg_maintrup} if $|Y_1|\geq (30d)^{30d}\cdot K_j$;
\item Run algorithm \textbf{MAINTRUP}, \textbf{Option 2} and apply Proposition \ref{prop.alg_maintrup}, if $|X_0^+(A)|\leq (30d)^{30d}\cdot K_j$ and $|Y_1|\leq (30d)^{30d}\cdot K_j$.
\end{itemize} In any case it is easy to verify that (\ref{eq.extra_case_1}) is true, so the proof is complete.
\end{proof}

\appendix
\section{Some auxiliary results}\label{app.aux} In this appendix, we first prove some  basic properties of the extended (E-) and truncated (T- or $(\Lambda,\Gamma)$-) dynamics, in Proposition \ref{prop.dynamics_E_T_property}.
\begin{proof}[Proof of Proposition \ref{prop.dynamics_E_T_property}] We need to prove the following statements (1)--(5). Note that the flow property for $\Sc_N^E(t)$ in (\ref{eq.semigroup}) follows from the same property for $\Hc_N^E(t)$ and the definition (\ref{eq.defS^M}) of $\Sc_N^E(t)$, since\[\sum_{\vz^0:\Hc_N^E(t+s)\vz^0=\vz}f(\vz^0)=\sum_{\vz^1:\Hc_N^E(t)\vz^1=\vz}\,\sum_{\vz^0:\Hc_N^E(s)\vz^0=\vz^1}f(\vz^0);\]
 also the equality $\Sc_N^E(t)f=\Sc_N(t)f$ in Proposition \ref{prop.dynamics_E_T_property} (\ref{it.dynamics_E_T_property_2}) follows from the previous statement concerning $\Hc_N^E(t)$ and $\Hc_N(t)$, in the same way as in the proof of (\ref{eq.S_trunc_equal_to_S}) in Proposition \ref{prop.dynamics_comparison}.
\begin{enumerate}
\item For E- and T-dynamics, the initial configurations in which two collisions happen at the same time, or an infinite number of collisions happen in a finite time interval, or some collision happens at a fixed time, form a Lebesgue zero set.
\item Apart from Lebesgue zero sets, the mapping $\Hc_N^E(t)$ defined by the E-dynamics (Definition \ref{def.E_dynamics}, Proposition \ref{prop.dynamics_E_T_property} (\ref{it.dynamics_E_T_property_2})) satisfies the flow property in \eqref{eq.semigroup}.
\item If the initial configuration $\vz^0\in \Dc_N$, then the E-dynamics and the O-dynamics coincide.
\item In the $(\Lambda,\Gamma)$-dynamics, the size and recollision number of any cluster $\Mb_j$ of the topological reduction $\Mb$ always satisfy that $|\Mb_j|\leq \Lambda$ and  $\rho(\Mb_j)\leq\Gamma$.
\item The equality (\ref{eq.L1pres}) holds for both E- and T-dynamics.
\item For O-, E- and T-dynamics, the total number of collisions has an absolute upper bound that depends only on $d$ and $N$.
\end{enumerate}

\textbf{Proof of (1).} This follows from the same proof as in \cite{Alexander}.

\textbf{Proof of (2)--(3).} Notice that by Definition \ref{def.E_dynamics}, the motion of the particles starting from any time $t$ depends only on the configuration of the system \emph{at time} $t$ but \emph{not on the collision history before time $t$}. The flow property in (2) then follows. Moreover, by Definition \ref{def.E_dynamics}, we know that if the distance between any two particles $i$ and $j$ ever becomes $\geq\varepsilon$ along the E-dynamics, then their distance will remain $\geq\varepsilon$ for the rest of time; therefore, if the distance between any two particles is $>\varepsilon$ initially, then all these distances will remain $\geq\varepsilon$ forever, so the E-dynamics coincide with the O-dynamics, which proves (3).

\textbf{Proof of (4).} This follows from the definition of $(\Lambda,\Gamma)$-truncated dynamics (Definition \ref{def.T_dynamics}): it is easy to see that any collision that is allowed cannot cause $|\Mb_j|>\Lambda$ or $\rho(\Mb_j)>\Gamma$, due to the definition of the condition $\mathtt{Pro}$ in (\ref{eq.trunc_dym}).

\textbf{Proof of (5).} We prove it for $\Sc_N^E$, and the $\Sc_N^{\Lambda,\Gamma}$ case is the same. By listing all the collisions and time ordering them, we can divide $\Rb^{2dN}$ into a countable union of disjoint open sets $(O_\alpha)_{\alpha\in\Nb}$, where on each $O_\alpha$ the E-dynamics contain only a fixed set of collisions between fixed particles with a fixed time ordering. The mapping $\Hc_N^E(t)$, when restricted to $O_\alpha$, is injective and volume preserving; we denote it by $\Hc^t$. Then
\begin{equation}
\int_{\Rb^{2dN}} (\Sc_N^{E}(t)f)(\vz)\,\mathrm{d}\vz=\int_{\Rb^{2dN}} \bigg(\sum_{\vz^0:\Hc_N^{E}(t)(\vz^0)=\vz}f(\vz^0)\bigg)\,\mathrm{d}\vz=\sum_{\alpha}\int_{\Rb^{2dN}} \bigg(\sum_{\vz^0\in O_\alpha:\Hc^t(\vz^0)=\vz}f(\vz^0)\bigg)\,\mathrm{d}\vz.
\end{equation}
The above sum in $\vz^0$ contains at most one element, and is nonzero if and only if $\vz\in\Hc^t(O_\alpha)$. By the above, the mapping $\Hc^t$ is  a volume preserving bijection between $O_\alpha$ and $\Hc^t(O_\alpha)$. Therefore,
\begin{equation}\label{eq.eq.proof_volume_pres_5}
\int \Sc_N^{E}(t)f(\vz)\,\mathrm{d}\vz=\sum_{\alpha}\int_{\Hc^t(O_\alpha)}f((\Hc^{t})^{-1}\vz)\,\mathrm{d}\vz=\sum_\alpha\int_{O_\alpha}f(\vz^0)\,\mathrm{d}\vz^0=\int_{\Rb^{2dN}} f(\vz^0)\,\mathrm{d}\vz^0.
\end{equation}

\textbf{Proof of (6).} For the O-dynamics this follows from \cite{BFK98}; in fact the upper bound is explicit $(32N^{3/2})^{N^2}$. For the $(\Lambda,\Gamma)$-truncated dynamics this follows directly from (5). Consider now the E-dynamics. We know that in the E-dynamics, for each pair of particles, their distance remains $<\varepsilon$ until some time when their distance becomes $\varepsilon$, and their distance remains $\geq\varepsilon$ afterwards. There are at most $\binom{q}{2}$ such transition times, and in any time interval $I=(t',t'')$ between each two consecutive transition times $t'$ and $t''$, the set of particle pairs with distance $<\varepsilon$, say $K$, remains fixed.

Note that, see \cite{BFK98}, the O-dynamics can be viewed as the billiard motion on the billiard table
\[\Dc:=\bigcap_{1\leq i<j\leq N}\Dc_{ij};\quad \Dc_{ij}=\big\{\vz\in \Rb^{2dN}:|x_i-x_j|\geq\varepsilon\big\}.\] The proof in \cite{BFK98} shows that the number of boundary reflections of this billiard ball is at most $(32N^{3/2})^{N^2}$, by establishing some geometric properties of $\Dc$. 

Now the same result in \cite{BFK98} still holds, if one replaces $\Dc$ by $\Dc':=\bigcap_{(i,j)\in K^c}\Dc_{ij}$; indeed $\Dc'$ is just the intersection of a sub-collection of the full collection of balls defining $\Dc$, and it can be checked that the same proof in \cite{BFK98} carries over. This then gives an upper bound on the number of collisions \emph{in the E-dynamics within any interval $I$}, because this dynamics exactly corresponds to the billiard motion on the billiard table $\Dc'$ (the billiard is also restricted to the sets $\Dc_{ij}^c$ for $(i,j)\in K$, but the boundaries of these sets are not involved in any reflections, as the billiard ball never reach the boundary of any of these sets within time interval $I$). The desired result for the E-dynamics then follows with the absolute upper bound $G_0(d,N)=\binom{N}{2}\cdot (32N^{3/2})^{N^2}$. 
\end{proof}
Finally, we state and prove a result concerning the preservation of the norm $\mathrm{Bol}^\beta$ used in Theorem \ref{th.main}.
\begin{proposition}
\label{prop.pres_of_decay} Recall the norm $\mathrm{Bol}^\beta$ defined in (\ref{eq.boltzmann_decay_2}). Suppose $n=n(t,x,v)$ is a positive solution to the Boltzmann equation (\ref{eq.boltzmann}), and let $\beta=\beta(t)$ be a fixed, strictly positive, strictly decreasing function with $\beta(0)\geq \beta>0$ and $\beta'(t)\leq -c<0$ for some constant $c$. Assume $n$ satisfies that
\begin{equation}\label{bound1}\sup_{t\in[0,T]}\big\|e^{\beta(t)|v|^2}n(t)\big\|_{L^\infty}\leq A<\infty,
\end{equation} and that
\begin{equation}\label{bound2}\|n(0)\|_{\mathrm{Bol}^{\beta(0)}}+\|\nabla_xn(0)\|_{\mathrm{Bol}^{\beta(0)}}\leq B<\infty,
\end{equation} then for any $t\in[0,T]$ we have
\begin{equation}\label{bound3}\|n(t)\|_{\mathrm{Bol}^{\beta(t)}}+\|\nabla_xn(t)\|_{\mathrm{Bol}^{\beta(t)}}\lesssim B,
\end{equation} where the constants in $\lesssim$ may depend on $(A,\beta,c,T)$.
\end{proposition}
\begin{proof}
First we prove (\ref{bound3}) for $n(t)$. Define $\|n(t)\|_{\mathrm{Bol}^{\beta(t)}}:=M_0(t)$, by (\ref{eq.boltzmann}) we have
\begin{equation}\label{bol2}
\begin{aligned}n(t,x,v)&=n(0,x-tv,v)+\int_0^t\int_{\Rb^d}\int_{\Sb^{d-1}}\big((v-v_1)\cdot\omega\big)_+\cdot\big\{n(s,x-(t-s)v,v')n(s,x-(t-s)v,v_1')\\&-n(s,x-(t-s)v,v)n(s,x-(t-s)v,v_1)\big\}\,\mathrm{d}\omega\mathrm{d}v_1\mathrm{d}s.
\end{aligned}
\end{equation} We focus on the Duhamel term (denote it by $n^{(1)}$). Using that $|v|^2+|v_1|^2=|v'|^2+|v_1'|^2$, we get that
\begin{equation}\label{bol3}
\begin{aligned}e^{\beta(t)|v|^2}n^{(1)}&=\int_0^t\int_{\Rb^d}e^{-\beta(t)|v_1|^2}\int_{\Sb^{d-1}}\big((v-v_1)\cdot\omega\big)_+\cdot e^{(\beta(t)-\beta(s))(|v|^2+|v_1|^2)}\\
&\times\big\{e^{\beta(s)|v'|^2}n(s,x-(t-s)v,v')\cdot e^{\beta(s)|v_1'|^2}n(s,x-(t-s)v,v_1')\\&\quad- e^{\beta(s)|v|^2}n(s,x-(t-s)v,v)\cdot e^{\beta(s)|v_1|^2}n(s,x-(t-s)v,v_1)\big\}\,\mathrm{d}\omega\mathrm{d}v_1\mathrm{d}s.
\end{aligned}
\end{equation} Note that $\beta(t)\geq\beta$ and $\beta(t)-\beta(s)\leq -c(t-s)$, it is easy to see that
\[\big((v-v_1)\cdot\omega\big)_+\cdot e^{(\beta(t)-\beta(s))(|v|^2+|v_1|^2)/2}\lesssim_{\beta,c}(t-s)^{-1/2}.\] Let $K\in 2^\Nb$ be such that $(t-s)\max(|v|,|v_1|)\sim(t-s)\max(|v'|,|v_1'|)\sim K$ (with $\sim$ replaced by $\lesssim $ if $K=1$), then from (\ref{bol3}) we have
\begin{equation}\label{bol4}
\begin{aligned}e^{\beta(t)|v|^2}n^{(1)}&\lesssim_{\beta,c,T}\sum_K e^{-\theta K^2}\int_0^t(t-s)^{-1/2}\int_{\Rb^d}e^{-\beta|v_1|^2}\mathbbm{1}_{(t-s)\max(|v|,|v_1|)\lesssim K}\\
&\times\int_{\Sb^{d-1}}\big\{e^{\beta(s)|v'|^2}n(s,x-(t-s)v,v')\cdot e^{\beta(s)|v_1'|^2}n(s,x-(t-s)v,v_1')\\&\quad+ e^{\beta(s)|v|^2}n(s,x-(t-s)v,v)\cdot e^{\beta(s)|v_1|^2}n(s,x-(t-s)v,v_1)\big\}\,\mathrm{d}\omega\mathrm{d}v_1\mathrm{d}s,
\end{aligned}
\end{equation} where $\theta\gtrsim_{\beta,c T}1$ is a positive constant. In (\ref{bol4}), if $k\in\Zb^d$ is fixed and $|x-k|\leq 1$, then for $y=x-(t-s)v$ we have $|y-k|\lesssim K$, therefore we have
\begin{equation}\label{bol5}e^{\beta(s)|v'|^2}n(s,x-(t-s)v,v')\leq\sum_{k'\in\Zb^d,|k'-k|\lesssim K}G_{k'}(s),\end{equation} for some $G_{k'}(s)$ satisfying $\sum_{k'}G_{k'}(s)\lesssim\|n(s)\|_{\mathrm{Bol}^{\beta(s)}}\leq M_0(s)$. The other three terms satisfy the same estimate; moreover, by (\ref{bound1}) we also have $e^{\beta(s)|v'|^2}n(s,x-(t-s)v,v')\leq A$ and the same for the other three terms. In summary, for every $|x-k|\leq 1$ and $v\in\Rb^d$, we have
\begin{equation}
\begin{aligned}e^{\beta(t)|v|^2}n^{(1)}(x,v)&\lesssim_{\beta,c,T}\sum_K e^{-\theta K^2}\int_0^t(t-s)^{-1/2}\int_{\Rb^d} e^{-\beta|v_1|^2}\int_{\Sb^{d-1}}A\cdot\sum_{|k'-k|\leq K} G_{k'}(s)\,\mathrm{d}\omega\mathrm{d}v_1\mathrm{d}s.
\end{aligned}
\end{equation} Upon taking supremum in $(x,v)$ and summing in $k$ and $K$, and doing the same for the transport term $n(0,x-tv,v)$, we obtain that
\[M_0(t)\lesssim_{A,\beta,c,T}B+\int_0^t (t-s)^{-1/2}M_0(s)\,\mathrm{d}s.\] By Gronwall, this implies $M_0(t)\lesssim_{A,\beta,c,T}B$ for any $t\in[0,T]$. Finally, for the proof of the $\nabla_x$ derivative, just notice that, by (\ref{eq.boltzmann}) we have $(\partial_t+v\cdot\nabla_x)(\nabla_xn)=\Qc(n,\nabla_xn)+\Qc(\nabla_xn,n)$, so the same result holds by repeating the proof above.
\end{proof}

\section{A dictionary between molecules and physical trajectories}\label{app.trans}

\begin{table}[H]
\centering
\begin{tabular}{!{\vrule width 1pt} c | c !{\vrule width 1pt} c |c !{\vrule width 1pt}}
\noalign{\hrule height 1pt}
 \textbf{Molecule terms} &\textbf{Physical picture} & \textbf{Molecule terms} & \textbf{Physical picture}\\
\noalign{\hrule height 1pt}
C/O atoms & Collisions/overlaps & Particle line & \makecell{Full trajectory of one particle\\(will be identified with particles)}\\
\hline
Layers $\ell'$ & Time intervals $[(\ell'-1)\tau,\ell'\tau]$ & Serial & \makecell{In the same particle line\\ (trajectory of the same particle)}\\
\hline
Edges & Linear pieces (of trajectories) & \makecell{Start/finish\\ layers in $[\underline{\ell}:\ell]$} & \makecell{First/last layer where the\\ particle occurs in the dynamics}\\
\hline
Bonds/ends &\makecell{Linear pieces between two\\ consecutive collisions or overlaps/\\first and last linear pieces} & \makecell{Finish layer $\ell+1$} & \makecell{(Root) particles occurring in the \\time $t=\ell\tau$ correlation/cumulant}\\
\hline
Empty ends & \makecell{Trajectories with\\ no collisions and overlaps} & Lifespan & \makecell{Time layers where the particle \\is present in $\Mb$ (excluding $\ell+1$)}\\
\hline
Parent/child & \makecell{Next/previous collision or \\overlap of an involved particle} & Initial links & \makecell{Initial closeness of particles\\ governed by initial cumulants} \\
\hline

Top/bottom edges & \makecell{Post/pre collisional configurations} & Ov-segments & \makecell{Linear pieces of trajectories\\ with possible overlaps in between}\\
\noalign{\hrule height 1pt}
\end{tabular}
\caption{Correspondence between molecule terms and physical picture. The notations are as in Definition \ref{def.molecule}. There is some special convention about finish layer being $\ell+1$ (recall $\ell$ is the highest layer).}
\label{tab.trans}
\end{table}
\section{Terms and definitions}\label{app.glossary}
\begin{table}[H]
\centering
\begin{tabular}{!{\vrule width 1pt} c | c !{\vrule width 1pt} c |c !{\vrule width 1pt}}
\noalign{\hrule height 1pt}
 \textbf{Notions} &\textbf{Location} & \textbf{Notions} & \textbf{Location}\\
\noalign{\hrule height 1pt}
O dynamics & Def \ref{def.hard_sphere} & Sets $\Ec_*$, $\Ec_{\mathrm{end}}$ & Def. \ref{def.sets_molecule}\\
\hline
Flow map $\Hc_N$, operator $\Sc_N$ & Def. \ref{def.hard_sphere} & Admissible trajectory, pre-col. config. & Def. \ref{def.admissible}\\
\hline
Density $W_N$, correlation $f_s$ & Def. \ref{def.grand_canon} & Topological reduction & Def. \ref{def.top_reduction}\\
\hline
The $\mathrm{Bol}^\beta$ norm & Eq. (\ref{eq.boltzmann_decay_2}) & Cluster, cluster graph, value $\rho(A)$ & Def. \ref{def.cluster}, \ref{def.recollision_number}\\
\hline
Constants $\Lf$, $\tau$ & Def. \ref{def.notation} (\ref{it.time_layer}) & Modified dynamics (associated $\Hc$ and $\Sc$) & Def. \ref{def.modified_dynamics}\\
\hline
Constants $C^*$, $C_j^*$ & Def. \ref{def.notation} (\ref{it.defC*})  & Prescribed dynamics & Def. \ref{def.molecule_truncated_dynamics}\\
\hline
Constants $A_\ell$, $\Lambda_\ell$, $\beta_\ell$, $\theta_\ell$, $\upsilon$ & Def. \ref{def.notation} (\ref{it.defnot3})  & Sets $\Fc_{\vLambda_\ell}$, $\Fc_{\vLambda_\ell}^{\mathrm{err}}$, $\Tc_{\Lambda_\ell}$, $\Tc_{\Lambda_\ell}^{\mathrm{err}}$, $\Fc_{\vLambda_\ell}^{\mathrm{trc.err}}$ & Def. \ref{def.set_T_F}\\
\hline
Molecule, ov-seg., particle lines & Def. \ref{def.molecule}  & $\mathbbm{1}_\Mb$, $\mathbbm{1}_\Lc^\varepsilon$, $(\Sc\circ\mathbbm{1})_\Mb$, $|(\Sc\circ\mathbbm{1})_\Mb|$, $\Ic\Nc_\Mb$, $|\Ic\Nc_\Mb|$ & Def. \ref{def.associated_op_nonlocal}, \ref{def.associated_int}\\
\hline
Descendant, ancestor & Def. \ref{def.molecule_order}  & $\widetilde{f}_s$, $f_s^{\mathrm{err}}$, $f^\Ac$, $E_H$, $\mathrm{Err}$, $\mathrm{Err}^j$ & Prop. \ref{prop.cumulant_formula}\\
\noalign{\hrule height 1pt}
\end{tabular}
\caption{Terms and definitions (Sections \ref{sec.intro}--\ref{sec.formula_cumulant}).}
\end{table}

\begin{table}[H]
\centering
\begin{tabular}{!{\vrule width 1pt} c | c !{\vrule width 1pt} c |c !{\vrule width 1pt}}
\noalign{\hrule height 1pt}
 \textbf{Notions} &\textbf{Location} & \textbf{Notions} & \textbf{Location}\\
\noalign{\hrule height 1pt}
E, T dynamics  & Def. \ref{def.E_dynamics}--\ref{def.T_dynamics} & Values $s_{\ell'}$, $\Rf_{\ell'}$, $\rho$ & Def. \ref{def.parameter_rho_old}\\
\hline
Sets $\Dc_N^{\Lambda,\Gamma}$, $\Ec_N^{\Lambda,\Gamma}$ &Def. \ref{def.truncated_domain} & $\vz_\Ec$, $\vt_\Mb$, $\Dirac_\nf$, $\Dc$, operator $\Ic_\Mb$ & Def. \ref{def.associated_op}\\
\hline
Correlation $\widetilde{f}_s$, $f_s^{\mathrm{err}}$ &Def. \ref{def.ftrunc}  & $Q_\Mb$, $\Ic_\Mb(Q_\Mb)$ & Prop. \ref{prop.local_int}\\
\hline
Partial dynamics  & Def. \ref{def.partial_dynamics}& Deleting, degree, sub-molecule&Def. \ref{def.delete}\\
\hline
$\mathbbm{1}_\Mb^{\Lambda,\Gamma}$, $\mathbbm{1}_{\Mb\not\sim\Mb'}$, $\mathbbm{1}_{(\Mb,e)\not\sim(\Mb',e')}$  & Def. \ref{def.indicator_cluster} & Cutting&Def. \ref{def.cutting}\\
\hline
$\mathbbm{1}_\Mb^{\Lambda,\Gamma}$ (O-atoms), $(\Sc\circ\mathbbm{1})_{\Mb}^{\mathrm{Pen}}$, $\Fs_\Lambda$, $\Fs_\Lambda^{\mathrm{err}}$ & Def. \ref{def.indicator_general} & Regular  & Def. \ref{def.reg}\\
\hline
Molecule sets $\Fc_{\Lambda_\ell}$, $\Fc_{\Lambda_\ell}^{\mathrm{err}}$  & Def. \ref{def.set_F_single_layer} & Ov-adjacent, ov-components etc. & Def. \ref{def.ov_connect}\\
\hline
Quantities $\widetilde{f}_\Mb$, $\widetilde{f}_{\Mb,\mathrm{upp}}$  & Prop. \ref{prop.molecule_representation} & Splitting, sub-cases, restrictions $\Sc_j$ & Def. \ref{def.splitting}\\
\hline
Connectivity $r(\Mb)\searrow H$  & Def. \ref{def.r_conn_H} & Elementary molecules & Def. \ref{def.elementary}\\
\hline
Quantities $\Ic\Nc_{\Mb,H}^{\mathrm{Pen}}$, $|\Ic\Nc_{\Mb,H}|$ & Eq. (\ref{eq.associated_int_notsim_single}) & Matching variables under cutting & Def. \ref{def.associated_vars_op}\\
\hline
Initial links $\Lc$, $\mathbbm{1}_\Lc^{\varepsilon}$, $\upsilon=3^{-d-1}$ & Prop. \ref{prop.initial_cumulant} & Cutting order $\prec_{\mathrm{cut}}$ & Def. \ref{def.ordered_disjoint}\\
\noalign{\hrule height 1pt}
\end{tabular}
\caption{Terms and definitions (Sections \ref{sec.truncation_large_molecule}--\ref{sec.cutting}).}
\end{table}

\begin{table}[H]
\centering
\begin{tabular}{!{\vrule width 1pt} c | c !{\vrule width 1pt} c |c !{\vrule width 1pt}}
\noalign{\hrule height 1pt}
 \textbf{Notions} &\textbf{Location} & \textbf{Notions} & \textbf{Location}\\
\noalign{\hrule height 1pt}
Good, normal, bad molecules& Def. \ref{def.good_normal}&Algorithm \textbf{MAINUD} (toy)&Def. \ref{def.toy1+_alg}\\
\hline
$\#_{\mathrm{good}}$, $\#_{\mathrm{bad}}$ etc.&Prop. \ref{prop.comb_est} &Toy models II \& III&Def. \ref{def.toy2}, \ref{def.toy3}\\
\hline
Quantity $\Bf$&Prop. \ref{prop.weight}&Algorithm \textbf{3COMPUP} (toy)&Def. \ref{def.toy2_alg}\\
\hline
Sets $\Sf$, $\Fc$, $\Gc$&Prop. \ref{prop.volume}&Algorithm \textbf{2CONNUP} (toy)&Def. \ref{def.toy3_alg}\\
\hline
New initial links $\widetilde{\Lc}$, $\mathbbm{1}_{\widetilde{\Lc}}^{\varepsilon+}$&Prop. \ref{prop.init_link}&Layer selection (toy)& Def. \ref{def.layer_select_toy}\\
\hline
Interval layers $\Mb_{[\ell_1:\ell_2]}$ etc.&Def. \ref{def.layer_interval}&Algorithm \textbf{UP}, (MONO) property&Def. \ref{def.alg_up}, Prop. \ref{prop.alg_up}\\
\hline
2-layer toy model&Def. \ref{def.simplified}&(PL) argument (in proof)&Prop. \ref{prop.case1}\\
\hline
Toy models I \& I+&Def. \ref{def.toy}, \ref{def.toy1+}& Strong degenerate, primitive&Def. \ref{def.strdeg}\\
\hline
Proper&Def. \ref{def.toy_proper}&\textbf{Choices 1--3}&Prop. \ref{prop.case2}\\
\hline
Algorithm \textbf{UP} (toy)& Def. \ref{def.up_toy}&Layer refinement, $\zeta=(\ell',k)$, $\Mb_\zeta^T$&Def. \ref{def.layer_refine}\\
\noalign{\hrule height 1pt}
\end{tabular}
\caption{Terms and definitions (Sections \ref{sec.treat_integral}--\ref{sec.layer_refine}).}
\end{table}

\begin{table}[H]
\centering
\begin{tabular}{!{\vrule width 1pt} c | c !{\vrule width 1pt} c |c !{\vrule width 1pt}}
\noalign{\hrule height 1pt}
 \textbf{Notions} &\textbf{Location} & \textbf{Notions} & \textbf{Location}\\
\noalign{\hrule height 1pt}
$\ell_1$, $\zeta_U=(\ell_1,\gamma)$, $\zeta_D$, $\zeta_0$&Def. \ref{def.layer_select}&Ov-proper&Def. \ref{def.ov_proper}\\
\hline
Set $S$, $(\rho',\rho'', \rho''')$&Def. \ref{def.layer_select}&Function \textbf{SELECT}&Def. \ref{def.func_select}\\
\hline
Sets $Y$, $\Mb_{\zeta_U}^T$, $\Mb_{(\zeta_D:\zeta_U)}^T$, $\Mb_{\zeta_D}^T$&Prop. \ref{prop.layer_select}&\textbf{MAINUD}, $X_0$, $X_1$&Def. \ref{def.alg_maincr}\\
\hline
Notions $\Vb$, $G$&Def. \ref{def.layer_cutting}&Ov-components $\#_{\mathrm{ovcp}}$&Prop. \ref{prop.alg_maincr}\\
\hline
\textbf{Options 1--2}, $\Mb_U$, $\Mb_D$, $\Mb_{UD}$&Prop. \ref{prop.layer_cutting}&\textbf{Stages 1--6}&Sec. \ref{sec.finish}\\
\hline
UD connections, $\#_{\mathrm{UD}}$&Prop. \ref{prop.layer_cutting}&$\overline{\Mb}_{\mathrm{2conn}}^2$, $\overline{Z}_0$, $\#_{\mathrm{ovcp}(\overline{Z}_0)}$, $\#_{\mathrm{ovcp}(X_1)}$&Sec. \ref{sec.finish}\\
\hline
Weak degenerate&Def. \ref{def.weadeg}&$\Mf_n$, $\Mf_n^*$, $\pb_{\mathrm{tr}}$, $\Mb^{\mathrm{twi}}$&Def. \ref{def.fa_notion}, Prop. \ref{prop.cancel}\\
\hline
\textbf{2CONNUP}, $\Mb_{\mathrm{2conn}}$, $\Mb_{\mathrm{2conn}}^{1,2}$&Def. \ref{def.alg_2connup}&Transversal set&Def. \ref{def.trans}\\
\hline
Full components $\#_{\mathrm{fucp}}$&Prop. \ref{prop.alg_2connup}&\textbf{TRANSDN}, \textbf{SELECT2}&Def. \ref{def.alg_transup}, \ref{def.func_select2}\\
\hline
\textbf{3COMPUP}, $\#_{\mathrm{3ovcp}}$&Def. \ref{def.3comp_alg}&\textbf{MAINTRUP}&Def. \ref{def.alg_maintrup}\\
\noalign{\hrule height 1pt}
\end{tabular}
\caption{Terms and definitions (Sections \ref{sec.layer_select}--\ref{sec.error}).}
\end{table}


\begin{thebibliography}{99}
\bibitem{Alexander} R.K. Alexander. The infinite hard sphere system. \emph{Ph.D. Thesis}, Department of Mathematics, University of California at Berkeley (1975).
\bibitem{BPS13} A. V. Boblylev, M. Pulvirenti, C. Saffirio. From Particle Systems to the Landau Equation: A Consistency Result. Commun. Math. Phys. 319, 683--702 (2013). 

\bibitem{BGS16} T. Bodineau, I. Gallagher, L. Saint-Raymond. The Brownian motion as the limit of a deterministic system of hard spheres. \emph{Inventiones mathematicae}, 1-61 (2016).

\bibitem{BGS17} T. Bodineau, I. Gallagher, L. Saint-Raymond. From hard sphere dynamics to the Stokes-Fourier equations: an $L^2$ analysis of the Boltzmann-Grad limit. \emph{Annals of PDE}, 3 (2017).

\bibitem{BGS18} T. Bodineau, I. Gallagher, L. Saint-Raymond. Derivation of an Ornstein-Uhlenbeck process for a massive particle in a rarefied gas of particles. \emph{Ann. IHP} 19, Issue 6, 1647--1709 (2018)

\bibitem{BGSS18} T. Bodineau, I. Gallagher, L. Saint-Raymond and S. Simonella. One-sided convergence in the Boltzmann-Grad limit. \emph{Ann. Fac. Sci. Toulouse Math.} (6)27 (2018), no.5, 985--1022.
\bibitem{BGSS22_2} T. Bodineau, I. Gallagher, L. Saint-Raymond and S. Simonella. Cluster expansion for a dilute hard sphere gas dynamics. \emph{J. Math. Phys.} 63 (2022), no. 7, Paper No. 073301, 26 pp.
\bibitem{BGSS20} T. Bodineau, I. Gallagher, L. Saint-Raymond and S. Simonella. Statistical dynamics of a hard sphere gas: fluctuating Boltzmann equation and large deviations. \emph{Ann. of Math. (2)} 198 (2023), no. 3, 1047--1201.
\bibitem{BGSS20_2} T. Bodineau, I. Gallagher, L. Saint-Raymond and S. Simonella. Long-time correlations for a hard sphere gas at equilibrium. \emph{Comm. Pure Appl. Math.} 76 (2023), no. 12, 3852--3911.
\bibitem{BGSS23} T. Bodineau, I. Gallagher, L. Saint-Raymond and S. Simonella. Dynamics of dilute gases: a statistical approach. \emph{EMS Press, Berlin}, 2023, 750--795.
\bibitem{BGSS22} T. Bodineau, I. Gallagher, L. Saint-Raymond and S. Simonella. Long-time derivation at equilibrium of the fluctuating Boltzmann equation. \emph{Ann. Probab. 52} (2024), no. 1, 217--295.
\bibitem{Bol72} L. Boltzmann, Weitere Studien uber das Warme gleichgenicht unfer Gasmolakular. Sitzungs-berichte der Akademie der Wissenschaften 66 (1872), 275--370. Translation: Further studies on the thermal equilibrium of gas molecules, in Kinetic Theory 2, 88--174, Ed. S.G. Brush, Pergamon, Oxford (1966).
\bibitem{Bol02} L. Boltzmann, \emph{Lecons sur la th\'{e}orie des gaz}, Gauthier-Villars (Paris, 1902-1905). R\'{e}-\'{e}dition Jacques Gabay, 1987.
\bibitem{BGHS19} T. Buckmaster, P. Germain, Z. Hani and J. Shatah. Onset of the wave turbulence description of the long-time behavior of the nonlinear Schr\"{o}dinger equation. \emph{Invent. Math.} 225 (2021), 787--855. 
\bibitem{BFK98} D. Burago, S. Ferleger and A. Kononenko. Uniform Estimates on the Number of Collisions in Semi-Dispersing Billiards. \emph{Ann. of Math.} 147 (1998), no. 3, 695--708.
\bibitem{CG19} C. Collot and P. Germain. On the derivation of the homogeneous kinetic wave equation. arXiv:1912.10368.
\bibitem{CG20} C. Collot and P. Germain. Derivation of the homogeneous kinetic wave equation: longer time scales. arXiv:2007.03508.
\bibitem{DEL89} A. de Masi, R. Esposito and J. L. Lebowitz. Incompressible Navier-stokes and Euler limits of the Boltzmann equation. \emph{Comm. Pure Appl. Math.} 42 (1989), no.8, 1189--1214.
\bibitem{DH19} Y. Deng and Z. Hani. On the derivation of the wave kinetic equation for NLS. \emph{Forum of Math. Pi.} 9 (2021), e6.
\bibitem{DH21} Y. Deng and Z. Hani. Full derivation of the wave kinetic equation. \emph{Invent. Math.} 233 (2023), no. 2, 543--724.
\bibitem{DH21_2} Y. Deng and Z. Hani. Propagation of chaos and higher order statistics in wave kinetic theory. \emph{J. Eur. Math. Soc. (JEMS)}, to appear.
\bibitem{DH23} Y. Deng and Z. Hani. Derivation of the wave kinetic equation: full range of scaling laws. arXiv:2301.07063.
\bibitem{DH23_2} Y. Deng and Z. Hani. Long time justification of wave turbulence theory. arXiv:2311.10082.
\bibitem{DHM25} Y. Deng, Z. Hani and X. Ma. Hilbert's sixth problem: derivation of fluid equations via Boltzmann's kinetic theory. arXiv: 2503.01800.
\bibitem{ESY08} L. Erd\"{o}s, M. Salmhofer and H-T. Yau. Quantum diffusion of the random Schr\"{o}dinger evolution in the scaling limit. \emph{Acta Math.} 200 (2008), no. 2, 211--277.
\bibitem{EY00} L. Erd\"{o}s and H-T. Yau. Linear Boltzmann equation as the weak coupling limit of a random Schr\"{o}dinger equation. \emph{Comm. Pure Appl. Math.} 53 (2000), 667--735.
\bibitem{GST14} I. Gallagher, L. Saint-Raymond and B. Texier. From Newton to Boltzmann: hard spheres and short-range potentials. \emph{Zurich Advanced Lectures in Mathematics Series}, vol. 18. EMS, Lewes (2014).
\bibitem{GT20} I. Gallagher and I. Tristani. On the convergence of smooth solutions from Boltzmann to Navier-Stokes. \emph{Ann. Henri Lebesgue} 3 (2020), 561--614.
\bibitem{GS04} F. Golse and L. Saint-Raymond. The Navier-Stokes limit of the Boltzmann equation for bounded collision kernels. \emph{Invent. Math.} 155 (2004), 81--161.
\bibitem{GS23} N. Guillen and L. Silvestre. The Landau equation does not blow up. arXiv:2311.09420, 2023.
\bibitem{Has62} K. Hasselmann. On the nonlinear energy transfer in a gravity wave spectrum. Part 1. \emph{J. Fluid Mech.} 12 (1962), 481--500.
\bibitem{Has63} K. Hasselmann. On the nonlinear energy transfer in a gravity wave spectrum. Part 2. \emph{J. Fluid Mech.} 15 (1963), 273--281.
\bibitem{H01} D. Hilbert. Mathematical Problems. \emph{Bulletin of the American Mathematical Society} 8 (1901), 437--479.
\bibitem{IP86} Reinhard Illner and Mario Pulvirenti. Global validity of the Boltzmann equation for a two-dimensional rare gas in vacuum. \emph{Comm. Math. Phys.} 105 (1986), no. 2, 189--203.
\bibitem{IP89} R. Illner and M. Pulvirenti. Global validity of the Boltzmann equation for two- and three-dimensional rare gas in vacuum: Erratum and improved result. \emph{Comm. Math. Phys.} 121 (1989), no. 1, 143--146.
\bibitem{ISV24} C. Imbert, L. Silvestre and C. Villani. On the monotonicity of the Fisher information for the Boltzmann equation. arXiv:2409.01183, 2024.
\bibitem{Nor28} L. Nordheim. On the kinetic method in the new statistics and application in the electron theory of conductivity. \emph{Proc. Roy. Soc. London Ser. A} 119 (1928), 689--698.
\bibitem{Kac56} Mark Kac. Foundations of Kinetic Theory. \emph{Proceedings of The third Berkeley symposium on mathematical statistics and probability} 3 (1956), 171--197.
\bibitem{King75} F. King, BBGKY hierarchy for positive potentials, Ph.D. dissertation, Dept. Mathematics, Univ. California, Berkeley, 1975.
\bibitem{Lan75} O.E. Lanford, Time evolution of large classical systems, \emph{Lect. Notes in Physics} 38, J. Moser ed., 1--111, Springer Verlag (1975).
\bibitem{LS11} J. Lukkarinen and H. Spohn. Weakly nonlinear Schr\"{o}dinger equation with random initial data. \emph{Invent. Math.} 183 (2011), 79--188.


\bibitem{Ma} X. Ma. Almost sharp wave kinetic theory of multidimensional KdV type equations with $d\geq 3$. arXiv:2204.06148, 2022. 

\bibitem{MS11} J. Marklof and A.  Strombergsson. The Boltzmann-Grad limit of the periodic Lorentz gas. \emph{Annals of Mathematics} 174 (2011), 225--298.

\bibitem{MT12} K. Matthies and F. Theil. A Semigroup Approach to the Justification of Kinetic Theory. \emph{SIAM Journal on Mathematical Analysis}, Vol. 44, Iss. 6 (2012).

\bibitem{MM13} S. Mischler and C. Mouhot. Kac's program in kinetic theory. \emph{Invent. Math. 193}, no. 1, 1--147 (2013).
\bibitem{Naz11} S. Nazarenko. Wave turbulence. \emph{Lecture Notes in Physics} Vol. 825 (2011), Springer, Heidelberg.
\bibitem{OVY93} S. Olla, S.R.S. Varadhan, and H.-T. Yau. Hydrodynamical limit for a Hamiltonian system with weak noise. \emph{Comm. Math. Phys. 155}, no. 3, 523--560 (1993).
\bibitem{Pei29} R. E. Peierls. Zur kinetischen Theorie der W\"{a}rmeleitung in Kristallen. \emph{Annalen Physik} 3 (1929), 1055--1101.
\bibitem{PS17} M. Pulvirenti and S. Simonella. The Boltzmann-Grad limit of a hard sphere system: analysis of the correlation error. \emph{Inventiones Math.} 207 (2017), 1135--1127.
\bibitem{PS21} M. Pulvirenti, S. Simonella. A Brief Introduction to the Scaling Limits and Effective Equations in Kinetic Theory. \emph{In: Albi, G., Merino-Aceituno, S., Nota, A., Zanella, M. (eds) Trails in Kinetic Theory. SEMA SIMAI Springer Series}, vol 25. Springer, Cham 2021.
\bibitem{QY98} J. Quastel and H.T. Yau. Lattice Gases, Large Deviations, and the Incompressible Navier-Stokes Equations. \emph{Ann. Math.} 148 (1998), no. 1, 51--108.
\bibitem{Sai09} L. Saint-Raymond. Hydrodynamic limits of the Boltzmann equation. Lecture Notes in Mathematics. Springer-Verlag Berlin Heidelberg 2009.
\bibitem{SpohnBook} H. Spohn, Large scale dynamics of interacting particles, \emph{Texts and Monographs in Physics}, vol. 174, Springer, 1991, xi+342 pages.  
\bibitem{Spo77} H. Spohn. Derivation of the transport equation for electrons moving through random impurities. \emph{J. Statist. Phys.}
17 (1977), no.6, 385--412.
\bibitem{UU33} E. Uehling and G. Uhlenbeck. Transport Phenomena in Einstein-Bose and Fermi-Dirac Gases. \emph{Phys. Rev.} 43 (1933), 552--561.
\bibitem{UT06} S. Ukai and T. Yang, Mathematical theory of the Boltzmann equation, \emph{Lecture Notes Series, vol. 8}, Liu Bie Ju Center for Mathematical Sciences, City University of
Hong-Kong, 2006.
\bibitem{Villani} C. Villani, A review of mathematical topics in collisional kinetic theory. \emph{Handbook of Math. Fluid Dyn.}, 1 (2002), pp. 71--305.
\bibitem{Zak65} V. E. Zakharov. Weak turbulence in media with decay spectrum. \emph{Zh. Priklad. Tech. Fiz.} 4 (1965), 5--39.
\end{thebibliography}
\end{document}